\documentclass[letterpaper, 11pt]{article}

\usepackage[T1]{fontenc}
\usepackage[english]{babel}
\usepackage{tikz}
\usepackage{hyperref}
\usepackage{amssymb}
\usepackage{amsthm}
\usepackage{mathtools}
\usepackage{array}
\usepackage[english]{cleveref}
\usepackage{bm}
\usepackage{ifthen}
\usepackage[affil-it]{authblk}
\usepackage{thm-restate}
\usepackage[margin=1in]{geometry}
\usepackage{longtable}
\usepackage{ifthen}
\usepackage{caption}
\usepackage{pdfpages}
\usepackage{multirow}
\usepackage{arydshln}
\usepackage{stmaryrd}
\usepackage{version}
\usepackage{multirow}
\usepackage{float}
\usepackage{etoolbox}
\usepackage{footmisc}
\usepackage{marvosym}
\usepackage{csquotes}
\usepackage{enumitem}

\newboolean{arxiv}
\setboolean{arxiv}{true}

\numberwithin{equation}{section}
\numberwithin{table}{section}
\numberwithin{figure}{section}

\Crefname{equation}{Property}{Properties}

\newcommand{\bit}{\mathfrak{b}}
\newcommand{\nsb}{\nu}
\newcommand{\nub}{\mu}

\newcommand{\bitset}{\mathcal{B}}

\newcommand{\flipset}[3]{\mathrm{Flip}(#1,#2\ifthenelse{\equal{#3}{}}{)}{,#3)}}
\newcommand{\unflipset}[3]{\mathrm{Unflip}(#1,#2\ifthenelse{\equal{#3}{}}{)}{,#3)}}
\newcommand{\flips}[3]{\mathrm{fl}(#1,#2\ifthenelse{\equal{#3}{}}{)}{,#3)}}
\newcommand{\flipmax}[3]{\mathrm{mfn}(#1,#2\ifthenelse{\equal{#3}{}}{)}{,#3)}}
\newcommand{\unflipmax}[3]{\mathrm{mufn}(#1,#2\ifthenelse{\equal{#3}{}}{)}{,#3)}}
\newcommand{\lastflip}[3]{\mathrm{lfn}(#1,#2\ifthenelse{\equal{#3}{}}{)}{,#3)}}
\newcommand{\lastunflip}[3]{\mathrm{lufn}(#1,#2\ifthenelse{\equal{#3}{}}{)}{,#3)}}
\newcommand{\e}{\varepsilon}
\newcommand{\scalar}[1]{\left\langle#1\right\rangle}
\newcommand{\minsig}[1]{m^{\sigma}_{#1}}
\newcommand{\minnegsig}[1]{m^{\neg\sigma}_{#1}}

\newcommand{\minsige}[1]{m^{\sigmae}_{#1}}
\newcommand{\minnegsige}[1]{m^{\neg\sigmae}_{#1}}
\newcommand{\Pref}[1]{Property~(\ref{property: #1})}

\newcommand{\maxocc}{\mathfrak{m}}

\newcommand{\indbit}{\beta}

\newcommand{\occrec}{\phi}
\newcommand{\incorrect}[1]{\mathfrak{I}^{#1}}
\newcommand{\floor}[1]{\left\lfloor#1\right\rfloor}
\newcommand{\ceil}[1]{\left\lceil#1\right\rceil}
\newcommand{\abs}[1]{\left|#1\right|}

\newcommand{\sigmabar}{{\bar \sigma}}
\newcommand{\nsigmabar}{\neg\sigmabar}
\newcommand{\sigmae}{\sigma\hspace{-1.75pt}e}
\newcommand{\sigmaebar}{\bar\sigma[e]}
\newcommand{\nsigmaebar}{\neg\sigmaebar}
\newcommand{\canstrat}{\sigma_{\bit}}
\newcommand{\canstratbar}{{\overline{\sigma}_{\bit}}}
\newcommand{\ncanstratbar}{\neg\canstratbar}

\newcommand{\valu}{\Xi}
\newcommand{\valustar}{{\hat \valu}}

\newcommand{\eps}{\varepsilon}
\newcommand{\rew}[1]{\langle#1\rangle}
\newcommand{\ubracket}[1]{\left\llbracket#1\right\rrbracket}

\newcolumntype{C}[1]{>{\centering\let\newline\\\arraybackslash\hspace{0pt}}m{#1}}

\newcommand{\B}{\mathbb{B}}
\newcommand{\G}{\mathbb{G}}
\newcommand{\E}{\mathbb{E}}
\newcommand{\D}{\mathbb{D}}
\renewcommand{\S}{\mathbb{S}}
\renewcommand{\P}{\mathrm{P}}
\newcommand{\M}{\mathrm{M}}

\newcommand{\id}{\bm{1}}

\usetikzlibrary{decorations.pathmorphing}
\usetikzlibrary{decorations.markings}

\tikzstyle{zero}=[circle, draw, align=center, scale=0.5, inner sep=0pt, minimum size=1.25cm]
\tikzstyle{one}=[rectangle, draw, align=center, inner sep=0pt, scale=0.5, minimum size=1.25cm]
\tikzstyle{label}=[rectangle, draw, align=center, inner sep=2pt, scale=0.55, fill=gray!33!white]

\def\x{1.75}
\def\y{1}

\definecolor{orange}{RGB}{255,127,0}
\definecolor{blue}{RGB}{0,127,255}
\definecolor{playblue}{RGB}{0,127,255}
\definecolor{ranred}{RGB}{255,0,0}
\newcommand{\boldall}[1]{\ifmmode\mathbf{#1}\else\textbf{\boldmath{#1}}\fi}

\newcommand{\reach}[1]{\Lambda_{#1}}
\newcommand{\sigdiff}[2]{\Delta(#1,#2)}
\newcommand{\priority}[1]{\Omega(#1)}
\newcommand{\relbit}[1]{\nub^{#1}}
\newcommand{\applied}[2]{\mathfrak{A}_{#1}^{#2}}

\newtheorem{theorem}{Theorem}[section]
\newtheorem{lemma}[theorem]{Lemma}

\newtheorem{corollary}[theorem]{Corollary}

\newtheoremstyle{footnotelemma}
{3pt}
{3pt}
{\itshape}
{}
{\bfseries}
{$(\star)$.}
{.5em}
{}

\theoremstyle{footnotelemma}

\theoremstyle{definition}
\newtheorem{definition}[theorem]{Definition}

\newenvironment{property}[1]{\begin{equation}\tag{#1}\label{property: #1}}{\end{equation}}

\newcommand\blfootnote[1]{%
  \begingroup
  \renewcommand\thefootnote{}\footnote{#1}%
  \addtocounter{footnote}{-1}%
  \endgroup
}

\begin{document}

\title{An Exponential Lower Bound for Zadeh's pivot rule\thanks{The work of Yann Disser and Alexander V. Hopp is supported by the `Excellence Initiative' of the German Federal and State Governments and the Graduate School~CE at TU~Darmstadt.}}

\author[1,\Letter]{Yann Disser}

\author[2,\Letter]{Oliver Friedmann}

\author[1,\Letter]{Alexander V. Hopp}

\affil[1]{Department of Mathematics, TU Darmstadt, Germany}
\affil[2]{Department of Informatics, LMU Munich, Germany}

\maketitle

\begin{abstract}
The question whether the Simplex Algorithm admits an efficient pivot rule remains one of the most important open questions in discrete optimization.
While many natural, deterministic pivot rules are known to yield exponential running times, the random-facet rule was shown to have a subexponential running time.
For a long time, Zadeh's rule remained the most prominent candidate for the first deterministic pivot rule with subexponential running time.
We present a lower bound construction that shows that Zadeh's rule is in fact exponential in the worst case.
Our construction is based on a close relation to the Strategy Improvement Algorithm for Parity Games and the Policy Iteration Algorithm for Markov Decision Processes, and we also obtain exponential lower bounds for Zadeh's rule in these contexts.\blfootnote{\hspace*{-1em}\textsuperscript{\Letter}disser@mathematik.tu-darmstadt.de, oliver.friedmann@googlemail.com, hopp@gsc.tu-darmstadt.de}
\end{abstract}

\section{Introduction}

\setcounter{footnote}{0}

The quest for discovering the best pivot rule for the Simplex Algorithm~\cite{Dantzig1963} remains one of the most important challenges in discrete optimization.
In particular, while several other \emph{weakly} polynomial algorithms for solving Linear Programs have been proposed in the past~\cite{Khachiyan1980,Karmarkar1984,BertsimasVempala/04,KelnerSpielman/06,DunaganVempala/07}, no fully ``combinatorial'' algorithm with strongly polynomial running time is known to date -- in fact, the question whether such an algorithm exists is contained in Smale's list of 18 mathematical problems for the century, among other famous unsolved problems like the Riemann hypothesis and the P versus NP problem~\cite{Smale1998}.
The Simplex Algorithm is inherently combinatorial and may yield a strongly polynomial algorithm if a suitable pivot rule exists.
The question what theoretical worst-case running time can be achieved with a pivot rule for the Simplex Algorithm is closely related to the question what the largest possible (combinatorial) diameter of a polytope is, and, in particular, to the weak Hirsch conjecture that states that the diameter is polynomially bounded~\cite{Dantzig1963,Santos/12,Todd/14}.

For various natural pivot rules, exponential worst-case examples were found soon after the Simplex Algorithm was proposed~\cite{Klee1972,Avis1978,Goldfarb1979}.
These examples are highly imbalanced in the sense that they cause some improving directions to be selected by the pivot rule only rarely, while others are selected often.
Randomized pivot rules were proposed as a way to average out the behavior of the Simplex Algorithm and to thus avoid imbalanced behavior.
The hope that this may lead to a better worst-case performance was met when subexponential upper bounds were eventually established for the random-facet pivot rule~\cite{Kalai1992,MatousekSharirWelzl/96,Hansen2015}.
Other promising candidates for efficient pivot rules were deterministic ``pseudo-random'' rules that balance out the behavior of the algorithm explicitly by considering all past decisions in each step, instead of obliviously deciding for an improvement independently in each step.
The two most prominent examples of such pivot rules are Cunningham's rule~\cite{Cunningham1979} that fixes an order of all possible improvement directions at the start and, in each step, picks the next improving direction in this order in round robin fashion, and Zadeh's rule~\cite{Zadeh1980} that picks an improving direction chosen least often so far in each step.
By design, bad examples are much more difficult to construct for these more balanced pivoting rules, and it took more than 30 years until the first lower bounds were established.
Eventually, a subexponential lower bound was shown for the random-facet rule~\cite{Friedmann2011Sub,Hansen2012,FriedmannHansenZwickErrata/14} and the random-edge rule~\cite{MatousekSzabo/06,Friedmann2011Sub}.
Most recently, a subexponential bound was shown for Zadeh's rule using an artificial tie-breaking rule~\cite{Friedmann2011,DisserHopp/19}, and an exponential bound for Cunningham's rule~\cite{Avis2017}.
An exponential lower bound for Zadeh's rule using a tie-breaking rule that is implemented as an ordered list is known on Acyclic Unique Sink Orientations~\cite{Antonis/17}, but it is unclear whether the corresponding construction can be realized as a Linear Program.
This means that Zadeh's rule remained the only promising candidate for a deterministic pivot rule to match the subexponential running time of the random-facet rule.

Local search algorithms similar to the Simplex Algorithm are important in other domains like V\"oge and Jurdzi\'nski's Strategy Improvement Algorithm for Parity Games~\cite{VoegeJurdzinski2000} and Howard's Policy Iteration Algorithm for Markov Decision Processes~\cite{Howard1960}.
Much like the Simplex Algorithm, these algorithms rely on a pivot rule that determines which local improvement to perform in each step.
And much like for the Simplex Algorithm, many natural deterministic pivot rules for these algorithms have been shown to be exponential~\cite{Friedmann2009,Fearnley2010Exp,Friedmann2011exp,Avis2017}, while a subexponential bound has been shown for the random-facet rule~\cite{Kalai1992,MatousekSharirWelzl/96,Kalai1997,FriedmannHansenZwick/11,FriedmannHansenZwickErrata/14}.
Again, Zadeh's rule remained as a promising candidate for a deterministic subexponential pivot~rule.

\vspace*{-1em}
\paragraph{Our results and techniques.} In this paper, we give the first exponential lower bound for Zadeh's pivot rule for the Strategy Improvement Algorithm for Parity Games, for the Policy Iteration Algorithm for Markov Decision Processes, and for the Simplex Algorithm.
This closes a long-standing open problem by ruling out Zadeh's pivot rule as a candidate for a deterministic, subexponential pivot rule in each of these three areas (up to tie-breaking).
Our result joins a recent series of lower bound results for different pivot rules.\footnote{Several of these results were highly valued by the community, e.g., the best paper award at STOC 2011 was awarded for the subexponential lower bound the random-facet rule~\cite{Friedmann2011Sub}, the Tucker Prize 2012 was awarded for the subexponential lower bound for Zadeh's rule~\cite{Friedmann2011}, and the Kleene Award 2009 was awarded for the exponential lower bound for the Strategy Improvement Algorithm~\cite{Friedmann2009}.} 
We note that while the lower bound for the Simplex Algorithm is arguably our most important result, the lower bounds for Parity Games and Markov Decision Processes are important in their own right and complement previous results in these areas~\cite{Friedmann2009,Fearnley2010Exp,Friedmann2011exp,FriedmannHansenZwick/11}.

Our lower bound construction is based on the technique used in~\cite{Friedmann2011,Avis2017} (among others).
In particular, we construct a Parity Game that forces the Strategy Improvement Algorithm to emulate a binary counter by enumerating strategies corresponding to the natural numbers~$0$ to~$2^{n-1}$.
The construction is then converted into a Markov Decision Process that behaves similarly (but not identically) regarding the Policy Iteration Algorithm.
Finally, using a well-known transformation, the Markov Decision Process can be turned into a Linear Program for which the Simplex Algorithm mimics the behavior of the Policy Iteration Algorithm.
We remark that we use an artificial, but systematic and polynomial time computable, tie-breaking rule for the pivot step whenever Zadeh's rule does not yield a unique improvement direction.
Importantly, while the tie-breaking rule is carefully crafted to simplify the analysis, conceptually, our construction is not based on exploiting the tie-breaking rule.
Note that it cannot be avoided to fix a tie-breaking rule when analyzing Zadeh's pivot rule, in the sense that, for every Markov Decision Process of size~$n$, a tie-breaking rule tailored to this instance exists, such that the policy iteration algorithm takes at most~$n$ steps~\cite[Corollary~4.79]{FriedmannThesis}.

Roughly speaking, much like the subexponential construction in~\cite{Friedmann2011}, our construction consists of multiple levels, one for each bit of the counter. 
The subexponential construction requires each level to connect to the level of the least significant set bit of the currently represented number, which yields a quadratic number~$m$ of edges in the construction, which in turn leads to a lower bound of~$2^{\Omega(n)} = 2^{\Omega(\sqrt{m})}$, i.e., a subexponential bound in the size~$\Theta(m)$ of the construction.
In contrast, our construction only needs each level to connect to one of the first two levels, depending on whether the currently represented number is even or odd.
Very roughly, this is the key idea of our result, since it allows us to reduce the size of the construction to~$\Theta(n)$, which leads to an exponential lower bound.
However, to make this change possible, many other technical details have to be addressed, and, in particular, we are no longer able to carry the construction for Parity Games over as-is to Markov Decision Processes.

A main challenge when constructing a lower bound for Zadeh's rule is to keep track not only of the exact sets of improving directions in each step, but also of the exact number of times every improving direction was selected in the past.
In contrast, the exponential lower bound construction for Cunningham's rule~\cite{Avis2017} ``only'' needs to keep track of the next improving direction in the fixed cyclic order.
As a consequence, the full proof of our result is very technical, because it requires us to consider \emph{all} possible improvements in every step, and there are many transitional steps between configurations representing natural numbers.
In this extended abstract we give an exact description of our construction and an outline of our proof. 
A complete and detailed proof can be found in the full version~\cite{FriedmannDisserHopp2019Arxiv}.
Importantly, our construction has been implemented and tested empirically for consistency with our formal treatment using the \textsc{PGSolver} library \cite{FriedmannLange/PGSolver}.
The resulting animations of the execution for~$n=3$ resp.~$n=4$, which take 160 resp. 466 steps, are available online~\cite{FriedmannImplemetation}.

\section{Parity Games and Strategy Improvement} \label{section: Parity Games}

A \emph{Parity Game (PG)} is a two player game that is played on a directed graph where every vertex has at least one outgoing edge.
Formally, it is defined as a tuple $G=(V_0,V_1,E,\Omega)$, where $(V_0\cup V_1,E)$ is a directed graph and $\Omega\colon V_0\cup V_1\to\mathbb{N}$ is the \emph{priority function}.
The set $V_p$ is the set of vertices of player $p\in\{0,1\}$ and the set $E_p\coloneqq\{(v,w)\in E\colon v\in V_p\}$ is the set of edges of player $p\in\{0,1\}$.
For convenience, we define $V\coloneqq V_0\cup V_1$.
A \emph{play} in $G$ is an infinite walk through the graph.
The \emph{winner} of a play is determined by the highest priority that occurs infinitely often along the walk.
If this priority is even, player 0 wins, otherwise, player 1 wins.

Formally, a play in $G$ can be described by a pair of \emph{strategies}.
A strategy for player $p$ is a function that chooses one outgoing edge for each vertex of player $p$.
To be precise, a \emph{(deterministic positional) strategy} for player $p$ is a function $\sigma\colon V_p\to V$ that selects for each vertex $v\in E_p$ a target vertex $\sigma(v)$ such that $(v,\sigma(v))\in E_p$ for all $v\in V_p$.
Throughout this paper we only consider deterministic positional strategies and henceforth simply refer to them as \emph{strategies}.
Two strategies $\sigma,\tau$ for players~$0,1$ and a starting vertex $v$ then define a unique play starting at $v$ with the corresponding walk being determined by the strategies of the players.
A play can thus fully be described by a tuple $(\sigma,\tau,v)$ and is denoted by $\pi_{\sigma,\tau,v}$.
A player~$0$ strategy $\sigma$ is \emph{winning} for player~$0$ at vertex $v$, if player $0$ is the winner of every game $\pi_{\sigma,\tau,v}$, regardless of $\tau$.
Winning strategies for player $1$ are defined analogously.
One of the most important theorems in the theory of parity games is that, for every starting vertex, there is always a winning strategy for exactly one of the two players.
The computational problem of \emph{solving} a parity game is to find corresponding partition of $V$.

\begin{theorem}[e.g. \cite{Kuesters2002,Fearnley2015}] \label{thm: PG Are Postionally Determined}
In every parity game, $V$ can be partitioned into \emph{winning sets} $(W_0,W_1)$, where player~$p$ has a positional winning strategy for all $v\in W_p$.
\end{theorem}

\subsection{Vertex Valuations, the Strategy Improvement Algorithm and Sink Games} \label{section: VV SIA and SG}

We now discuss the Strategy Improvement Algorithm of V\"oge and Jurdzi\'nski \cite{VoegeJurdzinski2000} and its theoretical background.
We discuss the concept of \emph{vertex valuations} and define a special class of games that our construction belongs to, called \emph{sink games} and define vertex valuations for this class of games.
We refer to \cite{Avis2017,FriedmannThesis} for a more in-depth and general discussion of these topics.

Fix a pair~$\sigma,\tau$ of strategies for players~$0,1$, respectively.
The idea of vertex valuations is to assign a valuation to every $v\in V$ that encodes how ``profitable'' vertex $v$ is for player~0.
By defining a suitable pre-order on these valuations, this enables us to compare the valuations of vertices and ``improve'' the strategy~$\sigma$ by changing the target $\sigma(v)$ of a vertex $v$ to a more ``profitable'' vertex $w\neq\sigma(v)$ with $(v,w)\in E$.
Since there are only finitely many strategies and vertices, improving the strategy of player~$0$ terminates at some point, resulting in a so-called \emph{optimal strategy} for player~$0$.
It is known (e.g. \cite{VoegeJurdzinski2000,FriedmannThesis}) that an optimal strategy can then be used to determine the winning sets~$W_0,W_1$ of the parity game and thus solve the game.

Formally, vertex valuations are given as a totally ordered set $(U,\preceq)$.
For every pair of strategies~$\sigma,\tau$, we are given a function $\valu_{\sigma,\tau}\colon V\to U$ assigning vertex valuations to vertices.
Since $U$ is totally ordered, this induces an ordering of the vertices for fixed strategies $\sigma,\tau$.
To eliminate the dependency on the player 1 strategy, we define the vertex valuation of $v$ with respect to $\sigma$ by $\valu_{\sigma}(v)\coloneqq\min_{\prec}\valu_{\sigma,\tau}(v)$ where the minimum is taken over all player 1 strategies $\tau$.
Formally, if $\valu_{\sigma,\tau}(\tau(v))\preceq\valu_{\sigma}(v)$ for all $(v,u)\in E_1$, then the player 1 strategy $\tau$ is called \emph{counterstrategy} for $\sigma$.
It is well-known that counterstrategies exist and can be computed efficiently \cite{VoegeJurdzinski2000}.
For a strategy~$\sigma$, an arbitrary but fixed counterstrategy is denoted by $\tau^{\sigma}$.

We can extend this ordering to a partial ordering of strategies by defining $\sigma\unlhd\sigma'$ if and only if $\valu_{\sigma}(v)\preceq\valu_{\sigma'}(v)$ for all $v\in V$.
We write $\sigma\lhd\sigma'$ if $\sigma\unlhd\sigma'$ and $\sigma\neq\sigma'$.
Given a strategy~$\sigma$, a strategies $\sigma'$ with $\sigma\lhd\sigma'$ can be obtained by applying \emph{improving switches}.
Intuitively, an improving switch is an edge such that including $e$ in $\sigma$ improves the strategy with respect to $\unlhd$.
Formally, let $e=(v,u)\in E_0$ and $\sigma(v)\neq u$.
We define $\sigmae$ via $\sigmae(v')\coloneqq\sigma(v)$ if $v'\neq v$ and $\sigmae(v)\coloneqq w$.
The edge $e$ is \emph{improving} for $\sigma$ if $\sigma\lhd\sigmae$ and we denote the set of improving switches for $\sigma$ by $I_{\sigma}$.

The Strategy Improvement Algorithm now operates as follows.
Given a initial strategy $\iota$, apply improving switches until a strategy $\sigma^*$ with $I_{\sigma^*}=\emptyset$ is reached.
Such a strategy is called \emph{optimal} and a strategy is optimal if and only if $\sigma\lhd\sigma^*$ for all player 0 strategies \cite{VoegeJurdzinski2000}.
The running time of this algorithm highly depends on the order in which improving switches are applied - a point which we discuss later in more detail.

This terminology allows us to introduce a special class of Parity Games, called \emph{sink games}.
This class allows for an easy definition of the vertex valuations as discussed after the definition.

\begin{definition} \label{definition: Sink Game}
A Parity Game $G=(V_0,V_1,E,\Omega)$ together with an initial player 0 strategy $\iota$ is a \emph{Sink Game} if the following two statements hold.
\begin{enumerate}
	\item There is a vertex $t$ with $(t,t)\in E$ and $\priority{t}=1$ reachable from all vertices.
		In addition, $\priority{v}>\priority{t}$ for all $v\in V\setminus\{t\}$.
		This unique vertex $t$ is called the \emph{sink} of the sink game.
	\item For each player $0$ strategy $\sigma$ with $\iota\unlhd\sigma$ and each vertex $v$, every play $\pi_{\sigma,\tau^{\sigma},v}$ ends in $t$.
\end{enumerate}
\end{definition}

Let $G=(V_0,V_1,E,\Omega)$ and $\iota$ define a Sink Game.
To simplify the presentation, assume that $\Omega$ is injective.
Since $G$ is a Sink Game, every play $\pi_{\sigma,\tau,v}$ in $G$ can be represented as the walk $\pi_{\sigma,\tau,v}=v,v_2,\dots,v_k,(t)^{\infty}$.
In particular, a play can be identified with its \emph{path component} $v,v_2,\dots,v_k$.
Now, defining $\valu_{\sigma}(v)$ as the path component of $\pi_{\sigma,\tau,v}$ is a well-studied choice of vertex valuations.
To give a total ordering of the vertex valuations, it thus suffices to give a ordering of all subsets of $V$.

Let $M,N\subseteq V, M\neq N$.
Intuitively, $M$ is better than $N$ for player 0 if it contains a vertex with large even priority not contained in $N$, or if it there is a vertex with large odd priority contained in~$N$ but not in $M$.
Formally, $v\in M\Delta N$ is called \emph{most significant difference of $M$ and $N$} if $\priority{v}>\priority{w}$ for all $w\in M\Delta N, w\neq v$.
The most significant difference $M$ and $N$ is denoted by $\sigdiff{M}{N}$ and allows us to define an ordering $\prec$ on the subsets of $V$.
For $M,N\subset V, M\neq N$ we define \[M\prec N \Longleftrightarrow [\sigdiff{M}{N}\in N\wedge \sigdiff{M}{N}\text{ is even}]\vee[\sigdiff{M}{N}\in M\wedge\sigdiff{M}{N}\text{ is odd}].\]
Note that $\prec$ is a total ordering as we assume $\Omega$ to be injective.
We mention here that injectivity is not necessary - it suffices if the most significant difference of any two vertex valuations is unique.

The following theorem summarizes the most important aspects related to parity games, vertex valuations and improving switches.
Note that the construction of vertex valuations given here is a simplified version of the general concept of vertex valuations used for parity games.
It is, however, in accordance with the general construction and we refer to \cite{FriedmannThesis} for a more detailed discussion.

\begin{theorem}[\cite{VoegeJurdzinski2000}]
Let $G=(V_0,V_1,E,\Omega)$ be a sink game and $\sigma$ be a player 0 strategy.
\begin{enumerate}
	\item The vertex valuations of a player 0 strategy are polynomial-time computable.
	\item There is an optimal player 0 strategy $\sigma^*$ with respect to the ordering $\lhd$.
	\item If $I_{\sigma}=\emptyset$, then $\sigma$ is optimal.
	\item We have $I_{\sigma}=\{(v,w)\in E_0\colon \valu_{\sigma}(\sigma(v))\lhd\valu_{\sigma}(w)\}$ and $\sigma\lhd\sigmae$ for all $e\in I_{\sigma}$.
	\item Given an optimal player 0 strategy, the winning sets $W_0$ and $W_1$ of player 0 and player 1 can be computed in polynomial time.
\end{enumerate}
\end{theorem}
\section{Lower Bound Construction} \label{section: Lower Bound Construction}

\begin{figure}[ht]
\centering
\begin{tikzpicture}[scale=0.6]
\def\x{3}
\def\y{1}

\foreach \i \in in {1,...,4}{
	\draw (0,1.5*\i) rectangle (\x,\y+1.5*\i);
	\node at (-1.15,{0.5*\y+1.5*\i}) {Level $\i$};
	\node[draw,circle, inner sep=1pt, scale=0.7] (F\i0) at (1,{0.5*\y+1.5*\i}) {$F_{\i,0}$};
	\node[draw,circle, inner sep=1pt, scale=0.7] (F\i0) at (2,{0.5*\y+1.5*\i}) {$F_{\i,1}$};
}

\begin{scope}[xshift=1.75*\x cm]
\def\x{3}
\def\y{1}
\foreach \i \in in {1,...,4}{
	\draw (0,1.5*\i) rectangle (\x,\y+1.5*\i);
}
\node at (-0.25,{0.5*\y+1.5*4}) {1};
\draw[fill=red, line width=3pt] (1,{0.5*\y+1.5*4}) circle (0.33);
\draw[fill=none, dashed, line width=3pt] (2,{0.5*\y+1.5*4}) circle (0.33);
\node at (-0.25,{0.5*\y+1.5*3}) {0};
\draw[fill=none, dashed, line width=3pt] (1,{0.5*\y+1.5*3}) circle (0.33);
\draw[fill=none, line width=3pt] (2,{0.5*\y+1.5*3}) circle (0.33);
\node at (-0.25,{0.5*\y+1.5*2}) {1};
\draw[fill=red, line width=3pt] (1,{0.5*\y+1.5*2}) circle (0.33);
\draw[fill=none, dashed, line width=3pt] (2,{0.5*\y+1.5*2}) circle (0.33);
\node at (-0.25,{0.5*\y+1.5*1}) {1};
\draw[fill=none, dashed, line width=3pt] (1,{0.5*\y+1.5*1}) circle (0.33);
\draw[fill=red, line width=3pt] (2,{0.5*\y+1.5*1}) circle (0.33);
\end{scope}

\begin{scope}[xshift=3.25*\x cm]
\def\x{3}
\def\y{1}
\foreach \i \in in {1,...,4}{
	\draw (0,1.5*\i) rectangle (\x,\y+1.5*\i);
}
\node at (-0.25,{0.5*\y+1.5*4}) {0};
\draw[fill=none, line width=3pt] (1,{0.5*\y+1.5*4}) circle (0.33);
\draw[fill=none, dashed, line width=3pt] (2,{0.5*\y+1.5*4}) circle (0.33);
\node at (-0.25,{0.5*\y+1.5*3}) {0};
\draw[fill=none, line width=3pt] (1,{0.5*\y+1.5*3}) circle (0.33);
\draw[fill=none, dashed, line width=3pt] (2,{0.5*\y+1.5*3}) circle (0.33);
\node at (-0.25,{0.5*\y+1.5*2}) {1};
\draw[fill=red, line width=3pt] (1,{0.5*\y+1.5*2}) circle (0.33);
\draw[fill=none, dashed, line width=3pt] (2,{0.5*\y+1.5*2}) circle (0.33);
\node at (-0.25,{0.5*\y+1.5*1}) {1};
\draw[fill=none, dashed, line width=3pt] (1,{0.5*\y+1.5*1}) circle (0.33);
\draw[fill=red, line width=3pt] (2,{0.5*\y+1.5*1}) circle (0.33);

\draw[fill=red, line width=2pt] (1.25*\x,{0.5*\y+1.5*1}) circle (0.22) node[right=5pt] {Closed \& active};
\draw[fill=none, line width=2pt] (1.25*\x,{0.5*\y+1.5*1.4}) circle (0.22) node[right=5pt] {Open \& active};
\draw[fill=none, dashed, line width=2pt] (1.25*\x,{0.5*\y+1.5*1.8}) circle (0.22) node[right=5pt] {Open \& inactive};
\end{scope}
\end{tikzpicture}
\caption{Visualization of the general structure of the binary counter for $n=4$.
Level~1 encodes the least significant bit and level 4 encodes the most significant bit.
The left picture shows the cycle centers and their distribution in the levels.
The two figures on the right give examples for settings of the cycles representing the numbers 11 and 3, respectively.
} \label{figure: Intuitive Idea}
\end{figure}
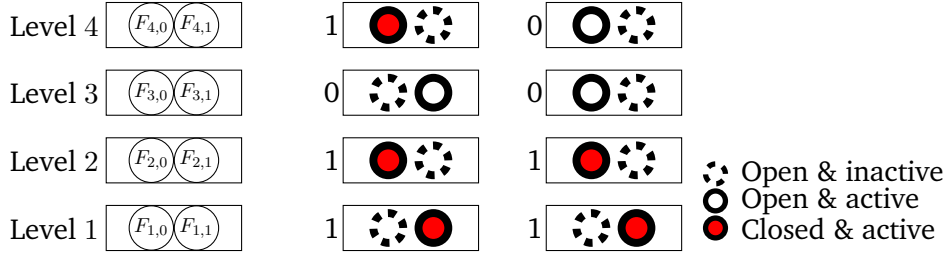

In this section, we describe a PG $S_n=(V_0,V_1,E,\Omega)$ such that the Strategy Improvement Algorithm performs at least $2^n$ iterations when using Zadeh's pivot rule and a specific tie-breaking rule.
Before giving a formal definition, we give a high-level intuition of the main idea of the construction.
A simplified visualization of the construction is given in \Cref{figure: Intuitive Idea}.

The key idea is that $S_n$ simulates an $n$-digit binary counter.
We thus introduce notation related to binary counting.
It will be convenient to consider counter configurations with more than $n$ bits, where unused bits are zero. 
In particular, we always interpret bit $n+1$ as 0.
Formally, we denote the set of $n$-bit configurations by $\bitset_n\coloneqq \{\bit\in\{0,1\}^{\infty}\colon\bit_i=0\quad\forall i>n\}$.
We start with index one, hence a counter configuration $\bit \in \bitset_n$ is a tuple $(\bit_{n},\ldots,\bit_1)$.
Here, $\bit_1$ is the least and $\bit_n$ is the most significant bit. 
The \emph{integer value} of $\bit \in \bitset_n$ is~$\sum_{i=1}^n\bit_i2^{i-1}$.
We identify the integer value of $\bit$ with its counter configuration and use the natural ordering of $\mathbb{N}$ to order counter configurations.
For $\bit\in\bitset_n,\bit\neq0$, we define $\nsb(\bit) \coloneqq \min\{i\in\{1,\dots,n\}\colon \bit_i = 1\}$ to be the \emph{least significant set bit of~$\bit$}.

The PG $S_n$ consists of $n$ (nearly) identical levels and each level encodes one bit of the counter.
Certain strategies and corresponding counterstrategies in $S_n$ are then interpreted as binary numbers.
If the Strategy Improvement Algorithm enumerates at least one player 0 strategy per $\bit\in\bitset_n$ before finding the optimal strategy, it enumerates at least $2^n$ strategies.
Since the game has size linear in $n$, this then establishes the exponential lower bound.

The main challenge is to obey Zadeh's pivot rule as it forces the algorithm to only use improving switches used least often during the execution.
Intuitively, a counter obeying this rule needs to switch bits in a ``balanced'' way.
However, counting from 0 to $2^n-1$ in binary does not switch individual bits equally often.
For example, the least significant bit is switched every time and the most significant bit is switched only once.
The key idea to overcome this obstacle is to have a substructure in each level that contains two gadgets.
These gadgets are called \emph{cycle centers}.
In every iteration of the algorithm, only one of the cycle centers is interpreted as encoding the bit of the current level.
This enables us to perform operations within the other cycle center without loosing the interpretation of the bit being equal 0 to 1.
This is achieved by an alternating encoding of the bit by the two cycle centers.

We now provide more details.
Consider some level $i$, some $\bit\in\bitset_n$ and denote the cycle centers of level~$i$ by $F_{i,0}$ and $F_{i,1}$.
One of them now encodes $\bit_i$.
Which of them represents $\bit_i$ depends on~$\bit_{i+1}$, since we always consider~$F_{i,\bit_{i+1}}$ to encode $\bit_i$.
This cycle center is called the \emph{active} cycle center of level $i$, while~$F_{i,1-\bit_{i+1}}$ is called \emph{inactive}.
A cycle center can additionally be \emph{closed} or \emph{open}.
These terms are used to formalize when a bit is interpreted as 0 or 1.
To be precise, $\bit_i$  is interpreted as~1 if and only if $F_{i,\bit_{i+1}}$ is closed.
In this way, cycle centers encode binary numbers.
Since bit $i+1$ switches every second time bit $i$ switches, counting from 0 to $2^n-1$ in binary then results in an alternating and balanced usage of both cycle centers of any level as required by Zadeh's pivot rule.

We now describe the construction of a parity game that implements this idea in detail.
Fix some $n\in\mathbb{N}$.
The vertex sets~$V_0,V_1$ of the underlying graph are composed as follows:
\begin{align*}
V_0\coloneqq & \{b_i,g_i\colon i\in\{1,\dots,n\}\} \cup \{d_{i,j,k},e_{i,j,k},s_{i,j} \colon i\in\{1,\dots,n-j\}, j,k\in\{0,1\}\} \\
V_1\coloneqq & \{F_{i,j} \colon i\in\{1,\dots, n-j\}, j\in\{0,1\}\} \cup \{h_{i,j} \colon i\in\{1,\dots, n-j\}, j\in\{0,1\}\}\cup \{t\}
\end{align*}

\begin{figure}[ht]
\centering
\begin{tikzpicture}[scale=0.8]

	\def\x{1.5}
	\def\y{0.75}
	\tikzstyle{zero}=[circle, draw, align=center, scale=0.4, inner sep=0pt, minimum size=1.1cm]
	\tikzstyle{one}=[rectangle, draw, align=center, inner sep=0pt, scale=0.4, minimum size=1.1cm]

	\node[zero, dashed] (PairSelector2) at  (0*\x,\y) {$b_2$};
	\node[zero, dashed] (Selector1) at  (0*\x,3*\y) {$g_1$};
	\node[zero, dashed] (PairSelector1) at  (0,4*\y) {$b_1$};
	\node[one] (t) at (0*\x,15*\y) {$t$\\1};
	\draw[->] (t) to[out=60, in=120, looseness=3] (t);
	
	\node[zero] (g1) at (0,0) {$g_1$\\$11$};
	\node[one] (F10) at (-3*\x,2*\y) {$F_{1,0}$\\6};
	\node[one] (F12) at (3*\x,2*\y) {$F_{1,1}$\\4};
	\node[zero] (d101) at (-2*\x,3*\y) {$d_{1,0,1}$\\3};
	\node[zero] (e101) at (-\x,3*\y) {$e_{1,0,1}$\\3};
	\node[zero] (d111) at (2*\x,3*\y) {$d_{1,1,1}$\\3};
	\node[zero] (e111) at (1*\x,3*\y) {$e_{1,1,1}$\\3};
	\node[zero] (d100) at (-2*\x,1*\y) {$d_{1,0,0}$\\3};
	\node[zero] (e100) at (-\x,1*\y) {$e_{1,0,0}$\\3};
	\node[zero] (d110) at (2*\x,1*\y) {$d_{1,1,0}$\\3};
	\node[zero] (e110) at (\x,1*\y) {$e_{1,1,0}$\\3};
	\node[one] (h10) at (-3*\x,4*\y) {$h_{1,0}$\\$12$};
	\node[zero] (b1) at (-4*\x,0) {$b_1$\\3};
	\node[zero] (s10) at (-2*\x,4*\y) {$s_{1,0}$\\10};
	\node[zero] (s11) at (2*\x,4*\y) {$s_{1,1}$\\8};
	\node[one] (h11) at (3*\x,4*\y) {$h_{1,1}$\\$12$};	
	
	\draw[dashed] (-4.5*\x,4.5*\y)--(3.5*\x,4.5*\y);
	
	\draw[->] (e101)-- (PairSelector2);
	\draw[->, red, thick] (e101)--(Selector1);
	
	\draw[->] (e111)-- (PairSelector2);
	\draw[->, red, thick] (e111)--  (Selector1);
	
	\draw[->] (e100)-- (PairSelector2);
	\draw[->, red, thick] (e100)-- (Selector1);

	\draw[->] (e110)--(PairSelector2);
	\draw[->, red, thick] (e110)--(Selector1);
	
	\draw[->] (s10)--(h10);
	\draw[->, red, thick] (s10)--(PairSelector1);
	
	\draw[->, red, thick, dashed] (s11)--(h11);
	\draw[->] (s11)--(PairSelector1);
	
\begin{scope}[yshift=5*\y cm]
	\node[zero, dashed] (PairSelector1) at  (0,4*\y) {$b_1$};
	\node[zero, dashed] (PairSelector2) at  (0*\x,\y) {$b_2$};
	\node[zero, dashed] (Selector1) at  (0*\x,3*\y) {$g_1$};

	\node[zero] (g2) at (0,0) {$g_2$\\$13$};
	\node[one] (F20) at (-3*\x,2*\y) {$F_{2,0}$\\6};
	\node[one] (F22) at (3*\x,2*\y) {$F_{2,1}$\\4};
	\node[zero] (d201) at (-2*\x,3*\y) {$d_{2,0,1}$\\3};
	\node[zero] (e201) at (-\x,3*\y) {$e_{2,0,1}$\\3};
	\node[zero] (d211) at (2*\x,3*\y) {$d_{2,1,1}$\\3};
	\node[zero] (e211) at (1*\x,3*\y) {$e_{2,1,1}$\\3};
	\node[zero] (d200) at (-2*\x,1*\y) {$d_{2,0,0}$\\3};
	\node[zero] (e200) at (-\x,1*\y) {$e_{2,0,0}$\\3};
	\node[zero] (d210) at (2*\x,1*\y) {$d_{2,1,0}$\\3};
	\node[zero] (e210) at (\x,1*\y) {$e_{2,1,0}$\\3};
	\node[one] (h20) at (-3*\x,4*\y) {$h_{2,0}$\\$14$};
	\node[zero] (b2) at (-4*\x,0) {$b_2$\\3};
	\node[zero] (s20) at (-2*\x,4*\y) {$s_{2,0}$\\10};
	\node[zero] (s21) at (2*\x,4*\y) {$s_{2,1}$\\8};
	\node[one] (h21) at (3*\x,4*\y) {$h_{2,1}$\\$14$};	
	
	\draw[dashed] (-4.5*\x,4.5*\y)--(3.5*\x,4.5*\y);
	
	\draw[->, red, thick] (d201) to [out=270-30, in=0+30] (F20);
	\draw[->] (d201)-- (e201);
	
	\draw[->] (e201)-- (PairSelector2);
	\draw[->, red, thick] (e201)--(Selector1);
	
	\draw[->] (d211) to [out=270+30, in=180-30] (F22);
	\draw[->, red, thick] (d211)-- (e211);
	
	\draw[->] (e211)-- (PairSelector2);
	\draw[->, red, thick] (e211)--  (Selector1);
	
	\draw[->, red, thick] (d200) to[out=90+30, in=0-30]  (F20);
	\draw[->] (d200)--(e200);
	
	\draw[->] (e200)-- (PairSelector2);
	\draw[->, red, thick] (e200)-- (Selector1);
	
	\draw[->] (d210) to[out=90-30, in=180+30]  (F22);
	\draw[->, red, thick] (d210)-- (e210);
	
	\draw[->] (e210)--(PairSelector2);
	\draw[->, red, thick] (e210)--(Selector1);
	
	\draw[->, red, thick, dashed] (g2) to[out=180-15, in=270] (F20);
	\draw[->] (g2)to[out=0+15, in=270]  (F22);
	
	\draw[->, red, thick] (b2)--(g2);
	
	\draw[->] (F20) to[out=90-30, in=180+30] (d201);
	\draw[->] (F20) to[out=270+30, in=180-30] (d200);
	\draw[->, blue, thick, dashed] (F20) to[out=90, in=180+30] (s20);
	
	\draw[->] (F22) to[out=90+30, in=0-30] (d211);
	\draw[->, blue, thick] (F22)to[out=270-30, in=0+30]  (d210);
	\draw[->] (F22) to[out=90, in=0-30](s21);
	
	\draw[->, red, thick, dashed] (s20)--(h20);
	\draw[->] (s20)--(PairSelector1);
	
	\draw[->] (s21)--(h21);
	\draw[->, red, thick] (s21)--(PairSelector1);
	\draw[->, blue, thick, dashed] (h20) --++(-\x,6*\y)--(t);
\end{scope}

\begin{scope}[yshift=10*\y cm]

	\node[zero, dashed] (PairSelector1) at  (0,4*\y) {$b_1$};
	\node[zero, dashed] (PairSelector2) at  (0*\x,\y) {$b_2$};
	\node[zero, dashed] (Selector1) at  (0*\x,3*\y) {$g_1$};

	\node[zero] (g3) at (0,0) {$g_3$\\$15$};
	\node[one] (F30) at (-3*\x,2*\y) {$F_{3,0}$\\6};
	\node[zero] (d301) at (-2*\x,3*\y) {$d_{3,0,1}$\\3};
	\node[zero] (e301) at (-\x,3*\y) {$e_{3,0,1}$\\3};
	\node[zero] (d300) at (-2*\x,1*\y) {$d_{3,0,0}$\\3};
	\node[zero] (e300) at (-\x,1*\y) {$e_{3,0,0}$\\3};
	\node[one] (h30) at (-3*\x,4*\y) {$h_{3,0}$\\$16$};
	\node[zero] (b3) at (-4*\x,0) {$b_3$\\3};
	\node[zero] (s30) at (-2*\x,4*\y) {$s_{3,0}$\\10};
	
	\draw[dashed] (-4.5*\x,4.5*\y)--(3.5*\x,4.5*\y);
	
	\draw[->] (d301) to [out=270-30, in=0+30] (F30);
	\draw[->, red, thick] (d301)-- (e301);
	
	\draw[->] (e301)-- (PairSelector2);
	\draw[->, red, thick] (e301)--(Selector1);
	
	\draw[->] (d300) to[out=90+30, in=0-30]  (F30);
	\draw[->, red, thick] (d300)--(e300);
	
	\draw[->] (e300)-- (PairSelector2);
	\draw[->, red, thick] (e300)-- (Selector1);
	
	\draw[->, red, thick] (g3) to[out=180-15, in=270] (F30);
	
	\draw[->] (b3)--(g3);
	
	\draw[->] (F30) to[out=90-30, in=180+30] (d301);
	\draw[->, blue, thick] (F30) to[out=270+30, in=180-30] (d300);
	\draw[->] (F30) to[out=90, in=180+30] (s30);
	
	\draw[->, red, thick] (s30)--(h30);
	\draw[->] (s30)--(PairSelector1);
	\draw[->, blue, thick] (h30) to[out=20,in=180] (t);
	
	\draw[->, red, thick, blue] (b3) --++(0,5*\y)--(t);
	
\end{scope}

	\draw[->] (d101) to [out=270-30, in=0+30] (F10);
	\draw[->, red, thick] (d101)-- (e101);
	
	\draw[->, red, thick] (d111) to [out=270+30, in=180-30] (F12);
	\draw[->] (d111)-- (e111);
	
	\draw[->] (d100) to[out=90+30, in=0-30]  (F10);
	\draw[->, red, thick] (d100)--(e100);
	
	\draw[->, red, thick] (d110) to[out=90-30, in=180+30]  (F12);
	\draw[->] (d110)-- (e110);
	
	\draw[->] (g1) to[out=180-15, in=270] (F10);
	\draw[->, red, thick, dashed] (g1)to[out=0+15, in=270]  (F12);
	
	\draw[->, red, thick, dashed] (b1)--(g1);
	\draw[->] (b1)--(b2);
	\draw[->] (b2)--(b3);
	
	\draw[->] (F10) to[out=90-30, in=180+30] (d101);
	\draw[->, blue, thick] (F10) to[out=270+30, in=180-30] (d100);
	\draw[->, blue] (F10) to[out=90, in=180+30] (s10);
	
	\draw[->] (F12) to[out=90+30, in=0-30] (d111);
	\draw[->] (F12)to[out=270-30, in=0+30]  (d110);
	\draw[->, blue, thick, dashed] (F12) to[out=90, in=0-30](s11);
	
	\draw[->, blue, thick, dashed] (h11) to[out=150,in=0] (g2);
	\draw[->, blue, thick] (h21) to[out=150,in=0] (g3);
	\draw[->, blue, thick] (h10)--(b3);
\end{tikzpicture}
\caption{The graph $S_3$ together with a canonical strategy representing the number 3 in the graph $S_3$.
The dashed copies  of the vertices $g_1,b_1$ and $b_2$ all refer to the corresponding vertices of levels 1 and 2.
Red edges belong to the strategy of player 0, blue edges belong to the counterstrategy of player 1.
The dashed edges indicate the spinal path.}
\label{figure: Representing 5}
\end{figure}
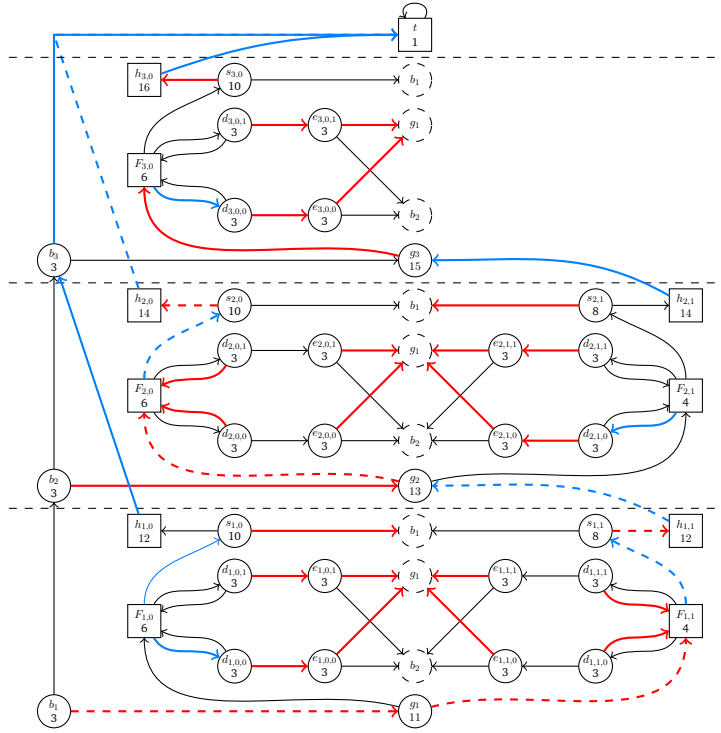

The priorities of the vertices and their sets of outgoing edges are given by \Cref{table: parity game lower bound edges}.
Note that every vertex $v\in V_0$ has at most two outgoing edges. 
For convenience of notation, we henceforth identify the node names~$b_i$ and~$g_i$ for $i > n$ with $t$.
The graph can be separated into $n$ levels, where the levels $i<n-1$ are structurally identical and the levels $n-1$ and $n$ differ slightly from the other levels.
The $i$-th level is shown in \Cref{figure: Level i}, the complete graph $S_3$ is shown in \Cref{figure: Representing 5}.

\begin{figure}[ht]
\centering
\begin{tikzpicture}[scale=0.8]
	\node[zero, dashed] (PairSelector2) at  (0*\x,\y) {$b_2$};
	\node[zero, dashed] (Selector1) at  (0*\x,3*\y) {$g_1$};
	\node[zero, dashed] (PairSelector1) at  (0,4*\y) {$b_1$};
	\node[zero, dashed] (SelectorPlus1) at (-4*\x,5.25*\y) {$b_{i+1}$};
	\node[zero, dashed] (SelectorPlus2) at (-3*\x,5.25*\y) {$b_{i+2}$};
	\node[zero, dashed] (PairSelectorPlus1) at  (0,5.25*\y) {$g_{i+1}$};
	
	\node[zero] (PairSelector) at (0,0) {$g_i$\\$2i+9$};
	\node[one] (CycleCenter1) at (-3*\x,2*\y) {$F_{i,0}$\\6};
	\node[one] (CycleCenter2) at (3*\x,2*\y) {$F_{i,1}$\\4};
	\node[zero] (CycleNodeX1) at (-2*\x,3*\y) {$d_{i,0,1}$\\3};
	\node[zero] (CycleNodeX1T) at (-\x,3*\y) {$e_{i,0,1}$\\3};
	\node[zero] (CycleNodeX2) at (2*\x,3*\y) {$d_{i,1,1}$\\3};
	\node[zero] (CycleNodeX2T) at (1*\x,3*\y) {$e_{i,1,1}$\\3};
	\node[zero] (CycleNodeY1) at (-2*\x,1*\y) {$d_{i,0,0}$\\3};
	\node[zero] (CycleNodeY1T) at (-\x,1*\y) {$e_{i,0,0}$\\3};
	\node[zero] (CycleNodeY2) at (2*\x,1*\y) {$d_{i,1,0}$\\3};
	\node[zero] (CycleNodeY2T) at (\x,1*\y) {$e_{i,1,0}$\\3};
	\node[one] (UpperSelector1) at (-3*\x,4*\y) {$h_{i,0}$\\$2i+10$};
	\node[zero] (Selector) at (-4*\x,0) {$b_i$\\3};
	\node[zero] (UpDown1) at (-2*\x,4*\y) {$s_{i,0}$\\10};
	\node[zero] (UpDown2) at (2*\x,4*\y) {$s_{i,1}$\\8};
	\node[one] (Helper) at (3*\x,4*\y) {$h_{i,1}$\\$2i+10$};	
	
	\draw[dashed] (-4.5*\x,4.5*\y)--(3.5*\x,4.5*\y);
	
	\draw[->] (CycleNodeX1) to [out=270-30, in=0+30] (CycleCenter1);
	\draw[->] (CycleNodeX1)-- (CycleNodeX1T);
	
	\draw[->] (CycleNodeX1T)--(PairSelector2);
	\draw[->] (CycleNodeX1T)--(Selector1);
	
	\draw[->] (CycleNodeX2) to [out=270+30, in=180-30] (CycleCenter2);
	\draw[->] (CycleNodeX2)--  (CycleNodeX2T);
		
	\draw[->] (CycleNodeX2T)-- (PairSelector2);
	\draw[->] (CycleNodeX2T)-- (Selector1);

	\draw[->] (CycleNodeY1) to[out=90+30, in=0-30] (CycleCenter1);
	\draw[->] (CycleNodeY1)--(CycleNodeY1T);
		
	\draw[->] (CycleNodeY1T)-- (PairSelector2);
	\draw[->] (CycleNodeY1T)-- (Selector1);
	
	\draw[->] (CycleNodeY2) to[out=90-30, in=180+30] (CycleCenter2);
	\draw[->] (CycleNodeY2)-- (CycleNodeY2T);
	
	\draw[->] (CycleNodeY2T)--(PairSelector2);
	\draw[->] (CycleNodeY2T)--(Selector1);
	
	\draw[->] (PairSelector) to[out=180-15, in=270] (CycleCenter1);
	\draw[->] (PairSelector)to[out=0+15, in=270]  (CycleCenter2);
		
	\draw[->] (Selector)-- (PairSelector);
	\draw[->] (Selector)--(SelectorPlus1);
%
	\draw[->] (CycleCenter1) to[out=90-30, in=180+30] (CycleNodeX1);
	\draw[->] (CycleCenter1) to[out=270+30, in=180-30] (CycleNodeY1);
	\draw[->] (CycleCenter1) to[out=90, in=180+30] (UpDown1);
	
	\draw[->] (CycleCenter2) to[out=90+30, in=0-30] (CycleNodeX2);
	\draw[->] (CycleCenter2)to[out=270-30, in=0+30]  (CycleNodeY2);
	\draw[->] (CycleCenter2) to[out=90, in=0-30](UpDown2);
	
	\draw[->] (UpDown1)--(UpperSelector1);
	\draw[->] (UpDown1)--(PairSelector1);
	
	\draw[->] (UpDown2)--(Helper);
	\draw[->] (UpDown2)--(PairSelector1);
	
%
	
	\draw[->] (UpperSelector1)--(SelectorPlus2);
	
	\draw[->] (Helper) to[out=90+45, in=0] (PairSelectorPlus1);
\end{tikzpicture}
\caption{Level $i$ of $S_n$ for $i\in\{1,\dots,n-2\}$.
Circular vertices are player~0 vertices, rectangular vertices are player~1 vertices.
Labels below vertex names denote their priorities.
Dashed vertices do not (necessarily) belong to level $i$.}
\label{figure: Level i}
\end{figure}
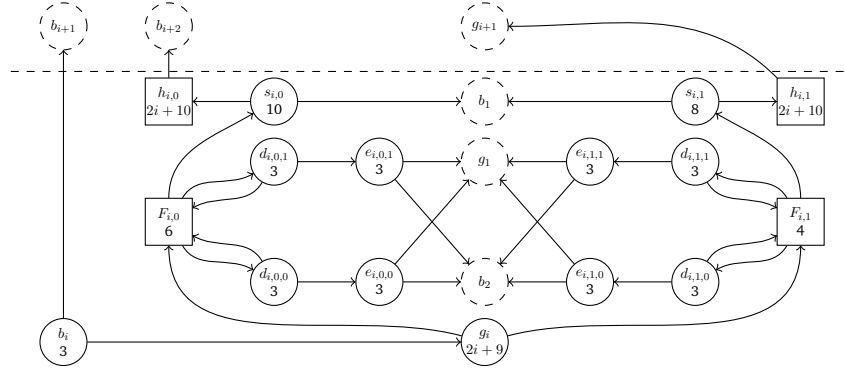

\begin{table}
\centering
\begin{tabular}{cc}
\begin{tabular}[ht]{|c||c|c|}  \hline
  Vertex 				& Successors 				& Priority 			\\  \hline  \hline
  $b_i$     		& $g_i$,     $b_{i+1}$      & $3$  				\\
  $e_{i,j,k}$ 		& $b_2$,     $g_1$     		& $3$  				\\
  $g_i$     		& $F_{i,0}$(, $F_{i,1}$)  	& $2i+9$  	\\
  $h_{i,0}$ 		& $b_{i+2}$                 & $2i+10$  	\\  \hline
\end{tabular}
&
\begin{tabular}{|c||c|c|}\hline
  Vertex 			&Successors 								& Priority \\  \hline  \hline
  $s_{i,j}$ 	& $h_{i,j}$, $b_1$                			& $10-2\cdot j$  \\
  $d_{i,j,k}$ 	& $F_{i,j}$, $e_{i,j,k}$     				& $3$  \\
  $F_{i,j}$ 	& $d_{i,j,0}$, $d_{i,j,1}$, $s_{i,j}$ 		& $6-2\cdot j$  \\
  $h_{i,1}$ 	& $g_{i+1}$                       			& $2i+10$  \\ \hline
\end{tabular}
\end{tabular}
\caption{Edges and vertex priorities of $S_n$.}
\label{table: parity game lower bound edges}
\end{table}

The general idea of the construction is the following.
Certain pairs of player 0 strategies $\sigma$ and counterstrategies $\tau^{\sigma}$ are interpreted as \emph{representing} a number $\bit\in\bitset_n$.
Such a pair of strategies induces a path starting at $b_1$ and ending at $t$, traversing the levels $i\in\{1,\dots,n\}$ with $\bit_i=1$ while ignoring levels with~$\bit_i=0$.
This path is called the \emph{spinal path} with respect to~$\bit\in\bitset_n$.
Ignoring and including levels in the spinal path is controlled by the \emph{entry vertex} $b_i$ of each level $i\in\{1,\dots,n\}$.
To be precise, when $\bit$ is represented, the entry vertex of level $i$ is intended to point towards the \emph{selector vertex}~$g_i$ of level $i$ if and only if $\bit_i=1$.
Otherwise, i.e., when $\bit_i=0$, level $i$ is ignored and the entry vertex $b_i$ points towards the entry vertex of the next level.


Consider a level $i\in\{1,\dots,n-2\}$.
Attached to the selector vertex $g_i$ are the \emph{cycle centers}~$F_{i,0}$ and~$F_{i,1}$ of level $i$.
As described at the beginning of this section, these player 1 vertices are the main structures used for interpreting whether the bit $i$ is equal to one.
They alternate in encoding bit~$i$.
As discussed before, this is achieved by interpreting the \emph{active} cycle center $F_{i,\bit_{i+1}}$ as encoding~$\bit_i$ while the \emph{inactive} cycle center $F_{i,1-\bit_{i+1}}$ does not interfere with the encoding.
This enables us to manipulate the inactive part of a level without loosing the encoded value of $\bit_i$.
To this end, the selector vertex $g_i$ is used to ensure that the active cycle center is contained in the spinal path.

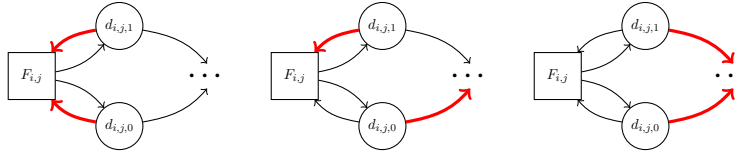
\begin{figure}[ht]
\centering
\begin{tikzpicture}[scale=0.666]

\node[one] (CCycle) at (-3*\x,0) {$F_{i,j}$};
\foreach \z in {-\x,2*\x,5*\x}{
	\node at (\z,0) {\dots};
}
\node[one] (PCCycle) at (0,0) {$F_{i,j}$};
\node[one] (OCycle) at (3*\x,0) {$F_{i,j}$};

\node[zero] (CCVertex1) at (-2*\x,\y) {$d_{i,j,1}$};
\node[zero] (CCVertex0) at (-2*\x,-\y) {$d_{i,j,0}$};

\node[zero] (PCVertex1) at (\x,\y) {$d_{i,j,1}$};
\node[zero] (PCVertex0) at (\x,-\y) {$d_{i,j,0}$};

\node[zero] (OVertex1) at (4*\x,\y) {$d_{i,j,1}$};
\node[zero] (OVertex0) at (4*\x,-\y) {$d_{i,j,0}$};

\draw[->] (CCycle) to[out=0+10,in=270-40] (CCVertex1);
\draw[->, very thick, red] (CCVertex1) to[out=180+10,in=90-40] (CCycle);
\draw[->] (CCVertex1) to[out=0-10,in=90+40] ++(\x,-0.75*\y);
\draw[->] (CCycle) to[out=0-10, in=90+40] (CCVertex0);
\draw[->, very thick, red] (CCVertex0) to[out=180-10, in=270+40] (CCycle); 
\draw[->] (CCVertex0) to[out=0+10,in=270-40] ++(\x,0.75*\y);

\draw[->] (PCCycle) to[out=0+10,in=270-40] (PCVertex1);
\draw[->, very thick, red] (PCVertex1) to[out=180+10,in=90-40] (PCCycle);
\draw[->] (PCVertex1) to[out=0-10,in=90+40] ++(\x,-0.75*\y);
\draw[->] (PCCycle) to[out=0-10, in=90+40] (PCVertex0);
\draw[->] (PCVertex0) to[out=180-10, in=270+40] (PCCycle);  
\draw[->, very thick, red] (PCVertex0) to[out=0+10,in=270-40] ++(\x,0.75*\y);

\draw[->] (OCycle) to[out=0+10,in=270-40] (OVertex1);
\draw[->] (OVertex1) to[out=180+10,in=90-40] (OCycle);
\draw[->, very thick, red] (OVertex1) to[out=0-10,in=90+40] ++(\x,-0.75*\y);
\draw[->] (OCycle) to[out=0-10, in=90+40] (OVertex0);
\draw[->] (OVertex0) to[out=180-10, in=270+40] (OCycle);  
\draw[->, very thick, red] (OVertex0) to[out=0+10,in=270-40] ++(\x,0.75*\y);

\end{tikzpicture}
\caption{A closed, a halfopen and an open cycle center.
Edges of the strategy of player 0 are marked in red.
Choices of player 1 are not depicted.
Dots represent that the cycle vertices point to some unspecified vertex.}
\label{figure: Cycle configurations}
\end{figure}

As discussed previously, a cycle center $F_{i,j}$ can have different configurations.
To be precise, it can be \emph{closed, halfopen,} or \emph{open}.
The configuration of $F_{i,j}$ is defined via the \emph{cycle vertices} $d_{i,j,0}$ and~$d_{i,j,1}$ of the cycle center and the two \emph{cycle edges} $(d_{i,j,0},F_{i,j})$ and $(d_{i,j,1},F_{i,j})$.
More precisely, $F_{i,j}$ is \emph{closed} with respect to a player 0 strategy $\sigma$ if both cycle vertices point towards the cycle center, i.e., when $\sigma(d_{i,j,0})=\sigma(d_{i,j,1})=F_{i,j}$.
If this is the case for exactly one of the two edges, the cycle center $F_{i,j}$ is called \emph{halfopen}.
A cycle that is neither closed nor halfopen is called \emph{open}.
An example of the different configurations is given in \Cref{figure: Cycle configurations}.

In addition, the cycle center is connected to its \emph{upper selection vertex} $s_{i,j}$.
It connects the cycle center $F_{i,j}$ with the first level via $(s_{i,j},b_1)$ and with either level $i+1$ or $i+2$ via $(s_{i,j},h_{i,j})$ via the respective edge $(h_{i,0},b_{i+2})$ or $(h_{i,1},g_{i+1})$ (depending on $j$).
This vertex is thus central in allowing~$F_{i,j}$ to get access to either the beginning of the spinal path or the next level of the spinal path.

We next discuss the cycle vertices.
If their cycle centers are not closed, these vertices still need to be able to access the spinal path.
The valuation of vertices along this path is usually very high and it is almost always very profitable for player 0 vertices to get access to this path.
Since the cycle vertices cannot obtain access via the cycle center (as this would, by definition, close the cycle center) they need to ``escape'' the level in another way.
This is handled by the \emph{escape vertices} $e_{i,j,0}$ and $e_{i,j,1}$.
The escape vertices are used to connect the levels with higher indices to the first two levels and thus enable each vertex to access the spinal path.
To be precise, they are connected with the entry vertex of level 2 and the selector vertex of level~1.
In principle, the escape vertices will point towards~$g_1$ when the least significant set bit of the currently represented number has the index 1 and towards~$b_2$ otherwise.

We now formalize the idea of a strategy encoding a binary number by defining the notion of a \emph{canonical strategy}.
Note that the definition also includes some aspects that are purely technical, i.e., solely required for some proofs, and do not have an immediate intuitive interpretation.

\begin{definition}\label{definition: Canonical Strategy}
Let $\bit\in\bitset_n$.
A player 0 strategy $\sigma$ for the Parity Game $S_n$	 is called \emph{canonical strategy for $\bit$} if it has the following properties.
\begin{enumerate}[parsep=-2pt]
	\item All escape vertices point to $g_1$ if $\bit_1=1$ and to $b_2$ if $\bit_1=0$.
	\item The following hold for all levels $i\in\{1,\dots,n\}$ with $\bit_i=1$:
	\begin{enumerate}[parsep=-1pt]
		\item Level $i$ needs to be accessible, i.e., $\sigma(b_i)=g_i$.
		\item The cycle center $F_{i,\bit_{i+1}}$ needs to be closed while $F_{i,1-\bit_{i+1}}$ must not be closed.			
		\item The selector vertex of level $i$ needs to select the active cycle center, i.e., $\sigma(g_i)=F_{i,\bit_{i+1}}$.
	\end{enumerate}	
	\item The following hold for all levels $i\in\{1,\dots,n\}$ with $\bit_i=0$:
	\begin{enumerate}[parsep=-1pt]
		\item Level $i$ must not be accessible and needs to be ``avoided'', i.e., $\sigma(b_i)=b_{i+1}$.
		\item The cycle center $F_{i,\bit_{i+1}}$ must not be closed.
		\item If the cycle center $F_{i,1-\bit_{i+1}}$ is closed, then $\sigma(g_i)=F_{i,1-\bit_{i+1}}$.
		\item If none of the cycle centers $F_{i,0},F_{i,1}$ is closed, then $\sigma(g_i)=F_{i,0}$.		
	\end{enumerate}
	\item Let $\bit_{i+1}=0$.
		Then, level $i+1$ is not accessible from level $i$, i.e., $\sigma(s_{i,0})=h_{i,0}$ and $\sigma(s_{i,1})=b_1$.
	\item Let $\bit_{i+1}=1$.
		Then, level $i+1$ is accessible from level $i$, i.e., $\sigma(s_{i,0})=b_1$ and $\sigma(s_{i,1})=h_{i,1}$. 
	\item Both cycle centers of level $\nsb(\bit+1)$ are open. 
\end{enumerate}
\end{definition}

We use $\canstrat$ to denote a canonical strategy for $\bit\in\bitset_n$.
A canonical strategy representing $(0,1,1)$ in $S_3$ is shown in \Cref{figure: Representing 5}.

As mentioned before, the main structure that is used to determine whether a bit is interpreted as being set are the cycle centers.
In fact, any configuration of the cycle centers can be interpreted as an encoded number in the following way.

\begin{definition} \label{definition: Representative and Induced States}
Let $\sigma$ be a player 0 strategy for $S_n$.
Then, the \emph{induced bit state} $\indbit^\sigma=(\indbit_n^\sigma,\dots,\indbit_1^\sigma)$ is defined inductively as follows.
We define $\indbit_n^{\sigma}=1$ if and only if $\sigma(d_{n,0,0})=\sigma(d_{n,0,1})=F_{n,0}$ and $\indbit_i^\sigma=1$ if and only if $\sigma(d_{i,\indbit_{i+1}^{\sigma},0})=\sigma(d_{i,\indbit_{i+1}^{\sigma},1})=F_{i,\indbit_{i+1}^{\sigma}}$ for $i<n$
\end{definition}

This definition is in accordance with our interpretation of encoding a number as $\indbit^{\canstrat}=\bit$ if $\canstrat$ is a canonical strategy for $\bit$.
%


\section{Lower Bound for Policy Iteration on MDPs} \label{section: Markov Decision Process}

In this section we discuss the Markov Decision Process (MDP) that is constructed analogously to the PG $S_n$.
We discuss how this MDP allows the construction of a Linear Program (LP) such that the results obtained for the MDP carry over to the LP formulation.
The main idea is to replace player~1 by the ``random player'' and to choose the probabilities in such a way that applying improving switches in the MDP behaves nearly the same way as in the PG.
Note that we continue to use the same language for valuations, strategies and so on in MDP context, although other notions (like policy instead of strategy) are more common.

We give a brief introduction to the theory of MDPs (see also~\cite{Avis2017}).
Similarly to a PG, an MDP is formally defined by its underlying graph $(V_0,V_R,E,r,p)$.
Here,~$V_0$ is the set of vertices controlled by player 0 and $V_R$ is the set of randomization vertices.
We let $V\coloneqq V_0\cup V_R$.
For $p\in\{0,R\}$, we define $E_{p}\coloneqq\{(v,w)\colon v\in E_p\}$.
The set~$E_0$ then corresponds to possible choices that player~0 can make, and each such choice is assigned a \emph{reward} by the reward function $r\colon E_0\to\mathbb{R}$.
The set $E_R$ corresponds to probabilistic transitions and transition probabilities are specified by the function $p\colon E_R\to[0,1]$ with $\sum_{u\colon(v,u)\in E_R}p(v,u)=1$.

As for $S_n$, a \emph{(player 0) strategy} is a function $\sigma\colon V_0\to V$ that selects for each vertex $v\in V_0$ a target corresponding to an edge, i.e., such that $(v,\sigma(v))\in E_0$.
There are several computational tasks that can be investigated for MDPs.
They are typically described via an \emph{objective}.
We consider the \emph{expected total reward objective} for MDPs which can be formulated using vertex valuations in the following sense.
Given an MDP, we define the \emph{vertex valuations} $\valu_{\sigma}^\M(*)$ with respect to a strategy $\sigma$ as the solution (if it exists) of the following set of equations: 
\[\valu_{\sigma}^\M(u)\coloneqq\begin{cases}r(u,\sigma(u))+\valu_{\sigma}^\M(\sigma(u)), &u\in V_0,\\\sum\limits_{v\colon(u,v)\in E_R}p(u,v)\cdot\valu_{\sigma}^\M(v), &u\in V_R.\end{cases}\]
We also impose the condition that the values sum up to 0 on each irreducible recurrent class of the Markov chain defined by $\sigma$, yielding uniqueness \cite{Avis2017}.
We intentionally use very similar notation as for vertex valuations in the context of Parity Games since this allows for a unified treatment.


We now discuss the \emph{Policy Iteration Algorithm} and refer to \cite{Howard1960} for further details.
Similar to the Strategy Improvement Algorithm for PGs, this algorithm starts with some initial policy $\iota=\sigma_0$.
In each step $i$, it generates a strategy $\sigma_i$ by changing the target vertex $\sigma_{i-1}(v)$ of some vertex $v\in V_0$ to some vertex~$w$ with $\valu_{\sigma}^M(w)>\valu_{\sigma}^M(\sigma_{i-1}(v))$.
For an arbitrary strategy $\sigma$, such an edge $(v,w)\in E_0$ with $w\neq\sigma(v)$ but $\valu_{\sigma}^\M(w)>\valu_{\sigma}^\M(\sigma(v))$ is called \emph{improving switch} and the set of improving switches is denoted by $I_{\sigma}$.
The term \emph{optimal strategy} is defined as in PG context.
In particular, a strategy~$\sigma$ is optimal if and only if $I_{\sigma}=\emptyset$.
Moreover, applying an improving switch cannot decrease the valuation of any vertex.
That is, if $e=(v,w)\in I_{\sigma}$ and $\sigmae$ denotes the strategy obtained after applying $e$ to $\sigma$, then $\valu_{\sigmae}^\M(v')\geq\valu_{\sigma}^\M(v')$ for all $v'\in V$ and $\valu_{\sigmae}^\M(v)>\valu_{\sigma}^\M(v)$.
Since there are only finitely many strategies, the algorithm thus generates a finite sequence $\sigma_0,\sigma_1,\dots,\sigma_N$ with $I_{\sigma_N}=\emptyset$.

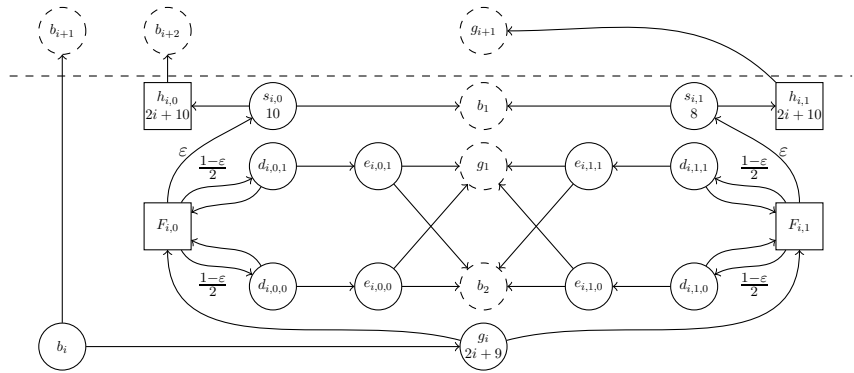
\begin{figure}[ht]
\begin{center}
\begin{tikzpicture}[scale=0.8]
	\node[zero, dashed] (PairSelector2) at  (0*\x,\y) {$b_2$};
	\node[zero, dashed] (Selector1) at  (0*\x,3*\y) {$g_1$};
	\node[zero, dashed] (PairSelector1) at  (0,4*\y) {$b_1$};
	\node[zero, dashed] (SelectorPlus1) at (-4*\x,5.25*\y) {$b_{i+1}$};
	\node[zero, dashed] (SelectorPlus2) at (-3*\x,5.25*\y) {$b_{i+2}$};
	\node[zero, dashed] (PairSelectorPlus1) at  (0,5.25*\y) {$g_{i+1}$};
	
	\node[zero] (PairSelector) at (0,0) {$g_i$\\$2i+9$};
	\node[one] (CycleCenter1) at (-3*\x,2*\y) {$F_{i,0}$};
	\node[one] (CycleCenter2) at (3*\x,2*\y) {$F_{i,1}$};
	\node[zero] (CycleNodeX1) at (-2*\x,3*\y) {$d_{i,0,1}$};
	\node[zero] (CycleNodeX1T) at (-\x,3*\y) {$e_{i,0,1}$};
	\node[zero] (CycleNodeX2) at (2*\x,3*\y) {$d_{i,1,1}$};
	\node[zero] (CycleNodeX2T) at (1*\x,3*\y) {$e_{i,1,1}$};
	\node[zero] (CycleNodeY1) at (-2*\x,1*\y) {$d_{i,0,0}$};
	\node[zero] (CycleNodeY1T) at (-\x,1*\y) {$e_{i,0,0}$};
	\node[zero] (CycleNodeY2) at (2*\x,1*\y) {$d_{i,1,0}$};
	\node[zero] (CycleNodeY2T) at (\x,1*\y) {$e_{i,1,0}$};
	\node[one] (UpperSelector1) at (-3*\x,4*\y) {$h_{i,0}$\\$2i+10$};
	\node[zero] (Selector) at (-4*\x,0) {$b_i$};
	\node[zero] (UpDown1) at (-2*\x,4*\y) {$s_{i,0}$\\$10$};
	\node[zero] (UpDown2) at (2*\x,4*\y) {$s_{i,1}$\\$8$};
	\node[one] (Helper) at (3*\x,4*\y) {$h_{i,1}$\\$2i+10$};	
	
	\draw[dashed] (-4.5*\x,4.5*\y)--(3.5*\x,4.5*\y);
	
	\draw[->] (CycleNodeX1) to [out=270-30, in=0+30] (CycleCenter1);
	\draw[->] (CycleNodeX1)-- (CycleNodeX1T);
	
	\draw[->] (CycleNodeX1T)--(PairSelector2);
	\draw[->] (CycleNodeX1T)--(Selector1);
	
	\draw[->] (CycleNodeX2) to [out=270+30, in=180-30] (CycleCenter2);
	\draw[->] (CycleNodeX2)--  (CycleNodeX2T);
		
	\draw[->] (CycleNodeX2T)-- (PairSelector2);
	\draw[->] (CycleNodeX2T)-- (Selector1);

	\draw[->] (CycleNodeY1) to[out=90+30, in=0-30] (CycleCenter1);
	\draw[->] (CycleNodeY1)--(CycleNodeY1T);
		
	\draw[->] (CycleNodeY1T)-- (PairSelector2);
	\draw[->] (CycleNodeY1T)-- (Selector1);
	
	\draw[->] (CycleNodeY2) to[out=90-30, in=180+30] (CycleCenter2);
	\draw[->] (CycleNodeY2)-- (CycleNodeY2T);
	
	\draw[->] (CycleNodeY2T)--(PairSelector2);
	\draw[->] (CycleNodeY2T)--(Selector1);
	
	\draw[->] (PairSelector) to[out=180-15, in=270] (CycleCenter1);
	\draw[->] (PairSelector)to[out=0+15, in=270]  (CycleCenter2);
		
	\draw[->] (Selector)-- (PairSelector);
	\draw[->] (Selector)--(SelectorPlus1);
	\draw[->] (CycleCenter1) to[out=90-30, in=180+30] node[above, scale=0.666] {$\frac{1-\eps}{2}$} (CycleNodeX1);
	\draw[->] (CycleCenter1) to[out=270+30, in=180-30] node[below, scale=0.666] {$\frac{1-\eps}{2}$} (CycleNodeY1);
	\draw[->] (CycleCenter1) to[out=90, in=180+30] node[left, scale=0.666] {$\eps$}(UpDown1);
	
	\draw[->] (CycleCenter2) to[out=90+30, in=0-30] node[above, scale=0.666] {$\frac{1-\eps}{2}$} (CycleNodeX2);
	\draw[->] (CycleCenter2)to[out=270-30, in=0+30] node[below, scale=0.666] {$\frac{1-\eps}{2}$} (CycleNodeY2);
	\draw[->] (CycleCenter2) to[out=90, in=0-30] node[right, scale=0.666] {$\eps$} (UpDown2);
	
	\draw[->] (UpDown1)--(UpperSelector1);
	\draw[->] (UpDown1)--(PairSelector1);
	
	\draw[->] (UpDown2)--(Helper);
	\draw[->] (UpDown2)--(PairSelector1);
	
	\draw[->] (UpperSelector1)--(SelectorPlus2);
	
	\draw[->] (Helper) to[out=90+45, in=0] (PairSelectorPlus1);
\end{tikzpicture}
\end{center}
\caption{Level $i$ of the MDP.
Circular vertices are vertices of the player, rectangular vertices are randomization vertices.
Numbers below vertex names, if present, encode priorities $\Omega$.
If a vertex has priority $\priority{v}$, then a reward of $\rew{v}\coloneqq(-N)^{\priority{v}}$ is associated with every edge leaving this vertex.}\label{figure: Level I Of MDP}
\end{figure}

We now discuss how the counter introduced in \Cref{section: Lower Bound Construction} is altered to obtain an MDP $M_n$.
A sketch of level~$i$ of $M_n$ can be found in \Cref{figure: Level I Of MDP}.
First, all player 1 vertices are replaced by randomization vertices.
This is a common technique used for obtaining MDPs that behave similar to given PGs and was used before (e.g., \cite{Fearnley2010Exp,Avis2017}).
While the ideas used in the transformations are similar, there is no standard reduction from PGs to MDPs preserving all properties.

In our construction, all cycle centers $F_{i,j}$ and all vertices $h_{i,j}$ are now randomization vertices.
As vertices of the type $h_{i,j}$ have only one outgoing edge, the probability of this edge is set to 1.
For defining the probabilities of the cycle edges, we introduce a small parameter $\eps>0$ and defer its exact to later.
The idea is to use $\eps$ to make the probabilities of edges $(F_{i,j},s_{i,j})$ very small by setting $p(F_{i,j},s_{i,j})=\eps$ and $p(F_{i,j},d_{i,j,k})=\frac{1-\eps}{2}$ for $k\in\{0,1\}$.
Then, the valuation of $s_{i,j}$ can only contribute significantly to the valuation of $F_{i,j}$ if the cycle center is closed.
If the cycle center is not closed, then the contribution of this vertex can often be neglected.
However, there are situations in which even this very low contribution has a significant impact on the valuation of the cycle center.
For example, if $F_{i,0}$ and $F_{i,1}$ are both open for $\sigma$, then $\valu_{\sigma}^\M(F_{i,0})>\valu_{\sigma}^\M(F_{i,1})$ if and only if $\valu_{\sigma}^\M(s_{i,0})>\valu_{\sigma}^\M(s_{i,1})$.
This sometimes results in a different behavior of the MDP when compared to the PG.
We discuss this later in more detail.

Second, all player 0 vertices remain player 0 vertices. 
Each player 0 vertex is assigned the same priority as in $S_n$.
This priority is now used to define the rewards of the edges leaving a vertex.
More precisely, if we denote the priority of $v\in V_0$ by $\priority{v}$, then we define the reward of any edge leaving $v$ as $\rew{v}\coloneqq(-N)^{\priority{v}}$, where $N\geq 7n$ is a large and fixed parameter.
Note that the reward of an edge thus only depends on its starting vertex.
The reward function that is defined in that way then has the effect that vertices with an even priority are profitable while vertices with an odd priority are not profitable.
In addition, the profitability of a vertex is better (resp.~worse) the higher its priority is.
By choosing a sufficiently large parameter $N$, it is also ensured that rewards are sufficiently separated.
For example, the profitability of some vertex $v$ with even priority cannot be dominated by traversing many vertices with lower but odd priorities.
In principle, this ensures that the MDP behaves very similarly to the PG.

Having introduced the parameter $N$, we now fix the parameter $\eps$ such that $\eps<(N^{2n+11})^{-1}$.
Note that both parameters can be encoded by a polynomial number of bits with respect to the parameter $n$.
By defining the reward of the edge $(t,t)$ as 0, this completely describes the MDP.

We now provide more details on the aspects where the PG and the MDP differ.
One of the main differences between the PG and the MDP are canonical strategies.
Consider a strategy $\sigma$ representing some $\bit\in\bitset_n$, some level $i$ and the two cycle centers $F_{i,0}, F_{i,1}$.
In PG context, both vertices have an even priority and the priority of $F_{i,0}$ is larger than the priority of $F_{i,1}$.
Thus, if both cycle centers escape the level, the valuation of $F_{i,0}$ is better than the valuation of $F_{i,1}$.
Consequently, if $\sigma(g_i)\neq F_{i,0}$, then $(g_i,F_{i,0})$ is improving for $\sigma$.
In some sense, this can be interpreted as the PG ``preferring'' $F_{i,0}$ over $F_{i,1}$.
A similar, but not the same phenomenon occurs in MDP context.
If both cycle centers $F_{i,0}$ and $F_{i,1}$ are in the same ``state'', then the valuation of the two upper selection vertices $s_{i,0}, s_{i,1}$ determines which cycle center has the better valuation.
It turns out that the valuation of $s_{i,\bit_{i+1}}$ is typically better than the valuation of $s_{i,1-\bit_{i+1}}$.
It is in particular not true that the valuation of $s_{i,0}$ is typically better than the valuation of $s_{i,1}$.
Hence, the MDP ``prefers'' vertices $F_{i,\bit_{i+1}}$ over vertices $F_{i,1-\bit_{i+1}}$.
We thus adjust the definition of a canonical strategy in MDP context in the following way.

\begin{definition} \label{definition: Canonical Strategy MDP}
Let $\bit\in\bitset_n$.
A player 0 strategy $\sigma$ for the MDP $M_n$ is called \emph{canonical strategy for $\bit$} if it has the properties defined in \Cref{definition: Canonical Strategy} where Property~3.(d) is replaced by the following:
If none of the cycle centers $F_{i,0}, F_{i,1}$ is closed, then $\sigma(g_i)=F_{i,\bit_{i+1}}$.
\end{definition}

\section{Lower Bound for the Simplex Algorithm and Linear Programs} \label{section: Linear Programs}
%

Following the arguments of \cite{Friedmann2011,Avis2017}, we now discuss how the MDP can be transformed into an LP such that the results obtained for the Policy Iteration Algorithm be transferred to the Simplex Algorithm.
This transformation makes use of the \emph{unichain condition}.
This condition (see \cite{Puterman2005}) states that the Markov Chain obtained from each strategy $\sigma$ has a single irreducible recurrent class.
Unfortunately, the MDP constructed previously does not fulfill the unichain condition.
As we prove in \Cref{lemma: Sink Game}, it however fulfills a \emph{weak} version of the unichain condition.
This weak version states that the \emph{optimal} policy has a single irreducible recurrent class and does not demand this to be true for every strategy.
We later argue why this implies that the same LP which can be obtained by transforming an MDP fulfilling the unichain condition can be used.

We thus return to the discussion for MDPs fulfilling the unichain condition. 
Optimal policies for MDPs fulfilling this condition can be found by solving the following Linear Program:
\begin{align} \tag{P}
\begin{aligned}
\max &\sum_{(u,v)\in E_0} r(u,v)\cdot x(u,v)\\
\text{s.t.} &\sum_{(u,v)\in E}x(u,v)-\sum_{\substack{(v,w)\in E_0\\ (w,u)\in E_R}}p(w,u)\cdot x(v,w)=1 &&\forall u\in V_0\\
&x(u,v)\geq 0 &&\forall(u,v)\in E_0
\end{aligned}
\end{align}

The variable $x(u,v)$ for $(u,v)\in E_0$ represents the probability (or frequency) of using the edge~$(u,v)$.
The constraints of (P) ensure that the probability of entering a vertex $u$ is equal to the probability of exiting $u$.
It is not difficult to see that the basic feasible solutions of~(P) correspond directly to strategies of the MDP, see eg. \cite{Avis2017}.
For each strategy $\sigma$ we can define a feasible setting of the variables $x(u,v)$ with $(u,v)\in E_0$ such that $x(u,v)>0$ only if $\sigma(u)=v$.
Conversely, for every basic feasible solution of (P), we can define a corresponding policy $\sigma$.
It is well-known that the policy corresponding to an optimal basic feasible solution of (P) is an optimal policy for the MDP (see, e.g., \cite{Puterman2005,Avis2017}).

As mentioned, our MDP only fulfills the weak unichain condition.
If the MDP is provided an initial strategy that has the same single irreducible recurrent class as the optimal policy, then the same Linear Program introduced above can be used \cite{FriedmannThesis}.
This follows since all considered basic feasible solutions will have the same irreducible recurrent class by monotonicity.
We refer to~\cite{Tijms2004} for more details.

\section{Lower Bound Proof} \label{section: Lower Bound Proof}

\subsection{The approach and basic definitions}

In this section we outline the proof for the exponential lower bound on the running time of the Strategy Improvement resp. Policy Iteration Algorithm using Zadeh's pivot rule and a strategy-based tie-breaking rule.
We discuss the following key components separately before combining them into our main result.\footnote{Formal proofs can be found in the appendices.}
In the following, we use the notation $G_n$ to simultaneously refer to $S_n$ and $M_n$.
If a statement or definition only holds for either $S_n$ or $M_n$, we explicitly state this.

\begin{enumerate}
	\item We first define an initial strategy $\iota$ such that the pair $(G_n,\iota)$ defines a sink game in PG context resp. has the weak unichain condition in MDP context.
		We also formalize the idea of counting how often an edge has been applied as improving switch.
	\item We then state and discuss the tie-breaking rule.
		Together with the initial strategy, this completely describes the application of the improving switches performed by the Strategy Improvement resp. Policy Iteration Algorithm.
		Further statements, proofs and explanations that are provided in the appendices thus only serve to prove that the algorithms and the tie-breaking rule indeed behave as intended.
	\item We then focus on a single transition from a canonical strategy $\canstrat$ to the next canonical strategy $\sigma_{\bit+1}$. 
		During such a transition, many improving switches need to be applied and thus many intermediate strategies need to be considered.
		These strategies are divided into five \emph{phases}, depending on the configuration of $G_n$ induced by the encountered strategies.
	\item To prove that the tie-breaking rule indeed proceeds along the described phases, we need to specify how often player 0 edges are applied as improving switches, which is formalized by an \emph{occurrence record}.
		We explicitly describe the occurrence records for canonical strategies.
	\item Finally, we combine the previous aspects to prove that applying the respective algorithms with Zadeh's pivot rule and our tie-breaking rule yields an exponential number of iterations.
\end{enumerate}

We begin by providing the initial strategy $\iota$ for $G_n$.
In principle, the initial strategy is a canonical strategy for $0$ in the sense of \Cref{definition: Canonical Strategy} resp. \ref{definition: Canonical Strategy MDP}.

\begin{definition} \label{definition: Initial Strategy}
The initial player $0$ strategy $\iota\colon V_0\mapsto V$ is defined as follows:
\begin{center}
\begin{tabular}{c||c|c|c|c|c|c|c|c|c|c|c|c|}
$v$ 		&$b_i (i<n)$	&$b_n$	&$g_i$		&$d_{i,j,k}$	&$e_{i,j,k}$	&$s_{i,0}$	&$s_{i,1} (i<n)$\\\hline
$\iota(v)$	&$b_{i+1}$		&$t$	&$F_{i,0}$	&$e_{i,j,k}$	&$b_2$			&$h_{i,0}$	&$b_1$ 
\end{tabular}
\end{center}
\end{definition}

We further introduce the notion of a \emph{reachable} strategy.
A strategy $\sigma'$ is reachable from some strategy $\sigma$ if it can be produced by the Strategy Improvement Algorithm starting from $\sigma$ and applying a finite number of improving switches.
Note that the notion of reachability does not depend on the pivot rule or the tie-breaking rule and that every strategy calculated by the Strategy Improvement resp. Policy Iteration Algorithm is reachable by definition.

\begin{definition} \label{definition: Reachable Strategies}
Let $\sigma$ be a player 0 strategy for $G_n$.
The set of all strategies that can be obtained from $\sigma$ by applying an arbitrary sequence of improving switches is denoted by $\reach{\sigma}$.
A strategy $\sigma'$ is \emph{reachable from $\sigma$} if $\sigma'\in\reach{\sigma}$.
\end{definition}

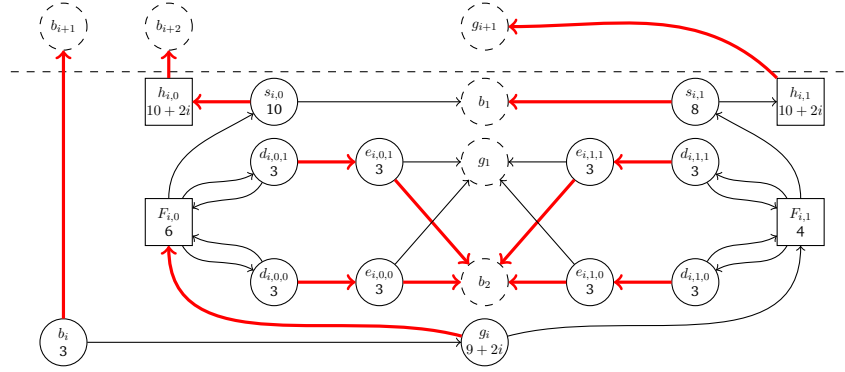
\begin{figure}
\centering
\begin{tikzpicture}[scale=0.8]
	
	\node[zero, dashed] (PairSelector2) at  (0*\x,\y) {$b_2$};
	\node[zero, dashed] (Selector1) at  (0*\x,3*\y) {$g_1$};
	\node[zero, dashed] (PairSelector1) at  (0,4*\y) {$b_1$};
	\node[zero, dashed] (SelectorPlus1) at (-4*\x,5.25*\y) {$b_{i+1}$};
	\node[zero, dashed] (SelectorPlus2) at (-3*\x,5.25*\y) {$b_{i+2}$};
	\node[zero, dashed] (PairSelectorPlus1) at  (0,5.25*\y) {$g_{i+1}$};
	
	\node[zero] (PairSelector) at (0,0) {$g_i$\\$9+2i$};
	\node[one] (CycleCenter1) at (-3*\x,2*\y) {$F_{i,0}$\\6};
	\node[one] (CycleCenter2) at (3*\x,2*\y) {$F_{i,1}$\\4};
	\node[zero] (CycleNodeX1) at (-2*\x,3*\y) {$d_{i,0,1}$\\3};
	\node[zero] (CycleNodeX1T) at (-\x,3*\y) {$e_{i,0,1}$\\3};
	\node[zero] (CycleNodeX2) at (2*\x,3*\y) {$d_{i,1,1}$\\3};
	\node[zero] (CycleNodeX2T) at (1*\x,3*\y) {$e_{i,1,1}$\\3};
	\node[zero] (CycleNodeY1) at (-2*\x,1*\y) {$d_{i,0,0}$\\3};
	\node[zero] (CycleNodeY1T) at (-\x,1*\y) {$e_{i,0,0}$\\3};
	\node[zero] (CycleNodeY2) at (2*\x,1*\y) {$d_{i,1,0}$\\3};
	\node[zero] (CycleNodeY2T) at (\x,1*\y) {$e_{i,1,0}$\\3};
	\node[one] (UpperSelector1) at (-3*\x,4*\y) {$h_{i,0}$\\$10+2i$};
	\node[zero] (Selector) at (-4*\x,0) {$b_i$\\3};
	\node[zero] (UpDown1) at (-2*\x,4*\y) {$s_{i,0}$\\10};
	\node[zero] (UpDown2) at (2*\x,4*\y) {$s_{i,1}$\\8};
	\node[one] (Helper) at (3*\x,4*\y) {$h_{i,1}$\\$10+2i$};

	\tikzstyle{label}=[rectangle, draw, align=center, inner sep=3pt, scale=0.55, fill=gray!33!white]
	
	\draw[dashed] (-4.5*\x,4.5*\y)--(3.5*\x,4.5*\y);
	
	\draw[->] (CycleNodeX1) to [out=270-30, in=0+30] (CycleCenter1);
	\draw[->, red, very thick] (CycleNodeX1)-- (CycleNodeX1T);
	
	\draw[->, red, very thick] (CycleNodeX1T)-- (PairSelector2);
	\draw[->] (CycleNodeX1T)--(Selector1);
	
	\draw[->] (CycleNodeX2) to [out=270+30, in=180-30] (CycleCenter2);
	\draw[->, red, very thick] (CycleNodeX2)-- (CycleNodeX2T);
	
	\draw[->, red, very thick] (CycleNodeX2T)-- (PairSelector2);
	\draw[->] (CycleNodeX2T)--  (Selector1);
	
	\draw[->] (CycleNodeY1) to[out=90+30, in=0-30] (CycleCenter1);
	\draw[->, red, very thick] (CycleNodeY1)--(CycleNodeY1T);
	
	\draw[->, red, very thick] (CycleNodeY1T)-- (PairSelector2);
	\draw[->] (CycleNodeY1T)-- (Selector1);
	
	\draw[->] (CycleNodeY2) to[out=90-30, in=180+30]  (CycleCenter2);
	\draw[->, red, very thick] (CycleNodeY2)-- (CycleNodeY2T);
	
	\draw[->, red, very thick] (CycleNodeY2T)--(PairSelector2);
	\draw[->] (CycleNodeY2T)--(Selector1);
	
	\draw[->, red, very thick] (PairSelector) to[out=180-15, in=270, pos=0.66] (CycleCenter1);
	\draw[->] (PairSelector)to[out=0+15, in=270, pos=0.66]  (CycleCenter2);

	\draw[->] (Selector)-- (PairSelector);
	\draw[->, red, very thick] (Selector)--(SelectorPlus1);
	
	\draw[->] (CycleCenter1) to[out=90-30, in=180+30] (CycleNodeX1);
	\draw[->] (CycleCenter1) to[out=270+30, in=180-30] (CycleNodeY1);
	\draw[->] (CycleCenter1) to[out=90, in=180+30] (UpDown1);
	
	\draw[->] (CycleCenter2) to[out=90+30, in=0-30] (CycleNodeX2);
	\draw[->] (CycleCenter2)to[out=270-30, in=0+30]  (CycleNodeY2);
	\draw[->] (CycleCenter2) to[out=90, in=0-30](UpDown2);
	
	\draw[->, red, very thick] (UpDown1)--(UpperSelector1);
	\draw[->] (UpDown1)--(PairSelector1);
	
	\draw[->] (UpDown2)--(Helper);
	\draw[->, red, very thick] (UpDown2)--(PairSelector1);
	
	\draw[->, red, very thick] (UpperSelector1)--(SelectorPlus2);
	
	\draw[->, red, very thick] (Helper) to[out=90+45, in=0] (PairSelectorPlus1);
\end{tikzpicture}
\caption{The initial strategy $\iota$  (red edges) in level $i$ for $i\in\{1,\dots,n-2\}$.}
\label{figure: Initial Strategy}
\end{figure}

Note that reachability is a transitive property and that we include $\sigma\in\reach{\sigma}$ for convenience.
The $i$-th level of the initial strategy is shown in \Cref{figure: Initial Strategy}.
The initial strategy is chosen such that~$G_n$ and $\iota$ define a sink game in PG context resp. have the weak unichain condition in MDP context.

\begin{restatable}{lemma}{SinkGame} \label{lemma: Sink Game}
For all $n\in\mathbb{N}$, the game $G_n$ and the initial player 0 strategy $\iota$ define a Sink Game with sink~$t$ in PG context, resp. have the weak unichain condition in MDP context.
\end{restatable}

As Zadeh's pivot rule is a memorizing pivot rule, the algorithms needs to maintain information about how often edges have been applied as improving switches.
During the execution of the algorithms, we thus maintain an \emph{occurrence record} $\occrec^{\sigma}\colon E_0\mapsto\mathbb{R}$ that specifies how often an improving switch was applied since the beginning of the algorithms.
Formally, we define $\occrec^{\iota}(e)\coloneqq 0$ for every edge $e\in E_0$, i.e., the occurrence record with respect to the initial strategy is equal to 0.
Then, whenever the algorithms applies an edge $e$, the occurrence record of $e$ is increased by 1.

\subsection{The Tie-Breaking Rule} \label{section: Tie Breaking Rule}

We now discuss the tie-breaking rule.
It specifies which edge to if there are multiple improving switches that minimize the occurrence record for the current strategy.
The tie-breaking rule is in principle implemented as an ordering of the set $E_0$ and depends on the current strategy~$\sigma$ as well as the occurrence records.
Whenever the algorithms has to break ties, it then chooses the first edge according to this ordering.
There is, however, one exception from this behavior.
During the transition from the canonical strategy representing 1 towards the canonical strategy representing~2, one improving switch~$e$ (which we do not specify yet) has to be applied earlier than during other transitions.
The reason is that the occurrence records of several edges, including $e$, are equal to zero at this point in time.
Thus, the unmodified tie-breaking rule decides which of these switches is applied and would not choose to apply $e$.
Since the algorithms produces an unwanted behavior if $e$ is not applied at this specific point in time, we need to handle this situation explicitly.
Fortunately, in later iterations, $e$ turns out to be the unique improving switch minimizing the occurrence record whenever it has to be applied, so this special treatment is not necessary later.

It turns out that it is not necessary to give a complete ordering of $E^0$.
In fact, it is sufficient to describe a pre-order of $E^0$ as any linear extension of this pre-order can be used.

Let $\sigma$ be a player 0 strategy for $G_n$.
Henceforth, we use the symbol $*$ as a wildcard.
More precisely, when using the symbol~$*$, this means any suitable index or vertex (depending on the context) can be inserted for~$*$ such that the corresponding edge exists.
For example, the set $\{(e_{*,*,*},*)\}$ would then denote the set of all edges starting in escape vertices. 
Using this notation, we define the following sets of edges.
\begin{itemize}
	\item $\G\coloneqq\{(g_{i},F_{i,*})\}$ is the set of all edges leaving selector vertices.
	\item $\E^0\coloneqq\{(e_{i,j,k},*)\colon \sigma(d_{i,j,k})\neq F_{i,j}\}$  is the set of edges leaving escape vertices whose cycle vertices do not point towards their cycle center.
			Similarly, $\E^1\coloneqq \{(e_{i,j,k},*)\colon \sigma(d_{i,j,k})=F_{i,j}\}$ is the set of edges leaving escape vertices whose cycle vertices point towards their cycle center.
	\item $\D^1\coloneqq\{(d_{*,*,*},F_{*,*})\}$  is the set of cycle edges an $\D^0\coloneqq\{(d_{*,*,*},e_{*,*,*})\}$ is the set of the other edges leaving cycle vertices. 
	\item $\B^0\coloneqq\bigcup_{i=1}^{n-1}\{(b_i,b_{i+1})\}\cup\{(b_n,t)\}$ is the set of all edges between entry vertices.	
		The set $\B^1\coloneqq\{(b_*,g_*)\}$ of all edges leaving entry vertices and entering selection vertices is defined analogously and $\B\coloneqq\B^0\cup\B^1$ is the set of all edges leaving entry vertices.
	\item $\S\coloneqq\{(s_{*,*},*)\}$ is the set of all edges leaving upper selection vertices.
\end{itemize}

We next define two pre-orders based on these sets.
However, we need to define finer pre-orders for the sets $\E^0,\E^1,\S$ and $\D^1$ first.

Informally, the pre-order on $\E^0$ forces the algorithms to favor switches of higher levels and to favor $(e_{i,0,k},*)$ over $(e_{i,1,k},*)$ in PG context and $(e_{i,\indbit^{\sigma}_{i+1},k},*)$ over $(e_{i,1-\indbit^{\sigma}_{i+1},k},*)$ in MDP context.
For a formal description let $(e_{i,j,x},*),(e_{k,l,y},*)\in\E^0$.
In PG context, we define $(e_{i,j,x},*)\prec_\sigma (e_{k,l,y},*)$ if either $i>k$, or $i=k$ and $j<l$.
In MDP context, we define $(e_{i,j,x},*)\prec_\sigma (e_{k,l,y},*)$ if either $i>k$, or $i=k$ and $j=\indbit^{\sigma}_{i+1}$.

Similarly, the pre-order on $\S$ also forces the algorithm to favor switches of higher levels.
Thus, for $(s_{i,j},*),(s_{k,l},*)\in\S$, we define $(s_{i,j},*)\prec_{\sigma} (s_{k,l},*)$ if $i>k$.

We now describe the pre-order for $\E^1$.
Let $(e_{i,j,x},*),(e_{k,l,y},*)\in \E^1$.
\begin{enumerate}
	\item The first criterion encodes that switches contained in higher levels are applied first.
		Thus, if $i>k$, then $(e_{i,j,x},*)\prec_{\sigma}(e_{k,l,y},*)$.
	\item If $i=k$, then we consider the states of the cycle centers $F_{i,j}$ and $F_{k,l}=F_{i,1-j}$.
		If exactly one cycle center of level $i$ is closed, then the improving switches within this cycle center are applied first.
	\item Consider the case where $i=k$ but no cycle center of level $i$ is closed.
		Let $t^{\rightarrow}\coloneqq b_2$ if $\nsb>1$ and $t^{\rightarrow}\coloneqq g_1$ if $\nsb=1$.
		If there is exactly one halfopen cycle center escaping to $t^{\rightarrow}$ in level $i$, then switches within this cycle center have to be applied first.
	\item Assume that none of the prior criteria applied.
		This includes the case where both cycle centers are in the same state, and $i=k$ holds in this case.
		Then, the order of application depends on whether we consider PG or MDP context.
		In PG context, improving switches within $F_{i,0}$ are applied first.
		In MDP context, improving switches within $F_{i,\indbit^{\sigma}_{i+1}}$ are applied first.
\end{enumerate}

We next give a pre-order for $\D_1$.
Let $(d_{i,j,x},F_{i,j}),(d_{k,l,y},F_{k,l})\in\D^1$.
\begin{enumerate}
	\item The first criterion states that improving switches that are part of open cycles are applied first.
		We thus define $(d_{i,j,x},F_{i,j})\prec_\sigma (d_{k,l,y},F_{k,l})$ if $\sigma(d_{k,l,1-y})=F_{k,l}$ but $\sigma(d_{i,j,1-x})\neq F_{i,j}$.
	\item The second criterion states the following.
		Among all halfopen cycles, improving switches contained in cycle centers such that the bit of the level the cycle center is part of is equal to zero are applied first. 
		If the first criterion does not apply, we thus define $(d_{i,j,x},F_{i,j})\prec_\sigma (d_{k,l,y},F_{k,l})$ if $\indbit_k^{\sigma}>\indbit_i^{\sigma}$.
	\item The third criterion states that among all partially closed cycles, improving switches inside cycles contained in lower levels are applied first.
		If none of the first two criteria apply, we thus define $(d_{i,j,x},F_{i,j})\prec_\sigma (d_{k,l,y},F_{k,l})$ if $k>i$.
	\item The fourth criterion states that improving switches within the active cycle center are applied first within one level.
		If none of the previous criteria apply, we thus define $(d_{i,j,x},F_{i,j})\prec_\sigma (d_{k,l,y},F_{k,l})$ if $\indbit^\sigma_{k+1}\neq l$ and $\indbit^{\sigma}_{i+1}= j$.
	\item The last criterion states that edges with last index equal to zero are preferred within one cycle center.
		That is, if none of the previous criteria apply, we define $(d_{i,j,x},F_{i,j})\prec_{\sigma}(d_{k,l,y},F_{k,l})$ if $x<y$.
		If this criterion does not apply either, the edges are incomparable.
\end{enumerate}

We now define the pre-order $\prec_{\sigma}$ and the-breaking rule, implemented by an ordering of $E_0$.

\begin{definition} \label{definition: Tie-Breaking exponential}
Let $\sigma$ be a player 0 strategy for $G_n$ and $\occrec^{\sigma}\colon E^0\to\mathbb{N}_0$ be an occurrence record. 
We define the pre-order $\prec_\sigma$ on $E_0$ by defining the set-based pre-order 
\[\G\prec_\sigma\D^0\prec_\sigma\E^1\prec_\sigma\B\prec_\sigma\S\prec_\sigma\E^0\prec_\sigma\D^1\]where the sets $\E^0,\E^1,\S$ and $\D^1$ are additionally pre-ordered as described before.
We extend the pre-order to an arbitrary but fixed total ordering and denote the corresponding order also by $\prec_\sigma$.
We define the following tie-breaking rule:
\emph{Let $I_\sigma^{\min}$ denote the set of improving switches with respect to $\sigma$ that minimize the occurrence record.
Apply the first improving switch contained in~$I_\sigma^{\min}$ with respect to the ordering $\prec_\sigma$ with the following exception:
If $\occrec^{\sigma}(b_1,b_2)=\occrec^{\sigma}(s_{1,1},h_{1,1})=0$, then apply $(s_{1,1},h_{1,1})$ instead of $(b_1,b_2)$.}
\end{definition}

We usually just use the notation $\prec$ to denote the ordering if it is clear from the context which strategy is considered and whether $\prec$ is defined via the default or special pre-order.

\begin{restatable}{lemma}{TieBreakingPolynomial} \label{lemma: Tie Breaking Polynomial}
Given a strategy $\sigma\in\reach{\iota}$ and an occurrence record $\occrec^{\sigma}\colon E^0\to\mathbb{N}_0$, the tie-breaking rule can be evaluated in polynomial time.
\end{restatable}

\subsection{The Phases of a Transition and the Application of Improving Switches} \label{section: Phases}

As explained earlier, the goal is to prove that Zadeh's pivot rule and tour tie-breaking rule enumerates at least one strategy per number $\bit\in\bitset_n$.
This is proven in an inductive fashion.
That, is we prove that given a canonical strategy $\canstrat$ for $\bit\in\bitset_n$, the algorithms eventually calculate a canonical strategy $\sigma_{\bit+1}$ for $\bit+1$.
This process is called a \emph{transition} and each transition is partitioned into up to five phases.
In each phase, a different ``task'' is performed in order to obtain the strategy $\sigma_{\bit+1}$.
These tasks are, for example, the opening and closing of cycle centers, updating the escape vertices or adjusting some of the selection vertices.

\begin{figure}[ht]
\centering
\begin{tikzpicture}
\def\x{2.3}
\def\y{0.7}
\draw[white] (-0*\x,-2.4*\y) rectangle (5.75*\x, 2.4*\y);
\tikzset{every node/.style={draw, rectangle, minimum size=1cm, align=center, scale=0.75}}
\node (pib) at (0,0) {$\canstrat$};
\node (pib+1) at (5.75*\x,0) {$\sigma_{\bit+1}$};

\node[dashed, scale=0.8] (p1t) at (1*\x, 1.75*\y) {Balance \\Occurrence Records};
\node[dashed, scale=0.8] (p1m) at (1*\x,0) {Close $F_{\nsb, \bit_{\nsb+1}}$};
\node[dashed, scale=0.8, fill=gray!30!white] (p1b) at (1*\x,-1.75*\y) {Update \\selector vertices};
\node[dashed, scale=0.8] (p2t) at (2.25*\x,0.75*\y) {Update \\$b_i$ for $i<\nsb$};
\node[dashed, scale=0.8] (p2b) at (2.25*\x,-0.75*\y) {Update \\upper selection vertices};
\node[dashed, scale=0.8] (p3t) at (3.5*\x,0.75*\y) {Open cycle centers};
\node[dashed, scale=0.8] (p3tt) at (3.5*\x, 2.25*\y) {Update $b_{\nsb}$};
\node[dashed, scale=0.8] (p3m) at (3.5*\x,-0.75*\y) {Update \\escape vertices};
\node[dashed, scale=0.8, fill=yellow!30!white] (p3b) at (3.5*\x,-2.25*\y) {Update \\ upper selection vertices};
\node[dashed, scale=0.8] (p5t) at (4.75*\x,1.75*\y) {Balance \\occurrence records};
\node[dashed, scale=0.8] (p5m) at (4.75*\x, 0) {Update \\escape vertices};
\node[dashed, scale=0.8, fill=gray!30!white] (p5b) at (4.75*\x, -1.75*\y) {Update \\selector vertices};

\draw[->, thick] (pib) to[out=45, in=180] (p1t);

\draw[->, thick] (p1t)--(p1m);
\draw[<->, thick, gray] (p1t) to[out={180+22.5}, in={180-22.5}] (p1b); 

\draw[->, thick, blue] (p1m)--++(0.6125*\x,0)--++(0,0.75*\y)--(p2t);
\draw[->, thick, red] (p1m)--++(0.6125*\x,0)--++(0,1.5*\y)--++(2*0.65*\x,0)--++(0,-1.5*\y)--(p3m);
\draw[thick, black] (p1m)--++(0.6125*\x,0)--++(0,\y);
\draw[<->, thick, blue] (p2t)--(p2b);

\draw[->, thick,blue] (p2b)--(p3m); 
\draw[<->, thick, blue] (p3m)--(p3b);
\draw[<->, thick] (p3t)--(p3m);

\draw[->, thick] (p3tt)--(p5t);
\draw[->, red, thick] (p3t)--(p3tt);
\draw[<->, thick] (p5t)--(p5m);
\draw[<->, thick, gray] (p5m)--(p5b);
\draw[->, thick, blue] (p3b)--++(0.625*\x,0)--++(0,4*\y)--(p3tt);

\draw[->, thick] (p5t) to[out=-15, in={180-30}](pib+1);
\draw[->, thick] (p5m)--(pib+1);
\draw[->, thick, gray] (p5b) to[out=15, in={180+30}] (pib+1);

\draw[<-, thick] (p1m)--(p1b);
\draw[<->, thick, gray] (p5t) to[out={180+22.5}, in={180-22.5}] (p5b); 
\draw[<->, thick] (p5m)--(p5b);
\draw[->, thick, gray] (p5b) to[out=15, in={180+30}] (pib+1);
\end{tikzpicture}
\caption[Sketch of the phases performed by the algorithm in $M_n$]{Sketch of some of the tasks performed during the different phases by the strategy improvement algorithm in $M_n$.
Each box marks one task that has to be fulfilled, and each vertical set of boxes corresponds to one of the phases $1,2,3$ and $5$ (from left to right).
Gray boxes and edges represent tasks that are not performed during all transitions.
Red/blue edges represent that the corresponding boxes are only relevant if $\bit+1$ is odd/even.
The yellow box represents the task that is phase $4$ in $S_n$ but part of phase $3$ in $M_n$.}
\label{figure: Sketch Of Phases}
\end{figure}
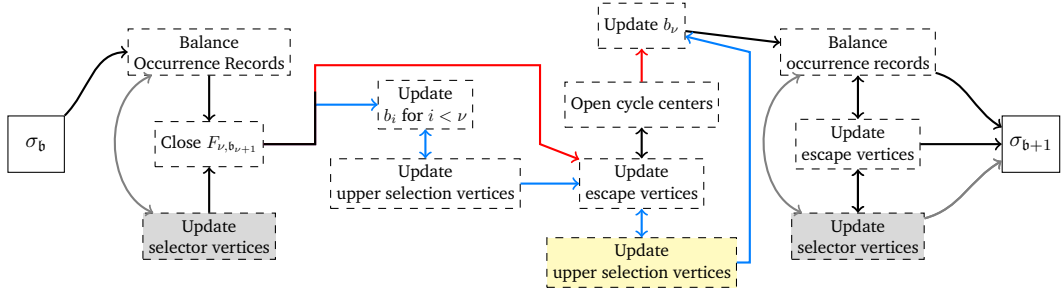

Depending on whether we consider PG or MDP context and $\nsb(\bit+1)$, there can be 3,4 or 5 different phases.
Phases 1,3 and 5 always take place while Phase 2 only occurs if $\nsb(\bit+1)>1$, as it updates the target vertices of some selection vertices $s_{i,j}$ with $i<\nsb(\bit+1)$.
The same holds for Phase 4, although this phase only exists in PG context.
In MDP context, we apply the corresponding switches already in Phase 3 and there is no separate Phase 4.
A very simplified and schematic sketch of the different phases and their interaction in the Markov decision process $M_n$ is given in \Cref{figure: Sketch Of Phases}.
Consider the canonical strategy $\canstrat$ for some $\bit\in\bitset_n$ and let $\nsb\coloneqq\ell(\bit+1)$.

We now give a detailed description of the phases.
\begin{enumerate}
	\item During phase $1$, several cycle vertices switch towards their cycle centers.
		The primary purpose of this phase is that the strategy $\sigma$ obtained after the application of the final switch represents $\bit+1$.
		A secondary purpose is to balance the occurrence records of the cycle edges by applying additional improving switches.
		This application might close inactive cycle centers $F_{i,1-\bit_{i+1}}$ and can thus make edges $(g_i,F_{i,1-\bit_{i+1}})$ improving.
		As balancing occurrence records might close additional cycle centers, more edges of the type $(g_*,F_{*,*})$ can become improving and are then applied.	
		The final improving switch applied during phase~$1$ closes the cycle center $F_{\nsb,(\bit+1)_{\nsb+1}}$.
		Depending on the parity of $\bit+1$, either phase $2$ (if $\bit+1\bmod2=0$) or $3$ (if $\bit+1\bmod2=1$) begins.
		
	\item In phase $2$, the upper selection vertices $s_{i,j}$ for $i\in[\nsb-1]$ and $j=(\bit+1)_{i+1}$ change their targets from $b_1$ to $h_{i,j}$.
		This is necessary as the induced bit state is now representing $\bit+1$, so the $\nsb$ least significant bits changed.
		Furthermore, the entry vertices $b_i$ of these levels switch towards $b_{i+1}$ (with the exception of $b_1$).
		Since $\bit_i=\bit_{i+1}$ for all $i\neq 1$ if $\nsb=1$, these operations only need to be performed if $\nsb>1$.
		
	\item Phase $3$ is partly responsible for applying improving switches involving escape vertices.
		Since the parity of $\bit$ and $\bit+1$ is not the same,  all escape vertices have to change their targets.
		During phase $3$, exactly the escape vertices $e_{i,j,k}$ whose cycle vertex $d_{i,j,k}$ points to the cycle center~$F_{i,j}$ change their targets.
		In addition, exactly these cycle vertices then also change their targets to $e_{i,j,k}$ unless the cycle center~$F_{i,j}$ is closed and active.
		This enables the application of $(d_{i,j,k},F_{i,j})$ which is necessary to balance the occurrence record of the cycle edges.

		At the end of this phase either $(b_1,g_1)$ (if $\nsb=1$) or $(b_1,b_2)$ (if $\nsb>1$) is applied.
		In~$M_n$, the improving switches of phase $4$ are also applied during phase $3$.
		
	\item During phase $4$, the remaining upper selection vertices $s_{i,j}$ with $ i\in[\nsb-1]$ and $j=1-(\bit+1)_{i+1}$ are updated by changing their targets to $b_1$.
		These updates are necessary to allow the cycle centers and cycle vertices to access the spinal path.
		Similarly to phase $2$, these switches are only applied if $\nsb>1$.
	
	\item During phase $5$, the remaining improving switches involving escape vertices are applied.
		Moreover, some of the cycle edges $(d_{i,j,k},F_{i,j})$ that have a very low occurrence record are also applied in order to increase their occurrence record.
		In some sense, the switches \enquote{catch up} to the other edges that have been applied more often.
		This application might close some inactive cycle centers $F_{i,j}$ and consequently make the corresponding edge $(g_i,F_{i,j})$ improving.
		This switch is then also applied.
		Phase $5$ ends once the set of improving switches only contains edges of the type $(d_{i,j,k},F_{i,j})$. 
\end{enumerate}
To give the formal definition of the phases, we require additional notation for describing strategies.
In particular, we encode the choices of~$\sigma$ by using integers.
For this purpose, we introduce a function $\sigmabar$ that will be extended  later to describe more complex configurations of $G_n$.
At this point, we only provide the first layer of complexity by defining $\sigmabar(v)$ for all $v\in V_0$ using \Cref{table: Sigmabar for vertices}.
In principle, $\sigmabar$ is used to abbreviate boolean expressions. 
These expressions are either true (i.e., equal to $1$) or false (i.e., equal to $0$).
For example, $\sigmabar(b_i)$ denotes the boolean expression $\sigma(b_i)=g_i$, so $\sigmabar(b_i)=1$ if and only if $\sigma(b_i)=g_i$.

\begin{table}[ht] 
\centering
\begin{tabular}{c||c|c|c|c|c|c}
Symbol 							&$\sigmabar(b_i)$	&$\sigmabar(s_{i,j})$			&$\sigmabar(g_i)$				&$\sigmabar(d_{i,j,k})$			&$\sigmabar(e_{i,j,k})$\\\hline
Encoded expression	&$\sigma(b_i)=g_i$	&$\sigma(s_{i,j})=h_{i,j}$	&$\sigma(g_i)=F_{i,1}$	&$\sigma(d_{i,j,k})=F_{i,j}$	&$\sigma(e_{i,j,k})=b_2$
\end{tabular}
\caption[The function $\sigmabar$ in the exponential construction.]{Definition of the function $\sigmabar$ for the player vertices and a strategy $\sigma$ in $G_n$.}\label{table: Sigmabar for vertices}
\end{table}

For convenience, we define $\sigmabar(t)\coloneqq 0$.
Since every player vertex has outdegree of at most two, the value of $\sigmabar(v)$ is in bijection to $\sigma(v)$.
We can thus in particular use $\nsigmabar(v)$ to denote $\sigmabar(v)=0$.
For convenience of notation, the precedence level of \enquote{$=$} and \enquote{$\neq$} is higher than the precedence level of $\wedge$ and $\vee$.
That is, an expression $x\wedge y=z$ is interpreted as $x\wedge (y=z)$.

Using this notation, we now introduce a strategy-based parameter $\relbit{\sigma}\in[n+1]$.
This parameter is called \emph{the next relevant bit} of the strategy $\sigma$.
Before defining this parameter formally, we briefly explain its importance and how it can be interpreted.

One of the central concepts of $G_n$ is that the two cycle centers of a fixed level alternate in representing bit $i$.
Consequently, the selector vertex $g_i$ of level $i$ needs to select the correct cycle center. 
Moreover, $b_i$ should point to $g_i$ if and only if bit $i$ is equal to $1$ (see \Cref{definition: Canonical Strategy} resp. \ref{definition: Canonical Strategy MDP}).
This in particular implies that the selector vertex $g_{i-1}$ of level $i-1$ needs to be in accordance with the entry vertex of level $b_i$ if bit $i-1$ is equal to $1$.
More precisely, it should not happen that $\sigma(b_i)=g_i, \sigma(b_{i+1})=g_{i+1}$ and $\sigma(g_i)=F_{i,0}$.
However, it cannot be guaranteed that this does not happen for some intermediate strategies encountered during $\canstrat\to\sigma_{\bit+1}$.
Such a configuration is then an indicator that \emph{some} operations have to be performed in the levels $i$ and $i+1$.
This is captured by the parameter $\relbit{\sigma}$ as it is defined as the lowest level higher than any level that is set \enquote{incorrectly} in that sense.
If there are no such levels, then $\relbit{\sigma}$ is the lowest level $i$ with $\sigma(b_i)=b_{i+1}$.
The parameter can thus be interpreted as an indicator encoding where \enquote{work needs to be done next}.
Formally, it is defined as follows.

\begin{definition}[Next relevant bit] \label{definition: Next relevant bit}
Let $\sigma\in\reach{\sigma_0}$.
The \emph{set of incorrect levels}\index{incorrect level} is defined as $\incorrect{\sigma}\coloneqq\{i\in[n]:\sigmabar(b_i)\wedge\sigmabar(g_i)\neq\sigmabar(b_{i+1})\}$.
The \emph{next relevant bit}\index{next relevant bit} $\relbit{\sigma}$ of the strategy $\sigma$ is \[\relbit{\sigma}\coloneqq\begin{cases}
	\min \{i>\max\{i'\in\incorrect{\sigma}\}:\sigmabar(b_i)\wedge\sigmabar(g_i)=\sigmabar(b_{i+1})\}\cup\{n\}, &\text{if }\incorrect{\sigma}\neq\emptyset,\\
	\min \{i\in[n+1]:\sigma(b_i)=b_{i+1}\}, &\text{if }\incorrect{\sigma}=\emptyset.
\end{cases}\]
\end{definition}

The next relevant bit now enables us to give a formal definition of the phases.
These phases are described by a set of \emph{properties} that the strategies of the corresponding phase have to fulfill.
Before we list and explain the properties that are used for defining the phases, we extend the function $\sigmabar$.
The definition given in \Cref{table: Sigmabar for vertices} enables us to use the function $\sigmabar$ to describe the state of individual vertices.
It is however convenient to also describe more complex configurations by encoding them as boolean expressions.
An example for such a configurations is the setting of the cycle centers.
We thus extend the notation of $\sigmabar$ and refer to \Cref{table: Sigmabar complex} for an overview over the complete definition of the function.

\begin{table}[ht]
\centering
\begin{tabular}{|c|c|}\hline
Symbol									&Encoded expression		\\\hline
$\sigmabar(b_i)$			&$\sigma(b_i)=g_i$\\
$\sigmabar(s_{i,j})$		&$\sigma(s_{i,j})=h_{i,j}$\\
$\sigmabar(g_i)$			&$\sigma(g_i)=F_{i,1}$\\
$\sigmabar(d_{i,j,k})$	&$\sigma(d_{i,j,k})=F_{i,j}$\\
$\sigmabar(e_{i,j,k})$	&$\sigma(e_{i,j,k})=b_2$\\\hline
\end{tabular}
\hspace*{1em}
\begin{tabular}{|c|c|}\hline
Symbol									&Encoded expression		\\\hline
$\sigmabar(s_i)$			&$\sigmabar(s_{i\sigmabar(g_i)})$\\
$\sigmabar(d_{i,j})$		&$\sigmabar(d_{i,j,0})\wedge\sigmabar(d_{i,j,1})$\\
$\sigmabar(d_i)$			&$\sigmabar(d_{i,\sigmabar(g_i)})$\\
$\sigmabar(eg_{i,j})$	&$\bigvee_{k\in\{0,1\}}[\nsigmabar(d_{i,j,k})\wedge\nsigmabar(e_{i,j,k})]$\\
$\sigmabar(eb_{i,j})$	&$\bigvee_{k\in\{0,1\}}[\nsigmabar(d_{i,j,k})\wedge\sigmabar(e_{i,j,k})]$\\
$\sigmabar(eg_i)$			&$\sigmabar(eg_{i,\sigmabar(g_i)})$\\
$\sigmabar(eb_i)$			&$\sigmabar(eb_{i,\sigmabar(g_i)})$\\\hline
\end{tabular}
\caption[Extension and full definition of the function $\sigmabar$.]{Extension and full definition of the function $\sigmabar$ given in \Cref{table: Sigmabar for vertices} to describe more complex configurations.
Here, $\nsigmabar(v)$ is the logical negation of $\sigmabar(v)$.} \label{table: Sigmabar complex}
\end{table}

Formally, a strategy belongs to one of the five phases if it has a certain set of properties.
These properties can be partitioned into several categories depending on the vertices or terms that are involved.
The properties might also depend on one or more parameters like a level or a cycle center.

Consider some fixed strategy $\sigma, \bit\in\bitset_n$ and let $\nsb\coloneqq\ell(\bit+1)$ denote the least significant set bit of $\bit+1$.
The first three properties are related to the \textbf{E}ntry \textbf{V}ertices.
\Pref{EV1}$_i$ states that the entry vertex of level $i$ should point to $g_i$ if and only if the the active (with respect to the induced bit state) cycle center $F_{i,,\indbit^{\sigma}_{i+1}}$ is closed, so
\begin{property}{EV1}\sigmabar(b_i)=\sigmabar(d_{i,\indbit^{\sigma}_{i+1}}).\end{property}

Similarly, \Pref{EV2}$_i$ states that $\sigma(b_i)=g_i$ implies that the selector vertex $g_i$ of level~$i$ should point to the corresponding active cycle center, so
\begin{property}{EV2}\sigmabar(b_i)\implies\sigmabar(g_i)=\indbit^{\sigma}_{i+1}.\end{property}

\Pref{EV3}$_i$ states that $\sigma(b_i)=g_i$ implies that the inactive cycle center is not closed, so
\begin{property}{EV3}\sigmabar(b_i)\implies\nsigmabar(d_{i,1-\indbit^{\sigma}_{i+1}}).\end{property}

This property is a good example for a property that will be violated during specific phases as several inactive cycle centers will be closed when the induced bit state changes from $\bit$ to $\bit+1$.

The next five properties are all related to the \textbf{ESC}ape vertices $e_{*,*,*}$.
\Pref{ESC1} states that the escape vertices are set \enquote{correctly}, that is, as they should be set for a canonical strategy representing $\indbit$, so
\begin{property}{ESC1}[\indbit^{\sigma}_1=0\implies \sigma(e_{*,*,*})=b_2]\wedge[\indbit^{\sigma}_1=1\implies\sigma(e_{*,*,*})=g_1].\end{property}
\Pref{ESC2} states that all escape vertices point to $g_1$, so
\begin{property}{ESC2}\sigma(e_{*,*,*})=g_1.\end{property}
Although this property seems redundant due to \Pref{ESC1}, it is crucial for properly defining the second phase.

The next three properties are used to describe the access of $F_{i,j}$ to the vertices $g_1$ and~$b_2$ via the escape vertices.
More precisely, they state whether $F_{i,j}$ has access to only $g_1$ (\Pref{ESC3}$_{i,j}$), only $b_2$ (\Pref{ESC4}$_{i,j}$) or to both of these vertices (\Pref{ESC5}$_{i,j}$).
We mention here that \Pref{ESC3} is technically not used for the definition of the phases, but as it will be used within several proofs and fits the other properties related to the escape vertices, we already provide it here.
Formally, these properties are given via
\begin{property}{ESC3}\sigmabar(eg_{i,j})\wedge\nsigmabar(eb_{i,j}),\end{property}\vspace*{-1em}
\begin{property}{ESC4}\sigmabar(eb_{i,j})\wedge\nsigmabar(eg_{i,j}),\end{property}\vspace*{-1em}
\begin{property}{ESC5}\sigmabar(eb_{i,j})\wedge\sigmabar(eg_{i,j}).\end{property}

The next  properties are concerned with the \textbf{U}pper \textbf{S}election \textbf{V}ertices $s_{*,*}$.
\Pref{USV1}$_i$ states that both upper selection vertices of level $i$ are set \enquote{correctly} with respect to the induced bit state, while \Pref{USV3}$_i$ states that both of these vertices are set incorrectly.
\Pref{USV2}$_{i,j}$ simply states $\sigma(s_{i,j})=h_{i,j}$ and will be used to identify strategies for which the upper selection vertices of lower levels need to be updated since the induced bit state changed.
\begin{property}{USV1}\sigma(s_{i,\indbit^{\sigma}_{i+1}})=h_{i,\indbit^{\sigma}_{i+1}}\wedge\sigma(s_{i,1-\indbit^{\sigma}_{i+1}})=b_1\end{property}\vspace*{-1em}
\begin{property}{USV2}\sigma(s_{i,j})=h_{i,j}\end{property}\vspace*{-1em}
\begin{property}{USV3}\sigma(s_{i,\indbit^{\sigma}_{i+1}})=b_1\wedge\sigma(s_{i,1-\indbit^{\sigma}_{i+1}})=h_{i,1-\indbit^{\sigma}_{i+1}}\end{property}

The next two properties are related to the \textbf{C}ycle \textbf{C}enters.
\Pref{CC1}$_i$ states that at least one cycle center of level $i$ has to be open or halfopen if $i<\relbit{\sigma}$.
\Pref{CC2} states that the active cycle center of level $\nsb=\ell(\bit+1)$ is closed and that the selector vertex of level $\nsb$ chooses the correct cycle center with respect to $\bit+1$, so
\begin{property}{CC1}i<\relbit{\sigma}\implies\nsigmabar(d_{i,0})\vee\nsigmabar(d_{i,1}),\end{property}\vspace*{-1em}
\begin{property}{CC2}\sigmabar(d_{\nsb})\wedge\sigmabar(g_{\nsb})=(\bit+1)_{\nsb+1}\end{property}

The following two properties are related to the \textbf{S}elector \textbf{V}ertices and are unique for either the \textbf{M}arkov decision process $M_n$  (\Pref{SVM})or the sink \textbf{G}ame $S_n$ (\Pref{SVG}).
They are related to the setting of selector vertices if the represented bit is equal to $0$.
According to \Cref{definition: Canonical Strategy} resp. \ref{definition: Canonical Strategy MDP}, the cycle center chosen by $g_i$ is fixed in this case, see condition~3.(d).
The two properties (\ref{property: SVM}) and (\ref{property: SVG}) now state that the selector vertex can only choose the other cycle center $F_{i,1}$ resp. $F_{i,1-\indbit^{\sigma}_{i+1}}$ if this cycle center is closed, so
\begin{property}{SVM}\sigmabar(g_i)=1-\indbit^{\sigma}_{i+1}\implies\sigmabar(d_{i,1-\indbit^{\sigma}_{i+1}}),\end{property}\vspace*{-1em}
\begin{property}{SVG}\sigmabar(g_i)=1\implies\sigmabar(d_{i,1}).\end{property}

The final two properties are related to the next \textbf{REL}evant bit $\relbit{\sigma}$ defined in \Cref{definition: Next relevant bit}.
\Pref{REL1} states that the set of incorrect levels is empty.
This in particular implies that $\relbit{\sigma}=\min\{i\in[n+1]:\sigma(b_i)=b_{i+1}\}$.
\Pref{REL2} states that this parameter is equal to the least significant set bit of the bit state induced by $\sigma$, so
\begin{property}{REL1}\nexists i:\sigma(b_{i-1})=g_{i-1}\wedge\sigmabar(b_i)\neq\sigmabar(g_{i-1}),\end{property}\vspace*{-1em}
\begin{property}{REL2}\relbit{\sigma}=\ell(\indbit^{\sigma}).\end{property}

Together with the induced bit state, these properties are now used to formally define the phases.
This is done by providing a table where each row corresponds to one of the properties and each column to one phase.
In addition, there are some special conditions that have to be fulfilled during some phases that cannot be phrased as a simple property.

\begin{definition}[Phase-$k$-strategy] \label{definition: Phase-$k$-strategy}
Let $\sigma$ be a strategy for $G_n$ and $k\in[5]$.
Then, $\sigma$ is a \emph{phase-$k$-strategy}\index{phase-$k$-strategy} if it has the properties of the corresponding column of \Cref{table: Definition of Phases} and fulfills the corresponding special conditions.
\end{definition}

\begin{table}
\footnotesize
\centering
\begin{tabular}{|c||c|c|c|c|c|}\hline%
Property						&Phase 1	&Phase 2				&Phase 3												&Phase 4										&Phase 5 \\\hline\hline
(\ref{property: EV1})$_{i}$	&$i\in[n]$	&$i>\relbit{\sigma}$	&$i>1$													&$i\in[n]$										&$i\in[n]$\\
(\ref{property: EV2})$_{i}$	&$i\in[n]$	&$i\geq\relbit{\sigma}$	&$i>1$													&$i\in[n]$										&$i\in[n]$\\	
(\ref{property: EV3})$_{i}$	&$i\in[n]\setminus\{\nsb\}$	&$i>\relbit{\sigma}$	&$i>1,i\neq\relbit{\sigma}$								&$i\in[n]$										&$i\in[n]$\\				
(\ref{property: USV1})$_{i}$	&$i\in[n]$	&$i\geq\relbit{\sigma}$	&$i\geq\relbit{\sigma}$									&$i\geq\nsb$									&$i\in[n]$\\
(\ref{property: USV2})$_{i,j}$	&-			&$(i,1-\indbit_{i+1})\colon i<\relbit{\sigma}$				&$(i,*)\colon i<\relbit{\sigma}$	&$(i,\indbit_{i+1})\colon i<\nsb$		&-\\
(\ref{property: ESC1})			&True		&-						&-														&-												&False*\\
(\ref{property: ESC2})			&-			&True					&-														&-												&-\\
(\ref{property: ESC4})$_{i,j}$	&-			&-						&-														&$(i,j)\in S_1$											&- \\
(\ref{property: ESC5})$_{i,j}$	&-			&-						&-														&$(i,j)\in S_2$											&-\\
(\ref{property: REL1})			&True		&-						&-														&True											&True\\
(\ref{property: REL2})			&-			&True					&True													&False											&False\\
(\ref{property: CC1})$_{i}$	&$i\in[n]$	&$i\in[n]$				&$i\in[n]$												&$i\in[n]$										&$i\in[n]$ \\			
(\ref{property: CC2})				&-			&True$^{\dagger}$						&True$^{\dagger}$														&True											&True\\
(\ref{property: SVM})$_{i}$/(\ref{property: SVG})$_{i}$		&$i\in[n]$	&$i\in[n]$				&-														&-												&-*\\
$\indbit=$				&$\bit$		&$\bit+1$				&$\bit+1$												&$\bit+1$										&$\bit+1$\\\hline										
Special	&Phase 2: 	&\multicolumn{4}{l|}{$\exists i<\relbit{\sigma}\colon$(\ref{property: USV3})$_{i}\wedge\neg$(\ref{property: EV2})$_{i}\wedge\neg$(\ref{property: EV3})$_{i}$}\\
		&Phase 2,3:	&\multicolumn{4}{l|}{$\dagger$A phase-2- resp. phase-3-strategy without \Pref{CC2} is}\\
		&								&\multicolumn{4}{l|}{called \emph{pseudo} phase-2- resp. phase-3-strategy.} \\		
		&Phase 4: 	&\multicolumn{4}{l|}{$\exists i<\nsb(\bit+1)\colon$(\ref{property: USV2})$_{i,1-\indbit_{i+1}}$}\\
		&Phase 5:	&\multicolumn{4}{l|}{*If $\sigma$ has \Pref{ESC1} and there is an index $i$ such that $\sigma$ does}\\
		&						&\multicolumn{4}{l|}{not have \Pref{SVM}$_{i}$ \textbackslash(\ref{property: SVG})$_{i}$, it is defined as a phase-5-strategy} \\\hline
\end{tabular}

\vspace*{0.5em}

\begin{tabular}{|c|ccccc}\hline
		$S_1=$ 	&\multicolumn{5}{l|}{$\{(i,1-\indbit_{i+1})\colon i\in [\nsb-1]\}\cup \{(i,1-\indbit_{i+1})\colon i\in\{\nsb,\dots,m-1\}\wedge\indbit_i=0\} $}\\
					&\multicolumn{5}{l|}{\hspace*{1em}$\cup\begin{cases}\emptyset, &\exists k\in\mathbb{N}\colon \bit+1=2^k\\\{(\nsb,1-\indbit_{\nsb+1})\}, &\nexists k\in\mathbb{N}\colon \bit+1=2^k \end{cases} $}\\\hline
$S_2=$		&\multicolumn{5}{l|}{$\{(i,\indbit^{\sigma}_{i+1})\colon i\in[\nsb-1]\}\cup\{(i,1-\indbit_{i+1})\colon i\in\{\nsb+1,\dots,m\}\wedge\indbit_i=1\} $}\\
					&\multicolumn{5}{l|}{\hspace*{1em}$\cup\{(i,\indbit_{i+1})\colon i\in\{\nsb,\dots,m-1\}\wedge\indbit_i=0\}\cup \{(i,k)\colon i>m, k\in\{0,1\}\}$}\\
					&\multicolumn{5}{l|}{\hspace*{1em}$\cup\begin{cases}\{(\nsb, 1)\}, &\exists k\in\mathbb{N}\colon \bit+1=2^k\\\emptyset, &\nexists k\in\mathbb{N}\colon \bit+1=2^k \end{cases} $}\\\hline
$S_3=$		&\multicolumn{5}{l|}{$\{(i,1-\indbit_{i+1})\colon i\in[u]\}\cup\{(i,1-\indbit_{i+1})\colon i\in\{u+1,\dots,m\}\wedge\indbit_{i}=1\}\cup$}\\
					&\multicolumn{5}{l|}{\hspace*{1em}$\cup\{(i,\indbit_{i+1})\colon i\in\{u+1,\dots,m-1\}\wedge\indbit_i=0\}$}\\
					&\multicolumn{5}{l|}{\hspace*{1em}$\cup\{(i,k)\colon i>m, k\in\{0,1\}\}\cup\{(u,\indbit_{u+1})\}$}\\\hline
$S_4=$		&\multicolumn{5}{l|}{$\{(i,1-\indbit_{i+1})\colon i\in\{u+1,\dots,m-1\}\wedge\indbit_{i}=0\}$}\\\hline
\end{tabular}
\caption[Definition of the phases in the exponential construction.]{Definition of the phases for a strategy $\sigma$ and a number $\bit\in\bitset_n$.
	The entries show for which indices the strategy has the corresponding property resp.~whether the strategy has the property at all.
	Expressions of the type \enquote{$i\in[n]$} or similar are meant as \enquote{$\forall i\in[n]$}
	A '-' signifies that it is not specified whether $\sigma$ has the corresponding property.
	The last row contains further properties used for the definition of the phases.
	The lower table contains all sets used for the definition of the phases and two additional sets $S_3,S_4$ that are necessary for later proofs.
	We use the abbreviations $\nsb\coloneqq\ell(\bit+1), m\coloneqq\max\{i\colon\indbit_i=1\}$ and $u\coloneqq\min\{i\colon\indbit_i=0\}$.} \label{table: Definition of Phases}
\end{table}

\subsection{The Occurrence Records} \label{section: Occurrence records}

We next describe the actual occurrence records that occur when applying the Strategy Iteration resp. Policy Iteration Algorithm.
To do so, we need to introduce notation related to binary counting.

The number of applications of specific edges in level $i$ as improving switches depends on the last time the corresponding cycle centers were closed or how often they were closed.
We thus define $\flips{\bit}{i}$ as the number of numbers smaller than $\bit$ with least significant set bit equal having index $i$.
To quantify how often a specific cycle center was closed, we introduce the \emph{maximal flip number} and the \emph{maximal unflip number}.
Let $\bit\in\bitset_n, i\in\{1,\dots,n\}$ and $j\in\{0,1\}$.
Then, we define the \emph{maximal flip number} $\flipmax{\bit}{i}{\{(i+1,j)\}}$ as the largest $\tilde{b}\leq\bit$ with $\nsb(\tilde{\bit})=i$ and $\tilde{\bit}_{i+1}=j$.
Similarly, we define the \emph{maximal unflip number} $\unflipmax{\bit}{i}{\{(i+1,j)\}}$ as the largest $\tilde{\bit}\leq\bit$ with $\tilde{\bit}_1=\dots=\tilde{\bit}_i=0$ and $\tilde{\bit}_{i+1}=j$.
If there are no such numbers, then $\flipmax{\bit}{i}{\{(i+1,j)\}}\coloneqq0, \unflipmax{\bit}{i}{\{(i+1,j)\}}\coloneqq 0$.
If we do not impose the condition that bit $i+1$ needs to be equal to $j$ then we omit the term in the notation., i.e., $\flipmax{\bit}{i}{}=\max(\{0\}\cup\{\bit'\leq\bit\colon\nsb(\bit')=i\})$ and $\unflipmax{\bit}{i}{}$ is defined analogously.

These notations enable us to properly describe the occurrence records.
We however do not describe the occurrence record for every strategy $\sigma$ produced by the Strategy Improvement resp. Policy Iteration Algorithm.
Instead, we only give a description of the occurrence records for canonical strategies.
When discussing the application of the improving switches, we later prove the following:
Assuming that the occurrence records are described correctly for $\canstrat$, they are also described correctly for $\sigma_{\bit+1}$ when improving switches are applied according to Zadeh's pivot rule and our tie-breaking rule.

\begin{table}[ht]
\small
\centering
\renewcommand{\arraystretch}{1.5}
\begin{tabular}{|C{1.59cm}||c||c|c|}\hline
Edge $e$			&$\occrec^{\canstrat}(e)$					&Edge $e$				&$\occrec^{\canstrat}(e)$\\\hline
$(e_{i,j,k},g_1)$	&$\ceil{\frac{\bit}{2}}$					&$(b_i,g_i)$			&$\flips{\bit}{i}{}$	\\
$(e_{i,j,k},b_2)$	&$\floor{\frac{\bit}{2}}$					&$(b_i,b_{i+1})$		&$\flips{\bit}{i}{}-\bit_i$	\\
$(s_{i,j},h_{i,j})$	&$\flips{\bit}{i+1}{}-(1-j)\cdot\bit_{i+1}$	&$(d_{i,j,k},e_{i,j,k})$&$\leq\begin{cases}\occrec^{\canstrat}(e_{i,j,k},g_1), &\bit_1=0\\\occrec^{\canstrat}(e_{i,j,k},b_2), &\bit_1=1\end{cases}$\\
$(s_{i,j},b_1)$		&$\flips{\bit}{i+1}{}-j\cdot\bit_{i+1}$		&$(g_i,F_{i,j})$			&$\leq\min\limits_{k\in\{0,1\}}\occrec^{\canstrat}(d_{i,j,k},F_{i,j})$\\\hline
\end{tabular}

\vspace*{0.5em}

\begin{tabular}{|c||c|c|}\hline
Condition						&$\occrec^{\canstrat}(d_{i,j,k},F_{i,j})$ 										&Tolerance	\\\hline
$\bit_i=1\wedge\bit_{i+1}=j$	&$\ceil{\frac{\flipmax{\bit}{i}{\{(i+1,j)\}}+1-k}{2}}$							&0			\\
$\bit_i=0\vee\bit_{i+1}\neq j$	&$\min\left (\floor{\frac{\bit+1-k}{2}},\ell^{\bit}(i,j,k)+t_{\bit}\right)$		&$t_{\bit}\in\begin{cases}\{0\},&i=1\vee\bit_i=1\\\{0,1\},&i\neq1\wedge\bit_1=0\\\{-1,0,1\},&i\neq1\wedge\bit_1=1\end{cases}$\\\hline\hline
\multicolumn{3}{|c|}{$\ell^{\bit}(i,j,k)\coloneqq\ceil{\frac{\flipmax{\bit}{i}{\{(i+1,j)\}+1-k}}{2}}+\bit-\id_{j=0}\flipmax{\bit}{i+1}{}-\id_{j=1}\unflipmax{\bit}{i+1}{}$}\\\hline
\end{tabular}
\caption{Occurrence records for the canonical strategy $\canstrat$.
Depending on the edge, we either give the exact occurrence record, give an upper bound, or we give the occurrence record up to a certain tolerance.
A parameter $t_{\bit}$ fulfilling the assumptions for the case $\bit_i=0\vee\bit_{i+1}\neq j$ is called \emph{feasible} for $\bit$.}  \label{table: Occurrence Records}
\end{table}

\begin{restatable}{theorem}{TableDescribesOR} \label{theorem: Table Describes OR}
Let $\canstrat$ be a canonical strategy for $\bit\in\bitset_n$ and assume that improving switches are applied as described in \Cref{section: Phases}.
Then \Cref{table: Occurrence Records} describes the occurrence records of the player 0 edges with respect to $\canstrat$.
\end{restatable}

We now give some intuition for the occurrence records of parts of \Cref{table: Occurrence Records}.
As the occurrence records of most of the edges are much more complicated to explain, we omit an intuitive description of their occurrence records here.
Let $\canstrat\in\reach{\iota}$ be a canonical strategy for $\bit\in\bitset_n$.

Consider some edge $(b_i,g_i)$.
This edge is applied as an improving switch whenever bit $i$ switches from 0 to 1.
That is, it is applied if and only if we transition towards some $\bit'\in\bitset_n$ with $\nsb(\bit')=i$ and $\bit'\leq\bit$.
Therefore, $\occrec^{\canstrat}(b_i,g_i)=\flips{\bit}{i}{}$.
Now consider $(b_i,b_{i+1})$.
This edge is only applied as an improving switch when bit $i$ switches from 1 to 0.
This can however only happen if bit $i$ switched from 0 to 1 earlier.
That is, applying $(b_i,b_{i+1})$ can only happen when $(b_i,g_i)$ was applied before.
Also, we can only apply the switch $(b_i,g_i)$ again after bit $i$ has been switched back to 0 again, i.e., after $(b_i,b_{i+1})$ was applied.
Consequently, $\occrec^{\canstrat}(b_i,b_{i+1})=\occrec^{\canstrat}(b_i,g_i)-\bit_i=\flips{\bit}{i}{}-\bit_i$.

Next, consider some edge $(s_{i,j},h_{i,j})$ and fix $j=1$ for now.
This edge is applied as an improving switch if and only if bit $i+1$ switches from 0 to~1.
Hence, as discussed before, $\occrec^{\canstrat}(s_{i,j},h_{i,j})=\flips{\bit}{i+1}{}$.
Now let $j=0$.
The switch $(s_{i,0},h_{i,0})$ is applied whenever bit $i+1$ switches from 1 to 0.
This requires the bit to have switched from 0 to 1 before.
Therefore,  $\occrec^{\canstrat}(s_{i,0},h_{i,0})=\occrec^{\canstrat}(s_{i,1},h_{i,1})-\bit_{i+1}=\flips{\bit}{i+1}{}-\bit_{i+1}$.
Further note that the switch $(s_{i,j},b_1)$ is applied in the same transitions in which the switch $(s_{i,1-j},h_{i,1-j})$ is applied.
Hence, $\occrec^{\canstrat}(s_{i,j},h_{i,j})=\flips{\bit}{i+1}{}-(1-j)\cdot\bit_{i+1}$ and $\occrec^{\canstrat}(s_{i,j},b_1)=\flips{\bit}{i+1}{}-j\cdot\bit_{i+1}$.

Finally consider some edge $(e_{i,j,k},g_1)$.
This edge is applied as improving switch whenever the first bit switches from 0 to 1.
Since $0$ is even, this happens once for every odd numbers smaller than or equal to $\bit$, i.e., $\ceil{\frac{\bit}{2}}$ times.
Since the switch $(e_{i,j,k},b_2)$ is applied during each transition in which the switch $(e_{i,j,k},g_1)$ is not applied, we have $\occrec(e_{i,j,k},g_1)=\bit-\ceil{\frac{\bit}{2}}= \floor{\frac{\bit}{2}}$ as $\bit\in\mathbb{N}$.

\subsection{Proving the lower bound} \label{section: Proving Lower Bound}

We now describe how the exponential lower bound for Zadeh's pivot rule for the strategy improvement algorithm is proven, implying the same bound for the simplex algorithm applied to linear programs and similar algorithms.
This section does however not contain all of the formal details.
Instead, it presents the core concepts and aspects of the proofs as well as the most important ideas and arguments that are used to derive the lower bound.
Of course, this approach does not suffice to give a formal proof.
Incorporating both the intuitive ideas and the \enquote{core} of our proof as well as the technical details in the same chapter would make it extremely hard to understand the approach and statements.
All of the necessary formalism is thus introduced, proven and discussed in detail in the appendices.

We prove that applying improving switches to $G_n$ using Zadeh's pivot rule and our tie-breaking rule requires an exponential number of iterations using an inductive argument.
Assume we are given a canonical strategy $\canstrat$ for $\bit\in\bitset_n$ that has some helpful additional properties as well as an occurrence record as described by \Cref{table: Occurrence Records}.
We prove that the application of improving switches eventually yields a canonical strategy $\sigma_{\bit+1}$ for $\bit+1$ that has the same additional properties and whose occurrence record is also described by \Cref{table: Occurrence Records} when interpreted for $\bit+1$.
It is then sufficient to prove that the initial strategy~$\sigma_0$ has these properties already and that  $\occrec^{\sigma_{0}}$ is described by \Cref{table: Occurrence Records}, then the exponential lower bound on the number of iterations follows immediately.
Since $G_n$ has a linear number of vertices and edges and the priorities, rewards and probabilities can be encoded using a polynomial number of bits, this implies an exponential lower bound for the respective algorithms.
We begin explaining the proof by introducing the mentioned set of properties.

\subsubsection{The canonical properties and basic statements}
The additional properties of the canonical strategies are called \emph{canonical properties}.
Two of these are straight-forward.
They state that \Cref{table: Occurrence Records} correctly describes the occurrence records for $\canstrat$ and that each improving switch was applied at most once per previous transition $\sigma_{\bit'}\to\sigma_{\bit'+1}$ with $\bit'<\bit$.
The remaining properties are more involved and introduced in more detail.
An overview over these properties can be found in \Cref{table: Occurrence Records Properties}.

Let $\bit\in\bitset_n$.
Consider some cycle center $F_{i,j}$ with $\bit_i=0\vee\bit_{i+1}\neq j$.
Then, $F_{i,j}$ should not be closed for $\canstrat$.
It can however still happen that $\canstrat(d_{i,j,k})=F_{i,j}$ for some $k\in\{0,1\}$.
\Pref{OR1} states that this can only happen if the occurrence record of $(d_{i,j,k},F_{i,j})$ is sufficiently low.
Formally, for a general strategy $\sigma$, the property is defined as 
\begin{property}{OR1}\sigma(d_{i,j,k})=F_{i,j}\wedge(\bit_i=0\vee\bit_{i+1}\neq j)\implies \occrec^{\sigma}(d_{i,j,k},F_{i,j})<\floor{\frac{\bit+1}{2}}.\end{property}

\Pref{OR2} characterizes under which circumstances the parameter $t_{\bit}$ used for describing the occurrence records of cycle edges (see \Cref{table: Occurrence Records}) is equal to $1$.
More precisely, it states that the parameter is equal to $1$ if and only if the cycle vertex points towards the cycle center.
Of course, this statement is only valid if either $\indbit^{\sigma}_i=0$ or $\indbit^{\sigma}_{i+1}\neq j$ since the parameter is only then relevant for describing the occurrence record.
Formally, 
\begin{property}{OR2}\indbit^{\sigma}_i=0\vee\indbit^{\sigma}_{i+1}\neq j\implies(\occrec^{\sigma}(d_{i,j,k},F_{i,j})=\ell^{\bit}(i,j,k)+1 \Longleftrightarrow \sigma(d_{i,j,k})=F_{i,j}). \end{property}

Analogously, \Pref{OR3} gives a characterization regarding the cases in which the parameter is equal to $-1$.
This characterization is more involved as it also depends on the exact value of the occurrence record of a cycle edge.
It  states that the parameter can only be $-1$ without being equal to $\floor{(\bit+1-k)/2}$ if and only if (i) $\bit$ is odd, (ii) $\bit+1$ is not a power of $2$, (iii) $i=\ell(\bit+1)$, (iv) $j\neq\bit_{i+1}$ and (v) $k=0$.
Formally, 
\begin{property}{OR3}\begin{split}
&\occrec^{\sigma}(d_{i,j,k},F_{i,j})=\ell^{\bit}(i,j,k)-1\,\wedge\,\occrec^{\sigma}(d_{i,j,k},F_{i,j})\neq\floor{\frac{\bit+1-k}{2}}\\
\Longleftrightarrow&\bit\bmod2=1\wedge\nexists l\in\mathbb{N}:\bit+1=2^l\wedge i=\ell(\bit+1)\wedge j\neq\bit_{i+1}\wedge k=0. 
\end{split}\end{property}

The final property states that the occurrence record of any cycle edge with $\sigma(d_{i,j,k})\neq F_{i,j}$ is relatively high.
Formally,
\begin{property}{OR4} \sigma(d_{i,j,k})\neq F_{i,j}\implies\occrec^{\sigma}(d_{i,j,k},F_{i,j})\in \left \{ \floor{\frac{\bit+1}{2}}-1, \floor{\frac{\bit+1}{2}}\right\}.\end{property}

\begin{table}[ht]
\small
\centering
\renewcommand{\arraystretch}{1.5}
\begin{tabular}{|c|l|}\hline
(\ref{property: OR1})$_{i,j,k}$											& $\sigma(d_{i,j,k})=F_{i,j}\wedge(\bit_i=0\vee\bit_{i+1}\neq j)\implies \occrec^{\sigma}(d_{i,j,k},F_{i,j})<\floor{\frac{\bit+1}{2}}$\\\hline
(\ref{property: OR2})$_{i,j,k}$											& $\indbit^{\sigma}_i=0\vee\indbit^{\sigma}_{i+1}\neq j\implies(\occrec^{\sigma}(d_{i,j,k},F_{i,j})=\ell^{\bit}(i,j,k)+1 \Longleftrightarrow \sigma(d_{i,j,k})=F_{i,j})$\\\hline
\multirow{2}{*}{(\ref{property: OR3})$_{i,j,k}$} 	& $\occrec^{\sigma}(d_{i,j,k},F_{i,j})=\ell^{\bit}(i,j,k)-1\,\wedge\,\occrec^{\sigma}(d_{i,j,k},F_{i,j})\neq\floor{\frac{\bit+1-k}{2}}$\\
																											& \hspace*{2em}$\Longleftrightarrow\bit\bmod2=1\wedge \nexists l\in\mathbb{N}:\bit+1=2^l\wedge i=\ell(\bit+1)\wedge j\neq\bit_{i+1}\wedge k=0$ \\\hline
(\ref{property: OR4})$_{i,j,k}$											&$\sigma(d_{i,j,k})\neq F_{i,j}\implies\occrec^{\sigma}(d_{i,j,k},F_{i,j})\in \left \{ \floor{\frac{\bit+1}{2}}-1, \floor{\frac{\bit+1}{2}}\right\}$\\\hline
\end{tabular}
\caption{The additional properties of canonical strategies.} \label{table: Occurrence Records Properties}
\end{table}

This now allows to define the canonical properties formally.

\begin{definition}[Canonical properties] \label{definition: Canonical properties}
Let $\sigma\in\reach{\sigma}$ be a strategy for $G_n$.
Then, $\sigma$ has the \emph{canonical properties}\index{canonical properties} if 
\begin{enumerate}
	\item the occurrence records $\occrec^{\sigma}$ are described correctly by \Cref{table: Occurrence Records},
	\item $\sigma$ has Properties (\ref{property: OR1})$_{*,*,*}$ to (\ref{property: OR4})$_{*,*,*}$ and
	\item any improving switch was applied at most once per previous transition between canonical strategies.
\end{enumerate}
\end{definition}

We begin our arguments by explicitly determining the set of improving switches for canonical strategies. 
We also extend our notation of transitions.
Let $\sigma,\sigma'\in\reach{\sigma_0}$ denote two strategies and let $\sigma'$ be reached after $\sigma$.
We denote the sequence of strategies calculated by the algorithm while transitioning from $\sigma$ to $\sigma'$ by $\sigma\to\sigma'$.
In addition, the sequence of actually applied improving switches is denoted by $\applied{\sigma}{\sigma'}$.
Throughout this section, let $\bit\in\bitset_n$ be fixed and $\nsb\coloneqq\ell(\bit+1)$.
As most of the statements presented here, its proof can be found in the appendices.

\begin{restatable}{lemma}{ImprovingSetsOfCanonicalStrategies} \label{lemma: Improving sets of canonical strategies}
Let $\canstrat\in\reach{\sigma_0}$ be a canonical strategy for $\bit\in\bitset_n$.
Then, $\canstrat$ is a phase-$1$-strategy for $\bit$ and $I_{\canstrat}=\{(d_{i,j,k},F_{i,j}):\canstrat(d_{i,j,k})\neq F_{i,j}\}$.
\end{restatable}

As our proof is inductive, we need to give a basis for the induction.

\begin{lemma} \label{lemma: Properties of initial strategy}
The initial strategy $\sigma_0$ is a canonical strategy for $\bit=0$ and has all canonical properties. 
\end{lemma}

\begin{proof}
As no improving switch was applied yet and it is easy to verify that $\sigma_0$ is a canonical strategy for $0$, it suffices to prove that $\sigma_0$ has Properties (\ref{property: OR1})$_{*,*,*}$ to (\ref{property: OR4})$_{*,*,*}$.

Let $i\in[n]$ and $j,k\in\{0,1\}$.
First,~$\sigma_0$ has \Pref{OR1}$_{i,j,k}$ as $\sigma_0(d_{i,j,k})=e_{i,j,k}$.
In addition, $0=\occrec^{\sigma_0}(d_{i,j,k},F_{i,j})<1\leq\ell^{\bit}(i,j,k)+1$, so~$\sigma_0$ has \Pref{OR2}$_{i,j,k}$.
Moreover, $\occrec^{\sigma_0}(d_{i,j,k},F_{i,j})=0=\floor{(1-k)/2}=\floor{(\bit+1-k)/2},$ hence the premise of \Pref{OR3}$_{i,j,k}$ is incorrect.
Thus, $\sigma_0$ has \Pref{OR3}$_{i,j,k}$.
Since it is immediate that $\sigma_0$ has \Pref{OR4}$_{i,j,k}$, the statement follows.
\end{proof}

We now discuss how the main statement is proven in more detail.
Consider a canonical strategy $\canstrat$ for $\bit$. 
We prove that applying improving switches according to Zadeh's pivot rule and the tie-breaking rule given in \Cref{definition: Tie-Breaking exponential} produces a specific phase-$k$-strategy for $\bit$ for every $k\in[5]$.
These strategies are the first phase-$k$-strategies for~$\bit$ that the algorithm reaches and have several properties that allow  us to simplify the proofs.
These properties are summarized in \Cref{table: Properties at start of phase}.
We furthermore explicitly state the improving switches with respect to these \enquote{initial phase-$k$-strategies} in \Cref{table: Switches at start of phase}.
Both tables distinguish whether the number $\bit+1$ is even or odd and specify whether certain entries are only valid for $S_n$ resp. $M_n$.

\begin{table}
\centering
\footnotesize
\renewcommand{\arraystretch}{1.33}
\begin{tabular}{|c||c|c|c|}\hline
\multirow{2}{*}{Phase}		&\multirow{2}{*}{$\nsb=1$}	&\multicolumn{2}{c|}{$\nsb>1$}\\\cdashline{3-4}
							&							&$S_n$			&$M_n$			\\\hline\hline
1							&\multicolumn{3}{c|}{Canonical strategy for $\bit$ having the canonical properties}\\\hline
\multirow{5}{*}{2}			&\multirow{5}{*}{-}			&\multicolumn{2}{c|}{$\sigma(d_{i,j,k})\neq F_{i,j}\Rightarrow\occrec^{\sigma}(d_{i,j,k},F_{i,j})=\floor{\frac{\bit+1}{2}}$}\\
							&							&\multicolumn{2}{c|}{$(g_i,F_{i,j})\in\applied{\canstrat}{\sigma}\hspace*{-1pt}\Rightarrow\hspace*{-1pt}[\bit_i=0\wedge\bit_{i+1}\neq j]\vee i=\nsb$ and $F_{i,j}$ is closed}\\
							&							&\multicolumn{2}{c|}{$\applied{\canstrat}{\sigma}\subseteq\mathbb{D}^1\cup\mathbb{G}$}\\
							&							&\multicolumn{2}{c|}{$\sigma(g_{\nsb})=F_{\nsb,\indbit_{\nsb+1}^{\sigma}}$ and $\sigma(g_i)=F_{i,1-\indbit^{\sigma}_{i+1}}$ for all $i<\nsb$}\\
							&							&\multicolumn{2}{c|}{$i<\nsb\Rightarrow \sigmabar(d_{i})$ and \Pref{USV3}$_i$}\\\hline
\multirow{3}{*}{3}			&\multicolumn{3}{c|}{$\sigma(d_{i,j,k})\neq F_{i,j}\Rightarrow\occrec^{\sigma}(d_{i,j,k},F_{i,j})=\occrec^{\canstrat}(d_{i,j,k},F_{i,j})=\floor{\frac{\bit+1}{2}}$}\\
							&\multicolumn{3}{c|}{$\sigma(s_{i,*})=h_{i,*}, \sigma(g_i)=F_{i,1-\indbit^{\sigma}_{i+1}}$ and $\sigmabar(d_i)$ for all $i<\nsb$ as well as $\sigma(g_{\nsb})=F_{\nsb,\bit_{\nsb+1}}$}\\
							&\multicolumn{3}{c|}{$\applied{\canstrat}{\sigma^3}\subseteq\mathbb{D}^1\cup\mathbb{G}\cup\mathbb{S}\cup\mathbb{B}$ and $(g_i,F_{i,j})\in\applied{\canstrat}{\sigma}\Rightarrow[\bit_i=0\wedge\bit_{i+1}\neq j]\vee i=\nsb$ and $F_{i,j}$ is closed}\\\hline
\multirow{8}{*}{4}			&\multirow{9}{*}{-}	&$\relbit{\sigma}=\min\{i\colon\indbit^{\sigma}_i=0\}$																	&\multirow{8}{*}{-}\\
							&					&$\sigma(d_{i,j,*})=F_{i,j}\Leftrightarrow\indbit^{\sigma}_{i}=1\wedge(\bit+1)_{i+1}=j$													&\\
							&					&$(i,j)\in S_1\Rightarrow\sigmabar(eb_{i,j})\wedge\nsigmabar(eg_{i,j})$																			&\\
							&					&$(i,j)\in S_2\Rightarrow\sigmabar(eb_{i,j})\wedge\sigmabar(eg_{i,j})$																			&\\
							&					&$(g_i,F_{i,j})\in\applied{\canstrat}{\sigma}\Rightarrow[\bit_i=0\wedge\bit_{i+1}\neq j]\vee i=\nsb$									&\\
							&					&$(g_i,F_{i,j})\in\applied{\canstrat}{\sigma}\Rightarrow$ $F_{i,j}$ is closed															&\\	
							&					&$\sigma(e_{i,j,k})=b_2\Rightarrow\occrec^{\sigma}(d_{i,j,k},F_{i,j})=\occrec^{\canstrat}(d_{i,j,k},F_{i,j})$	&\\
							&					&$\sigma(e_{i,j,k})=b_2\Rightarrow\occrec^{\sigma}(d_{i,j,k},F_{i,j})=\floor{\frac{\bit+1}{2}}$	&\\\hline
\multirow{7}{*}{5}			&\multicolumn{3}{c|}{$\relbit{\sigma}=\min\{i\colon\indbit^{\sigma}_{i+1}=0\}$}		\\
													&\multicolumn{3}{c|}{$\sigma(d_{i,j,k})=F_{i,j}\Leftrightarrow\indbit^{\sigma}_{i}=1\wedge\indbit^{\sigma}_{i+1}=j$}		\\
													&\multicolumn{3}{c|}{$(g_i,F_{i,j})\in\applied{\canstrat}{\sigma}\Rightarrow[\bit_i=0\wedge\bit_{i+1}\neq j]\vee i=\nsb$}		\\
													&\multicolumn{3}{c|}{$\sigma(e_{i,j,k})=t^{\rightarrow}\Rightarrow\occrec^{\sigma}(d_{i,j,k},F_{i,j})=\occrec^{\canstrat}(d_{i,j,k},F_{i,j})=\floor{\frac{\bit+1}{2}}$}\\\cdashline{2-4}
													&\multirow{2}{*}{$(i,j)\in S_3\Rightarrow\sigmabar(eg_{i,j})\wedge\sigmabar(eb_{i,j})$}	&\multicolumn{2}{c|}{$(i,j)\in S_1\Rightarrow\sigmabar(eb_{i,j})\wedge\nsigmabar(eg_{i,j})$}\\
													&\multirow{2}{*}{$(i,j)\in S_4\Rightarrow\sigmabar(eg_{i,j})\wedge\nsigmabar(eb_{i,j})$}	&\multicolumn{2}{c|}{$(i,j)\in S_2\Rightarrow\sigmabar(eb_{i,j})\wedge\sigmabar(eg_{i,j})$}\\
													&		&\multicolumn{2}{c|}{$i<\nsb\Rightarrow\sigmabar(g_i)=1-\indbit^{\sigma}_{i+1}$}\\\hline
1							&\multicolumn{3}{c|}{Canonical strategy for $\bit+1$ having the canonical properties}\\\hline
\end{tabular}
\caption[Properties of specific phase-$k$-strategies.]{Properties of specific phase-$k$-strategies.
To simplify notation, let $t^{\rightarrow}\coloneqq g_1$ if $\nsb=1$ and $t^{\rightarrow}\coloneqq b_2$ if $\nsb>1$.
A '-' signifies that the corresponding combination does not occur during the execution of the algorithm.} \label{table: Properties at start of phase}
\end{table}

\begin{table}[ht]
\centering
\footnotesize
\renewcommand{\arraystretch}{1.5}
\begin{tabular}{|c|c|c|c|}
\hline%
\multirow{2}{*}{Phase}		&\multirow{2}{*}{$\nsb=1$}	&\multicolumn{2}{c|}{$\nsb>1$}	\\\cdashline{3-4}
							&							&$S_n$			&$M_n$				\\\hline	
1							&\multicolumn{3}{c|}{$\mathfrak{D}^{\sigma}\coloneqq\{(d_{i,j,k},F_{i,j})\colon\sigma(d_{i,j,k})\neq F_{i,j}\}$} \\\hline
2							&-								&\multicolumn{2}{c|}{$\mathfrak{D}^{\sigma}\cup\{(b_{\nsb},g_{\nsb}), (s_{\nsb-1,1},h_{\nsb-1,1})\}$}\\\hline
3							&$\mathfrak{D}^{\sigma}\cup\{(b_1,g_1)\}\cup\{(e_{*,*,*},g_1)\}$	&\multicolumn{2}{c|}{$\mathfrak{D}^{\sigma}\cup\{(b_1,b_2)\}\cup\{(e_{*,*,*},b_2)\}$}\\\hline
4							&-							& $\mathfrak{E}^{\sigma}\cup X_0\cup X_1\cup\{(s_{\nsb-1,0},b_1)\}\cup\{(s_{i,1},b_1)\colon i\leq\nsb-2\}$					&-\\\hline
5							&$\displaystyle\mathfrak{E}^{\sigma}\cup\hspace*{-1.5pt}\bigcup^{m-1}_{\substack{i=\relbit{\sigma}+1\\\indbit^{\sigma}_i=0}}\hspace*{-1.5pt}\{(d_{i,1-\indbit^{\sigma}_{i+1},*},F_{i,1-\indbit^{\sigma}_{i+1}})\}$	&\multicolumn{2}{C{8cm}|}{$\displaystyle\mathfrak{E}^{\sigma} \cup \bigcup^{\nsb-1}_{i=1}\{d_{i,1-\indbit^{\sigma}_{i+1},*}, F_{i,-1-\indbit^{\sigma}_{i+1}})\}\cup X_0\cup X_1$}\\\hline
1							&\multicolumn{3}{c|}{$\mathfrak{D}^{\sigma}$}\\\hline\hline
\multicolumn{4}{|c|}{\multirow{4}{*}{$		\displaystyle X_k\coloneqq\begin{cases}
			\emptyset, &\bit+1\text{ is a power of two},\\
			\{(d_{\nsb,1-\indbit^{\sigma}_{\nsb},k},F_{\nsb,1-\indbit^{\sigma}_{\nsb}})\}\cup\bigcup\limits_{\substack{i=\nsb+1\\\indbit_{i}=0}}^{m-1}\{(d_{i,1-\indbit^{\sigma}_{i+1},k},F_{i,1-\indbit^{\sigma}_{i+1}})\}, &\text{otherwise.}\end{cases}
$}}\\\multicolumn{4}{|c|}{}\\\multicolumn{4}{|c|}{}\\\multicolumn{4}{|c|}{}\\\hline
\end{tabular}
\caption[The improving switches at the beginning of different phases.]{The improving switches at the beginning of the different phases.
We define $m\coloneqq\max\{i\colon\sigma(b_i)=g_i\}$ and $\mathfrak{E}^{\sigma}\coloneqq\{(d_{i,j,k},F_{i,j}),(e_{i,j,k},b_2)\colon\sigma(e_{i,j,k})=g_1\}$ if $\nsb>1$.
Analogously, we let $\mathfrak{E}^{\sigma}\coloneqq\{(d_{i,j,k},F_{i,j}),(e_{i,j,k},g_1)\colon\sigma(e_{i,j,k})=b_2\}$ if $\nsb=1$.
We do not interpret $1$ as a power of two.} \label{table: Switches at start of phase}
\end{table}

\subsubsection{Detailed application of the improving switches}

We now discuss the application of the improving switches during the individual phases more formally and link this description to the previously given tables.
We illustrate this procedure by providing sketches of the first phase-$k$-strategies that are calculated by the algorithm in the sink game $S_3$ when transitioning from $\sigma_3$ to $\sigma_4$.
As $4$ is an even number, we have $\nsb>1$ in this example.
The canonical strategy $\sigma_3$  in the sink game $S_3$ shown given in \Cref{figure: Start Phase 1}.

\begin{figure}[ht]
\centering
\begin{tikzpicture}[scale=0.33]

	\def\x{2.25}
	\def\y{1.55}
	\tikzset{
		player/.style={circle, draw, inner sep=1pt, scale=0.4, fill=playblue, text=playblue, minimum size=2.5em},
		random/.style={draw, inner sep=1pt,  scale=0.4, fill=ranred, text=ranred, minimum size=2.5em},
	} 
	
	\node[player, dashed,align=center] (PairSelector2) at  (0*\x,\y) {$b_2$};
	\node[player, dashed,align=center] (Selector1) at  (0*\x,3*\y) {$g_1$};
	\node[player, dashed,align=center] (PairSelector1) at  (0,4*\y) {$b_1$};
	\node[player,align=center] (t) at (0*\x,15*\y) {$t$};
	\draw[->,playblue, very thick ] (t) to[out=60, in=120, looseness=5] (t);
	
	\node[player,align=center] (g1) at (0,0) {$g_1$};
	\node[random,align=center] (F10) at (-3*\x,2*\y) {$F_{1,0}$};
	\node[random,align=center] (F12) at (3*\x,2*\y) {$F_{1,1}$};
	\node[player,align=center] (d101) at (-2*\x,3*\y) {$d_{1,0,1}$};
	\node[player,align=center] (e101) at (-\x,3*\y) {$e_{1,0,1}$};
	\node[player,align=center] (d111) at (2*\x,3*\y) {$d_{1,1,1}$};
	\node[player,align=center] (e111) at (1*\x,3*\y) {$e_{1,1,1}$};
	\node[player,align=center] (d100) at (-2*\x,1*\y) {$d_{1,0,0}$};
	\node[player,align=center] (e100) at (-\x,1*\y) {$e_{1,0,0}$};
	\node[player,align=center] (d110) at (2*\x,1*\y) {$d_{1,1,0}$};
	\node[player,align=center] (e110) at (\x,1*\y) {$e_{1,1,0}$};
	\node[player,align=center] (h10) at (-3*\x,4*\y) {$h_{1,0}$};
	\node[player,align=center] (b1) at (-4*\x,0) {$b_1$};
	\node[player,align=center] (s10) at (-2*\x,4*\y) {$s_{1,0}$};
	\node[player,align=center] (s11) at (2*\x,4*\y) {$s_{1,1}$};
	\node[player,align=center] (h11) at (3*\x,4*\y) {$h_{1,1}$};	
	
	\draw[dotted,  ] (-4.5*\x,4.5*\y)--(3.5*\x,4.5*\y);
	
	\draw[->] (e101)-- (PairSelector2);
	\draw[->, playblue, very thick] (e101)--(Selector1);
	
	\draw[->] (e111)-- (PairSelector2);
	\draw[->, playblue, very thick] (e111)--  (Selector1);
	
	\draw[->] (e100)-- (PairSelector2);
	\draw[->, playblue, very thick] (e100)-- (Selector1);

	\draw[->] (e110)--(PairSelector2);
	\draw[->,playblue, very thick] (e110)--(Selector1);
	
	\draw[->] (s10)--(h10);
	\draw[->, playblue, very thick ] (s10)--(PairSelector1);
	
	\draw[->, playblue, very thick] (s11)--(h11);
	\draw[->] (s11)--(PairSelector1);
	
\begin{scope}[yshift=5*\y cm]
	\node[player, dashed,align=center] (PairSelector1) at  (0,4*\y) {$b_1$};
	\node[player, dashed,align=center] (PairSelector2) at  (0*\x,\y) {$b_2$};
	\node[player, dashed,align=center] (Selector1) at  (0*\x,3*\y) {$g_1$};

	\node[player,align=center] (g2) at (0,0) {$g_2$\\$13$};
	\node[random,align=center] (F20) at (-3*\x,2*\y) {$F_{2,0}$};
	\node[random,align=center] (F22) at (3*\x,2*\y) {$F_{2,1}$};
	\node[player,align=center] (d201) at (-2*\x,3*\y) {$d_{2,0,1}$};
	\node[player,align=center] (e201) at (-\x,3*\y) {$e_{2,0,1}$};
	\node[player,align=center] (d211) at (2*\x,3*\y) {$d_{2,1,1}$};
	\node[player,align=center] (e211) at (1*\x,3*\y) {$e_{2,1,1}$};
	\node[player,align=center] (d200) at (-2*\x,1*\y) {$d_{2,0,0}$};
	\node[player,align=center] (e200) at (-\x,1*\y) {$e_{2,0,0}$};
	\node[player,align=center] (d210) at (2*\x,1*\y) {$d_{2,1,0}$};
	\node[player,align=center] (e210) at (\x,1*\y) {$e_{2,1,0}$};
	\node[player,align=center] (h20) at (-3*\x,4*\y) {$h_{2,0}$};
	\node[player,align=center] (b2) at (-4*\x,0) {$b_2$};
	\node[player,align=center] (s20) at (-2*\x,4*\y) {$s_{2,0}$};
	\node[player,align=center] (s21) at (2*\x,4*\y) {$s_{2,1}$};
	\node[player,align=center] (h21) at (3*\x,4*\y) {$h_{2,1}$};	
	
	\draw[dotted,  ] (-4.5*\x,4.5*\y)--(3.5*\x,4.5*\y);
	
	\draw[->, playblue, very thick ] (d201) to [out=270-30, in=0+30] (F20);
	\draw[->] (d201)-- (e201);
	
	\draw[->] (e201)-- (PairSelector2);
	\draw[->, playblue, very thick] (e201)--(Selector1);
	
	\draw[->, green!70!black, very thick] (d211) to [out=270+30, in=180-30] (F22);
	\draw[->, playblue, very thick  ] (d211)-- (e211);
	
	\draw[->] (e211)-- (PairSelector2);
	\draw[->,playblue, very thick] (e211)--  (Selector1);
	
	\draw[->, playblue, very thick] (d200) to[out=90+30, in=0-30]  (F20);
	\draw[->] (d200)--(e200);
	
	\draw[->] (e200)-- (PairSelector2);
	\draw[->,playblue, very thick] (e200)-- (Selector1);
	
	\draw[->, green!70!black, very thick] (d210) to[out=90-30, in=180+30]  (F22);
	\draw[->, playblue, very thick ] (d210)-- (e210);
	
	\draw[->] (e210)--(PairSelector2);
	\draw[->,playblue, very thick] (e210)--(Selector1);
	
	\draw[->, playblue, very thick] (g2) to[out=180-15, in=270] (F20);
	\draw[->] (g2)to[out=0+15, in=270]  (F22);
	
	\draw[->,playblue , very thick  ] (b2)--(g2);
	
	\draw[->] (F20) to[out=90-30, in=180+30] (d201);
	\draw[->] (F20) to[out=270+30, in=180-30] (d200);
	\draw[->, ranred, very thick] (F20) to[out=90, in=180+30] (s20);
	
	\draw[->] (F22) to[out=90+30, in=0-30] (d211);
	\draw[->, ranred, very thick ] (F22)to[out=270-30, in=0+30]  (d210);
	\draw[->] (F22) to[out=90, in=0-30](s21);
	
	\draw[->,playblue, very thick] (s20)--(h20);
	\draw[->] (s20)--(PairSelector1);
	
	\draw[->] (s21)--(h21);
	\draw[->, playblue, very thick] (s21)--(PairSelector1);
	\draw[->, playblue, very thick] (h20) --++(-\x,6*\y)--(t);
\end{scope}

\begin{scope}[yshift=10*\y cm]

	\node[player, dashed,align=center] (PairSelector1) at  (0,4*\y) {$b_1$};
	\node[player, dashed,align=center] (PairSelector2) at  (0*\x,\y) {$b_2$};
	\node[player, dashed,align=center] (Selector1) at  (0*\x,3*\y) {$g_1$};

	\node[player,align=center] (g3) at (0,0) {$g_3$};
	\node[random,align=center] (F30) at (-3*\x,2*\y) {$F_{3,0}$};
	\node[player,align=center] (d301) at (-2*\x,3*\y) {$d_{3,0,1}$};
	\node[player,align=center] (e301) at (-\x,3*\y) {$e_{3,0,1}$};
	\node[player,align=center] (d300) at (-2*\x,1*\y) {$d_{3,0,0}$};
	\node[player,align=center] (e300) at (-\x,1*\y) {$e_{3,0,0}$};
	\node[player,align=center] (h30) at (-3*\x,4*\y) {$h_{3,0}$};
	\node[player,align=center] (b3) at (-4*\x,0) {$b_3$};
	\node[player,align=center] (s30) at (-2*\x,4*\y) {$s_{3,0}$};
	
	\draw[dotted,  ] (-4.5*\x,4.5*\y)--(3.5*\x,4.5*\y);
	
	\draw[->, green!70!black, very thick] (d301) to [out=270-30, in=0+30] (F30);
	\draw[->, playblue, very thick] (d301)-- (e301);
	
	\draw[->] (e301)-- (PairSelector2);
	\draw[->,  playblue, very thick] (e301)--(Selector1);
	
	\draw[->,  green!70!black, very thick] (d300) to[out=90+30, in=0-30]  (F30);
	\draw[->, playblue, very thick] (d300)--(e300);
	
	\draw[->] (e300)-- (PairSelector2);
	\draw[->,  playblue, very thick] (e300)-- (Selector1);
	
	\draw[->,  playblue, very thick ] (g3) to[out=180-15, in=270] (F30);
	
	\draw[->] (b3)--(g3);
	
	\draw[->] (F30) to[out=90-30, in=180+30] (d301);
	\draw[->, ,  ] (F30) to[out=270+30, in=180-30] (d300);
	\draw[->, ranred, very thick] (F30) to[out=90, in=180+30] (s30);
	
	\draw[->,  playblue, very thick] (s30)--(h30);
	\draw[->] (s30)--(PairSelector1);
	\draw[->,  playblue, very thick] (h30) to[out=20,in=180] (t);
	
	\draw[->,very thick, playblue ] (b3) --++(0,5*\y)--(t);
	
\end{scope}

	\draw[->, green!70!black, very thick] (d101) to [out=270-30, in=0+30] (F10);
	\draw[->,  playblue, very thick] (d101)-- (e101);
	
	\draw[->, playblue, very thick ] (d111) to [out=270+30, in=180-30] (F12);
	\draw[->] (d111)-- (e111);
	
	\draw[->,  green!70!black, very thick] (d100) to[out=90+30, in=0-30]  (F10);
	\draw[->, playblue, very thick  ] (d100)--(e100);
	
	\draw[->,  playblue, very thick] (d110) to[out=90-30, in=180+30]  (F12);
	\draw[->] (d110)-- (e110);
	
	\draw[->] (g1) to[out=180-15, in=270] (F10);
	\draw[->,  playblue, very thick ] (g1)to[out=0+15, in=270]  (F12);
	
	\draw[->,  playblue, very thick] (b1)--(g1);
	\draw[->] (b1)--(b2);
	\draw[->] (b2)--(b3);
	
	\draw[->] (F10) to[out=90-30, in=180+30] (d101);
	\draw[->, ranred, very thick ] (F10) to[out=270+30, in=180-30] (d100);
	\draw[->] (F10) to[out=90, in=180+30] (s10);
	
	\draw[->] (F12) to[out=90+30, in=0-30] (d111);
	\draw[->] (F12)to[out=270-30, in=0+30]  (d110);
	\draw[->, ranred, very thick] (F12) to[out=90, in=0-30](s11);
	
	\draw[->, playblue, very thick] (h11) to[out=150,in=0] (g2);
	\draw[->, playblue, very thick ] (h21) to[out=150,in=0] (g3);
	\draw[->,playblue ,very thick  ] (h10)--(b3);
\end{tikzpicture}
\caption[The canonical strategy $\sigma_3$ in $S_3$]{The canonical strategy $\sigma_3$ in the sink game $S_3$.
Blue edges represent choices of player $0$, red choices represent choices of player $1$, and green edges represent improving switches. 
For simplification, we omit the labels here.} \label{figure: Start Phase 1}
\end{figure}

By \Cref{table: Switches at start of phase}, the set of improving switches is given by all edges $(d_{i,j,k},F_{i,j})$ with $\canstrat(d_{i,j,k})\neq F_{i,j}$ at the beginning of phase $1$.
During phase $1$, all improving switches $e=(d_{i,j,k},F_{i,j})\in I_{\canstrat}$ with $\occrec^{\canstrat}(e)=\floor{(\bit+1)/2}-1$ are applied first as these minimize the occurrence records, cf. \Pref{OR4}.
This might close an inactive cycle center $F_{i,1-\bit_{i+1}}$, making the edge $(g_i,F_{i,1-\bit_{i+1}})$ improving if $\bit_i=0$ and $\canstrat(g_i)\neq F_{i,1-\bit_{i+1}}$.
Since the occurrence record of this edge is then smaller than the occurrence record of the corresponding cycle edges by \Cref{table: Occurrence Records}, such a switch is then applied immediately.
After all improving switches $e$ with $\occrec^{\canstrat}(e)=\floor{(\bit+1)/2}-1$ are applied, the algorithm switches edges $e=(d_{i,j,k},F_{i,j})$ with $\occrec^{\canstrat}(e)=\floor{(\bit+1)/2}$ until there are no open cycle centers anymore.
This behavior is enforced by the tie-breaking rule which then ensures that the cycle center $F_{\nsb,\bit_{\nsb+1}}$ is closed next.
As before, this might \enquote{unlock} the improving switch $(g_{\nsb},\bit_{\nsb+1})$.
By \Cref{table: Occurrence Records}, this switch minimizes the occurrence record among all improving switches and is thus applied which concludes phase $1$.
In any case, phase $2$ then begins if $\nsb>1$, and phase~$3$ begins if $\nsb=1$.

\begin{lemma} \label{lemma: Reaching phase 2}
Let $\canstrat\in\reach{\sigma_0}$ be a canonical strategy for $\bit\in\bitset_n$ with ${\ell(\bit+1)>1}$ having the canonical properties.
After applying finitely many improving switches, the strategy improvement algorithm produces a phase-$2$-strategy $\sigma^{(2)}$ for $\bit$ as described by the corresponding rows of \Cref{table: Properties at start of phase,table: Switches at start of phase}.
\end{lemma}

A schematic example of the phase-$2$-strategy $\sigma^{(2)}$reached when transitioning from $\sigma_3$ to~$\sigma_4$ in the sink game $S_3$ is given in \Cref{figure: Start Phase 2}.

\begin{figure}[ht]
\centering
\begin{tikzpicture}[scale=0.33]

	\def\x{2.25}
	\def\y{1.55}
	\tikzset{
		player/.style={circle, draw, inner sep=1pt, scale=0.4, fill=playblue, text=playblue, minimum size=2.5em},
		random/.style={draw, inner sep=1pt,  scale=0.4, fill=ranred, text=ranred, minimum size=2.5em},
	} 
	
	\node[player, dashed,align=center] (PairSelector2) at  (0*\x,\y) {$b_2$};
	\node[player, dashed,align=center] (Selector1) at  (0*\x,3*\y) {$g_1$};
	\node[player, dashed,align=center] (PairSelector1) at  (0,4*\y) {$b_1$};
	\node[player,align=center] (t) at (0*\x,15*\y) {$t$};
	\draw[->,playblue, very thick ] (t) to[out=60, in=120, looseness=5] (t);
	
	\node[player,align=center] (g1) at (0,0) {$g_1$};
	\node[random,align=center] (F10) at (-3*\x,2*\y) {$F_{1,0}$};
	\node[random,align=center] (F12) at (3*\x,2*\y) {$F_{1,1}$};
	\node[player,align=center] (d101) at (-2*\x,3*\y) {$d_{1,0,1}$};
	\node[player,align=center] (e101) at (-\x,3*\y) {$e_{1,0,1}$};
	\node[player,align=center] (d111) at (2*\x,3*\y) {$d_{1,1,1}$};
	\node[player,align=center] (e111) at (1*\x,3*\y) {$e_{1,1,1}$};
	\node[player,align=center] (d100) at (-2*\x,1*\y) {$d_{1,0,0}$};
	\node[player,align=center] (e100) at (-\x,1*\y) {$e_{1,0,0}$};
	\node[player,align=center] (d110) at (2*\x,1*\y) {$d_{1,1,0}$};
	\node[player,align=center] (e110) at (\x,1*\y) {$e_{1,1,0}$};
	\node[player,align=center] (h10) at (-3*\x,4*\y) {$h_{1,0}$};
	\node[player,align=center] (b1) at (-4*\x,0) {$b_1$};
	\node[player,align=center] (s10) at (-2*\x,4*\y) {$s_{1,0}$};
	\node[player,align=center] (s11) at (2*\x,4*\y) {$s_{1,1}$};
	\node[player,align=center] (h11) at (3*\x,4*\y) {$h_{1,1}$};	
	
	\draw[dotted,  ] (-4.5*\x,4.5*\y)--(3.5*\x,4.5*\y);
	
	\draw[->] (e101)-- (PairSelector2);
	\draw[->, playblue, very thick] (e101)--(Selector1);
	
	\draw[->] (e111)-- (PairSelector2);
	\draw[->, playblue, very thick] (e111)--  (Selector1);
	
	\draw[->] (e100)-- (PairSelector2);
	\draw[->, playblue, very thick] (e100)-- (Selector1);

	\draw[->] (e110)--(PairSelector2);
	\draw[->,playblue, very thick] (e110)--(Selector1);
	
	\draw[->] (s10)--(h10);
	\draw[->, playblue, very thick ] (s10)--(PairSelector1);
	
	\draw[->, playblue, very thick] (s11)--(h11);
	\draw[->] (s11)--(PairSelector1);
	
\begin{scope}[yshift=5*\y cm]
	\node[player, dashed,align=center] (PairSelector1) at  (0,4*\y) {$b_1$};
	\node[player, dashed,align=center] (PairSelector2) at  (0*\x,\y) {$b_2$};
	\node[player, dashed,align=center] (Selector1) at  (0*\x,3*\y) {$g_1$};

	\node[player,align=center] (g2) at (0,0) {$g_2$\\$13$};
	\node[random,align=center] (F20) at (-3*\x,2*\y) {$F_{2,0}$};
	\node[random,align=center] (F22) at (3*\x,2*\y) {$F_{2,1}$};
	\node[player,align=center] (d201) at (-2*\x,3*\y) {$d_{2,0,1}$};
	\node[player,align=center] (e201) at (-\x,3*\y) {$e_{2,0,1}$};
	\node[player,align=center] (d211) at (2*\x,3*\y) {$d_{2,1,1}$};
	\node[player,align=center] (e211) at (1*\x,3*\y) {$e_{2,1,1}$};
	\node[player,align=center] (d200) at (-2*\x,1*\y) {$d_{2,0,0}$};
	\node[player,align=center] (e200) at (-\x,1*\y) {$e_{2,0,0}$};
	\node[player,align=center] (d210) at (2*\x,1*\y) {$d_{2,1,0}$};
	\node[player,align=center] (e210) at (\x,1*\y) {$e_{2,1,0}$};
	\node[player,align=center] (h20) at (-3*\x,4*\y) {$h_{2,0}$};
	\node[player,align=center] (b2) at (-4*\x,0) {$b_2$};
	\node[player,align=center] (s20) at (-2*\x,4*\y) {$s_{2,0}$};
	\node[player,align=center] (s21) at (2*\x,4*\y) {$s_{2,1}$};
	\node[player,align=center] (h21) at (3*\x,4*\y) {$h_{2,1}$};	
	
	\draw[dotted,  ] (-4.5*\x,4.5*\y)--(3.5*\x,4.5*\y);
	
	\draw[->, playblue, very thick ] (d201) to [out=270-30, in=0+30] (F20);
	\draw[->] (d201)-- (e201);
	
	\draw[->] (e201)-- (PairSelector2);
	\draw[->, playblue, very thick] (e201)--(Selector1);
	
	\draw[->, green!70!black, very thick] (d211) to [out=270+30, in=180-30] (F22);
	\draw[->, playblue, very thick  ] (d211)-- (e211);
	
	\draw[->] (e211)-- (PairSelector2);
	\draw[->,playblue, very thick] (e211)--  (Selector1);
	
	\draw[->, playblue, very thick] (d200) to[out=90+30, in=0-30]  (F20);
	\draw[->] (d200)--(e200);
	
	\draw[->] (e200)-- (PairSelector2);
	\draw[->,playblue, very thick] (e200)-- (Selector1);
	
	\draw[->,playblue, very thick] (d210) to[out=90-30, in=180+30]  (F22);
	\draw[->, ,  ] (d210)-- (e210);
	
	\draw[->] (e210)--(PairSelector2);
	\draw[->,playblue, very thick] (e210)--(Selector1);
	
	\draw[->, playblue, very thick] (g2) to[out=180-15, in=270] (F20);
	\draw[->] (g2)to[out=0+15, in=270]  (F22);
	
	\draw[->,playblue , very thick  ] (b2)--(g2);
	
	\draw[->] (F20) to[out=90-30, in=180+30] (d201);
	\draw[->] (F20) to[out=270+30, in=180-30] (d200);
	\draw[->, ranred, very thick] (F20) to[out=90, in=180+30] (s20);
	
	\draw[->, ranred, very thick] (F22) to[out=90+30, in=0-30] (d211);
	\draw[->, ,  ] (F22)to[out=270-30, in=0+30]  (d210);
	\draw[->] (F22) to[out=90, in=0-30](s21);
	
	\draw[->,playblue, very thick] (s20)--(h20);
	\draw[->] (s20)--(PairSelector1);
	
	\draw[->, green!70!black, very thick] (s21)--(h21);
	\draw[->, playblue, very thick] (s21)--(PairSelector1);
	\draw[->, playblue, very thick] (h20) --++(-\x,6*\y)--(t);
\end{scope}

\begin{scope}[yshift=10*\y cm]

	\node[player, dashed,align=center] (PairSelector1) at  (0,4*\y) {$b_1$};
	\node[player, dashed,align=center] (PairSelector2) at  (0*\x,\y) {$b_2$};
	\node[player, dashed,align=center] (Selector1) at  (0*\x,3*\y) {$g_1$};

	\node[player,align=center] (g3) at (0,0) {$g_3$};
	\node[random,align=center] (F30) at (-3*\x,2*\y) {$F_{3,0}$};
	\node[player,align=center] (d301) at (-2*\x,3*\y) {$d_{3,0,1}$};
	\node[player,align=center] (e301) at (-\x,3*\y) {$e_{3,0,1}$};
	\node[player,align=center] (d300) at (-2*\x,1*\y) {$d_{3,0,0}$};
	\node[player,align=center] (e300) at (-\x,1*\y) {$e_{3,0,0}$};
	\node[player,align=center] (h30) at (-3*\x,4*\y) {$h_{3,0}$};
	\node[player,align=center] (b3) at (-4*\x,0) {$b_3$};
	\node[player,align=center] (s30) at (-2*\x,4*\y) {$s_{3,0}$};
	
	\draw[dotted,  ] (-4.5*\x,4.5*\y)--(3.5*\x,4.5*\y);
	
	\draw[->, playblue, very thick] (d301) to [out=270-30, in=0+30] (F30);
	\draw[-> ] (d301)-- (e301);
	
	\draw[->] (e301)-- (PairSelector2);
	\draw[->,  playblue, very thick] (e301)--(Selector1);
	
	\draw[->,  playblue, very thick] (d300) to[out=90+30, in=0-30]  (F30);
	\draw[->, ,  ] (d300)--(e300);
	
	\draw[->] (e300)-- (PairSelector2);
	\draw[->,  playblue, very thick] (e300)-- (Selector1);
	
	\draw[->,  playblue, very thick ] (g3) to[out=180-15, in=270] (F30);
	
	\draw[->, green!70!black, very thick] (b3)--(g3);
	
	\draw[->] (F30) to[out=90-30, in=180+30] (d301);
	\draw[->, ,  ] (F30) to[out=270+30, in=180-30] (d300);
	\draw[->, ranred, very thick] (F30) to[out=90, in=180+30] (s30);
	
	\draw[->,  playblue, very thick] (s30)--(h30);
	\draw[->] (s30)--(PairSelector1);
	\draw[->,  playblue, very thick] (h30) to[out=20,in=180] (t);
	
	\draw[->, playblue, very thick  ] (b3) --++(0,5*\y)--(t);
	
\end{scope}

	\draw[->, green!70!black, very thick] (d101) to [out=270-30, in=0+30] (F10);
	\draw[->,  playblue, very thick] (d101)-- (e101);
	
	\draw[->, playblue, very thick ] (d111) to [out=270+30, in=180-30] (F12);
	\draw[->] (d111)-- (e111);
	
	\draw[->,  playblue, very thick] (d100) to[out=90+30, in=0-30]  (F10);
	\draw[->, ,  ] (d100)--(e100);
	
	\draw[->,  playblue, very thick] (d110) to[out=90-30, in=180+30]  (F12);
	\draw[->] (d110)-- (e110);
	
	\draw[->] (g1) to[out=180-15, in=270] (F10);
	\draw[->,  playblue, very thick ] (g1)to[out=0+15, in=270]  (F12);
	
	\draw[->,  playblue, very thick] (b1)--(g1);
	\draw[->] (b1)--(b2);
	\draw[->] (b2)--(b3);
	
	\draw[->, ranred, very thick] (F10) to[out=90-30, in=180+30] (d101);
	\draw[->, ,  ] (F10) to[out=270+30, in=180-30] (d100);
	\draw[->] (F10) to[out=90, in=180+30] (s10);
	
	\draw[->] (F12) to[out=90+30, in=0-30] (d111);
	\draw[->] (F12)to[out=270-30, in=0+30]  (d110);
	\draw[->, ranred, very thick] (F12) to[out=90, in=0-30](s11);
	
	\draw[->, playblue, very thick] (h11) to[out=150,in=0] (g2);
	\draw[->, playblue, very thick ] (h21) to[out=150,in=0] (g3);
	\draw[->,playblue ,very thick  ] (h10)--(b3);
\end{tikzpicture}
\caption[Example of a phase-$2$-strategy]{The phase-$2$-strategy $\sigma^{(2)}$ calculated when transitioning from $\sigma_3$ to $\sigma_4$ in the sink game $S_3$.
Blue edges represent choices of player $0$, red choices represent choices of player $1$, and green edges represent improving switches. 
For simplification, we omit the labels here.} \label{figure: Start Phase 2}
\end{figure}

We now consider the case $\nsb>1$.
Closing the cycle center $F_{\nsb,\bit_{\nsb+1}}$ at the end of phase~$1$ changes the induced bit state from $\bit$ to $\bit+1$.
This implies that the targets of the entry vertices contained in levels $i\leq\nsb$ and of all upper selection vertices contained in levels $i<\nsb$ need to be changed accordingly.
This is reflected by the set of improving switches containing the edges $(b_{\nsb},g_{\nsb})$ and $(s_{\nsb-1,1},h_{\nsb-1,1})$.
Moreover, cycle edges that were improving for $\canstrat$ but were not applied during phase $1$ remain improving, see \Cref{table: Switches at start of phase}.
These switches were not applied as they have an occurrence record of $\floor{(\bit+1)/2}$, and this large occurrence record guarantees that the algorithm will not apply these switches.
The algorithm thus applies the improving switches involving the entry vertices $b_2$ through $b_{\nsb}$ as well as the improving switches $(s_{i,(\bit+1)_{i+1}},h_{i,(\bit+1)_{i+1}})$ next. 
These switches are  applied until $(b_1,b_2)$ becomes improving.
This switch is however not applied yet, and the algorithm reaches phase $3$.
Since none of these switches needs to be applied if $\nsb=1$, the algorithm directly produces a phase-$3$-strategy after phase $1$ if $\nsb=1$.

\begin{lemma} \label{lemma: Reaching phase 3}
Let $\canstrat\in\reach{\sigma_0}$ be a canonical strategy for $\bit\in\bitset_n$ having the canonical properties.
After applying finitely many improving switches, the strategy improvement algorithm produces a phase-$3$-strategy $\sigma^{(3)}$ for $\bit$ as described by the corresponding rows of \Cref{table: Properties at start of phase,table: Switches at start of phase}.
\end{lemma}

A schematic example of the phase-$3$-strategy~$\sigma^{(3)}$ that is reached when transitioning from $\sigma_3$ to $\sigma_4$ in the sink game $S_3$ is given in \Cref{figure: Start Phase 3}.
Note that this is an example for the case $\nsb>1$ as $4$ is an even number.

\begin{figure}[ht]
\centering
\begin{tikzpicture}[scale=0.33]

	\def\x{2.25}
	\def\y{1.55}
	\tikzset{
		player/.style={circle, draw, inner sep=1pt, scale=0.4, fill=playblue, text=playblue, minimum size=2.5em},
		random/.style={draw, inner sep=1pt,  scale=0.4, fill=ranred, text=ranred, minimum size=2.5em},
	} 
	
	\node[player, dashed,align=center] (PairSelector2) at  (0*\x,\y) {$b_2$};
	\node[player, dashed,align=center] (Selector1) at  (0*\x,3*\y) {$g_1$};
	\node[player, dashed,align=center] (PairSelector1) at  (0,4*\y) {$b_1$};
	\node[player,align=center] (t) at (0*\x,15*\y) {$t$};
	\draw[->,playblue, very thick ] (t) to[out=60, in=120, looseness=5] (t);
	
	\node[player,align=center] (g1) at (0,0) {$g_1$};
	\node[random,align=center] (F10) at (-3*\x,2*\y) {$F_{1,0}$};
	\node[random,align=center] (F12) at (3*\x,2*\y) {$F_{1,1}$};
	\node[player,align=center] (d101) at (-2*\x,3*\y) {$d_{1,0,1}$};
	\node[player,align=center] (e101) at (-\x,3*\y) {$e_{1,0,1}$};
	\node[player,align=center] (d111) at (2*\x,3*\y) {$d_{1,1,1}$};
	\node[player,align=center] (e111) at (1*\x,3*\y) {$e_{1,1,1}$};
	\node[player,align=center] (d100) at (-2*\x,1*\y) {$d_{1,0,0}$};
	\node[player,align=center] (e100) at (-\x,1*\y) {$e_{1,0,0}$};
	\node[player,align=center] (d110) at (2*\x,1*\y) {$d_{1,1,0}$};
	\node[player,align=center] (e110) at (\x,1*\y) {$e_{1,1,0}$};
	\node[player,align=center] (h10) at (-3*\x,4*\y) {$h_{1,0}$};
	\node[player,align=center] (b1) at (-4*\x,0) {$b_1$};
	\node[player,align=center] (s10) at (-2*\x,4*\y) {$s_{1,0}$};
	\node[player,align=center] (s11) at (2*\x,4*\y) {$s_{1,1}$};
	\node[player,align=center] (h11) at (3*\x,4*\y) {$h_{1,1}$};	
	
	\draw[dotted,  ] (-4.5*\x,4.5*\y)--(3.5*\x,4.5*\y);
	
	\draw[->,green!70!black, very thick] (e101)-- (PairSelector2);
	\draw[->, playblue, very thick] (e101)--(Selector1);
	
	\draw[->,green!70!black, very thick] (e111)-- (PairSelector2);
	\draw[->, playblue, very thick] (e111)--  (Selector1);
	
	\draw[->,green!70!black, very thick] (e100)-- (PairSelector2);
	\draw[->, playblue, very thick] (e100)-- (Selector1);

	\draw[->,green!70!black, very thick] (e110)--(PairSelector2);
	\draw[->,playblue, very thick] (e110)--(Selector1);
	
	\draw[->, playblue, very thick] (s10)--(h10);
	\draw[->,] (s10)--(PairSelector1);
	
	\draw[->, playblue, very thick] (s11)--(h11);
	\draw[->] (s11)--(PairSelector1);
	
\begin{scope}[yshift=5*\y cm]
	\node[player, dashed,align=center] (PairSelector1) at  (0,4*\y) {$b_1$};
	\node[player, dashed,align=center] (PairSelector2) at  (0*\x,\y) {$b_2$};
	\node[player, dashed,align=center] (Selector1) at  (0*\x,3*\y) {$g_1$};

	\node[player,align=center] (g2) at (0,0) {$g_2$\\$13$};
	\node[random,align=center] (F20) at (-3*\x,2*\y) {$F_{2,0}$};
	\node[random,align=center] (F22) at (3*\x,2*\y) {$F_{2,1}$};
	\node[player,align=center] (d201) at (-2*\x,3*\y) {$d_{2,0,1}$};
	\node[player,align=center] (e201) at (-\x,3*\y) {$e_{2,0,1}$};
	\node[player,align=center] (d211) at (2*\x,3*\y) {$d_{2,1,1}$};
	\node[player,align=center] (e211) at (1*\x,3*\y) {$e_{2,1,1}$};
	\node[player,align=center] (d200) at (-2*\x,1*\y) {$d_{2,0,0}$};
	\node[player,align=center] (e200) at (-\x,1*\y) {$e_{2,0,0}$};
	\node[player,align=center] (d210) at (2*\x,1*\y) {$d_{2,1,0}$};
	\node[player,align=center] (e210) at (\x,1*\y) {$e_{2,1,0}$};
	\node[player,align=center] (h20) at (-3*\x,4*\y) {$h_{2,0}$};
	\node[player,align=center] (b2) at (-4*\x,0) {$b_2$};
	\node[player,align=center] (s20) at (-2*\x,4*\y) {$s_{2,0}$};
	\node[player,align=center] (s21) at (2*\x,4*\y) {$s_{2,1}$};
	\node[player,align=center] (h21) at (3*\x,4*\y) {$h_{2,1}$};	
	
	\draw[dotted,  ] (-4.5*\x,4.5*\y)--(3.5*\x,4.5*\y);
	
	\draw[->, playblue, very thick ] (d201) to [out=270-30, in=0+30] (F20);
	\draw[->] (d201)-- (e201);
	
	\draw[->,green!70!black, very thick] (e201)-- (PairSelector2);
	\draw[->, playblue, very thick] (e201)--(Selector1);
	
	\draw[->, green!70!black, very thick] (d211) to [out=270+30, in=180-30] (F22);
	\draw[->, playblue, very thick  ] (d211)-- (e211);
	
	\draw[->,green!70!black, very thick] (e211)-- (PairSelector2);
	\draw[->,playblue, very thick] (e211)--  (Selector1);
	
	\draw[->, playblue, very thick] (d200) to[out=90+30, in=0-30]  (F20);
	\draw[->] (d200)--(e200);
	
	\draw[->,green!70!black, very thick] (e200)-- (PairSelector2);
	\draw[->,playblue, very thick] (e200)-- (Selector1);
	
	\draw[->,playblue, very thick] (d210) to[out=90-30, in=180+30]  (F22);
	\draw[->, ,  ] (d210)-- (e210);
	
	\draw[->,green!70!black, very thick] (e210)--(PairSelector2);
	\draw[->,playblue, very thick] (e210)--(Selector1);
	
	\draw[->, playblue, very thick] (g2) to[out=180-15, in=270] (F20);
	\draw[->] (g2)to[out=0+15, in=270]  (F22);
	
	\draw[->] (b2)--(g2);
	
	\draw[->] (F20) to[out=90-30, in=180+30] (d201);
	\draw[->] (F20) to[out=270+30, in=180-30] (d200);
	\draw[->, ranred, very thick] (F20) to[out=90, in=180+30] (s20);
	
	\draw[->, ranred, very thick] (F22) to[out=90+30, in=0-30] (d211);
	\draw[->, ,  ] (F22)to[out=270-30, in=0+30]  (d210);
	\draw[->] (F22) to[out=90, in=0-30](s21);
	
	\draw[->,playblue, very thick] (s20)--(h20);
	\draw[->] (s20)--(PairSelector1);
	
	\draw[->, playblue, very thick] (s21)--(h21);
	\draw[->] (s21)--(PairSelector1);
	\draw[->, playblue, very thick] (h20) --++(-\x,6*\y)--(t);
\end{scope}

\begin{scope}[yshift=10*\y cm]

	\node[player, dashed,align=center] (PairSelector1) at  (0,4*\y) {$b_1$};
	\node[player, dashed,align=center] (PairSelector2) at  (0*\x,\y) {$b_2$};
	\node[player, dashed,align=center] (Selector1) at  (0*\x,3*\y) {$g_1$};

	\node[player,align=center] (g3) at (0,0) {$g_3$};
	\node[random,align=center] (F30) at (-3*\x,2*\y) {$F_{3,0}$};
	\node[player,align=center] (d301) at (-2*\x,3*\y) {$d_{3,0,1}$};
	\node[player,align=center] (e301) at (-\x,3*\y) {$e_{3,0,1}$};
	\node[player,align=center] (d300) at (-2*\x,1*\y) {$d_{3,0,0}$};
	\node[player,align=center] (e300) at (-\x,1*\y) {$e_{3,0,0}$};
	\node[player,align=center] (h30) at (-3*\x,4*\y) {$h_{3,0}$};
	\node[player,align=center] (b3) at (-4*\x,0) {$b_3$};
	\node[player,align=center] (s30) at (-2*\x,4*\y) {$s_{3,0}$};
	
	\draw[dotted,  ] (-4.5*\x,4.5*\y)--(3.5*\x,4.5*\y);
	
	\draw[->, playblue, very thick] (d301) to [out=270-30, in=0+30] (F30);
	\draw[->] (d301)-- (e301);
	
	\draw[->,green!70!black, very thick] (e301)-- (PairSelector2);
	\draw[->,  playblue, very thick] (e301)--(Selector1);
	
	\draw[->,  playblue, very thick] (d300) to[out=90+30, in=0-30]  (F30);
	\draw[->, ,  ] (d300)--(e300);
	
	\draw[->,green!70!black, very thick] (e300)-- (PairSelector2);
	\draw[->,  playblue, very thick] (e300)-- (Selector1);
	
	\draw[->,  playblue, very thick ] (g3) to[out=180-15, in=270] (F30);
	
	\draw[->, playblue, very thick] (b3)--(g3);
	
	\draw[->] (F30) to[out=90-30, in=180+30] (d301);
	\draw[->, ,  ] (F30) to[out=270+30, in=180-30] (d300);
	\draw[->, ranred, very thick] (F30) to[out=90, in=180+30] (s30);
	
	\draw[->,  playblue, very thick] (s30)--(h30);
	\draw[->] (s30)--(PairSelector1);
	\draw[->,  playblue, very thick] (h30) to[out=20,in=180] (t);
	
	\draw[->, ,  ] (b3) --++(0,5*\y)--(t);
	
\end{scope}

	\draw[->, green!70!black, very thick] (d101) to [out=270-30, in=0+30] (F10);
	\draw[->,  playblue, very thick] (d101)-- (e101);
	
	\draw[->, playblue, very thick ] (d111) to [out=270+30, in=180-30] (F12);
	\draw[->] (d111)-- (e111);
	
	\draw[->,  playblue, very thick] (d100) to[out=90+30, in=0-30]  (F10);
	\draw[->, ,  ] (d100)--(e100);
	
	\draw[->,  playblue, very thick] (d110) to[out=90-30, in=180+30]  (F12);
	\draw[->] (d110)-- (e110);
	
	\draw[->] (g1) to[out=180-15, in=270] (F10);
	\draw[->,  playblue, very thick ] (g1)to[out=0+15, in=270]  (F12);
	
	\draw[->,  playblue, very thick] (b1)--(g1);
	\draw[->, green!70!black, very thick] (b1)--(b2);
	\draw[->, playblue, very thick] (b2)--(b3);
	
	\draw[->, ranred, very thick] (F10) to[out=90-30, in=180+30] (d101);
	\draw[->, ,  ] (F10) to[out=270+30, in=180-30] (d100);
	\draw[->] (F10) to[out=90, in=180+30] (s10);
	
	\draw[->] (F12) to[out=90+30, in=0-30] (d111);
	\draw[->] (F12)to[out=270-30, in=0+30]  (d110);
	\draw[->, ranred, very thick] (F12) to[out=90, in=0-30](s11);
	
	\draw[->, playblue, very thick ] (h11) to[out=150,in=0] (g2);
	\draw[->, playblue, very thick ] (h21) to[out=150,in=0] (g3);
	\draw[->, playblue, very thick ] (h10)--(b3);
\end{tikzpicture}
\caption[Example of a phase-$3$-strategy]{The phase-$3$-strategy $\sigma^{(3)}$ calculated when transitioning from $\sigma_3$ to $\sigma_4$ in the sink game $S_3$.
Blue edges represent choices of player $0$, red choices represent choices of player $1$, and green edges represent improving switches. 
For simplification, we omit the labels here.} \label{figure: Start Phase 3}
\end{figure}

When phase $3$ begins, all edges $(e_{*,*,*}, g_1)$ resp. $(e_{*,*,*},b_2)$ become improving, depending on $\nsb$.
The reason is that closing $F_{\nsb,\bit_{\nsb+1}}$ resp. closing this cycle center and updating the spinal path in the levels $1$ to $\nsb$ significantly increases the valuation of $g_1$ resp. $b_2$.
Since the improving switches of the form $(d_{*,*,*},F_{*,*})$ still maximize the occurrence records, the switches involving the escape vertices are applied next.
As all of them have the same occurrence record, all improving switches $(e_{i,j,k},*)$ with $\sigma^{(3)}(d_{i,j,k})=F_{i,j}$ are applied due to the tie-breaking rule.
If the corresponding cycle center is not closed and active, then this application unlockss the improving switch $(d_{i,j,k},e_{i,j,k})$ as this edge allows the cycle vertex to gain access to the very profitable spinal path.

At this point, there is a major difference between the behavior of $S_n$ and $M_n$.
In $S_n$, the application of a switch $(d_{i,j,k},e_{i,j,k})$ does not change the valuation of its cycle center~$F_{i,j}$.
The reason is that player $1$ controls the cycle center and can then react by choosing vertex $d_{i,j,1-k}$, yielding the same valuation as before.
This is also true if the cycle center was closed, since player $1$ chooses the upper selection vertices in both cases.
In $M_n$, the application of the switch $(d_{i,j,k},e_{i,j,k})$ however has an immediate consequence regarding the valuation of $F_{i,j}$.
As the valuation of the cycle center is (roughly) the arithmetic mean of the valuation of its cycle vertices, the increase of the valuations of $d_{i,j,k}$ also increases the valuation of $F_{i,j}$.
This then makes the cycle center $F_{i,j}$ profitable since it grants access to the spinal path.
Most importantly, it enables the upper selection vertex $s_{i,j}$ to use this access by switching to $b_1$ as the path starting in $b_1$ then leads to $F_{i,j}$.
The reason is that the exact way that improving switches are applied in this phase since the tie-breaking rule dictates that switches of higher levels are applied prior to switches of lower levels.

In summary, all switches $(e_{i,j,k},*)$ with $\sigma^{(3)}(d_{i,j,k})=F_{i,j}$ are applied during phase~$3$.
If $F_{i,j}$ is not closed and active, this makes $(d_{i,j,k},e_{i,j,k})$ improving, and this switch is applied next.
In $M_n$, this might also make switches $(s_{i,j},b_1)$ improving if $i<\nsb$ and $j=1-(\bit+1)_{i+1}$, and this switch is then also applied immediately.
Phase $3$ ends with the application of $(b_1,g_1)$ if $\nsb=1$ and $(b_1,b_2)$ if $\nsb>1$.
Depending on whether we consider $M_n$ or $S_n$ and depending on $\nsb$, we then either obtain a phase-$4$-strategy or a phase-$5$-strategy.

\begin{lemma} \label{lemma: Reaching phase 4 or 5}
Let $\canstrat\in\reach{\sigma_0}$ be a canonical strategy for $\bit\in\bitset_n$ having the canonical properties.
After applying finitely many improving switches, the strategy improvement algorithm produces a strategy $\sigma$ with the following properties.
If $\nsb>1$, then $\sigma$ is a phase-$4$-strategy for $\bit$ in $S_n$ and a phase-$5$-strategy for $\bit$ in $M_n$.
If $\nsb=1$, then $\sigma$ is a phase-$5$-strategy for $\bit$.
In any case, $\sigma$ is described by the corresponding rows of  \Cref{table: Properties at start of phase,table: Switches at start of phase}.
\end{lemma}

In the sink game $S_3$ that we consider as an example in this section, the algorithm thus produces a phase-$4$-strategy when transitioning from $\sigma_3$ to $\sigma_4$ which is visualized in \Cref{figure: Start Phase 4}.

\begin{figure}[ht]
\centering
\begin{tikzpicture}[scale=0.33]

	\def\x{2.25}
	\def\y{1.55}
	\tikzset{
		player/.style={circle, draw, inner sep=1pt, scale=0.4, fill=playblue, text=playblue, minimum size=2.5em},
		random/.style={draw, inner sep=1pt,  scale=0.4, fill=ranred, text=ranred, minimum size=2.5em},
	} 
	
	\node[player, dashed,align=center] (PairSelector2) at  (0*\x,\y) {$b_2$};
	\node[player, dashed,align=center] (Selector1) at  (0*\x,3*\y) {$g_1$};
	\node[player, dashed,align=center] (PairSelector1) at  (0,4*\y) {$b_1$};
	\node[player,align=center] (t) at (0*\x,15*\y) {$t$};
	\draw[->,playblue, very thick ] (t) to[out=60, in=120, looseness=5] (t);
	
	\node[player,align=center] (g1) at (0,0) {$g_1$};
	\node[random,align=center] (F10) at (-3*\x,2*\y) {$F_{1,0}$};
	\node[random,align=center] (F12) at (3*\x,2*\y) {$F_{1,1}$};
	\node[player,align=center] (d101) at (-2*\x,3*\y) {$d_{1,0,1}$};
	\node[player,align=center] (e101) at (-\x,3*\y) {$e_{1,0,1}$};
	\node[player,align=center] (d111) at (2*\x,3*\y) {$d_{1,1,1}$};
	\node[player,align=center] (e111) at (1*\x,3*\y) {$e_{1,1,1}$};
	\node[player,align=center] (d100) at (-2*\x,1*\y) {$d_{1,0,0}$};
	\node[player,align=center] (e100) at (-\x,1*\y) {$e_{1,0,0}$};
	\node[player,align=center] (d110) at (2*\x,1*\y) {$d_{1,1,0}$};
	\node[player,align=center] (e110) at (\x,1*\y) {$e_{1,1,0}$};
	\node[player,align=center] (h10) at (-3*\x,4*\y) {$h_{1,0}$};
	\node[player,align=center] (b1) at (-4*\x,0) {$b_1$};
	\node[player,align=center] (s10) at (-2*\x,4*\y) {$s_{1,0}$};
	\node[player,align=center] (s11) at (2*\x,4*\y) {$s_{1,1}$};
	\node[player,align=center] (h11) at (3*\x,4*\y) {$h_{1,1}$};	
	
	\draw[dotted,  ] (-4.5*\x,4.5*\y)--(3.5*\x,4.5*\y);
	
	\draw[->,green!70!black, very thick] (e101)-- (PairSelector2);
	\draw[->, playblue, very thick] (e101)--(Selector1);
	
	\draw[->,playblue, very thick] (e111)-- (PairSelector2);
	\draw[->, ] (e111)--  (Selector1);
	
	\draw[->,playblue, very thick] (e100)-- (PairSelector2);
	\draw[->] (e100)-- (Selector1);

	\draw[->,playblue, very thick] (e110)--(PairSelector2);
	\draw[->] (e110)--(Selector1);
	
	\draw[->, playblue, very thick] (s10)--(h10);
	\draw[->,] (s10)--(PairSelector1);
	
	\draw[->, playblue, very thick] (s11)--(h11);
	\draw[->, green!70!black, very thick] (s11)--(PairSelector1);
	
\begin{scope}[yshift=5*\y cm]
	\node[player, dashed,align=center] (PairSelector1) at  (0,4*\y) {$b_1$};
	\node[player, dashed,align=center] (PairSelector2) at  (0*\x,\y) {$b_2$};
	\node[player, dashed,align=center] (Selector1) at  (0*\x,3*\y) {$g_1$};

	\node[player,align=center] (g2) at (0,0) {$g_2$\\$13$};
	\node[random,align=center] (F20) at (-3*\x,2*\y) {$F_{2,0}$};
	\node[random,align=center] (F22) at (3*\x,2*\y) {$F_{2,1}$};
	\node[player,align=center] (d201) at (-2*\x,3*\y) {$d_{2,0,1}$};
	\node[player,align=center] (e201) at (-\x,3*\y) {$e_{2,0,1}$};
	\node[player,align=center] (d211) at (2*\x,3*\y) {$d_{2,1,1}$};
	\node[player,align=center] (e211) at (1*\x,3*\y) {$e_{2,1,1}$};
	\node[player,align=center] (d200) at (-2*\x,1*\y) {$d_{2,0,0}$};
	\node[player,align=center] (e200) at (-\x,1*\y) {$e_{2,0,0}$};
	\node[player,align=center] (d210) at (2*\x,1*\y) {$d_{2,1,0}$};
	\node[player,align=center] (e210) at (\x,1*\y) {$e_{2,1,0}$};
	\node[player,align=center] (h20) at (-3*\x,4*\y) {$h_{2,0}$};
	\node[player,align=center] (b2) at (-4*\x,0) {$b_2$};
	\node[player,align=center] (s20) at (-2*\x,4*\y) {$s_{2,0}$};
	\node[player,align=center] (s21) at (2*\x,4*\y) {$s_{2,1}$};
	\node[player,align=center] (h21) at (3*\x,4*\y) {$h_{2,1}$};	
	
	\draw[dotted,  ] (-4.5*\x,4.5*\y)--(3.5*\x,4.5*\y);
	
	\draw[->] (d201) to [out=270-30, in=0+30] (F20);
	\draw[->, playblue, very thick ] (d201)-- (e201);
	
	\draw[->,playblue, very thick] (e201)-- (PairSelector2);
	\draw[->] (e201)--(Selector1);
	
	\draw[->, green!70!black, very thick] (d211) to [out=270+30, in=180-30] (F22);
	\draw[->, playblue, very thick  ] (d211)-- (e211);
	
	\draw[->,green!70!black, very thick] (e211)-- (PairSelector2);
	\draw[->,playblue, very thick] (e211)--  (Selector1);
	
	\draw[->] (d200) to[out=90+30, in=0-30]  (F20);
	\draw[->, playblue, very thick] (d200)--(e200);
	
	\draw[->,playblue, very thick] (e200)-- (PairSelector2);
	\draw[->] (e200)-- (Selector1);
	
	\draw[->] (d210) to[out=90-30, in=180+30]  (F22);
	\draw[->,playblue, very thick] (d210)-- (e210);
	
	\draw[->,playblue, very thick] (e210)--(PairSelector2);
	\draw[->] (e210)--(Selector1);
	
	\draw[->, playblue, very thick] (g2) to[out=180-15, in=270] (F20);
	\draw[->] (g2)to[out=0+15, in=270]  (F22);
	
	\draw[->] (b2)--(g2);
	
	\draw[->] (F20) to[out=90-30, in=180+30] (d201);
	\draw[->] (F20) to[out=270+30, in=180-30] (d200);
	\draw[->, ranred, very thick] (F20) to[out=90, in=180+30] (s20);
	
	\draw[->, ranred, very thick] (F22) to[out=90+30, in=0-30] (d211);
	\draw[->, ,  ] (F22)to[out=270-30, in=0+30]  (d210);
	\draw[->] (F22) to[out=90, in=0-30](s21);
	
	\draw[->,playblue, very thick] (s20)--(h20);
	\draw[->, green!70!black, very thick] (s20)--(PairSelector1);
	
	\draw[->, playblue, very thick] (s21)--(h21);
	\draw[->] (s21)--(PairSelector1);
	\draw[->, playblue, very thick] (h20) --++(-\x,6*\y)--(t);
\end{scope}

\begin{scope}[yshift=10*\y cm]

	\node[player, dashed,align=center] (PairSelector1) at  (0,4*\y) {$b_1$};
	\node[player, dashed,align=center] (PairSelector2) at  (0*\x,\y) {$b_2$};
	\node[player, dashed,align=center] (Selector1) at  (0*\x,3*\y) {$g_1$};

	\node[player,align=center] (g3) at (0,0) {$g_3$};
	\node[random,align=center] (F30) at (-3*\x,2*\y) {$F_{3,0}$};
	\node[player,align=center] (d301) at (-2*\x,3*\y) {$d_{3,0,1}$};
	\node[player,align=center] (e301) at (-\x,3*\y) {$e_{3,0,1}$};
	\node[player,align=center] (d300) at (-2*\x,1*\y) {$d_{3,0,0}$};
	\node[player,align=center] (e300) at (-\x,1*\y) {$e_{3,0,0}$};
	\node[player,align=center] (h30) at (-3*\x,4*\y) {$h_{3,0}$};
	\node[player,align=center] (b3) at (-4*\x,0) {$b_3$};
	\node[player,align=center] (s30) at (-2*\x,4*\y) {$s_{3,0}$};
	
	\draw[dotted,  ] (-4.5*\x,4.5*\y)--(3.5*\x,4.5*\y);
	
	\draw[->, playblue, very thick] (d301) to [out=270-30, in=0+30] (F30);
	\draw[->] (d301)-- (e301);
	
	\draw[->,playblue, very thick] (e301)-- (PairSelector2);
	\draw[->] (e301)--(Selector1);
	
	\draw[->,  playblue, very thick] (d300) to[out=90+30, in=0-30]  (F30);
	\draw[->, ,  ] (d300)--(e300);
	
	\draw[->, playblue, very thick] (e300)-- (PairSelector2);
	\draw[->] (e300)-- (Selector1);
	
	\draw[->,  playblue, very thick ] (g3) to[out=180-15, in=270] (F30);
	
	\draw[->, playblue, very thick] (b3)--(g3);
	
	\draw[->] (F30) to[out=90-30, in=180+30] (d301);
	\draw[->, ,  ] (F30) to[out=270+30, in=180-30] (d300);
	\draw[->, ranred, very thick] (F30) to[out=90, in=180+30] (s30);
	
	\draw[->,  playblue, very thick] (s30)--(h30);
	\draw[->] (s30)--(PairSelector1);
	\draw[->,  playblue, very thick] (h30) to[out=20,in=180] (t);
	
	\draw[->, ,  ] (b3) --++(0,5*\y)--(t);
	
\end{scope}

	\draw[->, green!70!black, very thick] (d101) to [out=270-30, in=0+30] (F10);
	\draw[->,  playblue, very thick] (d101)-- (e101);
	
	\draw[-> ] (d111) to [out=270+30, in=180-30] (F12);
	\draw[->, playblue, very thick] (d111)-- (e111);
	
	\draw[->] (d100) to[out=90+30, in=0-30]  (F10);
	\draw[->, playblue, very thick] (d100)--(e100);
	
	\draw[->] (d110) to[out=90-30, in=180+30]  (F12);
	\draw[->, playblue, very thick] (d110)-- (e110);
	
	\draw[->] (g1) to[out=180-15, in=270] (F10);
	\draw[->,  playblue, very thick ] (g1)to[out=0+15, in=270]  (F12);
	
	\draw[->] (b1)--(g1);
	\draw[->, playblue, very thick] (b1)--(b2);
	\draw[->, playblue, very thick] (b2)--(b3);
	
	\draw[->, ranred, very thick] (F10) to[out=90-30, in=180+30] (d101);
	\draw[->, ,  ] (F10) to[out=270+30, in=180-30] (d100);
	\draw[->] (F10) to[out=90, in=180+30] (s10);
	
	\draw[->] (F12) to[out=90+30, in=0-30] (d111);
	\draw[->] (F12)to[out=270-30, in=0+30]  (d110);
	\draw[->, ranred, very thick] (F12) to[out=90, in=0-30](s11);
	
	\draw[->, playblue, very thick ] (h11) to[out=150,in=0] (g2);
	\draw[->, playblue, very thick ] (h21) to[out=150,in=0] (g3);
	\draw[->, playblue, very thick ] (h10)--(b3);
\end{tikzpicture}
\caption[Example of a phase-$4$-strategy]{The phase-$4$-strategy $\sigma^{(4)}$ calculated when transitioning from $\sigma_3$ to $\sigma_4$ in the sink game $S_3$.
Blue edges represent choices of player $0$, red choices represent choices of player $1$, and green edges represent improving switches. 
For simplification, we omit the labels here.} \label{figure: Start Phase 4}
\end{figure}

Consider the case that there is a phase $4$.
This only happens in~$S_n$ and if~$\nsb>1$.
During this phase, the improving switches $(s_{i,j},b_1)$ with $i<\nsb$ and $j=1-(\bit+1)_{i+1}$, which were already applied during phase $3$ in $M_n$, are applied.
The reason that these switches only become improving now in $S_n$ is that the application of $(b_1,b_2)$ at the end of phase~$3$ significantly increases the valuation of $b_1$ as this vertex then enables accessing the spinal path.
The switches are then applied from higher levels to lower levels, so the final improving switch applied in this phase is $(s_{1,1-(\bit+1)_2},b_1)$, resulting in a phase-$5$-strategy.
Most importantly, the algorithm eventually produces a phase-$5$-strategy in any case.

\begin{lemma} \label{lemma: Reaching phase 5}
Let $\canstrat\in\reach{\sigma_0}$ be a canonical strategy for $\bit\in\bitset_n$ having the canonical properties.
After applying finitely many improving switches, the strategy improvement algorithm produces a phase-$5$-strategy $\sigma^{(5)}$ for $\bit$ as described by the corresponding rows of \Cref{table: Properties at start of phase,table: Switches at start of phase}.
\end{lemma}

A schematic example of the phase-$3$-strategy~$\sigma^{(3)}$ that is reached when transitioning from $\sigma_3$ to $\sigma_4$ in the sink game $S_3$ is given in \Cref{figure: Start Phase 5}.
\begin{figure}[ht]
\centering
\begin{tikzpicture}[scale=0.33]

	\def\x{2.25}
	\def\y{1.55}
	\tikzset{
		player/.style={circle, draw, inner sep=1pt, scale=0.4, fill=playblue, text=playblue, minimum size=2.5em},
		random/.style={draw, inner sep=1pt,  scale=0.4, fill=ranred, text=ranred, minimum size=2.5em},
	} 
	
	\node[player, dashed,align=center] (PairSelector2) at  (0*\x,\y) {$b_2$};
	\node[player, dashed,align=center] (Selector1) at  (0*\x,3*\y) {$g_1$};
	\node[player, dashed,align=center] (PairSelector1) at  (0,4*\y) {$b_1$};
	\node[player,align=center] (t) at (0*\x,15*\y) {$t$};
	\draw[->,playblue, very thick ] (t) to[out=60, in=120, looseness=5] (t);
	
	\node[player,align=center] (g1) at (0,0) {$g_1$};
	\node[random,align=center] (F10) at (-3*\x,2*\y) {$F_{1,0}$};
	\node[random,align=center] (F12) at (3*\x,2*\y) {$F_{1,1}$};
	\node[player,align=center] (d101) at (-2*\x,3*\y) {$d_{1,0,1}$};
	\node[player,align=center] (e101) at (-\x,3*\y) {$e_{1,0,1}$};
	\node[player,align=center] (d111) at (2*\x,3*\y) {$d_{1,1,1}$};
	\node[player,align=center] (e111) at (1*\x,3*\y) {$e_{1,1,1}$};
	\node[player,align=center] (d100) at (-2*\x,1*\y) {$d_{1,0,0}$};
	\node[player,align=center] (e100) at (-\x,1*\y) {$e_{1,0,0}$};
	\node[player,align=center] (d110) at (2*\x,1*\y) {$d_{1,1,0}$};
	\node[player,align=center] (e110) at (\x,1*\y) {$e_{1,1,0}$};
	\node[player,align=center] (h10) at (-3*\x,4*\y) {$h_{1,0}$};
	\node[player,align=center] (b1) at (-4*\x,0) {$b_1$};
	\node[player,align=center] (s10) at (-2*\x,4*\y) {$s_{1,0}$};
	\node[player,align=center] (s11) at (2*\x,4*\y) {$s_{1,1}$};
	\node[player,align=center] (h11) at (3*\x,4*\y) {$h_{1,1}$};	
	
	\draw[dotted,  ] (-4.5*\x,4.5*\y)--(3.5*\x,4.5*\y);
	
	\draw[->,green!70!black, very thick] (e101)-- (PairSelector2);
	\draw[->, playblue, very thick] (e101)--(Selector1);
	
	\draw[->,playblue, very thick] (e111)-- (PairSelector2);
	\draw[->, ] (e111)--  (Selector1);
	
	\draw[->,playblue, very thick] (e100)-- (PairSelector2);
	\draw[->] (e100)-- (Selector1);

	\draw[->,playblue, very thick] (e110)--(PairSelector2);
	\draw[->] (e110)--(Selector1);
	
	\draw[->, playblue, very thick] (s10)--(h10);
	\draw[->,] (s10)--(PairSelector1);
	
	\draw[->] (s11)--(h11);
	\draw[->, playblue, very thick] (s11)--(PairSelector1);
	
\begin{scope}[yshift=5*\y cm]
	\node[player, dashed,align=center] (PairSelector1) at  (0,4*\y) {$b_1$};
	\node[player, dashed,align=center] (PairSelector2) at  (0*\x,\y) {$b_2$};
	\node[player, dashed,align=center] (Selector1) at  (0*\x,3*\y) {$g_1$};

	\node[player,align=center] (g2) at (0,0) {$g_2$\\$13$};
	\node[random,align=center] (F20) at (-3*\x,2*\y) {$F_{2,0}$};
	\node[random,align=center] (F22) at (3*\x,2*\y) {$F_{2,1}$};
	\node[player,align=center] (d201) at (-2*\x,3*\y) {$d_{2,0,1}$};
	\node[player,align=center] (e201) at (-\x,3*\y) {$e_{2,0,1}$};
	\node[player,align=center] (d211) at (2*\x,3*\y) {$d_{2,1,1}$};
	\node[player,align=center] (e211) at (1*\x,3*\y) {$e_{2,1,1}$};
	\node[player,align=center] (d200) at (-2*\x,1*\y) {$d_{2,0,0}$};
	\node[player,align=center] (e200) at (-\x,1*\y) {$e_{2,0,0}$};
	\node[player,align=center] (d210) at (2*\x,1*\y) {$d_{2,1,0}$};
	\node[player,align=center] (e210) at (\x,1*\y) {$e_{2,1,0}$};
	\node[player,align=center] (h20) at (-3*\x,4*\y) {$h_{2,0}$};
	\node[player,align=center] (b2) at (-4*\x,0) {$b_2$};
	\node[player,align=center] (s20) at (-2*\x,4*\y) {$s_{2,0}$};
	\node[player,align=center] (s21) at (2*\x,4*\y) {$s_{2,1}$};
	\node[player,align=center] (h21) at (3*\x,4*\y) {$h_{2,1}$};	
	
	\draw[dotted,  ] (-4.5*\x,4.5*\y)--(3.5*\x,4.5*\y);
	
	\draw[->, green!70!black, very thick] (d201) to [out=270-30, in=0+30] (F20);
	\draw[->, playblue, very thick ] (d201)-- (e201);
	
	\draw[->,playblue, very thick] (e201)-- (PairSelector2);
	\draw[->] (e201)--(Selector1);
	
	\draw[->, green!70!black, very thick] (d211) to [out=270+30, in=180-30] (F22);
	\draw[->, playblue, very thick  ] (d211)-- (e211);
	
	\draw[->,green!70!black, very thick] (e211)-- (PairSelector2);
	\draw[->,playblue, very thick] (e211)--  (Selector1);
	
	\draw[->, green!70!black, very thick] (d200) to[out=90+30, in=0-30]  (F20);
	\draw[->, playblue, very thick] (d200)--(e200);
	
	\draw[->,playblue, very thick] (e200)-- (PairSelector2);
	\draw[->] (e200)-- (Selector1);
	
	\draw[->] (d210) to[out=90-30, in=180+30]  (F22);
	\draw[->,playblue, very thick] (d210)-- (e210);
	
	\draw[->,playblue, very thick] (e210)--(PairSelector2);
	\draw[->] (e210)--(Selector1);
	
	\draw[->, playblue, very thick] (g2) to[out=180-15, in=270] (F20);
	\draw[->] (g2)to[out=0+15, in=270]  (F22);
	
	\draw[->] (b2)--(g2);
	
	\draw[->] (F20) to[out=90-30, in=180+30] (d201);
	\draw[->, ranred, very thick] (F20) to[out=270+30, in=180-30] (d200);
	\draw[->] (F20) to[out=90, in=180+30] (s20);
	
	\draw[->, ranred, very thick] (F22) to[out=90+30, in=0-30] (d211);
	\draw[->, ,  ] (F22)to[out=270-30, in=0+30]  (d210);
	\draw[->] (F22) to[out=90, in=0-30](s21);
	
	\draw[->] (s20)--(h20);
	\draw[->,playblue, very thick] (s20)--(PairSelector1);
	
	\draw[->, playblue, very thick] (s21)--(h21);
	\draw[->] (s21)--(PairSelector1);
	\draw[->, playblue, very thick] (h20) --++(-\x,6*\y)--(t);
\end{scope}

\begin{scope}[yshift=10*\y cm]

	\node[player, dashed,align=center] (PairSelector1) at  (0,4*\y) {$b_1$};
	\node[player, dashed,align=center] (PairSelector2) at  (0*\x,\y) {$b_2$};
	\node[player, dashed,align=center] (Selector1) at  (0*\x,3*\y) {$g_1$};

	\node[player,align=center] (g3) at (0,0) {$g_3$};
	\node[random,align=center] (F30) at (-3*\x,2*\y) {$F_{3,0}$};
	\node[player,align=center] (d301) at (-2*\x,3*\y) {$d_{3,0,1}$};
	\node[player,align=center] (e301) at (-\x,3*\y) {$e_{3,0,1}$};
	\node[player,align=center] (d300) at (-2*\x,1*\y) {$d_{3,0,0}$};
	\node[player,align=center] (e300) at (-\x,1*\y) {$e_{3,0,0}$};
	\node[player,align=center] (h30) at (-3*\x,4*\y) {$h_{3,0}$};
	\node[player,align=center] (b3) at (-4*\x,0) {$b_3$};
	\node[player,align=center] (s30) at (-2*\x,4*\y) {$s_{3,0}$};
	
	\draw[dotted,  ] (-4.5*\x,4.5*\y)--(3.5*\x,4.5*\y);
	
	\draw[->, playblue, very thick] (d301) to [out=270-30, in=0+30] (F30);
	\draw[->] (d301)-- (e301);
	
	\draw[->,playblue, very thick] (e301)-- (PairSelector2);
	\draw[->] (e301)--(Selector1);
	
	\draw[->,  playblue, very thick] (d300) to[out=90+30, in=0-30]  (F30);
	\draw[->, ,  ] (d300)--(e300);
	
	\draw[->, playblue, very thick] (e300)-- (PairSelector2);
	\draw[->] (e300)-- (Selector1);
	
	\draw[->,  playblue, very thick ] (g3) to[out=180-15, in=270] (F30);
	
	\draw[->, playblue, very thick] (b3)--(g3);
	
	\draw[->] (F30) to[out=90-30, in=180+30] (d301);
	\draw[->, ,  ] (F30) to[out=270+30, in=180-30] (d300);
	\draw[->, ranred, very thick] (F30) to[out=90, in=180+30] (s30);
	
	\draw[->,  playblue, very thick] (s30)--(h30);
	\draw[->] (s30)--(PairSelector1);
	\draw[->,  playblue, very thick] (h30) to[out=20,in=180] (t);
	
	\draw[->, ,  ] (b3) --++(0,5*\y)--(t);
	
\end{scope}

	\draw[->, green!70!black, very thick] (d101) to [out=270-30, in=0+30] (F10);
	\draw[->,  playblue, very thick] (d101)-- (e101);
	
	\draw[->, green!70!black, very thick ] (d111) to [out=270+30, in=180-30] (F12);
	\draw[->, playblue, very thick] (d111)-- (e111);
	
	\draw[->] (d100) to[out=90+30, in=0-30]  (F10);
	\draw[->, playblue, very thick] (d100)--(e100);
	
	\draw[->, green!70!black, very thick] (d110) to[out=90-30, in=180+30]  (F12);
	\draw[->, playblue, very thick] (d110)-- (e110);
	
	\draw[->] (g1) to[out=180-15, in=270] (F10);
	\draw[->,  playblue, very thick ] (g1)to[out=0+15, in=270]  (F12);
	
	\draw[->] (b1)--(g1);
	\draw[->, playblue, very thick] (b1)--(b2);
	\draw[->, playblue, very thick] (b2)--(b3);
	
	\draw[->, ranred, very thick] (F10) to[out=90-30, in=180+30] (d101);
	\draw[->, ,  ] (F10) to[out=270+30, in=180-30] (d100);
	\draw[->] (F10) to[out=90, in=180+30] (s10);
	
	\draw[->] (F12) to[out=90+30, in=0-30] (d111);
	\draw[->, ranred, very thick] (F12)to[out=270-30, in=0+30]  (d110);
	\draw[->] (F12) to[out=90, in=0-30](s11);
	
	\draw[->, playblue, very thick ] (h11) to[out=150,in=0] (g2);
	\draw[->, playblue, very thick ] (h21) to[out=150,in=0] (g3);
	\draw[->, playblue, very thick ] (h10)--(b3);
\end{tikzpicture}
\caption[Example of a phase-$5$-strategy]{The phase-$5$-strategy $\sigma^{(5)}$ calculated when transitioning from $\sigma_3$ to $\sigma_4$ in the sink game $S_3$.
Blue edges represent choices of player $0$, red choices represent choices of player $1$, and green edges represent improving switches. 
For simplification, we omit the labels here.} \label{figure: Start Phase 5}
\end{figure}

We now discuss phase $5$.
During this phase, the remaining improving switches of the type $(e_{*,*,*},*)$ are applied.
Applying such a switch then forces every cycle center to point towards the spinal path.
But this implies that the valuation of $F_{i,j}$ increases significantly, making cycle edges improving again.
Several of these cycle edges may now have very low occurrence records as their cycle center was closed for a large number of iterations, and are thus applied immediately after being unlocked.
Similarly to phase $1$, this can make edges $(g_i,F_{i,j})$ improving, and these edges are applied immediately if they become improving.
After all switches $(e_{*,*,*},*)$ as well as possible switches $(d_{i,j,k},F_{i,j})$ with low occurrence records and corresponding switches $(g_i,F_{i,j})$ are applied, this yields a canonical strategy for $\bit+1$.
Then, phase $1$ begins again.

All of this is formalized by the following two statements.
Note that we use the expression of the \enquote{next feasible row} as certain phases may not be present in  certain cases.
Thus,  \enquote{the next row} may not always be accurate.

\begin{restatable}{lemma}{ReachingNextPhase} \label{lemma: Reaching next phase}
Let $\canstrat\in\reach{\sigma_0}$ be a canonical strategy for $\bit\in\bitset_n$ having the canonical properties.
Let $\sigma$ be a strategy obtained by applying a sequence of improving switches to~$\canstrat$.
Let $\sigma$ and $I_{\sigma}$ have the properties of row $k$ of \Cref{table: Properties at start of phase} and \ref{table: Switches at start of phase} for some $k\in[5]$.
Then, applying improving switches according to Zadeh's pivot rule and the tie-breaking rule of \Cref{definition: Tie-Breaking exponential} produces a strategy $\sigma'$ that is described by the next feasible rows of \Cref{table: Properties at start of phase,table: Switches at start of phase}.
\end{restatable}

This lemma then enables us to prove the two main theorems.

\begin{restatable}{theorem}{ReachingCanStrat} \label{theorem: Reaching canonical strategy}
Let $\canstrat\in\reach{\sigma_0}$ be a canonical strategy for $\bit\in\bitset_n$ having the canonical properties.
After applying finitely many improving switches according to Zadeh's pivot rule and the tie-breaking rule of \Cref{definition: Tie-Breaking exponential}, the strategy improvement algorithm calculates a strategy $\sigma_{\bit+1}$ with the following properties.
\begin{enumerate}
	\item $I_{\sigma_{\bit+1}}=\{(d_{i,j,k},F_{i,j}):\sigma_{\bit+1}(d_{i,j,k})\neq F_{i,j}\}.$
	\item The occurrence records are described by \Cref{table: Occurrence Records} when interpreted for $\bit+1$.
	\item $\sigma_{\bit+1}$ is a canonical strategy for $\bit+1$ and has Properties (\ref{property: OR1})$_{*,*,*}$ to (\ref{property: OR4})$_{*,*,*}$.
	\item When transitioning from $\canstrat$ to $\sigma_{\bit+1}$, every improving switch is applied at most once.
\end{enumerate}
In particular, $\sigma_{\bit+1}$ has the canonical properties.
\end{restatable}

A schematic example of the canonical strategy~$\sigma_4$ that is reached when transitioning from $\sigma_3$ to $\sigma_4$ in the sink game $S_3$ is given in \Cref{figure: Final Canstrat}.
\begin{figure}[ht]
\centering
\begin{tikzpicture}[scale=0.33]

	\def\x{2.25}
	\def\y{1.55}
	\tikzset{
		player/.style={circle, draw, inner sep=1pt, scale=0.4, fill=playblue, text=playblue, minimum size=2.5em},
		random/.style={draw, inner sep=1pt,  scale=0.4, fill=ranred, text=ranred, minimum size=2.5em},
	} 
	
	\node[player, dashed,align=center] (PairSelector2) at  (0*\x,\y) {$b_2$};
	\node[player, dashed,align=center] (Selector1) at  (0*\x,3*\y) {$g_1$};
	\node[player, dashed,align=center] (PairSelector1) at  (0,4*\y) {$b_1$};
	\node[player,align=center] (t) at (0*\x,15*\y) {$t$};
	\draw[->,playblue, very thick ] (t) to[out=60, in=120, looseness=5] (t);
	
	\node[player,align=center] (g1) at (0,0) {$g_1$};
	\node[random,align=center] (F10) at (-3*\x,2*\y) {$F_{1,0}$};
	\node[random,align=center] (F12) at (3*\x,2*\y) {$F_{1,1}$};
	\node[player,align=center] (d101) at (-2*\x,3*\y) {$d_{1,0,1}$};
	\node[player,align=center] (e101) at (-\x,3*\y) {$e_{1,0,1}$};
	\node[player,align=center] (d111) at (2*\x,3*\y) {$d_{1,1,1}$};
	\node[player,align=center] (e111) at (1*\x,3*\y) {$e_{1,1,1}$};
	\node[player,align=center] (d100) at (-2*\x,1*\y) {$d_{1,0,0}$};
	\node[player,align=center] (e100) at (-\x,1*\y) {$e_{1,0,0}$};
	\node[player,align=center] (d110) at (2*\x,1*\y) {$d_{1,1,0}$};
	\node[player,align=center] (e110) at (\x,1*\y) {$e_{1,1,0}$};
	\node[player,align=center] (h10) at (-3*\x,4*\y) {$h_{1,0}$};
	\node[player,align=center] (b1) at (-4*\x,0) {$b_1$};
	\node[player,align=center] (s10) at (-2*\x,4*\y) {$s_{1,0}$};
	\node[player,align=center] (s11) at (2*\x,4*\y) {$s_{1,1}$};
	\node[player,align=center] (h11) at (3*\x,4*\y) {$h_{1,1}$};	
	
	\draw[dotted,  ] (-4.5*\x,4.5*\y)--(3.5*\x,4.5*\y);
	
	\draw[->,playblue, very thick] (e101)-- (PairSelector2);
	\draw[->] (e101)--(Selector1);
	
	\draw[->,playblue, very thick] (e111)-- (PairSelector2);
	\draw[->, ] (e111)--  (Selector1);
	
	\draw[->,playblue, very thick] (e100)-- (PairSelector2);
	\draw[->] (e100)-- (Selector1);

	\draw[->,playblue, very thick] (e110)--(PairSelector2);
	\draw[->] (e110)--(Selector1);
	
	\draw[->, playblue, very thick] (s10)--(h10);
	\draw[->,] (s10)--(PairSelector1);
	
	\draw[->] (s11)--(h11);
	\draw[->, playblue, very thick] (s11)--(PairSelector1);
	
\begin{scope}[yshift=5*\y cm]
	\node[player, dashed,align=center] (PairSelector1) at  (0,4*\y) {$b_1$};
	\node[player, dashed,align=center] (PairSelector2) at  (0*\x,\y) {$b_2$};
	\node[player, dashed,align=center] (Selector1) at  (0*\x,3*\y) {$g_1$};

	\node[player,align=center] (g2) at (0,0) {$g_2$\\$13$};
	\node[random,align=center] (F20) at (-3*\x,2*\y) {$F_{2,0}$};
	\node[random,align=center] (F22) at (3*\x,2*\y) {$F_{2,1}$};
	\node[player,align=center] (d201) at (-2*\x,3*\y) {$d_{2,0,1}$};
	\node[player,align=center] (e201) at (-\x,3*\y) {$e_{2,0,1}$};
	\node[player,align=center] (d211) at (2*\x,3*\y) {$d_{2,1,1}$};
	\node[player,align=center] (e211) at (1*\x,3*\y) {$e_{2,1,1}$};
	\node[player,align=center] (d200) at (-2*\x,1*\y) {$d_{2,0,0}$};
	\node[player,align=center] (e200) at (-\x,1*\y) {$e_{2,0,0}$};
	\node[player,align=center] (d210) at (2*\x,1*\y) {$d_{2,1,0}$};
	\node[player,align=center] (e210) at (\x,1*\y) {$e_{2,1,0}$};
	\node[player,align=center] (h20) at (-3*\x,4*\y) {$h_{2,0}$};
	\node[player,align=center] (b2) at (-4*\x,0) {$b_2$};
	\node[player,align=center] (s20) at (-2*\x,4*\y) {$s_{2,0}$};
	\node[player,align=center] (s21) at (2*\x,4*\y) {$s_{2,1}$};
	\node[player,align=center] (h21) at (3*\x,4*\y) {$h_{2,1}$};	
	
	\draw[dotted,  ] (-4.5*\x,4.5*\y)--(3.5*\x,4.5*\y);
	
	\draw[->, green!70!black, very thick] (d201) to [out=270-30, in=0+30] (F20);
	\draw[->, playblue, very thick ] (d201)-- (e201);
	
	\draw[->,playblue, very thick] (e201)-- (PairSelector2);
	\draw[->] (e201)--(Selector1);
	
	\draw[->, green!70!black, very thick] (d211) to [out=270+30, in=180-30] (F22);
	\draw[->, playblue, very thick  ] (d211)-- (e211);
	
	\draw[->,playblue, very thick] (e211)-- (PairSelector2);
	\draw[->] (e211)--  (Selector1);
	
	\draw[->, green!70!black, very thick] (d200) to[out=90+30, in=0-30]  (F20);
	\draw[->, playblue, very thick] (d200)--(e200);
	
	\draw[->,playblue, very thick] (e200)-- (PairSelector2);
	\draw[->] (e200)-- (Selector1);
	
	\draw[->, green!70!black, very thick] (d210) to[out=90-30, in=180+30]  (F22);
	\draw[->,playblue, very thick] (d210)-- (e210);
	
	\draw[->,playblue, very thick] (e210)--(PairSelector2);
	\draw[->] (e210)--(Selector1);
	
	\draw[->, playblue, very thick] (g2) to[out=180-15, in=270] (F20);
	\draw[->] (g2)to[out=0+15, in=270]  (F22);
	
	\draw[->] (b2)--(g2);
	
	\draw[->] (F20) to[out=90-30, in=180+30] (d201);
	\draw[->, ranred, very thick] (F20) to[out=270+30, in=180-30] (d200);
	\draw[->] (F20) to[out=90, in=180+30] (s20);
	
	\draw[->, ranred, very thick] (F22) to[out=90+30, in=0-30] (d211);
	\draw[->, ,  ] (F22)to[out=270-30, in=0+30]  (d210);
	\draw[->] (F22) to[out=90, in=0-30](s21);
	
	\draw[->] (s20)--(h20);
	\draw[->,playblue, very thick] (s20)--(PairSelector1);
	
	\draw[->, playblue, very thick] (s21)--(h21);
	\draw[->] (s21)--(PairSelector1);
	\draw[->, playblue, very thick] (h20) --++(-\x,6*\y)--(t);
\end{scope}

\begin{scope}[yshift=10*\y cm]

	\node[player, dashed,align=center] (PairSelector1) at  (0,4*\y) {$b_1$};
	\node[player, dashed,align=center] (PairSelector2) at  (0*\x,\y) {$b_2$};
	\node[player, dashed,align=center] (Selector1) at  (0*\x,3*\y) {$g_1$};

	\node[player,align=center] (g3) at (0,0) {$g_3$};
	\node[random,align=center] (F30) at (-3*\x,2*\y) {$F_{3,0}$};
	\node[player,align=center] (d301) at (-2*\x,3*\y) {$d_{3,0,1}$};
	\node[player,align=center] (e301) at (-\x,3*\y) {$e_{3,0,1}$};
	\node[player,align=center] (d300) at (-2*\x,1*\y) {$d_{3,0,0}$};
	\node[player,align=center] (e300) at (-\x,1*\y) {$e_{3,0,0}$};
	\node[player,align=center] (h30) at (-3*\x,4*\y) {$h_{3,0}$};
	\node[player,align=center] (b3) at (-4*\x,0) {$b_3$};
	\node[player,align=center] (s30) at (-2*\x,4*\y) {$s_{3,0}$};
	
	\draw[dotted,  ] (-4.5*\x,4.5*\y)--(3.5*\x,4.5*\y);
	
	\draw[->, playblue, very thick] (d301) to [out=270-30, in=0+30] (F30);
	\draw[->] (d301)-- (e301);
	
	\draw[->,playblue, very thick] (e301)-- (PairSelector2);
	\draw[->] (e301)--(Selector1);
	
	\draw[->,  playblue, very thick] (d300) to[out=90+30, in=0-30]  (F30);
	\draw[->, ,  ] (d300)--(e300);
	
	\draw[->, playblue, very thick] (e300)-- (PairSelector2);
	\draw[->] (e300)-- (Selector1);
	
	\draw[->,  playblue, very thick ] (g3) to[out=180-15, in=270] (F30);
	
	\draw[->, playblue, very thick] (b3)--(g3);
	
	\draw[->] (F30) to[out=90-30, in=180+30] (d301);
	\draw[->, ,  ] (F30) to[out=270+30, in=180-30] (d300);
	\draw[->, ranred, very thick] (F30) to[out=90, in=180+30] (s30);
	
	\draw[->,  playblue, very thick] (s30)--(h30);
	\draw[->] (s30)--(PairSelector1);
	\draw[->,  playblue, very thick] (h30) to[out=20,in=180] (t);
	
	\draw[->, ,  ] (b3) --++(0,5*\y)--(t);
	
\end{scope}

	\draw[->, green!70!black, very thick] (d101) to [out=270-30, in=0+30] (F10);
	\draw[->,  playblue, very thick] (d101)-- (e101);
	
	\draw[->, green!70!black, very thick ] (d111) to [out=270+30, in=180-30] (F12);
	\draw[->, playblue, very thick] (d111)-- (e111);
	
	\draw[->, green!70!black, very thick] (d100) to[out=90+30, in=0-30]  (F10);
	\draw[->, playblue, very thick] (d100)--(e100);
	
	\draw[->, green!70!black, very thick] (d110) to[out=90-30, in=180+30]  (F12);
	\draw[->, playblue, very thick] (d110)-- (e110);
	
	\draw[->,  playblue, very thick ] (g1) to[out=180-15, in=270] (F10);
	\draw[->] (g1)to[out=0+15, in=270]  (F12);
	
	\draw[->] (b1)--(g1);
	\draw[->, playblue, very thick] (b1)--(b2);
	\draw[->, playblue, very thick] (b2)--(b3);
	
	\draw[->] (F10) to[out=90-30, in=180+30] (d101);
	\draw[->, ranred, very thick ] (F10) to[out=270+30, in=180-30] (d100);
	\draw[->] (F10) to[out=90, in=180+30] (s10);
	
	\draw[->] (F12) to[out=90+30, in=0-30] (d111);
	\draw[->, ranred, very thick] (F12)to[out=270-30, in=0+30]  (d110);
	\draw[->] (F12) to[out=90, in=0-30](s11);
	
	\draw[->, playblue, very thick ] (h11) to[out=150,in=0] (g2);
	\draw[->, playblue, very thick ] (h21) to[out=150,in=0] (g3);
	\draw[->, playblue, very thick ] (h10)--(b3);
\end{tikzpicture}
\caption[The canonical strategy $\sigma_4$ in $S_3$]{The canonical strategy $\sigma_4$ calculated when transitioning from $\sigma_3$ to $\sigma_4$ in the sink game $S_3$.
Blue edges represent choices of player $0$, red choices represent choices of player $1$, and green edges represent improving switches. 
For simplification, we omit the labels here.} \label{figure: Final Canstrat}
\end{figure}

\begin{restatable}{theorem}{MainTheorem} \label{theorem: Main theorem}
Using Zadeh's pivot rule and the tie-breaking rule of \Cref{definition: Tie-Breaking exponential} when applying
\begin{enumerate} 
	\item the strategy improvement algorithm of \cite{VoegeJurdzinski2000} to $S_n$,
	\item the policy iteration algorithm of \cite{Howard1960} to $M_n$,
	\item the simplex algorithm of \cite{Dantzig1951} to the linear program induced by $M_n$
\end{enumerate}
to the game $G_n$/induced linear program, all of size $\mathcal{O}(n)$, requires at least $2^{n}$ iterations for finding the optimal strategy/solution when using $\sigma_0$ as initial strategy.
\end{restatable}

This concludes our informal proof of the exponential lower bound.
We mention here that some of our findings were implemented by Oliver Friedmann.
More precisely, the sink game~$S_n$ was implemented using the \textsc{PGSolver} library \cite{FriedmannLange/PGSolver}, and the strategy improvement algorithm using Zadeh's pivot rule and our tie-breaking rule was applied to examples of small sizes.
We manually verified that all of our results are correct for constructions with up to $10$ levels. 
The full execution of the strategy improvement algorithm for games with~$3$ and~$4$ levels are available online \cite{FriedmannImplemetation}.

\section{Conclusion}

In this paper, we have shown that Zadeh's pivot rule has an exponential worst-case running time in PG, MDP, and Simplex context.
Together with previous results, we now have a complete picture regarding the worst-case performance of all traditional pivot rules (up to tie-breaking).
This means that new pivot rules will have to be introduced and analyzed in order to make further progress towards the question whether an efficient pivot rule exists.
In particular, addressing the following questions might be the next step:
\begin{enumerate}
	\item Can the known lower bound constructions for history-based pivot rules be generalized to eliminate a larger class of pivot rules? 
	\item Can we devise new natural pivot rules that systematically rule out the general ideas behind the recent lower bound constructions for history-based pivot rules?
\end{enumerate}

Another approach to eliminating pivot rules is the analysis of their inherent complexity.
For example, it was shown recently that predicting the behavior of the Simplex Algorithm with Dantzig's original pivot rule is computationally intractable~\cite{DisserSkutella/18,Fearnley2015}.
It would be interesting to investigate whether these results carry over to Zadeh's pivot rule, and, in particular, whether a unified approach emerges that encompasses a broader class of pivot rules with a single construction.

\bibliographystyle{alpha}
\bibliography{00Main.bib}

\newpage

\appendix

\section{Vertex Valuations and Well-Behaved Strategies} \label{section: Vertex Valuations}

We begin by developing characterizations of the vertex valuations.
This requires us to analyze the strategies calculated by the strategy improvement algorithm in more detail.
As it will turn out, all strategies that the algorithm produces have a certain set of properties.
The strategies are thus \enquote{well-behaved}, and the properties represent the way the algorithm interacts with the instance.
These properties drastically simplify the proofs but it is tedious to prove that every strategy is well-behaved.
We thus prove that (i) the initial strategy is well-behaved and (ii) whenever the algorithm applies an improving switch to a well-behaved strategy, the resulting strategy is well-behaved.
In particular, these two statements imply that any strategy calculated by the algorithm is well-behaved.
This concept was not mentioned previously since it is solely used for proving our results and  most of them have no clear intuitive explanation.
We explicitly encoded the properties defining well-behaved properties in the implementation of the sink game provided by Oliver Friedmann, verifying that the produced strategies are indeed all well-behaved.

The properties are summarized in \Cref{table: Well behaved properties}.
We briefly discuss the properties next and introduce an additional set of parameters and additional notation.
The parameters are abbreviations that denote the first level in which $\sigmabar(x_{*})$ is either true or false for $x_*\in\{b_*,s_*,g_*\}$.
More precisely, we define $\minsig{x}\coloneqq\min(\{i\in[n]:\sigmabar(x_i)\}\cup\{n+1\})$ as well as $\minnegsig{x}\coloneqq\min(\{i\in[n]:\neg\sigmabar(x_i)\}\cup\{n+1\})$ where $x\in\{b,s,g\}$.
Furthermore, for a level $i\in[n]$ and a strategy $\sigma$ for $G_n$, we refer to the cycle center $F_{i,\sigmabar(g_i)}$, that is, the cycle center chosen by the selector vertex $g_i$, as the \emph{chosen cycle center of level $i$}. \index{chosen cycle center}\index{cycle center!chosen|see{chosen cycle center}}
Note that the chosen cycle center and the active cycle center of level $i$ do not necessarily coincide.

\subsection{Properties of well-behaved strategies}
We now introduce the properties that all strategies produced by the algorithm have.
The abbreviations used to refer to the properties are typically related to the vertices that they are related to, and similar or closely related properties have similar abbreviations.
Note that several of the properties are in fact implications.
For these properties, we thus demand that the full implication is true.
In all of the following, let $i\in[n], j,k\in\{0,1\}$ be suitable indices.

Let $\sigma$ be a strategy for $G_n$.
Consider a level $i\geq\relbit{\sigma}$ such that $\sigma(b_i)=g_i$. 
Then, a well-behaved strategy has $\sigma(s_{i,j})=h_{i,j}$ where $j=\sigmabar(g_i)$.
Intuitively, this states that in levels above $\relbit{\sigma}$ that represent a bit equal to one, the upper selection vertex is set correctly.
Formally, \begin{property}{S1}i\geq\relbit{\sigma}\wedge\sigma(b_i)=g_i\implies \sigma(s_{i,\sigmabar(g_i)})=h_{i,\sigmabar(g_i)}.\end{property}\vspace*{-1em}

Let $i<\relbit{\sigma}$.
Assume that either $\sigma(b_2)=g_2$ and $i>1$ or that the chosen cycle center of level $i$ is closed.
Then, we demand that the upper selection vertex $s_{i,\sigmabar(g_i)}$ points to $h_{i,\sigmabar(g_i)}$.
Formally, \begin{property}{S2}i < \relbit{\sigma} \wedge ((\sigma(b_2)=g_2 \wedge i > 1) \vee  \sigmabar(d_i))\implies\sigmabar(s_i).\end{property}\vspace*{-1em}

Let $i<\relbit{\sigma}-1$ and  $\sigma(b_i)=\sigma(b_{i+1})$.
Then, we demand that $\sigma(b_{i+1})=b_{i+2}$.
Intuitively, this encodes that entry vertices below level $\relbit{\sigma}-1$ are reset \enquote{from top to bottom}.
Formally, \begin{property} {B1}i<\relbit{\sigma}-1\wedge\sigma(b_{i})=b_{i+1}\implies\sigma(b_{i+1})=b_{i+2}.\end{property}\vspace*{-1em}

Assume that $\relbit{\sigma}\neq 1$ and that the entry vertex of level $\relbit{\sigma}-1$ points towards the next entry vertex.
Then, we demand that the entry vertex of level $\relbit{\sigma}$ points towards its selector vertex.
Formally, \begin{property} {B2} \relbit{\sigma}\neq 1\wedge\sigma(b_{\relbit{\sigma}-1})=b_{\relbit{\sigma}}\implies\sigma(b_{\relbit{\sigma}})=g_{\relbit{\sigma}}.\end{property}\vspace*{-1em}

Consider a level $i$ such that $\sigma(s_{i,1})=h_{i,1}$.
Further assume $\sigma(b_{i+1})=b_{i+2}$.
Then the values of $\sigmabar(g_{i+1})$ and $\sigmabar(b_{i+2})$ do not coincide for well-behaved strategies.
Formally, \begin{property} {B3} \sigma(s_{i,1})=h_{i,1} \wedge \sigma(b_{i+1})=b_{i+2}\implies \sigmabar(g_{i+1}) \neq\sigmabar(b_{i+2}).\end{property}\vspace*{-1em}

Consider some level $i<\relbit{\sigma}$.
For well-behaved strategies, $F_{i,1}$ is the chosen cycle center of level $i$ if and only if $i\neq\relbit{\sigma}-1$.
This encodes the statuses of the selector vertices of levels below $\relbit{\sigma}$.
Formally,  \begin{property} {BR1} i<\relbit{\sigma}\implies [\sigma(g_i)=F_{i,1} \iff i  \neq\relbit{\sigma}- 1].\end{property}\vspace*{-1em}

We demand that the chosen cycle center of a level $i<\relbit{\sigma}$ does not escape towards $g_1$.
Formally, \begin{property} {BR2}i<\relbit{\sigma}\implies\nsigmabar(eg_{i,\sigmabar(g_i)}).\end{property}	\vspace*{-1em}

Consider a level $i$ such that $\sigma(b_i)=g_i$. 
Let either $i>1, \relbit{\sigma}=1$ or let $\sigma(b_2)=g_2$ be equivalent to $\relbit{\sigma}>2$.
In any of these cases, this implies that the cycle center $F_{i,\sigmabar(g_i)}$ is closed for well-behaved strategies.
Intuitively, this gives a list of situations in which the active cycle center of a level corresponding to a bit that should be equal to 1 is already closed.
Formally, \begin{property} {D1}\sigma(b_i)=g_i \wedge (i > 1 \vee \relbit{\sigma}= 1 \vee  (\sigma(b_2)=g_2 \iff \relbit{\sigma}> 2))\implies \sigmabar(d_i).\end{property}\vspace*{-1em}

Let $\sigma(b_2)=g_2$ and $i\in\{2,\dots,\relbit{\sigma}-1\}$.
Then, we demand that the chosen cycle center of level $i$ is closed.
Formally, \begin{property} {D2}\sigma(b_2)=g_2 \wedge (2 \leq i < \relbit{\sigma})\implies\sigmabar(d_i).\end{property}\vspace*{-1em}

Let $\relbit{\sigma}=1, \minsig{b} \leq \minnegsig{s}, \minnegsig{g}$ and $G_n=S_n$.
Then, we demand that the chosen cycle center of level $1$  does not escape towards $b_2$.
Formally, \begin{property} {MNS1}\relbit{\sigma} = 1 \wedge \minsig{b} \leq \minnegsig{s}, \minnegsig{g}\wedge G_n=S_n\implies \nsigmabar(eb_1).\end{property}\vspace*{-1em}

Let $\relbit{\sigma}=1$ and consider a level $i$ such that $i < \minnegsig{g} <\minnegsig{s}, \minsig{b}$. 
Further assume that $G_n=S_n$ implies $\nsigmabar(b_{\minnegsig{g}+1})$.
For a well-behaved strategy, the chosen cycle center of level~$i$ does not escape towards $b_2$.
Formally, \begin{property} {MNS2}\relbit{\sigma} = 1 \wedge i < \minnegsig{g} < \minnegsig{s}, \minsig{b}\wedge[G_n=S_n\implies\nsigmabar(b_{\minnegsig{g}+1})]\implies\neg \sigmabar(eb_i).\end{property}\vspace*{-1em}

Let $\relbit{\sigma}=1$ as well as $i<\minnegsig{s}\leq \minnegsig{g}<\minsig{b}$ and $G_n=M_n$.
Then, we demand that the chosen cycle center of level $i$ is closed.
Formally, \begin{property} {MNS3}\relbit{\sigma}=1	\wedge i <\minnegsig{s}\leq\minnegsig{g}<\minsig{b}\wedge G_n=M_n \implies \sigmabar(d_i).\end{property}\vspace*{-1em}

Assume $\relbit{\sigma}=1$ as well as $\minnegsig{s}\leq\minnegsig{g}<\minsig{b}$.
Then, we demand that the chosen cycle center of level $\minnegsig{g}$ escapes to $b_2$ but not to $g_1$.
Formally, \begin{property} {MNS4}\relbit{\sigma}=1 \wedge\minnegsig{s}\leq\minnegsig{g}<\minsig{b} \implies \sigmabar(eb_{\minnegsig{s}})\wedge\nsigmabar(eg_{\minnegsig{s}}).\end{property}\vspace*{-1em}

Assume $\relbit{\sigma}=1$ as well as $i<\minnegsig{s}<\minsig{b}\leq\minnegsig{g}$ and $G_n=M_n$.
Then, we demand that the chosen cycle center of level $i$ is closed.
Formally, \begin{property} {MNS5}\relbit{\sigma}=1\wedge i<\minnegsig{s}<\minsig{b}\leq\minnegsig{g}\wedge G_n=M_n\implies \sigmabar(d_i).\end{property}\vspace*{-1em}

Assume $\relbit{\sigma}=1$ as well as $\minnegsig{s}<\minsig{b}\leq\minnegsig{g}$.
Then, we demand that the chosen cycle center of level $\minnegsig{s}$ escapes to $b_2$ but not to $g_1$.
Formally, \begin{property} {MNS6}\relbit{\sigma}=1\wedge\minnegsig{s}<\minsig{b}\leq\minnegsig{g}\implies\sigmabar(eb_{\minnegsig{s}})\wedge\nsigmabar(eg_{\minnegsig{s}}).\end{property}\vspace*{-1em}

Consider a cycle center $F_{i,j}$ escaping only to $g_1$.
Let $\relbit{\sigma}=1$.
Then, we demand that the upper selection vertex corresponding to the chosen cycle center of level $i$ escapes to $b_1$.
Formally, \begin{property} {EG1}\sigmabar(eg_{i,j}) \wedge \neg \sigmabar(eb_{i,j}) \wedge \relbit{\sigma} = 1\implies\nsigmabar(s_{i,j}).\end{property}\vspace*{-1em}

Consider a cycle center $F_{i,j}$ escaping only to $g_1$.
Let $\relbit{\sigma}=1$.
Then, we demand that the chosen cycle center of level 1 is closed.
Formally, \begin{property} {EG2}\sigmabar(eg_{i,j}) \wedge \neg \sigmabar(eb_{i,j}) \wedge \relbit{\sigma} = 1\implies\sigmabar(d_1).\end{property}\vspace*{-1em}

Consider a cycle center $F_{i,j}$ escaping only to $g_1$.
Then, we demand that the upper selection vertex corresponding to the chosen cycle center of level $1$ points towards the next level.
Formally, \begin{property} {EG3}\sigmabar(eg_{i,j}) \wedge \neg \sigmabar(eb_{i,j})\implies\sigmabar(s_1).\end{property}\vspace*{-1em}

Consider a cycle center $F_{i,j}$ escaping only to $g_1$.
Let $\relbit{\sigma}=1$.
We demand that this implies $\sigmabar(g_1)=\sigmabar(b_2)$ for well-behaved strategies.
Formally, \begin{property} {EG4}	\sigmabar(eg_{i,j}) \wedge \neg \sigmabar(eb_{i,j}) \wedge \relbit{\sigma} = 1\implies \sigmabar(g_1) = \sigmabar(b_2).\end{property}\vspace*{-1em}

Consider a cycle center $F_{i,j}$ escaping only to $g_1$.
Let $\relbit{\sigma}\neq 1$ and assume that the upper selection vertex of $F_{i,j}$ escapes towards $b_1$.
Then, we demand that $\sigmabar(b_{i+1})=j$, i.e., the entry vertex of the next level is set correctly with respect to level $i$.
Formally, \begin{property} {EG5}\sigmabar(eg_{i,j}) \wedge \neg \sigmabar(eb_{i,j}) \wedge \relbit{\sigma}\neq 1 \wedge \sigmabar(s_{i,j})\implies\sigmabar(b_{i+1}) = j.\end{property}\vspace*{-1em}

Consider a cycle center $F_{i,j}$ escaping only to $b_2$.
Let $\sigma(b_1)=b_2$.
Then $\sigmabar(b_{i+1})\neq j$.
Formally, \begin{property} {EB1}\sigmabar(eb_{i,j}) \wedge \neg\sigmabar(eg_{i,j}) \wedge \sigma(b_1)=g_1\implies\sigmabar(b_{i+1}) \neq j.	\end{property}\vspace*{-1em}

Consider a cycle center $F_{i,0}$ escaping only to $b_2$.
Let $\sigma(s_{i,0})=h_{i,0}$ and $\sigma(b_1)=g_1$.
Then, we demand that $\relbit{\sigma}=i+1$.
Formally, \begin{property} {EB2} \sigmabar(eb_{i,0}) \wedge \neg\sigmabar(eg_{i,0}) \wedge  \sigma(b_1)=g_1 \wedge \sigma(s_{i,0})=h_{i,0}\implies\relbit{\sigma} = i+1.\end{property}\vspace*{-1em}

Consider a cycle center $F_{i,j}$ escaping only to $b_2$.
Let $\sigma(s_{i,j})=h_{i,j}, \sigma(b_1)=g_1$ and $i>1$.
We then demand that the entry vertex of level $2$ does not grant access to this level.
Formally, \begin{property} {EB3}\sigmabar(eb_{i,j}) \wedge \neg\sigmabar(eg_{i,j})  \wedge \sigmabar(s_{i,j})=h_{i,j} \wedge i > 1\wedge  \sigma(b_1)=g_1\implies \sigma(b_2)=b_3.\end{property}\vspace*{-1em}

Consider a cycle center $F_{i,1}$ escaping only to $b_2$.
Let $\sigma(s_{i,1})=h_{i,1}$ and $\sigma(b_1)=g_1$.
Then, we demand that $\relbit{\sigma}>i+1$.
Formally, \begin{property} {EB4}\sigmabar(eb_{i,1}) \wedge \neg\sigmabar(eg_{i,1}) \wedge \sigma(s_{i,1})=h_{i,1} \wedge  \sigma(b_1)=g_1\implies\relbit{\sigma}> i+1.\end{property}\vspace*{-1em}

Consider a cycle center $F_{i,j}$ escaping only to $b_2$.
Let $\sigma(b_1)=g_1$.
Then, we demand that the entry vertex of level $\relbit{\sigma}$ grants access to this level.
Formally, \begin{property} {EB5}\sigmabar(eb_{i,j}) \wedge \neg\sigmabar(eg_{i,j}) \wedge \sigma(b_1)=g_1\implies \sigmabar(b_{\relbit{\sigma}})=g_{\relbit{\sigma}}.\end{property}	\vspace*{-1em}

Consider a cycle center $F_{i,j}$ escaping only to $b_2$.
Let $\relbit{\sigma}>2$.
Then, we demand that the entry vertex of level $2$ does not grant access to this level.
Formally, \begin{property} {EB6}\sigmabar(eb_{i,j}) \wedge \neg\sigmabar(eg_{i,j}) \wedge \relbit{\sigma}> 2\implies\sigma(b_2)=b_3.\end{property}\vspace*{-1em}

Consider a cycle center $F_{i,j}$ that can escape towards both $g_1$ and $b_2$.
Further assume that $\sigma(s_{i,j})=h_{i,j}$.
Then, we demand that $\sigma(b_{i+1})=j$, so the entry vertex of level $i+1$ is set in accordance with the upper selection vertex of level $i$.
Formally, \begin{property} {EBG1}\sigmabar(eb_{i,j}) \wedge \sigmabar(eg_{i,j}) \wedge \sigma(s_{i,j})=h_{i,j}\implies\sigmabar(b_{i+1})=j.\end{property}\vspace*{-1em}

Consider a cycle center $F_{i,j}$ that can escape towards both $g_1$ and $b_2$.
Further assume that $\sigmabar(g_1)=\sigmabar(b_2)$.		
Then, we demand that the upper selection vertex corresponding to the chosen cycle center of level 1 points towards $h_{1,\sigmabar(g_1)}$.
Formally, \begin{property} {EBG2}\sigmabar(eb_{i,j}) \wedge \sigmabar(eg_{i,j})\wedge \sigmabar(g_1) = \sigmabar(b_2)\implies\sigmabar(s_1).\end{property}\vspace*{-1em}

Consider a cycle center $F_{i,j}$ that can escape towards both $g_1$ and $b_2$.
Further assume that $\sigmabar(g_1)=\sigmabar(b_2)$.
Then the chosen cycle center of level 1 has to be closed for well-behaved strategies.
Formally, \begin{property} {EBG3}\sigmabar(eb_{i,j}) \wedge \sigmabar(eg_{i,j}) \wedge \sigmabar(g_1) = \sigmabar(b_2)\implies\sigmabar(d_1).\end{property}	\vspace*{-1em}

Consider a cycle center $F_{i,j}$ that can escape towards both $b_2$ and $g_1$.
Let $F_{1,0}$ be the chosen cycle center of level 1 and $\sigma(b_2)=g_2$.
Then, we demand that $\relbit{\sigma}\leq 2$.
Formally, \begin{property} {EBG4}\sigmabar(eb_{i,j}) \wedge \sigmabar(eg_{i,j}) \wedge \sigma(g_1)=F_{1,0} \wedge \sigma(b_2)=g_2\implies\relbit{\sigma}\leq 2.\end{property}\vspace*{-1em}

Consider a cycle center $F_{i,j}$ that can escape towards both $b_2$ and $g_1$.
Further assume that $\sigmabar(g_1)\neq\sigmabar(b_2)$.
Then, we demand that $\relbit{\sigma}\neq2$.
Formally, \begin{property} {EBG5}\sigmabar(eb_{i,j}) \wedge \sigmabar(eg_{i,j}) \wedge \sigmabar(g_1)\wedge \neg\sigmabar(b_2)\implies\relbit{\sigma}\neq 2.\end{property}\vspace*{-1em}

If the cycle center of level $n$ is closed, then $\sigma(b_n)=g_n$ or $\sigma(b_1)=g_1$ has to hold for wee-behaved strategies.
Formally, \begin{property} {DN1}\sigmabar(d_n)\implies\sigmabar(b_n)\vee\sigmabar(b_1).\end{property}

Let the cycle center of level $n$ is closed or let $F_{i,1}$ be the chosen cycle center for all $i\in[n-1]$.
Then, we demand that there is some level $i\in[n]$ such that $\sigma(b_i)=g_i$.
Formally, \begin{property} {DN2}\sigmabar(d_n)\vee\minnegsig{g}=n\implies\exists i\in[n]\colon\sigmabar(b_i).\end{property}\vspace*{-1em}

\begin{table}
\centering
\small
\renewcommand{\arraystretch}{1.15}
\begin{tabular}{|c|l|l|}
\hline
   & Premise & Conclusion \\\hline
  (\ref{property: S1}) 	& $i \geq \relbit{\sigma} \wedge \sigma(b_i)=g_i$ & $\sigmabar(s_{i})$ \\
  (\ref{property: S2}) 	& $i < \relbit{\sigma} \wedge ((\sigma(b_2)=g_2 \wedge i > 1) \vee  \sigmabar(d_i)\vee \sigma(b_1)=b_2)$ & $\sigmabar(s_i)$ \\\hdashline
  (\ref{property: B1}) 	& $i < \relbit{\sigma}- 1 \wedge \sigma(b_i)=b_{i+1}$ &  $\sigma(b_{i+1})=b_{i+2}$ \\  
  (\ref{property: B2})	&  $\relbit{\sigma}\neq 1 \wedge \sigma(b_{\relbit{\sigma}-1})=b_{\relbit{\sigma}}$ & $\sigma(b_{\relbit{\sigma}})=g_{\relbit{\sigma}}$ \\
  (\ref{property: B3}) & $\sigma(s_{i,1})=h_{i,1} \wedge \sigma(b_{i+1})=b_{i+2}$ & $\sigmabar(g_{i+1}) \neq\sigmabar(b_{i+2})$ \\\hdashline
  (\ref{property: BR1}) & $i<\relbit{\sigma}$ & $\sigma(g_i)=F_{i,1} \hspace*{-2.3pt}\iff \hspace*{-2.3pt}i\neq\relbit{\sigma}- 1$ \\
  (\ref{property: BR2})	&$i<\relbit{\sigma}$	&$\nsigmabar(eg_{i,\sigmabar(g_i)})$\\\hdashline
  (\ref{property: D1}) & $\sigma(b_i)=g_i \wedge (i > 1 \vee \relbit{\sigma}= 1 \vee  (\sigma(b_2)=g_2 \iff \relbit{\sigma}> 2))$ & $\sigmabar(d_i)$ \\
  (\ref{property: D2}) & $\sigma(b_2)=g_2 \wedge (2 \leq i < \relbit{\sigma})$ &  $\sigmabar(d_i)$ \\\hdashline
  (\ref{property: MNS1}) & $\relbit{\sigma} = 1 \wedge \minsig{b} \leq \minnegsig{s}, \minnegsig{g}\wedge G_n=S_n$&$\nsigmabar(eb_1)$ \\
  (\ref{property: MNS2}) & $\relbit{\sigma} = 1 \wedge i < \minnegsig{g} < \minnegsig{s}, \minsig{b}\wedge[G_n=S_n\implies\nsigmabar(b_{\minnegsig{g}+1})]$& $\neg \sigmabar(eb_i)$ \\
  (\ref{property: MNS3}) & $\relbit{\sigma} = 1 \wedge i<\minnegsig{s}\leq\minnegsig{g}<\minsig{b}\wedge G_n=M_n$ &  $\sigmabar(d_i)$\\
  (\ref{property: MNS4}) & $\relbit{\sigma} = 1 \wedge \minnegsig{s}\leq\minnegsig{g}<\minsig{b}$ & $\sigmabar(eb_{\minnegsig{s}})\wedge\nsigmabar(eg_{\minnegsig{s}})$\\
  (\ref{property: MNS5}) & $\relbit{\sigma} = 1 \wedge i < \minnegsig{s}<\minsig{b}\leq\minnegsig{g}\wedge G_n=M_n$ & $\sigmabar(d_i)$\\
  (\ref{property: MNS6}) & $\relbit{\sigma} = 1 \wedge \minnegsig{s}<\minsig{b}\leq\minnegsig{g}$ & $\sigmabar(eb_{\minnegsig{s}})\wedge\nsigmabar(eg_{\minnegsig{s}})$ \\\hdashline
  (\ref{property: EG1}) & $\sigmabar(eg_{i,j}) \wedge \neg \sigmabar(eb_{i,j}) \wedge \relbit{\sigma} = 1$ & $\neg \sigmabar(s_{i,j})$ \\
  (\ref{property: EG2}) & $\sigmabar(eg_{i,j}) \wedge \neg \sigmabar(eb_{i,j}) \wedge \relbit{\sigma} = 1$ & $\sigmabar(d_1)$ \\
  (\ref{property: EG3}) & $\sigmabar(eg_{i,j}) \wedge \neg \sigmabar(eb_{i,j})$ &  $\sigmabar(s_1)$ \\
  (\ref{property: EG4}) & $\sigmabar(eg_{i,j}) \wedge \neg \sigmabar(eb_{i,j}) \wedge \relbit{\sigma} = 1$ & $\sigmabar(g_1) = \sigmabar(b_2)$ \\
  (\ref{property: EG5}) & $\sigmabar(eg_{i,j}) \wedge \neg\sigmabar(eb_{i,j}) \wedge \relbit{\sigma}\neq 1 \wedge \sigmabar(s_{i,j})$ & $\sigmabar(b_{i+1}) = j$ \\\hdashline
  (\ref{property: EB1}) & $\sigmabar(eb_{i,j}) \wedge \neg\sigmabar(eg_{i,j}) \wedge \sigma(b_1)=g_1$ & $\sigmabar(b_{i+1}) \neq j$ \\
  (\ref{property: EB2}) & $\sigmabar(eb_{i,0}) \wedge \neg\sigmabar(eg_{i,0}) \wedge  \sigma(b_1)=g_1 \wedge \sigma(s_{i,0})=h_{i,0}$ & $\relbit{\sigma} = i+1$ \\
  (\ref{property: EB3}) & $\sigmabar(eb_{i,j}) \wedge \neg\sigmabar(eg_{i,j})  \wedge \sigma(s_{i,j})=h_{i,j} \wedge i > 1\wedge  \sigma(b_1)=g_1$ &  $\sigma(b_2)=b_3$  \\
  (\ref{property: EB4}) & $\sigmabar(eb_{i,1}) \wedge \neg\sigmabar(eg_{i,1}) \wedge \sigma(s_{i,1})=h_{i,1} \wedge  \sigma(b_1)=g_1$ & $\relbit{\sigma}> i+1$ \\
  (\ref{property: EB5}) & $\sigmabar(eb_{i,j}) \wedge \neg\sigmabar(eg_{i,j}) \wedge \sigma(b_1)=g_1$ & $\sigma(b_{\relbit{\sigma}})=g_{\relbit{\sigma}}$ \\
  (\ref{property: EB6}) & $\sigmabar(eb_{i,j}) \wedge \neg\sigmabar(eg_{i,j}) \wedge \relbit{\sigma}> 2$ & $\sigma(b_2)=b_3$ \\\hdashline
  (\ref{property: EBG1}) & $\sigmabar(eb_{i,j}) \wedge \sigmabar(eg_{i,j}) \wedge \sigma(s_{i,j})=h_{i,j}$ & $\sigmabar(b_{i+1})=j$ \\
  (\ref{property: EBG2}) & $\sigmabar(eb_{i,j}) \wedge \sigmabar(eg_{i,j})\wedge \sigmabar(g_1) = \sigmabar(b_2)$ & $\sigmabar(s_1)$ \\
  (\ref{property: EBG3}) & $\sigmabar(eb_{i,j}) \wedge \sigmabar(eg_{i,j}) \wedge \sigmabar(g_1) = \sigmabar(b_2)$ & $\sigmabar(d_1)$ \\
  (\ref{property: EBG4}) & $\sigmabar(eb_{i,j}) \wedge \sigmabar(eg_{i,j}) \wedge \sigma(g_1)=F_{1,0} \wedge \sigma(b_2)=g_2$ & $\relbit{\sigma}\leq 2$ \\
  (\ref{property: EBG5}) & $\sigmabar(eb_{i,j}) \wedge \sigmabar(eg_{i,j}) \wedge \sigmabar(g_1)\wedge \sigma(b_2)=b_3$ & $\relbit{\sigma}\neq 2$ \\ \hdashline
  (\ref{property: DN1})	&$\sigmabar(d_n)$	&$\sigma(b_n)=g_n\vee\sigma(b_1)=g_1$ \\
  (\ref{property: DN2})	&$\sigmabar(d_n)\vee\minnegsig{g}=n$	&$\exists i\in[n]:\sigma(b_i)=g_i$\\  \hline
  \end{tabular}
  \caption[Properties that all calculated strategies have.]{Properties that all calculated strategies have.
  A strategy that has all of these properties is called well-behaved.}  \label{table: Well behaved properties}
\end{table}

As mentioned previously, the abbreviations of the properties summarized in \Cref{table: Well behaved properties} are chosen according to configurations of $G_n$ or individual vertices.
An explanation of these names is given in \Cref{table: Names of well-behaved properties}.

\begin{table}[ht]
\centering
\begin{tabular}{|c|c|}\hline
Abbreviation &Explanation: Property involves... \\\hline
(S*)		&upper \textbf{S}election vertex\\
(B*)		&entry vertices \bm{$b_*$}\\
(BR*)	&levels \textbf{B}elow the next \textbf{R}elevant bit $\relbit{\sigma}$\\
(D*)		&cycle vertices $\bm{d}_{*,*,*}$\\
(MNS*)	&some parameter $\bm{\minnegsig{}}$\\
(EG*)	&cycle centers \textbf{E}scaping only to \bm{$g_1$}\\
(EB*)	&cycle centers \textbf{E}scaping only to \bm{$b_2$}\\
(EBG*)	&cycle centers \textbf{E}scaping only \bm{$b_2$} and \bm{$g_1$}\\
(DN*)	&the cycle vertices $\bm{d}_{*,*,*}$ of level \bm{$n$}.\\\hline
\end{tabular}
\caption[Abbreviations used for defining the properties summarized in \Cref{table: Well behaved properties}]{Explanation of the abbreviations used for defining the properties of \Cref{table: Well behaved properties}.} \label{table: Names of well-behaved properties}
\end{table}

\begin{definition}[Well-behaved strategy] \label{definition: Well-behaved strategy}
A strategy $\sigma$ for $G_n$ is \emph{well-behaved}\index{well-behaved}\index{strategy!well-behaved|see{well-behaved}} if it has all properties of \Cref{table: Well behaved properties}.
\end{definition}

\subsection{Derived properties of well-behaved strategies}
We begin by giving results related to the next relevant bit $\relbit{\sigma}$ of a strategy $\sigma$.
We first show that  its definition can be simplified for well-behaved strategies.
We recall here that the set of incorrect levels is defined as $\incorrect{\sigma}\coloneqq\{i\in[n]:\sigmabar(b_{i})\wedge\sigmabar(g_{i})\neq\sigmabar(b_{i+1})\}.$
Note that we interpret expressions of the form $x\wedge x=y$ as $x\wedge (x=y)$, so  the precedence level of \enquote{$=$} and \enquote{$\neq$} is higher than the precedence level of $\wedge$ and $\vee$.
Also, we assume the parameter $n\in\mathbb{N}$ to be sufficiently large and in particular larger than $3$.

\begin{lemma} \label{lemma: Relbit For Well Behaved Strategies}
Let $\sigma\in\reach{\sigma_0}$ have Properties (\ref{property: B2}) and (\ref{property: BR1}) and $\incorrect{\sigma}\neq\emptyset$.
Then there is an index $i>\max\{i'\in\incorrect{\sigma}\}$ with $\sigmabar(b_{i})$ and $\sigmabar(g_i)=\sigmabar(b_{i+1})$.
As a consequence, for arbitrary well-behaved strategies $\sigma\in\reach{\sigma_0}$ it holds that \[\relbit{\sigma}=\begin{cases}\min\{i>\max\{i'\in\incorrect{\sigma}\}:\sigmabar(b_{i})\wedge\sigmabar(g_{i})=\sigmabar(b_{i+1})\}, &\text{if } \incorrect{\sigma}\neq\emptyset,\\
\min\{i\in[n+1]:\sigma(b_{i})=b_{i+1}\}, &\text{if }\incorrect{\sigma}=\emptyset.\end{cases}\]
In particular, $\relbit{\sigma}=n+1$ implies $\incorrect{\sigma}=\emptyset$ for well-behaved strategies $\sigma$.
\end{lemma}

\begin{proof}
Let $\incorrect{\sigma}\neq\emptyset$.
By construction and since we interpret $t$ as $b_{n+1}$, it follows that $\sigmabar(g_n)=0=\sigmabar(b_{n+1})$.
This implies $\max\{i'\in\incorrect{\sigma}\}\leq n-1$ and that indices larger than this maximum exist.
For the sake of a contradiction, assume that there was no $i>\max\{i'\in\incorrect{\sigma}\}$ with $\sigmabar(b_i)$ and $\sigmabar(g_i)=\sigmabar(b_{i+1})$.
Then, $\relbit{\sigma}=n$ by definition and in particular $\relbit{\sigma}\neq 1$.
By the definition of~$\incorrect{\sigma}$ this implies $\sigma(b_n)=t$.
Now let $\max\{i'\in \incorrect{\sigma}\}\neq n-1$. 
Then, \Pref{BR1} and $\relbit{\sigma}=n$ imply $\sigma(g_{n-1})=F_{n-1,0}$.
In particular, $\sigmabar(g_{n-1})=\sigmabar(b_n)$, implying $\sigma(b_{n-1})=b_n$ since we assume $\max\{i'\in \incorrect{\sigma}\}\neq n-1$.
Consequently, $\sigma(b_{\relbit{\sigma}})=g_{\relbit{\sigma}}$ by \Pref{B2} as $\relbit{\sigma}\neq 1$.
But this is a contradiction to $\sigma(b_n)=t$.
Now assume $\max\{i'\in\incorrect{\sigma}\}=n-1$.
Then, $\sigmabar(g_{n-1})\neq\sigmabar(b_n)=0$ by the definition of $\incorrect{\sigma}$.
Thus, $\sigma(g_{n-1})=F_{n-1,1}$.
But, since $n-1=\relbit{\sigma}-1$, we also have $\sigma(g_{n-1})=F_{n-1,0}$ by \Pref{BR1} which is a contradiction.

Hence there is an index $i>\max\{i'\in\incorrect{\sigma}\}$ with $\sigmabar(b_{i})\wedge\sigmabar(g_{i})=\sigmabar(b_{i+1})$.
\end{proof}

Whenever discussing $\relbit{\sigma}$ for a well-behaved strategy $\sigma$, we implicitly use \Cref{lemma: Relbit For Well Behaved Strategies} without explicitly mentioning it.
We now prove that $\sigma(b_1)=b_2$ is equivalent to $\relbit{\sigma}=1$ for well-behaved strategies and deduce similar helpful statements related to $\relbit{\sigma}$.

\begin{lemma} \label{lemma: b1 iff relbit}
Let $\sigma\in\reach{\sigma_0}$ have \Pref{B1} and let $\relbit{\sigma}\neq 1$.
Then $\sigma(b_1)=g_1$.
Consequently, if $\sigma$ has \Pref{B1}, then $\relbit{\sigma}=1$ is equivalent to $\sigma(b_1)=b_2$.
\end{lemma}

\begin{proof}
By the definition of $\relbit{\sigma}$, it holds that $\relbit{\sigma}=1$ implies $\sigma(b_1)=b_2$.
It thus suffices to prove the first part of the statement, so let $\relbit{\sigma}\neq 1$, implying $\relbit{\sigma}>1$.

Let $\incorrect{\sigma}=\emptyset$, implying $\relbit{\sigma}=\min\{i\in[n+1]:\sigma(b_i)=b_{i+1}\}$.
Since $\relbit{\sigma}>1$, it needs to hold that $\sigma(b_1)=g_1$ since the minimum would be attained for $i=1$ otherwise.

Let $\incorrect{\sigma}\neq\emptyset$. 
Then $\relbit{\sigma}=\min\{i'>\max\{i\in\incorrect{\sigma}\}:\sigmabar(b_{i'})\wedge\sigmabar(g_{i'})=\sigmabar(b_{i'+1})\}$.
If $\relbit{\sigma}=2$, then $\max\{i\in\incorrect{\sigma}\}=1$, implying $\sigma(b_1)=g_1$ by the definition of $\incorrect{\sigma}$.
If $\relbit{\sigma}>2$, the contraposition of \Pref{B1} states \[\sigma(b_{i+1})=g_{i+1}\implies [i\geq\relbit{\sigma}-1\vee\sigma(b_i)=g_i].\]
Let $m\coloneqq\max\{i\in\incorrect{\sigma}\}$.
Then, by definition, $m<\relbit{\sigma}$ and $\sigma(b_m)=g_m$.
We thus either have $m-1\geq\relbit{\sigma}-1$ or $\sigma(b_{m-1})=g_{m-1}$.
Since $m-1\geq\relbit{\sigma}-1$ contradicts $m<\relbit{\sigma}$, we have $\sigma(b_{m-1})=g_{m-1}$.
The argument can now be applied iteratively, implying $\sigma(b_1)=g_1$.
\end{proof}

The following statement now shows a deep connection between the choice of $b_{\relbit{\sigma}}$ and the set of incorrect levels if the strategy $\sigma$ has certain properties.

\begin{lemma} \label{lemma: Traits of Relbit}
Let $\sigma\in\reach{\sigma_0}$ have Properties (\ref{property: B1}), (\ref{property: B2}) and (\ref{property: BR1}).
\begin{enumerate}
	\item Let $\incorrect{\sigma}\neq\emptyset$.
		Then $\sigma(b_{i})=g_{i}$ for all $i\leq\max\{i'\in\incorrect{\sigma}\}$ and $\sigma(b_{\relbit{\sigma}})=g_{\relbit{\sigma}}$.
	\item Let $\incorrect{\sigma}=\emptyset$.
		Then $\sigma(b_{i})=g_{i}$ for all $i<\relbit{\sigma}$ and $\sigma(b_{\relbit{\sigma}})=b_{\relbit{\sigma}+1}$.
		In addition, $\relbit{\sigma}>1$ implies $\sigma(b_2)=g_2\Longleftrightarrow\relbit{\sigma}>2$.
\end{enumerate}
Consequently, $\incorrect{\sigma}=\emptyset$ if and only if $\sigma(b_{\relbit{\sigma}})=b_{\relbit{\sigma}+1}$. 
\end{lemma}

\begin{proof}
The last statements follows directly from the first two, so only these are proven.
\begin{enumerate}
	\item Since $\incorrect{\sigma}\neq\emptyset$ implies $\relbit{\sigma}\neq 1$, the first statement follows by the same arguments used in the proof of \Cref{lemma: b1 iff relbit}.
			The second statement follows directly from \Cref{lemma: Relbit For Well Behaved Strategies} as $\sigma$ has \Pref{B2} and  \Pref{BR1}.
	\item The first statement from $\relbit{\sigma}=\min\{i\in[n+1]\colon\sigma(b_i)=b_{i+1}\}$ in this case.
		The second statement follows directly since $\relbit{\sigma}>1$ implies $\sigma(b_1)=g_1$ in this case. \qedhere
\end{enumerate}
\end{proof}

The second statement of \Cref{lemma: Traits of Relbit} yields the following corollary for well-behaved strategies by \Pref{D1}.
This corollary allows us to simplify several proofs regarding valuations of vertices in $M_n$ later.

\begin{corollary} \label{corollary: Simplified MDP Valuation}
Let $\sigma$ be a well-behaved strategy with $\incorrect{\sigma}=\emptyset$.
Then $i<\relbit{\sigma}$ implies $\sigmabar(d_{i})$.
\end{corollary}

Before discussing canonical strategies in general, we provide one more lemma that significantly simplifies several proofs.
It is closely related to properties of the type (MNS*) and proves that several of their assumptions imply useful statements.

\begin{lemma} \label{lemma: Config implied by Aeb}
Let $\sigma\in\reach{\sigma_0}$ have Properties (\ref{property: B1}) and (\ref{property: B3}) and let $\relbit{\sigma}=1$.
\begin{enumerate}
	\item If $\minsig{b}\leq\minnegsig{s},\minnegsig{g}$, then $\minsig{b}=2$.
	\item If $\minnegsig{g}<\minnegsig{s},\minsig{b}$ and $\minnegsig{g}>1$, then $\minnegsig{g}+1=\minsig{b}$.
\end{enumerate} 
\end{lemma}

\begin{proof}
Note that $\relbit{\sigma}=1$ implies $\sigma(b_1)=b_2$.
\begin{enumerate}
	\item By $\sigma(b_1)=b_2$, we have $\minsig{b}\geq 2$.
		Assume $\minsig{b}>2$ and let $i\coloneqq\minsig{b}-2$.
		Then $i<\minnegsig{g}$, implying $\sigma(g_i)=F_{i,1}$.
		In addition, $i<\minnegsig{s}$, implying $\sigma(s_{i,1})=h_{i,1}$.
		Since $i+1<\minsig{b}$, also $\sigma(b_{i+1})=b_{i+2}$.
		Consequently, by \Pref{B3}, $\sigmabar(g_{i+1})\neq\sigmabar(b_{i+2})=1$ since $i+2=\minsig{b}$.
		But this implies $\sigma(g_{i+1})=F_{i,0}$, contradicting $i+1\leq\minnegsig{g}$.
	\item Let $i\coloneqq\minnegsig{g}-1$.
		Then $\sigma(g_i)=F_{i,1}$, implying $\sigma(s_{i,1})=h_{i,1}$ by the choice of $i$.
		Furthermore, $\sigma(b_{i+1})=b_{i+2}$ as $i+1=\minnegsig{g}<\minsig{b}$.
		Consequently, by \Pref{B3}, $0=\sigmabar(g_{i+1})\neq\sigmabar(b_{i+2})$, implying $\sigmabar(b_{i+2})=1$ and thus $\minnegsig{g}+1=i+2=\minsig{b}$.\qedhere
\end{enumerate}
\end{proof}

\subsection{The framework for the vertex valuations}
We now discuss the general framework used for describing and characterizing the vertex valuations.
Whenever referring to valuations, we henceforth add one of three possible upper indices.
If we consider the valuations exclusively in the sink game $S_n$ resp. the Markov decision process $M_n$, we include an upper index $\P$ resp. $\M$.
If the arguments or statements apply to both $S_n$ and $M_n$, then we use the general wildcard symbol $*$.

For most proofs, we do not consider the \emph{real} valuations as described in the main part of this paper but a \enquote{reduced} version, referred to as $\valustar$.
In $S_n$, the motivation for considering reduced valuations is that the game is constructed in such a way that the most significant difference between two vertex valuations will always be unique and typically have a priority larger than six.
That is, vertices of priority three or four will rarely ever be relevant when comparing valuations.
They can thus be ignored in most cases, simplifying the vertex valuations.
In $M_n$, the reduced valuations are motivated differently.
Consider some cycle center $F_{i,j}$.
Its valuation is equal to \[\e\valu_{\sigma}^\M(s_{i,j})+\frac{1-\e}{2}\valu_{\sigma}^\M(d_{i,j,0})+\frac{1-\e}{2}\valu_{\sigma}^\M(d_{i,j,1}).\]
Intuitively, if $F_{i,j}$ is not closed, then the contribution of $s_{i,j}$ to the valuation of $F_{i,j}$ is very likely to be negligible.
However, if $F_{i,j}$ is closed, then $\e\valu_{\sigma}^\M(F_{i,j})=\e\valu_{\sigma}^\M(s_{i,j})$, so $\valu_{\sigma}^\M(F_{i,j})=\valu_{\sigma}^\M(s_{i,j})$ for every $\e>0$.
Thus, defining $\valustar_{\sigma}^\M$ as the limit of $\valu_{\sigma}^\M$ for $\e\to 0$ yields an easier way of calculating valuations as it eliminates terms of order $o(1)$.
There are however several cases in which the real valuations $\valu_{\sigma}^\M$ need to be considered since $\e\valu_{\sigma}^\M(s_{i,j})$ is not always negligible.
This motivation justifies the following definition.

\begin{definition}[Reduced valuation]
Let $v\in V$ and $\sigma\in\reach{\sigma_0}$.
The \emph{reduced valuation}\index{reduced valuation} of~$v$ with respect to $\sigma$ in $S_n$ is $\valustar_\sigma^\P(v)\coloneqq\valu_\sigma^\P(v)\setminus\{v'\in\valu_\sigma^\P(v)\colon\Omega(v')\in\{3,4,6\}\}.$\\
\noindent
In $M_n$, the \emph{reduced valuation} of $v\in V$ with respect to  $\sigma$ is $\valustar_{\sigma}^\M(v)\coloneqq\lim_{\e\to0}\valu_{\sigma}^\M(v).$
\end{definition}

We now introduce a unified notation for reduced valuations.
This enables us to perform several calculations and arguments for $S_n$ and $M_n$ simultaneously.
Since vertex valuations in $S_n$ are sets of vertices, we begin by arguing that valuations in $M_n$ can also be described as sets of vertices, although they are usually defined via edges.

Since $M_n$ is weakly unichain, the reduced valuation of a vertex is typically a path ending in the vertex $t$ with probability $1$ after finitely many steps.
The only exception are cycle centers escaping to both $b_2$ and $g_1$ which we discuss later.
By construction, the reward of any edge leaving a vertex $v$ is $\scalar{v}\coloneqq (-N)^{\Omega(v)}$, where $\Omega(v)$ denotes the priority of~$v$.
If the reduced valuation of a vertex $v$ corresponds to a path $P$ ending in $t$, then the total reward collected along the edges of $P$ in $M_n$  can thus be expressed as $\sum_{v\in P}\scalar{v}$.
This argument does however not apply to cycle centers.
The reduced valuation of a cycle center~$F_{i,j}$ might depend on both the reduced valuations of $g_1$ and $b_2$.
This is the case if $F_{i,j}$ escapes towards both of these vertices using its cycle vertices and corresponding escape vertices.
In this case, the reduced valuation of $F_{i,j}$ is the arithmetic mean of the reduced valuations of $g_1$ and $b_2$.
Since $M_n$ is weakly unichain and by the definition of the reduced valuations, the reduced valuations of these two vertices are disjoint paths ending in $t$.
In particular, the previous interpretation can be applied to both of these vertices and it can thus be extended naturally to the cycle center~$F_{i,j}$.
In summary, the reduced valuation of any vertex can be interpreted as either a single path or a union of two disjoint paths leading to $t$.

\Cref{table: Valuation abbreviations} introduces a unified notation that can be used for discussing vertex valuations in both $S_n$ and $M_n$ simultaneously.
In addition, it defines several subsets of vertices that turn out to be useful when describing vertex valuations.
These sets are, for example, all vertices in a level $i$ that contribute to the valuations of the vertices, or sets that will typically be part of several valuations.
Although the sets contained in this table formally depend on the current strategy, we do not include an index denoting this strategy as it will always be clear from the context.
To simplify this notation, we write $\scalar{v_1,v_2,\dots,v_k}$ to denote $\sum_{i\in[k]}\scalar{v_i}$ for arbitrary sets $\{v_1,\dots,v_k\}$ of vertices.

{\renewcommand{\arraystretch}{2.75}
\small
\begin{table}[ht]
\centering
\begin{tabular}{|l|l|}\hline
$W_i^\P\coloneqq\{g_i,s_{i,\sigmabar(g_i)},h_{i,\sigmabar(g_i)}\}$ 	{\color{white}$\displaystyle\bigcup_1^2$}		&$W_i^\M\coloneqq\scalar{g_i,s_{i,\sigmabar(g_i)},h_{i,\sigmabar(g_i)}}$ 																					\\\hdashline
$L_{i,\ell}^\P\coloneqq\displaystyle\bigcup_{i'=i}^{\ell}\{W_{i'}^P:\sigma(b_{i'})=g_{i'}\}$			&$L_{i,\ell}^\M\coloneqq\displaystyle\sum_{i'=i}^{\ell}\{W_{i'}^\M:\sigma(b_{i'})=g_{i'}\}$																\\\hdashline
$R_{i,\ell}^\P\coloneqq\displaystyle\hspace*{-1pt}\bigcup_{i'=i}^{\relbit{\sigma}-1}\hspace*{-2pt}W_{i'}^P\cup\hspace*{-5pt}\bigcup_{i'=\relbit{\sigma}+1}^{\ell}\hspace*{-5pt}\{W_{i'}^P:\sigma(b_{i'})=g_{i'}\}$	&$R_{i,\ell}^\M\coloneqq\displaystyle\sum_{i'=i}^{\relbit{\sigma}-1}W_{i'}^\M+\hspace*{-5pt}\sum_{i'=\relbit{\sigma}+1}^{\ell}\hspace*{-5pt}\{W_{i'}^\M:\sigma(b_{i'})=g_{i'}\}$	\\\hdashline
$B_{i,\ell}^\P\coloneqq\begin{cases} R_{i,\ell}^P		&\text{if } i<\relbit{\sigma} \text{ and } \sigma(b_i)=g_i\\
						L_{i,\ell}^P		&\text{otherwise} \end{cases}$	
	&$B_{i,\ell}^\M\coloneqq\begin{cases} R_{i,\ell}^\M		&\text{if } i<\relbit{\sigma} \text{ and } \sigma(b_i)=g_i\\
						L_{i,\ell}^\M		&\text{otherwise} \end{cases}$\\\hline
\end{tabular}

\bigskip

\begin{tabular}{|c||c|c|c|c|c|}\hline
unified notation	&$\oplus/\bigoplus_*^*$	&$\ubracket{\cdot}$	&$\prec/\succ$	&0				&$W\subset\valustar_{\sigma}^*(\cdot)$.\\\hline
corr. notation for  $S_n$			&$\cup/\bigcup_*^*$		&$\{\cdot\}$		&$\lhd/\rhd$	&$\emptyset$	&All $w\in W$ are contained in $\valustar_{\sigma}^\P(\cdot)$	\\\hdashline
corr. notation for $M_n$			&$+/\sum_*^*$			&$\scalar{\cdot}$		&$</>$			&0				&All $w\in W$ are summands of $\valustar_{\sigma}^\M(\cdot)$.\\\hline
\end{tabular}
\caption[Notation used for unified arguments and vertex valuations.]{Abbreviations and notation used for unified arguments and vertex valuations.
We also define $L_{i}^*\coloneqq L_{i,n}^*, R_i^*\coloneqq R_{i,n}^*, B_i^*\coloneqq B_{i,n}^*$} \label{table: Valuation abbreviations}
\end{table}
}

It is immediate that $\valustar_{\sigma}^\P(v)\lhd\valustar_{\sigma}^\P(w)$ implies $\valu_{\sigma}^\P(v)\lhd\valu_{\sigma}^\P(w)$. 
This is not however not completely obvious for $M_n$ since it is not clear how much we \enquote{lose} by using the reduced instead of the real valuation.
However, as shown by the following lemma, we only lose a negligible amount of $o(1)$.
Hence, if $\valustar_{\sigma}^\M(v)>\valustar_{\sigma}^\M(w)$ and if the difference between the two terms is sufficiently large, then we can deduce $\valu_{\sigma}^\M(v)>\valu_{\sigma}^\M(w)$.

\begin{restatable}{lemma}{BooundingValuesForMDPs} \label{lemma: Bounding Values For MDPs}
Let $P=\{g_*,s_{*,*},h_{i*,*}\}$ be the set of vertices with priorities in $M_n$.
Let $S,S',P\subseteq P$ be non-empty subsets, let $\sum(S)\coloneqq\sum_{v\in S}\rew{v}$ and define $\sum(S')$ analogously.
\begin{enumerate}
	\item $|\sum(S)|<N^{2n+11}$ and $\e\cdot|\sum(S)|<1$ for every subset $S\subseteq P$, and
	\item $|\max_{v\in S}\rew{v}|<|\max_{v\in S'}\rew{v}|$ if and only if $|\sum(S)|<|\sum(S')|$ .
\end{enumerate}
\end{restatable}

Since $N=7n$ is larger than the number of vertices with priorities, \Cref{lemma: Bounding Values For MDPs} implies that we can use reduced valuations in $M_n$ in the following way.

\begin{corollary} \label{corollary: Reduced Valuations For MDPs}
Let $\sigma\in\reach{\sigma_0}$.
Then $\valustar_{\sigma}^\M(w)>\valustar_{\sigma}^\M(v)$ implies $\valu_{\sigma}^\M(w)>\valu_{\sigma}^\M(v)$.
\end{corollary}

\begin{proof}
Reduced valuations can be represented as sums of powers of $N$.
If the reduced valuations of two vertices differ, then they thus differ by terms of order at least $N$.
However, by \Cref{lemma: Bounding Values For MDPs}, the reduced valuation of a vertex and its real valuation only differ by terms of order $o(1)$ since the difference between the real and the reduced valuation of a vertex is always an expression of the type $\e\cdot|\sum(S)|$ for some subset $S$.
Consequently, if the reduced valuation of $v$ is larger than the reduced valuation of $w$, the same is true for the real valuations.
\end{proof}

It is possible that $\valustar_{\sigma}^\M(w)=\valustar_{\sigma}^\M(v)$ but $\valu_{\sigma}^\M(w)\neq\valu_{\sigma}^\M(v)$.
This case can occur if there are two cycle centers $F_{i,0}$ and $F_{i,1}$ which are in the same state.
Then, the valuations of the corresponding upper selection vertices decide which of these two vertices has the better valuation.
Since the influence of these vertices is however neglected when considering the reduced valuation, we need to investigate the real valuations in such a case.

Before characterizing the vertex valuations, we state the following general statements regarding the terms of \Cref{table: Valuation abbreviations}.
We will not always refer to this lemma when we use it as it is used in nearly all calculations.
However, we want to especially underline the last statement, as this formalizes the intuition that traversing a single level $i$ completely is more beneficial then traversing all levels below level $i$.

\begin{restatable}{lemma}{VVLemma} \label{lemma: VV Lemma}
Let $\sigma\in\reach{\sigma_0}$ be well-behaved.
\begin{enumerate}
	\item Let $\sigma(b_{\relbit{\sigma}})=b_{\relbit{\sigma}+1}$.
		Then $L_i^*\preceq R_i^*$ for all $i\in[n]$ and $L_i^*\prec R_j^*$ for $j<i\leq\relbit{\sigma}$.
	\item Let $\sigma(b_{\relbit{\sigma}})=g_{\relbit{\sigma}}$.
		Then $L_i^*\succeq R_i^*$ for all $i\in[n]$ and $L_i\succ R_j^*$ for $i\leq\relbit{\sigma}$ and $j\in[n]$ and $L_i\oplus\ubracket{g_j}\succ R_{j}$ for $i\leq\relbit{\sigma}$ and $j<\relbit{\sigma}$.
	\item Let $i\geq\relbit{\sigma}>j$.
		Then $R_j^*\prec\ubracket{s_{i,j},h_{i,j}}\oplus L_{i+1}^*$.
	\item For all $i\in[n]$, it holds that $\ubracket{g_i,s_{i,*},h_{i,*}}\succ\bigoplus_{i'<i}W_{i'}^*$ and $L_1^*\prec\ubracket{s_{i,j},h_{i,j}}\oplus L_{i+1}^*$.
\end{enumerate}
\end{restatable}

\subsection{Characterizing vertex valuations}
The remainder of this section is dedicated to explicitly determine the vertex valuations for well-behaved strategies.
Most of the proofs are very technical and are thus deferred to \Cref{appendix: Proofs Exponential}.
We however also provide some proofs here in the main part to show how these statements are proven.
We begin by discussing the valuations of the entry vertices $b_i$ for $i>1$ and of selector vertices $g_i$ when $i<\relbit{\sigma}$ and $\sigma(b_2)=g_2$.

\begin{lemma} \label{lemma: Valuation of b and g if i>1}
Let $\sigma\in\reach{\sigma_0}$ be well-behaved and $i>1$.
Then $\valustar^*_\sigma(b_i)=B_i^*$ and $i<\relbit{\sigma}$ and $\sigma(b_2)=g_2$ imply $\valustar^*_\sigma(g_i)=R_i^*$.
\end{lemma}

\begin{proof}
We prove both statements by backwards induction on $i$ and begin with the first statement.
Let $i=n$ and $\sigma(b_n)\neq g_n$.
Then $\valustar_{\sigma}^*(b_n)=0$.
Since $B_n^*=L_n^*=0$, the statement follows.
Now let $\sigma(b_n)=g_n$.
Then, by \Pref{D1}, $\sigmabar(d_n)$ and $\sigmabar(s_i)$ by \Pref{S2} since $\relbit{\sigma}\leq n$ by \Cref{lemma: Relbit For Well Behaved Strategies}.
Hence $\valustar_{\sigma}^*(b_n)=W_n^*$.
Since $W_n^*=B_n^*$ in this case, the statement follows for both $S_n$ and $M_n$. 

Now let $i<n, i>1,$ and assume that the statement holds for $i+1$.
We show that it holds for $i$ as well.
This part of the proof uses the second statement of the lemma \emph{directly}, i.e., in a non-inductive way.
Since we use the first statement when proving the second \emph{inductively}, the induction is correct.
We distinguish several cases.
\begin{itemize}
	\item Let $\sigma(b_i)=b_{i+1}$ and $\relbit{\sigma}=\min\{i'\colon \sigma(b_{i'})=b_{i'+1}\}$.
		Then, $\relbit{\sigma}\leq i$ and we show $\valustar_\sigma^*(b_i)=L_i^*$.
		By the definition of $\relbit{\sigma}$, there is no $i'\in[n]$ such that $\sigma(b_{i'})=g_{i'}$ and $\sigmabar(b_{i'+1})\neq\sigmabar(g_{i'})$.
		Since the statement holds if $\sigma(b_{i'})=b_{i'+1}$ for all $i'>i$ consider the smallest $i'>i$ with $\sigma(b_{i'})=g_{i'}$.
		Since $i'>i>1$ we obtain $\valustar_\sigma^*(b_{i'})=B_{i'}^*$ by the induction hypotheses.
		Furthermore, $B_{i'}^*=L_{i'}^*$ by $\relbit{\sigma}\leq i<i'$.
		By the choice of~$i'$ we have $\valustar_\sigma^*(b_i)=\valustar_\sigma^*(b_{i+1})=\dots=\valustar_\sigma^*(b_{i'})=L_{i'}^*$ as well as $L_{i'}^*=L_{i'-1}^*=\dots=L_i^*$.
		As $\relbit{\sigma}\leq i$ implies $L_i^*=B_i^*$, we thus have $\valustar_\sigma^*(b_i)=B_i^*.$
	\item  Let $\sigma(b_i)=b_{i+1}$ and $\relbit{\sigma}=\min\{i'>\max\{i\in\incorrect{\sigma}\}\colon\sigmabar(b_{i'})\wedge\sigmabar(g_{i'})=\sigmabar(b_{i'+1})\}$.
		Assume $i\geq\relbit{\sigma}$.
		Then $\valustar_{\sigma}^*(b_i)=\valustar_\sigma^*(b_{i+1})=B_{i+1}^*$ by the induction hypotheses and $B_{i+1}^*=L_{i+1}^*$ by the choice of $i$.
		Since $B_i^*=L_i^*$ and $L_i^*=L_{i+1}^*$ as $\sigma(b_i)=b_{i+1}$, the statement follows.
		Hence assume $i<\relbit{\sigma}$.
		Since $\sigma(b_i)=b_{i+1}$ and since $\sigma$ is well-behaved, \Pref{B1} yields $\sigma(b_{i'})=b_{i'+1}$ for all $i'\in\{i,\dots,\relbit{\sigma}-1\}$.
		By the induction hypothesis, we thus have $\valustar_{\sigma}^*(b_{i'})=B_{i'}^*=L_{i'}^*$ for these indices.
		In particular, $\valustar_{\sigma}^*(b_{i+1})=L_{i+1}^*$.
		Note that this also holds for the case $i+1=\relbit{\sigma}$.
		Since $\sigma(b_{i})=b_{i+1}$ implies $L_i^*=L_{i+1}^*$ and $\valustar_{\sigma}^*(b_i)=\valustar_{\sigma}^*(b_{i+1})$, this then yields $\valustar_{\sigma}^*(b_{i})=L_{i+1}^*=L_i^*=B_i^*$.
	\item Let $\sigma(b_i)=g_i$ and $i\geq\relbit{\sigma}$.
		As before, the induction hypothesis yields $\valustar_\sigma^*(b_{i'})=L_{i'}^*$ for all $i'>i$.
		Let $j\coloneqq\sigmabar(g_i)$.
		Then, $F_{i,j}$ is closed by \Pref{D1} since $i>1$.
		Thus  $\valustar_{\sigma}^*(F_{i,j})=\valustar_{\sigma}^*(s_{i,j})$.
		In addition, $\sigma(s_{i,j})=h_{i,j}$ by \Pref{S1}.
		Since $i\geq\relbit{\sigma}$ and $\sigma(b_i)=g_i$ we have $j=\sigmabar(g_i)=\sigmabar(b_{i+1})$.
		By construction we thus have \[\valustar_\sigma^*(b_i)=\valustar_\sigma^*(g_i)=W_i^*\oplus\valustar_\sigma^*(b_{i+1})=W_i^*\oplus L_{i+1}^*=L_i^*=B_i^*.\]
	\item Finally, let $\sigma(b_i)=g_i$ and $i<\relbit{\sigma}$.
		Using the contraposition of \Pref{B1} we obtain that either $i-1\geq \relbit{\sigma}-1$ or $\sigma(b_{i-1})=g_{i-1}$.
		Since $i-1\geq\relbit{\sigma}-1$ contradicts $i<\relbit{\sigma}$, it follows that $\sigma(b_{i-1})=g_{i-1}$.
		Applying this statement inductively then yields $\sigma(b_2)=g_2$.
		Using the second statement of this lemma \emph{directly} we then obtain $\valustar_\sigma^*(b_i)=\valustar_\sigma^*(g_i)=R_i^*=B_i^*.$
\end{itemize}

We now show that $i<\relbit{\sigma}$ and $\sigma(b_2)=g_2$ imply $\valustar^*_\sigma(g_i)=R_i^*$.
This proof uses the first statement in an inductive way.
The statement is shown by backwards induction on $i$, so let $i=\relbit{\sigma}-1$.
Then $\sigma(g_i)=F_{i,0}$ by \Pref{BR1} and $\valustar_{\sigma}^*(F_{i,0})=\valustar_{\sigma}^*(s_{i,0})$ since $F_{i,0}$ is closed by \Pref{D2}.
Also, $\sigma(s_{i,0})=h_{i,0}$ by \Pref{S2}.
By construction, and using the first statement inductively we obtain \[\valustar_\sigma^*(g_i)=W_i^*\oplus\valustar_\sigma^*(b_{i+2})=W_i^*\oplus B_{i+2}^*=W_{\relbit{\sigma}-1}^*\oplus B_{\relbit{\sigma}+1}^*=W_{\relbit{\sigma}-1}^*\oplus L_{\relbit{\sigma}+1}^*=R_{\relbit{\sigma}-1}^*.\]

Let $i<\relbit{\sigma}-1$. 
By Properties (\ref{property: BR1}) and (\ref{property: D2}), $\sigma(g_i)=F_{i,1}$ and $\valustar_{\sigma}^*(F_{i,1})=\valustar_{\sigma}^*(s_{i,1})$.
By \Pref{S2}, also $\sigma(s_{i,1})=h_{i,1}$. 
By construction, the induction hypotheses thus yields $\valustar_\sigma^*(g_i)=W_i^*\oplus\valustar_\sigma^*(g_{i+1})=R_{i+1}^*\cup W_i^*=R_i^*=B_i^*.$
\end{proof}

The next lemma shows how the valuation of $g_i$ might change if the additional requirements used in the second statement of \Cref{lemma: Valuation of b and g if i>1} are not met resp.~if $i=1$.
As its proof is rather involved, requires some case distinctions but uses again a backwards induction and \ref{lemma: Valuation of b and g if i>1}, its proof is deferred to the appendix.

\begin{restatable}{lemma}{ValuationOfGIfLevelSmall} \label{lemma: Valuation of g if level small}
Let $\sigma\in\reach{\sigma_0}$ be well-behaved and $i<\relbit{\sigma}$.
Then $\valustar_\sigma^\P(g_i)=R_i^\P$ and\[\valustar_\sigma^\M(g_i)=\begin{cases} 
B_2^\M+\displaystyle\sum_{j=i}^{k-1}W_j^\M+\rew{g_k},&\text{ if }k\coloneqq\min\{k\geq i\colon\neg\sigmabar(d_k)\}<\relbit{\sigma}\\
\valustar_{\sigma}^\M(g_i)=R_i^\M,&\text{ otherwise}.\end{cases}\]
\end{restatable}

This lemma can now be used to generalize the first statement of \Cref{lemma: Valuation of b and g if i>1}.

\begin{lemma} \label{lemma: Valuation of b}
Let $\sigma\in\reach{\sigma_0}$ be well-behaved.
Then $\valustar_\sigma(b_i)^\P=B_i^\P$ for all $i\in[n]$ and $\valustar_{\sigma}^\M(b_i)=B_i^\M$ for all $i\in\{2,\dots,n\}$.
Furthermore, \[\valustar_{\sigma}^\M(b_1)=\begin{cases}B_2^\M+\displaystyle\sum_{j=1}^{k-1}W_j^\M+\rew{g_k}, &\text{ if }k\coloneqq\min\{i\geq1\colon\nsigmabar(d_i)\}<\relbit{\sigma},\\
	B_1^\M, &\text{ otherwise}.\end{cases}\]
\end{lemma}

\begin{proof}
The case $i>1$ follows by \Cref{lemma: Valuation of b and g if i>1}.
It therefore suffices to consider the case $i=1$.
Let $\sigma(b_1)=b_2$.
Then $\relbit{\sigma}=1$ by \Cref{lemma: b1 iff relbit}.
Therefore, by \Cref{lemma: Valuation of b and g if i>1}, we have $\valustar_\sigma^*(b_2)=B_2^*=L_2^*$. 
Since $\sigma(b_1)=b_2$ implies $L_2^*=L_1^*$ we thus obtain $\valustar_\sigma^*(b_1)=\valustar_\sigma^*(b_2)=B_2^*=L_2^*=L_1^*=B_1^*.$

Assume $\sigma(b_1)=g_1$.
Then $\relbit{\sigma}>1$ by \Cref{lemma: b1 iff relbit}.
Consider the case $G_n=S_n$ first.
Then, $\valustar_\sigma^\P(g_1)=R_1^\P$ by \Cref{lemma: Valuation of g if level small}.
Hence, since $i=1<\relbit{\sigma}$ and $\sigma(b_1)=g_1$ it holds that $B_1^\P=R_1^\P$.
Thus $\valustar_\sigma^\P(b_1)=\valustar_\sigma^\P(g_1)=R_1^\P=B_1^\P.$

Consider the case $G_n=M_n$ next.
If $\valustar_\sigma^\M(g_1)=R_1^\M$, the statement follows by the same arguments used for the case $G_n=S_n$.
Hence let $k\coloneqq\min\{i\geq1\colon\neg\sigmabar(d_i)\}<\relbit{\sigma}$.
Then, since $\sigma(b_1)=g_1$ implies $\valustar_{\sigma}^\M(b_1)=\valustar_{\sigma}^\M(g_1)$, \Cref{lemma: Valuation of g if level small} implies  the statement.
\end{proof}

We thus completely characterized the valuation of all vertices $b_i$.
The next vertex valuation we discuss is the valuation of $g_1$ for the special case of $\relbit{\sigma}=1$.
As the vertex valuation of this vertex is rather complex and the proof requires several case distinctions, we defer it to \Cref{appendix: Proofs Exponential}.

As always, we identify $b_i$ for $i>n$ with $t$ for convenience of notation.

\begin{restatable}{lemma}{ValuationOfGOne} \label{lemma: Valuation of g1}
Let $\relbit{\sigma}=1$ and $m\coloneqq \min\{\minnegsig{g},\minnegsig{s}\}$.
Then \[\valustar_\sigma^*(g_1)=\begin{cases}
\rew{g_1}+\valustar_{\sigma}^{\M}(b_2), 																												&\text{if }\minsig{b}\leq\minnegsig{s},\minnegsig{g}\wedge G_n=M_n\wedge\nsigmabar(d_1),\\
W_1^*\oplus\valustar_{\sigma}^*(b_2), 																													&\text{if }	\minsig{b}\leq\minnegsig{s},\minnegsig{g},\\
																																																		&\phantom{if }\wedge(G_n=S_n\vee[G_n=M_n\wedge\sigmabar(d_1)]),\\
\smash{\bigoplus\limits_{i'=1}^{m}W_{i'}^*\oplus\valustar_\sigma^*(b_{\minnegsig{g}+2})}	&\text{if }\minnegsig{g}<\minnegsig{s},\minsig{b},\\
																																																		&\phantom{if }\wedge[(\sigmabar(b_{\minnegsig{g}+1})\wedge G_n=S_n)\vee\nsigmabar(eb_{\minnegsig{g}})],\\
\bigoplus\limits_{i'=1}^{m-1}W_{i'}^*\oplus \ubracket{g_m}\oplus\valustar_\sigma^*(b_2) 		&\text{otherwise}.
\end{cases}\]
\end{restatable}

The next vertex valuation that we investigate in detail is the valuation of the vertices~$F_{i,j}$, i.e., of the cycle centers.
For these vertices we need to distinguish between the sink game $S_n$ and the Markov decision process $M_n$.
We begin with case $G_n=M_n$ as the corresponding statement follows directly from the definition of $\valustar_{\sigma}^\M$.

\begin{lemma} \label{lemma: Exact Behavior Of Random Vertex}
Let $G_n=M_n$.
Let $\sigma\in\reach{\sigma_0}$ be well-behaved and $i\in[n], j\in\{0,1\}$.
Then \[\valustar_{\sigma}^\M(F_{i,j})=\begin{cases}
\valustar_{\sigma}^\M(s_{i,j}),												&\text{if } \sigmabar(d_{i,j}),	\\
\valustar_{\sigma}^\M(g_1),													&\text{if } \sigmabar(eg_{i,j})\wedge\neg\sigmabar(eb_{i,j}),	\\
\valustar_{\sigma}^\M(b_2),													&\text{if } \sigmabar(eb_{i,j})\wedge\neg\sigmabar(eg_{i,j}),	\\
\frac{1}{2}\valustar_{\sigma}^\M(g_1)+\frac{1}{2}\valustar_{\sigma}^\M(b_2),	&\text{if } \sigmabar(eg_{i,j})\wedge\sigmabar(eb_{i,j}).
\end{cases}\]
\end{lemma}

The exact behavior of player $1$ in the sink game $S_n$ requires a more sophisticated analysis.
The reason is that the behavior of player $1$ very much depends on the configuration of the complete counter and the exact setting of several vertices in different levels.
In particular, depending on the setting of the cycle vertices and the upper selection vertex of a cycle center, the valuations of $d_{i,j,0}, d_{i,j,1}$ $s_{i,j}$ can be completely different.
Consequently, player~$1$  can theoretically choose from up to three different valuations.
As the player always minimizes the valuation of the vertices, this requires us to analyze and compare a lot of valuations exactly.
We thus do not provide its proof here but defer it to \Cref{appendix: Proofs Exponential}.

\begin{restatable}{lemma}{ExactBehaviorOfCounterstrategy} \label{lemma: Exact Behavior Of Counterstrategy}
Let $G_n=S_n$.
Let $\sigma\in\reach{\sigma_0}$ be well-behaved and $i\in[n], j\in\{0,1\}$.
Then \[\valustar_\sigma^\P(F_{i,j})=\begin{cases}
\valustar_\sigma^\P(s_{i,j}),				&\text{if }\sigmabar(d_{i,j}),	\\
\{s_{i,j}\}\cup\valustar_\sigma^\P(b_2),	&\text{if }\sigmabar(eg_{i,j})\wedge\neg\sigmabar(eb_{i,j})\wedge\relbit{\sigma}=1,	\\
\valustar_\sigma^\P(g_1),		&\text{if }\sigmabar(eg_{i,j})\wedge\neg\sigmabar(eb_{i,j})\wedge\relbit{\sigma}\neq 1,	\\
\valustar_\sigma^\P(b_2),		&\text{if }\sigmabar(eb_{i,j})\wedge\neg\sigmabar(eg_{i,j})\wedge\relbit{\sigma}=1\\
															&\phantom{if }\wedge (\neg\sigmabar(s_{i,j})\vee\sigmabar(b_{i+1})=j),	\\
\valustar_\sigma^\P(s_{i,j}),	&\text{if }\sigmabar(eb_{i,j})\wedge\neg\sigmabar(eg_{i,j})\\
															&\phantom{if }\wedge(\relbit{\sigma}\neq 1\vee (\sigmabar(s_{i,j})\wedge \sigmabar(b_{i+1})\neq j)) ,\\
\valustar_\sigma^\P(g_1),		&\text{if }\sigmabar(eb_{i,j})\wedge\sigmabar(eg_{i,j})\wedge\sigmabar(g_1)\neq\sigmabar(b_2),	\\
\valustar_\sigma^\P(b_2),		&\text{if }\sigmabar(eb_{i,j})\wedge\sigmabar(eg_{i,j})\wedge\sigmabar(g_1)=\sigmabar(b_2).
\end{cases}\]
\end{restatable}

This exact characterization of the valuations of the cycle centers can be used to determine the exact valuations of all selector vertices.
We begin by considering the case $G_n=M_n$.

\begin{corollary} \label{corollary: Complete Valuation Of Selection Vertices MDP}
Let $G_n=M_n$.
Let $\sigma\in\reach{\sigma_0}$ be well-behaved, $i\in[n]$  and define \[\lambda_i^\M\coloneqq\min\{\ell\geq i\colon\sigma(b_\ell)=g_\ell\vee\sigma(g_\ell)=F_{\ell,0}\vee\sigma(s_{\ell,\sigmabar(g_\ell)})=b_1\vee\neg\sigmabar(d_{\ell})\}.\]
Then $\valustar_{\sigma}^\M(g_i)=\sum_{\ell=i}^{\lambda-1}W_{\ell}^\M+\valustar_{\sigma}^\M(g_{\lambda})$ where $\lambda\coloneqq\lambda_i^\M$ and \[
\valustar_{\sigma}^\M(g_{\lambda})=\begin{cases} 
\valustar_{\sigma}^\M(b_{\lambda}),																&\text{if }\sigmabar(b_{\lambda}),\\
\rew{g_{\lambda}}+\frac{1}{2}\valustar_{\sigma}^\M(g_1)+\frac{1}{2}\valustar_{\sigma}^\M(b_2),	&\text{if }\neg\sigmabar(b_{\lambda})\wedge\sigmabar(eg_{\lambda})\wedge\sigmabar(eb_{\lambda}),\\
\rew{g_{\lambda}}+\valustar_{\sigma}^\M(g_1),													&\text{if }\neg\sigmabar(b_{\lambda})\wedge\sigmabar(eg_{\lambda})\wedge\neg\sigmabar(eb_{\lambda}),\\
\rew{g_{\lambda}}+\valustar_{\sigma}^\M(b_2),													&\text{if }\neg\sigmabar(b_{\lambda})\wedge\neg\sigmabar(eg_{\lambda})\wedge\sigmabar(eb_{\lambda}),\\
\rew{g_{\lambda},s_{\lambda, \sigmabar(g_{\lambda})}}+\valustar_{\sigma}^\M(b_1),	&\text{if }\neg\sigmabar(b_{\lambda})\wedge\sigmabar(d_{\lambda})\wedge\neg\sigmabar(s_{\lambda}),\\
W_{\lambda}^\M+\valustar_{\sigma}^\M(b_{\lambda+2}),											&\text{otherwise}.
\end{cases}\]
\end{corollary}

\begin{proof}
To simplify notation let $\lambda\coloneqq\lambda_i^\M$.
By the definition of $\lambda$, for all $\ell\in\{i,\dots,\lambda-1\}$, it holds that $\sigma(g_{\ell})=F_{\ell,1}, \sigmabar(d_{\ell,1})$ and $\sigma(s_{\ell,1})=h_{\ell,1}$.
This implies the first part of the statement as this yields $\valustar_{\sigma}^\M(g_{\ell})=W_{\ell}^{\M}+\valustar_{\sigma}^\M(g_{\ell+1})$ for each such index.

Thus consider $\valustar_{\sigma}^\M(g_{\lambda})$.
The first four cases follow immediately resp. by \Cref{lemma: Exact Behavior Of Random Vertex}. 
Consider the case $\neg\sigmabar(b_{\lambda})\wedge\sigmabar(d_{\lambda})\wedge\neg\sigmabar(s_{\lambda})$.
Then $\sigma(s_{\lambda,\sigmabar(g_{\lambda})})=b_1$ and the statement follows from $\valustar_{\sigma}^\M(F_{\lambda,\sigmabar(g_{\lambda})})=\valustar_{\sigma}^\M(s_{\lambda,\sigmabar(g_{\lambda})})$.
Hence consider the ``otherwise'' case, implying $\neg\sigmabar(b_{\lambda_i})\wedge\sigmabar(d_{\lambda})\wedge\sigmabar(s_{\lambda})$.
By the definition of $\lambda$, this yields  $\sigma(g_{\lambda_i})=F_{\lambda,0}$.
\end{proof}

We now prove the corresponding statement for the case $G_n=S_n$.

\begin{corollary} \label{corollary: Complete Valuation Of Selection Vertices PG}
Let $G_n=S_n$.
Let $\sigma\in\reach{\sigma_0}$ be well-behaved, $i\in[n]$ and define \[\lambda_i^\P\coloneqq\min\{\ell \geq i \colon \sigma(b_\ell)=g_\ell\vee\sigma(g_\ell)=F_{\ell,0}\vee \sigma(s_{\ell,\sigmabar(g_\ell)})=b_1\vee \sigma(b_{\ell+1})=g_{\ell+1}\}.\]
Then $\valustar_\sigma^\P(g_i)=\bigcup_{i'=i}^{\lambda_i-1}W_{i'}^\P\cup\valustar_\sigma^\P(g_{\lambda_i})$, where $\lambda\coloneqq\lambda_i^\P$ and
{
\[\valustar_\sigma^\P(g_{\lambda}) =\begin{cases}    	
	\valustar_\sigma^\P(b_{\lambda}) 																&\text{if } \sigma(b_{\lambda})=g_{\lambda}, \\    	
    \{g_{\lambda}\} \cup \valustar_\sigma^\P(g_1) 											&\text{if } \nsigmabar(b_{\lambda})\wedge\sigmabar(eg_{\lambda}) \wedge \neg\sigmabar(eb_{\lambda}) \wedge \relbit{\sigma}\neq 1,\\
    \{g_{\lambda}\} \cup \valustar_\sigma^\P(b_2) 											&\text{if } \nsigmabar(b_{\lambda})\wedge\sigmabar(eb_{\lambda}) \wedge \neg\sigmabar(eg_{\lambda}) \wedge \relbit{\sigma}=1\\
    																																					&\phantom{\text{if } }\wedge(\neg \sigmabar(s_{\lambda})\vee\sigmabar(b_{\lambda+1})=\sigmabar(g_{\lambda})) ,\\
    \{g_{\lambda}\} \cup \valustar_\sigma^\P(g_1) 											&\text{if } \nsigmabar(b_{\lambda})\wedge\sigmabar(eb_{\lambda}) \wedge \sigmabar(eg_{\lambda}) \wedge \sigmabar(g_1) \neq \sigmabar(b_2), \\
    \{g_{\lambda}\} \cup \valustar_\sigma^\P(b_2) 											&\text{if } \nsigmabar(b_{\lambda})\wedge\sigmabar(eb_{\lambda}) \wedge  \sigmabar(eg_{\lambda}) \wedge \sigmabar(g_1) = \sigmabar(b_2), \\    	
	\{g_{\lambda},s_{\lambda,\sigmabar(g_{\lambda})}\}\cup\valustar_\sigma^\P(b_1)	&\text{if none of the above and } \sigma(s_{\lambda,\sigmabar(g_{\lambda})})=b_1,\\	
	W_{\lambda}^\P\cup\valustar_\sigma^\P(b_{\lambda+2})										&\text{if none of the above and } \sigma(g_{\lambda})=F_{\lambda,0},\\
	W_{\lambda}^\P\cup\valustar_\sigma^\P(b_{\lambda+1})										&\text{otherwise}.		
\end{cases}\]}
\end{corollary}

\begin{proof}
	Let $\lambda\coloneqq\lambda_i^\P$ and $\ell\in\{i,\dots,\lambda-1\}$.
	We prove $\valustar_{\sigma}^\P(F_{\ell,\sigmabar(g_{\ell})})=\valustar_{\sigma}^\P(s_{\ell,\sigmabar(g_{\ell})})$.
	Since $\ell<\lambda$, it follows that $\sigma(b_{\ell})=b_{\ell+1}, \sigma(g_{\ell})=F_{\ell,1}$, $\sigma(s_{\ell,1})=h_{\ell,1}$ and $\sigma(b_{\ell+1})=b_{\ell+2}$.
	We show that this implies that none of the cases 3,4,6 and 7 of \Cref{lemma: Exact Behavior Of Counterstrategy} can occur.
	
	If the conditions of the third case were true, then $\sigma(s_{\ell,1})=b_1$ by \Pref{EG1}, contradicting $\sigma(s_{\ell,1})=h_{\ell,1}$.	
	If the conditions of the fourth case were true, then $\neg\sigmabar(s_{\ell,1})\vee\sigmabar(b_{\ell+1})=1$.
	But, since $\sigma(s_{\ell,1})=h_{\ell,1}$ and $\sigma(b_{\ell+1})=b_{\ell+2}$, this cannot hold.
	If the conditions of the sixth or seventh case were true, then $\sigmabar(eb_{\ell,1})\wedge\sigmabar(eg_{\ell,1})$.
	But then, since $\sigma(s_{\ell,1})=h_{\ell,1}$, \Pref{EBG1} implies $\sigma(b_{\ell+1})=g_{\ell+1}$, contradicting $\sigma(b_{\ell+1})=b_{\ell+2}$.
	
	Hence $\valustar_{\sigma}^\P(F_{\ell,\sigmabar(g_{\ell})})=\valustar_{\sigma}^\P(s_{\ell,\sigmabar(g_{\ell})})$.
	Since also $\sigmabar(g_{\ell})=1$ and $\sigma(s_{\ell,1})=h_{\ell,1}$, this implies the first part of the statement.
	It thus remains to investigate $\valustar_{\sigma}^\P(g_{\lambda})$.

	The first five statements follow directly from \Cref{lemma: Exact Behavior Of Counterstrategy}.
	Thus consider the sixth.
	Since none of the five previous cases must hold, one of the following holds:
	\begin{enumerate}
		\item $\sigma(b_{\lambda})=b_{\lambda+1}\wedge\neg\sigmabar(eb_{\lambda})\wedge[\neg\sigmabar(eg_{\lambda})\vee\relbit{\sigma}=1]$
		\item $\sigma(b_{\lambda})=b_{\lambda+1}\wedge\neg\sigmabar(eg_{\lambda})\wedge[\neg\sigmabar(eb_{\lambda})\vee\relbit{\sigma}\neq 1\vee(\sigmabar(s_{\lambda})\wedge\sigmabar(b_{\lambda+1})\neq\sigmabar(g_{\lambda}))]$.
	\end{enumerate}
	We now consider these two cases together with the assumption $\sigma(s_{\lambda, \sigmabar(g_{\lambda})})=b_1$.
	It again suffices to show \smash{$\valustar_{\sigma}^\P(F_{\lambda, \sigmabar(g_{\lambda})})=\valustar_{\sigma}^\P(s_{\lambda, \sigmabar(g_{\lambda})})$.}
	\begin{enumerate}
		\item If $\neg\sigmabar(eb_{\lambda})\wedge\neg\sigmabar(eg_{\lambda})$, then $\sigmabar(d_{\lambda})$.
			Consequently, by \Cref{lemma: Exact Behavior Of Counterstrategy}, the statement follows.
			Otherwise we have $\neg\sigmabar(eb_{\lambda})\wedge\sigmabar(eg_{\lambda})\wedge\relbit{\sigma}=1$.
			But then, the conditions of the second case of \Cref{lemma: Exact Behavior Of Counterstrategy} hold and the statement follows again.
		\item As previously, the statement follows if $\neg\sigmabar(eg_{\lambda})\wedge\neg\sigmabar(eb_{\lambda})$.
			Since $\sigma(s_{\lambda,\sigmabar(g_i)})=b_1$ by assumption, the conditions can thus only be fulfilled if $\neg\sigmabar(eg_{\lambda})\wedge\sigmabar(eb_{\lambda})\wedge\relbit{\sigma}\neq1$.
			But then, the conditions of case 5 of \Cref{lemma: Exact Behavior Of Counterstrategy} hold and the statement follows. 
	\end{enumerate}
	
	Next consider the seventh case.
	Then $\sigmabar(g_{\lambda})=0$.
	Note that $\sigma(s_{\lambda, 0})=h_{\lambda,0}$ holds by assumption and that it again suffices to show $\valustar_{\sigma}^\P(F_{\lambda,0})=\valustar_{\sigma}^\P(s_{\lambda,0})$.
	We thus again investigate the two cases mentioned before together with the assumption $\sigma(s_{\lambda,0})=h_{\lambda,0}$ and $\sigmabar(g_{\lambda})=0$.
	\begin{enumerate}
		\item Here, the same arguments used before can be applied again.
		\item If either $\neg\sigmabar(eg_{\lambda})\wedge\neg\sigmabar(eb_{\lambda})$ or $\neg\sigmabar(eg_{\lambda})\wedge\sigmabar(eb_{\lambda})\wedge\relbit{\sigma}\neq 1$, then the statement follows by the previously given arguments.
			Hence consider the case $\neg\sigmabar(eg_{\lambda})\wedge\sigmabar(eb_{\lambda})\wedge\relbit{\sigma}=1\wedge(\sigmabar(s_{\lambda})\wedge\sigmabar(b_{\lambda+1})\neq\sigmabar(g_{\lambda}))$.
			But then, the conditions of the fifth case of \Cref{lemma: Exact Behavior Of Counterstrategy} are fulfilled again, implying the statement.
	\end{enumerate}
	Finally, consider the eighth case.
	We then have $\sigma(b_{\lambda})=b_{\lambda+1}, \sigma(s_{\lambda,1})=h_{\lambda,1}$ and $\sigmabar(g_{\lambda})=1$.
	By the definition of $\lambda$, we thus need to have $\sigma(b_{\lambda+1})=g_{\lambda+1}$.
	It hence suffices to prove that $\valustar_{\sigma}^\P(F_{\lambda,1})=\valustar_{\sigma}^\P(s_{\lambda,1})$.
	This however follows by the same arguments used in the last case.
\end{proof}

To conclude the characterization of the vertex valuations, we state one additional lemma.
This lemma allows us to simplify the evaluation of the valuation of the cycle centers under certain conditions without having to check the conditions of \Cref{lemma: Exact Behavior Of Random Vertex} resp. \ref{lemma: Exact Behavior Of Counterstrategy}.
It will in particular be used when analyzing cycle centers during phase~$1$.

\begin{lemma} \label{lemma: Cycle centers in Phase One}
Let $\sigma\in\reach{\sigma_0}$ be well-behaved, let $i\in[n], j\in\{0,1\}$ and consider the cycle center $F_{i,j}$.
Assume that $\sigma$ has Properties (\ref{property: ESC1}),  (\ref{property: EV1})$_1$ and (\ref{property: USV1})$_{i}$.
It then holds that $\valustar_{\sigma}^*(F_{i,j})=\valustar_{\sigma}^*(s_{i,j})$ if $\sigmabar(d_{i,j})$ and $\valustar_{\sigma}^*(F_{i,j})=\valustar_{\sigma}^*(b_1)$ otherwise.
\end{lemma}

\begin{proof}
If $\sigmabar(d_{i,j})$, then the statement follows from \Cref{lemma: Exact Behavior Of Random Vertex} resp. \ref{lemma: Exact Behavior Of Counterstrategy}.
Hence consider the case that~$F_{i,j}$ is not closed and let $G_n=S_n$.
We show that the conditions of either the first or the fourth case of \Cref{lemma: Exact Behavior Of Counterstrategy} are fulfilled and that the corresponding valuations can be expressed as $\valustar_{\sigma}^\P(b_1)$.

By \Pref{ESC1}, the last two cases of \Cref{lemma: Exact Behavior Of Counterstrategy} cannot occur.
Let, for the sake of contradiction, the conditions of the second case be fulfilled, i.e., $\sigmabar(eg_{i,j}),\nsigmabar(eb_{i,j})$ and $\relbit{\sigma}=1$.
Then $\sigma(b_1)=b_2$.
By \Pref{EV1}$_1$ and the definition of the induced bit state, this implies $\indbit_1=0$.
Hence, by \Pref{ESC1}, $\sigma(e_{*,*,*})=b_2$.
Since $F_{i,j}$ is not closed, this implies that there is at least one $k\in\{0,1\}$ such that $\sigma(d_{i,j,k})=e_{i,j,k}$ and $\sigma(e_{i,j,k})=b_2$.
But then $\sigmabar(eb_{i,j})$ contradicting $\nsigmabar(eb_{i,j})$.

Now let, for the sake of contradiction, the conditions of the fifth case of \Cref{lemma: Exact Behavior Of Counterstrategy} be fulfilled, i.e., $\sigmabar(eb_{i,j}),\nsigmabar(eg_{i,j})$ and either $\relbit{\sigma}\neq 1$ or $\sigmabar(s_{i,j})\wedge\sigmabar(b_{i+1})\neq j$.
If $\relbit{\sigma}\neq 1$, we can deduce $\sigmabar(eg_{i,j})$ by the same arguments used for the second case, again resulting in a contradiction.
Thus let $\sigma(s_{i,j})=h_{i,j}$ and $\sigmabar(b_{i+1})\neq j$.
Then, by \Pref{USV1}$_i$, $j=\indbit_{i+1}$.
But then, the other condition states $\sigmabar(b_{i+1})\neq\indbit_{i+1}$ which is a contradiction since $\indbit_{i+1}=\sigmabar(d_{i,\bit_{i+1}})$ by definition.

Consider the third case of \Cref{lemma: Exact Behavior Of Counterstrategy}.
Then, $\relbit{\sigma}\neq 1$ implies $\sigma(b_1)=g_1$ by \Cref{lemma: b1 iff relbit}.
Thus, $\valustar_{\sigma}^\P(b_1)=\valustar_{\sigma}^\P(g_1)=\valustar_{\sigma}^\P(F_{i,j})$.
Consider the fourth case of \Cref{lemma: Exact Behavior Of Counterstrategy}.
Then $\relbit{\sigma}=1$,  implying $\valustar_{\sigma}^\P(b_1)=\valustar_{\sigma}^\P(b_2)=\valustar_{\sigma}^\P(F_{i,j})$ as $\sigma(b_1)=b_2$.

Now let $G_n=M_n$ and let $F_{i,j}$ not be closed.
By \Pref{ESC1}, $F_{i,j}$ cannot escape towards both $g_1$ and $b_2$.
By \Pref{EV1}$_1$, $\indbit_1=1$ if and only if$\sigma(b_1)=g_1$.
The statement thus follows since \Pref{ESC1} implies that $F_{i,j}$ escapes to $g_1$ if $\indbit_1=1$ and to $b_2$ if $\indbit_1=0$.
\end{proof}

This concludes our general results on the vertex valuations.
As it will turn out that every strategy calculated by the strategy improvement resp. policy iteration algorithm is well-behaved.
Consequently, these characterization can be applied to all strategies that are considered in the following sections.

\section{The Application of Individual Improving Switches} \label{section: Improving switches technical}

This section contains technical details related to the application of individual improving switches and the different phases of a single transition.
We consider a fixed number $\bit\in\bitset_n$.
Before analyzing the single phases, we develop some general statements that are either not related to a single phase or are used repeatedly in the upcoming proofs.
Henceforth, $\bit\in\bitset_n$ is a fixed number and $\nsb\coloneqq\ell(\bit+1)$ denotes the least significant set bit of $\bit+1$.
We further define the abbreviation $\sum(\bit,i)\coloneqq\sum_{l<i}\bit_l\cdot 2^{l-1}$.

Most of the proofs of this section are deferred to \Cref{appendix: Proofs Exponential}.

\subsection{Basic statements and statements independent of the phases} 

The first lemma enables us to compare the valuations of cycle centers in~$M_n$ for several well-behaved strategies calculated by the strategy improvement algorithm.

\begin{restatable}{lemma}{BothCCOpenForMDP} \label{lemma: Both CC Open For MDP}
Let $G_n=M_n$.
Let $\sigma\in\reach{\sigma_0}$ be a well-behaved phase-$k$-strategy for some $\bit\in\bitset_n$ having \Pref{USV1}$_i$ and \Pref{EV1}$_{i+1}$ for some $i\in[n]$ where $k\in[5]$.
If $F_{i,0}$ and $F_{i,1}$ are in the same state and if either $i\geq\nsb$ or~$\sigma$ has \Pref{REL1}, then $\valu_{\sigma}^\M(F_{i,\indbit_{i+1}})>\valu_{\sigma}^\M(F_{i,1-\indbit_{i+1}})$.
\end{restatable}

As mentioned at the beginning of \Cref{section: Vertex Valuations}, formal lemmas describing the applications of the improving switches are proven only for well-behaved strategies.
We will always prove that the strategies obtained by the application of improving switches remain well-behaved and prove that canonical strategies are well-behaved.
Consequently, \emph{all} strategies calculated by the strategy improvement are well-behaved.

As a basis for these arguments, we prove that canonical strategies are well-behaved.

\begin{lemma} \label{lemma: Canonical strategies are well-behaved}
Let $\canstrat$ be a canonical strategy for some $\bit\in\bitset_{n}$.
Then $\canstrat$ is well-behaved.
\end{lemma}

\begin{proof}
Let $\sigma\coloneqq\canstrat$ and let $i\in[n]$ such that $\sigma(b_{i})=g_{i}$.
Then, by the definition of a canonical strategy, we have $\bit_{i}=1$, implying $\sigma(g_{i})=F_{i,\bit_{i+1}}$.
Hence, $\sigmabar(g_{i})=\bit_{i+1}=\sigmabar(b_{i+1})$.
Thus $\incorrect{\sigma}=\emptyset$, implying 
\begin{equation} \label{equation: Relbit for canonical strategy}
\relbit{\sigma}=\min\{i\in[n+1]\colon\sigma(b_{i})=b_{i+1}\}.
\end{equation}
We prove that $\sigma$ has all properties of \Cref{table: Well behaved properties}.
We investigate each property and show that either its premise is false or that both the premise and the conclusion are true.
\begin{itemize}[leftmargin=1.75cm]
	\item[(\ref{property: S1})] Let $i\geq\relbit{\sigma}$ with $\sigma(b_i)=g_i$.
		Then, $\sigmabar(g_i)=\bit_{i+1}$ and $\sigmabar(s_{i,\bit_{i+1}})=1$, hence $\sigmabar(s_{i}).$
	\item[(\ref{property: S2})] Let $i<\relbit{\sigma}$.
		Then, by (\ref{equation: Relbit for canonical strategy}), $\sigma(b_i)=g_i$.
		By \Pref{S1}, this implies $\sigmabar(s_i)$.
	\item[(\ref{property: B1})] Let $i<\relbit{\sigma}-1$.
		Then $\sigma(b_i)=g_i$ by (\ref{equation: Relbit for canonical strategy}), hence the premise is false.
	\item[(\ref{property: B2})] By (\ref{equation: Relbit for canonical strategy}), we have $\sigma(b_{\relbit{\sigma}-1})=g_{\relbit{\sigma}-1}$, so the premise is false.
	\item[(\ref{property: B3})] Let $i\in[n]$ with $\sigma(b_{i+1})=b_{i+2}$.
		This implies $\bit_{i+1}=0$, so $\sigma(s_{i,1})=b_1$.
	\item[(\ref{property: BR1})] Let $i<\relbit{\sigma}$.
		This implies $\sigma(b_i)=g_i$ and $\sigma(g_i)=F_{i,\bit_{i+1}}$.
		Assume $\bit_{i+1}=0$.
		We then have $\relbit{\sigma}\leq i+1$ by (\ref{equation: Relbit for canonical strategy}), implying $\relbit{\sigma}=i+1$. 
		But then $\sigma(g_i)=F_{i,0}$ and $i=\relbit{\sigma}-1$.
		Now assume $\bit_{i+1}=1$.
		We then have $\relbit{\sigma}>i$ and $\relbit{\sigma}\neq i+1$, so $\relbit{\sigma}>i+1$.
		But then, $\sigma(g_i)=F_{i,1}$ and $i<\relbit{\sigma}-1$ and in particular $i\neq\relbit{\sigma}-1$.
	\item[(\ref{property: BR2})] Since $\incorrect{\sigma}=\emptyset$, $i<\relbit{\sigma}$ implies $\sigmabar(d_i)$ by \Cref{corollary: Simplified MDP Valuation}, so $\nsigmabar(eg_{i,\sigmabar(eg_{i})})$.
	\item[(\ref{property: D1})] Since $\sigma(b_i)=g_i$ implies $\sigmabar(g_i)=\bit_{i+1}$ and that $F_{i,\bit_{i+1}}$ is closed, both premise and conclusion are true.
	\item[(\ref{property: D2})] This follows by the same arguments used in the last case since $\sigma(b_i)=g_i$ for any $i<\relbit{\sigma}$ by (\ref{equation: Relbit for canonical strategy}).
	\item[(\ref{property: MNS1})] By \Cref{lemma: Config implied by Aeb}, the premise implies $\minsig{b}=2$, hence $\sigma(b_2)=g_2$ and $\sigma(b_1)=b_2$.
		Consequently, $\sigma(g_1)=F_{i,1}$ as $\minsig{b}\leq\minnegsig{g},\minnegsig{s}$, contradicting the definition of a canonical strategy if $G_n=S_n$.
	\item[(\ref{property: MNS2})] Assuming that there was some index $i<\minnegsig{g}<\minnegsig{s},\minsig{b}$ and $\relbit{\sigma}=1$ implies $\sigma(b_1)=b_2$ as well as $\sigma(g_1)=F_{1,1}$, contradicting the definition of a canonical strategy for the case $G_n=S_n$.
		If $G_n=M_n$, then we need to have $\bit_2=1$, implying $\minsig{b}=2$.
		This is however a contradiction to $1<\minnegsig{g}<\minsig{b}$.
	\item[(\ref{property: MNS3})] Let $\relbit{\sigma}=1$ and assume there was some $i<\minnegsig{s}\leq\minnegsig{g}<\minsig{b}$.
		As $G_n=M_n$ by assumption, $\sigma(s_{1,1})=h_{1,1}$ then implies $\bit_2=1$ and thus $\minsig{b}=2$, contradicting $1<\minnegsig{s}<\minsig{b}$.
	\item[(\ref{property: MNS4})] Let $\relbit{\sigma}=1$, assume $\minnegsig{s}\leq\minnegsig{g}<\minsig{b}$ and let $i\coloneqq\minnegsig{s}$.
		We prove that the premise either yields a contradiction or implies $\sigmabar(eb_{i})\wedge\nsigmabar(eg_{i})$.
		Since $\relbit{\sigma}=1$ implies $\sigma(b_1)=b_2$ and thus $\bit_1=0$, the definition of a canonical strategy implies that it suffices to prove that $F_{i,\sigmabar(g_i)}$ is not closed.
		This follows from the definition of a canonical strategy if $i=1$, so let $i>1$.
		Then $1<i=\minnegsig{s}\leq\minnegsig{g}$, hence $\sigma(g_1)=F_{1,1}$.
		This however contradicts the definition of a canonical strategy if $G_n=S_n$.
		If $G_n=M_n$, then this implies $\bit_2=1$ and hence $\sigma(b_2)=g_2$, thus $\minsig{b}=2$.
		But this contradicts the premise $\minnegsig{g}<\minsig{b}$.
	\item[(\ref{property: MNS5})] Let $\relbit{\sigma}=1$.
		We show that there is no $i<\minnegsig{s}<\minsig{b}\leq\minnegsig{g}$.
		Assume there was such an index $i$, implying that $1<\minnegsig{s}<\minsig{b}\leq\minnegsig{g}$.
		Thus, $\sigma(g_1)=F_{1,1}$ and $\sigma(s_{1,1})=h_{1,1}$, implying $\bit_2=1$.
		But then $\minsig{b}=2$, contradicting $1<\minnegsig{s}<\minsig{b}$.
	\item[(\ref{property: MNS6})] Let $\relbit{\sigma}=1$.
		If $\minnegsig{s}\neq 1$, the same arguments used when discussing \Pref{MNS5} can be applied.
		However, for $\minnegsig{s}=1$, the statement follows since both cycle centers of level 1 are open and since these cycle centers escape to $b_2$.
	\item[(\ref{property: EG1})] By $\relbit{\sigma}=1$, we have $\sigma(b_1)=b_2$, implying $\bit_1=0$.
		Thus, any cycle center which is not closed escapes towards $b_2$ by definition, hence the premise is incorrect.
	\item[(\ref{property: EG2})] Follows by the same arguments used in the last case.
	\item[(\ref{property: EG3})] Assume that there is some cycle center escaping towards $g_1$.
		Then $\bit_1=1$.
		This implies $\sigma(b_1)=g_1$ and by the same arguments used earlier in this proof, this implies $\sigmabar(s_1)$.
	\item[(\ref{property: EG4})] This follows by the same arguments used when discussing Properties (\ref{property: EG1}) and (\ref{property: EG2}).	
	\item[(\ref{property: EG5})] It is easy to see that $\sigmabar(s_{i,j})$ implies $\sigmabar(b_{i+1})=j$.
	\item[(EB*)] Every premise of any of the properties (EB*) contains $\sigma(b_1)=g_1$.
		Hence, we always have $\bit\bmod2=1$, implying that no cycle center can escape towards $b_2$.
		But this implies that the premise any of these properties is false.
	\item[(EBG*)] By the definition of a canonical strategy, no cycle center can escape towards both $b_2$ and $g_1$.
	\item[(\ref{property: DN1})] By the definition of a canonical strategy, $\sigmabar(d_n)$ holds if and only if $\sigma(b_n)=g_n$. 
		Hence both the premise and the conclusion are true.
	\item[(\ref{property: DN2})] If $\sigmabar(d_n)$ the statement follows analogously as in the last case.
		Hence assume $\minnegsig{g}=n$.
		Then $\sigma(g_i)=F_{i,1}$ for all $i<n$, so in particular, $\sigma(g_1)=F_{1,1}$.
		But, by the definition of a canonical strategy, this immediately implies $\bit_1=1$ and $\sigma(b_1)=g_1$  if $G_n=S_n$ and $\bit_2$ and thus $\sigma(b_2)=g_2$ if $G_n=M_n$. \qedhere
\end{itemize}
\end{proof}

Our goal is to prove \Cref{lemma: Improving sets of canonical strategies} next as this statement describes the set of improving switches of canonical strategies.
Before doing so, we analyze the terms used in \Cref{table: Occurrence Records} to describe the occurrence records of the cycle vertices in more detail.
As the proofs of the following lemmas are rather technical, they are deferred to \Cref{appendix: Proofs Exponential}.

\begin{restatable}{lemma}{NumericsOfEll} \label{lemma: Numerics Of Ell}
Let $\bit\in\bitset_n$.
If $\id_{j=0}\lastflip{\bit}{i+1}{}+\id_{j=1}\lastunflip{\bit}{i+1}{}=0$ for  $i\in[n], j\in\{0,1\}$, then $\ell^{\bit}(i,j,k)\geq\bit$ for $,k\in\{0,1\}$.
Otherwise, the following hold:
\begin{center}
\renewcommand{\arraystretch}{1.5}
	\begin{tabular}{|c||c|c|c|} \hline
		Setting of bits				&$\bit_i=1\wedge\bit_{i+1}=1-j$						&$\bit_i=0\wedge\bit_{i+1}=j$						&$\bit_i=0\wedge\bit_{i+1}=1-j$\\\hdashline
		\multirow{2}{*}{$\ell^{\bit}(i,j,k)=$}	&\multirow{2}{*}{$\ceil{\frac{\bit+\sum(\bit,i)+1-k}{2}}$}	&\multirow{2}{*}{$\ceil{\frac{\bit+2^{i-1}+\sum(\bit,i)+1-k}{2}}$}	&\multirow{2}{*}{$\ceil{\frac{\bit-2^{i-1}+\sum(\bit,i)+1-k}{2}}$}\\&&&\\\hline
	\end{tabular}
\end{center}
\end{restatable}

\begin{restatable}{lemma}{ProgressForUnimportantCC} \label{lemma: Progress for unimportant CC}
Let $\bit\in\bitset_n$ and $i\in[n]$ and $ j\in\{0,1\}$ such that $\bit_i=0$ or $\bit_{i+1}\neq j$.
Then, \[\id_{j=0}\lastflip{\bit}{i+1}{}-\id_{j=1}\lastunflip{\bit}{i+1}{}=\id_{j=0}\lastflip{\bit+1}{i+1}{}-\id_{j=1}\lastunflip{\bit+1}{i+1}{}.\]
Moreover, if $i\neq\nsb$, then $\ell^{\bit}(i,j,k)+1=\ell^{\bit+1}(i,j,k)$.
\end{restatable}

\begin{restatable}{lemma}{OccurrenceRecordsCycleVertices} \label{lemma: Occurrence Records Cycle Vertices}
Let $\canstrat$ be a canonical strategy for $\bit$ such that its occurrence records are described by \Cref{table: Occurrence Records}.
Assume that $\canstrat$ has Properties (\ref{property: OR1})$_{*,*,*}$ to (\ref{property: OR4})$_{*,*,*}$.
Then, the following hold.
\begin{enumerate}
	\item Let $i\in[n] $ and $j\in \{0,1\}$ and assume that either $\bit_i=0$ or $\bit_{i+1}\neq j$.
		Then, it holds that $\occrec^{\canstrat}(d_{i,j,*},F_{i,j})\leq\floor{(\bit+1)/2}$.
	\item Let $j\coloneqq\bit_{\nsb+1}$.
		Then, $\occrec^{\canstrat}(d_{\nsb,j,0},F_{\nsb,j})=\floor{(\bit+1)/2}$.
		In addition, $\nsb=1$ implies $\occrec^{\canstrat}(d_{\nsb,j,1},F_{\nsb,j})=\floor{(\bit+1)/2}$ and $\nsb>1$ implies $\occrec^{\sigma}(d_{\nsb,j,1},F_{\nu,j})=\floor{(\bit+1)/2}-1$.
	\item If $i=1$, then $\canstrat(d_{1,1-\bit_{2},*})\neq F_{1,1-\bit_{2}}$ and $\occrec^{\canstrat}(d_{1,1-\bit_{2},0},F_{1,1-\bit_{2}})=\floor{(\bit+1)/2}$.
\end{enumerate}
\end{restatable}

\begin{restatable}{lemma}{NumericsOfOR} \label{lemma: Numerics Of OR}
Let $\bit\in\bitset_n$ and $i\in[n]$.
It holds that $\flips{\bit}{i}{}=\floor{(\bit+2^{i-1})/2^{i}}$ and $\flips{\bit+1}{i}{}=\flips{\bit}{i}{}+\id_{i=\nsb}$.
In addition, for indices $i_1,i_2\in[n]$ with $i_1<i_2$ and $\bit\geq2^{i_1-1}$ imply $\flips{\bit}{i_1}{}>\flips{\bit}{i_2}{}$.
Furthermore, if $k\coloneqq \frac{\bit+1}{2^{\nsb-1}}$ and $x\in[\nsb-1]$, then $\flips{\bit}{\nsb-x}{}=k\cdot 2^{x-1}$.
\end{restatable}

Now all general statements required for the upcoming proofs and statements are in place.
We begin by analyzing phase $1$, or, more precisely, the statements that prove the application of improving switches until reaching phase $2$.

We first restate \Cref{lemma: Improving sets of canonical strategies} and provide its formal proof.

\ImprovingSetsOfCanonicalStrategies*

\begin{proof}
It is easy to verify that canonical strategies are phase-$1$-strategies.
To simplify notation, let $\sigma\coloneqq\canstrat$ and $\mathfrak{D}^{\sigma}\coloneqq\{(d_{i,j,k},F_{i,j}):\sigma(d_{i,j,k})\neq F_{i,j}\}$.
It then suffices to show $I_{\sigma}=\mathfrak{D}^{\sigma}$.
We thus have to prove that $\sigma(d_{i,j,k})\neq F_{i,j}$ implies $\valu_{\sigma}^*(F_{i,j})>\valu_{\sigma}^*(e_{i,j,k})$ and that there are no other improving switches.

Let $e=(d_{i,j,k},F_{i,j})$ with $\sigma(d_{i,j,k})\neq F_{i,j}$.
By \Cref{lemma: Canonical strategies are well-behaved}, $\sigma$ is well-behaved, and the results of \Cref{section: Vertex Valuations} can be applied.
By \Cref{lemma: Cycle centers in Phase One}, $\valustar_{\sigma}^*(F_{i,j})=\valustar_{\sigma}^*(b_1)$.
Consider the case $G_n=S_n$.
Then, $\valu_{\sigma}^\P(F_{i,j})=\{F_{i,j},d_{i,j,k'},e_{i,j,k'}\}\cup\valu_{\sigma}^\P(\sigma(e_{i,j,k'}))$ for some $k'\in\{0,1\}$, implying $\valu_{\sigma}^\P(e_{i,j,*})\lhd\valu_{\sigma}^\P(F_{i,j})$ since $\sigma(e_{i,j,0})=\sigma(e_{i,j,1})$ by \Pref{ESC1}.
Now let $G_n=M_n$.
By \Pref{ESC1} and \Pref{EV1}$_1$, $\valu_{\sigma}^\M(e_{*,*,*})=\valu_{\sigma}^\M(b_1)$ for all escape vertices $e_{*,*,*}$.
Thus, $\valu_{\sigma}^\M(F_{i,j})=(1-\e)\valu_{\sigma}(b_1)+\e\valu_{\sigma}^\M(s_{i,j})$ and it suffices to prove $\valu_{\sigma}^\M(s_{i,j})>\valu_{\sigma}^\M(b_1)$.
If $\sigma(s_{i,j})=b_1$, then this follows immediately from $\valu_{\sigma}^\M(s_{i,j})=\rew{s_{i,j}}+\valu_{\sigma}^\M(b_1)>\valu_{\sigma}^\M(b_1)$.
Thus assume $\sigma(s_{i,j})=h_{i,j}$.
Then, by \Pref{USV1}$_{i}$, $j=\indbit_{i+1}$ and $\valustar_{\sigma}^\M(s_{i,j})=\rew{s_{i,j},h_{i,j}}+\valustar_{\sigma}^\M(b_{i+1})$ by \Pref{EV1}$_{i+1}$.
Hence $F_{i,j}$ is the active cycle center of level $i$.
Since it is not closed by assumption, we thus have $\bit_i=0$ by \Pref{EV1}$_i$ and $i\geq\relbit{\sigma}$ by \Pref{REL1}, implying $\valustar_{\sigma}^\M(b_{i+1})=L_{i+1}^\M$.
Furthermore, by \Cref{lemma: Valuation of b} and \Cref{corollary: Simplified MDP Valuation}, $\valustar_{\sigma}^\M(b_1)=B_1^\M$.
If $B_1^\M=L_1^\M$, then \Cref{lemma: VV Lemma}~(4.) implies $L_1^\M<\rew{s_{i,j},h_{i,j}}+L_{i+1}^\M$. 
If $B_1^\M=R_1^\M$, then $i\geq\relbit{\sigma}$ and \Cref{lemma: VV Lemma}~(3.) yields $R_1^\M<\rew{s_{i,j},h_{i,j}}+L_{i+1}^\M$. 
Hence $\valustar_{\sigma}^\M(s_{i,j})>\valustar_{\sigma}^\M(b_1)$, implying $\valustar_{\sigma}^\M(F_{i,j})>\valustar_{\sigma}^\M(e_{i,j,*})$.
Consequently, $\sigma(d_{i,j,k})\neq F_{i,j}$ implies $\valu_{\sigma}^*(F_{i,j})>\valu_{\sigma}^*(\sigma(d_{i,j,k}))$ in both $S_n$ and $M_n$, so $(d_{i,j,k},F_{i,j})\in I_{\sigma}$.
It remains to show that there are no other improving switches.

We first show that there is no improving switch $e=(b_i,*)$.
Let $i\in[n]$ and $\sigma(b_i)=g_i$.
We need to show $\valu_{\sigma}^*(b_{i+1})\preceq\valu_{\sigma}^*(g_i)$.
Since $\sigma(b_i)=g_i$ implies $\valustar_{\sigma}^*(b_i)=\valustar_{\sigma}^*(g_i)$, it suffices to show $\valustar_{\sigma}^*(b_{i+1})\prec\valustar_{\sigma}^*(b_i)$.
By \Cref{lemma: Valuation of b} and \Cref{corollary: Simplified MDP Valuation}, we have $\valustar_{\sigma}^*(b_{i+1})=B_{i+1}^*$ and $\valustar_{\sigma}^*(b_i)=B_i^*$.
Assume $i<\relbit{\sigma}$.
Then $B_i^*=R_i^*$ and the statement follows directly if $B_{i+1}^*=R_{i+1}^*$.
If $B_{i+1}^*=L_{i+1}^*$, we need to have $i+1=\relbit{\sigma}$.
But then $\sigma(b_{i+1})=b_{i+2}$ and thus $L_{i+1}^*\prec R_i^*=W_i^*\oplus L_{i+1}^*$.
Thus assume $i\geq\relbit{\sigma}$.
Then $B_{i+1}^*=L_{i+1}^*$ and $B_i^*=L_i^*$ and the statement follows by $\sigma(b_i)=g_i$.

Now let $\sigma(b_i)=b_{i+1}$.
We prove $\valustar_{\sigma}^*(g_i)\prec\valustar_{\sigma}^*(b_{i+1})$.
Note that $\bit_i=0$ implies $\relbit{\sigma}\leq i$, so $\valustar_{\sigma}^*(b_{i+1})=L_{i+1}^*=L_i^*$.
We use \Cref{corollary: Complete Valuation Of Selection Vertices MDP} resp. \Cref{corollary: Complete Valuation Of Selection Vertices PG}  to compute the valuation of~$g_i$.
We thus need to evaluate $\lambda_i^\M$ resp.~$\lambda_i^\P$.
If $\sigma(g_i)=F_{i,0}$, we have $\lambda_i^*=i$. 
If $\sigma(g_i)=F_{i,1}\wedge\sigma(b_{i+1})=b_{i+2}$, we have $\sigmabar(g_i)=1\neq0=\bit_{i+1}$.
Thus, by the definition of a canonical strategy, $\sigma(s_{i,\sigmabar(g_i)})=b_1$, implying $\lambda_i^*=i$.
Hence assume $\sigma(g_i)=F_{i,1}\wedge\sigma(b_{i+1})=g_{i+1}$.
If $G_n=S_n$, then this implies $\lambda_i^\P=i$ by the definition of $\lambda_i^\P$.
If $G_n=M_n$, then, $\lambda_i^\M=i$ follows since we need to have $\ncanstratbar(d_i)$ due to $\sigma(b_i)=b_{i+1},\sigma(g_i)=F_{i,1}$ and $\sigma(b_{i+1})=g_{i+1}$.
Hence, $\lambda_i^*=i$ in both cases.

Let $G_n=S_n$ and consider the different cases listed in \Cref{corollary: Complete Valuation Of Selection Vertices PG} describing the vertex valuations for selection vertices in $S_n$.
In order to show the statement we distinguish the cases listed in that corollary.
Note that the first case cannot occur.
\begin{itemize}
	\item Let $\sigmabar(eg_{i}), \nsigmabar(eb_i)$ and $\relbit{\sigma}\neq 1$.
		Then, $\valustar_{\sigma}^\P(g_i)=\{g_i\}\cup\valustar_{\sigma}^\P(g_1)$ by \Cref{corollary: Complete Valuation Of Selection Vertices PG}.
		Since $\relbit{\sigma}\neq 1$ implies $1<\relbit{\sigma}$,  $\valustar_{\sigma}^\P(g_1)=R_1^\P$ by \Cref{lemma: Valuation of g if level small}.
		Hence, since $i\geq\relbit{\sigma}$, \begin{align*}
		\valustar_{\sigma}^\P(g_i)&=\{g_i\}\cup R_1^\P\lhd\bigcup_{i'\geq i}\{W_{i'}^\P\colon\sigma(b_{i'})=g_{i'}\}=L_{i}^\P.
		\end{align*}
	\item Let $\sigmabar(eb_i),\nsigmabar(eg_i),\relbit{\sigma}=1$ and $(\nsigmabar(s_{i,j})\vee\sigmabar(b_{i+1})=\sigmabar(g_i))$.
		Then, by \Cref{corollary: Complete Valuation Of Selection Vertices PG}, $\valustar_{\sigma}^\P(g_i)=\{g_i\}\cup\valustar_{\sigma}^\P(b_2)$.
		Since $\relbit{\sigma}=1$ we have $\valustar_{\sigma}^\P(b_2)=L_2^\P$.
		Thus, since $\sigma(b_i)=b_{i+1}$, we obtain $\valustar_{\sigma}^\P(g_i)=\{g_{i}\}\cup\valustar_{\sigma}^\P(b_2)=\{g_i\}\cup L_2^\P\lhd L_{i+1}^\P.$
\end{itemize}

\vspace*{-.75em}

	This covers the first three cases.	The fourth and fifth case cannot occur since they require a cycle center to escape towards both $g_1$ and $b_2$. 
	
\begin{itemize}	
	\item Let the conditions of case six be fulfilled.
		Then $\valustar_{\sigma}^\P(g_i)=\{g_i,s_{i,j}\}\cup\valustar_{\sigma}^\P(b_1)$ where $j=\sigmabar(g_i)$.
		If $\valustar_{\sigma}^\P(b_1)=L_1^\P$, then $\valustar_{\sigma}^\P(g_i)=\{g_i,s_{i,j}\}\cup L_1^\P\lhd L_{i+1}^\P.$
		If $\valustar_{\sigma}^\P(b_1)=R_1^\P$, then the statement follows by the same calculations used in the first case.
	\item Let the conditions of case seven be fulfilled.
		Then $\sigma(g_i)=F_{i,0}$.
		It is easy to verify that we then have $\sigma(b_{i})=b_{i+1}$ and either $\nsigmabar(eb_{i})\wedge[\nsigmabar(eg_{i})\vee\relbit{\sigma}=1]$  or $\nsigmabar(eg_{i})\wedge[\nsigmabar(eb_{i})\vee\relbit{\sigma}\neq 1\vee\sigmabar(b_{i+1})\nsigmabar(g_{i})]$.
		If $\nsigmabar(eb_{i})\wedge\nsigmabar(eg_{i})$, then $\sigmabar(d_i)$.
		But then, $\sigma(b_{i})=b_{i+1}$ implies $\sigmabar(g_i)\neq\sigmabar(b_{i+1})$.
		Hence, \Pref{USV1}$_i$ implies $\sigma(s_{i,0})=b_1$, contradicting the currently considered case.
		
		Thus consider the case $\nsigmabar(eb_i)\wedge\sigmabar(eg_i)\wedge\relbit{\sigma}=1$ next.
		Then, since $\relbit{\sigma}=1$ and since~$\sigma$ is well-behaved, $\bit_1=0$.
		But $\sigmabar(eg_i)$ implies $\bit_1=1$ which is a contradiction.
		
		Next, consider the case $\nsigmabar(eg_i)\wedge\sigmabar(eb_i)\wedge\relbit{\sigma}\neq 1$.
		As before, $\relbit{\sigma}\neq 1$ implies $\bit_1=1$ whereas $\sigma(eb_i)$ implies $\bit_1=0$, again resulting in a contradiction.
		
		Thus, consider the case $\nsigmabar(eg_i)\wedge\sigmabar(eb_i)\wedge\relbit{\sigma}=1\wedge\sigmabar(b_{i+1})\neq\sigmabar(g_i)$.
		Then, since $\sigma(g_i)=F_{i,0}$, we have $\sigma(b_{i+1})=g_{i+1}$.
		Since $\relbit{\sigma}=1$ implies $\valustar_{\sigma}^\P(b_{i+1})=L_{i+1}^\P$, we thus have	$\valustar_{\sigma}^\P(g_i)=W_i^\P\cup\valustar_{\sigma}^\P(b_{i+2})\lhd W_{i+1}^\P\cup\valustar_{\sigma}^\P(b_{i+2})=\valustar_{\sigma}^\P(b_{i+1}).$ 
	\item Case eight can only occur if $F_{i,\bit_{i+1}}$ is closed, contradicting $\bit_i=0$.
\end{itemize}

\vspace*{0.75em}

Let $G_n=M_n$ and consider \Cref{corollary: Complete Valuation Of Selection Vertices MDP}.
As before, we distinguish between the different cases listed in this corollary.
The first two cases cannot occur due to $\sigma(b_i)=b_{i+1}$ resp. \Pref{ESC1}.
Consider the third case, implying $\sigma(b_1)=g_1$ and consequently $\valustar_{\sigma}^\M(g_1)=\valustar_{\sigma}^\M(b_1)=R_1^\M$.
Using $i\geq\relbit{\sigma}$ and $\valustar_{\sigma}^\M(g_i)=\rew{g_i}+\valustar_{\sigma}^\M(g_1)$ we obtain \begin{align*}
	\valustar_{\sigma}^\M(g_i)&=\rew{g_1}+R_1^\M<\sum_{\ell\geq i}\{W_{\ell}^\M\colon\sigma(b_\ell)=g_{\ell}\}=L_i^\M.
\end{align*}
Consider the fourth case.
Then $\valustar_{\sigma}^\M(g_i)=\rew{g_i}+\valustar_{\sigma}^\M(b_2)$ and $\valustar_{\sigma}^\M(b_2)=L_2^\M$.
Thus, by $\sigma(b_i)=b_{i+1}$, we obtain $\valustar_{\sigma}^\M(g_i)=\rew{g_i}+\valustar_{\sigma}^\M(b_2)=\rew{g_i}+L_2^\M< L_{i+1}^\M.$
Consider the fifth case, implying $\valustar_{\sigma}^\M(g_i)=\rew{g_i,s_{i,\sigmabar(g_i)}}+\valustar_{\sigma}^\M(b_1)$.
Then, the statement follows analogously to the third case if $\sigma(b_1)=g_1$ and analogously to the fourth case if $\sigma(b_1)=b_2$.
The sixth case requires that the active cycle center of level $i$ is closed, contradicting $\bit_i=0$ resp. \Pref{EV1}$_i$.
Therefore there are no improving switches $e=(b_i,*)$.

Now consider some  $g_i$ with $ i\in[n-1]$ since $\sigma(g_n)=F_{n,0}$ for every $\sigma$ by construction.
First assume $\bit_i=0$.
Then, by \Cref{definition: Canonical Strategy} resp. \ref{definition: Canonical Strategy MDP}, $F_{i,\bit_{i+1}}$ is not closed.
Assume that $F_{i,1-\bit_{i+1}}$ is not closed either.
Then, by \Pref{ESC1} and \Pref{REL1}, $\relbit{\sigma}=1$ implies $\sigmabar(eb_{i,j})\wedge\nsigmabar(eg_{i,j})$ and $\relbit{\sigma}\neq 1$ implies $\sigmabar(eg_{i,j})\wedge\nsigmabar(eb_{i,j})$ for both $j\in\{0,1\}$.
Let $G_n=S_n$.
Then, by \Cref{lemma: Cycle centers in Phase One}, both cycle centers of level $i$ escape towards the same vertex via some escape vertex.
Since $\Omega(F_{i,0})=6$ and $\Omega(F_{i,1})=4$, this implies $\valu_{\sigma}^\P(F_{i,0})\rhd\valu_{\sigma}^\P(F_{i,1})$.
Thus, $(g_i,F_{i,1-\sigmabar(g_i)})\notin I_{\sigma}$ as $\sigma(g_i)=F_{i,0}$ by \Cref{definition: Canonical Strategy}.
Let $G_n=M_n$.
Then, $\valu_{\sigma}^\M(F_{i,\bit_{i+1}})>\valu_{\sigma}^\M(F_{i,1-\bit_{i+1}})$ by \Cref{lemma: Both CC Open For MDP}, also implying the statement since $\sigma(g_i)=F_{i,\bit_{i+1}}$ by \Cref{definition: Canonical Strategy MDP}.
Thus consider the case that $F_{i,1-\bit_{i+1}}$ is closed.
Then $\sigmabar(g_i)=1-\bit_{i+1}$ by \Cref{definition: Canonical Strategy} resp. \ref{definition: Canonical Strategy MDP}.
Since $F_{i,\bit_{i+1}}$ is not closed, \Cref{lemma: Cycle centers in Phase One} implies $\valustar_{\sigma}^*(F_{i,\bit_{i+1}})=\valustar_{\sigma}^*(b_1)$.
The statement thus follows since \Pref{USV1}$_i$ implies $\valustar_{\sigma}^*(F_{i,1-\bit_{i+1}})=\ubracket{s_{i,1-\bit_{i+1}}}\oplus\valustar_{\sigma}^*(b_1)$.

Let $\bit_i=1$, implying $\sigma(g_i)=F_{i,\bit_{i+1}}, \sigmabar(d_{i,\bit_{i+1}})$ and $\valustar_{\sigma}^*(F_{i,\bit_{i+1}})=\valustar_{\sigma}^*(s_{i,\bit_{i+1}})$.
By the definition of a canonical strategy, $\valustar_{\sigma}^*(s_{i,\bit_{i+1}})=\ubracket{s_{i,\bit_{i+1}},h_{i,\bit_{i+1}}}\oplus\valustar_{\sigma}^*(b_{i+1})$ since $\sigma(b_{i+1})=g_{i+1}$ if and only if $\bit_{i+1}=1$, and $F_{i,1-\bit_{i+1}}$ is not closed.
Hence, by \Cref{lemma: Cycle centers in Phase One}, $\valustar_{\sigma}^*(F_{i,1-\bit_{i+1}})=\valustar_{\sigma}^*(b_1)$.
It thus suffices to show $\valustar_{\sigma}^*(b_1)\prec\valustar_{\sigma}^*(s_{i,\bit_{i+1}})$.
This however follows immediately since $\sigma(b_i)=g_i$ implies $\valustar_{\sigma}^*(s_{i,\bit_{i+1}})\subseteq\valustar_{\sigma}^*(b_1)$. 

Next consider some escape vertex $e_{i,j,k}$ with $i\in[n], j,k\in\{0,1\}$ and let $\bit$ be even.
Then $\sigma(e_{i,j,k})=b_2$, so we prove $\valu_{\sigma}^*(g_1)\preceq\valu_{\sigma}^*(b_2)$.
Since $\bit_1=0$, we have $\sigma(b_1)=b_2$ by \Pref{EV1}$_1$.
Since we however already proved that $(b_1,g_1)\notin I_{\sigma}$, we need to have $\valu_{\sigma}^*(g_1)\preceq\valu_{\sigma}^*(b_2)$.
Now let $\bit\bmod2=1$, implying $\sigma(b_1)=g_1$.
In this case, $\sigma(e_{i,j,k})=g_1$, and $\valu_{\sigma}^*(b_2)\preceq\valu_{\sigma}^*(g_1)$ follows since $(b_1,b_2)\notin I_{\sigma}$.

Consider some upper selection vertex $s_{i,j}$ with $i\in[n]$ and  $j=\bit_{i+1}$.
Then $\sigma(s_{i,j})=h_{i,j}$, so we prove $\valu_{\sigma}^*(b_1)\preceq\valu_{\sigma}^*(h_{i,j})$.
By \Pref{EV1}$_{i+1}$, we have  $\valustar_{\sigma}^*(h_{i,j})=\ubracket{h_{i,j}}\oplus\valustar_{\sigma}^*(b_{i+1})$.
There are two cases. 
If $\bit_i=0$, then we have $h_{i,j}\notin\valustar_{\sigma}^*(b_1)$.
If $\bit_{i}=1$, then we have $g_i\in\valustar_{\sigma}^*(b_1)$.
However, this implies $\valustar_{\sigma}^*(b_1)\prec\valustar_{\sigma}^*(h_{i,j})$ in either case since $\valustar_{\sigma}^*(b_{i+1})\subseteq\valustar_{\sigma}^*(b_1)$.
Now let $j\neq\bit_{i+1}$.
In this case we prove $\valu_{\sigma}^*(h_{i,j})\preceq\valu_{\sigma}^*(b_1)$.
Consider the case $j=0$ first.
Then $\valustar_{\sigma}^*(h_{i,j})=\ubracket{h_{i,j}}\oplus\valustar_{\sigma}^*(b_{i+2})$, so $W_{i+1}^*\notin\valustar_{\sigma}^*(h_{i,j})$.
In particular we then have $\bit_{i+1}=1$, implying $W_{i+1}\subseteq\valustar_{\sigma}(b_1)$.
We thus have $\valustar_{\sigma}^*(h_{i,j})\preceq\valustar_{\sigma}^*(b_1)$.
Similarly, if $j=1$, we have $g_{i+1}\in\valustar_{\sigma}^*(h_{i,j})$ and $g_{i+1}\notin\valustar_{\sigma}^*(b_1)$, implying the statement.
\end{proof}

We now begin with our discussion of the application of individual improving switches.
This is organized as follows.
For each phase, we provide a table that contains a summary of most of the statements related to the corresponding phase.
Each row of such a table is then proven by an individual lemma.
Both the tables and the proofs are very technical, and it is not obvious why the strategies the strategy improvement algorithm produces have the corresponding properties.
We thus defer most of the proofs to \Cref{appendix: Proofs Exponential}.
We do not discuss the execution of the corresponding algorithm in all technical details here, but provide lemmas summarizing several of the more technical lemmas.
These lemmas then also relate the application of the individual switches to the occurrence records given in \Cref{table: Occurrence Records}.
The formal and exact description of the application of the algorithm is then given in \Cref{section: Final Formal Proofs} where the results of this section will be applied.

We refer to \Cref{figure: Start Phase 1} through \ref{figure: Final Canstrat} for visualizations of the strategies at the beginning of the different phases in the graph $S_3$.
To simplify the description of the improving switches, we define $\mathfrak{D}^{\sigma}\coloneqq\{(d_{i,j,k},F_{i,j}):\sigma(d_{i,j,k})\neq F_{i,j}\}$ as in \Cref{table: Switches at start of phase}.

\subsection{Improving switches of phase 1} 

In this phase, cycle edges $(d_{*,*,*},F_{*,*})$ and edges $(g_*.F_{*,*})$ are applied.
As explained previously, we provide an overview describing the application of individual switches during phase $1$ in \Cref{table: Phase 1 Switches}.
We interpret each row of this table stating that if a strategy $\sigma$ fulfills the given conditions, applying the given switch $e$ results in a strategy $\sigmae$ that has the claimed properties.
For convenience, conditions specifying the improving switch, resp. the level or cycle center corresponding to the switch, are contained in the second column.
Note that we also include one improving switch that technically belongs to phase~2.
This is included as \Cref{table: Phase 1 Switches} then contains all statements necessary to prove that applying improving switches to $\canstrat$ yields the phase-$2$-strategy that is described in \Cref{table: Properties at start of phase,table: Switches at start of phase}.

\begin{table}[ht]
\centering
\footnotesize
\begin{tabular}{|c|c|l|}\hline
Conditions for $\sigma$														&Switch $e$																&Properties of $\sigmae$				\\\hline\hline
$F_{i,j}$ is open and $I_{\sigma}=\mathfrak{D}^{\sigma}$ 					&$(d_{i,j,k},F_{i,j})$													&Phase-$1$-strategy for $\bit$ and $I_{\sigmae}=\mathfrak{D}^{\sigmae}$ \\\hdashline

$G_n=S_n$ and $I_{\sigma}=\mathfrak{D}^{\sigma}$							&$(d_{i,1-\bit_{i+1},k},F_{i,1-\bit_{i+1}})$							&Phase-$1$-strategy for $\bit$\\
$\sigma(g_i)=F_{i,1-\bit_{i+1}}$											&$i\neq 1$																&$I_{\sigmae}=\mathfrak{D}^{\sigmae}$\\\hdashline

$I_{\sigma}=\mathfrak{D}^{\sigma}$	and $\sigma(g_i)=F_{i,\bit_{i+1}}$		&$(d_{i,1-\bit_{i+1},k},F_{i,1-\bit_{i+1}})$							&Phase-$1$-strategy for $\bit$\\
$\sigma(d_{i,1-\bit_{i+1},1-k})=F_{i,1-\bit_{i+1}}$							&$\bit_i=0$																&$I_{\sigmae}=\mathfrak{D}^{\sigmae}\cup\{(g_i,F_{i,1-\bit_{i+1}})\}$\\\hdashline

$I_{\sigma}=\mathfrak{D}^{\sigma}\cup\{(g_{i},F_{i,1-\bit_{i+1}})\}$		&$(g_{i},F_{i,1-\bit_{i+1}})$											&Phase-$1$-strategy for $\bit$\\
$F_{i,j}$ is closed															&$i\neq 1 \wedge \bit_i=0$												&$I_{\sigmae}=I_{\sigma}\setminus\{e\}=\mathfrak{D}^{\sigmae}$\\\hdashline

																			&\multirow{9}{*}{$(d_{\nsb,\bit_{\nsb+1},k},F_{\nsb,\bit_{\nsb+1}})$}	&$\nsb=1\Rightarrow$ Phase-$3$-strategy for $\bit$\\
																			&																		&$\nsb=1\wedge\sigma(g_{\nsb})=F_{\nsb,\bit_{\nsb+1}}$ imply\\
																			&																		&\hspace*{0.75em}$I_{\sigmae}=\mathfrak{D}^{\sigmae}\cup\{(b_{1},g_1)\}\cup\{(e_{*,*,*},g_1)\}$\\
\multirow{2}{*}{$I_{\sigma}=\mathfrak{D}^{\sigma}$}											&																		&$\nsb>1\Rightarrow$ Phase-$2$-strategy for $\bit$\\
\multirow{2}{*}{$\sigma(d_{\nsb,\bit_{\nsb+1},1-k})=F_{\nsb,\bit_{\nsb+1}}$}					&																		&$\nsb>1\wedge\sigma(g_{\nsb})=F_{\nsb,\bit_{\nsb+1}}$ imply\\
																			&																		&\hspace*{0.75em}$I_{\sigmae}=\mathfrak{D}^{\sigmae}\cup\{(b_{\nsb},g_{\nsb}), (s_{\nsb-1,1},h_{\nsb-1,1})\}$\\																																					
																			&																		&$\sigma(g_{\nsb})\neq F_{\nsb,\bit_{\nsb+1}}$ implies\\
																			&																		&\hspace*{0.75em}$I_{\sigmae}=\mathfrak{D}^{\sigmae}\cup\{(g_{\nsb},F_{\nsb,\bit_{\nsb+1}})\}$\\
																			&																		&\hspace*{0.75em}Pseudo phase-$2$- resp. phase-$3$-strategy\\\hdashline
Pseudo phase-$2$-strategy and $\nsb>1$									&\multirow{2}{*}{$(g_{\nsb},F_{\nsb,\bit_{\nsb+1}})$}					&Phase-$2$-strategy for $\bit$\\
$I_{\sigma}=\mathfrak{D}^{\sigma}\cup\{(g_{\nsb},F_{\nsb,\bit_{\nsb+1}})\}$	&																		&$I_{\sigmae}=\mathfrak{D}^{\sigmae}\cup\{(b_{\nsb},g_{\nsb}), (s_{\nsb-1,1},h_{\nsb-1,1})\}$\\\hline
\end{tabular}
\caption[Improving switches applied during phase 1.]{Improving switches applied during phase 1.
For convenience, we always assume $\sigma\in\reach{\sigma_0}$ and that $\sigma$ is a phase-$1$-strategy for $\bit$ if not stated otherwise.
We thus also always have $\sigmae\in\reach{\sigma_0}$.} \label{table: Phase 1 Switches}
\end{table}

The first lemma shows that performing switches at cycle vertices that do not close any cycle centers does not create any new improving switches and does not make existing switches unimproving.

\begin{restatable}[First row of \Cref{table: Phase 1 Switches}]{lemma}{NorClosingCycleCenterInPhaseOne} \label{lemma: Not closing cycle center in phase one}
Let $\sigma\in\reach{\sigma_0}$ be a well-behaved phase-$1$-strategy for $\bit\in\bitset_n$ with $I_{\sigma}=\mathfrak{D}^{\sigma}$.
Let $i\in[n], j,k\in\{0,1\}$ such that $e\coloneqq(d_{i,j,k},F_{i,j})\in I_{\sigma}$ and $\sigma(d_{i,j,1-k})\neq F_{i,j}$.
Then $\sigmae$ is a well-behaved phase-$1$-strategy for $\bit$ with $\sigmae\in\reach{\sigma_0}$ and $I_{\sigmae}=\mathfrak{D}^{\sigmae}$.
\end{restatable}

The next lemma describes what happens when the inactive cycle center $F_{i,1-\indbit_{i+1}^{\sigma}}$ is closed under the assumption that the selector vertex of level $i$ points towards this cycle center.
This happens when cycle centers with a low occurrence record have to \enquote{catch up}.
We exclude level $1$ here since the edges of the cycle centers in this level switch sufficiently often.
Consequently, this behavior does not occur for $i=1$.
Also, we only need to consider this for $G_n=S_n$ since it cannot happen that $g_i$ points towards $F_{i,1-\indbit_{i+1}^{\sigma}}$ if $G_n=M_n$.

\begin{restatable}[Second row of \Cref{table: Phase 1 Switches}]{lemma}{ClosingWhenSelectorIsPointing} \label{lemma: Closing when selector is pointing}
Let $G_n=S_n$.
Let $\sigma\in\reach{\sigma_0}$ be a well-behaved phase-$1$-strategy for $\bit\in\bitset_n$ with $I_{\sigma}=\mathfrak{D}^{\sigma}$.
Let $i\in[n], j,k\in\{0,1\}$ such that $e\coloneqq(d_{i,j,k},F_{i,j})\in I_{\sigma}$ and $\sigma(d_{i,j,1-k})=F_{i,j}, i\neq 1, j\neq\bit_{i+1}$ as well as $\sigma(g_{i})=F_{i,j}$.
Then $\sigmae$ is a well-behaved phase-$1$-strategy for $\bit$ with $I_{\sigmae}=\mathfrak{D}^{\sigmae}$ and $\sigmae\in\reach{\sigma_0}$.
\end{restatable}

The next lemma describes what happens when the inactive cycle center $F_{i,1-\indbit_{i+1}^{\sigma}}$ of some level $i\in[n-1]$ is closed under the assumption that the selector vertex of level $i$ does \emph{not} point towards that cycle center.
In this case, the valuation of $F_{i,1-\indbit^{\sigma}_{i+1}}$ increases significantly, making the switch $(g_i,F_{i,1-\indbit^{\sigma}_{i+1}})$ improving.

\begin{restatable}[Third row of \Cref{table: Phase 1 Switches}]{lemma}{ClosingWhenSelectorIsNotPointing} \label{lemma: Closing when selector is not pointing}
Let $\sigma\in\reach{\sigma_0}$ be a well-behaved phase-$1$-strategy for $\bit\in\bitset_{n}$ with $I_{\sigma}=\mathfrak{D}^{\sigma}$.
Let $i\in[n-1], j,k\in\{0,1\}$ such that $e\coloneqq(d_{i,j,k},F_{i,j})\in I_{\sigma}$ and $\sigma(d_{i,j,1-k})=F_{i,j}, j=1-\indbit_{i+1}^{\sigma},\sigma(b_i)=b_{i+1}$ and $\sigma(g_i)=F_{i,1-j}$.
Then $\sigmae$ is a well-behaved phase-$1$-strategy for $\bit$ with $\sigmae\in\reach{\sigma_0}$ and $I_{\sigmae}=\mathfrak{D}^{\sigmae}\cup\{(g_{i},F_{i,j})\}$.
\end{restatable}

It can thus happen that improving switches $(g_i,F_{i,j})$ are created.
We prove that applying this switch again yields a strategy $\sigma$ with $I_{\sigma}=\mathfrak{D}^{\sigma}$.

\begin{restatable}[Fourth row of \Cref{table: Phase 1 Switches}]{lemma}{SwitchingSelectorInPhaseOne} \label{lemma: Switching Selector in Phase One}
Let $\sigma\in\reach{\sigma_0}$ be a well-behaved phase-$1$-strategy for $\bit\in\bitset_n$ with $I_{\sigma}=\mathfrak{D}^{\sigma}\cup\{(g_{i},F_{i,1-\bit_{i+1}})\}$ for some index $i\in[n-1]$.
Let $e\coloneqq(g_i,F_{i,1-\bit_{i+1}})\in I_{\sigma}$ and $\bit_i=0, i\neq 1$ and $\sigmabar(d_{i,j})$.
Then $\sigmae$ is a well-behaved phase-$1$-strategy for $\bit$ with $I_{\sigmae}=I_{\sigma}\setminus\{e\}$.
\end{restatable}

This now allows us to formalize the application of the first set of improving switches that are applied during phase $1$.

\begin{restatable}{lemma}{PhaseOneLowOR} \label{lemma: Phase 1 Low OR}
Let $\sigma\in\reach{\canstrat}$ be a well-behaved phase-$1$-strategy for $\bit$ with $I_{\sigma}=\mathfrak{D}^{\sigma}$.
Let $\canstrat\in\reach{\sigma_0}$ and let $\canstrat$ have the canonical properties.
Let $i\in[n], j,k\in\{0,1\}$ such that $e\coloneqq(d_{i,j,k},F_{i,j})\in I_{\sigma},I_{\canstrat}$ with $\occrec^{\sigma}(e)=\occrec^{\canstrat}(e)=\floor{(\bit+1)/2}-1$.
Then $\sigmae$ is a well-behaved phase-$1$-strategy for $\bit$ with $\sigmae\in\reach{\sigma_0}$.
Furthermore, $\sigma(d_{i,j,1-k})=F_{i,j}, j\neq\bit_{i+1}, \sigma(g_i)=F_{i,1-j}$ and $\sigma(b_i)\neq g_i$ imply $I_{\sigmae}=(I_{\sigma}\setminus\{e\})\cup\{(g_i,F_{i,j})\}$.
Otherwise, $I_{\sigmae}=I_{\sigma}\setminus\{e\}$.
In addition, the occurrence record of $e$ with respect to $\sigmae$ is described correctly by \Cref{table: Occurrence Records} when interpreted for $\bit+1$.
\end{restatable}

The next lemma now formalizes the last row of \Cref{table: Phase 1 Switches}.
It describes what happens when the cycle center $F_{\nsb,\bit_{\nsb+1}}$ is closed, concluding phase $1$.

\begin{restatable}[Fifth row of \Cref{table: Phase 1 Switches}]{lemma}{ClosingActiveCC} \label{lemma: Closing active CC}
Let $\sigma\in\reach{\sigma_0}$ be a well-behaved phase-$1$-strategy for $\bit\in\bitset_n$ and $I_{\sigma}=\mathfrak{D}^{\sigma}$.
Let $\nsb\coloneqq\ell(\bit+1)$ and $j\coloneqq\bit_{\nsb+1}$.
Let $e\coloneqq(d_{\nsb,j,k},F_{\nsb,j})\in I_{\sigma}$ and $\sigma(d_{\nsb,j,1-k})=F_{\nsb,j}$ for some $k\in\{0,1\}$.
The following statements hold.
\begin{enumerate}
	\item $\indbit^{\sigmae}=\bit+1$.
	\item $\sigmae$ has Properties (\ref{property: EV1})$_{i}$ and (\ref{property: EV3})$_{i}$ for all $i>\nsb$.
		It also has \Pref{EV2}$_{i}$ and \Pref{USV1}$_{i}$ for all $i\geq\nsb$ as well as \Pref{REL1}, and $\relbit{\sigmae}=\relbit{\sigma}=\nsb$.
	\item $\sigmae$ is well-behaved and $\sigmae\in\reach{\sigma_0}$.
	\item If $\nsb=1$, then $\sigmae$ is a phase-$3$-strategy for $\bit$.
		If $\sigma(g_\nsb)=F_{\nsb,j}$, then it holds that  $I_{\sigmae}=\mathfrak{D}^{\sigmae}\cup\{(b_1,g_1)\}\cup\{(e_{*,*,*},g_1)\}.$
		If $\sigma(g_\nsb)\neq F_{\nsb,j}$, then $I_{\sigmae}=\mathfrak{D}^{\sigmae}\cup\{(g_{\nsb},F_{\nsb,j})\}$ and $\sigmae$ is a pseudo phase-$3$-strategy.
	\item If $\nsb>1$, then $\sigmae$ is a phase-$2$-strategy for $\bit$.
		If $\sigma(g_\nsb)=F_{\nsb,j}$, then it holds that $I_{\sigmae}=\mathfrak{D}^{\sigmae}\cup\{(b_\nsb,g_\nsb)\}\cup\{(s_{\nsb-1,1},h_{\nsb-1,1})\}.$
		If $\sigma(g_\nsb)\neq F_{\nsb,j}$, then $I_{\sigmae}=\mathfrak{D}^{\sigmae}\cup\{(g_{\nsb},F_{\nsb,j})\}$ and $\sigmae$ is a pseudo phase-$2$-strategy.
\end{enumerate}
\end{restatable}

The final statement contained in \Cref{table: Phase 1 Switches} does technically not belong to Phase 1.
It considers the case that $\sigma(g_\nsb)\neq F_{\nsb,\bit_{\nsb+1}}$ when the cycle center $F_{\nsb,\bit_{\nsb+1}}$ is closed.
We show that applying $(g_i,F_{\nsb,\bit_{\nsb+1}})$ then results in the same strategy that would be achieved if $\sigma(g_\nsb)=F_{\nsb,j}$ already held.

\begin{restatable}[Sixth row of \Cref{table: Phase 1 Switches}]{lemma}{TransitionToPhaseTwo} \label{lemma: Transition to Phase Two resp Three}
Let $\sigma\in\reach{\sigma_0}$ be a well-behaved pseudo phase-$2$-strategy for $\bit\in\bitset_n$ with $\nsb>1$.
Let $e\coloneqq(g_{\nsb},F_{\nsb,\bit_{\nsb+1}})$ and $I_{\sigma}=\mathfrak{D}^{\sigma}\cup\{(g_{\nsb},F_{\nsb,\bit_{\nsb+1}})\}.$
Assume that $\sigma$ has \Pref{REL1}.
Then $\sigmae$ is a well-behaved phase-$2$-strategy for $\bit$ with $\sigmae\in\reach{\sigma_0}$ and $I_{\sigmae}=\mathfrak{D}^{\sigmae}\cup\{(b_\nsb,g_\nsb), (s_{\nsb-1,1},h_{\nsb-1,1})\}$.
\end{restatable}

This concludes our discussion of the application of improving switches that potentially yield a phase-$2$-strategy for $\bit$ as described by \Cref{table: Properties at start of phase,table: Switches at start of phase}.
We next provide the lemmas necessary for proving that the strategy improvement algorithm reaches a phase-$3$-strategy regardless of whether we have $G_n=S_n$ or $G_n=M_n$ and of the parity of $\bit$.
This is done by investigating the improving switches of phase $2$ as well as proving how a \enquote{real} phase-$3$-strategy can be obtained by the respective algorithm only yields a pseudo phase-$3$-strategy at the end of phase $1$.

\subsection{Improving switches of phase 2} 
During phase~$2$, the entry vertices $b_i$ of levels $i\in\{2,\dots,\nsb\}$ and the upper selection vertices $s_{i,(\bit+1)_{i+1}}$ of levels $i\leq\nsb-1$ are updated.
We again provide an overview describing the application of individual improving switches during phase~$2$ as well as the application of the switch $(g_{\nsb},F_{\nsb,(\bit+1)_{\nsb+1}})$ if the algorithm produces a pseudo phase-$3$-strategy.

{\begin{table}[ht]
\centering
\footnotesize
\renewcommand{\arraystretch}{1.05}
\begin{tabular}{|c|c|l|}\hline
Properties of $\sigma$																																																	&Switch $e$								&Properties of $\sigmae$				\\\hline\hline
																																																			&\multirow{5}{*}{$(b_{\nsb},g_{\nsb})$}		&Phase-$2$-strategy for $\bit$\\
$I_{\sigma}=\mathfrak{D}^{\sigma}\cup\{(b_{\nsb},g_{\nsb}),(s_{\nsb-1,1},h_{\nsb-1,1})\}$		&																								&$\nsb\neq 2$ implies\\
\Pref{REL1}																																													&																								&\hspace*{0.75em}$I_{\sigmae}=(I_{\sigma}\setminus\{e\})\cup\{(b_{\nsb-1},b_{\nsb}),(s_{\nsb-2,0},h_{\nsb-2,0})\}$\\
\Pref{USV3}$_i$ for all $i<\nsb$																																		&																								&$\nsb=2$ implies\\
																																																			&																								&\hspace*{0.75em}$I_{\sigmae}=\mathfrak{D}^{\sigmae}\hspace*{-0.2pt}\cup\hspace*{-0.2pt}\{(b_1,b_2),(s_{1,1},h_{1,1})\}\hspace*{-0.2pt}\cup\hspace*{-0.2pt}\{(e_{*,*,*},b_2)\}$\\\hdashline

$i'<\relbit{\sigma}\Rightarrow F_{i',\sigmabar(g_{i'})}$ is closed												&$(s_{i,j},h_{i,j})$									&$i\neq 1\Rightarrow$ Phase-$2$-strategy for $\bit$\\
\Pref{USV3}$_{i'}$ for all $i'\leq i$																	&$i<\relbit{\sigma}$									&$i=1\Rightarrow$ Phase-$3$-strategy for $\bit$\\
Properties~(\ref{property: EV1})$_{\relbit{\sigma}}$, (\ref{property: EV1})$_{i+1}$	&$j=\indbit^{\sigma}_{i+1}$								&$I_{\sigmae}=I_{\sigma}\setminus\{e\}$\\\hdashline

\multirow{2}{*}{$i'<\relbit{\sigma}\Rightarrow F_{i',\sigmabar(g_{i'})}$ is closed}												&												&Phase-$2$-strategy for $\bit$\\
\multirow{2}{*}{\Pref{USV3}$_{i'}$ for all $i'\leq i$}																												&$(b_i,b_{i+1})$				&$i\neq 2$ implies \\
\multirow{2}{*}{$i'>i\Rightarrow$ Properties (\ref{property: EV1})$_{i'}\wedge$(\ref{property: EV2})$_{i'}$}			&	$i>1$								&\hspace*{0.75em}{$I_{\sigmae}=(I_{\sigma}\setminus\{e\})\cup\{(b_{i-1},b_{i}),(s_{i-2,0},h_{i-2,0})\}$}\\
\multirow{2}{*}{$i'>i,i'\neq\relbit{\sigma}\Rightarrow$\Pref{EV3}$_{i'}$}																&$i<\relbit{\sigma}$	&$i=2$ implies\\
																																																								&												&\hspace*{0.75em}{$ I_{\sigmae}=(I_{\sigma}\setminus\{e\})\cup\{(b_1,b_2)\}\cup\{(e_{*,*,*},b_2)\}$}\\\hdashline

Pseudo phase-$3$-strategy and $\nsb=1$																		&\multirow{2}{*}{$(g_{\nsb},F_{\nsb,\bit_{\nsb+1}})$}	&Phase-$3$-strategy for $\bit$		\\
$I_{\sigma}=\mathfrak{D}^{\sigma}\cup\{(g_{\nsb},F_{\nsb,\bit_{\nsb+1}})\}$										&														&$I_{\sigmae}=\mathfrak{D}^{\sigmae}\cup\{(b_1,g_1)\}\cup\{(e_{*,*,*},g_1)\}$\\\hline	
\end{tabular}
\caption[Improving switches applied during phase $2$.]{Improving switches applied during phase $2$.
For convenience, we always assume $\sigma\in\reach{\sigma_0}$, that $\sigma$ is a phase-$2$-strategy for $\bit$ and that $\nsb>1$ if not stated otherwise.
We thus also always have $\sigmae\in\reach{\sigma_0}$.
We also include one application here that technically belongs to phase $3$.} \label{table: Phase 2 Switches}
\end{table}
}

We now formalize and prove the statements summarized in \Cref{table: Phase 2 Switches}.
We begin by describing the application of $(b_{\nsb},g_{\nsb})$. 

\begin{restatable}[First row of \Cref{table: Phase 2 Switches}]{lemma}{BeginningOfPhaseTwo} \label{lemma: Beginning of Phase Two}
Let $\sigma\in\reach{\sigma_0}$ be a well-behaved phase-$2$-strategy for $\bit\in\bitset_n$ with $\nsb>1$.
Let $I_{\sigma}=\mathfrak{D}^{\sigma}\cup\{(b_{\nsb},g_{\nsb}),(s_{\nsb-1,1},h_{\nsb-1,1})\}.$
Let $\sigma$ have \Pref{REL1} as well as \Pref{USV3}$_{i}$  for all $i<\nsb$.
Let $e\coloneqq (b_\nsb,g_\nsb)$.
Then, $\sigmae$ is a well-behaved phase-$2$-strategy for $\bit$ with $\sigmae\in\reach{\sigma_0}$. 
In addition ,$\nsb\neq 2$ implies \[I_{\sigmae}=\mathfrak{D}^{\sigmae}\cup\{(b_{\nsb-1},b_\nsb),(s_{\nsb-1,1},h_{\nsb-1,1}),(s_{\nsb-2,0},h_{\nsb-2,0})\}\] if $\nsb\neq 2$ and $\nsb=2$ implies \[I_{\sigmae}=\mathfrak{D}^{\sigmae}\cup\{(b_1,b_2),(s_{1,1},h_{1,1})\}\cup\{(e_{*,*,*},b_2)\}.\]
\end{restatable}

The following lemma describes the application of switches $(s_{i,j},h_{i,j})$ for $i\in[\relbit{\sigma}-1]$ and $j=\indbit^{\sigma}_{i+1}$.
Depending on whether $i\neq1$ or $i=1$, applying this switch might conclude phase $2$ and thus lead to a phase-$3$-strategy for $\bit$.
As the following lemma describes a strategy that is obtained after the application of several improving switches during phase~$2$, we include several additional assumptions that encode the application of these previously applied switches.

\begin{restatable}[Second row of \Cref{table: Phase 2 Switches}]{lemma}{UpperSelectionVerticesInPhaseTwo} \label{lemma: Upper Selection Vertices in Phase Two}
Let $\sigma\in\reach{\sigma_0}$ be a well-behaved phase-$2$-strategy for some $\bit\in\bitset_n$ with $\nsb>1$.
Assume that $\sigmabar(d_{i'})=1$ for all $i'<\relbit{\sigma}$ and that $e=(s_{i,j},h_{i,j})\in I_{\sigma}$ for some $i\in[\relbit{\sigma}-1]$ where $j\coloneqq\indbit^{\sigma}_{i+1}$.
Further assume that $\sigma$ has \Pref{USV3}$_{i'}$ for all $i'\leq i$.
Also, assume that $\sigma$ has Properties (\ref{property: EV1})$_{\relbit{\sigma}}$ and (\ref{property: EV1})$_{i+1}$.
If $i\neq 1$, then $\sigmae$ is a well-behaved phase-$2$-strategy for $\bit$.
If $i=1$, then $\sigmae$ is a well-behaved phase-$3$-strategy for $\bit$.
In either case, $I_{\sigmae}=I_{\sigma}\setminus\{e\}$. 
\end{restatable}

The following lemma describes the application of an improving switch $(b_i,b_{i+1})$ for levels $i\in\{2,\dots,\nsb-1\}$ during phase~$2$.

\begin{restatable}[Third row of \Cref{table: Phase 2 Switches}]{lemma}{ResettingEntryVertices} \label{lemma: Resetting Entry Vertices}
Let $\sigma\in\reach{\sigma_0}$ be a well-behaved phase-$2$-strategy for $\bit\in\bitset_n$ with $\nsb>1$.
Assume that $\sigmabar(d_{i'})=1$ for all $i'<\relbit{\sigma}$ and $e=(b_i,b_{i+1})\in I_{\sigma}$ for some $i\in\{2,\dots,\relbit{\sigma}-1\}$.
In addition, assume that $\sigma$ has \Pref{USV3}$_{i'}$ for all $i'<i$, \Pref{EV1}$_{i'}$ and \Pref{EV2}$_{i'}$ for all $i'>i$ as well as \Pref{EV3}$_{i'}$ for all $i'>i, i'\neq\relbit{\sigma}$.

Then $\sigmae$ is a well-behaved phase-$2$-strategy for $\bit$.
Furthermore, $i\neq2$ implies \[I_{\sigmae}=(I_{\sigma}\setminus\{e\})\cup\{(b_{i-1},b_i),(s_{i-2,0},h_{i-2,0})\}\] and $i=2$ implies $I_{\sigmae}=(I_{\sigma}\setminus\{e\})\cup\{(b_1,b_2)\}\cup\{(e_{*,*,*},b_2)\}.$
\end{restatable}

This concludes our overview related to the improving switches applied during phase~$2$.
The next lemma considers a special case that can occur at the beginning of phase $3$.
Although we closed the cycle center $F_{\nsb,(\bit+1)_{\nsb+1}}$ at the end of phase 1, it is not guaranteed that the selection vertex of level $\nsb$ points towards this cycle center if $\nsb=1$.
That is, it is not guaranteed that we immediately obtain a \enquote{proper} phase-$3$-strategy.
Such a strategy is then called \emph{pseudo} phase-$3$-strategy.
If the first phase-$3$-strategy is a pseudo phase-$3$-strategy, then the improving switch $(g_{\nsb},F_{\nsb,(\bit+1)_{\nsb-1}})$ will be applied immediately at the beginning of phase $3$.
The lemma thus describes the last row of \Cref{table: Phase 2 Switches}.

\begin{restatable}[Last row of \Cref{table: Phase 2 Switches}]{lemma}{PossibleBeginningOfPhaseThree} \label{lemma: Possible Beginning of Phase 3}
Let $\sigma\in\reach{\sigma_0}$ be a well-behaved pseudo phase-$3$-strategy for some $\bit\in\bitset_n$ with $\nsb=1$.
Let $I_{\sigma}=\mathfrak{D}^{\sigma}\cup\{(g_{\nsb},F_{\nsb,\bit_{\nsb+1}})\}$ and $e\coloneqq(g_\nsb,F_{\nsb,\bit_{\nsb+1}})$.
Then $\sigmae$ is a well-behaved phase-$3$-strategy for $\bit$ with $\sigmae\in\reach{\sigma_0}$ and \[I_{\sigmae}=(I_{\sigma}\setminus\{e\})\cup\{(b_1,g_1)\}\cup\{(e_{*,*,*},g_1)\}.\]
\end{restatable}

These are all lemmas necessary for describing phase $2$.
We consider the statements related to the application of the improving switches during phase $3$ next.

\subsection{Improving switches of phase 3}

We now discuss the application of improving switches during phase 3, which highly depends on whether we have $G_n=S_n$ or $G_n=M_n$ and on the least significant set bit of $\bit+1$.
As usual, we provide an overview describing the application of individual improving switches during phase~$3$.
To simplify and unify the arguments, we define $t^{\rightarrow}\coloneqq b_2$ if $\nsb>1$ and $t^{\rightarrow}\coloneqq g_1$ if $\nsb=1$.
Similarly, let $t^{\leftarrow}\coloneqq g_1$ if $\nsb>1$ and $t^{\leftarrow}\coloneqq b_2$ if $\nsb=1$.
We furthermore define $\mathfrak{E}^{\sigma}\coloneqq\{(d_{i,j,k},F_{i,j}), (e_{i,j,k},t^{\rightarrow})\colon\sigma(e_{i,j,k})=t^{\leftarrow}\}$.

{\begin{table}
\footnotesize
\centering
\begin{tabular}{|c|c|l|}\hline
Properties of $\sigma$																								&Switch $e$											&Properties of $\sigmae$				\\\hline\hline
																																										&\multirow{6}{*}{$(e_{i,j,k},t^{\rightarrow})$}		&Phase-$3$-strategy for $\bit$\\
\multirow{2}{*}{If $G_n=S_n$: \Pref{USV2}$_{i',*} \forall i'<\relbit{\sigma}$}	&																											&$\sigma(d_{i,j,1-k})=e_{i,j,1-k}$\\
\multirow{2}{*}{If $G_n=M_n$: $\sigma(s_{i',*})=b_1$ implies}									&																															&\hspace*{0.333em}$\vee[\sigma(d_{i,j,1-k})=F_{i,j}\wedge j\neq\indbit^{\sigma}_{i+1}]$ imply \\
\multirow{2}{*}{$\sigmabar(eb_{i',*})\wedge\nsigmabar(eg_{i',*})\mid \forall i'<\relbit{\sigma}$}	&																											&\hspace*{1em}$I_{\sigmae}=(I_{\sigma}\setminus\{e\})\cup\{(d_{i,j,k},e_{i,j,k})\}$\\
																																										&																											&$\sigma(d_{i,j,1-k})=F_{i,j}\wedge j=\indbit^{\sigma}_{i+1}$ imply \\
																																										&																											&\hspace*{1em}$I_{\sigmae}=I_{\sigma}\setminus \{e\}$\\\hdashline

$G_n=S_n$																														&														&\multirow{2}{*}{Phase-$3$-strategy for $\bit$}\\
$\sigmabar(d_{i,j})\implies j\neq\indbit^{\sigma}_{i+1}$					&$(d_{i,j,k},e_{i,j,k})$			&\multirow{2}{*}{$I_{\sigmae}=I_{\sigma}\setminus\{e\}$}\\
$\sigma(e_{i,j,k})=t^{\rightarrow}$																	&														&\\\hdashline
			
\multirow{2}{*}{$G_n=M_n$}																													&$(d_{i,j,k},e_{i,j,k})$							&\multirow{2}{*}{Phase-$3$-strategy for $\bit$}\\
\multirow{2}{*}{$\sigma(e_{i,j,k})=t^{\rightarrow}$}																&$\indbit^{\sigma}_i=1$					&\multirow{2}{*}{$I_{\sigmae}=I_{\sigma}\setminus\{e\}$}\\
																								
																																												&$j=\indbit^{\sigma}_{i+1}$			&\\\hdashline		
																								
$G_n=M_n$ and $\sigma(g_i)=F_{i,1-j}$																								&$(d_{i,j,k},e_{i,j,k})$							&\multirow{2}{*}{Phase-$3$-strategy for $\bit$}\\
$F_{i,j}$ is $t^{\leftarrow}$-halfopen, $F_{i,1-j}$ is $t^{\rightarrow}$-open 				&$\indbit^{\sigma}_i=0$								&\multirow{2}{*}{$I_{\sigmae}=I_{\sigma}\setminus\{e\}$}\\
$\sigma(e_{i,j,k})=t^{\rightarrow}$																					&$j=\indbit^{\sigma}_{i+1}$							&\\\hdashline

$G_n=M_n$																										&$(s_{i,j},b_1)$									&\multirow{2}{*}{Phase-$3$-strategy for $\bit$}\\
$\nsb>1$																											&$i<\nsb$											&\multirow{2}{*}{$I_{\sigmae}=(I_{\sigma}\setminus\{e\})$}\\
$\sigmabar(eb_{i,j})\wedge\nsigmabar(eg_{i,j})$																		&$j=1-\indbit^{\sigma}_{i+1}$						&\\\hdashline

$G_n=M_n$ and $\sigma(e_{i,j,k})=t^{\rightarrow}$																	&\multirow{2}{*}{$(d_{i,j,k},e_{i,j,k})$}			&\multirow{2}{*}{Phase-$3$-strategy for $\bit$}\\
$F_{i,j}$ is $t^{\rightarrow}$-halfopen																		&\multirow{2}{*}{$j=1-\indbit^{\sigma}_{i+1}$}		&\multirow{2}{*}{$I_{\sigmae}=(I_{\sigma}\setminus\{e\})$}\\
$\indbit^{\sigma}_i=0\Rightarrow[\sigma(g_i)=F_{i,j}\wedge F_{i,1-j}$ is $t^{\leftarrow}$-halfopen$]$ 					&													&\\\hdashline

$G_n=S_n$, $\nsb>1$ and $I_{\sigma}=\mathfrak{E}^{\sigma}\cup\{(b_1,b_2)\}$																&													&\\
$\sigma(d_{i,j,k})=F_{i,j}\Leftrightarrow\indbit^{\sigma}_{i}=1\wedge\indbit^{\sigma}_{i+1}=j$			&													&Phase-$4$-strategy for $\bit$ with $\relbit{\sigmae}=1$\\
(\ref{property: ESC4})$_{i,j}$ for all $(i,j)\in S_1$																												&$(b_1,b_2)$						&$I_{\sigmae}=(I_{\sigma}\setminus\{e\})\cup\{(s_{\nsb-1,0},b_1)\}$\\
(\ref{property: ESC5})$_{i,j}$ for all $(i,j)\in S_2$																												&													&\hspace*{0.75em}$\cup\{(s_{i,1},b_1)\colon i\leq\nsb-2\}\cup X_0\cup X_1$\\
$i<\relbit{\sigma}\Rightarrow\sigma(s_{i,*})=h_{i,*}$																										&													&\\\hdashline

$G_n=M_n$,  $\nsb>1$ and $I_{\sigma}=\mathfrak{E}^{\sigma}\cup\{(b_1,b_2)\}$														&\multirow{4}{*}{$(b_1,b_2)$}						&\multirow{2}{*}{Phase-$5$-strategy for $\bit$ with $\relbit{\sigmae}=1$}\\
$\sigma(d_{i,j,k})=F_{i,j}\Leftrightarrow\indbit^{\sigma}_i=1\wedge\indbit^{\sigma}_{i+1}=j$						&													&\multirow{2}{*}{$I_{\sigmae}=(I_{\sigma}\setminus\{e\})\cup X_0\cup X_1$}\\
(\ref{property: ESC4})$_{i,j}$ for all $(i,j)\in S_1$		&													&\multirow{3}{*}{\hspace*{0.75em}$\cup\{(d_{i,-1,\indbit^{\sigma}_{i+1},*},F_{i,1-\indbit^{\sigma}_{i+1}})\colon i<\nsb\}$}\\
(\ref{property: ESC5})$_{i,j}$ for all $(i,j)\in S_2$&&\\
\Cref{property: USV1}$_i$ for all $i\in[n]$																&													&\\\hdashline

																																																				&																		&Phase-$5$-strategy with $\relbit{\sigmae}=u$\\
$\nsb=1$ and $I_{\sigma}=\mathfrak{E}^{\sigma}\cup\{(b_1,g_1)\}$																&																		&$I_{\sigmae}=(I_{\sigma}\setminus\{e\})$\\
$\sigma(d_{i,j,k})=F_{i,j}\Leftrightarrow\indbit^{\sigma}_i=1\wedge\indbit^{\sigma}_{i+1}=j$		&\multirow{2}{*}{$(b_1,g_1)$}		&\multirow{4}{*}{\hspace*{0.75em}$\displaystyle\cup\bigcup_{\substack{i'=u+1\\\indbit^{\sigma}_i=0}}^{m-1}\{(d_{i,1-\indbit^{\sigma}_{i+1},*},F_{i,1-\indbit^{\sigma}_{i+1}})\}$}\\
(\ref{property: ESC3})$_{i,j}$ for all $(i,j)\in S_4$																										&																		&\\
(\ref{property: ESC5})$_{i,j}$ for all $(i,j)\in S_3$																										&																		&\\
																																																				&																		&\\\hdashline
\end{tabular}

\caption[Improving switches applied during phase 3.]{Improving switches applied during phase 3.
For convenience, we always assume $\sigma\in\reach{\sigma_0}$ and that $\sigma$ is a phase-$3$-strategy for $\bit$ if not stated otherwise, implying that $\sigmae\in\reach{\sigma_0}$.
The definition of the sets $X_k, S_i$ can be found in \Cref{table: Definition of Phases,table: Switches at start of phase}.} \label{table: Phase 3 Switches}
\end{table}
}

There are also additional statements describing the application of the improving switches during phase $3$.
These statements are however more involved and cannot be stated in the same way the statements contained in \Cref{table: Phase 3 Switches} can be described.
We defer these statements and their discussion for the moment and begin with a lemma characterizing the vertex valuations for several phase-$3$-strategies.
As its proof is rather short and yields some interesting insights regarding phase-$3$-strategies, it is also given directly here and not deferred to the appendix.

\begin{lemma} \label{lemma: Valuations in Phase Three}
Let $\sigma\in\reach{\sigma_0}$ be a well-behaved phase-$3$-strategy for $\bit\in\bitset_n$.
\begin{enumerate}
	\item If $\nsb=1$, then $\valustar_{\sigma}^*(b_2)=L_2^*$ and $\valustar_{\sigma}^*(g_1)=W_1^*\oplus\valustar_{\sigma}^*(b_2)$, so in particular $\valu_{\sigma}^*(g_1)\succ\valu_{\sigma}^*(b_2)$.
	\item If $\nsb>1$, then $\valustar_{\sigma}^*(b_2)=L_2^*$ and $\valu_{\sigma}^*(b_2)\succ\valu_{\sigma}^*(g_1)\oplus\ubracket{s_{i,j}}$ for $i\in[n], j\in\{0,1\}$, so in particular $\valustar_{\sigma}^*(b_2)\succ\valustar_{\sigma}^*(g_1)$.
\end{enumerate}
\end{lemma}

\begin{proof}
Let $\nsb=1$.
Since $\sigma$ is a phase-$3$-strategy, this implies $\valustar_{\sigma}^*(b_2)=L_2^*$ as $\relbit{\sigma}=1$.
Let $G_n=S_n$.
By \Pref{EV1}$_2$ and \Pref{CC2}, $\sigmabar(g_1)=\sigmabar(b_2)=\indbit^{\sigma}_2$.
Thus, $\indbit^{\sigma}_2=0$ implies $\sigma(g_1)=F_{1,0}$ and $\indbit_2=1$ implies $\sigma(b_2)=g_2$.
In either case, $\lambda_1^\P=1$.
We now investigate which of the cases of \Cref{corollary: Complete Valuation Of Selection Vertices PG} can occur and prove that $\valustar_{\sigma}^\P(g_1)=W_1^\P\cup\valustar_{\sigma}^\P(b_2)$ for the respective cases.
The first case cannot occur since $\relbit{\sigma}=1$ implies $\sigma(b_1)=b_2$.
The second up to the fifth case cannot occur since \Pref{REL2} and \Pref{CC2} imply $\sigmabar(d_1)$.
In addition, \Pref{CC2} implies $\sigmabar(g_1)=\indbit_2$.
Thus, by \Pref{USV1}$_1$, $\sigmabar(s_{1,\sigmabar(g_1)})=h_{1,\sigmabar(g_1)}$, so the conditions of the sixth case of \Cref{corollary: Complete Valuation Of Selection Vertices PG} cannot hold.
Consequently, the conditions of one of the last two cases of \Cref{corollary: Complete Valuation Of Selection Vertices PG} are fulfilled.
However, since $\sigma(b_2)=b_3$ if $\sigma(g_1)=F_{1,0}$ by \Pref{CC2} and \Pref{EV1}$_2$, the statement follows in either case.

Let $G_n=M_n$ and consider \Cref{corollary: Complete Valuation Of Selection Vertices MDP}.
If $\indbit_2=0$, then $\sigma(g_1)=F_{1,0}$ by \Pref{CC2}, implying $\lambda_1^\M=1$.
Since $\sigmabar(d_1)\wedge\sigmabar(s_1)$ as shown previously,it follows that $\valustar_{\sigma}^\M(g_1)=W_{1}^\M+\valustar_{\sigma}^\M(b_2)$ since the conditions of the last case are fulfilled.
If $\indbit^{\sigmae}_2=1$, then $\lambda_1^\M=2$ by \Pref{EV2}$_2$.
Consequently, $\valustar_{\sigmae}^\M(g_1)=W_1^\M+\valustar_{\sigmae}^\M(b_2)$ since the conditions of the first case are fulfilled.

This concludes the case $\nsb=1$, hence assume $\nsb>1$.
We prove $\valustar_{\sigma}^*(b_2)=L_2^*$ first.
If $\sigma(b_2)=b_3$, this follows by definition.
Hence assume $\sigma(b_2)=g_2$.
Then, by \Pref{EV1}$_2$ and \Pref{CC2}, $\sigma(g_1)=F_{1,0}$.
In addition, $\nsb=\relbit{\sigma}>1$ implies $\sigma(b_1)=g_1$.
Consequently, $1\in\incorrect{\sigma}$.
Since $\sigma$ has \Pref{EV1}$_i$ and \Pref{EV2}$_i$ for all $i>1$, no other index can be contained in $\incorrect{\sigma}$.
But this implies $\incorrect{\sigma}=\{1\}$ and thus, since $\sigma(b_2)=g_2$, $\relbit{\sigma}=2$, implying $\valustar_{\sigma}^*(b_2)=L_2^*$ as claimed.

We now prove that $\nsb>1$ implies $\valu_{\sigma}^*(b_2)\succ\valu_{\sigma}^*(g_1)\oplus\ubracket{s_{i,j}}$ for $i\in[n], j\in\{0,1\}$.
Since $\nsb>1$ implies $\sigma(b_1)=g_1$, this implies that $\valustar_{\sigma}^*(g_1)=\valustar_{\sigma}^*(b_1)$.
Furthermore, by $1<\relbit{\sigma}$ and $\sigma(b_1)=g_1$, \Cref{lemma: Valuation of b} and $B_2^*=L_2^*$ imply that either $\valustar_{\sigma}^*(b_1)=R_1^*$ or $G_n=M_n$ and $\valustar_{\sigma}^\M(b_1)=\rew{g_k}+\sum_{j=1}^{k-1}W_j^\M+L_2^\M$ where $k=\min\{i\geq 1\colon\nsigmabar(d_i)\}<\relbit{\sigma}$.
In the second case the statement follows directly, in the first it follows by \Cref{lemma: VV Lemma} since $\incorrect{\sigma}\neq\emptyset$ implies $\sigma(b_{\relbit{\sigma}})=g_{\relbit{\sigma}}$ by \Cref{lemma: Traits of Relbit}.
\end{proof}

We now begin with the lemmas describing phase $3$.
The first lemma describes the application of $(e_{i,j,k},g_1)$ resp. $(e_{i,j,k},b_2)$ for the case that $\sigma(d_{i,j,k})=F_{i,j}$.
As all of the following lemmas, this lemma contains some conditions encoding the behavior of the strategy improvement algorithm and the application of previous improving switches.
Since phase $3$ is not exactly identical for $S_n$ and $M_n$, there are also some conditions distinguishing between the two.

\begin{restatable}[First row of \Cref{table: Phase 3 Switches}]{lemma}{EscapeVerticesPhaseThree} \label{lemma: Escape Vertices Phase Three}
Let $\sigma\in\reach{\sigma_0}$ be a well-behaved phase-$3$-strategy for $\bit\in\bitset_n$.
Let $i\in[n],j,k\in\{0,1\}$ such that $(e_{i,j,k},t^\to)\in I_{\sigma}$ and $\sigma(d_{i,j,k})=F_{i,j}$.
Further assume the following.
\begin{enumerate}
	\item If $G_n=S_n$, then, $\sigma$ has \Pref{USV2}$_{i',j'}$ for all $i'<\relbit{\sigma}, j'\in\{0,1\}$.
	\item If $G_n=M_n$, then, $\sigma(s_{i',j'})=b_1$ implies $\sigmabar(eb_{i',j'})\wedge\nsigmabar(eg_{i',j'})$ for all $i'<\relbit{\sigma}$ and $ j'\in\{0,1\}$.
\end{enumerate}
Then $\sigmae$ is a well-behaved phase-$3$-strategy for $\bit$ with $\sigmae\in\reach{\sigma_0}$.
If $\sigma(d_{i,j,1-k})=e_{i,j,1-k}$ or $[\sigma(d_{i,j,1-k})=F_{i,j}$ and $j\neq\indbit^{\sigma}_{i+1}]$, then $I_{\sigmae}=(I_{\sigma}\setminus\{e\})\cup\{(d_{i,j,k},e_{i,j,k})\}.$
If $\sigma(d_{i,j,1-k})=F_{i,j}$ and $j=\indbit^{\sigma}_{i+1}$, then $I_{\sigmae}=I_{\sigma}\setminus\{e\}$.
\end{restatable}

We now want to describe the application of improving switches $(d_{i,j,k}, e_{i,j,k})$ in phase~$3$.
The main challenge regarding these switches is that there are several different cases that need to be considered when a switch of this type is applied.
We thus provide several individual lemmas that are combined later to give a lemma summarizing the application of these switches.
We first show that we always obtain a well-behaved phase-$3$-strategy.

\begin{restatable}{lemma}{PhaseThreeWellBehaved} \label{lemma: Phase 3 Well-Behaved}
Let $\sigma\in\reach{\sigma_0}$ be a well-behaved phase-$3$-strategy for $\bit$.
Let $i\in[n]$ and $j,k\in\{0,1\}$ such that  $\sigma(e_{i,j,k})=t^{\rightarrow}$ and $e\coloneqq(d_{i,j,k},e_{i,j,k})\in I_{\sigma}$.
Let $\sigma(d_{i,j,1-k})=e_{i,j,1-k}$ or $[\sigma(d_{i,j,1-k})=F_{i,j}$ and $j\neq\indbit^{\sigma}_{i+1}]$.
Then $\sigmae$ is a well-behaved phase-$3$-strategy for $\bit$ with $\sigmae\in\reach{\sigma_0}$.
\end{restatable}

\Cref{lemma: Phase 3 Well-Behaved} justifies to omit the upper index when referring to the induced bit state.
For a well-behaved phase-$3$-strategy $\sigma$ with $e\in I_{\sigma}$, we thus define $\indbit\coloneqq\indbit^{\sigma}=\indbit^{\sigmae}=\bit+1$ while we discuss phase $3$.

We now describe the application of switches $(d_{*,*,*},e_{*,*,*})$.
While this application is not hard to describe in $S_n$, it is very complex in $M_n$.
The reason is that applying these switches \emph{always} has an influence on the valuation of the cycle centers in $M_n$.
Thus, we need to carefully investigate the application of these switches and need to pay heavy attention to the exact order of application.

We begin with the application of an improving switch $(d_{*,*,*},e_{*,*,*})$ during phase 3 in $S_n$.
This lemma is significantly easier than the lemmas for the case $G_n=M_n$ as the valuation of the cycle center $F_{i,j}$ does not change when applying $(d_{i,j,k}, e_{i,j,k}), i\in[n], j,k\in\{0,1\}$  in $S_n$.

\begin{restatable}[Second row of \Cref{table: Phase 3 Switches}]{lemma}{PGCycleVerticesPhaseThree} \label{lemma: PG Cycle Vertices Phase Three}
Let $G_n=S_n$.
Let $\sigma\in\reach{\sigma_0}$ be a well-behaved phase-$3$-strategy for $\bit\in\bitset_n$.
Let $i\in[n], j,k\in\{0,1\}$ such that $e\coloneqq(d_{i,j,k},e_{i,j,k})\in I_{\sigma}$ and $\sigma(e_{i,j,k})=t^{\rightarrow}$.
Further assume that $\sigmabar(d_{i,j})$ implies $j\neq\indbit^{\sigma}_{i+1}$.
Then $\sigmae$ is a well-behaved phase-$3$-strategy for $\bit$ with $\sigmae\in\reach{\sigma_0}$ and $I_{\sigmae}=I_{\sigma}\setminus\{e\}$.
\end{restatable}

We now focus on the case $G_n=M_n$.
The next lemma describes the application of switches $(d_{i,j,k},e_{i,j,k})$ where $i\in[n], j,k\in\{0,1\}$ within levels $i$ with $\indbit^{\sigma}_i=1$.
We skip the upper index $\M$ to denote that we have $G_n=M_n$ since we exclusively consider this case.

\begin{restatable}[Third row of \Cref{table: Phase 3 Switches}]{lemma}{MDPPhaseThreeLevelIsSet} \label{lemma: MDP Phase 3 Level Is Set}
Let $G_n=M_n$ and let $\sigma\in\reach{\sigma_0}$ be a well-behaved phase-$3$-strategy for $\bit$.
Let $i\in[n]$ with $\indbit_i=1$ and let $j\coloneqq 1-\indbit_{i+1}$.
Let $e\coloneqq(d_{i,j,k},e_{i,j,k})\in I_{\sigma}$ and $\sigma(e_{i,j,k})=t^{\rightarrow}$ for some $k\in\{0,1\}$.
Then $\sigmae$ is a well-behaved phase-$3$-strategy for $\bit$ with $\sigmae\in\reach{\sigma_0}$ and $I_{\sigmae}=I_{\sigma}\setminus\{e\}$.
\end{restatable}

The next lemma describes the application of an improving switch $(d_{i,j,k},e_{i,j,k})$ within a $t^{\rightarrow}$-open cycle center.

\begin{restatable}[Fourth row of \Cref{table: Phase 3 Switches}]{lemma}{ImprovingSwitchinOtherCCPhaseThree} \label{lemma: Improving Switch in Other CC Phase Three}
Let $G_n=M_n$.
Let $\sigma\in\reach{\sigma_0}$ be a well-behaved phase-$3$-strategy for $\bit\in\bitset_n$.
Let $i\in[n]$ with $\indbit^{\sigma}_{i}=0, j\coloneqq\indbit^{\sigma}_{i+1}$ and let $F_{i,j}$ be $t^{\leftarrow}$-halfopen.
Let $F_{i,1-j}$ be $t^{\rightarrow}$-open and $\sigma(g_i)=F_{i,1-j}$.
Let $e\coloneqq(d_{i,j,k},e_{i,j,k})\in I_{\sigma}$ and $\sigma(e_{i,j,k})=t^{\rightarrow}$ with $k\in\{0,1\}$.
Then $\sigmae$ is a well-behaved phase-$3$-strategy for $\bit$ with $\sigmae\in\reach{\sigma_0}$ and $I_{\sigmae}=I_{\sigma}\setminus\{e\}$.
\end{restatable}

The next lemma describes the application of improving switches within levels $i\in[n]$ in which no cycle center is closed at the beginning of phase $3$.
The first case describes the first improving switch that is applied in such a level.
This switch is applied in the cycle center $F_{i,\sigmabar(g_i)}$ to avoid the creation of an additional improving switch at the selector vertex.
The second case describes the second improving switch that is then applied in the cycle center $F_{i,1-\sigmabar(g_i)}$.
The statement of this lemma is not included in \Cref{table: Phase 3 Switches} as it is to involved and does not fit the framework of the lemmas summarized there.

\begin{restatable}{lemma}{NoCCClosedPhaseThree} \label{lemma: No CC Closed Phase Three}
Let $G_n=M_n$.
Let $\sigma\in\reach{\sigma_0}$ be a well-behaved phase-$3$-strategy for $\bit\in\bitset_n$.
Let $i\geq\relbit{\sigma}+1$ and assume $\sigmabar(g_i)=\indbit^{\sigma}_{i+1}$.
\begin{enumerate}
	\item If both cycle centers of level $i$ are $t^{\leftarrow}$-halfopen, then let $j\coloneqq\sigmabar(g_i)$.
	\item If $F_{i,\indbit^{\sigma}_{i+1}}$ is mixed and $F_{i,1-\indbit^{\sigma}_{i+1}}$ is $t^{\leftarrow}$-halfopen, then let $j\coloneqq1-\sigmabar(g_i)$.
\end{enumerate}
In any case, assume $e\coloneqq (d_{i,j,k},e_{i,j,k})\in I_{\sigma}$ and $\sigma(e_{i,j,k})=t^{\rightarrow}$ for $k\in\{0,1\}$.
Then, $\sigmae$ is a well-behaved phase-$3$-strategy for $\bit$ with $\sigmae\in\reach{\sigma_0}$ and $I_{\sigmae}=I_{\sigma}\setminus\{e\}$.
\end{restatable}

The next lemma describes the application of a switch $(d_{i,*,*},e_{i,*,*})$ within a closed but inactive cycle center for the case that $\indbit_i=0$.
The lemma requires that the strategy $\sigma$ fulfills several rather complicated assumptions.
As usual, these assumptions somehow \enquote{encode} the order of application of the improving switches.

\begin{restatable}{lemma}{MDPPhaseThreeOpenClosedCycleCenter} \label{lemma: MDP Phase 3 Open Closed Cycle Center}
Let $G_n=M_n$.
Let $\sigma$ be a well-behaved phase-$3$-strategy for $\bit\in\bitset_n$ with {$\sigma\in\reach{\sigma_0}$}.
Let $i\in[n]$ and $j\coloneqq 1-\indbit^{\sigma}_{i+1}$.
Let $e\coloneqq(d_{i,j,k},e_{i,j,k})\in I_{\sigma}$ and $\sigma(e_{i,j,k})=t^{\rightarrow}$ for some $k\in\{0,1\}$. 
Further assume that there is no other triple of indices $i',j',k'$ with $(d_{i',j',k'},e_{i',j',k'})\in I_{\sigma}$, that $F_{i,j}$ is closed and that $\sigma$ fulfills the following assumptions:
\begin{enumerate}
	\item If $\indbit^{\sigma}_i=0$, then $\sigma(g_i)=F_{i,j}$ and $F_{i,1-j}$ is $t^{\leftarrow}$-halfopen. 
	\item $i<\relbit{\sigma}$ implies [$\sigma(s_{i,j})=h_{i,j}$ and $\sigma(s_{i',j'})=h_{i',j'}\wedge\sigmabar(d_{i'})$ for all $i'<i, j'\in\{0,1\}$] and that the cycle center $F_{i',1-\sigmabar(g_{i'})}$ is $t^{\leftarrow}$-halfopen for all $i'<i$.
		In addition, $i<\relbit{\sigma}-1$ implies $\sigmabar(eb_{i+1})$.
	\item $i'>i$ implies $\sigma(s_{i,1-\indbit^{\sigma}_{i'+1}})=b_1$. 
	\item $i'>i$ and $\indbit^{\sigma}_{i'}=0$ imply that either [$\sigmabar(g_{i'})=\indbit^{\sigma}_{i'+1}$ and $F_{i,0},F_{i,1}$ are mixed] or [$\sigmabar(g_{i'})=1-\indbit^{\sigma}_{i'+1}$, $F_{i',1-\indbit^{\sigma}_{i'+1}}$ is $t^{\rightarrow}$-open and $F_{i',\indbit^{\sigma}_{i'+1}}$ is mixed] and 
	\item $i'>i$ and $\indbit^{\sigma}_{i'}=1$ imply that $F_{i',1-\indbit^{\sigma}_{i'+1}}$ is either mixed or $t^{\rightarrow}$-open. 
\end{enumerate}
Then $\sigmae$ is a well-behaved phase-$3$-strategy for $\bit$ with $\sigmae\in\reach{\sigma_0}$ and $I_{\sigmae}=I_{\sigma}\setminus\{e\}$ if $i\geq\relbit{\sigma}$ and $I_{\sigmae}=[I_{\sigma}\cup\{(s_{i,j},b_1)\}]\setminus\{e\}$ if $i<\relbit{\sigma}$.
\end{restatable}

The next lemma describes the application of an improving switch $(s_{i,j}, b_1)$ that might be unlocked by the application of a switch $(d_{i,j,k}, e_{i,j,k})$.
In $M_n$, these switches are already implied during phase $3$ while they are applied during phase $4$ in $S_n$.
Thus, the following lemma only considers $M_n$.

\begin{restatable}[Fifth row of \Cref{table: Phase 3 Switches}]{lemma}{USVInPhaseThree} \label{lemma: USV in Phase Three}
Let $G_n=M_n$.
Let $\sigma\in\reach{\sigma_0}$ be a well-behaved phase-$3$-strategy for $\bit\in\bitset_n$ with $\nsb>1$.
Let $i<\relbit{\sigma}, j=1-\indbit^{\sigma}_{i+1}$ and $e\coloneqq(s_{i,j},b_1)\in I_{\sigma}$.
Further assume $\sigmabar(eb_{i,j})\wedge\nsigmabar(eg_{i,j})$.
Then $\sigmae$ is a well-behaved phase-$3$-strategy for $\bit$ with $I_{\sigmae}=I_{\sigma}\setminus\{e\}$ and $\sigmae\in\reach{\sigma_0}$.
\end{restatable}

The next lemma now describes the application of the second improving switch of the kind $(d_{i,j,k},e_{i,j,k})$ within a cycle center that was closed in phase 1.

\begin{restatable}[Sixth row of \Cref{table: Phase 3 Switches}]{lemma}{OpenCompletelyInPhaseThree} \label{lemma: Open Completely In Phase Three}
Let $G_n=M_n$.
Let $\sigma\in\reach{\sigma_0}$ be a well-behaved phase-$3$-strategy for $\bit\in\bitset_n$.
Let $i\in[n]$ and $j\coloneqq 1-\indbit^{\sigma}_{i+1}$.
Let $F_{i,j}$ be $t^{\rightarrow}$-halfopen and assume that $\indbit^{\sigma}_{i}=0$ implies that $F_{i,1-j}$ is $t^{\leftarrow}$-halfopen as well as $\sigma(g_i)=F_{i,j}$.
Let $e\coloneqq (d_{i,j,k},e_{i,j,k})\in I_{\sigma}$ and $\sigma(e_{i,j,k})=t^{\rightarrow}$ for $k\in\{0,1\}$.
Then $\sigmae$ is a well-behaved phase-$3$-strategy for $\bit$ with $\sigmae\in\reach{\sigma_0}$ and $I_{\sigmae}=I_{\sigma}\setminus\{e\}$. 
\end{restatable}

This concludes the discussion of the application of switches $(d_{*,*,*},e_{*,*,*})$.
The next lemma describes the end of phase 3 in $S_n$ for $\nsb> 1$.
In contrast to the Markov decision process $M_n$, none of the switches $(s_{i,1-\indbit_{i+1}},b_1)$ with $i<\relbit{\sigma}$ is applied during phase $3$.
In $S_n$, these switches only become improving after applying the switch $(b_1,b_2)$.
This then starts phase $4$ and the beginning of this phase is described by the following lemma.
We refer to \Cref{table: Definition of Phases} resp. \Cref{table: Switches at start of phase} for the definition of the sets $S_1$ and $S_2$ resp. $X_k$ that are used in the statement.

\begin{restatable}[Seventh row of \Cref{table: Phase 3 Switches}]{lemma}{StartOfPhaseFourValu} \label{lemma: Start Of Phase Four Valu}
Let $G_n=S_n$.
Let $\sigma\in\reach{\sigma_0}$ be a well-behaved phase-$3$-strategy for $\bit\in\bitset_n$ with $\nsb>1$.
Let \[I_{\sigma}=\{(b_1,b_2)\}\cup\{(d_{i,j,k},F_{i,j}),(e_{i,j,k},b_2):\sigma(e_{i,j,k})=g_1\}\] and $\sigma(d_{i,j,k})=F_{i,j}\Leftrightarrow\indbit^{\sigma}_{i}=1\wedge\indbit^{\sigma}_{i+1}=j$ for all $i\in[n], j,k\in\{0,1\}$.
Assume that $\sigma$ has \Pref{ESC4}$_{i,j}$ for all $(i,j)\in S_1$ and \Pref{ESC5}$_{i,j}$ for all $(i,j)\in S_2$.
Further assume that $\sigma(s_{i,j})=h_{i,j}$ for all $i<\nsb, j\in\{0,1\}$.
Let $e\coloneqq(b_1,b_2)$ and $m\coloneqq\max\{i\colon \indbit^{\sigma}_i=1\}$.
Then $\sigmae$ is a well-behaved phase-$4$-strategy for $\bit$ with $\relbit{\sigmae}=1$ and \[I_{\sigmae}=(I_{\sigma}\setminus\{e\})\cup\{(s_{\nsb-1,0},b_1)\}\cup\{(s_{i,1},b_1)\colon i\leq\nsb-2\}\cup X_0\cup X_1\] where $X_k$ is defined as in \Cref{table: Switches at start of phase}.
\end{restatable}

As mentioned earlier, there is no phase $4$ if $G_n=M_n$, even for $\nsb>1$.
Hence, the application of the improving switch $(b_1,b_2)$ directly yields a phase-$5$-strategy if all of the switches $(s_{i,j},b_1)$ have been applied before.

\begin{restatable}[Eighth row of \Cref{table: Phase 3 Switches}]{lemma}{PhaseThreeToFiveMDP} \label{lemma: Phase Three To Five MDP}
Let $G_n=M_n$.
Let $\sigma\in\reach{\sigma_0}$ be a well-behaved phase-$3$-strategy for $\bit\in\bitset_n$ with $\nsb>1$.
Let \[I_{\sigma}=\{(b_1,b_2)\}\cup\{(d_{i,j,k},F_{i,j}), (e_{i,j,k},b_2)\colon\sigma(e_{i,j,k})=g_1\}.\]
Let $\sigma$ have \Pref{USV1}$_{i}$ for all $i\in[n]$ and let $\sigma(d_{i,j,k})=F_{i,j}\Leftrightarrow\indbit^{\sigma}_{i}=1\wedge\indbit^{\sigma}_{i+1}=j$ for all $i\in[n], j,k\in\{0,1\}$.
Let $\sigma$ have \Pref{ESC4}$_{i,j}$ for all $(i,j)\in S_1$ and \Pref{ESC5}$_{i,j}$ for all $(i,j)\in S_2$.
Further assume that $e\coloneqq(b_1,b_2)\in I_{\sigma}$ and let $m\coloneqq\max\{i\colon\indbit^{\sigma}_i=1\}$.
Then, $\sigmae$ is a well-behaved phase-$5$-strategy for $\bit$ with $\relbit{\sigmae}=1$ and \[I_{\sigmae}=(I_{\sigma}\setminus\{e\})\cup\{(d_{i,1-\indbit^{\sigma}_{i+1},k},F_{i,1-\indbit^{\sigma}_{i+1}})\colon i<\nsb\}\cup X_0\cup X_1\]where $X_k$ is defined as in \Cref{table: Switches at start of phase}.
\end{restatable}

The next lemma now describes the direct transition from phase 3 to phase 5 for $\bit\in\bitset_n$ with $\nsb=1$.
In this case, there is no need to distinguish whether $G_n=S_n$ or $G_n=M_n$.

\begin{restatable}[Last row of \Cref{table: Phase 3 Switches}]{lemma}{PhaseThreeToFive} \label{lemma: Phase Three To Five}
Let $\sigma\in\reach{\sigma_0}$ be a well-behaved phase-$3$-strategy for $\bit\in\bitset_n$ with $\nsb=1$.
Let $I_{\sigma}=\{(b_1,g_1)\}\cup\{(d_{i,j,k},F_{i,j}),(e_{i,j,k},g_1)\colon\sigma(e_{i,j,k})=b_2\}$ and assume that $\sigma$ has \Pref{ESC5}$_{i,j}$ for all $(i,j)\in S_3$ and \Pref{ESC3}$_{i,j}$ for all $(i,j)\in S_4$.
Let $\sigma(d_{i,j,k})=F_{i,j}\Leftrightarrow\indbit^{\sigma}_i=1\wedge\indbit^{\sigma}_{i+1}=j$ for all $i\in[n], j,k\in\{0,1\}$.
Let $e\coloneqq (b_1,g_1)$ and define $ m\coloneqq\max\{i\colon\indbit^{\sigma}_i=1\}$ and $u\coloneqq\min\{i\colon\indbit^{\sigma}_i=0\}$.
Then $\sigmae$ is a well-behaved phase-$5$-strategy for $\bit$ with $\relbit{\sigmae}=u, \sigmae\in\reach{\sigma_0}$ and \[I_{\sigmae}=(I_{\sigma}\setminus\{e\})\cup\bigcup_{\substack{i'=u+1\\\indbit^{\sigma}_{i}=0}}^{m-1}\{(d_{i,1-\indbit^{\sigma}_{i+1},0},F_{i,1-\indbit^{\sigma}_{i+1}}),(d_{i,1-\indbit^{\sigma}_{i+1},1},F_{i,1-\indbit^{\sigma}_{i+1}})\}.\]
\end{restatable}

This concludes our discussion of the application of the improving switches of phase $3$.
We now discuss the improving switches that are applied during phase $4$ if this phase is present.

\subsection{Improving switches of phase 4}

As explained earlier,  it is still necessary to apply improving switches $(s_{*,*},b_1)$ in $S_n$ if $\nsb>1$.
These switches are applied during phase $4$.
Since these switches are the only switches that are applied during phase $4$, we do not provide a table summarizing the application of improving switches during this phase.
Instead, we provide the following single lemma.

\begin{restatable}{lemma}{PhaseFourComplete} \label{lemma: Phase Four Complete}
Let $G_n=S_n$.
Let $\sigma\in\reach{\sigma_0}$ be a well-behaved phase-$4$-strategy for $\bit\in\bitset_n$ with $\nsb>1$.
Assume that there is an index $i<\nsb$ such that $e\coloneqq (s_{i,j},b_1)\in I_{\sigma}$ where $j\coloneqq 1-\indbit^{\sigma}_{i+1}$.
Further assume the following:
\begin{enumerate}
	\item $\sigma$ has \Pref{USV1}$_{i'}$ for all $i'>i$.
	\item For all $i',j',k'$, it holds that $\sigma(d_{i',j',k'})=F_{i',j'}$ if and only if $\indbit^{\sigma}_{i'}=1\wedge\indbit^{\sigma}_{i'+1}=j'$.
	\item $i'<\nsb$ implies $\sigmabar(g_{i'})=1-\indbit^{\sigma}_{i'+1}$.
	\item $i'<i$ implies $\sigma(s_{i',*})=h_{i',*}$.
\end{enumerate}
If there is an index $i'<i$ such that $(s_{i',1-\indbit^{\sigma}_{i'+1}},b_1)\in I_{\sigma}$, then $\sigmae$ is a well-behaved phase-$4$-strategy for $\bit$.
Otherwise, it is a well-behaved phase-$5$-strategy for $\bit$.
In either case, it holds that $I_{\sigmae}=(I_{\sigma}\setminus\{e\})\cup\{(d_{i,j,0},F_{i,j}), (d_{i,j,1},F_{i,j})\}$.
\end{restatable}

\subsection{Improving switches of phase 5}

We now discuss the improving switches which are applied during phase 5.
As usual, we provide a table that contains one row per \enquote{type} of improving switch and provide a statement for each row of that table.
This overview is given by \Cref{table: Phase 5 Switches}.
There is one more complex statement describing the application of improving switches of the type $(e_{*,*,*}, g_1)$ resp. $(e_{*,*,*}, b_2)$ during phase $5$.
Due to its complexity, this statement is not contained in \Cref{table: Phase 5 Switches}.

We begin by providing the lemma describing the application of the improving switches involving the escape vertices.
As usual, we define $t^{\rightarrow}\coloneqq g_1\wedge t^{\leftarrow}\coloneqq b_2$ if $\nsb=1$ and $t^{\leftarrow}\coloneqq b_2\wedge t^{\rightarrow}\coloneqq g_1$ if $\nsb>1$.

\begin{restatable}{lemma}{PhaseFiveEscapeEasy} \label{lemma: Phase Five Escape Easy}
Let $\sigma\in\reach{\sigma_0}$ be a well-behaved phase-$5$-strategy for $\bit\in\bitset_n$.
Let $i\in[n]$ and $j,k\in\{0,1\}$ with $e\coloneqq (e_{i,j,k},t^{\rightarrow})\in I_{\sigma}$ and $\sigmabar(eb_{i,j})\wedge\sigmabar(eg_{i,j})$.
Furthermore assume that $G_n=S_n$ implies \[j=1\wedge\nsb>1\implies\nsigmabar(eg_{i,1-j})\quad\text{ and }\quad j=1\wedge\nsb=1\implies\nsigmabar(eb_{i,1-j}).\]
Similarly, assume that $G_n=M_n$ implies \[j=1-\indbit^{\sigma}_{i+1}\wedge\nsb>1\implies\nsigmabar(eg_{i,1-j})\quad\text{ and }\quad j=1-\indbit^{\sigma}_{i+1}\wedge\nsb=1\implies\nsigmabar(eb_{i,1-j}).\]
Moreover, assume that $\nsb=2$ implies $\sigma(g_1)=F_{1,0}$ if $G_n=S_n$.
Then the following hold.
\begin{enumerate}
	\item If there are indices $(i',j',k')\neq(i,j,k)$ with $(e_{i',j',k'},t^{\rightarrow})\in I_{\sigma}$ or if there is an index~$i'$ such that $\sigma$ does not have \Pref{SVG}$_{i'}$/(\ref{property: SVM})$_{i'}$, then $\sigmae$ is a phase-$5$-strategy for~$\bit$.
	\item The strategy $\sigmae$ is well-behaved.
	\item If there are no indies $(i',j',k')\neq(i,j,k)$ with $(e_{i',j',k'},t^{\rightarrow})\in I_{\sigma}$ and if $\sigma$ has \Pref{SVG}$_{i'}$ resp. (\ref{property: SVM})$_{i'}$ for all $i'[n]$, then $\sigmae$ is a phase-$1$-strategy for $\bit+1$.
	\item If $G_n=S_n$, then \[(g_i,F_{i,j})\in I_{\sigmae}\Longleftrightarrow\indbit^{\sigma}_i=0\wedge\sigmaebar(g_i)=1\wedge j=0\wedge\begin{cases}\sigmabar(eb_{i,1-j}), &\nsb>1\\\sigmabar(eg_{i,1-j}), &\nsb=1 \end{cases}.\]
		If $G_n=M_n$, then \[(g_i,F_{i,j})\in I_{\sigmae}\Longleftrightarrow\indbit^{\sigma}_i=0\wedge\sigmaebar(g_i)=1-\indbit^{\sigma}_{i+1}\wedge j=\indbit^{\sigma}_{i+1}\wedge\begin{cases}\sigmabar(eb_{i,1-j}), &\nsb>1\\\sigmabar(eg_{i,1-j}), &\nsb=1 \end{cases}.\]
		If the corresponding conditions are fulfilled, then \[I_{\sigmae}=(I_{\sigma}\setminus\{e\})\cup\{(d_{i,j,1-k},F_{i,j}),(g_i,F_{i,j}))\}.\]
		Otherwise, $I_{\sigmae}=(I_{\sigma}\setminus\{e\})\cup\{(d_{i,j,1-k},F_{i,j})\}$.
\end{enumerate}
\end{restatable}

{\begin{table}[ht]
\footnotesize
\centering
\begin{tabular}{|c|c|l|}\hline
Properties of $\sigma$																				&Switch $e$											&Properties of $\sigmae$				\\\hline\hline
\multirow{2}{*}{$\sigma(b_i)=b_{i+1}$}												&$(d_{i,j,k},F_{i,j})$														&\multirow{2}{*}{Phase-$5$-strategy for $\bit$}\\
\multirow{2}{*}{$\sigmabar(g_i)=1-\indbit^{\sigma}_{i+1}$}	&$i\neq 1 $																		&\multirow{2}{*}{$I_{\sigmae}=I_{\sigma}\setminus\{e\}$}\\
																																	&	$j=1-\indbit^{\sigma}_{i+1}$								&\\\hdashline																			
																																									
$\indbit^{\sigma}_i=0$																																	&\multirow{5}{*}{$(g_i,F_{i,j})$}					&\multirow{2}{*}{$\sigma(e_{*,*,*})=t^{\rightarrow}\wedge$(\ref{property: SVG})$_{i'}$/(\ref{property: SVM})$_{i'}$ $\forall i'\in[n]$}\\
$\nsb=1\Rightarrow\sigmabar(eg_{i,j})\wedge\nsigmabar(eg_{i,j})$									&																						&\multirow{2}{*}{\hspace*{0.75em}$\Rightarrow$ Phase-$1$-strategy for $\bit+1$}\\
$\nsb>1\Rightarrow\sigmabar(eb_{i,j})\wedge\nsigmabar(eg_{i,j})$									&																						&\multirow{2}{*}{Otherwise phase-$5$-strategy for $\bit$}\\
$\relbit{\sigma}=1$ implies 																															&																						&\multirow{2}{*}{$I_{\sigmae}=I_{\sigma}\setminus\{e\}$}\\
$[i'\geq i\Rightarrow\sigmabar(d_{i',*})\vee(\sigmabar(eb_{i',*})\wedge\nsigmabar(eg_{i',*})]$		&															&\\\hline
\end{tabular}
\caption[Improving switches applied during phase 5.]{Improving switches applied during phase 5.
For convenience, we always assume $\sigma\in\reach{\sigma_0}$ and that $\sigma$ is a phase-$5$-strategy for $\bit$ .
Note that we thus also always have $\sigmae\in\reach{\sigma_0}$.} \label{table: Phase 5 Switches}
\end{table}
}

The following lemma describes the application of improving switches that involve cycle vertices during phase $5$.

\begin{restatable}[First row of \Cref{table: Phase 5 Switches}]{lemma}{PhaseFiveCycleVertices} \label{lemma: Phase Five Cycle Vertices}
Let $\sigma\in\reach{\sigma_0}$ be a well-behaved phase-$5$-strategy for $\bit\in\bitset_n$.
Let $i\in[n], j=1-\indbit^{\sigma}_{i+1}, k\in\{0,1\}$ with $e\coloneqq(d_{i,j,k},F_{i,j})\in I_{\sigma}$ and assume $\sigma(b_i)=b_{i+1}, \sigmabar(g_i)=1-\indbit^{\sigma}_{i+1}$ and $i\neq 1$.
Then $\sigmae$ is a well-behaved Phase-$5$-strategy for $\bit$ with $\sigmae\in\reach{\sigma_0}$ and $I_{\sigmae}=I_{\sigma}\setminus\{e\}$.
\end{restatable}

The next lemma concludes our discussion on the application of the improving switches and the corresponding transition through the phases.
It describes the application of switches involving selector vertices during phase 5.

\begin{restatable}[Second row of of \Cref{table: Phase 5 Switches}]{lemma}{PhaseFiveSelectorVertices} \label{lemma: Phase Five Selector Vertices}
Let $\sigma\in\reach{\sigma_0}$ be a well-behaved phase-$5$-strategy for $\bit\in\bitset_n$.
Let $i\in[n], j\in\{0,1\}$ with $e\coloneqq (g_i,F_{i,j})\in I_{\sigma}$ and $\indbit^{\sigma}_i=0$.
Assume that $\nsb=1$ implies $\sigmabar(eg_{i,j})\wedge\nsigmabar(eb_{i,j})$ and that $\nsb>1$ implies $\sigmabar(eb_{i,j})\wedge\nsigmabar(eg_{i,j})$.
Further assume that $\relbit{\sigma}=1$ implies that for any $i'\geq i$ and $j'\in\{0,1\}$, either $\sigmabar(d_{i',j'})$ or $\sigmabar(eb_{i',j'})\wedge\nsigmabar(eg_{i',j'})$.
If $\sigma(e_{i',j',k'})=t^{\rightarrow}$ for all $i'\in[n],j',k'\in\{0,1\}$ and if $\sigmae$ has \Pref{SVG}$_{i'}$/(\ref{property: SVM})$_{i'}$ for all $i'\in[n]$, then $\sigmae$ is a phase-$1$-strategy for $\bit+1$.
Otherwise it is a phase-$5$-strategy for $\bit$.
In either case, $\sigmae$ is well-behaved and $I_{\sigmae}=I_{\sigma}\setminus\{e\}$.
\end{restatable}

This concludes our discussion of the lemmas describing the exact application of individual improving switches in $G_n$.
The next section now applies the results of this section to provide formal proofs of the statements of \Cref{section: Lower Bound Proof}.

\section{Proving the Main Statements} \label{section: Final Formal Proofs}

In this section, we provide the formal proofs for the statements given in \Cref{section: Lower Bound Proof}.
Before providing these proofs, we briefly explain how this section is organized.
We begin by considering a canonical strategy $\canstrat$ for some $\bit\in\bitset_n$ that has the canonical properties.
For each of the $k$ phases, we prove that applying improving switches according to Zadeh's rule yields a phase-$k$-strategy $\sigma^{(k)}$ as described by \Cref{table: Properties at start of phase,table: Switches at start of phase}.
This is done by considering the phases one after another.
At the end, we prove that applying the improving switches of phase $5$ to $\sigma^{(5)}$ yields a canonical strategy $\sigma_{\bit+1}$ for $\bit+1$.

When proving these statements, we typically immediately prove that the occurrence record of an edge $e$ is described correctly by \Cref{table: Occurrence Records} when interpreted for $\bit+1$ after its application.
The only kind of edges for which this is not proven immediately are edges of the form $(g_*,F_{*,*})$.
The reason is that we need to analyze more than a single transition to properly analyze the occurrence records of these edges.
Consequently, we defer the discussion of the occurrence records of these edges to the end of this section.

When proving the statements of this section, we often state smaller statements within the proofs as \textbf{Claims}.
Using claims allows us to hide several more technical aspects on the macroscopic level, making the important proofs shorter and thus easier to comprehend.
The proofs of all claims can however be found in \Cref{appendix: Proofs Exponential}.

Consider some fixed $\bit\in\bitset_n$ , let $\nsb\coloneqq\ell(\bit+1)$ denote the least significant set bit of $\bit+1$ and let $\canstrat\in\reach{\sigma_0}$ be a canonical strategy for $\bit$ that has the canonical conditions.
As a a reminder, for two strategies $\sigma,\sigma'$ with $\sigma'\in\reach{\sigma}$, the set $\applied{\sigma}{\sigma'}$ describes the set of improving switches applied by the strategy improvement resp. policy iteration algorithm during the transition $\sigma\to\sigma'$.
Also, we define the parameter $\maxocc\coloneqq\floor{(\bit+1)/2}$ as this quantity will often be used when analyzing the occurrence records as it serves as a natural upper bound and is the $\mathfrak{m}$aximum occurrence record that improving switches have.

\subsection{Reaching a phase-2-strategy}

We begin by proving that Zadeh's pivot rule together with the tie-breaking rule given in \Cref{definition: Tie-Breaking exponential} yields a phase-$2$-strategy $\sigma^{(2)}\in\reach{\sigma_0}$ as described by the corresponding rows of \Cref{table: Properties at start of phase,table: Switches at start of phase}.
That is, we  provide the proof of a slightly extended version of \Cref{lemma: Reaching phase 2}.
This extension states that $\sigma^{(2)}$ is well-behaved.
This was not included in the original statement as the term \enquote{well-behaved} was only introduced in the appendices.

\begin{lemma}[Extended version of \Cref{lemma: Reaching phase 2}] \label{lemma: Extended Reaching phase 2}
Let $\canstrat\in\reach{\sigma_0}$ be a canonical strategy for $\bit\in\bitset_n$ with $\nsb=\ell(\bit+1)>1$ having the canonical properties.
After applying finitely many improving switches, the strategy improvement algorithm produces a well-behaved phase-$2$-strategy $\sigma^{(2)}$ for $\bit$ as described by the corresponding rows of \Cref{table: Properties at start of phase,table: Switches at start of phase}.
\end{lemma}

\begin{proof}
By \Cref{lemma: Canonical strategies are well-behaved}, $\canstrat$ is well-behaved.
Let $j\coloneqq\bit_{\nsb+1}=(\bit+1)_{\nsb+1}$.
Since $\canstrat$ is a canonical strategy, we have $\canstrat(d_{\nsb,j,*})\neq F_{\nsb,j}$.
Moreover, $I_{\canstrat}=\mathfrak{D}^{\canstrat}$ and $\canstrat$ is a phase-1-strategy for~$\bit$ by \Cref{lemma: Improving sets of canonical strategies}.
In particular, $(d_{\nsb,j,*},F_{\nsb,j})\in I_{\canstrat}$.
By \Cref{lemma: Occurrence Records Cycle Vertices}, $(d_{\nsb,j,0},F_{\nsb,j})$ maximizes the occurrence record among all improving switches and $\occrec^{\canstrat}(d_{\nsb,j,1},F_{\nsb,j})=\maxocc-1$.
By \Pref{OR4}$_{*,*,*}$, $I_{\canstrat}$ can be partitioned into $I_{\canstrat}=I_{\canstrat}^{<\maxocc}\cup I_{\canstrat}^{\maxocc}$ where $e\in I_{\canstrat}^{<\maxocc}$ if $\occrec^{\canstrat}(e)=\maxocc-1$ and $e\in I_{\canstrat}^{\maxocc}$ if $\occrec^{\canstrat}(e)=\maxocc$.
If $I_{\canstrat}^{<\maxocc}\neq\emptyset$, then a switch contained in this set is applied first as the \textsc{LeastEntered} pivot rule always chooses an improving switch minimizing the occurrence record.
By applying \Cref{lemma: Phase 1 Low OR} iteratively, the algorithm applies switches $e\in I_{\canstrat}^{<\maxocc}$ until it either reaches a strategy $\sigma$ with $I^{<\maxocc}_{\sigma}=\emptyset$ or until an edge $(g_i,F_{i,j'})$ with $j'\neq j$ becomes improving.
By \Cref{lemma: Phase 1 Low OR},  \Cref{table: Occurrence Records} (interpreted for $\bit+1$) correctly describes the occurrence record of all switches applied in the process.

\begin{restatable}{claim}{SelectorVerticesInPhaseOne} \label{claim: Selector Vertices in Phase One}
If an edge $(g_i,F_{i,j'})$ with $i\in[n]$ and $j'\neq\bit_{\nsb+1}$ becomes improving during the application of improving switches contained in $I^{<\maxocc}$, then it is applied immediately.
Its application is described by row 4 of \Cref{table: Phase 1 Switches}.
\end{restatable}

Such an edge $(g_i,F_{i,j'})$ is only applied if $F_{i,j'}$ was closed by the previous applications.
This implies that either $(d_{i,j',0},F_{i,j'}), (d_{i,j',1}, F_{i,j'})\in I^{<\maxocc}_{\canstrat}$ or $\sigma(d_{i,j',1-k})=\canstrat(d_{i,j',1-k})=F_{i,j'}$.
The first case can only happen for $i=\nsb$ and $j'=1-\bit_{\nsb+1}$.
Thus, if a switch $(g_i,F_{i,j'})$ is applied, then either $i=\nsb$ or $\canstrat(d_{i,j',1-k})=F_{i,j'}$ and $(d_{i,j',k},F_{i,j'})\in I^{<\maxocc}_{\canstrat}$.

The previous arguments can now be applied until a strategy is reached for which no edge has a \enquote{low} occurrence record.
Thus let $\sigma$ be a phase-1-strategy $\sigma$ with $I^{<\maxocc}_{\sigma}=\emptyset$ and $\mathbb{G}\cap I_{\sigma}=\emptyset$.
Further note that $\applied{\canstrat}{\sigma}\subseteq\mathbb{D}^1\cup\mathbb{G}$ and that $(g_i,F_{i,j})\in\applied{\canstrat}{\sigma}$ implies $\bit_i=0\wedge\bit_{i+1}\neq j$ and that the previous arguments hold independent of $\nsb$.

We discuss improving switches contained in $I_{\canstrat}^{\maxocc}$ next.

\begin{restatable}{claim}{NoCycleCenterOpenPhaseOne} \label{claim: No Cycle Center Open Phase One}
Let $\nsb>1$ and let $\sigma$ denote the strategy obtained after applying all improving switches contained in $I_{\sigma}^{<\maxocc}$.
For all suitable indices $i\in[n], j'\in\{0,1\}$  it holds that $\sigma(d_{i,j',1})=F_{i,j'}$, implying that no cycle center is open for $\sigma$.
\end{restatable}

Let $\sigma\in\reach{\canstrat}$ be a phase-1-strategy with $I_{\sigma}=\{e=(d_{i,j,k},F_{i,j})\colon\occrec^{\sigma}(e)=\maxocc\}=\mathfrak{D}^{\sigma}$ as described previously.
We prove that $e\coloneqq(d_{\nsb,j,0},F_{\nsb,j})$ is applied next.
By \Cref{lemma: Occurrence Records Cycle Vertices}, the definition of a canonical strategy and since only edges with an occurrence record less than $\maxocc$ were applied so far, this implies $e\in I_{\sigma}$.
Since all improving switches have the same occurrence records, it is sufficient to show that no other improving switch is ranked lower by the tie-breaking rule.
By \Cref{claim: No Cycle Center Open Phase One}, there are no open cycle centers.
Hence, the ordering of the edges is based on the bits represented by the levels, the index of the levels and whether the cycle center is active.
To be precise, the first switch according to the tie-breaking rule is the improving switch contained in the active cycle center of the lowest level with a bit equal to 0.
This edge is precisely  $e=(d_{\nsb,j,0},F_{\nsb,j})$.

We now prove that the occurrence record of $e$ is described by \Cref{table: Occurrence Records} when interpreted for $\bit+1$ after the application.
Since~$F_{\nsb,j}$ is closed for $\sigma$ but was open for $\canstrat$, we prove \[\occrec^{\sigmae}(e)=\ceil{\frac{\lastflip{\bit+1}{\nsb}{\{(\nsb+1,j)\}}+1}{2}}.\]
By the definition of $\nsb$, it holds that $\bit+1=\lastflip{\bit+1}{\nsb}{\{(\nsb+1,j)\}}$.
The statement thus follows since $\maxocc+1=\ceil{(\bit+1+1)/2}$.

By row 5 of \Cref{table: Phase 1 Switches}, $\sigmae$ is a well-behaved (potentially pseudo) phase-2-strategy for~$\bit$.
If $\sigma(g_\nsb)\neq F_{\nsb,j}$, then $(g_{\nsb},F_{\nsb,j})$ minimizes the occurrence record among all improving switches.
Due to the tie-breaking rule, this switch is then applied next, and this application is formalized in row 6 of \Cref{table: Phase 1 Switches}.

Let $\sigma$ denote the strategy obtained after applying  $(g_{\nsb},F_{\nsb,j})$ if $\sigma(g_{\nsb})\neq F_{\nsb,j}$ resp. after applying $(d_{\nsb,j,0},F_{\nsb,j})$ if $\sigma(g_{\nsb})=F_{\nsb,j}$.
Then, by row~5 resp.~6 of \Cref{table: Phase 1 Switches}, \[I_{\sigma}=\mathfrak{D}^{\sigma}\cup\{(b_{\nsb},g_{\nsb})\}\cup\{(s_{\nsb-1,1},h_{\nsb-1,1})\}.\]
Furthermore, $\sigma$ has \Pref{USV3}$_i$ for all $i<\nsb$ as $\canstrat$ has \Pref{USV1}$_i$ and $\bit_i=1-(\bit+1)_i$ for $i\leq\nsb$. 
In addition, $\sigma(d_{i,j,k})\neq F_{i,j}$ implies $\occrec^{\sigma^{(2)}}(d_{i,j,k},F_{i,j})=\maxocc$ by \Cref{corollary: Switches With Low OR In Phase One}.
Moreover, since no improving switch $(d_{*,*,*},e_{*,*,*})$ was applied and $\bit_{i}=1-\indbit^{\sigma}_{i+1}$ for all $i<\nsb$, it holds that $\sigmabar(g_i)=1-\indbit^{\sigma}_{i+1}$ and $\sigmabar(d_{i,1-\indbit^{\sigma}_{i+1}})$ for all $i<\nsb$.
\end{proof}

We henceforth refer to the phase-2-strategy that is described by the corresponding rows of \Cref{table: Properties at start of phase,table: Switches at start of phase} and whose existence we just proved by $\sigma^{(2)}$.
When proving the existence of this strategy, we furthermore implicitly proved the following three corollaries.
We later show that the condition~$\nsb>1$ can be dropped in the first corollary.

\begin{corollary} \label{corollary: Selection Vertices In Phase One}
Let $\canstrat$ be a canonical strategy for $\bit$ having the canonical properties and $\nsb>1$.
Let $i\in[n]$ and $j\in\{0,1\}$.
Then, the edge $(g_i,F_{i,j})$ is applied as improving switch during phase $1$ if and only if $F_{i,j}$ is closed during phase~$1$, $\canstrat(g_i)\neq F_{i,j}$ and $i\neq\nsb$.
A cycle center can only be closed during phase~$1$ if either $i=\nsb$ or if there exists an index $k\in\{0,1\}$ such that  $\canstrat(d_{i,1-\bit_{i+1},k})=F_{i,1-\bit_{i+1}}, \occrec^{\canstrat}(d_{i,1-\bit_{i+1},1-k},F_{i,\bit_{i+1}})<\maxocc$ and $\canstrat(b_i)=b_{i+1}$.
\end{corollary}

\begin{corollary} \label{corollary: Switches With Low OR In Phase One}
Let $\canstrat$ be a canonical strategy for $\bit$ having the canonical properties and let $i\in[n],j,k\in\{0,1\}$ such that $\canstrat(d_{i,j,k})\neq F_{i,j}$.
If $\occrec^{\canstrat}(d_{i,j,k},F_{i,j})<\maxocc$, then $(d_{i,j,k},F_{i,j})$ is applied during phase 1.
\end{corollary}

\begin{corollary} \label{corollary: No Open CC For Phase Two}
No cycle center is open with respect to~$\sigma^{(2)}$.
\end{corollary}

\begin{corollary} \label{corollary: Improving Switches of Phase 1}
\Cref{table: Occurrence Records} correctly specifies the occurrence record of every improving switch applied during $\canstrat\to\sigma^{(2)}$ when interpreted for $\bit+1$, excluding switches $(g_*,F_{*,*})$.
In addition, each switch is applied at most once.
\end{corollary}

\subsection{Reaching a phase-3-strategy}

We now prove that the algorithm produces a phase-$3$-strategy by proving a slightly extended version of \Cref{lemma: Reaching phase 3}.
If $\nsb=1$, then this follows by analyzing phase 1 in a similar fashion as done when proving \Cref{lemma: Extended Reaching phase 2}.
It in fact turns out that nearly the identical arguments can be applied.
If $\nsb>1$, then we use that lemma to argue that we obtain a phase-2-strategy.
We then investigate phase 2 in detail and prove that we also obtain a phase-3-strategy.
The proof uses the statements summarized in \Cref{table: Phase 1 Switches,table: Phase 2 Switches}, and we refer to these tables and the corresponding statements in proofs of the related statements.

\begin{lemma}[Extended version of \Cref{lemma: Reaching phase 3}] \label{lemma: Extended Reaching phase 3}
Let $\canstrat\in\reach{\sigma_0}$ be a canonical strategy for $\bit\in\bitset_n$ having the canonical properties.
After applying a finite number of improving switches, the strategy improvement algorithm produces a well-behaved phase-3-strategy $\sigma^{(3)}\in\reach{\sigma_0}$ as described by the corresponding rows of \Cref{table: Properties at start of phase,table: Switches at start of phase}.
\end{lemma}

\begin{proof}
Consider the case $\nsb=1$ first.
As shown in the proof of \Cref{lemma: Extended Reaching phase 2}, the set $I_{\canstrat}$ can be partitioned into~$I_{\canstrat}^{<\maxocc}$ and~$I_{\canstrat}^{\maxocc}$.
Since \Cref{lemma: Phase 1 Low OR} also applies for $\nsb=1$, the same arguments imply that the algorithms calculate a phase-1-strategy $\sigma\in\reach{\sigma_0}$ with $I_{\sigma}=\{e=(d_{i,j,k},F_{i,j})\colon\occrec^{\sigma}(e)=\maxocc\}=\mathfrak{D}^{\sigma}$.
We can again deduce $\applied{\canstrat}{\sigma}\subseteq\mathbb{D}^1\cup\mathbb{G}$ and that $(g_i,F_{i,j})\in\applied{\canstrat}{\sigma}$ implies $\bit_i=0\wedge\bit_{i+1}\neq j$ or $i=\nsb$ for all $i\in[n], j\in\{0,1\}$.
We can further assume $(g_*,F_{*,*})\notin I_{\sigma}$.
Also, by \Cref{lemma: Phase 1 Low OR}, the occurrence records of edges $(d_{*,*,*},F_{*,*})\in\applied{\canstrat}{\sigma}$ is described by \Cref{table: Occurrence Records} when interpreted for $\bit+1$.

Since all improving switches now have the same occurrence records, their order of application depends on the tie-breaking rule.
Due to the first criterion, improving switches contained in open cycle centers are applied first.
Hence, a sequence of strategies is produced until a strategy without open cycle centers is reached.
All produced strategies are well-behaved phase-$1$-strategies for~$\bit$, reachable from $\sigma_0$ by row $1$ of \Cref{table: Phase 1 Switches}.
Also, by the tie-breaking rule, the edge $(d_{*,*,0},F_{*,*})$ is applied as improving switch in an open cycle center $F_{*,*}$.
By the same arguments used when proving \Cref{lemma: Extended Reaching phase 2}, the second switch of~$F_{\nu,\bit_{\nsb+1}}$ is applied next and, possibly, $(g_{\nsb},F_{\nsb,\bit_{\nsb+1}})$ is applied afterwards.

Let $\sigma^{(3)}$ denote the strategy obtained after closing the cycle center $F_{\nsb,\bit_{\nsb+1}}$ resp. after applying $(g_{\nsb},F_{\nsb,\bit_{\nsb+1}})$ if it becomes improving.

\begin{restatable}{claim}{OcccurrenceRecordOfCycleVerticesPhaseOne} \label{claim: Occurrence Records Of Cycle Vertices Phase One}
Let $i\in[n],j,k\in\{0,1\}$ such that $(d_{i,j,k},F_{i,j})\in\applied{\canstrat}{\sigma^{(3)}}$.
The occurrence records of $(d_{i,j,k},F_{i,j})$ with respect to $\sigma^{(3)}$ is specified by \Cref{table: Occurrence Records} when interpreted for $\bit+1$.
\end{restatable}

Note that the last row of \Cref{table: Phase 2 Switches} can be used to describe the application of $(g_{\nsb},F_{\nsb,\bit_{\nsb+1}})$.
Then, by row~5 of \Cref{table: Phase 1 Switches} resp. the last row of \Cref{table: Phase 2 Switches} and our previous arguments, \smash{$\sigma^{(3)}$} has all properties listed in the respective rows of \Cref{table: Properties at start of phase,table: Switches at start of phase}.
Furthermore, as we used the same arguments, \Cref{corollary: Selection Vertices In Phase One} is also valid for $\nsb=1$ and we can drop the assumption $\nsb>1$.

Consider the case $\nsb>1$, implying $\bit\geq 1$.
By \Cref{lemma: Extended Reaching phase 2}, applying improving switches to $\canstrat$ yields a phase-2-strategy $\sigma=\sigma^{(2)}$ for $\bit$ with $I_{\sigma}=\mathfrak{D}^{\sigma}\cup\{(b_\nsb,g_\nsb),(s_{\nsb-1,1},h_{\nsb-1,1})\}$ and $\sigma\in\reach{\sigma_0}$.
By \Cref{table: Occurrence Records} and \Cref{lemma: Numerics Of OR},  \[\occrec^{\sigma}(b_\nsb,g_\nsb)=\flips{\bit}{\nsb}{}=\occrec^{\sigma}(s_{\nsb-1,1},h_{\nsb-1,1})\quad\text{and}\quad\flips{\bit}{\nsb}{}=\floor{\frac{\bit+2^{\nsb-1}}{2^\nsb}}.\]
Since $\nsb>1$ and $\bit\geq 1$, this implies $\flips{\bit}{\nsb}{}\leq\lfloor(\bit+2)/4\rfloor\leq\maxocc$.
By \Cref{lemma: Extended Reaching phase 2}, any improving switch $(d_{*,*,*},F_{*,*})\in I_{\sigma}$ has an occurrence record of $\maxocc$.
Thus, by the tie-breaking rule, $(b_\nsb,g_\nsb)$ is applied next.
Let $\sigmae$ denote the strategy obtained after applying $(b_\nsb,g_\nsb)$.
It is easy to verify that $\sigma$ has the properties of row 1 of \Cref{table: Phase 2 Switches}.
Consequently,~$\sigmae$ is a phase-2-strategy for $\bit$ with $\sigmae\in\reach{\sigma_0}$.
By \Cref{lemma: Numerics Of OR}, $\occrec^{\sigmae}(b_{\nsb},g_{\nsb})=\flips{\bit}{\nsb}{}+1=\flips{\bit+1}{\nsb}{},$ so \Cref{table: Occurrence Records} describes the occurrence record of $(b_{\nsb},g_{\nsb})$ with respect to $\bit+1$.
The set of improving switches for $\sigmae$ now depends on $\nsb$, see row 1 of \Cref{table: Phase 2 Switches}.

Let $\nsb=2$.
Then $I_{\sigmae}=\mathfrak{D}^{\sigmae}\cup\{(b_1,b_2),(s_{1,1},h_{1,1})\}\cup\{(e_{*,*,*},b_2)\}$.
In this case, $(s_{1,1},h_{1,1})$ is applied next and its application yields the desired phase-$3$-strategy $\sigma^{(3)}$.

\begin{restatable}{claim}{USVPhaseTwoIfNSBTwo} \label{claim: USV Phase Two If NSB Two}		
Let $\nsb=2$ and consider the phase-2-strategy $\sigma$ obtained after the application of $(b_{\nsb}, g_{\nsb})$.
Then, the edge $(s_{1,1},h_{1,1})$ is applied next, and the obtained strategy is a well-behaved phase-$3$-strategy for $\bit$ described by the respective rows of \Cref{table: Properties at start of phase,table: Switches at start of phase}.
\end{restatable}

If $\nsb>1$, then we do not obtain the desired strategy yet and we have to consider a longer sequence of improving switches that are applied.
Thus, let $\nsb>2$, implying $\bit\neq 1$.
Then, the first row of \Cref{table: Phase 2 Switches} implies \[I_{\sigmae}=\mathfrak{D}^{\sigmae}\cup\{(b_{\nsb-1},b_{\nsb}),(s_{\nsb-1,1},h_{\nsb-1,1}),(s_{\nsb-2,0},h_{\nsb-2,0})\}.\]
By \Cref{table: Occurrence Records}, $\occrec^{\sigmae}(b_{\nsb-1},b_{\nsb})=\flips{\bit}{\nsb-1}{}-1$ and $\occrec^{\sigmae}(s_{\nsb-1,1},h_{\nsb-1,1})=\flips{\bit}{\nsb}{}$.
In addition, it holds that  $\occrec^{\sigmae}(s_{\nsb-2,0},h_{\nsb-2,0})=\flips{\bit}{\nsb-1}{}-1.$
Hence, both edges $(b_{\nsb-1},b_{\nsb})$ and $(s_{\nsb-2,0},h_{\nsb-2,0})$ minimize the occurrence record.
By the tie-breaking rule, the switch $e'\coloneqq(b_{\nsb-1},b_{\nsb})$ is now applied.
We show that the application of $e'$ can be described by row 3 of \Cref{table: Phase 2 Switches}.
We thus need to show the following:
\begin{itemize}
	\item \boldall{$\sigmaebar(d_{i'})$ for all $i'<\relbit{\sigmae}$:} This follows from \Cref{lemma: Extended Reaching phase 2} as no switch $(d_{*,*,*},e_{*,*,*})$ was applied during $\canstrat\to\sigma^{(2)}$ and no improving switch involving selector vertices was applied in a level $i'<\relbit{\sigmae}$.
	\item \boldall{$\sigmae$ has \Pref{USV3}$_{i'}$ for all $i'<\nsb-1$:} 
		Since no switch $(s_{i',*},*)$ was applied for $i'<\nsb-1$, this follows since $\canstrat$ has \Pref{USV1}$_{i'}$ for those indices.				
	\item \boldall{$\sigmae$ has \Pref{EV1}$_{i'}$ and (\ref{property: EV2})$_{i'}$ for all $i'>\nsb-1$ and (\ref{property: EV3})$_{i'}$ for all $i'>\nsb-1$ with $i'\neq\relbit{\sigmae}$:}
		Since $\relbit{\sigmae}-1=\nsb-1$ and $\sigmae$ is a phase-2-strategy for $\bit$, it suffices to prove that $\sigmae$ has \Pref{EV1}$_{\nsb}$ and \Pref{EV2}$_{\nsb}$.
		This however follows since the strategy in which $(b_{\nsb},g_{\nsb})$ was applied had \Pref{CC2}. 
\end{itemize}
		
For simplicity, we denote the strategy that is obtained by applying $e'$ to $\sigmae$ by $\sigma$.
By our previous arguments and row 3 of \Cref{table: Phase 2 Switches}, $\sigma$ is a well-behaved phase-2-strategy for $\bit$ that has \Pref{CC2} as well as Properties (\ref{property: EV1})$_i$ and (\ref{property: EV2})$_i$ for all $i\geq\nsb-1$ and \Pref{EV3}$_i$ for all $i>\nsb-1, i\neq\nsb$.
In addition, $\sigmabar(d_{i})$ for all $i<\nsb$ and $\sigma$ has \Pref{USV3}$_i$ for all $i<\nsb-1$ 
Furthermore, \Cref{lemma: Numerics Of OR} implies \[\occrec^{\sigma}(e)=\flips{\bit}{\nsb-1}{}-1+1=\flips{\bit}{\nsb-1}{}=\flips{\bit+1}{\nsb-1}{}-(\bit+1)_{\nsb-1},\] so \Cref{table: Occurrence Records} describes the occurrence record of $e$ with respect to $\bit+1$.		
By row 3 of \Cref{table: Phase 2 Switches}, $\nsb-1>2$ implies \[I_{\sigma}=\mathfrak{D}^{\sigma}\cup\{(s_{\nu-1,1},h_{\nu-1,1}),(s_{\nu-2,0},h_{\nu-2,0}),(b_{\nu-2},b_{\nu-1}),(s_{\nu-3,0},h_{\nu-3,0})\}.\]
Similarly, $\nsb-1=2$ implies \[I_{\sigma}=\mathfrak{D}^{\sigma}\cup\{(e_{i,j,k},b_2)\}\cup\{(b_1,b_2),(s_{2,1},h_{2,1}),(s_{1,0},h_{1,0})\}.\]
In both cases, $e\coloneqq(s_{\nsb-1,1},h_{\nsb-1,1})\in I_{\sigma}$ is applied next.

\begin{restatable}{claim}{SecondSelectionSwitchInPhaseTwo} \label{claim: Second Selection Switch in Phase Two}
After the application of $(b_{\nsb-1},b_{\nsb})$ in the case $\nsb>2$, the switch $e=(s_{\nsb-1,1},h_{\nsb-1,1})$ is applied next.
Its application can be described by row 2 of \Cref{table: Phase 2 Switches} and \Cref{table: Occurrence Records} specifies its occurrence record after the application correctly when interpreted for $\bit+1$.
\end{restatable}
		
Let $\nsb-1>2$.
We argue that applying improving switches according to Zadeh's pivot rule and our tie-breaking rule then results in a sequence of strategies such that we finally obtain a strategy $\sigma'$ with $I_{\sigma'}=\mathfrak{D}^{\sigma'}\cup\{(e_{*,*,*},b_2)\}\cup\{(b_1,b_2),(s_{1,0},h_{1,0})\}.$
Note that such a strategy is also obtained after the application of $(s_{\nsb-1,1},h_{\nsb-1,1})$ if $\nsb-1=2$.
For any $x\in\{2,\dots,\nsb-2\}$, \Cref{lemma: Numerics Of OR} implies \begin{equation} \label{equation: Phase 2 General}
	\begin{split}			
	\occrec^{\sigmae}(s_{\nsb-x,0},h_{\nsb-x,0})&<\occrec^{\sigmae}(b_{\nsb-x},b_{\nsb-(x+1)})=\occrec^{\sigmae}(s_{\nsb-(x-1),0},h_{\nsb-(x-1),0})\\
		&<\occrec^{\sigmae}(e_{*,*,*},b_2).
	\end{split}			
	\end{equation}
Thus, $(s_{\nsb-2,0},h_{\nsb-2,0})$ is applied next.
It is easy to verify that $\sigmae$ meets the requirements of row 2 of \Cref{table: Phase 2 Switches}, so it can be used to describe the application of $(s_{\nsb-2,0},h_{\nsb-2,0})$.
		
Let $\sigma'$ denote the strategy obtained.
Then $I_{\sigma'}=\mathfrak{D}^{\sigma}\cup\{(b_{\nsb-2},b_{\nsb-1}),(s_{\nsb-3,0},h_{\nsb-3,0})\}$.
Also, the occurrence record of $(s_{\nsb-2,0},h_{\nsb-2,0})$ is described by \Cref{table: Occurrence Records} when interpreting the table for $\bit+1$.
By Equation~(\ref{equation: Phase 2 General}) and the tie-breaking rule, $(b_{\nsb-2},b_{\nsb-1})$ is applied next.
Similar to the previous cases, it is easy to check that row 3 of \Cref{table: Phase 2 Switches} applies to this switch.
We thus obtain a strategy $\sigma$ such $\nsb-2\neq 2$ implies \[I_{\sigma}=\mathfrak{D}^{\sigma}\cup\{(s_{\nsb-3,0},h_{\nsb-3,0}),(b_{\nsb-3},b_{\nsb-2}),(s_{\nsb-4,0},h_{\nsb-4,0})\}\] and $\nsb-2=2$ implies \[I_{\sigma}=\mathfrak{D}^{\sigma}\cup\{(e_{*,*,*,},b_2)\}\cup\{(b_1,b_2),(s_{1,0},h_{1,0})\}.\]
In either case, a simple calculation implies that the occurrence record of $(b_{\nsb-2},b_{\nsb-1})$ is described by \Cref{table: Occurrence Records} interpreted for $\bit+1$.		

In the first case, we can now apply the same arguments again iteratively as Equation~(\ref{equation: Phase 2 General}) remains valid for $\sigma'$ and $x\in\{2,\dots,\nsb-3\}$.
After applying a finite number of improving switches we thus obtain a phase-$2$-strategy $\sigma\in\reach{\sigma_0}$ with \[I_{\sigma}=\mathfrak{D}^{\sigma}\cup\{(e_{*,*,*},b_2)\}\cup\{(b_1,b_2),(s_{1,0},h_{1,0})\}.\]
Furthermore,  $\sigma$ has Properties (\ref{property: EV1})$_i$, (\ref{property: EV2})$_i$ and (\ref{property: USV2})$_{i,\indbit_{i+1}}$ for all $i>1$ as well as \Pref{EV3}$_i$ for all $i>1, i\neq\relbit{\sigma}$ and \Pref{CC2}.		
In addition, $\sigmabar(g_i)=1-\indbit_{i+1}$ and $\sigmabar(d_{i,1-\indbit_{i+1}})$ for all $i<\nsb$ and the occurrence records of all edges applied so far (with the exception of switches $(g_*,F_{*,*})$) is described by \Cref{table: Occurrence Records} when being interpreted for $\bit+1$.
Note that all of this also holds if $\nsb-1=2$.		
		
Consequently, $\sigma$ meets the requirements of row 2 of \Cref{table: Phase 2 Switches}.
As $\nsb>2$, we have $\indbit_2=0$.
By \Cref{table: Occurrence Records},  \[\occrec^{\sigma}(s_{1,0},h_{1,0})=\flips{\bit}{2}{}-1<\flips{\bit}{1}{}-1=\occrec^{\sigma}(b_1,b_2)\] as well as \[\flips{\bit}{2}{}-1=\floor{(\bit+2)/4}-1<\floor{\bit/2}=\occrec^{\sigma}(e_{*,*,*},b_2).\]
Hence, the switch $e=(s_{1,0},h_{1,0})$ is applied next and by row 2 of \Cref{table: Phase 2 Switches}, $\sigma^{(3)}\coloneqq\sigmae$ is a phase-$3$-strategy for~$\bit$ with \[I_{\sigma^{(3)}}=\mathfrak{D}^{\sigma^{(3)}}\cup\{(e_{i,j,k},b_2)\}\cup\{(b_1,b_2)\}.\]
We thus obtain a strategy as described by the corresponding rows of \Cref{table: Properties at start of phase,table: Switches at start of phase}.
\end{proof}

We henceforth use $\sigma^{(3)}$ to refer to the phase-3-strategy described by \Cref{lemma: Extended Reaching phase 3}.
Note that we implicitly proved the following corollaries where the second follows by \Cref{corollary: Improving Switches of Phase 1}.

\begin{corollary} \label{corollary: No Open CC In Phase Three}
No cycle center is open with respect to $\sigma^{(3)}$.
\end{corollary}

\begin{corollary} \label{corollary: Improving Switches of Phase 1 and 2}
\Cref{table: Occurrence Records} specifies the occurrence record of every improving switch applied during $\canstrat\to\sigma^{(3)}$ when interpreted for $\bit+1$, excluding switches $(g_*,F_{*,*})$.
In addition, each such switch was applied once.
\end{corollary}

\subsection{Reaching a phase-4-strategy or a phase-5-strategy}

We now discuss the application of improving switches during phase 3, which highly depends on whether $G_n=M_n$ or $G_n=S_n$ and on the least significant set bit of $\bit+1$.
The next lemma now summarizes the application of improving switches during phase~3 and is a generalization of \Cref{lemma: Reaching phase 4 or 5}.
Depending on $G_n$ and $\nsb$, we then either obtain a phase-4-strategy or a phase-5-strategy for $\bit$.
As with the previous lemmas, this lemma is an extension of \Cref{lemma: Reaching phase 4 or 5}.
We also use the usual notation and define $t^{\rightarrow}\coloneqq b_2$ if $\nsb>1$ and $t^{\rightarrow}\coloneqq g_1$ if $\nsb=1$.
Similarly, let $t^{\leftarrow}\coloneqq g_1$ if $\nsb>1$ and $t^{\leftarrow}\coloneqq b_2$ if $\nsb=1$.

\begin{lemma}[Extended version of \Cref{lemma: Reaching phase 4 or 5}] \label{lemma: Extended Reaching phase 4 or 5}
Let $\canstrat\in\reach{\sigma_0}$ be a canonical strategy for $\bit\in\bitset_n$ having the canonical properties.
After applying finitely many improving switches, the strategy improvement algorithm produces a well-behaved strategy $\sigma$ with the following properties:
If $\nsb>1$, then $\sigma$ is a phase-$k$-strategy for $\bit$, where $k=4$ if $G_n=S_n$ and $k=5$ if $G_n=M_n$.
If $\nsb=1$, then $\sigma$ is a phase-5-strategy for~$\bit$.
In any case, $\sigma\in\reach{\sigma_0}$ and $\sigma$ is described by the corresponding rows of \Cref{table: Properties at start of phase,table: Switches at start of phase}.
\end{lemma}

Before proving this lemma, we provide an additional lemma that summarizes the application of switches of the type $(d_{*,*,*}, e_{*,*,*})$.
Its proof is omitted here and deferred to \Cref{appendix: Proofs Exponential}.

\begin{restatable}{lemma}{EasyCVSwitchesPhaseThree} \label{lemma: Easy CV Switches Phase Three}
Let $\sigma\in\reach{\sigma^{(3)}}$ be a well-behaved phase-3-strategy for $\bit$ obtained through the application of a sequence $\applied{\sigma^{(3)}}{\sigma}\subseteq\mathbb{E}^1\cup\mathbb{D}^0$ of improving switches.
Assume that the conditions of row~1 of \Cref{table: Phase 3 Switches} were fulfilled for each intermediate strategy $\sigma'$ of the transition $\sigma^{(3)}\to\sigma$.
Let $t^{\rightarrow}\coloneqq b_2$ if $\nsb>1$ and $t^{\rightarrow}\coloneqq g_1$ if $\nsb=1$.
Let $i\in[n],j,k\in\{0,1\}$ such that $e\coloneqq(d_{i,j,k}, e_{i,j,k})\in I_{\sigma}$ is applied next  and assume $\sigma(e_{i,j,k})=t^{\rightarrow}, \indbit^{\sigma}_{i}=0\vee\indbit^{\sigma}_{i+1}\neq j$ and $I_{\sigma}\cap\mathbb{D}^0=\{e\}$.
Further assume that either $i\geq\nsb$ or that we consider the case $G_n=S_n$.
Then $\sigmae$ is a phase-3-strategy for $\bit$ with $I_{\sigmae}=(I_{\sigma}\setminus\{e\})$.
\end{restatable}

This now enables us to prove \Cref{lemma: Extended Reaching phase 4 or 5}.

\begin{proof}[Proof of \Cref{lemma: Extended Reaching phase 4 or 5}]
By \Cref{lemma: Extended Reaching phase 3}, applying improving switches according to Zadeh's pivot rule and our tie-breaking rule yields a phase-3-strategy $\sigma^{(3)}\in\reach{\sigma_0}$ described by the corresponding rows of \Cref{table: Properties at start of phase,table: Switches at start of phase}.
As it simplifies the formal proof significantly, we begin by describing phase~3 informally.

For every cycle vertex $d_{*,*,*}$, it either holds that $\sigma^{(3)}(d_{*,*,*})=F_{*,*}$ or $\sigma^{(3)}(d_{*,*,*})=e_{*,*,*}$ and $(d_{*,*,*},F_{*,*})\in I_{\sigma^{(3)}}$.
It will turn out that only switches corresponding to cycle vertices of the first type are applied during phase 3.
Consider an arbitrary but fixed such cycle vertex $d_{i,j,k}$ for some suitable indices $i,j,k$.
Then, the switch $(e_{i,j,k},t^{\rightarrow})$ will be applied.
If $(\bit+1)_i=0$ or $(\bit+1)_{i+1}\neq j$, then $(d_{i,j,k},e_{i,j,k})$ becomes improving and is applied next.
This procedure then continues until all such improving switches have been applied.
During this procedure, it might happen that an edge $(s_{i',*},b_1)$ with $i'<\nsb$ becomes improving after applying some switch $(d_{*,*,*},e_{*,*,*})$ if $\nsb>1$ and $G_n=M_n$.
In this case, the corresponding switch is applied immediately.
Finally, $(b_1,b_2)$ resp. $(b_1,g_1)$ is applied, resulting in a phase-4-strategy if $\nsb>1$ and $G_n=S_n$ and in a phase-5-strategy otherwise.

We now formalize this behavior.
We first show that switches $(e_{*,*,*}, t^{\rightarrow})$ minimize the occurrence record among all improving switches.
Consider some indices $i,j,k$ such that $(d_{i,j,k},F_{i,j})\in I_{\sigma^{(3)}}$.
Then, $\occrec^{\sigma^{(3)}}(d_{i,j,k},F_{i,j})=\maxocc$ by \Cref{lemma: Extended Reaching phase 3} resp. \Cref{table: Properties at start of phase}.
If $\nsb>1$, then \[\occrec^{\sigma^{(3)}}(e_{i,j,k},b_2)=\floor{\frac{\bit}{2}}=\maxocc-1=\flips{\bit}{1}{}-\bit_1=\occrec^{\sigma^{(3)}}(b_1,b_2)\] by \Cref{table: Occurrence Records}.
Similarly, if $\nsb=1$, then $\occrec^{\sigma^{(3)}}(e_{i,j,k},g_1)=\occrec^{\sigma^{(3)}}(b_1,g_1).$
By the tie-breaking rule, a switch of the type $(e_{i',j',k'},t^{\rightarrow})$ with $\sigma^{(3)}(d_{i',j',k'})=F_{i',j'}$ for some suitable indices is thus applied next.
Since $\sigma^{(3)}(s_{i',*})=h_{i',*}$ for all $i'<\relbit{\sigma^{(3)}}$ by \Cref{lemma: Extended Reaching phase 3}, the statement of row 1 of \Cref{table: Phase 3 Switches} can be applied.

Let $i\in[n],j,k\in\{0,1\}$ denote the indices such that $e\coloneqq(e_{i,j,k},t^{\rightarrow})\in I_{\sigma}\cap\mathbb{E}^1$ is the switch that is applied next.
We prove that the characterization given in the first row of \Cref{table: Phase 3 Switches} implies $I_{\sigmae}=(I_{\sigma}\setminus\{e\})\cup\{(d_{i,j,k},e_{i,j,k})\}$ if $\indbit_{i}=0\vee\indbit_{i+1}\neq j$ and $I_{\sigmae}=I_{\sigma}\setminus\{e\}$ else.
As explained earlier, the strategy $\sigma$ fulfills the requirements of the first row of \Cref{table: Phase 3 Switches}.
Consider the strategy $\sigmae$.
By the first row of \Cref{table: Phase 3 Switches}, $(d_{i,j,k},e_{i,j,k})$ is improving for $\sigmae$ if and only if either $\sigma(d_{i,j,1-k})=e_{i,j,1-k}$ or $[\sigma(d_{i,j,1-k})=F_{i,j}$ and $j\neq\indbit_{i+1}]$.
It thus suffices to prove that $\indbit_i=0\vee\indbit_{i+1}\neq j$ is equivalent to the disjunction of these two conditions.
We do so by showing $\indbit_i=1\wedge\indbit_{i+1}=j\Leftrightarrow\sigma(d_{i,j,1-k})=F_{i,j}\wedge\indbit_{i+1}=j$. 
The direction ``$\Rightarrow$'' follows since the cycle center $F_{i,j}$ is then active and closed. 
The direction ``$\Leftarrow$'' follows since $\sigma(d_{i,j,1-k})=F_{i,j}$ implies that $F_{i,j}$ is closed as $e$ being improving for $\sigma$ implies $\sigma(d_{i,j,k})=F_{i,j}$.
But then, by the definition of $\indbit$ and the choice of $j$, $\indbit_i=1$.

Consequently, by the tie-breaking rule and row 1 of \Cref{table: Phase 3 Switches}, improving switches $(e_{*,*,*},t^{\rightarrow})\in \mathbb{E}^1$ are applied until a switch of this type with $i\in[n], j\in\{0,1\}$ and $\indbit_i=0\vee\indbit_{i+1}\neq j$ is applied.
The occurrence record of each applied switch is described by \Cref{table: Occurrence Records} when interpreted for $\bit+1$ since $\floor{\bit/2}+1=\maxocc$ if $\bit$ is odd and $\ceil{\bit/2}+1=\ceil{(\bit+1)/2}$ if $\bit$ is even.
By row~$1$ of \Cref{table: Phase 3 Switches}, $(d_{i,j,k},e_{i,j,k})$ now becomes improving.
As $(d_{i,j,k},e_{i,j,k})\notin\applied{\canstrat}{\sigma}$ and since switches of the type $(e_{*,*,*},t^{\rightarrow})$ minimize the occurrence record, \Cref{table: Occurrence Records} and the tie-breaking rule imply that $(d_{i,j,k},e_{i,j,k})$ is applied next.
In particular, an edge $(d_{i,j,k},e_{i,j,k})$ is applied immediately if it becomes improving and this requires that $(e_{i,j,k},t^{\rightarrow})$ was applied earlier.
Therefore, the application of improving switches $(e_{*,*,*},t^{\rightarrow})$ is described by row~1 of \Cref{table: Phase 3 Switches} and whenever an edge $(d_{*,*,*},e_{*,*,*})$ becomes improving, its application is described by \Cref{lemma: Easy CV Switches Phase Three}.
In particular, the occurrence record of all these edges is described by \Cref{table: Occurrence Records} when interpreted for $\bit+1$.

Let $G_n=S_n$.
Then, row 1 of \Cref{table: Phase 3 Switches} and \Cref{lemma: Easy CV Switches Phase Three} can be applied until reaching a strategy~$\sigma$ such that all improving switches $(e_{i,j,k},t^{\rightarrow})$ with $i\in[n],j,k\in\{0,1\}$ and $\sigma^{(3)}(d_{i,j,k})=F_{i,j}$ were applied.
Since a fixed improving switch $(d_{i,j,k},e_{i,j,k})$ was applied if and only if $\indbit_i=0\vee\indbit_{i+1}=j$, this implies that $\sigma(d_{i,j,k})=F_{i,j}$ is equivalent to $\indbit_i=1\wedge\indbit_{i+1}=j$ for all $i\in[n],j,k\in\{0,1\}$.
Consequently, every cycle center is closed or escapes towards~$t^{\rightarrow}$.
In addition, for suitable indices $i,j,k$, an edge $(d_{i,j,k},F_{i,j})$ is an improving switch exactly if the corresponding switch $(e_{i,j,k},t^{\rightarrow})$ was not applied.
Consequently, \[I_{\sigma}=\{(d_{i,j,k},F_{i,j}), (e_{i,j,k}, t^{\rightarrow})\colon\sigma(e_{i,j,k})=t^{\leftarrow}\}\cup\{(b_1,t^{\rightarrow})\}.\]
Now, as $\occrec^{\sigma}(b_1,t^{\rightarrow})=\occrec^{\sigma}(e_{*,*,*},t^{\rightarrow})$ and $\mathbb{E}^1=\emptyset$, the switch $e\coloneqq(b_1,t^{\rightarrow})$ is applied next due to the tie-breaking rule.
We prove that we can apply row 7 resp. 9 of \Cref{table: Phase 3 Switches}, implying the statement for the case $G_n=S_n$ and arbitrary $\nsb$ and for the case $G_n=M_n$ and $\nsb=1$.
The following claim shows that one of the key requirements for the application of the corresponding statements is fulfilled.

\begin{restatable}{claim}{SigmabarInPhaseThree} \label{claim: Sigmabar in Phase Three}
Let $\sigma$ denote the phase-$3$-strategy in which the improving switch $(b_1,t^{\rightarrow})$ should be applied next.
If $\nsb>1$, then $\sigmabar(eb_{i,j})\wedge\nsigmabar(eg_{i,j})$ for all $(i,j)\in S_1$ and, in addition, $\sigmabar(eb_{i,j})\wedge\sigmabar(eg_{i,j})$ for all $(i,j)\in S_2$.
If $\nsb=1$, then $\sigmabar(eg_{i,j})\wedge\nsigmabar(eb_{i,j})$ for all $(i,j)\in S_4$ and $\sigmabar(eb_{i,j})\wedge\sigmabar(eg_{i,j})$ for all $(i,j)\in S_3$.
\end{restatable}

In addition to the two statements of the claim, it holds that $\sigma(d_{i,j,*})=F_{i,j}$ if and only if $\indbit_i=1\wedge\indbit_{i+1}=j$ for all $i\in[n], j\in\{0,1\}.$.
Consequently, all requirements of row $7$ are met for the case $G_n=S_n$ and $\nsb>1$, implying that the application of $e=(b_1,b_2)$ yields a phase-4-strategy as described by the corresponding rows of \Cref{table: Properties at start of phase,table: Switches at start of phase}.
Analogously, all requirements of row $9$ are met for the case that $\nsb=1$, implying that the application of $e=(b_1,g_1)$ yields a phase-5-strategy as described by the corresponding rows of \Cref{table: Properties at start of phase,table: Switches at start of phase} in this case.

It remains to consider the  case $\nsb>1$ for $G_n=M_n$, implying $t^{\rightarrow}=b_2$ and $t^{\leftarrow}=g_1$.
Using the same argumentation as before, row 1 of \Cref{table: Phase 3 Switches} and \Cref{lemma: Easy CV Switches Phase Three} imply that improving switches within levels $i\geq\nsb$ are applied until we obtain a phase-3-strategy $\sigma$ for $\bit$ with \begin{align*}
I_{\sigma}=&\{(d_{i,j,k},F_{i,j})\colon i<\nsb\wedge\sigma(d_{i,j,k})\neq F_{i,j}\}\\
	&\cup\{(d_{i,j,k},F_{i,j}),(e_{i,j,k},b_2)\colon i\geq\nsb\wedge\sigma(e_{i,j,k})=g_1\}\cup\{(b_1,b_2)\}.
\end{align*}

As no cycle center in any level $i'<\nsb$ was opened yet, the switch $e=(e_{i,j,k},b_2)$ with $i=\nsb-1, j=1-\indbit_{i+1}$ and $k\in\{0,1\}$ is applied next.
Since $\sigma(d_{i,j,k})=F_{i,j}$, row~1 of \Cref{table: Phase 3 Switches} implies $I_{\sigmae}=(I_{\sigma}\setminus\{e\})\cup\{(d_{i,j,k},e_{i,j,k})\}$.
Due to the tie-breaking rule, $(d_{i,j,k},e_{i,j,k})$ is applied next.

\begin{restatable}{claim}{ApplicationCycleEdgesInPhaseThreeMDP} \label{claim: Application Cycle Edges in Phase Three MDP}
The strategy $\sigmae$ meets the five requirements of \Cref{lemma: MDP Phase 3 Open Closed Cycle Center} and the lemma thus describes the application of the improving switch $(d_{i,j,k},e_{i,j,k})$.
\end{restatable}

Therefore, applying $(d_{i,j,k},e_{i,j,k})$ yields a well-behaved phase-3-strategy $\sigma\in\reach{\sigma_0}$ for $\bit $ with $I_{\sigma}=(I_{\sigmae}\setminus\{(d_{i,j,k},e_{i,j,k})\})\cup\{(s_{i,j},b_1)\}$.
We prove $\occrec^{\sigma}(s_{i,j},b_1)<\occrec^{\sigma}(e_{i,j,k},b_2)=\floor{\bit/2}$, implying that $(s_{i,j},b_1)$ is applied next.
It is easy to verify that $(s_{i,j},b_1)\notin\applied{\canstrat}{\sigma}$.
Consequently, by \Cref{table: Occurrence Records} and as $i=\nsb-1$ and $j=1-\indbit_{i+1}=0$, \[\occrec^{\sigma}(s_{i,j},b_1)=\flips{\bit}{i+1}{}-j\cdot\bit_{i+1}=\flips{\bit}{\nsb}{}\leq\floor{\frac{\bit+2}{4}}<\floor{\frac{\bit+1}{2}}\] if $\bit\geq 3$ since $\nsb\geq 2$.
If $\bit_1=1$, then $(s_{i,j},b_1)$ is also the next switch applied as the tie-breaking rule then ranks $(s_{i,j},b_1)$ higher than any switch of the type $(e_{*,*,*},b_2)$.
Since $(e_{i,j,k},b_2), (d_{i,j,k},e_{i,j,k})\in\applied{\canstrat}{\sigma}$ and since the cycle center $F_{i,j}$ was closed when $(e_{i,j,k},b_2)$ was applied, we have $\sigmabar(eb_{i,j})\wedge\nsigmabar(eg_{i,j})$.
Therefore, the fifth row of \Cref{table: Phase 3 Switches} describes the application of $e=(s_{i,j},b_1)$.
Consequently, $\sigmae$ is a phase-3-strategy with $I_{\sigmae}=I_{\sigma}\setminus\{e\}$ and $\occrec^{\sigmae}(s_{i,j},b_1)=\flips{\bit}{\nsb}{}+1=\flips{\bit+1}{\nsb}{}$ by \Cref{lemma: Numerics Of OR}.
Thus, \Cref{table: Occurrence Records} describes the occurrence record of $(s_{i,j},b_1)$ when interpreted for $\bit+1$.
Since $F_{i,j}$ is $b_2$-halfopen for $\sigmae$ whereas $F_{i,1-j}$ is $g_1$-halfopen, $(e_{i,j,1-k},b_2)$ is applied next.
By the first row of \Cref{table: Phase 3 Switches}, this application unlocks $(d_{i,j,1-k},e_{i,j,1-k})$.
Using our previous arguments and observations, it is easy to verify that $(d_{i,j,1-k},e_{i,j,1-k})$ is applied next and that its application is described by the second-to-last row of \Cref{table: Phase 3 Switches}.
The tie-breaking rule then chooses to apply $(e_{i,1-j,k},b_2)\in\mathbb{E}^1$ next.
By row 1 of \Cref{table: Phase 3 Switches}, $(d_{i,1-j,k},e_{i,1-j,k})$ then becomes improving and is applied next.
Its application is described by row 5 of \Cref{table: Phase 3 Switches}.
After applying this switch, we then obtain a strategy $\sigma$ with \begin{align*}
	I_{\sigma}=&\{(d_{i,j,k},F_{i,j})\colon i<\nsb-1\wedge\sigma(d_{i,j,k})\neq F_{i,j}\}\\
		&\cup\{(d_{i,j,k},F_{i,j}),(e_{i,j,k},b_2)\colon i\geq\nsb-1\wedge\sigma(e_{i,j,k})=g_1\}\cup\{(b_1,b_2)\}.
\end{align*}
It is easy to verify that the same arguments can be applied iteratively as applying a switch $(s_{i',j'},b_1)$ with $i'<\nsb$ always requires to open the corresponding cycle center $F_{i',j'}$ first.
Thus, after finitely many iterations, we obtain a strategy $\sigma$ with \[I_{\sigma}=\{(d_{i,j,k},F_{i,j}), (e_{i,j,k}, b_2)\colon\sigma(e_{i,j,k})=g_1\}\cup\{(b_1,b_2)\}.\]
By the same arguments as for $G_n=S_n$, the conditions of the row $8$ of \Cref{table: Phase 3 Switches} are met, so we obtain a strategy as described by the corresponding rows of \Cref{table: Properties at start of phase,table: Switches at start of phase}.
\end{proof}

Note that we implicitly proved the following which follows from \Cref{corollary: Improving Switches of Phase 1 and 2}.

\begin{corollary} \label{corollary: Improving Switches of Phase 3 if there is a 4}
Let $\sigma^{(4)}$ be the phase-4-strategy calculated by the strategy improvement algorithm when starting with a canonical strategy $\canstrat$ having the canonical properties as described by \Cref{lemma: Extended Reaching phase 4 or 5}.
Then, \Cref{table: Occurrence Records} specifies the occurrence record of every improving switch applied during $\canstrat\to\sigma^{(4)}$ when interpreted for $\bit+1$, excluding switches $(g_*,F_{*,*})$.
In addition, each such switch was applied once.
\end{corollary}

As indicated by \Cref{lemma: Extended Reaching phase 4 or 5}, we do not always obtain a phase-5-strategy immediately after phase~3 as there might be improving switches involving selection vertices $s_{i,*}$ in levels $i<\nsb$ that still need to be applied if $G_n=S_n$.
We thus prove that we also reach a phase-$5$-strategy after applying these switches.
Consequently, we always reach a phase-$5$-strategy.
The following lemma generalizes \Cref{lemma: Reaching phase 5}.

\begin{lemma}[Extended version of \Cref{lemma: Reaching phase 5}] \label{lemma: Extended Reaching phase 5}
Let $\canstrat\in\reach{\sigma_0}$ be a canonical strategy for $\bit\in\bitset_n$ having the canonical properties.
After applying finitely many improving switches, the strategy improvement resp. policy iteration algorithm produces a well-behaved phase-5-strategy $\sigma^{(5)}\in\reach{\sigma_0}$ as described by the corresponding rows of \Cref{table: Properties at start of phase,table: Switches at start of phase}.
\end{lemma}

\begin{proof}
By \Cref{lemma: Extended Reaching phase 4 or 5}, it suffices to consider the case $G_n=S_n$ and $\nsb>1$.
The same lemma implies that the strategy improvement algorithm calculates a phase-4-strategy $\sigma$ for $\bit$ with $\sigma\in\reach{\sigma_0}$ and \begin{align*}
I_{\sigma}&=\{(d_{i,j,k},F_{i,j}),(e_{i,j,k},b_2)\colon\sigma(e_{i,j,k})=g_1\}\cup\{(s_{\nsb-1,0},b_1)\}\\
	&\cup\{(s_{i,1},b_1)\colon i\leq\nsb-2\}\cup X_0\cup X_1.
\end{align*}

\begin{restatable}{claim}{FirstSwitchPhaseFour} \label{claim: First Switch Phase Four}
Let $\sigma$ denote the first phase-$4$-strategy in $S_n$ for $\nsb>1$.
Then, the switch $(s_{\nsb-1,0},b_1)$ is applied next and the application of this switch is described by \Cref{lemma: Phase Four Complete}.
\end{restatable}

Consider the case $\nsb=2$ first.
Then, applying $e=(s_{1,0},b_1)$ yields a phase-5-strategy and $\occrec^{\sigmae}(e)=\flips{\bit}{\nsb}{}+1=\flips{\bit+1}{\nsb}{}$ by \Cref{lemma: Numerics Of OR}.
Hence, \Cref{table: Occurrence Records} describes the occurrence record of $e$ with respect to $\bit+1$.
In addition, we then have \begin{align*}
	I_{\sigmae}&=(I_{\sigma}\setminus\{e\})\cup\{(d_{1,0,0},F_{1,0}),(d_{1,0,1},F_{1,0})\})\\
		&=\{(d_{i,j,k},F_{i,j}),(e_{i,j,k},b_2)\colon\sigmae(e_{i,j,k})=g_1\}\cup \{(d_{i,1-\indbit_{i+1},*},F_{i,1-\indbit_{i+1}})\colon i\leq\nsb-1\}\\
		&\quad\cup X_0\cup X_1.
\end{align*}
Since $\sigmae$ is a phase-5-strategy, it has \Pref{REL1}, implying $\relbit{\sigmae}=u=\min\{i\colon\indbit_i=0\}$.
Thus, $\sigmae$ has all properties listed in the corresponding rows of \Cref{table: Properties at start of phase,table: Switches at start of phase}.

Before discussing the case $\nsb>2$, we discuss edges $(d_{i,j,*},F_{i,j})$ that become improving when a switch $(s_{i,j},b_1)$ with $i<\nsb$ and $j=1-\indbit_{i+1}$ is applied, see \Cref{lemma: Phase Four Complete}.
Since $i<\nsb$ implies $1-\indbit_{i+1}=\bit_{i+1}$, their cycle centers $F_{i,j}$ were closed for $\canstrat$.
Therefore, their occurrence records might be very low with respect to the current strategy $\sigma$.
However, their occurrence records are not ``too low'' in the sense that they interfere with the improving switches applied during phase 4.
More precisely, we prove that $i<\nsb$ and $j=\bit_{i+1}$ imply $\occrec^{\canstrat}(d_{i,j,k},F_{i,j})>\floor{(\bit+2)/4}-1$.
By \Cref{table: Occurrence Records}, \[\occrec^{\canstrat}(d_{i,j,k},F_{i,j})=\ceil{\frac{\lastflip{\bit}{i}{\{(i+1,j)\}}+1-k}{2}}.\]
Since $i<\nsb$, we have $\bit_1=\dots=\bit_i=1$ and, by the choice of $j$, $\bit_{i+1}=j$ and $\bit\geq2^{\nsb-1}-1$.
This implies $\lastflip{\bit}{i}{\{(i+1,j)\}}=\bit-\sum(\bit,i)=\bit-2^{i-1}+1$.
Thus \begin{align*}
	\occrec^{\canstrat}(d_{i,j,k},F_{i,j})&=\ceil{\frac{\bit-2^{i-1}+2-k}{2}}\geq\ceil{\frac{\bit-2^{i-1}+1}{2}}=\floor{\frac{2\bit-2^{i}+4}{4}}.
\end{align*}
Since $\floor{(\bit+2)/4}-1=\floor{(\bit-2)/4}$, it suffices to prove $2\bit-2^{i}+4-(\bit-2)>4$.
This follows as $i\leq\nsb-1$ implies \[2\bit-2^i+4-\bit+2=\bit-2^{i}+6\geq2^{\nsb-1}-1-2^i+6\geq2^{i}-2^{i}+5=5.\]

Let $\nsb>2$.
We obtain $\occrec^{\sigmae}(e)=\flips{\bit+1}{\nsb}{}$ as before.
Furthermore, \Cref{lemma: Phase Four Complete} yields\begin{align*}
I_{\sigmae}&=\{(d_{i,j,k},F_{i,j}),(e_{i,j,k},b_2)\colon\sigmae(e_{i,j,k})=g_1\}\\
	&\quad\cup\{(s_{i,1},b_1)\colon i\leq\nsb-2\}\cup\{(d_{\nsb-1,0,0},F_{\nsb-1,0}),(d_{\nsb-1,0,1},F_{\nsb-1,0})\}.
\end{align*}
We show that the switches $(s_{\nsb-2,1},b_1),\dots,(s_{1,1},b_1)$ are applied next and in this order.
To simplify notation, we denote the current strategy by $\sigma$.
By \Cref{table: Occurrence Records}, it holds that $\occrec^{\sigma}(s_{i,1},b_1)=\flips{\bit}{i+1}{}-1$ for all $i\leq\nsb-2$.
Hence $\occrec^{\sigma}(s_{\nsb-2,1},b_1)<\dots<\occrec^{\sigma}(s_{1,1},b_1)$ by \Cref{lemma: Numerics Of OR}.
It thus suffices to show that the occurrence record of $(s_{1,1},b_1)$ is smaller than the occurrence record of any switch improving for $\sigma$ and any improving switch that might be unlocked by applying some switch $(s_{i,1},b_1)$ for $i\leq\nsb-2$.

The second statement follows since $\occrec^{\canstrat}(s_{1,1},b_1)=\flips{\bit}{2}{}-1=\floor{(\bit+2)/4}-1$ and since the occurrence record of any edge that becomes improving is bounded by $\floor{(\bit+2)/4}$ as discussed earlier.
It thus suffices to show the first statement.

Let $e\coloneqq(d_{i,j,k},F_{i,j})\in I_{\sigma}$ with $i\in[n], j,k\in\{0,1\}$ and  $\sigma(e_{i,j,k})=g_1$.
By \Cref{lemma: Extended Reaching phase 3}  and \Cref{lemma: Numerics Of OR}, it then holds that $\occrec^{\sigma}(e)=\maxocc=\flips{\bit}{1}{}$.
In addition, $\nsb>2$ implies $\flips{\bit}{1}{}>\flips{\bit}{\nsb-1}{}$, hence $\occrec^{\sigma}(e)<\occrec^{\sigma}(s_{1,1},b_1)$ follows.
Next let $e\coloneqq(e_{i,j,k},b_2)\in I_{\sigma}$ with $i\in[n], j,k\in\{0,1\}$ and $\sigma(e_{i,j,k})=g_1$.
Then, since $\bit$ is odd, \Cref{table: Occurrence Records} implies  $\occrec^{\sigma}(e)=\floor{\bit/2}=\floor{(\bit+1)/2}-1=\flips{\bit}{1}{}-1$.
Consequently, we have $\occrec^{\sigma}(e)>\occrec^{\sigma}(s_{1,1},b_1)$.
If $\bit+1$ is not a power of two, we need to show this estimation for some more improving switches.
But this can be shown by easy calculations similar to the calculations necessary when discussing the application of $(s_{\nsb-1,0},b_1)$ which can be found in the proof of \Cref{claim: First Switch Phase Four} in \Cref{appendix: Proofs Exponential}.

Consequently, the switches $(s_{\nsb-1},b_1),\dots,(s_{1,1},b_1)$ are applied next, and they are applied in this order.
It is easy to verify that the requirements of \Cref{lemma: Phase Four Complete} are always met, so this lemma describes the application of these switches.
It is also easy to check that the occurrence records of these edges are described by \Cref{table: Occurrence Records} after applying them.
Let $\sigma$ denote the strategy obtained after applying $(s_{1,1},b_1)$.
Then $\sigma$ is a well-behaved phase-$5$-strategy for $\bit$ with $\sigma\in\reach{\sigma_0}$ and $\relbit{\sigma}=\min\{i\colon\indbit_i=0\}$.
This  further implies 
\begin{align*}
I_{\sigma}&=\{(d_{i,j,k},F_{i,j}),(e_{i,j,k},b_2)\colon\sigma(e_{i,j,k})=g_1\}\\
	&\quad\cup\{(d_{i,1-(\bit+1)_{i+1},*},F_{i,1-(\bit+1)_{i+1}})\colon i\leq\nsb-1\}\cup X_0\cup X_1.
\end{align*}
We observe that $\sigma(e_{i,j,k})=g_1$ still implies $\occrec^{\canstrat}(d_{i,j,k},F_{i,j})=\occrec^{\sigma}(d_{i,j,k},F_{i,j})=\maxocc$ for all indices $i\in[n],j,k\in\{0,1\}$ since the corresponding switches are improving since the end of phase~1.
Also, every improving switch was applied at most once and we proved that the occurrence record of every improving switch that was applied is described correctly by \Cref{table: Occurrence Records} when interpreted for $\bit+1$.
Since no improving switches involving cycle vertices were applied, $\sigma(d_{i,j,*})=F_{i,j}$ if and only if $(\bit+1)_i=1$ and $(\bit+1)_{i+1}=j$ where $i\in[n],j\in\{0,1\}$.
Hence, all conditions listed in the corresponding rows of \Cref{table: Properties at start of phase,table: Switches at start of phase} are fulfilled, proving the statement.
\end{proof}

We henceforth use $\sigma^{(5)}$ to refer to the phase-$5$-strategy described by \Cref{lemma: Extended Reaching phase 5}.
As before, we implicitly proved the following corollary which follows from \Cref{corollary: Improving Switches of Phase 1 and 2,corollary: Improving Switches of Phase 3 if there is a 4}.

\begin{corollary} \label{corollary: Improving Switches for Phase 5}
Let $\sigma^{(5)}$ be the phase-5-strategy calculated by the strategy improvement algorithm when starting with a canonical strategy $\canstrat$ having the canonical properties as described by \Cref{lemma: Extended Reaching phase 5}.
Then, \Cref{table: Occurrence Records} specifies the occurrence record of every improving switch applied during $\canstrat\to\sigma^{(5)}$ when interpreted for $\bit+1$, excluding switches $(g_*,F_{*,*})$.
In addition, each such switch was applied once.
\end{corollary}

\subsection{Reaching a canonical strategy part I: Everything but the occurrence records}

There are two major statements that we still have to prove.
First, we have to prove that applying improving switches to $\sigma^{(5)}$ yields a canonical strategy $\sigma_{\bit+1}$ for $\bit+1$ having the canonical properties.
Note that this implies \Cref{lemma: Reaching next phase}, stating that applying improving switches yields the strategies as described by \Cref{table: Properties at start of phase,table: Switches at start of phase}.
Second, we need to investigate the occurrence records of edges $(g_*,F_{*,*})$ which we ignored until now.

We begin by proving the first statement.
We also prove several smaller statements implicitly which will be used when proving that $\sigma_{\bit+1}$ has the canonical properties.

\begin{lemma} \label{lemma: Extended Reaching canonical strategy}
Let $\canstrat\in\reach{\sigma_0}$ be a canonical strategy for $\bit$ having the canonical properties.
Then, applying improving switches according to Zadeh's pivot rule and the tie-breaking rule produces a canonical strategy $\sigma_{\bit+1}\in\reach{\sigma_0}$ for $\bit+1$ with $I_{\sigma_{\bit+1}}=\mathfrak{D}^{\sigma_{\bit+1}}.$
\end{lemma}

\begin{proof}
By \Cref{lemma: Extended Reaching phase 5}, applying improving switches according to Zadeh's pivot rule and our tie-breaking rule yields a phase-$5$-strategy $\sigma\coloneqq\sigma^{(5)}$ for $\bit$ with $\sigma^{(5)}\in\reach{\sigma_0}$ and $\relbit{\sigma}=u=\min\{i\colon\indbit_i=0\}$.
Let $m\coloneqq\max\{i\colon\indbit_i=1\}$.

\noindent
\boldall{Consider the case \boldall{$\nsb=1$.}}
We begin by proving that the occurrence records of the improving switches are bounded by $\maxocc$.
We furthermore characterize the improving switches which will be applied next.

\begin{restatable}{claim}{BoundedORBeginningPhaseFive} \label{claim: Bounded OR Beginning Phase 5}
For all $e\in I_{\sigma}$, it holds that $\occrec^{\sigma}(e)\leq\maxocc$.
Let $e\in I_{\sigma}$ with $\occrec^{\sigma}(e)<\maxocc$.
Then, $e=(d_{i,j,k},F_{i,j})$ with $i\in\{u+1,\dots,m-1\}, j=1-\indbit^{i+1}, k\in\{0,1\}$ and $\canstrat(d_{i,j,k})=F_{i,j}$.
\end{restatable}

Thus, improving switches $(d_{i,j,k},F_{i,j})$ with $i\in\{u+1,\dots,m-1\}, \indbit_i=0, j=1-\indbit_{i+1}$ and $k\in\{0,1\}$ are applied first.
Let $e=(d_{i,j,k}, F_{i,j})$ denote such a switch with $\occrec^{\sigma}(e)<\maxocc$ minimizing the occurrence record.
Since $\canstrat(d_{i,j,k})=F_{i,j}$,  $e$ was not applied during phase~1, it follows that $\occrec^{\sigma}(e)=\occrec^{\canstrat}(e)=\ell^{\bit}(i,j,k)+1$.

\begin{restatable}{claim}{RowOneInPhaseFive} \label{claim: Row One in Phase Five}
Let $\sigma$ denote the phase-$5$-strategy at the beginning of phase $5$ for $\nsb=1$.
Let $i,\in[n], j,k\in\{0,1\}$ such that $e=(d_{i,j,k}, F_{i,j})\in I_{\sigma}$ and $\occrec^{\sigma}(e)<\maxocc$.
Row~1 of \Cref{table: Phase 5 Switches} can be applied to describe the application of $e$.
\end{restatable}

Thus, $\sigmae$ is a well-behaved phase-$5$-strategy for $\bit$ with $\sigmae\in\reach{\sigma_0}$ and $I_{\sigmae}=I_{\sigma}\setminus\{e\}$.
By \Cref{lemma: Progress for unimportant CC} and the choice of $i$ and $j$, it follows that $\ell^{\bit}(i,j,k)+1=\ell^{\bit+1}(i,j,k)$.
In particular, \[\occrec^{\sigmae}(e)=\ell^{\bit}(i,j,k)+1+1=\ell^{\bit+1}(i,j,k)+1\leq\floor{\frac{\bit+1}{2}}\leq\floor{\frac{(\bit+1)+1-k}{2}}.\]
Thus, by choosing the parameter $t_{\bit+1}=1$, which is feasible since $i\neq 1$, the occurrence record of $e$ is described by \Cref{table: Occurrence Records} when interpreted for $\bit+1$.

Now, the same arguments can be used for all improving switches $e'\in\mathbb{D}^1\cap I_{\sigma}$ with $\occrec^{\sigma}(e')<\maxocc$.
All of these switches are thus applied and their occurrence records are specified by \Cref{table: Occurrence Records} when interpreted for $\bit+1$.
After the application of these switches, we obtain a well-behaved phase-5-strategy $\sigma$ for $\bit$ with $\sigma\in\reach{\sigma_0}$ and \begin{align}\label{equation: Characterization NSB==1}\begin{split}
I_{\sigma}=&\{(d_{i,j,k},F_{i,j}),(e_{i,j,k},g_1)\colon\sigma(e_{i,j,k})=b_2\}\\
	&\cup\bigcup_{\substack{i=u+1\\\indbit_i=0}}^{m-1}\left\{e=(d_{i,1-\indbit_{i+1},*},F_{i,1-\indbit_{i+1}})\colon\occrec^{\sigma}(e)=\maxocc\right\}.
\end{split}
\end{align}
In particular, all improving switches have an occurrence record of $\maxocc$.
Thus, the tie-breaking rule now applies a switch of the type $(e_{*,*,*},g_1)$.
Let $i\in[n], j,k\in\{0,1\}$ such that $e\coloneqq(e_{i,j,k},g_1)$ is the next applied improving switch.

\begin{restatable}{claim}{EscapeVerticesPhaseFiveNSBOne} \label{claim: Escape Vertices Phase Five NSB One}
Let $\nsb=1$ and let $\sigma$ denote the strategy obtained after applying all improving switches with an occurrence record less than $\maxocc$ during phase $5$.
Then, \Cref{lemma: Phase Five Escape Easy} can be applied to describe the application of $e=(e_{i,j,k},g_1)$.
\end{restatable}

In fact, \Cref{claim: Escape Vertices Phase Five NSB One} can be applied for any improving switch of the type $(e_{*,*,*},g_1)$.
Furthermore, $\occrec^{\sigmae}(e)$ is specified by \Cref{table: Occurrence Records} when interpreted for $\bit +1$ as $\nsb=1$ implies $\ceil{\bit/2}+1=\ceil{(\bit+1)/2}$.
Depending on whether the conditions listed in the fourth case of \Cref{lemma: Phase Five Escape Easy} are fulfilled, either \begin{align*}
	I_{\sigmae}&=(I_{\sigmae}\setminus\{e\})\cup\{(d_{i,j,1-k},F_{i,j}), (g_{i},F_{i,j})\}\quad\text{or}\quad I_{\sigmae}=(I_{\sigmae}\setminus\{e\})\cup\{(d_{i,j,1-k},F_{i,j})\}.
\end{align*}
In particular, $\tilde{e}\coloneqq(d_{i,j,1-k},F_{i,j})$ becomes improving in either case.
As formalized by the following corollary,  $\tilde{e}$ has an occurrence record of at least $\maxocc$.
This corollary will be used in later arguments, hence it is not a claim as we use the term claim solely for statements that are only relevant within a single proof.
Nevertheless, its proof is deferred to \Cref{appendix: Proofs Exponential}.

\begin{restatable}{corollary}{FulfillingORFiveIfNSBIsOne} \label{corollary: Fulfilling ORFive If NSB Is One}
Let $\nsb=1$ and $i\in[n], j,k\in\{0,1\}$.
If the edge $\tilde{e}=(d_{i,j,1-k},F_{i,j})$ becomes improving during phase $5$ due to the application of $(e_{i,j,k},g_1)$, then the corresponding strategy has \Pref{OR4}$_{i,j,1-k}$.
\end{restatable}

Now, consider the case that $(g_i,F_{i,j})$ becomes improving when applying $(e_{i,j,k},g_1)$.
We prove that this implies $(g_i,F_{i,j})\notin\applied{\canstrat}{\sigmae}$.
The conditions stated in \Cref{lemma: Phase Five Escape Easy} imply that the switch was not applied previously in phase~$5$.
For the sake of a contradiction, assume that $(g_i,F_{i,j})$ was applied during phase~$1$ of the current transition.
Then, by \Cref{corollary: Selection Vertices In Phase One}, the cycle center $F_{i,j}$ was closed during phase~$1$.
Since  $(e_{i,j,k},g_1)$ was applied immediately before unlocking $(g_i,F_{i,j})$, we have $\occrec^{\canstrat}(d_{i,j,k},F_{i,j})=\maxocc$ by \Cref{lemma: Extended Reaching phase 5}.
However, by \Cref{corollary: Selection Vertices In Phase One}, a cycle center can only be closed during phase~1 if either $i=\nsb$ or if the occurrence record of both cycle edges is less than $\maxocc$.
We thus need to have $i=\nsb=1$.
But then $\indbit_i=1$, implying that $(g_i,F_{i,j})$ cannot become improving.
Hence, a switch $(g_i,F_{i,j})$ that is unlocked during phase 5 was not applied earlier in the same transition if $\nsb=1$.

Since $\occrec^{\sigmae}(g_i,F_{i,j})=\occrec^{\canstrat}(g_i,F_{i,j})$, we have $\occrec^{\sigmae}(g_i,F_{i,j})\leq\occrec^{\canstrat}(d_{i,j,k},F_{i,j})=\maxocc$ by \Cref{table: Occurrence Records}.
By \Cref{corollary: Fulfilling ORFive If NSB Is One}, $\occrec^{\sigmae}(d_{i,j,1-k},F_{i,j})\geq\maxocc$.
Therefore, the occurrence record of any improving switch except $(g_i,F_{i,j})$ is at least $\maxocc$.
Thus, $(g_i,F_{i,j})$ either uniquely minimizes the occurrence record or has the same occurrence record as all other improving switches.
Consequently, by the tie-breaking rule, $(g_i,F_{i,j})$ is applied next in either case.

We prove that row~2 of \Cref{table: Phase 5 Switches} applies to this switch.
Since $\nsb=1, \relbit{\sigmae}=u>1$ and $\indbit_i=0$, it suffices to prove $\sigmaebar(eg_{i,j})\wedge\nsigmaebar(eb_{i,j})$.
But this follows as we applied $(e_{i,j,k},g_1)$ earlier and since $F_{i,j}$ was mixed when this switch was applied.
Observe that the following corollary holds due to the conditions which specify when a switch $(g_i,F_{i,j})$ is unlocked, independent on $\nsb$.

\begin{corollary} \label{corollary: Selection Vertices Phase 5}
Let $\nsb=1$.
If an improving switch $(g_i,F_{i,j})$ is applied during phase $5$, then the resulting strategy has \Pref{SVG}$_i$/(\ref{property: SVM})$_{i}$.
\end{corollary}

Let $\sigma$ denote the strategy obtained after applying $(e_{i,j,k},g_1)$ (and potentially $(g_i,F_{i,j})$ if it became improving).
Assume that there is an improving switch of the type $(e_{*,*,*},g_1)\in I_{\sigma}$.
Then, by \Cref{lemma: Phase Five Escape Easy} resp. row~ 2 of \Cref{table: Phase 5 Switches}, $\sigma$ is a phase-5-Strategy for $\bit$.
By our previous discussion, the occurrence records of all improving switches are at least $\maxocc$.
Among all improving switches with an occurrence record of exactly $\maxocc$, the tie-breaking rule then decides which switch to apply.
There are two types of improving switches.
Each switch is either of the form $(d_{*,*,*},F_{*,*j})$ or of the form $(e_{i',j',k'},g_1)$ for indices $i'\in[n], j',k'\in\{0,1\}$ with $\sigma(d_{i',j',k'})=e_{i',j',k'}$.
Since every edge $(e_{*,*,*},g_1)$ minimizes the occurrence record among all improving switches, an edge of this type is chosen.
Let $(e_{i',j',k'},g_1)$ denote this switch.
Then, the same arguments used previously can be used again.
More precisely, \Cref{lemma: Phase Five Escape Easy} applies to this such a switch, making the edge $(d_{i',j',1-k'},F_{i',j'})$ and eventually also $(g_{i'},F_{i',j'})$ improving.
Also, \Cref{corollary: Fulfilling ORFive If NSB Is One,corollary: Selection Vertices Phase 5} apply to these switches and another switch of the form $(e_{*,*,*},g_1)$ is applied afterwards.
Thus, inductively, all remaining improving switches $(e_{*,*,*},g_1)$ are applied.

Let $\sigma$ denote the strategy that is reached before the last improving switch $(e_{*,*,*},g_1)$ is applied.
We argue that this switch is $e\coloneqq(e_{1,1-\indbit_2,k},g_1)$ for some $k\in\{0,1\}$ and that $\sigma$ has \Pref{SVG}$_{i}$/(\ref{property: SVM})$_{i}$ for all $i\in[n]$.
As the tie-breaking rule applies improving switches in higher levels first, it suffices to prove that there there is a $k\in\{0,1\}$ such that $e\in I_{\sigma^{(5)}}$.
This however follows from \Cref{lemma: Extended Reaching phase 5} as $\nsb=1$ implies $(1,\indbit_2)\in S_3$.
It remains to prove that $\sigma$ has \Pref{SVG}$_{i}$/(\ref{property: SVM})$_{i}$ for all $i\in[n]$.

\begin{restatable}{claim}{SVGSVMInPhaseFive} \label{claim: SVG SVM In Phase Five}
If $\nsb=1$, then the strategy $\sigma$ obtained before the application of the switch $e\coloneqq(e_{1,1-\indbit_2,k},g_1)$ has \Pref{SVG}$_{i}$/(\ref{property: SVM})$_{i}$ for all $i\in[n]$.
\end{restatable}

Thus, \Cref{lemma: Phase Five Escape Easy} applies to $e\coloneqq (e_{1,\indbit^{\sigma}_2,k},g_1)$.
Let $\sigma_{\bit+1}\coloneqq\sigmae$ denote the strategy obtained by applying $e$.
Then, as we assume that there are no further indices $(i',j',k')$ such that $(e_{i',j',k'},g_1)\in I_{\sigma_{\bit+1}}$, \Cref{lemma: Phase Five Escape Easy} implies that $\sigma_{\bit+1}$ is a phase-1-strategy for $\bit+1$ with $\sigma_{\bit+1}\in\reach{\sigma_0}$.
Since every edge was applied at most once during $\canstrat\to\sigma^{(5)}$ by \Cref{lemma: Extended Reaching phase 5} and since no edge applied during $\sigma^{(5)}\to\sigma_{\bit+1}$ was applied earlier, every edge was applied at most once as improving switch during $\canstrat\to\sigma_{\bit+1}$.
We furthermore implicitly proved the following corollary where the second statement follows from \Cref{corollary: Fulfilling ORFive If NSB Is One}.

\begin{corollary} \label{corollary: Cycle Edges Phase Five If NSB One}
Let $\nsb=1$ and let $\sigma_{\bit+1}$ denote the strategy obtained after the application of the final improving switch $(e_{*,*,*},g_1)$.
Let $i\in[n]$ and $j,k\in\{0,1\}$.
Then, $(d_{i,j,k},F_{i,j})\in\applied{\sigma^{(5)}}{\sigma_{\bit+1}}$ if and only if $\canstrat(d_{i,j,k})=F_{i,j}, \occrec^{\canstrat}(d_{i,j,k},F_{i,j})<\maxocc, i\in\{u+1,\dots,m-1\}, \indbit_i=0$ and $j=1-\indbit_{i+1}$.
In addition, $\sigma_{\bit+1}$ has \Pref{OR2}$_{i,j,k}$.
\end{corollary}

It remains to prove that $\sigma_{\bit+1}$ is a canonical strategy for $\bit+1$ with $I_{\sigma_{\bit+1}}=\mathfrak{D}^{\sigma_{\bit+1}}$.

We begin with the second statement.
This can be proven by using the characterization given in Equation~(\ref{equation: Characterization NSB==1}) and showing $I_{\sigma}\subseteq \mathfrak{D}^{\sigma_{\bit+1}}$ and $I_{\sigma}\supseteq \mathfrak{D}^{\sigma_{\bit+1}}$.

\begin{restatable}{claim}{ImprovingSwitchesForBPlusOneForNSBOne} \label{claim: Improving Switches For B plus One for NSB}
It holds that $I_{\sigma_{\bit+1}}=\{(d_{i,j,k},F_{i,j})\colon\sigma_{\bit+1}(d_{i,j,k})\neq F_{i,j}\}$.
\end{restatable}

To simplify notation, let $\sigma\coloneqq\sigma_{\bit+1}$.
We now prove that $\sigma$ is a canonical strategy for $\bit$, concluding the case $\nsb=1$.
Since $\sigma$ is a phase-1-strategy for $\bit+1$, it holds that $\bit+1=\indbit$.
Consider the conditions listed in \Cref{definition: Canonical Strategy} resp. \ref{definition: Canonical Strategy}.
Condition~1 is fulfilled since $\sigma(e_{*,*,*})=g_1$ and $\nsb=1$.
Condition~2(a) is fulfilled since $\indbit^{\sigma}_i=(\bit+1)_i=1$ implies $\sigma(b_i)=g_i$ by \Pref{EV1}$_i$ for every $i\in[n]$.
Consider condition~2(b) and let $i\in[n]$.
If $(\bit+1)_i=1$, then $F_{i,(\bit+1)_{i+1}}$ is closed by \Pref{EV1}$_i$.
We prove that $(\bit+1)_{i}=1$ implies that $F_{i,j}$ with $j\coloneqq1-(\bit+1)_{i+1}$ cannot be closed.

Consider $\sigma^{(5)}$ and let $k\in\{0,1\}$.
Then, $\sigma^{(5)}(d_{i,j,k})=F_{i,j}$ if and only if $\indbit^{\sigma^{(5)}}_i=1\wedge\indbit^{\sigma^{(5)}}_{i+1}=j$.
Hence, $\sigma^{(5)}(d_{i,j,0})\neq F_{i,j}$ and it suffices to show that $e\coloneqq(d_{i,j,0},F_{i,j})$ was not applied during $\sigma^{(5)}\to\sigma$.
By \Cref{corollary: Cycle Edges Phase Five If NSB One}, it suffices to show $\occrec^{\sigma^{(5)}}(e)\geq\maxocc$.
By \Cref{lemma: Numerics Of Ell}, it holds that $\ell^{\bit}(i,j,0)\geq\maxocc.$
Since $\nsb=1$,  \Pref{OR4}$_{i,j,0}$ implies $\occrec^{\canstrat}(e)\neq\ell^{\bit}(i,j,0)-1$, hence $\occrec^{\sigma^{(5)}}(e)\geq\occrec^{\canstrat}(e)\geq\maxocc$.
Consequently, condition~2(b) is fulfilled.
Condition~2(c) is fulfilled by $\indbit^{\sigma}=\bit+1$ and \Pref{EV2}$_*$.

Conditions~3(a) and 3(b) are fulfilled since~$\sigma$ has \Pref{EV1}$_*$.
Consider condition~3(c) and let $i\in[n]$.
We prove that $(\bit+1)_i=0, j=1-(\bit+1)_{i+1}$ and $\sigmabar(d_{i,j})$ imply $\sigmabar(g_{i})=F_{i,j}$.
Since $S_n$ is a sink game and $M_n$ is weakly unichain, $F_{i,j}$ being closed implies $\valustar_{\sigma}^*(F_{i,j})=\valustar_{\sigma}^*(s_{i,j})$.
Thus, $\valu_{\sigma}^*(F_{i,j})=\ubracket{s_{i,j}}\oplus\valu_{\sigma}^*(g_1)$ by the choice of $j$ and since $\nsb=1$.
As shown by \Cref{lemma: Exact Behavior Of Random Vertex,lemma: Exact Behavior Of Counterstrategy}, $\relbit{\sigma}\neq1, \sigmabar(eg_{i,1-j}), \neg\sigmabar(eb_{i,1-j})$ and $1-j=\indbit_{i+1}$ implies $\valu_{\sigma}^*(F_{i,1-j})=\{F_{i,1-j}, d_{i,1-j,k}, e_{i,1-j,k}, b_1\}\cup\valu_{\sigma}^*(g_1)$ for some $k\in\{0,1\}$.
But this implies $\sigma(g_i)=F_{i,j}$ since $(g_i,F_{i,1-j})\in I_{\sigma}$ otherwise, contradicting $I_{\sigma}=\{(d_{i,j,k},F_{i,j})\colon\sigma(d_{i,j,k})=F_{i,j}\}$.
Consider condition~3(d) and let $i\in[n]$ and let $j\coloneqq 0$ if $G_n=S_n$ and $j\coloneqq\indbit_{i+1}$ if $G_n=M_n$.		
It suffices to prove $\valu_{\sigma}^*(F_{i,j})\succ\valu_{\sigma}^*(F_{i,1-j})$ if none of the cycle centers are closed.
For $G_n=M_n$, this follows from \Cref{lemma: Both CC Open For MDP} or an easy calculation using $i\geq1=\nsb$.
For $G_n=S_n$, this follows from $\Omega(F_{i,0})>\Omega(F_{i,1})$ and since both priorities are even.		

Conditions~4 and 5 follow easily since $\sigma$ has \Pref{USV1}$_{*}$.
For condition~6, let $i\coloneqq\ell(\bit+2), j\coloneqq(\bit+1)_{i+1}$ and $k\in\{0,1\}$.
Since $\ell(\bit+1)=1$, we have $i\geq 2$ and $\bit_i=(\bit+1)_i=0$ as well as $\bit_{i+1}=(\bit+1)_{i+1}=j$.
We prove $\sigma(d_{i,j,k})\neq F_{i,j}$.
For the sake of a contradiction, let $\sigma(d_{i,j,k})=F_{i,j}$.
Then, by the choice of~$i$ and $j$ and \Cref{lemma: Extended Reaching phase 5}, it holds that$(d_{i,j,k},F_{i,j})\in\applied{\sigma^{(5)}}{\sigma}$.
Thus, by \Cref{corollary: Cycle Edges Phase Five If NSB One} and \Pref{OR2}$_{i,j,k}$, it holds that $\occrec^{\sigma^{(5)}}(d_{i,j,k},F_{i,j})<\maxocc$ and $\occrec^{\canstrat}(d_{i,j,k},F_{i,j})=\ell^{\bit}(i,j,k)+1$.
But, by \Cref{lemma: Numerics Of Ell}, we have \[\ell^{\bit}(i,j,k)=\ceil{\frac{\bit+2^{i-1}+\sum(\bit,i)+1-k}{2}}\geq\ceil{\frac{\bit+3-k}{2}}=\floor{\frac{\bit+2-k}{2}},\] which is a contradiction.
Hence, $\sigma(d_{i,j,k})\neq F_{i,j}$.

This concludes the case $\nsb=1$.
We now prove the same statements for the case $\nsb>1$.

\noindent
\boldall{Consider the case \boldall{$\nsb>1$.}}
Then, $\bit$ is odd and  $\maxocc=\floor{\bit/2}+1$.
By \Cref{lemma: Extended Reaching phase 5}, applying improving switches according to Zadeh's pivot rule and the tie-breaking rule given in \Cref{definition: Tie-Breaking exponential} yields a well-behaved phase-5-strategy $\sigma$ for $\bit$ with $\sigma\in\reach{\sigma_0}$ and $\relbit{\sigma}=u$.
In addition \begin{equation} \label{equation: IS Beginning Phase 5 For NSB>1}
\begin{split}
I_{\sigma} =\, &\{(d_{i,j,k},F_{i,j}),(e_{i,j,k},b_2)\colon\sigma(e_{i,j,k})=g_1\}\\
	&\cup\bigcup_{i=1}^{\nsb-1}\{(d_{i,1-\indbit_{i+1},*},F_{i,1-\indbit_{i+1}})\}\cup X_0\cup X_1,
\end{split}
\end{equation}
where $X_k$ is defined as in \Cref{table: Switches at start of phase}.

To deduce which improving switch is applied next, it is necessary to analyze their occurrence records.

\begin{restatable}{claim}{ORBeginningPhaseFiveForNSBLargerOne} \label{claim: OR Beginning Phase 5 For NSB>1}
Let $\nsb>1$.
The occurrence records of the improving switches with respect to the phase-$5$-strategy $\sigma$ described by \Cref{lemma: Extended Reaching phase 5} is described correctly by \Cref{table: Or Beginning Phase 5 For NSB>1}.
\end{restatable}

\begin{table}[ht]
\centering
\begin{tabular}{|c||c|c|c|}\hline
Switch $e$ 							&	$(d_{i,j,k}, F_{i,j})$				&$(e_{i,j,k},b_2)$			&$(d_{\nsb,1-\bit_{\nsb+1},k}, F_{\nsb,1-\bit_{\nsb+1}})$\\\hline
Condition 								&\multicolumn{2}{c|}{$\sigma(e_{i,j,k})=g_1$}	&--\\\hline
$\occrec^{\sigma}(e)$	&$=\maxocc$							&$=\maxocc-1$				&$=\maxocc$\\\hline
\end{tabular}

\bigskip

\begin{tabular}{|c||c|c|}\hline
Switch $e$ 										&	\multicolumn{2}{c|}{$(d_{i,j,k}, F_{i,j})$}																									\\\hline
\multirow{2}{*}{Condition} 		&\multicolumn{2}{c|}{$i\in\{\nsb+1,\dots,m\}, \bit_i=0, j=1-\bit_{i+1}, k\in\{0,1\}$}	\\\cline{2-3}
																	&$\qquad\canstrat(d_{i,j,k})=F_{i,j}\qquad$	&$\canstrat(d_{i,j,k})\neq F_{i,j}$					\\\hline
$\occrec^{\sigma}(e)$				&$\leq \maxocc-1 $							&$=\maxocc$				\\\hline
\end{tabular}

\bigskip

\begin{tabular}{|c||c|c|c|c|}\hline
Switch $e$ 										&	\multicolumn{4}{c|}{$(d_{i,j,k}, F_{i,j})$}																									\\\hline
\multirow{2}{*}{Condition} 	&\multicolumn{4}{c|}{$i\leq \nsb-1, j=1-\bit_{i+1}$}	\\\cline{2-5}
																&$i=1$	& $i=2$		&$i=3$		&$i>3$			\\\hline
$\occrec^{\sigma}(e)$				&$\maxocc$	&$=\maxocc-k$	&$=\maxocc-1-k$		&$<\maxocc-1$\\\hline
\end{tabular}
\caption[Occurrence records at the beginning of phase $5$ for $\nsb>1$]{Occurrence records of the improving switches at the beginning of phase $5$ for $\nsb>1$.} \label{table: Or Beginning Phase 5 For NSB>1}
\end{table}

We partition $I_{\sigma}$ into three subsets, based on their occurrence records.
An improving switch $e\in I_{\sigma}$ is called \begin{itemize}
	\item \emph{type 1 switch} if $\occrec^{\sigma}(e)=\maxocc$
	\item \emph{type 2 switch} if $\occrec^{\sigma}(e)=\maxocc-1$ and
	\item \emph{type 3 switch} if $\occrec^{\sigma}(e)<\maxocc-1$.
\end{itemize}

By Zadeh's pivot rule, type 3 switches are applied first, and we discuss the application of these switches next.


\begin{restatable}{claim}{ApplicationOfPhaseThreeSwitches} \label{claim: Application of Phase Three Switches}
Let $\nsb>1$ and consider the first phase-$5$-strategy.
The application of type $3$ switches is described by row 1 of \Cref{table: Phase 5 Switches}.
\end{restatable}

Let $i\in[n], j,k\in\{0,1\}$ and let $e=(d_{i,j,k},F_{i,j})$ denote the type $3$ switch that is applied next.
We show that \Cref{table: Occurrence Records} specifies the occurrence record of $e$ after its application when interpreted for $\bit+1$.
Consider the case  $i\in\{\nsb+1,\dots,m-1\}, \indbit_i=0, j=1-\indbit_{i+1}$ and $k\in\{0,1\}$ first.
Since $e$ is a type 3 switch, it holds that $\canstrat(d_{i,j,k})=F_{i,j}$, implying $\occrec^{\canstrat}(e)=\ell^{\bit}(i,j,k)+1$ by \Pref{OR2}$_{i,j,k}$.
Thus, the statement follows since $\ell^{\bit+1}(i,j,k)=\ell^{\bit}(i,j,k)+1$ by \Cref{lemma: Numerics Of OR}.
Now consider the case  $i\leq\nsb-1$.
Then, $F_{i,j}$ was closed with respect to $\canstrat$ and $j=\bit_{i+1}=1-\indbit_{i+1}$.
It is easy to verify that this implies $\occrec^{\canstrat}(e)=\ceil{(\bit-\sum(\bit,i)+1-k)/2}$.
Since $(\bit+1)_i=0\wedge(\bit+1)_{i+1}\neq j$ and the switch $e$ is applied, it suffices to prove $\ell^{\bit+1}(i,j,k)=\ceil{(\bit-\sum(\bit,i)+1-k)/2}$ as we can then choose $t^{\bit+1}=1$ as feasible parameter.
This however follows directly from \begin{align*}
	\ell^{\bit+1}(i,j,k)&=\ceil{\frac{\bit+1-2^{i-1}+\sum(\bit+1,i)+1-k}{2}}=\ceil{\frac{\bit+1-2^{i-1}+1-k}{2}}\\
		&=\ceil{\frac{\bit+1+\sum(\bit,i)-1+1-k}{2}}=\ceil{\frac{\bit-\sum(\bit,i)+1-k}{2}}.
\end{align*}
Note that we do not prove yet that choosing this parameter is in accordance with Properties (\ref{property: OR1})$_{*,*,*}$ to (\ref{property: OR4})$_{*,*,*}$.
Since $e$ is a type 3 switch, this furthermore implies  $\occrec^{\sigmae}(e)\leq\maxocc-1=\floor{(\bit+1+1)/2}-1$.
Hence, $\sigmae$ has \Pref{OR1}$_{i,j,k}$ and we have implicitly proven the following corollary.

\begin{corollary} \label{corollary: Improving Switches Of Type 3}
Let $\nsb>1$ and $i\in[n],j,k\in\{0,1\}$.
Every switch $e=(d_{i,j,k},F_{i,j})$ with $\occrec^{\canstrat}(e)<\maxocc-1$  is applied during phase 5, and the resulting strategy has \Pref{OR1}$_{i,j,k}$.
\end{corollary}

Now, the first row of \Cref{table: Phase 5 Switches} and the corresponding arguments can be applied for every improving switch of type~$3$.
Thus, we obtain a phase-$5$-strategy $\sigma\in\reach{\sigma_0}$ such that every improving switch is of type $1$ or $2$.
The next improving switch that is applied has an occurrence record of $\floor{(\bit+1)/2}-1$, i.e., it is of type $2$, so we discuss the application of these switches next.

Since any improving switch is either of the form $(d_{*,*,*},F_{*,*})$ or $(e_{*,*,*},b_2)$ and since the latter switches are of type~2, some improving switch $(e_{*,*,*},b_2)$ is applied next due to the tie-breaking rule.

\begin{restatable}{claim}{ApplicationOfTypeTwoSwitchesEscape} \label{claim: Application of Phase Three Switches Escape}
Let $\nsb>1$ and let $\sigma$ denote the strategy obtained after the application of all improving switches of type $3$ during phase $5$.
The application of type $2$ switches of the form $(e_{*,*,*},b_2)$ is described by row~$1$ of \Cref{lemma: Phase Five Escape Easy}.
\end{restatable}

Let $i\in[n], j,k\in\{0,1\}$ and let $e=(e_{i,j,k},b_2)$ denote the applied improving switch.
Then, \Cref{table: Occurrence Records} describes the occurrence record of $e$ after the application when interpreted for $\bit+1$ since $\occrec^{\sigmae}(e)=\occrec^{\canstrat}(e)+1=\floor{\bit/2}+1=\maxocc$.
Now, by \Cref{lemma: Phase Five Escape Easy}, $(d_{i,j,1-k},F_{i,j})\in I_{\sigmae}$ and the edge $(g_{i},F_{i,j})$ might become improving for $\sigmae$.
The strategy~$\sigmae$ is now either a phase-5-strategy for $\bit$ or a phase-1-strategy for $\bit+1$.
The following corollary which is proven in \Cref{appendix: Proofs Exponential} now describes the improving switch $(d_{i,j,1-k},F_{i,j})$ in more detail.

\begin{restatable}{corollary}{FulfillingORFiveIfNSBIsNotOne} \label{corollary: Fulfilling ORFive If NSB Is Not One}
Let $i\in[n], j,k\in\{0,1\}$ and let $\sigma$ denote the strategy obtained after the application of an improving switch $(e_{i,j,k},b_2)$ during phase $5$.
If $(d_{i,j,1-k},F_{i,j})\in I_{\sigma}$, then$\sigma$ has \Pref{OR4}$_{i,j,1-k}$ and it holds that $\min_{k'\in\{0,1\}}\occrec^{\canstrat}(d_{i,j,k'},F_{i,j})\leq\maxocc-1$.
\end{restatable}

We now use \Cref{corollary: Fulfilling ORFive If NSB Is Not One} to prove that $e\coloneqq (g_i,F_{i,j})$ is applied next if it becomes improving.
For simplicity, let $\sigma$ denote the current strategy that was obtained by applying an improving switch $(e_{i,j,*},b_2)$ according to \Cref{lemma: Phase Five Escape Easy}.

By the tie-breaking rule and \Cref{corollary: Fulfilling ORFive If NSB Is Not One}, it suffices to prove \begin{equation} \label{equation: OR of selector vertex bounded phase 5}
\occrec^{\sigma}(g_i,F_{i,j})\leq\floor{\frac{\bit+1}{2}}-1.
\end{equation}
Since \Cref{table: Occurrence Records} and \Cref{corollary: Fulfilling ORFive If NSB Is Not One} yield \[\occrec^{\canstrat}(g_i,F_{i,j})\leq\min_{k'\in\{0,1\}}\occrec^{\canstrat}(d_{i,j,k'},F_{i,j})\leq\floor{\frac{\bit+1}{2}}-1,\] it suffices to prove $(g_i,F_{i,j})\notin\applied{\canstrat}{\sigma}$.

\begin{restatable}{claim}{SelectionVerticesPhaseFiveNSBLargerOne} \label{claim: Selection Vertices Phase 5 if NSB>1}
Let $i\in[n], j,k\in\{0,1\}$ and let $\sigma$ denote the strategy obtained after the application of an improving switch $(e_{i,j,k},b_2)$ during phase $5$.
If $(g_i,F_{i,j})\in I_{\sigma}$, then  $(g_i,F_{i,j})\notin\applied{\canstrat}{\sigma}$.
\end{restatable}

Due to the tie-breaking rule, $(g_i,F_{i,j})$ is thus applied next.
We prove that row~2 of \Cref{table: Phase 3 Switches} applies to the application of $e$.

First, $\indbit_i=0$ follows from the conditions of \Cref{lemma: Phase Five Escape Easy}.
Second, $\sigmabar(eb_{i,j})\wedge\nsigmabar(eg_{i,j})$ follows as the cycle center $F_{i,j}$ was mixed earlier and since we just applied $(e_{i,j,k},b_2)$.
To prove that  $\sigmabar(d_{i',j'})\vee[\sigmabar(eb_{i',j'})\wedge\nsigmabar(eg_{i',j'})]$ holds for all $i'\geq i$ and $j\in\{0,1\}$, fix some $i'\geq i$ and $j'\in\{0,1\}$.
If $\indbit_{i'}=1\wedge\indbit_{i'+1}=j'$, then the statement follows from \Pref{EV1}$_{i'}$.
We may hence assume $\indbit_{i'}=0\vee j'\neq\indbit_{i+1}$ and that $F_{i',j'}$ is not closed.
Then, by \Cref{lemma: Extended Reaching phase 5}, either $\sigmabar(eb_{i',j'})\wedge\sigmabar(eg_{i',j'})$ or $\sigmabar(eb_{i',j'})\wedge\nsigmabar(eg_{i',j'})$.
Assume that the first case was true, implying $i'\neq i$.
Then, $\sigma(e_{i',j',k})=g_1$ and $\sigma(d_{i',j',k})=e_{i',j',k}$ for some $k\in\{0,1\}$.
This in particular implies $(e_{i',j',k},b_2)\in I_{\sigma}$.
But this is a contradiction to the fact that we apply improving switches according to the tie-breaking rule since $i'> i$ implies that the switch $(e_{i',j',k},b_2)$ is applied before the switch $(e_{i,j',k},b_2)$.

Hence, all requirements of the second row of \Cref{table: Phase 5 Switches} are met.
Further note that the strategy obtained after applying the switch has \Pref{SVG}$_i$/(\ref{property: SVM})$_{i}$ due to the conditions described in \Cref{lemma: Phase Five Escape Easy}.
In particular, \Cref{corollary: Selection Vertices Phase 5} also holds for $\nsb>1$.

After the application of $(e_{i,j,k},b_2)$ (or $(g_i,F_{i,j})$ if it becomes improving), the tie-breaking rule determines which switch is applied next.
Since $(d_{i,j,1-k},F_{i,j})$ has an occurrence record of at least $\maxocc-1$, another switch of the type $(e_{*,*,*},b_2)$ is applied.
But then, the same arguments used previously can be applied again.
That is, we can apply some switch $(e_{i',j',k'},b_2)$, making $(d_{i',j',1-k'},F_{i',j'})$ improving, and eventually making $(g_{i'},F_{i',j'})$ improving as well.
The switch $(g_{i'}, F_{i',j'})$ is applied immediately (if it becomes improving) whereas the other switch is not applied.
Then, inductively, all remaining switches of the form $(e_{*,*,*},b_2)$ are applied.

Let $\sigma$ denote the strategy that is reached after applying the final improving switch of the type $(e_{*,*,*},b_2)$.
We prove that $\sigma$ has Property (SV*)$_1$ if $(g_{1},F_{1,j})$ does not become improving and Property (SV*)$_i$ for all $i\geq 2$.
We first determine which is the last switch of the form $(e_{*,*,*},b_2)$ that will be applied.
It holds that $(1,\indbit_2)\in S_2$, implying $(e_{1,\indbit_2,k},b_2)\in I_{\sigma^{(5)}}$ for some $k\in\{0,1\}$ by \Cref{lemma: Extended Reaching phase 5}.
Due to the tie-breaking rule, this is thus the last switch of the form $(e_{*,*,*},b_2)$ that will be applied.
This might also unlock the corresponding improving switch $(g_1,F_{1,\indbit_2})$.
Let $\sigma$ denote the strategy obtained after the application of the switch $(e_{1,\indbit_2,k},b_2)$ resp. after the application of the switch $(g_1,F_{1,\indbit_2})$ if it becomes improving.

\begin{restatable}{claim}{PropertySVAtEndOfPhaseFive} \label{claim: Property SV at End of Phase Five}
Let $\nsb>1$.
The strategy $\sigma$ obtained after the application of the final improving switch of phase $5$ has Property (SV*)$_{i}$ for all $i\in[n]$.
\end{restatable}

Thus, by \Cref{lemma: Phase Five Escape Easy} resp. the row~2 of \Cref{table: Phase 5 Switches}, $\sigma$ is a well-behaved phase-1-strategy for $\bit+1$ with $\sigma\in\reach{\sigma_0}$.
It remains to show that $\sigma$ is a canonical strategy for $\bit+1$ with $I_{\sigma}=\{(d_{i,j,k},F_{i,j})\colon\sigma(d_{i,j,k})\neq F_{i,j}\}$.
This is formalized by the two following claims whose proofs can be found in \Cref{appendix: Proofs Exponential}.
The first statement is again shown by proving that the two sets are contained in each other.
The proof of the second statement is analogous to the corresponding statement for $\nsb=1$ and is thus deferred to the appendix.

\begin{restatable}{claim}{ImprovingSwitchesOfNewCanstrat} \label{claim: Improving Switches of new Canstrat}
Let $\sigma$ denote the strategy obtained after applying the final improving switch of phase~$5$ for $\nsb>1$.
Then $I_{\sigma}=\{(d_{i,j,k},F_{i,j})\colon\sigma(d_{i,j,k})\neq F_{i,j}\}$.
\end{restatable}

\begin{restatable}{claim}{NewCanstratIsCanstrat} \label{claim: New Canstrat is Canstrat}
Let $\sigma$ denote the strategy obtained after applying the final improving switch of phase~$5$ for $\nsb>1$.
Then $\sigma$ is a canonical strategy for $\bit+1$.
\end{restatable}

This concludes the case $\nsb>1$ and hence proves the statement.
\end{proof}

Using the previous similar corollaries of this type, it follows that we also implicitly proved the following corollary.

\begin{corollary} \label{corollary: Improving Switches for new canstrat}
Let $\sigma_{\bit+1}$ be the canonical strategy for $\bit+1$ calculated by the strategy improvement algorithm as described by \Cref{lemma: Extended Reaching canonical strategy}.
Then, \Cref{table: Occurrence Records} specifies the occurrence record of every improving switch applied until reaching $\sigma_{\bit+1}$, excluding switches $(g_*,F_{*,*})$, when interpreted for $\bit+1$.
In addition, each such switch was applied once.
\end{corollary}

It remains to prove that the canonical strategy $\sigma_{\bit+1}$ fulfills the canonical conditions and to investigate the occurrence records of edges of the type $(g_*,F_{*,*})$.
By \Cref{corollary: Improving Switches for new canstrat}, it suffices to prove that $\sigma_{\bit+1}$ has Properties (\ref{property: OR1})$_{*,*,*}$ to (\ref{property: OR4})$_{*,*,*}$ and that \Cref{table: Occurrence Records} specifies the occurrence records of all edges that were not applied during $\canstrat\to\sigma_{\bit+1}$.

We begin by investigating the canonical properties.
The following statement is required when discussing Properties (\ref{property: OR1})$_{*,*,*}$ to (\ref{property: OR4})$_{*,*,*}$.
It states that the occurrence record of the cycle edges of $F_{\ell(\bit+2),1-(\bit+2)}$ are large if $\bit$ is even and will be used repeatedly.
This is useful as the canonical properties with respect to $\sigma_{\bit+1}$ depend on $\bit+2$.
Its proof can be found in \Cref{appendix: Proofs Exponential}.

\begin{restatable}{lemma}{OccRecForBPlusTwo}\label{claim: OccRec for B Plus Two}
Let $\bit\in\bitset_n$ be even, $i\coloneqq\ell(\bit+2)$ and $j\coloneqq 1-(\bit+2)_{i+1}$.
If $\bit+2$ is a power of 2, then $\occrec^{\canstrat}(d_{i,j,*},F_{i,j})=\maxocc.$
Otherwise, $\occrec^{\canstrat}(d_{i,j,0},F_{i,j})=\floor{(\bit+1)/2}$ and $\occrec^{\canstrat}(d_{i,j,1},F_{i,j})=\maxocc-1.$
In any case, $\canstrat(d_{i,j,k})\neq F_{i,j}$ for both $k\in\{0,1\}$.
\end{restatable}

We now prove that $\sigma_{\bit+1}$ has Properties (\ref{property: OR1})$_{*,*,*}$ to (\ref{property: OR4})$_{*,*,*}$.

\begin{lemma} \label{lemma: Next CS has ORs}
Let $\sigma_{\bit+1}$ denote the canonical strategy calculated by the strategy improvement algorithm as described by \Cref{lemma: Extended Reaching canonical strategy}.
Then $\sigma_{\bit+1}$ has Properties (\ref{property: OR1})$_{*,*,*}$ to (\ref{property: OR4})$_{*,*,*}$.
\end{lemma}

\begin{proof}
To simplify notation, let $\sigma\coloneqq\sigma_{\bit+1}$.
We first prove that $\sigma$ has Properties (\ref{property: OR1})$_{*,*,*}$, (\ref{property: OR2})$_{*,*,*}$ and (\ref{property: OR4})$_{*,*,*}$ and discuss \Pref{OR3}$_{*,*,*}$ at the end.

\smallskip
\noindent
\boldall{Consider the case $\nsb>1$ first.}
Let $i\in[n], j,k\in\{0,1\}$ and consider \Pref{OR4}$_{i,j,k}$.
We prove that any improving switch has an occurrence record of either $\maxocc$ or $\maxocc-1$ as $\maxocc=\floor{(\bit+1+1)/2}$ due to $\nsb>1$.
Any $e\in I_\sigma$ was either improving for~$\sigma^{(5)}$ or became improving when transitioning from $\sigma^{(5)}$ to $\sigma$.			
As shown in the proof of \Cref{lemma: Extended Reaching canonical strategy}, all improving switches not applied during phase 5 had an occurrence record of at least $\maxocc-1$.
More precisely, this was shown implicitly when giving the characterization of the improving switches.
Also, the occurrence records of these edges are at most $\maxocc$, proving the statement for these edges.
For improving switches that were unlocked during phase 5, the statement follows by \Cref{corollary: Fulfilling ORFive If NSB Is Not One}.
Hence, $\sigma$ has \Pref{OR4}$_{i,j,k}$.

We prove that $\sigma$ has \Pref{OR2}$_{*,*,*}$ and \Pref{OR1}$_{*,*,*}$.
Consider some indices $i\in[n], j\in\{0,1\}$ with $\indbit_i=0\vee\indbit_{i+1}\neq j$ and let $k\in\{0,1\}$.
We prove \begin{equation} \label{equation: Proving OR properties 1}
 \sigma(d_{i,j,k})=F_{i,j}\Longleftrightarrow\occrec^{\sigma}(d_{i,j,k},F_{i,j})=\ell^{\bit+1}(i,j,k)+1.
 \end{equation}

Let $\sigma(d_{i,j,k})=F_{i,j}$.
Then, since $\sigma^{(5)}(d_{i,j,k})\neq F_{i,j}$ by the choice of $i$ and $j$ and \Cref{table: Properties at start of phase}, the switch was applied during $\sigma^{(5)}\to\sigma$.
Consequently, the edge $(d_{i,j,k},F_{i,j})$ was not applied as improving switch before phase 5 as switches are applied at most once by \Cref{corollary: Improving Switches for new canstrat}.
Thus, $\occrec^{\canstrat}(d_{i,j,k},F_{i,j})=\occrec^{\sigma^{(5)}}(d_{i,j,k},F_{i,j})<\maxocc-1.$
But this implies $\canstrat(d_{i,j,k})=F_{i,j}$ since the switch would have been applied in phase 1 otherwise.
Consequently, by \Cref{lemma: Numerics Of OR}, \[\occrec^{\sigma}(d_{i,j,k},F_{i,j})=\occrec^{\canstrat}(d_{i,j,k},F_{i,j})+1=\ell^{\bit}(i,j,k)+1+1=\ell^{\bit+1}(i,j,k)+1\leq\maxocc-1.\]
This implies both \enquote{$\Rightarrow$} of the equivalence (\ref{equation: Proving OR properties 1}) as well as \Pref{OR1}$_{i,j,k}$.

Now, let $\occrec^{\sigma}(d_{i,j,k},F_{i,j})=\ell^{\bit+1}(i,j,k)+1$.
We prove that this implies $\sigma(d_{i,j,k})=F_{i,j}$.
First, $\occrec^{\sigma}(d_{i,j,k},F_{i,j})=\ell^{\bit+1}(i,j,k)+1\leq\floor{(\bit+1+1-k)/2}$ implies $\ell^{\bit+1}(i,j,k)\leq\floor{(\bit-k)/2}$.
By \Cref{lemma: Numerics Of Ell}, this implies that $\indbit_{i+1}=1-j$.
Consider the case $\bit_i=0\wedge\bit_{i+1}\neq j$.
Then, we have $\occrec^{\canstrat}(d_{i,j,k},F_{i,j})=\min(\floor{(\bit+1-k)/2},\ell^{\bit}(i,j,k)+t_{\bit})$ for some $t_{\bit}$ feasible for~$\bit$.
Assume $\occrec^{\canstrat}(d_{i,j,k},F_{i,j})\neq \ell^{\bit}(i,j,k)+t_{\bit}$ for all feasible parameters and note that this implies $\occrec^{\canstrat}(d_{i,j,k},F_{i,j})=\floor{(\bit+1-k)/2}$.
Then $\occrec^{\sigma}(d_{i,j,k},F_{i,j})<\ell^{\bit}(i,j,k)+1$, implying \[\ell^{\bit+1}(i,j,k)=\ell^{\bit}(i,j,k)+1>\floor{\frac{\bit+1-k}{2}}+1=\floor{\frac{\bit+3-k}{2}}\geq\floor{\frac{\bit+1+1-k}{2}}\]which is a contradiction.
Consequently, $\occrec^{\canstrat}(d_{i,j,k},F_{i,j})=\ell^{\bit}(i,j,k)+t_{\bit}$ for some feasible~$t_{\bit}$.			
Assume $\occrec^{\sigma}(d_{i,j,k},F_{i,j})=\ell^{\bit}(i,j,k)$.
Then  \[\occrec^{\sigma}(d_{i,j,k},F_{i,j})=\ell^{\bit+1}(i,j,k)+1=\ell^{\bit}(i,j,k)+2=\occrec^{\canstrat}(d_{i,j,k},F_{i,j})+2,\] implying that the switch would have been applied twice during $\canstrat\to\sigma$.
This is a contradiction.
The same contradiction follows if we assume $\occrec^{\canstrat}(d_{i,j,k},F_{i,j})=\ell^{\bit}(i,j,k)-1$.
Hence, it holds that $\occrec^{\canstrat}(d_{i,j,k},F_{i,j})=\ell^{\bit}(i,j,k)+1$, implying $\canstrat(d_{i,j,k})=F_{i,j}$.
Since $\ell^{\bit}(i,j,k)=\ell^{\bit+1}(i,j,k)-1$, this also implies that the switch was indeed applied during the transition.
However, $\canstrat(d_{i,j,k})=F_{i,j}$ implies that the switch was not applied during phase~1 of that transition.
But then it must have been applied in phase 5, implying $\sigma(d_{i,j,k})=F_{i,j}$.

We now show that the same holds if $\bit_i=1$ and $\bit_{i+1}=j$, implying $i<\nsb$.
This then yields
\begin{align*}
	\occrec^{\canstrat}(d_{i,j,k},F_{i,j})&=\ceil{\frac{\lastflip{\bit}{i}{\{(i+1,j)\}}+1-k}{2}}=\ceil{\frac{\bit-2^{i-1}+1+1-k}{2}}\\
		&=\ceil{\frac{\bit+1-2^{i-1}+1-k}{2}}=\ell^{\bit+1}(i,j,k).
\end{align*}
Since $\occrec^{\sigma}(d_{i,j,k},F_{i,j})=\ell^{\bit+1}(i,j,k)+1$, this implies that the switch was applied during phase~5 of $\canstrat\to\sigma$.
Consequently, $\sigma(d_{i,j,k})=F_{i,j}$. 
This proves \enquote{$\Leftarrow$} and hence the equivalence (\ref{equation: Proving OR properties 1}).
Most importantly, $\sigma$ thus has \Pref{OR2}$_{i,j,k}$.

\smallskip
\noindent
\boldall{Now assume $\nsb=1$.}
Let $i\in[n],j,k\in\{0,1\}$ and consider \Pref{OR4}$_{i,j,k}$.
We prove that  $e\coloneqq(d_{i,j,k},F_{i,j})\in I_{\sigma}$ implies that $e$ has an occurrence record of $\floor{(\bit+2)/2}-1=\maxocc$ or $\floor{(\bit+2)/2}=\maxocc+1$.
It is easy to verify that for such an edge $e$, one of the following cases holds.
\begin{itemize}
	\item $e\in I_{\sigma'}$ for all $\sigma'\in\reach{\canstrat}$, i.e., the switch was improving during the complete transition.
		Then, $\occrec^{\sigma}(e)=\occrec^{\canstrat}(e)=\maxocc$ by \Cref{corollary: Switches With Low OR In Phase One}.
	\item There is a strategy $\sigma'\in\reach{\sigma^{(5)}}$ with $(d_{i,j,k},F_{i,j})\in I_{\sigma'}$ but $(d_{i,j,k},F_{i,j})\notin I_{\sigma^{(5)}}$.
		That is, the switch became improving during phase~5.
		Then, $\sigma$ has \Pref{OR4}$_{i,j,k}$ by \Cref{corollary: Fulfilling ORFive If NSB Is One}.
	\item The edge $e$ became an improving switch when applying $(b_1,g_1)$ at the end of phase $3$.
		Then $i\in\{u+1,\dots,m-1\}, j=1-\indbit_{i+1}$ and $\indbit_i=0$.
		Thus, by the characterization of~$I_{\sigma}$ given in the beginning of the proof of \Cref{lemma: Extended Reaching canonical strategy}, $\occrec^{\sigma}(d_{i,j,k},F_{i,j})=\maxocc$ since the switch would have been applied during phase $5$ otherwise.
\end{itemize}
Thus, $\sigma$ has \Pref{OR4}$_{i,j,k}$.

Now let $i\in[n]$ and $j\in\{0,1\}$ with $\indbit_i=0\vee\indbit_{i+1}\neq j$, let $k\in\{0,1\}$ and consider \Pref{OR2}$_{i,j,k}$.
Then $\sigma^{(5)}(d_{i,j,k})\neq F_{i,j}$ by \Cref{lemma: Extended Reaching phase 5}.
We again prove that $\sigma$ fulfills the equivalence (\ref{equation: Proving OR properties 1}) and that $\sigma$ has \Pref{OR1}$_{i,j,k}$ simultaneously.		

Let $\sigma(d_{i,j,k})=F_{i,j}$.
By \Cref{lemma: Extended Reaching phase 5}, $(d_{i,j,k},F_{i,j})\in\applied{\sigma^{(5)}}{\sigma}$.
Since improving switches are applied at most once per transition, this implies \[\occrec^{\sigma^{(5)}}(d_{i,j,k},F_{i,j})=\occrec^{\canstrat}(d_{i,j,k},F_{i,j})<\maxocc\] and $\canstrat(d_{i,j,k})=F_{i,j}$ by \Cref{corollary: Cycle Edges Phase Five If NSB One}.
Thus, $\occrec^{\canstrat}(d_{i,j,k},F_{i,j})=\ell^{\bit}(i,j,k)+1=\ell^{\bit+1}(i,j,k)$ by \Pref{OR2}$_{i,j,k}$ and \Cref{lemma: Numerics Of OR}.
Hence \begin{align*}
	\occrec^{\sigma}(d_{i,j,k},F_{i,j})&=\occrec^{\sigma^{(5)}}(d_{i,j,k},F_{i,j})+1=\ell^{\bit+1}(i,j,k)+1<\maxocc+1=\floor{\frac{\bit+2}{2}}
\end{align*} by integrality.
Thus, \enquote{$\Rightarrow$} as well as \Pref{OR1}$_{i,j,k}$ follow.

Let $\occrec^{\sigma}(d_{i,j,k},F_{i,j})=\ell^{\bit+1}(i,j,k)+1$.
By \Cref{lemma: Extended Reaching phase 5}. $\sigma^{(5)}(d_{i,j,k})=F_{i,j}$ if and only if $\indbit_i=1\wedge\indbit_{i+1}=j$.
It thus suffices to prove $(d_{i,j,k},F_{i,j})\in\applied{\sigma^{(5)}}{\sigma}$.
By \Cref{corollary: Cycle Edges Phase Five If NSB One}, we thus need to show that \begin{enumerate}
	\item $\occrec^{\canstrat}(d_{i,j,k},F_{i,j})<\maxocc\wedge\canstrat(d_{i,j,k})=F_{i,j}$,
	\item  $\indbit_i=0\wedge\indbit_{i+1}\neq j$ and 
	\item $i\in\{u+1,\dots,m-1\}$.
\end{enumerate}
Since \begin{equation} \label{equation: Proving OR properties second part equivalence}
\occrec^{\sigma}(d_{i,j,k},F_{i,j})=\ell^{\bit+1}(i,j,k)+1\leq\floor{\frac{\bit+1+1-k}{2}},
\end{equation}
 \Cref{lemma: Numerics Of Ell} implies that $\indbit_i=0\wedge\indbit_{i+1}=1-j$,.
Consequently since $\nsb=1$ implies that no bit switches from 1 to 0,  it follows that $\bit_i=0\wedge\bit_{i+1}=1-j$.
This implies that there is a feasible $t_{\bit}$ with $\occrec^{\canstrat}(d_{i,j,k},F_{i,j})=\min(\floor{(\bit+1-k)/2}, \ell^{\bit}(i,j,k)+t_{\bit})$.
Note that $t_{\bit}\neq -1$ due to the parity of $\bit$ and \Pref{OR3}$_{i,j,k}$.
We prove that $\occrec^{\canstrat}(d_{i,j,k},F_{i,j})=\ell^{\bit}(i,j,k)+1$ by ruling out the other possible cases.
\begin{itemize}
	\item Assume $\occrec^{\canstrat}(d_{i,j,k},F_{i,j})=\floor{(\bit+1-k)/2}$ and that neither $0$ nor $1$ are feasible parameters.
		As $i\neq\nsb$, this implies $\ell^{\bit+1}(i,j,k)=\ell^{\bit}(i,j,k)+1>\floor{(\bit+1-k)/2}$.
		But then $\ell^{\bit+1}(i,j,k)+1>\floor{(\bit+1+1-k)/2},$ contradicting \Cref{equation: Proving OR properties second part equivalence}.
	\item Next assume $\occrec^{\canstrat}(d_{i,j,k},F_{i,j})=\ell^{\bit}(i,j,k)$.
		Then, since $\ell^{\bit}(i,j,k)=\ell^{\bit+1}(i,j,k)-1$, the switch $(d_{i,j,k},F_{i,j})$ would have been switched twice during $\canstrat\to\sigma$.
		This is a contradiction.
\end{itemize}
Hence $\occrec^{\canstrat}(d_{i,j,k},F_{i,j})=\ell^{\bit}(i,j,k)+1$.
It remains to prove $i\in\{u+1,\dots,m-1\}$.
Since $i\geq m$ implies $\ell^{\bit}(i,j,k)\geq\bit$, this implies that we need to have $i<m$ as we have $\occrec^{\canstrat}(d_{i,j,k},F_{i,j})=\ell^{\bit}(i,j,k)+1<\maxocc$.
Also, assuming $i=u$ yields $\occrec^{\canstrat}(d_{i,j,k},F_{i,j})\geq\maxocc$ as discussed earlier.
Consequently, all of the three necessary conditions hold, so \Cref{corollary: Cycle Edges Phase Five If NSB One} implies the direction \enquote{$\Leftarrow$} of the equivalence (\ref{equation: Proving OR properties 1}).
Thus, $\sigma$ has Properties (\ref{property: OR1})$_{*,*,*}$, (\ref{property: OR2})$_{*,*,*}$ and (\ref{property: OR4})$_{*,*,*}$ if $\nsb=1$.

\smallskip
\noindent
\boldall{It remains to prove that $\sigma$ has \Pref{OR3}$_{*,*,*}$.}
By \Pref{OR3}$_{i,j,k}$,   $\occrec^{\sigma}(d_{i,j,k},F_{i,j})=\ell^{\bit+1}(i,j,k)-1\wedge\occrec^{\sigma}(d_{i,j,k},F_{i,j})\neq \floor{(\bit+1+1-k)/2}$ if and only if $\bit+1$ is odd, $\bit+2$ is not a power of 2, $i=\ell(\bit+2)$, $j\neq(\bit+2)_{i+1}$ and $k=0$.
We first prove the \enquote{if} part.
Since $\bit+1$ is odd, $\bit$ is even.
As $\bit+2$ is not a power of~2 by assumption, $\occrec^{\canstrat}(d_{i,j,0},F_{i,j})=\maxocc$ and $\occrec^{\canstrat}(d_{i,j,1},F_{i,j})=\maxocc-1$  as well $\canstrat(d_{i,j,k})\neq F_{i,j}$ for both $k\in\{0,1\}$ by \Cref{claim: OccRec for B Plus Two}.
Consider phase~1 of $\canstrat\to\sigma$.
Then, $(d_{i,j,1},F_{i,j})$ is applied in this phase by \Cref{corollary: Switches With Low OR In Phase One}.
Thus, by the tie breaking rule, $(d_{i,j,0},F_{i,j})$ is not applied during phase~$1$.
Since no switch with an occurrence record of $\maxocc$ is applied during phase 5, the switch is also not applied during phase 5.
Consequently, \[\occrec^{\sigma}(d_{i,j,0},F_{i,j})=\occrec^{\canstrat}(d_{i,j,0},F_{i,j})=\floor{\frac{\bit+1}{2}}=\floor{\frac{\bit+1+1}{2}}-1\] since $\bit+1$ is odd.
It remains to show $\ell^{\bit+1}(i,j,0)=\floor{(\bit+1+1)/2}$.
Since $\bit+1$ is odd, $\ell(\bit+2)\neq\nsb$ and $\bit_i=0$.
Hence, $\ell^{\bit+1}(i,j,0)=\ell^{\bit}(i,j,0)+1=\floor{\bit/2}+1=\floor{(\bit+1+1)/2}$ by \Cref{lemma: Numerics Of OR}.
Thus, the \enquote{if} part is fulfilled.

The \enquote{only if} part can be show using contraposition by dividing the proof into several small statements, each proving that one of the conditions is necessary.
We state all of the statements here and defer their proofs to \Cref{appendix: Proofs Exponential}.
More precisely, the following statements imply the \enquote{only if} part:

\begin{restatable}{claim}{ContrapositionORThree} \label{claim: Contraposition OR Three}
Let $i\in[n], j,k\in\{0,1\}$ and consider the two equations
\begin{align}
	\occrec^{\sigma}(e)&\neq\ell^{\bit+1}(i,j,k)-1, \label{equation: Occrec not ell-1}\\
	\occrec^{\sigma}(e)&=\floor{\frac{\bit+1+1-2}{2}}. \label{equation: Occrec is high}
\end{align}
\begin{enumerate}
	\item If $j=(\bit+2)_{i+1}$, then either \Cref{equation: Occrec not ell-1} or \Cref{equation: Occrec is high} holds.
		
	\item If $i\neq\ell(\bit+2)$ and $j\neq(\bit+2)_{i+1}$, then either \Cref{equation: Occrec not ell-1} or \Cref{equation: Occrec is high} holds.
		
	\item  If $\bit+1$ is even, $i=\ell(\bit+2)$ and $j\neq(\bit+2)_{i+1}$, then \Cref{equation: Occrec is high} holds.
		
	\item If $\bit+1$ is odd, $i=\ell(\bit+2), j= 1-(\bit+2)_{i+1}, k\in\{0,1\}$ and $\bit+2$ is a power of two, then \Cref{equation: Occrec is high} holds.
		
	\item If $\bit$ is even, $i=\ell(\bit+2), j\neq(\bit+2)_{i+1}, k=1$ and $\bit+2$ is not a power of two, then \Cref{equation: Occrec is high} holds.
\end{enumerate}
\end{restatable}

This show that $\sigma$ has \Pref{OR3}$_{*,*,*}$ and thus yields the statement.
\end{proof}

\subsection{Reaching a canonical strategy part II: The occurrence records}

It now remains to prove that \Cref{table: Occurrence Records} specifies the occurrence records with respect to the canonical strategy $\sigma_{\bit+1}$ for $\bit+1$ when it is interpreted for $\bit+1$.
This then implies that $\sigma_{\bit+1}$ has the canonical properties which can then be used to give inductive proofs of the main statements of \Cref{section: Lower Bound Proof}.

As in particular the investigation of edges of the type $(g_*,F_{*,*})$ is rather involved, we show two separate statements and consider all other edges first.

\begin{lemma} \label{lemma: Improving Switches for new canstrat}
Let $\sigma_{\bit+1}$ be the canonical strategy for $\bit+1$ calculated by the strategy improvement resp. policy iteration algorithm when starting with a canonical strategy $\canstrat$ having the canonical properties as described by \Cref{lemma: Extended Reaching canonical strategy}.
Then, \Cref{table: Occurrence Records} specifies the occurrence records of all edges $e\in E_0$ but edges of the type $(g_*,F_{*,*})$ that were applied during $\canstrat\to\sigma_{\bit+1}$.
\end{lemma}

\begin{proof}
There are two types of edges.
Each edge was either applied as improving switch when transitioning from $\canstrat$ to $\sigma_{\bit+1}$ or was not applied as an improving switch.
We already proved that \Cref{table: Occurrence Records} specifies the occurrence records of all improving switch that were applied, with the exception of switches $(g_*,F_{*,*})$.
It thus suffices to consider switches that were not applied when transitioning from $\canstrat$ to $\sigma_{\bit+1}$.
We thus identify edges that were not applied as improving switches and prove that their occurrence records are described by \Cref{table: Occurrence Records}.
To simplify notation, let $\sigma\coloneqq\sigma_{\bit+1}$.

Let $\nsb>1$ and let $i\in[n], j,k\in\{0,1\}$ be suitable indices.
We first prove the statement for all edges that are not of the type $(d_{*,*,*},F_{*,*})$.		
\begin{enumerate}
	\item Consider edges of the type $(b_i,*)$.
		Since $\nsb>1$, the edges $(b_i,b_{i+1})$ for $i\in[\nsb-1]$ as well as the edge  $(b_\nsb,g_{\nsb})$ were applied as improving switches.
		Let $e=(b_i,b_{i+1})$ and $i\geq\nsb$.
		Then $\occrec^{\sigma}(e)=\flips{\bit}{i}{}-\bit_i=\flips{\bit+1}{i}{}-(\bit+1)_{i}$ since either $\flips{\bit}{i}{}=\flips{\bit+1}{i}{}$ and $\bit_i=(\bit+1)_{i+1}$ (if $i>\nsb$) or $\flips{\bit+1}{i}{}=\flips{\bit}{i}{}+1$ and $\bit_i=0, (\bit+1)_i=1$ (if $i=\nsb$).
		Let $e=(b_i,g_i)$ for $i\neq \nsb$.
		Then, by \Cref{lemma: Numerics Of OR}, $\occrec^{\sigma}(e)=\flips{\bit}{i}{}=\flips{\bit+1}{i}{}$.
	\item Consider some edge $(g_i,F_{i,j})$ that was not applied during $\canstrat\to\sigma$.
		Then, the upper bound remains valid as it can only increase.
	\item Consider some vertex $s_{i,j}$.
		Since $\nsb>1$, the edges $(s_{\nsb-1,1},h_{\nsb-1,1}),(s_{\nsb-1,0},b_1)$ as well as the edges $(s_{i,0},h_{i,0}),(s_{i,1},b_1)$ for $i\in[\nsb-2]$ were switched.
		It thus suffices to consider indices $i\geq\nsb$.
		For these edges, the choice of $i$ implies \begin{align*}
			\occrec^{\sigma}(s_{i,j},b_1)&=\flips{\bit}{i+1}{}-j\cdot\bit_{i+1}=\flips{\bit+1}{i+1}{}-j\cdot(\bit+1)_{i+1},\\
			\occrec^{\sigma}(s_{i,j},h_{i,j})&=\flips{\bit}{i+1}{}-(1-j)\bit_{i+1}=\flips{\bit+1}{i+1}{}-(1-j)(\bit+1)_{i+1}.
		\end{align*}
	\item For $e=(e_{i,j,k},g_1)$, \Cref{table: Occurrence Records} implies $\occrec^{\sigma}(e_{i,j,k},g_1)=\ceil{\bit/2}=\ceil{(\bit+1)/2}$ since $\nsb>1$.
	\item For $e=(d_{i,j,k},e_{i,j,k})$, we need to prove $\occrec^{\sigma}(e)\leq\occrec^{\sigma}(e_{i,j,k},g_1)=\ceil{(\bit+1)/2}$ since $\bit$ is odd.
		This follows from $\occrec^{\sigma}(e)\leq\occrec^{\canstrat}(e)+1\leq\floor{\bit/2}+1=\floor{(\bit+2)/2}=\ceil{(\bit+1)/2}.$
\end{enumerate}
Let $i\in[n], j,k\in\{0,1\}$ and consider some edge $e=(d_{i,j,k},F_{i,j})$ that was not switched during $\canstrat\to\sigma$.
We distinguish the following cases.
\begin{enumerate}
	\item Let $(\bit_i=1\wedge\bit_{i+1}=j)$ and $((\bit+1)_{i}=1\wedge(\bit+1)_{i+1}=j)$.
		Then, since any intermediate strategy had \Pref{EV1}$_i$, the cycle  $F_{i,j}$ was always closed during $\canstrat\to\sigma$.
		Thus $i\neq\nsb$, implying $\lastflip{\bit}{i}{\{(i+1,j)\}}=\lastflip{\bit+1}{i}{\{(i+1,j)\}}$.
		Therefore, $\occrec^{\sigma}(e)$ is described by \Cref{table: Occurrence Records}.
	\item Let $(\bit_i=1\wedge\bit_{i+1}=j)$ and $(\bit+1)_{i}=0$, implying $i<\nsb$.
		Then bit $i+1$ also switched, so $(\bit+1)_i=0\wedge(\bit+1)_{i+1}\neq j$.
		Consequently, $e$ was not switched during phase~1 since $F_{i,j}$ was closed with respect to any intermediate strategy due to \Pref{EV1}$_i$.
		It is however possible that such a switch is applied during phase~5.
		Since $i\leq\nsb-1$, this switch is applied if and only if $\occrec^{\canstrat}(e)<\maxocc-1$.
		We may thus assume $\occrec^{\canstrat}(e)\geq\maxocc-1=\floor{(\bit-1)/2}$ and only need to consider the edge~$e$ if $\floor{(\lastflip{\bit}{i}{\{(i+1,j)\}}+2-k)/2}\geq\floor{(\bit-1)/2}.$
		This inequality holds if and only if one of the following three cases applies:
		\begin{itemize}
			\item $\lastflip{\bit}{i}{\{(i+1,j)\}}+2-k\geq\bit-1$.
			\item $\lastflip{\bit}{i}{\{(i+1,j)\}}+2-k$ is even and $\lastflip{\bit}{i}{\{(i+1,j)\}}+2-k=\bit-2$.
			\item $\lastflip{\bit}{i}{\{(i+1,j)\}}+2-k$ is odd and $\lastflip{\bit}{i}{\{(i+1,j)\}}+2-k=\bit$.
		\end{itemize}
		These assumptions can only hold if $i\in\{1,2\}\vee(i=3\wedge k=0)$.
		It thus suffices to consider three more cases.

		For $i=1$, it holds that \[\ell^{\bit+1}(i,j,k)=\ceil{\frac{\bit+1-1+1-k}{2}}=\ceil{\frac{\bit+1-k}{2}}=\occrec^{\sigma}(e).\]

		Similarly, for $i=2$, it holds that \[\ell^{\bit+1}(i,j,k)=\ceil{\frac{\bit+1-2+1-k}{2}}=\ceil{\frac{\bit-k}{2}}=\floor{\frac{\bit+1-k}{2}}=\occrec^{\sigma}(e).\]

		Finally, for $i=3$ and $k=0$,it holds that \[\ell^{\bit+1}(i,j,k)=\ceil{\frac{\bit+1-4+1-k}{2}}=\ceil{\frac{\bit-2-k}{2}}=\occrec^{\sigma}(e).\]

		Hence, the parameter $t_{\bit+1}=0$ can be chosen in all three cases.
	
	\item Let $(\bit_i=0\wedge\bit_{i+1}\neq j)$ and $((\bit+1)_i=0\wedge(\bit+1)_{i+1}\neq j)$, implying $i>\nsb$.
		First let $\id_{j=0}\lastflip{\bit}{i+1}{}+\id_{j=1}\lastunflip{\bit}{i+1}{}=0$.
		Then $\ell^{\bit}(i,j,k)\geq\bit$ by \Cref{lemma: Numerics Of Ell}, implying $\occrec^{\canstrat}(e)=\floor{(\bit+1-k)/2}$.
		Since $\bit$ is odd, $\floor{(\bit+1-1)/2}<\maxocc$.
		Hence, $(d_{i,j,1},F_{i,j})$ was applied during phase 1 of $\canstrat\to\sigma$ and $e=(d_{i,j,0},F_{i,j})\notin\applied{\canstrat}{\sigma}$.
		Thus, since $\ell^{\bit+1}(i,j,k)\geq\bit+1$ by the choice of $i$, choosing $t_{\bit+1}=0$ yields the desired characterization.
		
		Let $\id_{j=0}\lastflip{\bit}{i+1}{}+\id_{j=1}\lastunflip{\bit}{i+1}{}\neq 0$, implying $i<m=\max\{i\colon\indbit_i=1\}$.
		Using $i>\nsb\geq 2$, this yields \begin{align*}
			\ell^{\bit}(i,j,k)&=\ceil{\frac{\bit-2^{i-1}+\sum(\bit,i)+1-k}{2}}\\
				&=\ceil{\frac{\bit-2^{i-1}+\sum_{l=1}^{\nsb-1}2^{l-1}+\sum_{l=\nsb+1}^{i-1}\bit_l2^{l-1}+1-k}{2}}\\
				&\leq\ceil{\frac{\bit-2^{i-1}+2^{\nsb-1}-1+2^{i-1}-2^{\nsb}+1-k}{2}}=\ceil{\frac{\bit-2^{\nsb-1}-k}{2}}\\
				&\leq\ceil{\frac{\bit-2-k}{2}}=\floor{\frac{\bit-1-k}{2}}\leq\floor{\frac{\bit+1-k}{2}}-1.
		\end{align*}
		If $\canstrat(d_{i,j,k})\neq F_{i,j}$, this implies $\occrec^{\canstrat}(d_{i,j,k},F_{i,j})\leq\ell^{\bit}(i,j,k)\leq\floor{(\bit+1-k)/2}-1$.
		Then, by \Cref{corollary: Switches With Low OR In Phase One} the switch was applied during phase 1. 
		We may hence assume $\canstrat(d_{i,j,k})=F_{i,j}$, implying $\occrec^{\canstrat}(e)=\ell^{\bit}(i,j,k)+1\leq\floor{(\bit+1-k)/2}$ as well as $\occrec^{\canstrat}(e)\leq\maxocc-1$ by \Pref{OR1}$_{i,j,k}$.
		As we assume $e\notin\applied{\canstrat}{\sigma}$, it suffices to consider the case $\occrec^{\sigma}(e)=\occrec^{\canstrat}(e)=\floor{(\bit+1)/2}-1$ since $e$ is applied during phase 5 otherwise (see \Cref{corollary: Fulfilling ORFive If NSB Is Not One}).	
		Since $\ell^{\bit+1}(i,j,k)=\ell^{\bit}(i,j,k)+1$ by \Cref{lemma: Numerics Of OR}, choosing $t_{\bit+1}=0$ yields the desired characterization. 

	\item Let $\bit_i=0\wedge\bit_{i+1}\neq j$ and $(\bit+1)_i=1\wedge(\bit+1)_{i+1}\neq j$, so $i=\nsb$.
		The statement follows by the same argument used earlier if $\id_{j=0}\lastflip{\bit}{i+1}{}+\id_{j=1}\lastunflip{\bit}{i+1}{}=0$.
		Hence consider the case $\id_{j=0}\lastflip{\bit}{i+1}{}+\id_{j=1}\lastunflip{\bit}{i+1}{}\neq 0$.
		This implies that we have $\ell^{\bit}(i,j,k)=\floor{(\bit+1-k)/2}.$
		Since $\canstrat$ is a canonical strategy for $\bit$, we have $\canstrat(d_{i,j,k})\neq F_{i,j}$.
		If $\occrec^{\canstrat}(e)=\ell^{\bit}(i,j,k)$, then $\occrec^{\canstrat}(e)=\floor{(\bit+1-k)/2}$ and the same arguments used in the third case can be used to show the statement.
		If $\occrec^{\canstrat}(e)=\ell^{\bit}(i,j,k)-1$, then $\occrec^{\canstrat}(e)=\maxocc-1$ since we need to have  $k=0$ by \Pref{OR3}$_{i,j,k}$.
		But this implies that $e$ was switched during phase 1 and that we do not need to consider it here.
	\item Finally, let $\bit_i=0\wedge\bit_{i+1}=j$.
		It suffices to consider the case $(\bit+1)_i=0\wedge(\bit+1)_{i+1}=j$, implying $i>\nsb$.
		If $\id_{j=0}\lastflip{\bit}{i+1}{}+\id_{j=1}\lastunflip{\bit}{i+1}{}=0$, then the statement follows by the same arguments made earlier.
		Otherwise, we can also use the previous same arguments  since $\ell^{\bit}(i,j,k)>\floor{(\bit+1-k)/2}$ implies $\occrec^{\sigma}(e)=\floor{(\bit+1-k)/2}$.
\end{enumerate}

Now let $\nsb=1$ and $i\in[n],j,k\in\{0,1\}$.
We again begin by proving the statement for all edges that are not of the type $(d_{*,*,*},F_{*,*})$.
\begin{enumerate}
	\item Consider edges of the type $(b_i,*)$.
		Since $\nsb=1$, the only such edge that was applied was $(b_1,g_1)$.
		Let $e=(b_i,g_i)$ and $i\neq 1$.
		Then, $\occrec^{\sigma}(e)=\occrec^{\canstrat}(e)=\flips{\bit}{i}{}=\flips{\bit+1}{i}{}$ by \Cref{table: Occurrence Records} and \Cref{lemma: Numerics Of OR} as required.
		
		For $e=(b_i,b_{i+1})$ and $i\in[n]$, we have $\occrec^{\sigma}(e)=\occrec^{\canstrat}(e)=\flips{\bit}{i}{}-\bit_i.$
		If $i\neq 1$, then $\flips{\bit+1}{i}{}=\flips{\bit}{i}{}$ and $\bit_i=(\bit+1)_i$.
		If $i=1$, then $\flips{\bit+1}{i}{}=\flips{\bit}{i}{}+1$ and $\bit_i=0$ as well as $(\bit+1)_i=1$.
		In both cases, the occurrence record is described by \Cref{table: Occurrence Records}.
	\item Consider some edge $(g_i,F_{i,j})$ that was not applied.
		Then, the upper bound can only increase and thus remains valid.
	\item Consider some vertex $s_{i,j}$.
		Then, since $\nsb=1$, no edge $(s_{*,*},*)$ was switched.
		The statement then follows since $\flips{\bit}{i+1}{}-(1-j)\bit_{i+1}=\flips{\bit+1}{i+1}{}-(1-j)(\bit+1)_{i+1}$ and $\flips{\bit}{i+1}{}-j\cdot\bit_{i+1}=\flips{\bit+1}{i+1}{}-j(\bit+1)_{i+1}.$
	\item Consider some edge $e=(e_{i,j,k},b_2)$.
		Then, the statement follows directly as $\nsb=1$ implies $\occrec^{\sigma}(e)=\occrec^{\canstrat}(e)=\floor{\bit/2}=\maxocc.$
	\item Consider some edge of the type $e=(d_{i,j,k},e_{i,j,k})$ that was not applied.
		Then, $\nsb=1$ implies that $\occrec^{\sigma}(e)=\occrec^{\canstrat}(e)\leq\ceil{\bit/2}=\floor{(\bit+1)/2}=\floor{\bit/2}$.
		The upper bound is thus valid for $\sigma$ since $\occrec^{\canstrat}(e_{i,j,k},b_2)=\occrec^{\sigma}(e_{i,j,k},b_2)$,.		
\end{enumerate}
It remains to consider edges of the type $e=(d_{i,j,k},F_{i,j})$ that were not applied.
As the arguments used for proving this are similar to the ones used for the case $\nsb>1$, we defer the discussion of these edges to \Cref{appendix: Proofs Exponential}.
\begin{restatable}{claim}{OROfNotAppliedEdges}
Let $\nsb=1$.
The occurrence records of edges of the type $(d_{*,*,*}, F_{*,*})$ not applied during $\canstrat\to\sigma$ is described correctly by \Cref{table: Occurrence Records}.
\end{restatable}
\end{proof}

We now investigate the occurrence records of edges of the type $(g_*,F_{*,*})$ that were applied during $\canstrat\to\sigma$.
Determining the exact occurrence records of these edges is challenging as it is challenging to describe the exact conditions under which edges of these type become improving.
In particular, these conditions depend on whether we consider the sink game $S_n$ or the Markov decision process $M_n$, making it even harder to describe these in terms of the unified framework $G_n$.
For these reasons, we prove that the upper bound on the occurrence records of these edges given in \Cref{table: Occurrence Records} remains valid for $\sigma$ by an inductive argument as follows.
We begin by determining the exact set of conditions under which the application of an improving switch $(g_*,F_{*,*})$ might yield an occurrence record that could violate \Cref{table: Occurrence Records}.
That is, we explicitly determine properties that $\sigma$ needs to have for the upper bound to hold with equality.
The idea of the proof is then to show that these conditions imply that there is an earlier canonical strategy $\sigma'$  in which there was a slack between the upper bound of the occurrence record of $(g_*,F_{*,*})$ and the actual occurrence record of this edge.
We then prove that this slack is still present when the strategy $\sigma$ is reached, implying that the upper bound cannot hold with equality and remains valid.

For this reason, the proof itself is an inductive proof that requires us to consider up to 4 previous transitions and uses the statements of this section.
As we always excluded the occurrence records of edges of the type $(g_*,F_{*,*})$ in these statements, we can in fact use them within the following induction.
For example,  we heavily use that each improving switch is applied at most once per transition.
To keep the proof more readable and not refer to one of the previous lemmas in every second sentence, we do not always explicitly mention the lemma proving such a statement.

\begin{lemma} \label{lemma: OR of edges at Selector Vertices}
Let $\canstrat\in\reach{\sigma_0}$ be a canonical strategy for $\bit\in\bitset_n$ calculated by the strategy improvement resp. policy iteration algorithm having the canonical properties.
Then $\occrec^{\canstrat}(g_i,F_{i,j})\leq\min_{k\in\{0,1\}}\occrec^{\canstrat}(d_{i,j,k},F_{i,j})$.
In particular, \Cref{table: Occurrence Records} specifies the occurrence records of all edges of the type $(g_*,F_{*,*})$.
\end{lemma}

\begin{proof}
Let $i\in\{n\}, j\in\{0,1\}$ be fixed and let $e\coloneqq(g_i,F_{i,j})$ be an arbitrary but fixed edge of the type $(g_*,F_{*,*})$.
We prove the statement via induction on $\bit$.
We first consider the case $i\neq 1$ and discuss the case $i=1$ later.

We begin by investigating the first transition in which $e$ could have been applied.
Thus, let $\bit\leq 2^{i}\eqqcolon\mathfrak{t}$.
Then $\mathfrak{t}_{i+1}=1$ and for all $\mathfrak{d}\leq\mathfrak{t}$, it holds that $\mathfrak{d}_{i+1}=0$.
We prove that $e$ was applied at most once when transitioning from $\sigma_0$ to $\sigma_{\mathfrak{t}}$ and that this application can only happen during $\sigma_{\mathfrak{t}-1}\to\sigma_{\mathfrak{t}}$.
The statement then follows since the occurrence records of the cycle edges of $F_{i,j}$ are both at least one.

Since $\sigma_0(g_i)=F_{i,0}$, the switch $e$ cannot have been applied during phase 1 of any transition encountered during the sequence $\sigma_0\to\sigma_{\mathfrak{t}}$ as the choice of $\mathfrak{t}$ implies that there is no $\mathfrak{d}\leq\mathfrak{t}$ with $\mathfrak{d}_i=1\wedge\mathfrak{d}_{i+1}=1$.
It is also easy to show that this implies that it cannot happen that the cycle center $F_{i,j}$ was closed during phase 1 if $j=1-\mathfrak{d}_{i+1}$ as the occurrence record of the cycle edges is too low.
The switch $(g_i,F_{i,j})$ can thus only have been applied during some phase 5.
However, since $\sigma_0(g_i)=F_{i,0}$ and due to the choice of $\mathfrak{t}$, this can only happen when transitioning from $\sigma_{\mathfrak{t}-1}$ to $\sigma_{\mathfrak{t}}$.

Thus, the statement holds for all canonical strategies $\canstrat$ representing numbers $\bit\leq 2^{i}$.
Now, assume that it holds for all $\bit'<\bit$ for an arbitrary but fixed $\bit>2^{i}$.
We prove that the statement also holds for $\canstrat$.
Consider the strategy $\sigma_{\bit-1}$.
We begin by arguing that several cases do not need to be considered.

First of all, every improving switch is applied at most once in a single transition.
If it holds that $\min_{k\in\{0,1\}}\occrec^{\canstrat}(d_{i,j,k},F_{i,j})>\min_{k\in\{0,1\}}\occrec^{\sigma_{\bit-1}}(d_{i,j,k},F_{i,j})$, then the statement thus follows by the induction hypothesis.
We thus assume \begin{equation} \label{equation: Minima are the same}
\min_{k\in\{0,1\}}\occrec^{\canstrat}(d_{i,j,k},F_{i,j})=\min_{k\in\{0,1\}}\occrec^{\sigma_{\bit-1}}(d_{i,j,k},F_{i,j}).
\end{equation}
Similarly, if $e$ is not applied during $\sigma_{\bit-1}\to\canstrat$, then the statement also follows by the induction hypothesis.
We thus assume $e\in\applied{\sigma_{\bit-1}}{\canstrat}$.

These observations give first structural insights on $\bit-1$ and $\bit$.
First, if $\bit_i=1\wedge(\bit-1)_i=1$, then it is not possible to apply $e$ during $\sigma_{\bit-1}\to\canstrat$.
Second, if $\bit_i=1\wedge(\bit-1)_i=0$, then $i=\ell(\bit)$.
By \Cref{definition: Canonical Strategy} resp. \ref{definition: Canonical Strategy MDP}, both cycle centers of level $\ell(\bit)$ are open for $\sigma_{\bit-1}$.
Hence, \Cref{corollary: Selection Vertices In Phase One} implies that $F_{i,j}$ is closed during $\sigma_{\bit-1}\to\canstrat$ by applying both switches $(d_{i,j,0},F_{i,j})$ and $(d_{i,j,1},F_{i,j})$.
But then, Equation~(\ref{equation: Minima are the same}) is not fulfilled and the statement follows.
This implies that it suffices to consider the case $\bit_i=0$.

We now show that these three conditions imply that the occurrence record of $(d_{i,j,*},F_{i,j})$ is \enquote{large}.
To be precise, we prove that Equation~(\ref{equation: Minima are the same}), $e\in\applied{\sigma_{\bit-1}}{\canstrat}$ and $\bit_i=0$ imply \begin{equation}\label{equation: Desired Inequality}
\min_{k\in\{0,1\}}\occrec^{\sigma_{\bit-1}}(d_{i,j,k},F_{i,j})\geq\floor{\frac{\bit}{2}}-1.
\end{equation}
It then suffices to prove $\occrec^{\sigma_{\bit-1}}(g_i,F_{i,j})<\floor{\frac{\bit}{2}}-1$ to complete the proof.

It is easy but tedious to verify that these conditions either already imply the desired inequality or give additional structural insights.

\begin{restatable}{claim}{InequalityisMetOrOneEdgeIsSwitched} \label{claim: Inequality is met or one edge is switched}
Equation~(\ref{equation: Minima are the same}), $e=(g_i,F_{i,j})\in\applied{\sigma_{\bit-1}}{\canstrat}$ and $\bit_i=0$  either imply Inequality~(\ref{equation: Desired Inequality}) directly or that exactly one of the cycle edges of $F_{i,j}$ is switched during $\sigma_{\bit-1}\to\canstrat$.
\end{restatable}

Consequently, it suffices to consider the case that exactly one of the two cycle edges $(d_{i,j,*},F_{i,j})$ is applied during $\sigma_{\bit-1}\to\canstrat$.
However, by Equation~(\ref{equation: Minima are the same}), this implies that the the occurrence record of both edges $(d_{i,j,0},F_{i,j}), (d_{i,j,1}, F_{i,j})$ is the same with respect to $\sigma_{\bit-1}$.

Now assume that $F_{i,j}$ is open or halfopen for $\sigma_{\bit-1}$.
Then, there is at least one $k\in\{0,1\}$ with $\sigma_{\bit-1}(d_{i,j,k})\neq F_{i,j}$.
The statement thus follows since \Pref{OR4}$_{i,j,k}$ implies $\occrec^{\sigma_{\bit-1}}(d_{i,j,k},F_{i,j})\geq\floor{\bit/2}-1$.
Thus assume that $F_{i,j}$ is closed for $\sigma_{\bit-1}$.
This implies that either $(\bit-1)_i=1\wedge(\bit-1)_{i+1}=j$ or $(\bit-1)_i=0\wedge(\bit-1)_{i+1}\neq j$ holds.
In the first case, $\occrec^{\sigma_{\bit-1}}(d_{i,j,k},F_{i,j})=\floor{(\bit-2^{i-1}+2-k)/2}$  and $\ell(\bit)>1$ need to hold.
This implies that $\bit-2^{i-1}+2$ is even as we have $i\neq 1$ by assumption.
But this implies that we have $\occrec^{\sigma_{\bit-1}}(d_{i,j,1},F_{i,j})<\occrec^{\sigma_{\bit-1}}(d_{i,j,0},F_{i,j})$ which is a contradiction.

We thus need to have $(\bit-1)_i=0\wedge(\bit-1)_{i+1}\neq j$.
Then, since $F_{i,j}$ is closed, \Pref{OR2}$_{i,j,*}$ implies $\ell^{\bit-1}(i,j,0)=\ell^{\bit-1}(i,j,1)$ since \[\occrec^{\sigma_{\bit-1}}(d_{i,j,0},F_{i,j})=\ell^{\bit-1}(i,j,0)+1=\ell^{\bit-1}(i,j,1)+1=\occrec^{\sigma_{\bit-1}}(d_{i,j,1},F_{i,j}).\]
But this implies $\ceil{(\bit-1-2^{i-1}+\sum(\bit-1,i))/2}=\ceil{(\bit-2^{i-1}+\sum(\bit-1,i))/2}.$
Since $\bit-2^{i-1}+\sum(\bit-1,i)$ is always odd due to $i\neq 1$, this is however not possible.

This concludes this part of the proof. 
Since $\min_{k\in\{0,1\}}\occrec^{\sigma_{\bit-1}}(d_{i,j,k},F_{i,j})\geq\floor{\bit/2}-1,$ it suffices to prove $\occrec^{\sigma_{\bit-1}}(g_i,F_{i,j})<\floor{\bit/2}-1$ under the given conditions.

We begin by stating one more structural insight.
\begin{restatable}{claim}{NeedToHaveZero} \label{claim: Need to have 0}
Assume that Equation~(\ref{equation: Minima are the same}), $e=(g_i,F_{i,j})\in\applied{\sigma_{\bit-1}}{\canstrat}, \bit_i=0$ hold and  that exactly one of the two cycle edges $(d_{i,j,0},F_{i,j}), (d_{i,j,1}, F_{i,j})$ is applied during $\sigma_{\bit-1}\to\canstrat$.
Then $(\bit-1)_i=0$.
\end{restatable}

To simplify notation, we denote the binary number obtained by subtracting 1 from a binary number $(\bit'_n,\dots,\bit'_1)$ by $[\bit'_{n},\dots,\bit'_{1}]-1$.
Then, $\bit$ and $\bit-1$ can be represented as \begin{align*}
\bit&=(\bit_n,\dots,\bit_{i+1},\bm{0},\bit_{i-1},\dots,\bit_1), \\
\bit-1&=(\bit_n,\dots,\bit_{i+1},\bm{0},[\bit_{i-1},\dots,\bit_1]-1)
\end{align*} where bit $i$ is marked in \textbf{bold}.
The idea of the proof is now the following.
We define two smaller numbers that are relevant for the application of $(g_i,F_{i,j})$.
We use these numbers and the induction hypothesis to prove that even if $(g_i,F_{i,j})$ was applied during the maximum number of transitions, the claimed bound still holds.

We define \begin{align*}
	\bar{\bit}\coloneqq([\bit_n,\dots,\bit_{i+1}]-1,&\bm{1},1\dots,1) \\
	\widehat{\bit}\coloneqq([\bit_n,\dots,\bit_{i+1}]-1,&\bm{1},0,\dots,0)
\end{align*} where bit $i$ is again marked in \textbf{bold}.
These numbers are well-defined since $\bit\geq2^{i}$.

Consider $\widehat{\bit}$.
Let $\mathfrak{N}(\widehat{\bit}, \bit-1)$ denote the number of applications of $(g_i,F_{i,j})$ when transitioning from $\sigma_{\widehat{\bit}}$ to $\sigma_{\bit-1}$.
Then, since $\bit'_i=1$ for all $\bit'\in\{\widehat{\bit},\dots,\bar{\bit}\}$, we have $\mathfrak{N}(\widehat{\bit},\bit-1)=\mathfrak{N}(\bar{\bit}, \bit-1)$.
We thus can describe the occurrence record of $(g_i,F_{i,j})$ as \[\occrec^{\sigma_{\bit-1}}(g_i,F_{i,j})=\mathfrak{N}(0,\bit-1)=\mathfrak{N}(0,\widehat{\bit})+\mathfrak{N}(\widehat{\bit},\bit-1)=\occrec^{\sigma_{\widehat{\bit}}}(g_i,F_{i,j})+\mathfrak{N}(\bar{\bit},\bit-1).\]

Our goal is to bound the two terms on the right-hand side.
Using the induction hypothesis and that $\widehat{\bit}$ is even, it is easy to verify that the first term can  be bounded by $\lfloor\widehat{\bit}/2\rfloor$.
Since every improving switch is applied at most once per transition by \Cref{corollary: Improving Switches for new canstrat}, we have $\mathfrak{N}(\bar{\bit},\bit-1)\leq (\bit-1)-\bar{\bit}$.
However, this upper bound is not strong enough.
To improve this bound, we now distinguish between when exactly $(g_i,F_{i,j})$ is applied during $\sigma_{\bit-1}\to\canstrat$.
Note that we refer to even earlier transitions in the last statement of the following claim.

\begin{restatable}{claim}{SelectorSwitchAppliedInPhaseOne} \label{claim: Selector Switch Applied in Phase One}
Assume that Equation~(\ref{equation: Minima are the same}), $e=(g_i,F_{i,j})\in\applied{\sigma_{\bit-1}}{\canstrat}, \bit_i=0$ hold and  that exactly one of the two cycle edges $(d_{i,j,0},F_{i,j}), (d_{i,j,1}, F_{i,j})$ is applied during $\sigma_{\bit-1}\to\canstrat$.
If $(g_i,F_{i,j})$ is applied during phase 1 of $\sigma_{\bit-1}\to\canstrat$, then
\begin{enumerate}
	\item $\bit$ is even and $i\neq 2$,
	\item $\sum(\bit,i)=2^{i-1}-2$ and
	\item  if $(g_i,F_{i,j})\in\applied{\sigma_{\bit-2}}{\sigma_{\bit-1}}$, then $(g_i,F_{i,j})\notin\applied{\sigma_{\bit-3}}{\sigma_{\bit-2}}$.
\end{enumerate}
\end{restatable}


If $(g_i,F_{i,j})$ was applied during phase~1, then the last statement of \Cref{claim: Selector Switch Applied in Phase One} implies $\mathfrak{N}(\bar{\bit},\bit-1)\leq (\bit-1)-\bar{\bit}-1$.
Combining these results and using $\bar{\bit}=\bit-\sum(\bit,i)-1$ and $\widehat{\bit}=\bit-\sum(\bit,i)-2^{i-1}$ yields the statement as \begin{align*}
	\occrec^{\sigma_{\bit-1}}(g_i,F_{i,j})&=\occrec^{\sigma_{\widehat{\bit}}}(g_i,F_{i,j})+\mathfrak{N}(\bar{\bit},\bit-1)\\
		&\leq\floor{\frac{\widehat{\bit}}{2}}+(\bit-1)-\bar{\bit}-1=\floor{\frac{\widehat{\bit}}{2}}+\bit-\bar{\bit}-2\\
		&=\floor{\frac{\bit-\sum(\bit,i)-2^{i-1}}{2}}+\bit-\bit+\sum(\bit,i)+1-2\\
		&=\floor{\frac{\bit-2^{i-1}+\sum(\bit,i)}{2}}-1=\floor{\frac{\bit-2^{i-1}+2^{i-1}-2}{2}}-1\\
		&=\floor{\frac{\bit-2}{2}}-1=\floor{\frac{\bit}{2}}-2<\floor{\frac{\bit}{2}}-1.
\end{align*}
This proves the statement if $(g_i,F_{i,j})$ was applied during phase 1 of $\sigma_{\bit-1}\to\canstrat$.

Hence assume that $e=(g_i,F_{i,j})$ was applied during phase 5 of $\sigma_{\bit-1}\to\canstrat$.
Let $\sigma$ denote the phase-5-strategy in which $e$ is applied.
Then, $\sigmabar(g_i)=1-j$ needs to hold.
Consequently, either $\sigmabar_{\bit-1}(g_i)=1-j$ or $\sigmabar_{\bit-1}(g_i)=j$ and $(g_i,F_{i,1-j})\in\applied{\sigma_{\bit-1}}{\sigma}$.
We thus distinguish between these two cases.
In the first case, the following statement similar to \Cref{claim: Selector Switch Applied in Phase One} can be shown.

\begin{restatable}{claim}{SelectorSwitchAppliedInPhaseFive} \label{claim: Selector Switch Applied in Phase Five}
Assume that Equation~(\ref{equation: Minima are the same}), $e=(g_i,F_{i,j})\in\applied{\sigma_{\bit-1}}{\canstrat}, \bit_i=0$ hold and  that exactly one of the two cycle edges $(d_{i,j,0},F_{i,j}), (d_{i,j,1}, F_{i,j})$ is applied during $\sigma_{\bit-1}\to\canstrat$.
If $(g_i,F_{i,j})$ is applied during phase 5 of $\sigma_{\bit-1}\to\canstrat$ and $\sigmabar_{\bit-1}(g_i)=1-j$, then
\begin{enumerate}
	\item $i\neq 2$,
	\item $\sum(\bit,i)=2^{i-1}-2$ and
	\item  if $(g_i,F_{i,j})\in\applied{\sigma_{\bit-2}}{\sigma_{\bit-1}}$, then$(g_i,F_{i,j})\notin\applied{\sigma_{\bit-3}}{\sigma_{\bit-2}}$.
\end{enumerate}
\end{restatable}

The statement thus follows analogously.

Thus, assume $\sigmabar_{\bit-1}(g_i)=j$ and $(g_i,F_{i,1-j})\in\applied{\sigma_{\bit-1}}{\canstrat}$.
Since at most one improving switch involving a selection vertex is applied during phase 5, this implies that $(g_i,F_{i,1-j})$ was applied during phase 1 of $\sigma_{\bit-1}\to\canstrat$.
It could technically also be applied at the beginning of phase $2$ resp. $3$ when closing the final cycle center only creates a pseudo phase-$2$ resp. pseudo phase-$3$-strategy.
For clarity of presentation, we include this case in the second case and interpret the application of this improving switch as being part of phase $1$.
In particular, we thus have $1-j=1-\bit_{i+1}$ resp. $j=\bit_{i+1}$ as $(\bit-1)_i=\bit_i=0$ implies $i\neq\ell(\bit)$.
We prove that we need to have $i\neq 2$.

For the sake of a contradiction, assume $i=2$.
Then, since $(\bit-1)_i=0$, we have $\bit_1=1$ and $\bit-2=\bar{\bit}$.
Consequently, $1-j=1-\bit_{3}=(\bit-2)_3$.
Thus, the cycle center $F_{i,1-j}$ is active and closed with respect to $\sigma_{\bit-2}$.
As $\bit$ is odd, this implies \begin{align*}
	\occrec^{\sigma_{\bit-2}}(d_{i,1-j,k},F_{i,1-j})&=\floor{\frac{\lastflip{\bit-2}{i}{\{(i+1,(\bit-2)_3)\}}-k}{2}}+1\\
		&=\floor{\frac{\bit-2-1-k}{3}}+1=\floor{\frac{\bit-1-k}{2}}=\floor{\frac{\bit-1}{2}}-k.
\end{align*}
Since the cycle center is closed with respect to $\sigma_{\bit-2}$, none of these two edges is applied as improving switch during phase 1 of $\sigma_{\bit-2}\to\sigma_{\bit-1}$.
However, since $\ell(\bit-1)>1$, the switches are also not applied during phase 5 of that transition.
But this implies $\sigma_{\bit-1}(d_{i,1-j,k})\neq F_{i,1-j}$ for both $k\in\{0,1\}$, contradicting that $(g_i,F_{i,1-j})$ is applied during phase 1 of $\sigma_{\bit-1}\to\sigma_{\bit}$.
Note that this argument further proves that we cannot have $\bit-2=\bar{\bit}$.

We can thus assume $i>2$ and $\bit-3\geq\bar{\bit}$.
Since we apply $(g_i, F_{i,1-j})$ during phase 1 of $\sigma_{\bit-1}\to\canstrat$, we can use the same arguments used when proving \Cref{claim: Selector Switch Applied in Phase One} resp. \ref{claim: Selector Switch Applied in Phase Five} to prove $\sum(\bit,i)\leq 2^{i-1}-2$.
Similar to the last cases, we prove that there is at least one transition between $\bar{\bit}$ and $\bit-1$ in which the switch $(g_i,F_{i,j})$ is not applied.
As this follows if $(g_i,F_{i,j})\notin\applied{\sigma_{\bit-2}}{\sigma_{\bit-1}}$, assume $(g_i,F_{i,j})\in\applied{\sigma_{\bit-2}}{\sigma_{\bit-1}}$.

First, since $\bit-3\geq\bar{\bit}$ and since we apply $(g_i,F_{i,j})$ in phase 5 of $\sigma_{\bit-1}\to\canstrat$, it holds that $j=\bit_{i+1}=(\bit-1)_{i+1}=(\bit-2)_{i+1}$.
This implies that $(g_i,F_{i,j})$ was not applied during phase~1 of $\sigma_{\bit-3}\to\sigma_{\bit-2}$.
The reason is that this could only happen if $i=\ell(\bit-2)$, contradicting $(\bit-2)_i=0$, or if $(\bit-3)_i=0\wedge j\neq(\bit-3)_{i+1}$.
However, this then contradicts the previous identities regarding $j$.
Thus, $(g_i,F_{i,j})$ is applied during phase 5 of $\sigma_{\bit-3}\to\sigma_{\bit-2}$.

For the sake of a contradiction, assume $(\bit-3)_i=1$.
Then, $\bit-3=\bar{\bit}$, implying $(\bit-3)_{i+1}\neq j$.
This further implies $(\bit-3)_{i'}=1$ for all $i'\leq i$.
Then, since $F_{i,1-j}$ is closed with respect to $\sigma_{\bit-3}$ and $i\geq 3$, it is easy to calculate that we then need to have $\occrec^{\sigma_{\bit-3}}(d_{i,1-j,k},F_{i,1-j})<\floor{(\bit-3+1)/2}-1$.
But, since $\ell(\bit-2)>1$, this implies that both of these edges are applied at the beginning of phase~5 of $\sigma_{\bit-3}\to\sigma_{\bit-2}$.
Thus, $F_{i,1-j}$ is closed at the beginning of phase~5 of $\sigma_{\bit-3}\to\sigma_{\bit-2}$, contradicting the assumption that $(g_i,F_{i,j})$ is applied during phase 5 of $\sigma_{\bit-3}\to\sigma_{\bit-2}$, see \Cref{lemma: Phase Five Escape Easy}.

Thus assume $(\bit-3)_i=0\wedge(\bit-3)_{i+1}=j$.
This implies $\bit-4\geq\bar{\bit}$ and that the transition from $\sigma_{\bit-4}$ to $\sigma_{\bit-3}$ is thus part of the currently considered sequence of transitions.
Then, since $(g_i,F_{i,j})$ is applied during both $\sigma_{\bit-3}\to\sigma_{\bit-2}$ and $\sigma_{\bit-2}\to\sigma_{\bit-1}$, the improving switch $(g_i,F_{i,1-j})$ has to be applied in between.
This switch can only be applied during phase 1 of $\sigma_{\bit-2}\to\sigma_{\bit-1}$.
It is easy to see that this implies that there is a $k\in\{0,1\}$ such that $\occrec^{\sigma_{\bit-2}}(d_{i,1-j,k},F_{i,1-j})=\ell^{\bit-2}(i,1-j,k)+1\leq\floor{(\bit-1)/2}-1$ and $\occrec^{\sigma_{\bit-2}}(d_{i,1-j,1-k},F_{i,1-j})=\floor{(\bit-1)/2}-1$.
Using $\sum(\bit,i)\leq2^{i-1}-2$, it is then easy to verify that $\ell^{\bit-3}(i,1-j,k')\leq\floor{(\bit-k')/2}-3$  for $k'\in\{0,1\}$.
%

This implies that we need to have $\occrec^{\sigma_{\bit-3}}(d_{i,1-j,1},F_{i,1-j})=\ell^{\bit-3}(i,1-j,1)+1$ and that the edge $(d_{i,1-j,1}, F_{i,1-j})$ is applied as improving switch during phase 5 of $\sigma_{\bit-3}\to\sigma_{\bit-2}$.
We thus need to have $k=1$.
Consider $(d_{i,1-j,0}, F_{i,1-j})$.
Then, \[\ell^{\bit-3}(i,1-j,0)=\floor{\frac{\bit}{2}}-3=\floor{\frac{\bit-2}{2}}-2<\floor{\frac{\bit-2}{2}}-1=\floor{\frac{\bit-3+1}{2}}-1.\]
By \Pref{OR4}$_{i,1-j,0}$, we thus need to have $\sigma_{\bit-3}(d_{i,1-j,0})=F_{i,1-j}$.
In particular, it implies that $F_{i,1-j}$ is closed for $\sigma_{\bit-3}$ and thus $\sigma_{\bit-3}(g_i)=F_{i,1-j}$.

This now enables us to show that the edge $(g_i,F_{i,j})$ was not applied as improving switch during the transition $\sigma_{\bit-4}\to\sigma_{\bit-3}$.
Independent on whether $\bit-4=\bar{\bit}$ or $\bit-4\neq\bar{\bit}$, the switch was not applied during phase $5$ of $\sigma_{\bit-4}\to\sigma_{\bit-3}$ as we have $\sigma_{\bit-3}(g_i)=F_{i,1-j}$.
If $(\bit-4)_i=0\wedge (\bit-4)_{i+1}=j$, then it also follows directly that the switch was not applied during phase $1$ of that transition.
Thus consider the case $(\bit-4)_i=1 \wedge(\bit-4)_{i+1}\neq j$, implying $\bit-4=\bar{\bit}$.
But this immediately implies that the switch is not applied during phase 1.

This concludes the proof for the case that $(g_i,F_{i,j})$ was applied during phase 5 of $\sigma_{\bit-1}\to\canstrat$ and thus concludes the proof for $i\geq 2$.

It remains to consider the case $i=1$.
We prove the statement again via induction on~$\bit$.
It is easy to verify that the statement hold for both $G_n=S_n$ and $G_n=M_n$ for the canonical strategies $\sigma_0,\sigma_1,\sigma_2$.
Hence let $\bit>2$ and assume that the statement holds for all $\bit'<\bit$.
We show that the statement then also holds for $\bit$.

An improving switch $(g_1,F_{1,j})$ can only be applied during phase 1 if $1=\ell(\bit)$ and $j=\bit_2$.
Since we then have $(\bit-1)_i=0$, both edges $(d_{1,\bit_2,0},F_{1,\bit_2})$ and $(d_{1,\bit_2,1}, F_{1,\bit_2})$ are switched during phase 1 of $\sigma_{\bit-1}\to\canstrat$.
Thus, the statement follows by the induction hypothesis.

Thus consider the case $(\bit-1)_i=1$.
Then, a switch $(g_1,F_{1,j})$ can only be applied in phase 5.
Consider the case $G_n=S_n$ first.
Then, by the conditions of the application of such a switch in phase 5, we need to have $j=0$ and $\sigma_{\bit-1}(g_1)=1$.
Since $(\bit-1)_1=1$ this implies $\bit=(\dots,0,0), \bit-1=(\dots,1,1)$ and $\bit-2=(\dots,1,0)$.
It follows directly that $(g_1,F_{1,0})$ is not applied during the transition $\sigma_{\bit-2}\to\sigma_{\bit-1}$.
However, during phase 1 of both transitions $\sigma_{\bit-2}\to\sigma_{\bit-1}$, exactly one of the cycle edges of $F_{1,0}$ is switched.
Using the induction hypothesis, this implies the statement.
If $G_n=M_n$, then we then need to have $j=\bit_2$ and $\sigma_{\bit-1}(g_1)=1-\bit_2$.
If $j=0$, then the statement follows by the exact same arguments used for the case $G_n=S_n$.
If $j=1$, it follows by similar arguments.

\end{proof}

We can now prove the statements of \Cref{section: Lower Bound Proof}.
For convenience, we restate these statements before proving them.

We begin by showing that applying the improving switches according to Zadeh's pivot rule and the tie-breaking rule of \Cref{definition: Tie-Breaking exponential} yields the strategies described by \Cref{table: Properties at start of phase,table: Switches at start of phase}.

\ReachingNextPhase*

\begin{proof}
By \Cref{lemma: Extended Reaching canonical strategy}, applying improving switches to $\sigma$ produces a canonical strategy $\sigma_{\bit+1}$ for $\bit+1$.
When proving this lemma, we proved that the algorithm produces the intermediate strategies as described by the corresponding rows of \Cref{table: Properties at start of phase,table: Switches at start of phase}.
More precisely, this follows from  \Cref{lemma: Extended Reaching phase 2,,lemma: Extended Reaching phase 3,,lemma: Extended Reaching phase 4 or 5,,lemma: Extended Reaching phase 5}.

By \Cref{lemma: Extended Reaching canonical strategy}, $I_{\sigma_{\bit+1}}=\{(d_{i,j,k},F_{i,j}):\sigma_{\bit+1}\neq F_{i,j}\}$.
In particular, this set is described as specified by \Cref{table: Switches at start of phase}.
It thus remains to prove that $\sigma_{\bit+1}$ has the canonical properties.
More precisely, we prove the following three statements:
\begin{enumerate}
	 \item The occurrence records $\occrec^{\sigma_{\bit+1}}$ are described correctly by \Cref{table: Occurrence Records}:
	 	This follows from \Cref{lemma: Improving Switches for new canstrat,lemma: OR of edges at Selector Vertices}.
	 \item $\sigma_{\bit+1}$ has Properties (\ref{property: OR1})$_{*,*,*}$ to (\ref{property: OR4})$_{*,*,*}$:
	 	This follows from \Cref{lemma: Next CS has ORs}.
	 \item Any improving switch was applied at most once per previous transition between canonical strategies:
	 	This follows from \Cref{corollary: Improving Switches for new canstrat}.
\end{enumerate}
Thus, $\sigma_{\bit+1}$ is a canonical strategy for $\bit+1$ and has the canonical properties.
\end{proof}

This now immediately implies the following theorem of \Cref{section: Lower Bound Proof}, stating that applying improving switches to a canonical strategy for $\bit$  having the canonical properties produces such a strategy for $\bit+1$.

\ReachingCanStrat*

This now enables us to prove the remaining statements of \Cref{section: Lower Bound Proof} simultaneously.

\TableDescribesOR*

\MainTheorem*

\begin{proof}
By \Cref{lemma: Properties of initial strategy}, the initial strategy $\sigma_0$ is a canonical strategy representing $0$ having the canonical properties.
In addition, by \Cref{lemma: Sink Game}, it is a sink strategy for $S_n$ and a weak unichain policy for $M_n$.

By \Cref{theorem: Reaching canonical strategy}, applying improving switches to $\sigma_0$ yields a canonical strategy $\sigma_1$ representing the number $1$.
Also, the occurrence record of the edges is described correctly by \Cref{table: Occurrence Records} for $\sigma_1$.
In particular, \Cref{theorem: Reaching canonical strategy} can be applied to $\sigma_1$ again, yielding a canonical strategy $\sigma_2$ representing the number $2$.

This argument can now be applied iteratively.
Thus, applying improving switches according to Zadeh's pivot rule and the tie-breaking rule defined in \Cref{definition: Tie-Breaking exponential} produces the strategies $\sigma_0,\sigma_1,\dots,\sigma_{2^n-1}$.
By \Cref{theorem: Reaching canonical strategy}, \Cref{table: Occurrence Records} describes the occurrence records for each of these strategies, implying \Cref{theorem: Table Describes OR}. 
Since these are $2^n$ different strategies, since $G_n$ has size $\mathcal{O}(n)$ and all rewards and probabilities (for $G_n=M_n$) and priorities (for $G_n=S_n$) can be encoded using a polynomial number of bits, this implies the exponential lower bound for the strategy improvement and policy iteration algorithm.
Consequently, there is a linear program such that the simplex algorithm using the same pivot and tie-breaking rule requires the same number of iterations as it requires for $M_n$ and the lower bound also applies to the simplex algorithm.
\end{proof}
%
%
%

%
%

\newpage

\section{Formal Proofs} \label{appendix: Proofs Exponential}

Thisappendix contains all proofs of statements contained in the previous appendices.
We begin by providing proofs for the statements of \Cref{section: Vertex Valuations}.

\subsection{Omitted proofs of \Cref{section: Vertex Valuations}}

\BooundingValuesForMDPs*

\begin{proof}
The highest priority of any vertex is $2n+10$ and there are no more then $5n$ vertices with priorities.
Let $v,w$ be two vertices with $\Omega(v)>\Omega(w)$.
Then, by construction, $\rew{v}\geq N\cdot \rew{w}$.
In other words, if two vertices $v,w$ do not have the same priority, then the rewards associated with the vertices are apart by at least a factor of $N$.
Thus, for $S\subseteq P$, \[\abs{\sum(S)}\leq\abs{\sum(P)}\leq\abs{P}\cdot\abs{\max_{v\in P}\rew{v}}\leq \abs{P}\cdot N^{2n+10}=5n\cdot N^{2n+10}< N^{2n+11},\] implying the first statement since $\e=N^{-(2n+11)}$ by definition.

Let $S,S'\subseteq P$ be non-empty.
Let $\abs{\max_{v\in S}\rew{v}}<\abs{\max{_{v\in S'}}\rew{v}}$.
Then \[\abs{\sum(S)}\leq\abs{S}\cdot\abs{\max_{v\in S}\rew{v}}<5n\abs{\max_{v\in S}\rew{v}}\leq 5n\frac{\abs{\max_{v\in S'}\rew{v}}}{N}< \abs{\max_{v\in S'}\rew{v}}\leq\abs{\sum(S')}.\]
Now let $\abs{\sum(S)}<\abs{\sum(S')}$.
Then \begin{align*}
\abs{\max_{v\in S}\rew{v}}\leq \abs{\sum(S)}<\abs{\sum(S')}<5n\abs{\max_{v\in S'}\rew{v}},
\end{align*}
so $\abs{\max_{v\in S}\rew{v}}<5n\abs{\max_{v\in S'}\rew{v}}$.
Since $N\geq 7n$, this implies the statement.
\end{proof}

\VVLemma*

\begin{proof}
We prove the statements one after another.
\begin{enumerate}
	\item The first statement follows directly if $i\geq\relbit{\sigma}$ since this implies $R_i^*=L_i^*$.
		Thus assume $i<\relbit{\sigma}$.
		Then, $\sigma(b_{\relbit{\sigma}})=b_{\relbit{\sigma}+1}$ implies $L_i^*=\bigoplus_{\ell=i}^{\relbit{\sigma}-1}\{W_{\ell}^*\colon\sigma(b_{\ell})=g_{\ell}\}\oplus L_{\relbit{\sigma}+1}^*.$
		The first statement follows since \[\bigoplus_{\ell=i}^{\relbit{\sigma}-1}\{W_{\ell}^*\colon\sigma(b_{\ell})=g_{\ell}\}\preceq\bigoplus_{\ell=i}^{\relbit{\sigma}-1}W_{\ell}^*\quad\text{and}\quad R_i^*=\bigoplus_{\ell=i}^{\relbit{\sigma}-1}W_{\ell}^*\oplus L_{\relbit{\sigma}+1}.\]
		The second statement follows since $j<i\leq\relbit{\sigma}$ implies $\bigoplus_{\ell=i}^{\relbit{\sigma}-1}W_{\ell}^*\prec\bigoplus_{\ell=j}^{\relbit{\sigma}-1}W_{\ell}^*$.
	\item The first statement follows directly if $i>\relbit{\sigma}$ since this implies $R_i^*=L_i^*$.
		Thus assume $i\leq\relbit{\sigma}$.
		Since $\sigma(b_{\relbit{\sigma}})=g_{\relbit{\sigma}}$, the statement then follows since \begin{align*}
			L_i^*&=\bigoplus_{\ell\geq i}\{W_{\ell}^*\colon\sigma(b_{\ell})=g_{\ell}\}=\bigoplus_{\ell=i}^{\relbit{\sigma}-1}\{W_{\ell}^*\colon\sigma(b_{\ell})=g_{\ell}\}\oplus W_{\relbit{\sigma}}\oplus L_{\relbit{\sigma}+1}^*\\
				&\succ\bigoplus_{\ell=1}^{\relbit{\sigma}-1}W_{\ell}^*\oplus L_{\relbit{\sigma}+1}^*=R_1^*\succeq R_j^*.
		\end{align*}
		The same calculation implies the third statement as the estimations remain correct if $j<\relbit{\sigma}$.
	\item By $i\geq\relbit{\sigma}>j$, we have \begin{align*}
		R_j^*&=\bigoplus_{\ell=j}^{\relbit{\sigma}-1}W_{\ell}\oplus\bigoplus_{\ell\geq\relbit{\sigma}+1}\{W_{\ell}\colon\sigma(b_{\ell})=g_{\ell}\}\prec\bigoplus_{\ell=1}^{\relbit{\sigma}}W_{\ell}\oplus\bigoplus_{\ell\geq\relbit{\sigma}+1}\{W_{\ell}\colon\sigma(b_{\ell})=g_{\ell}\}\\
			&\prec\bigoplus_{\ell=1}^iW_{\ell}\oplus\bigoplus_{\ell\geq i+1}\{W_{\ell}\colon\sigma(b_{\ell})=g_{\ell}\}\prec\ubracket{s_{i,j},h_{i,j}}\oplus L_{i+1}^*.
	\end{align*}
	\item For $G_n=S_n$, the statement follows since the most significant difference is the vertex $h_{i*}$ and since the priority of this vertex is even.
		For $G_n=M_n$, this follows intuitively since the $\rew{h_{i,*}}$ has the largest exponent of all terms in the expression and since it is the only vertex that has this exponent.
		Thus, $\rew{h_{i,*}}$ is by a factor of $N$ larger than all other quantities in the given expression, and as $N$ is sufficiently large, this term dominates.
		Formally, this can be shown by an easy but tedious calculation. 
		
	  The second statement follows from the first. \qedhere
\end{enumerate}
\end{proof}

\ValuationOfGIfLevelSmall*

\begin{proof}
This statement is shown by backwards induction on $i$.
Let $i=\relbit{\sigma}-1$, implying $\sigma(g_i)=F_{i,0}$ by \Pref{BR1}.
\begin{itemize}
	\item Let $\sigmabar(d_i)$.
		Since $\sigma(s_{i,0})=h_{i,0}$ by \Pref{S2}, \Cref{lemma: Valuation of b and g if i>1} yields \[\valustar_\sigma^*(g_i)=W_i^*\oplus\valustar_\sigma^*(b_{i+2})=W_i^*\oplus B_{i+2}^*=W_i^*\oplus L_{i+2}^*=R_i^*.\]
	\item Let $\neg\sigmabar(d_i)$.
		Consider $G_n=S_n$ first.
		By \Pref{BR2}, either $\tau^{\sigma}(F_{i,0})=d_{i,0,k}$ where $\sigma(d_{i,0,k})=e_{i,0,k}$ and $\sigma(e_{i,0,k})=b_2$ for some $k\in\{0,1\}$ or $\tau^{\sigma}(F_{i,0})=s_{i,0}$.
		Since player~1 chooses $\tau^{\sigma}(F_{i,0})$ such that the valuation of $g_i$ is minimized we need to compare $B_2^\P\cup\{g_{i}\}$ (if player 1 chooses $d_{i,0,k}$) and $R_{i}^\P$ (if player 1 chooses $s_{i,0}$).
		Note that $\sigma(b_2)=b_3$ if $\relbit{\sigma}>2$ by \Pref{EB6} and that $\relbit{\sigma}=2$ implies $B_2^\P=L_2^\P$.
		We prove that this implies $R_i^\P\lhd L_2^\P\cup\{g_i\}$ and thus $\valustar_{\sigma}^\P(g_i)=R_i^\P$.
		As mentioned before, $\sigma(g_i)=F_{i,0}$.
		In addition, we have $\sigmabar(eb_{i,0})\wedge\nsigmabar(eg_{i,0})$.
		Thus, by \Pref{EB1}, we have $\sigmabar(b_{i+1})=\sigmabar(b_{\relbit{\sigma}})\neq 0$.
		Hence $\sigma(b_{\relbit{\sigma}})=g_{\relbit{\sigma}}$ and $\sigma(b_{i+1})=g_{i+1}$.
		The statement thus follows from \Cref{lemma: VV Lemma}~(2).
		Consider $G_n=M_n$ next.
		By the choice of $i$ and assumption, it holds that $i=\min\{k\geq i\colon\neg\sigmabar(d_k)\}<\relbit{\sigma}$.
		By \Pref{BR2}, this implies $\nsigmabar(eg_{i,0})$.
		Thus $\valustar_{\sigma}^\M(F_{i,j})=\valustar_{\sigma}^\M(b_2)$.
		Therefore, $\valustar_{\sigma}^\M(g_i)=\rew{g_i}+\valustar_\sigma^\M(b_2)=\rew{g_i}+B_2^\M$ by \Cref{lemma: Valuation of b and g if i>1}.
\end{itemize}

Now let $i<\relbit{\sigma}-1$, implying $\sigma(g_i)=F_{i,1}$ by \Pref{BR1} and $\relbit{\sigma}\geq 3$.
\begin{enumerate}
	\item Let $\sigmabar(d_i)$.
		Consider $G_n=S_n$ first.
		Then, $\tau^{\sigma}(F_{i,1})=s_{i,1}$ and $\sigma(s_{i,1})=h_{i,1}$ by \Pref{S2}.
		Using the induction hypotheses, this yields \[\valustar_\sigma^\P(g_i)=W_i^\P\cup\valustar_\sigma^\P(g_{i+1})=W_i^\P\cup R_{i+1}^\P=R_i^\P.\]
		If $G_n=M_n$, then the same property implies $\valustar_{\sigma}^\M(F_{i,1})=\valustar_{\sigma}^\M(s_{i,1})$ as well as $\sigma(s_{i,1})=h_{i,1}$.
		Thus $\valustar_\sigma^\M(g_i)=W_i^\M+\valustar_{\sigma}^\M(g_{i+1})$.
		Applying the induction hypotheses to $\valustar_{\sigma}^\M(g_{i+1})$ yields the result.

	\item Let $\nsigmabar(d_i)$.
		Consider $G_n=S_n$ first.
		By the same arguments used for $i=\relbit{\sigma}-1$, we need to show $R_i^\P\lhd B_2^\P\cup\{g_i\}$.
		By Properties (\ref{property: BR2}) and (\ref{property: EB2}) we thus have $\sigma(b_2)=b_3$ and $B_2=L_2$ as $\relbit{\sigma}\geq 3$.
		It thus suffices to prove $R_i^\P\lhd L_2^\P\cup\{g_i\}$.
		Let, for the sake of contradiction, $\relbit{\sigma}=\min\{i'\colon\sigma(b_{i'})=b_{i'+1}\}$.
		By \Pref{BR2} and $\neg\sigmabar(d_i)$, we have $\sigmabar(eb_{i,1})\wedge\neg\sigmabar(eg_{i,1})$.
		Thus, by \Pref{EB1}, $\sigma(b_{i+1})=b_{i+2}$.
		But this implies $\relbit{\sigma}\leq i+1$, contradicting $i<\relbit{\sigma}$.
		Hence $\incorrect{\sigma}\neq\emptyset$, implying $\sigma(b_{\relbit{\sigma}})=g_{\relbit{\sigma}}$ by \Cref{lemma: Traits of Relbit}.
		But then, the statement follows from \Cref{lemma: VV Lemma}~(2).

		Consider $G_n=M_n$ next.
		Then $i=\min\{k\geq i\colon\neg\sigmabar(d_k)\}<\relbit{\sigma}$.
		As in the case $i=\relbit{\sigma}-1$, we have $\valustar_{\sigma}^\M(F_{i,j})=\valustar_{\sigma}^\M(b_2)$.
		Therefore, \[\valustar_{\sigma}^\M(g_i)=\rew{g_i}+\valustar_{\sigma}^\M(b_2)=\rew{g_i}+B_2^\M\] by \Cref{lemma: Valuation of b and g if i>1}.
\end{enumerate}
\end{proof}

\ValuationOfGOne*

\begin{proof}
This statement is proven by distinguishing between several cases.
Most of the cases are proven by backwards induction, some are proven directly.
\begin{enumerate}
	\item \boldall{$\minsig{b}\leq\minnegsig{s},\minnegsig{g}, G_n=M_n$ and $\nsigmabar(d_1)$:}
		We prove that this implies that we have $\valustar_{\sigma}^\M(g_1)=\rew{g_1}+\valustar_{\sigma}^\M(b_2)$.
		By \Cref{lemma: Config implied by Aeb}, $\minsig{b}\leq\minnegsig{s},\minnegsig{g}$ implies $\minsig{b}=2$.
		Thus, $\sigma(b_2)=g_2, \sigma(g_1)=F_{1,1}$ and it suffices to prove $\sigmabar(eb_{1,1})\wedge\nsigmabar(eg_{1,1})$.
		By \Pref{EG1}, it cannot hold that $\sigmabar(eg_{1,1})\wedge\nsigmabar(eb_{1,1})$ as this would imply $\sigma(s_{1,1})=b_1$.
		This however contradicts $\sigma(s_{1,1})=h_{1,1}$ which follows from $1<\minnegsig{s}$ and $\sigma(g_1)=F_{1,1}$.
		By \Pref{EBG3}, we cannot have $\sigmabar(eb_{1,1})\wedge\sigmabar(eg_{1,1})$ as this would imply $\sigmabar(d_1)$, contradicting  the current assumptions.
		Thus, $\sigmabar(eb_{1,1})\wedge\nsigmabar(eg_{1,1})$.
		
	\item \boldall{$\minsig{b}\leq\minnegsig{s},\minnegsig{g}, G_n=M_n$  and $\sigmabar(d_1)$:}
		By \Cref{lemma: Config implied by Aeb}, it holds that $\minsig{b}=2$.
		This implies $\sigma(b_2)=g_2, \sigma(g_1)=F_{1,1}$ and $\sigma(s_{1,1})=h_{1,1}$.
		Thus, the chose cycle center of level $1$ is closed, implying $\valustar_{\sigma}^\M(g_1)=W_1^\M+\valustar_{\sigma}^\M(b_2)$.
		
	\item \boldall{$\minsig{b}\leq\minnegsig{s},\minnegsig{g}$ and $G_n=S_n$:}
		By the same argument used in the last case, it suffices to prove $\sigmabar(d_1)$.
		This however follows since $\nsigmabar(eb_1)$ by \Pref{MNS1} and $\nsigmabar(eg_1)$ by \Pref{EG1}. 

		\item \boldall{$\minnegsig{g}<\minnegsig{s},\minsig{b}\wedge\nsigmabar(eb_{\minnegsig{g}})\wedge[G_n=S_n\implies\nsigmabar(b_{\minnegsig{g}+1})]$:}
			We prove that \begin{equation} \label{equation: Case Four Valuation GOne Lemma}
				\valustar_{\sigma}^*(g_i)=\bigoplus_{i'=i}^{\minnegsig{g}}W_{i'}^*\oplus\valustar_\sigma^*(b_{\minnegsig{g}+2})
			\end{equation} for all $i\leq\minnegsig{g}$ by backwards induction.
			Let $i=\minnegsig{g}$.
			Then, by the choice of $i$, $\sigma(g_i)=F_{i,0}$ and $\sigma(s_{i,0})=h_{i,0}$.
			Since we assume $\nsigmabar(eb_{\minnegsig{g}})=\nsigmabar(eb_{i,0})$, $F_{i,0}$ does not escape towards $b_2$.
			In addition, by \Pref{EG1}, it cannot be the case that $\sigmabar(eg_{i,0})\wedge\nsigmabar(eb_{i,0})$ as this would imply $\sigma(s_{i,0})=b_1$.
			Hence $F_{i,0}$ is closed, so $\valustar_{\sigma}^*(F_{i,0})=\valustar_{\sigma}^*(s_{i,0})$.
			Consequently, $\valustar_\sigma^*(g_i)=W_i^*\oplus\valustar_\sigma^*(b_{i+2})$ which is exactly \Cref{equation: Case Four Valuation GOne Lemma} by the choice of $i$.
		
			Let $i<\minnegsig{g}$.
			By $i<\minnegsig{g}$, it holds that $\sigma(g_i)=F_{i,1}$.
			By \Pref{MNS2}, it also holds that $\neg\sigmabar(eb_i)$.
			Using \Pref{EG1} as before, we conclude that $F_{i,1}$ is closed, so $\valustar_\sigma^*(F_{i,1})=\valustar_\sigma^*(s_{i,1})$.
			Also, $\sigma(s_{i,1})=h_{i,1}$ since $i<\minnegsig{s}$.
			This implies $\valustar_\sigma^*(g_i)=W_i^*\oplus\valustar_\sigma^*(g_{i+1})$, so \Cref{equation: Case Four Valuation GOne Lemma} follows from the induction hypothesis.
		
		\item \boldall{$\minnegsig{g}<\minnegsig{s},\minsig{b}\wedge\sigmabar(b_{\minnegsig{g}+1})$ and $G_n=S_n$:}
			We prove that \begin{equation} \label{equation: Case Five Valuation GOne Lemma}
				\valustar_{\sigma}^\P(g_i)=\bigcup_{i'=i}^{\minnegsig{g}}W_{i'}^\P\cup\valustar_\sigma^\P(b_{\minnegsig{g}+2})
			\end{equation} for all $i\leq\minnegsig{g}$ by backwards induction.
			Let $i=\minnegsig{g}$.
			Then, by the choice of $i$, it holds that $\sigma(g_i)=F_{i,0}$ and $\sigma(s_{i,0})=h_{i,0}$.
			By \Pref{EG1}, it cannot be the case that $\sigmabar(eg_{i,0})\wedge\nsigmabar(eb_{i,0})$ as this would imply $\sigma(s_{i,0})=b_1$.
			By \Pref{EBG1}, it cannot be the case that $\sigmabar(eg_{i,0})\wedge\sigmabar(eb_{i,0})$ as this would imply $\sigmabar(b_{i+1})=0$, contradicting the assumption.
			In particular, this implies $\nsigmabar(eg_{i,0})$.
			Thus, depending on the choice of player $1$, either $\valustar_{\sigma}^\P(F_{i,0})=\valustar_{\sigma}^\P(b_2)$ or $\valustar_{\sigma}^\P(F_{i,0})=\valustar_{\sigma}^\P(s_{i,0})$.
			It is now easy to see that $\sigmabar(b_{\minnegsig{g}+1})$ implies $\valustar_{\sigma}^\P(s_{i,0})\lhd\valustar_{\sigma}^\P(b_2)$ and thus $\tau^{\sigma}(F_{i,0})=s_{i,0}$.
			Hence $\valustar_\sigma^\P(g_i)=W_i^\P\cup\valustar_\sigma^\P(b_{i+2})$ as required.
		
			Let $i<\minnegsig{g}$.
			By $i<\minnegsig{g}$, it follows that $\sigma(g_i)=F_{i,1}$ and consequently $\sigma(s_{i,1})=h_{i,1}$.
			Using Properties (\ref{property: EG1}), (\ref{property: EBG1}) and $i<\minsig{b}$, we can again conclude that $\nsigmabar(eg_{i,1})$.
			Consequently, either $\valustar_{\sigma}^\P(F_{i,1})=\valustar_{\sigma}^\P(b_2)$ or $\valustar_{\sigma}^\P(F_{i,1})=\valustar_{\sigma}^\P(s_{i,1})$.
			Using $\sigmabar(b_{\minsig{g}+1})$ and the induction hypotheses, it is again follows that$\valustar_{\sigma}^\P(s_{i,1})\lhd\valustar_{\sigma}^\P(b_2)$.
			Thus, the statement again follows by applying the induction hypotheses.

		\item \boldall{$\minnegsig{g}<\minnegsig{s},\minsig{b}\wedge\nsigmabar(b_{\minnegsig{g}+1})\wedge\sigmabar(eb_{\minnegsig{g}})$ and $G_n=S_n$:}
			Let, for the sake of contradiction, $\minnegsig{g}>1$.
			Then, by \Cref{lemma: Config implied by Aeb}, $\minsig{b}=\minnegsig{g}+1$ and in particular $\sigmabar(b_{\minnegsig{g}+1})$, contradicting the assumption.
			Thus, $\minnegsig{g}=1$.
			This implies $\sigma(g_1)=F_{1,0}$ and $\sigma(s_{1,0})=h_{1,0}$.
			Let, for the sake of contradiction, $\sigmabar(eb_{1,0})\wedge\sigmabar(eg_{1,0})$.
			Then, $\sigmabar(g_1)=\sigmabar(b_2)$ as $\nsigmabar(b_2)$ by assumption and $\minnegsig{g}=1$.
			But then, \Pref{EBG3} implies $\sigmabar(d_1)$ which is a contradiction.
			Consequently, $\sigmabar(eb_{i,0})\wedge\nsigmabar(eg_{i,0})$.
			Since \[\valustar_{\sigma}(b_2)\lhd\{s_{1,0},h_{1,0}\}\cup\valustar_{\sigma}(b_2)=\rew{s_{1,0},h_{1,0}}\cup\valustar_{\sigma}(b_3),\] this yields $\valustar_{\sigma}^\P(g_1)=\{g_1\}\cup\valustar_{\sigma}(b_{3})$.
		
	\item \boldall{$\minnegsig{g}<\minnegsig{s},\minsig{b}\wedge\sigmabar(eb_{\minnegsig{g}})$ and $G_n=M_n$:}
		We prove that \begin{equation} \label{equation: Case Seven Valuation GOne Lemma}
			\valustar_{\sigma}^\M(g_i)=\sum_{i'=i}^{\minnegsig{g}-1}W_{i}^\M+\rew{g_{\minnegsig{g}}}+\valustar_{\sigma}^\M(b_2)
		\end{equation} for all $i\leq\minnegsig{g}$ by backwards induction.
		Let $i=\minnegsig{g}$.
		Then, by construction, $\minnegsig{g}\neq n$ as $\minsig{b}\leq n$ .
		We prove $\nsigmabar(eg_{i,0})$.
		Assume otherwise, implying $\sigmabar(eb_{i,0})\wedge\sigmabar(eg_{i,0})$.
		Assume $\minnegsig{g}>1$.
		Then, by \Cref{lemma: Config implied by Aeb}, $\sigmabar(b_{\minnegsig{g}+1})=1$, contradicting \Pref{EBG1}. 
		Hence assume $\minnegsig{g}=1$.
		If $\sigmabar(b_2)=\sigmabar(g_1)$, then \Pref{EBG3} implies $\sigmabar(d_1)$, contradicting the assumption.
		If $\sigmabar(b_2)\neq\sigmabar(g_1)$, then $\sigmabar(b_2)=1$, again contradicting \Pref{EBG1}.
		Thus, $\nsigmabar(eg_{i,0})$ needs to hold.
		Consequently, as $\sigmabar(eb_{i,0})\wedge\nsigmabar(eg_{i,0})$, this yields $\valustar_{\sigma}^\M(g_{\minnegsig{g}})=\langle g_{\minnegsig{g}}\rangle+\valustar_{\sigmae}^\M(b_2)$.
		
		Let $i<\minnegsig{g}$.
		Then $\sigma(g_i)=F_{i,1}$ and $\sigma(s_{i,1})=h_{i,1}$.
		In addition, by \Pref{MNS2}, $\nsigmabar(eb_{i,1})$.
		Since $\sigmabar(eg_{i,1})$ would imply $\sigma(s_{i,1})=b_1$ by \Pref{EB1}, we have $\nsigmabar(eb_{i})\wedge\nsigmabar(eg_i)$, implying $\sigmabar(d_i)$.
		This implies $\valustar_{\sigma}^\M(g_i)=W_i^\M+\valustar_{\sigma}^\M(g_{i+1})$ and the statement then follows by using the induction hypotheses.
		
	\item \boldall{$\minnegsig{s}\leq\minnegsig{g}<\minsig{b}$ or $\minnegsig{s}<\minsig{b}\leq\minnegsig{g}$:}
		Let  $m\coloneqq\min\{\minnegsig{g},\minnegsig{s}\}$.
		We prove that \begin{equation} \label{equation: Case Eight Valuation GOne Lemma}
			\valustar_\sigma^*(g_i)=\bigoplus_{i'=i}^{m-1}W_{i'}^*\oplus\ubracket{g_m}\oplus\valustar_\sigma^*(b_2)
		\end{equation}	for all $i\leq m$. 
		Let $i=m$ and $j\coloneqq\sigmabar(g_i)$.
		We can assume $i=\minnegsig{s}$ in both cases, implying $i=\minnegsig{s}\leq n$ in both cases.
		By either \Pref{MNS4} or \Pref{MNS6}, we have $\sigmabar(eb_i)\wedge\nsigmabar(eg_i)$.
		This implies $\valustar_{\sigma}^\M(F_{i,j})=\valustar_{\sigma}^\M(b_2)$, hence the statement follows for $G_n=M_n$.
		Consider the case $G_n=S_n$.
		Since $i=\minnegsig{s}$, we have $\sigma(s_{i,j})=b_1$.
		Therefore, using $\sigma(b_1)=b_2$, we obtain $\valustar_\sigma^\P(s_{i,j})=\{s_{i,j}\}\cup\valustar_\sigma^\P(b_2)\rhd\valustar_\sigma^\P(b_2).$		
		Thus $\tau^{\sigma}(F_{i,j})=d_{i,j,k}$ and therefore $\valustar_\sigma^\P(g_m)=\{g_m\}\cup\valustar_\sigma^\P(b_2)$.
		
		Let $i<m$.
		Since $i<\minnegsig{s}\leq\minnegsig{g}$ in both cases, $\sigma(g_i)=F_{i,1}$ and $\sigma(s_{i,1})=h_{i,1}$.
		Let $G_n=S_n$.		
		By \Pref{EG1}, it cannot be the case that $\sigmabar(eg_{i,1})\wedge\nsigmabar(eb_{i,1})$ as this would imply $\sigma(s_{i,1})=b_1$.
		By \Pref{EBG1}, it cannot be the case that $\sigmabar(eb_{i,1})\wedge\sigmabar(eg_{i,1})$ as this would imply $\sigma(b_{i+1})=g_{i+1}$, contradicting the choice of~$i$.
		Hence, either $\sigmabar(d_{i,j})$ or $\sigmabar(eb_{i,1})\wedge\nsigmabar(eg_{i,1})$.
		We prove that $\tau^{\sigma}(F_{i,j})=s_{i,j}$ holds in any case.	
		It suffices to consider the second case as this follows directly in the first case.
		Since $\sigmabar(eb_{i,1})\wedge\nsigmabar(eg_{i,1})$, it suffices to prove $\valustar_{\sigma}^\P(s_{i,1})\lhd\valustar_{\sigma}^\P(b_2)$.
		This however follows by the induction hypotheses and $\valustar_{\sigma}^\P(s_{i,1})=\{s_{i,1},h_{i,1}\}\cup\valustar_{\sigma}^\P(g_{i+1})$.
		Thus, $\tau^{\sigma}(F_{i,j})=s_{i,j}$ for $G_n=S_n$ in any case.

		If $G_n=M_n$, then either \Pref{MNS3} or \Pref{MNS5} implies the cycle center $F_{i,1}$ is closed.
		Hence, using the induction hypotheses, \Cref{equation: Case Eight Valuation GOne Lemma} follows from $\valustar_{\sigma}^*(g_i)=W_i^*\oplus\valustar_{\sigma}^*(g_{i+1})$.
	\end{enumerate}
	Note that the cases listed here suffice, i.e., every possible relation between the three parameters $\minnegsig{s},\minnegsig{g}$ and $\minsig{b}$ is covered by exactly one of the cases.
\end{proof}

\ExactBehaviorOfCounterstrategy*

Proving this statement requires the following additional lemma.

\begin{lemma} \label{lemma: Can escape right lane}
Let $G_n=S_n$ and let $\sigma\in\reach{\sigma_0}$.
Let $i\in[n]$ such that $\sigmabar(eg_{i})$.
Then there is some $i'<i$ such that either $\sigma(g_{i'})=F_{i',0}$ or $\sigma(s_{i',\sigmabar(g_{i'})})=b_1$.
\end{lemma}

\begin{proof}
Let, for the sake of contradiction, $\sigma(g_{i'})=F_{i',1}$ and $\sigma(s_{i',\sigmabar(g_{i'})})=h_{i',\sigmabar(g_{i'})}$ for all $i'<i$.
Then, player 1 can create a cycle by setting $\tau^{\sigma}(F_{k,1})=s_{k,1}$ for all $k<i$ and $\tau(F_{i,\sigmabar(g_i)})=d_{i,\sigmabar(g_i),k}$ where $k$ is chosen such that the cycle center escapes towards $g_1$.
But this contradicts the fact that $S_n$ is a sink game.
\end{proof}

\begin{proof}[Proof of \Cref{lemma: Exact Behavior Of Counterstrategy}]
Consider a cycle center $F_{i,j}$.
We distinguish the following cases:
\begin{enumerate}
	\item \boldall{$\sigmabar(d_{i,j})$:}
		Then $F_{i,j}$ is closed and $\valustar_\sigma^\P(F_{i,j})=\valustar_\sigma^\P(s_{i,j})$ since $S_n$ is a sink game.
	\item \boldall{$\sigmabar(eg_{i,j}),\neg\sigmabar(eb_{i,j})$ and $\relbit{\sigma}=1$:}
		We prove $\valustar_\sigma^\P(F_{i,j})=\{s_{i,j}\}\cup\valustar_{\sigma}^\P(b_2).$
		Since $\nsigmabar(eb_{i,j})$, player 1 choose a cycle vertex escaping towards $g_1$ or $s_{i,j}$.
		As player~$1$ minimizes the vertex valuations, it suffices to prove $\valustar_\sigma^\P(s_{i,j})\lhd\valustar_\sigma^\P(g_1)$. 				
		\Pref{EG2} implies $\sigmabar(d_1)$ and hence $\valustar_\sigma^\P(g_1)=\{g_1\}\cup\valustar_\sigma^\P(s_{1,\sigmabar(g_1)}).$
		By \Pref{EG3}, also $\sigma(s_{1,\sigmabar(g_1)})=h_{1,\sigmabar(g_1)}$.
		Since $\sigmabar(g_1)=\sigmabar(b_2)$ by \Pref{EG4}, it then follows that $\valustar_{\sigma}^\P(h_{1,\sigmabar(g_1)})=\valustar_{\sigma}^\P(b_2)$.
		Hence $\valustar_\sigma^\P(g_1)=W_1^\P\cup\valustar_\sigma^\P(b_2)$.
		Furthermore,  $\valustar_\sigma^\P(s_{i,j})=\{s_{i,j}\}\cup\valustar_\sigma^\P(b_1)=\{s_{i,j}\}\cup\valustar_\sigma^\P(b_2)$ by \Pref{EG1} and $\relbit{\sigma}=1$, and consequently $\valustar_\sigma^\P(s_{i,j})\lhd\valustar_\sigma^\P(g_1).$
	\item \boldall{$\sigmabar(eg_{i,j}),\neg\sigmabar(eb_{i,j})$ and $\relbit{\sigma}\neq 1$:}
		We prove $\valustar_\sigma^\P(F_{i,j})=\valustar_\sigma^\P(g_1)$.
		Assume $\sigma(s_{i,j})=b_1$.
		Since $\sigma(b_1)=g_1$ by \Cref{lemma: b1 iff relbit}, \[\valustar_\sigma^\P(g_1)=\valustar_\sigma^\P(b_1)\lhd\{s_{i,j}\}\cup\valustar_\sigma^\P(b_1)=\{s_{i,j}\}\cup\valustar_\sigma^\P(g_1)=\valustar_\sigma^\P(s_{i,j}),\] implying $\tau^{\sigma}(F_{i,j})=s_{i,j}$ since $\neg\sigmabar(eb_{i,j})$.
		Assume $\sigma(s_{i,j})=h_{i,j}$.
		By \Pref{EG5}, it then holds that $j=\sigmabar(b_{i+1})$. 
		Since also $\sigma(b_1)=g_1$ and $\relbit{\sigma}\neq 1$, it suffices to prove $R_1^\P\lhd\{s_{i,j},h_{i,j}\}\cup\valustar_\sigma(b_{i+2-j})$.
		This can be shown by the following case distinction based on $j$ and the relation between $i+1$ and $\relbit{\sigma}$.
		\begin{enumerate}
			\item Let $j=1$ and $i+1<\relbit{\sigma}$.
				Then $\valustar_\sigma^\P(b_{i+2-j})=\valustar_\sigma^\P(b_{i+1})=R_{i+1}^\P$ since $\sigma(b_{i+1})=g_{i+1}$ by \Pref{EG5}.
				It thus suffices to prove $R_1^\P\lhd\{s_{i,j},h_{i,j}\}\cup R_{i+1}^\P$ which follows from $\sigdiff{R_{i+1}^\P\cup\{s_{i,j},h_{i,j}\}}{R_1^\P}=g_{i}$.
			\item Let$j=1$ and $i+1\geq\relbit{\sigma}$.
				Then $\valustar_\sigma^\P(b_{i+1})=L_{i+1}^\P$ and it suffices to prove $R_1^\P\lhd\{s_{i,j},h_{i,j}\}\cup L_{i+1}^\P$.
				This follows from \Cref{lemma: VV Lemma} if $i\geq\relbit{\sigma}$ and is easy to verify for $i+1=\relbit{\sigma}$.
			\item Let $j=0$ and $i+2\leq\relbit{\sigma}$.
				We first show $\valustar_\sigma^\P(b_{i+2})=L_{i+2}^\P$.
				If $i+2=\relbit{\sigma}$, then this follows by definition.
				Thus let $i+1<\relbit{\sigma}-1$.
				Since  $\sigma(s_{i,j})=h_{i,j}$, \Pref{EG5} implies $\sigma(b_{i+1})=b_{i+2}$.
				This implies $\sigma(b_{i+2})=b_{i+3}$ by \Pref{B1}.
				Consequently, also $\valustar_\sigma^\P(b_{i+2})=L_{i+2}^\P$ in this case.
				As usual, we have $\valustar_\sigma(b_1)=R_1^\P$ and thus prove $R_1^\P\lhd\{s_{i,j},h_{i,j}\}\cup L_{i+2}^\P$.
				But this follows from \Cref{lemma: VV Lemma} since $\sigma(b_{i+1})=b_{i+2}$ implies $L_{i+1}^\P=L_{i+2}^\P$ as follows:
				By $\relbit{\sigma}\geq i+2$ and $\sigma(b_{i+1})=b_{i+2}$, we have $\relbit{\sigma}\neq\min\{i'\in[n]\colon\sigma(b_{i'})=b_{i'+1}\}$.
				Hence $\incorrect{\sigma}\neq\emptyset$ and $\sigma(b_{\relbit{\sigma}})=g_{\relbit{\sigma}}$ by \Cref{lemma: Traits of Relbit}.				
				Thus $i+2\leq\relbit{\sigma}$ implies $R_1^\P\lhd L_{i+2}^\P\lhd L_{i+2}^\P\cup\{s_{i,j},h_{i,j}\}$ by \Cref{lemma: VV Lemma}.

			\item Let $j=0$ and $i+2>\relbit{\sigma}$.
				Then $\valustar_\sigma^\P(b_{i+2})=L_{i+2}^\P=L_{i+1}^\P$ since $\sigma(b_{i+1})=b_{i+2}$ by \Pref{EG5}.
				Thus, $R_1^\P\lhd\{s_{i,j},h_{i,j}\}\cup L_{i+1}^\P$ by \Cref{lemma: VV Lemma}.
		\end{enumerate}

		\item \boldall{$\sigmabar(eb_{i,j}),\neg\sigmabar(eg_{i,j}),\relbit{\sigma}=1$ and $(\sigma(s_{i,j})=b_1\vee\sigmabar(b_{i+1})=j)$:}
			We prove $\valustar_\sigma^\P(F_{i,j})=\valustar_\sigma^\P(b_2)$.
			Due to $\neg\sigmabar(eg_{i,j})$, it suffices to show $\valustar_\sigma^\P(b_2)\lhd\valustar_\sigma^\P(s_{i,j})$.
			If $\sigma(s_{i,j})=b_1$, then the statement follows directly since $\relbit{\sigma}=1$ implies $\sigma(b_1)=b_2)$.			
			Hence let $\sigma(s_{i,j})=h_{i,j}\wedge\sigmabar(b_{i+1})=j$.
			As $\relbit{\sigma}=1$, it holds that $\valustar_\sigma^\P(b_i)=L_i^\P$ for all $i\in[n]$.
			Hence, $\valustar_\sigma^\P(b_2)=L_2^\P$ and $\valustar_\sigma^\P(s_{i,j})=L_{i+2-j}^\P\cup\{s_{i,j},h_{i,j}\}$ since $\sigmabar(b_{i+1})=j$. 
			We thus prove $L_2^\P\lhd L_{i+2-j}^\P\cup\{s_{i,j},h_{i,j}\}$.
			It is sufficient to show $\sigma(b_i)=b_{i+1}$ since this implies $W_i\nsubseteq L_2$.
			For the sake of contradiction let $\sigma(b_i)=g_i$.
			Since $\relbit{\sigma}=1$, \Pref{D1} implies that the chosen cycle center of level $i$ is closed.
			However, $\sigmabar(b_{i+1})=j$ and $\incorrect{\sigma}=\emptyset$ then imply $\sigmabar(g_i)=j$.
			Hence this cycle center is $F_{i,j}$, contradicting $\sigmabar(eb_{i,j})$.
		\item \boldall{$\sigmabar(eb_{i,j}),\neg\sigmabar(eg_{i,j})$ and $\relbit{\sigma}\neq 1\vee(\sigma(s_{i,j})=h_{i,j}\wedge\sigmabar(b_{i+1})\neq j)$:}
			By \Pref{EB1}, we can assume $\sigmabar(b_{i+1})\neq j$ in either case.
			We prove that it holds that	$\valustar_\sigma^\P(F_{i,j})=\valustar_\sigma^\P(s_{i,j})$ by proving $\valustar_{\sigma}^{\P}(s_{i,j})\lhd\valustar_{\sigma}^{\P}(b_2)$.
			\begin{enumerate}
				\item Let $\sigma(s_{i,j})=h_{i,j}$ and $\sigma(b_1)=b_2$, implying $\relbit{\sigma}=1$.
				Let $j=0$.
				Since none of the cycle vertices of $F_{i,j}$ escapes to $g_1$ and $\sigma(s_{i,j})=\sigma(s_{i,0})=h_{i,0}$, we have to prove $\{s_{i,0},h_{i,0}\}\cup\valustar_\sigma^\P(L_{i+2})\lhd\valustar_\sigma^\P(L_2)$ since $\relbit{\sigma}=1$.
				This however follows as $\sigmabar(b_{i+1})\neq j=0$ implies $\sigma(b_{i+1})=g_{i+1}$. 
				Now let $j=1$.
				We then need to prove $\valustar_\sigma^\P(s_{i,1})\lhd\valustar_\sigma^\P(b_2)$.
				However, the exact valuation of~$s_{i,1}$ is not clear in this case and depends on several vertices of level 1 and $i+1$.
				To be precise we can have the following paths:
				
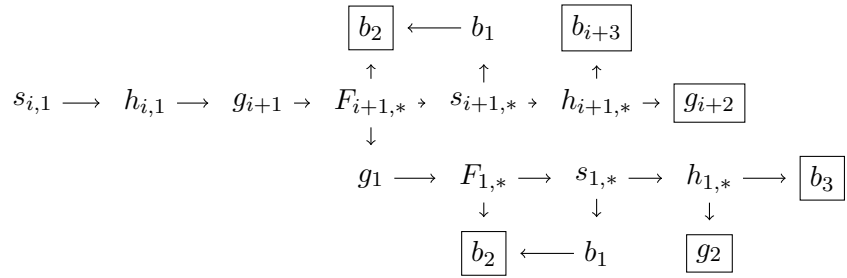
\begin{center}
	\begin{tikzpicture}
		\def\x{1.5}
		\def\y{1}
		\node (si1) at (0,0) {$s_{i,1}$};
		\node (hi1) at (\x,0) {$h_{i,1}$};
		\node (gi+1) at (2*\x,0) {$g_{i+1}$};
		\node (Fi+1*) at (3*\x,0) {$F_{i+1,\ast}$};
		\node[draw] (b22) at (3*\x,\y) {$b_2$};
		\node (g1) at (3*\x,-\y) {$g_1$};
		\node (F1*) at (4*\x,-\y) {$F_{1,\ast}$};
		\node[draw] (b21) at (4*\x,-2*\y) {$b_2$};
		\node (s1*) at (5*\x,-\y) {$s_{1,\ast}$};
		\node (b11) at (5*\x,-2*\y) {$b_1$};
		\node (h1*) at (6*\x,-\y) {$h_{1,\ast}$};
		\node[draw] (b3) at (7*\x,-\y) {$b_3$};
		\node[draw] (g2) at (6*\x,-2*\y) {$g_2$};
		\node (si+1*) at (4*\x,0) {$s_{i+1,\ast}$};
		\node (b12) at (4*\x,\y) {$b_1$};
		\node (hi+1*) at (5*\x,0) {$h_{i+1,\ast}$};
		\node[draw] (bi+3) at (5*\x,\y) {$b_{i+3}$};
		\node[draw] (gi+2) at (6*\x,0) {$g_{i+2}$};
		\draw[->, shorten >=1ex] (si1) edge (hi1) (hi1) edge (gi+1) (gi+1) edge (Fi+1*) (Fi+1*) edge (si+1*) (si+1*) edge (hi+1*) (hi+1*) edge (gi+2);
		\draw[->, shorten >=1ex] (si+1*) edge (b12) (b12) edge (b22);
		\draw[->, shorten >=1ex] (Fi+1*)--(b22);
		\draw[->, shorten >=1ex] (Fi+1*)--(g1) (g1) edge (F1*) (F1*) edge (s1*) (s1*) edge (h1*) (h1*) edge (b3);
		\draw[->, shorten >=1ex] (h1*)--(g2);
		\draw[->, shorten >=1ex] (s1*)--(b11) (b11) edge (b21);
		\draw[->, shorten >=1ex] (F1*)--(b21);		
		\draw[->, shorten >=1ex] (hi+1*)--(bi+3);			
	\end{tikzpicture}
\end{center}

				We show that $\valustar_\sigma^\P(s_{i,1})\lhd\valustar_\sigma^\P(b_2)$ holds for all marked ``endpoints'' that could be reached by $s_{i,1}$.
				In all cases, $j=1$ and $\sigmabar(b_{i+1})\neq j$ imply $\sigma(b_{i+1})=b_{i+2}$.
				\begin{itemize}
					\item \boldall{$b_2$:}
						Then, the statement follows directly as it is easy to verify that we then have $\sigdiff{\valustar_\sigma^\P(b_2)}{\valustar_\sigma^\P(s_{i,1})}=g_{i+1}$ in each possible case.
					\item \boldall{$b_{i+3}$:}
						Then $\valustar_\sigma^\P(s_{i,1})=\{s_{i,1},h_{i,1},g_{i+1},s_{i+1,0},h_{i+1,0}\}\cup\valustar_\sigma^\P(b_{i+3})$ and $\sigma(g_{i+1})=F_{i+1,0}$.
						Since $\sigma(s_{i,1})=h_{i,1}$, \Pref{B3} implies $\sigmabar(g_{i+1})\neq\sigmabar(b_{i+2})$, hence $\sigma(b_{i+2})=g_{i+2}$.	
						Since $\relbit{\sigma}=1$ implies $\valustar_\sigma^\P(b_2)=L_2^\P$ we therefore have $W_{i+2}^\P\subseteq \valustar_\sigma^\P(b_2)$. 
						This yields the statement.	
					\item \boldall{$b_3$:}
						Then $\valustar_\sigma^\P(s_{i,1})=\{s_{i,1},h_{i,1},g_{i+1}\}\cup W_1^\P\cup\valustar_\sigma^\P(b_3)$ and thus, \begin{align*}
							\valustar_\sigma^\P(s_{i,1})&=\{s_{i,1},h_{i,1},g_{i+1}\}\cup W_1^\P\cup L_3^\P\\
								&\unlhd\{s_{i,1},h_{i,1},g_{i+1}\}\cup W_1^\P\cup L_2^\P\lhd\valustar_\sigma^\P(b_2).								
						\end{align*}
					\item \boldall{$g_2$:} 
						Then, it holds that $\valustar_\sigma^\P(s_{i,1})=\{s_{i,1},h_{i,1},g_{i+1}\}\cup W_1^\P\cup\valustar_\sigma^\P(g_2)$ and $\valustar_\sigma^\P(b_2)=L_2^\P$. 
						As before we need to show $\valustar_\sigma^\P(s_{i,1})\lhd\valustar_\sigma^\P(b_2)$.
						Note that we can assume $i\geq 2$ since $S_n$ is a sink game and the valuation of $s_{i,1}$ contains a cycle for $i=1$.
						
						Let $\sigma(b_2)=g_2$.
						Then $\valustar_\sigma^\P(g_2)=\valustar_\sigma^\P(b_2)$.
						This implies $\sigdiff{\valustar_\sigma^\P(s_{i,1})}{\valustar_\sigma^\P(b_2)}=g_{i+1}$ since $i\geq 2$ and $W_{i+1}^\P\nsubseteq \valustar_{\sigma}^\P(b_2)$ due to $\sigma(b_{i+1})=b_{i+2}$ and $\relbit{\sigma}=1$.
						Since $g_{i+1}\subseteq\valustar_\sigma^\P(s_{i,1})$, this implies $\valustar_\sigma^\P(s_{i,1})\lhd\valustar_\sigma^\P(b_2)$.
						
						Thus let $\sigma(b_2)=b_3$, implying $\valustar_{\sigma}^{\P}(b_2)=\valustar_{\sigma}^{\P}(b_3)$ and let $k\coloneqq\sigmabar(g_2)$.
						Similar to the picture showing the ``directions'' to which the vertex $s_{i,1}$ can lead, there are several possibilities towards which vertex the path starting in $F_{2,k}$ leads.
						For all of the following cases, the main argument will be the following.
						No matter what choices are made in the lower levels and no matter how many levels the path starting in $g_2$ might traverse, the vertex $g_{i+1}$ contained in $\valustar_\sigma^\P(s_{i,1})$ will always ensure $\valustar_\sigma^\P(s_{i,1})\lhd\valustar_\sigma^\P(b_2)$.
						We distinguish the following cases.						
						\begin{enumerate}
							\item $F_{2,k}$ escapes towards $b_2$.
								Then $\valustar_\sigma^\P(g_2)=\{g_2\}\cup\valustar_\sigma^\P(b_2)$ and thus \[\valustar_\sigma^\P(s_{i,1})=\{s_{i,1},h_{i,1},g_{i+1},g_2\}\cup W_1^\P\cup\valustar_\sigma^\P(b_2)\lhd\valustar_\sigma^\P(b_2).\]
							\item $F_{2,k}$ escapes towards $g_1$.
								Depending on the configuration of level 1, the path can end in different vertices.
								As $S_n$ is a sink game, it cannot end in $g_1$ or $g_2$ since this would constitute a cycle.
								It can thus either end in $b_1,b_2$ or $b_3$.
								Since $\sigma(b_1)=b_2$ and $\sigma(b_2)=b_3$, we then have $\valustar_\sigma^\P(g_2)=\{g_2\}\cup W_1^\P\cup\valustar_\sigma^\P(b_3)$ in either case and thus \begin{align*}
									\valustar_\sigma^\P(s_{i,1})&=\{s_{i,1},h_{i,1},g_{i+1},g_2\}\cup W_1^\P\cup\valustar_\sigma^\P(b_2)\lhd\valustar_\sigma^\P(b_2).
								\end{align*}
							\item $F_{2,k}$ does not escape level 2 and $k=0$.
								In this case, $\tau(F_{2,k})=s_{2,0}$.
								If $\sigma(s_{2,0})=b_1$, then the statement follows by the same arguments used in the last case.
								Thus consider the case $\sigma(s_{2,0})=h_{2,0}$, implying  $\valustar_\sigma^\P(g_2)=W_2^\P\cup\valustar_\sigma^\P(b_4)$.
								Then, since $i\geq2$, \begin{align*}
									\valustar_\sigma^\P(s_{i,1})&=\{s_{i,1},h_{i,1},g_{i+1}\}\cup W_1^\P\cup W_2^\P\cup\valustar_\sigma^\P(b_4)\\
										&\lhd \valustar_\sigma^\P(b_4)\unlhd\valustar_\sigma^\P(b_2).
								\end{align*}
							\item $F_{2,k}$ does not escape level 2 and $k=1$.
								In this case we can use the exact same arguments to show that either $\valustar_\sigma(s_{i,1})\lhd\valustar_\sigma(b_2)$ or that the path reaches $g_4$.
								In fact, the same arguments can be used until vertex $g_{i-1}$ is reached.
								We now show that once this vertex is reached the inequality $\valustar_\sigma^\P(s_{i,1})\lhd\valustar_\sigma^\P(b_2)$ is fulfilled.
								Let $k'\coloneqq\sigmabar(g_{i-1})$.
								If $F_{i-1,k'}$ escapes towards $b_2$,  then $\valustar_\sigma^\P(s_{i,1})\lhd\valustar_\sigma^\P(b_2)$ follows from \[\valustar_\sigma^\P(g_2)=\bigcup_{i'=2}^{i-2}W_{i'}^\P\cup\{g_{i-1}\}\cup\valustar_\sigma^\P(b_2).\]
								If $F_{i-1,k'}$ escapes towards $b_1$ via $s_{i-1,k'}$, then the statement follows analogously since $\sigma(b_1)=b_2$.
								Thus assume that the cycle center escapes towards $g_1$.
								By the same arguments used before, it can be shown that level~1 needs to escape towards either $b_1,b_2$ or $b_3$.
								However, the same calculation used before can be applied in each of these cases.
								
								Next assume that the cycle center $F_{i-1,k'}$ does not escape level $i-1$ but traverses the level and reaches vertex $b_{i+1}$.
								Then, \[\valustar_\sigma^\P(g_2)=\bigcup_{i'=2}^{i-1}W_{i'}^\P\cup\valustar_\sigma^\P(b_{i+1})\]and $\valustar_\sigma^\P(s_{i,1})\lhd\valustar_\sigma^\P(b_{i+1})\unlhd\valustar_\sigma^\P(b_2).$
								The last case we need to consider is if level $i-1$ is traversed and $g_i$ is reached.
								In this case we need to have $\sigma(g_i)=F_{i,0}$ since player 1 could create a cycle otherwise, contradicting that $S_n$ is a sink game.
								If the cycle center $F_{i,0}$ escapes towards $g_1$ or $b_2$ the statement follows by the same arguments used before.
								If it reaches~$b_{i+2}$, then the statement follows from\begin{align*}
									\valustar_\sigma^\P(s_{i,1})&=\{s_{i,1},h_{i,1},g_{i+1}\}\cup\bigcup_{i'=1}^{i}W_{i'}^\P\cup\valustar_\sigma^\P(b_{i+2})\unlhd\valustar_\sigma^\P(b_2).
								\end{align*}
						\end{enumerate}
					\item \boldall{$g_{i+2}$:}
						This implies $\valustar_\sigma^\P(s_{i,1})=\{s_{i,1},h_{i,1}\}\cup W_{i+1}^\P\cup\valustar_\sigma^\P(g_{i+2})$ and we prove $\valustar_\sigma^\P(s_{i,1})\lhd\valustar_\sigma^\P(b_2)$.
						This case is organized similarly to the last case.
						We prove that the statement follows for all but one possible configurations of the levels $i+2$ to $n-1$.
						It then turns out that this missing configuration contradicts \Pref{DN1}.
						We distinguish the following cases.
						\begin{enumerate}
							\item Level $i+2$ escapes towards  $b_2$ via some cycle center $F_{i+2,*}$.
								Then $\valustar_\sigma^\P(g_{i+2})=\{g_{i+2}\}\cup\valustar_\sigma^\P(b_2)$ implies the statement.
							\item Level $i+2$ escapes towards $b_1$ via some upper selection vertex $s_{i+2,*}$.
								Then,  $\valustar_\sigma^\P(g_{i+2})=\{g_{i+2},s_{i+2,*}\}\cup\valustar_\sigma^\P(b_2)$ as $\sigma(b_1)=b_2$, implying the statement.
							\item Level $i+2$ is traversed completely and reaches $b_{i+4}$ directly.
								In this case, $\sigma(g_{i+2})=F_{i+2,0}$.
								Since $\sigma(s_{i+1,1})=h_{i+1,1}$ and $\sigma(b_{i+2})=b_{i+3}$ by \Pref{B3}, this implies $\sigma(b_{i+3})=g_{i+3}$.
								Thus, it holds that $\valustar_\sigma^\P(g_{i+2})=W_{i+2}^\P\cup\valustar_\sigma^\P(b_{i+4})$ which implies the statement due to $\sigma(b_{i+3})=g_{i+3}$.
							\item Level $i+2$ escapes towards $g_1$ via some cycle center $F_{i+2,*}$.
								Then, $\valustar_\sigma^\P(g_{i+2})=\{g_{i+2}\}\cup\valustar_\sigma^\P(g_1)$.
								Consider level 1.
								If  $F_{1,\sigmabar(g_1)}$ escapes towards $b_1,b_2$ or $b_3$, the statement follows by the same arguments used in the last two cases since $\valustar_\sigma^\P(b_3)\unlhd\valustar_\sigma^\P(b_2)$.
								Since~$S_n$ is a sink game, the cycle center cannot escape towards $g_1$.
								Thus assume that it escapes towards $g_2$.
								By the same arguments used previously, the statement either holds or level 3 is traversed and the path reaches $g_4$.
								We now iterate this argument until we reach a level $k<i+2$ such that either $\sigma(g_{k})=F_{k,0}$ or $\sigma(s_{k,\sigmabar(g_k)})=b_1$.
								Such a level exists by \Cref{lemma: Can escape right lane}.
								We only consider the second case here since the statement follows by calculations similar to the previous ones if $\sigma(g_k)=F_{k,0}$.
								Then \[\valustar_\sigma^\P(g_{i+2})=\{g_{i+2}\}\cup\bigcup_{i'=1}^{k-1}W_{i'}^\P\cup\{g_{k},s_{k,*}\}\cup\valustar_\sigma^\P(b_1),\] implying the statement since $k\leq i+1$.
							\item Level $i+2$ is traversed completely and reaches $g_{i+3}$.
								This implies that  $\sigma(s_{i+1,1})=h_{i+1,1}, \sigma(b_{i+2})=b_{i+3}$ and $\sigma(g_{i+2})=F_{i+2,1}$.
								Thus, by \Pref{B3}, also $\sigma(b_{i+3})=b_{i+4}$.
								We can therefore use the same arguments used before since $\valustar_\sigma^\P(g_{i+2})=W_{i+2}^\P\cup\valustar_\sigma^\P(g_{i+3})$.
								That is, the statement either holds or we reach the vertex $g_{n-1}$.
								If level $n-1$ escapes towards $b_1,b_2$ or $g_1$, the statement follows by the same arguments used for level $i+2$.
								We thus assume that level $n-1$ is traversed completely.
								Note that $\sigma(s_{n-2,1})=h_{n-2,1}$ and $\sigma(b_{n-1})=b_n$ (apply \Pref{B3} iteratively).
								Consider the case $\sigma(g_{n-1})=F_{n-1,0}$ first.								
								Then, $\sigma(b_n)=g_n$ by \Pref{B3} and thus \[\valustar_\sigma^\P(s_{i,1})=\{s_{i,1},h_{i,1}\}\cup\bigcup_{i'=i+1}^{n-1}W_{i'}^\P\lhd W_n^\P =\valustar_\sigma^\P(b_n)\unlhd\valustar_\sigma^\P(b_2).\]
								Consider the case $\sigma(g_{n-1})=F_{n-1,1}$ next.
								Since we assume that level $n-1$ does not escape towards one of the vertices $b_1,b_2$ or $g_1$, we traverse this level and reach~$g_n$.
								If level $n$ escapes towards $b_1,b_2$ or $g_1$ the statement follows as usual. 
								We thus assume that the level $n$ is traversed completely.
								We observe that the vertex $h_{n,1}$ has the highest even priority among all vertices in the parity game.
								Thus, player 1 would avoid this vertex if this was possible.
								We thus need to have $\sigmabar(d_n)$.
								But this is a contradiction to \Pref{DN1} since $\sigma(b_{n})=t$ by \Pref{B3} and $\sigma(b_1)=b_2$ by assumption.
						\end{enumerate}
				\end{itemize}
				\item Let $\sigma(s_{i,j})=h_{i,j}$ and $\sigma(b_1)=g_1$, implying $\relbit{\sigma}\neq 1$ by \Cref{lemma: b1 iff relbit}.
					\begin{enumerate}
						\item Let $j=0$ and $i=1$.
							Then, $\valustar_\sigma^\P(s_{i,j})=\{s_{1,0},h_{1,0}\}\cup\valustar_\sigma^\P(b_{3}).$
							By \Pref{EB2}, it holds that $\relbit{\sigma}=2$, implying $\valustar_\sigma^\P(b_2)=L_2^\P$ and $\valustar_\sigma^\P(b_3)=L_3^\P$.
							It thus suffices to prove $L_2^\P\rhd \{s_{1,0},h_{1,0}\}\cup L_3^\P$ which follows from \Pref{EB1} as this implies $\sigma(b_{i+1})\neq j$, so $\sigma(b_2)=g_2$.							
						\item Let $j=0$ and $i>1$.
							Then $\valustar_\sigma^\P(s_{i,j})=\{s_{i,0},h_{i,0}\}\cup\valustar_\sigma^\P(b_{i+2})$.
							By \Pref{EB2}, it follows that $\relbit{\sigma}=i+1$, implying $\valustar_\sigma^\P(b_{i+2})=L_{i+2}^\P$.
							In addition, $\sigma(b_{2})=b_3$ by \Pref{EB3}, implying $\valustar_\sigma^\P(b_2)=L_2^\P$.
							We thus prove $L_2^\P\rhd \{s_{i,0},h_{i,0}\}\cup L_{i+2}^\P$.
							This follows since \Pref{EB1} implies $\sigma(b_{i+1})=\sigma(b_{\relbit{\sigma}})\neq j=0$ and thus $\sigma(b_{\relbit{\sigma}})=g_{\relbit{\sigma}}$
						\item Let $j=1$.
							Then $\valustar_\sigma^\P(s_{i,j})=\{s_{i,1},h_{i,1}\}\cup\valustar_\sigma^\P(g_{i+1}).$
							By \Pref{EB4}, we have $i+1<\relbit{\sigma}$ and, by \Pref{EB1}, also $\sigma(b_{i+1})=b_{i+2}$.
							Therefore, $\valustar_\sigma^\P(g_{i+1})=R_{i+1}^\P$.
							This also implies $\relbit{\sigma}\neq\min\{i'\in[n]\colon\sigma(b_{i'})=b_{i'+1}\}$ since this would contradict $\relbit{\sigma}>i+1$.
							In particular it holds that $\incorrect{\sigma}\neq\emptyset$, implying $\sigma(b_{\relbit{\sigma}})=g_{\relbit{\sigma}}$ by \Cref{lemma: Traits of Relbit}.
							Since $\relbit{\sigma}>i+1\geq 2$, we can apply \Pref{EB6}, implying $\sigma(b_2)=b_3$.							
							Hence $\valustar_\sigma^\P(b_2)=L_2^\P$.
							It thus suffices to prove $L_2^{\P}\rhd \{s_{i,1},h_{i,1}\}\cup R_{i+1}^{\P}$ which follows from $\sigma(b_{\relbit{\sigma}})=g_{\relbit{\sigma}}$ and $i+1<\relbit{\sigma}$.
					\end{enumerate}
				\item Next let $\sigma(s_{i,j})=b_1$ and $\sigma(b_1)=g_1$, implying $\relbit{\sigma}\neq 1$.
					Consider the case $\relbit{\sigma}>2$ first.
					Then, by \Pref{EB6}, $\sigma(b_2)=b_3$, so $\valustar_\sigma^\P(b_2)=L_2^\P$ and $\relbit{\sigma}\neq \min\{i'[n]\colon\sigma(b_{i'})=b_{i'+1}\}$.
					Hence $\incorrect{\sigma}\neq\emptyset$ and thus $\sigma(b_{\relbit{\sigma}})=g_{\relbit{\sigma}}$.
					Since $\valustar_\sigma^\P(b_1)=R_1^\P$,  we prove $\valustar_\sigma^\P(s_{i,j})=\{s_{i,j}\}\cup R_1^\P \lhd L_2^\P$ which follows from $\sigma(b_{\relbit{\sigma}})=g_{\relbit{\sigma}}$.
					Now consider the case $\relbit{\sigma}=2$, implying $\valustar_\sigma^\P(b_2)=L_2^\P$.
					Since $\sigma(b_2)=g_2$ by \Pref{EB5}, the statement follows by the same arguments.
			\end{enumerate}
		\item \boldall{$\sigmabar(eb_{i,j}),\sigmabar(eg_{i,j})$ and $\sigmabar(g_1)\neq\sigmabar(b_2)$:}
			We prove $\valustar_\sigma^\P(F_{i,j})=\valustar_\sigma^\P(g_1)$.
			To simplify the proof, we show $\valustar_\sigma^\P(g_1)\lhd\valustar_\sigma^\P(b_2)$ first. 
			If $\valustar_{\sigma}^\P(F_{1,\sigmabar(g_1)})=\valustar_{\sigma}^\P(b_2)$, then the statement follows from $\valustar_\sigma^\P(g_1)=\{g_1\}\cup\valustar_{\sigma}^{\P}(F_{1,\sigmabar(g_1)})=\{g_1\}\cup\valustar_\sigma^\P(b_2)$.
			In addition, since $S_n$ is a sink game, the chosen cycle center of level $1$ cannot escape towards $g_1$ since this would close a cycle.
			If $\valustar_{\sigma}^{\P}(F_{i,\sigmabar(g_1)}=\valustar_{\sigma}^{\P}(g_1)$, then the claim also follows as player $1$ always minimizes the valuations and could choose vertex $b_2$ but prefers $g_1$.
			Thus let $\tau^{\sigma}(F_{1,\sigmabar(g_1)})=s_{1,\sigmabar(g_1)}$ and assume $\sigma(s_{1,\sigmabar(g_1)})=b_1$.
			Then, since $S_n$ is a sink game, we need to have $\sigma(b_1)=b_2$ since there would be a cycle otherwise.
			But then $\valustar_\sigma^\P(g_1)=\{g_1,s_{1,\sigmabar(g_1)}\}\cup\valustar_\sigma^\P(b_2)\lhd\valustar_\sigma^\P(b_2)$. 
			We can therefore assume $\sigma(s_{1,\sigmabar(g_1)})=h_{1,\sigmabar(g_1)}$ and distinguish two cases.
			\begin{itemize}
				\item Let $\sigma(g_1)=F_{1,0}$. 
					Then $\sigma(b_2)=g_2$.
					Therefore, by \Pref{EBG4}, $\relbit{\sigma}\leq 2$.
					Thus $\valustar_\sigma^\P(b_2)=L_2^\P$ and $\valustar_\sigma^\P(b_3)=L_3^\P$, hence \[\valustar_\sigma^\P(g_1)=W_1^\P\cup\valustar_\sigma^\P(b_3)=W_1^\P\cup L_3^\P\lhd W_2^\P\cup L_3^\P=L_2^\P=\valustar_\sigma^\P(b_2).\]
				\item Let $\sigma(g_1)=F_{1,1}$.
					Then, $\valustar_\sigma^\P(g_1)=W_1^\P\cup\valustar_\sigma^\P(g_2)$ and $\sigma(b_2)=b_3$.
					This implies $\relbit{\sigma}\neq 2$ by \Pref{EBG5}.
					Consider the case $\relbit{\sigma}>2$ first.
					\Cref{lemma: Valuation of g if level small} then implies $\valustar_\sigma^\P(b_2)=L_2^\P$ as well as $\valustar_\sigma^\P(g_2)=R_2^\P$.
					Also, $\sigma(b_2)=b_3$ and $\relbit{\sigma}>2$ together imply  $\relbit{\sigma}\neq\min\{i'\in[n]\colon\sigma(b_{i'})=b_{i'+1}\}$.
					Thus, $\incorrect{\sigma}\neq\emptyset$ and $\sigma(b_{\relbit{\sigma}})=g_{\relbit{\sigma}}$ by \Cref{lemma: Traits of Relbit}.
					Combining all of this then yields \begin{align*}
						\valustar_\sigma^\P(g_1)&=W_1^\P\cup R_2^\P\lhd W_{\relbit{\sigma}}^\P\cup L_{\relbit{\sigma}+1}=L_{\relbit{\sigma}}^{\P}\unlhd L_{2}^{\P}=\valustar_\sigma^\P(b_2).
					\end{align*}
					
					Now consider the case $\relbit{\sigma}=1$.
					Then again $\valustar_\sigma^\P(b_2)=L_2^\P$.
					We apply \Cref{lemma: Valuation of g1} to give the exact valuation of $g_1$.
					The case $\minsig{b}\leq\minnegsig{s},\minnegsig{g}$ cannot occur as this would imply $\sigma(b_2)=g_2$ by \Cref{lemma: Config implied by Aeb}.				
				
					Consider the case $\minnegsig{g}<\minnegsig{s},\minnegsig{b}$.
					As $\sigma(g_1)=F_{1,1}$, it holds that $\minnegsig{g}\neq1$.
					Let $i\coloneqq\minnegsig{g}$.
					Thus, by assumption, $\sigma(g_{i-1})=F_{i-1,1}$ and consequently $\sigma(s_{i-1,1})=h_{i-1,1}$ as well as $\sigma(b_{i})=b_{i+1}$.
					Thus, \Pref{B3} implies $0=\sigmabar(g_{i})\neq\sigmabar(b_{i+1})$, hence $\sigmabar(b_{i+1})=1$.
					But then, $\sigmabar(b_{\minnegsig{g}+1})$, so \Cref{lemma: Valuation of g1} yields \begin{align*}
						\valustar_{\sigma}(g_1)&=\bigcup_{i'=1}^{\minnegsig{g}}W_{i'}^\P\cup\valustar_{\sigma}^\P(b_{\minnegsig{g}+2})\\
							&\lhd W_{\minnegsig{g}+1}\cup\valustar_{\sigma}^\P(b_{\minnegsig{g}+2})=\valustar_{\sigma}^\P(b_{\minnegsig{g}+1})\unlhd\valustar_{\sigma}^\P(b_2).
					\end{align*}

			\end{itemize}
			Thus $\valustar_\sigma^\P(g_1)\lhd\valustar_\sigma^\P(b_2)$.
			We next prove that we have $\valustar_\sigma^\P(g_1)\lhd\valustar_\sigma^\P(s_{i,j})$, implying that player $1$ chooses to escape to $g_1$.
			\begin{enumerate}
				\item Let $\sigma(s_{i,j})=b_1$.
					If $\sigma(b_1)=g_1$,  then $\valustar_\sigma^\P(s_{i,j})=\{s_{i,j}\}\cup\valustar_\sigma^\P(g_1)$ implies the statement.
					If $\sigma(b_1)=b_2$ we have $\valustar_\sigma^\P(b_1)=\valustar_\sigma^\P(b_2)$.
					The statement then follows since $\valustar_\sigma^\P(g_1)\lhd\valustar_\sigma^\P(b_2)$ and $\valustar_\sigma^\P(s_{i,j})=\{s_{i,j}\}\cup\valustar_\sigma^\P(b_1)$.
				\item Let $\sigma(s_{i,j})=h_{i,j}$.
					Then, \Pref{EBG1} implies $\sigmabar(b_{i+1})=j$ and thus $\valustar_\sigma^\P(s_{i,j})=\{s_{i,j},h_{i,j}\}\cup\valustar_\sigma^\P(b_{i+1})$ .
					Let $\relbit{\sigma}=1$.
					Then, $\sigma(b_1)=b_2$, implying $\valustar_\sigma^\P(b_1)=\valustar_\sigma^\P(b_2)$ and $\valustar_\sigma^\P(b_{i+1})=L_{i+1}^\P$.
					Combining this with $\valustar_\sigma^\P(g_1)\lhd\valustar_\sigma^\P(b_2)$ yields the statement as	\[\valustar_\sigma^\P(g_1)\lhd\valustar_\sigma^\P(b_2)=L_2^\P=\lhd\{s_{i,j},h_{i,j}\}\cup L_{i+1}^\P=\valustar_\sigma^\P(s_{i,j}).\]
					Now let $\relbit{\sigma}\neq 1$, implying $\valustar_\sigma^\P(g_1)=R_1^\P$ by \Cref{lemma: Valuation of g if level small}.
					Consider the case $\relbit{\sigma}\geq i+1$ and $B_{i+1}^\P=L_{i+1}^\P$ first.
					Then, since $\sigma(b_{i+1})=b_{i+2}$ by $B_{i+1}^\P=L_{i+1}^\P$, we have $\relbit{\sigma}\neq i+1$.
					This implies $\relbit{\sigma}\neq\min(\{i'\colon\sigma(b_{i'})=b_{i'+1}\})$ and consequently $\sigma(b_{\relbit{\sigma}})=g_{\relbit{\sigma}}$.
					The statement then follows from \Cref{lemma: VV Lemma}~(2) and $\valustar_{\sigma}^{\P}(s_{i,j})\rhd L_{i+1}^{\P}$.
					If $\relbit{\sigma}\geq i+1$ and $B_{i+1}^\P=R_{i+1}^\P$ then the statement follows since $R_1^{\P}\lhd \{s_{i,j},h_{i,j}\}\cup R_{i+1}^{\P}$ in this case.
					If $\relbit{\sigma}<i+1$, then the statement follows from \Cref{lemma: VV Lemma}~(3).
			\end{enumerate}
		\item \boldall{$\sigmabar(eb_{i,j}),\sigmabar(eg_{i,j})$ and $\sigmabar(g_1)=\sigmabar(b_2)$:}
			We prove $\valustar_\sigma^\P(F_{i,j})=\valustar_\sigma^\P(b_2)$.
			Similar to the last case we prove $\valustar_\sigma^\P(b_2)\lhd\valustar_\sigma^\P(g_1)$ and $\valustar_\sigma^\P(b_2)\lhd\valustar_\sigma^\P(s_{i,j})$.
			
			The assumption $\sigmabar(g_1)=\sigmabar(b_2)$ implies $\sigma(h_{1,\sigmabar(g_1)})=\sigma(b_2)$.
			By \Pref{EBG3}, the chosen cycle center of level 1 is closed.
			In addition, $\sigmabar(s_1)$ by \Pref{EBG2}.
			Hence, $\valustar_\sigma^\P(g_1)=W_1^{\P}\cup\valustar_\sigma^\P(b_2)$, implying $\valustar_\sigma^\P(b_2)\lhd\valustar_\sigma^\P(g_1)$.
			
			It remains to show $\valustar_\sigma^\P(b_2)\lhd\valustar_\sigma^\P(s_{i,j})$.
			Let $\sigma(s_{i,j})=b_1$ first.
			If  $\sigma(b_1)=b_2$, then the statement follows from $\valustar_\sigma^\P(s_{i,j})=\{s_{i,j}\}\cup\valustar_\sigma^\P(b_2)$.
			Thus let $\sigma(b_1)=g_1$, implying $\valustar_\sigma^\P(s_{i,j})=\{s_{i,j}\}\cup\valustar_\sigma^\P(g_1)$.
			But this implies $\valustar_\sigma^\P(s_{i,j})\rhd\valustar_\sigma^\P(g_1)$ and consequently also $\valustar_\sigma^\P(s_{i,j})\rhd\valustar_\sigma^\P(b_2)$.
			Thus, let $\sigma(s_{i,j})=h_{i,j}$, implying $\sigmabar(b_{i+1})=j$ by \Pref{EBG1}.
			We distinguish two cases. 
			\begin{enumerate}
				\item Let $j=0$.
					Then $\valustar_\sigma^\P(s_{i,j})=\{s_{i,0},h_{i,0}\}\cup B_{i+2}^\P$.
					We first consider the case that $B_{i+2}^\P=L_{i+2}^\P$ and show $\{s_{i,0},h_{i,0}\}\cup L_{i+2}^\P\rhd L_2^\P,R_2^\P$ since this suffices to show $\valustar_\sigma(s_{i,j})\rhd\valustar_\sigma(b_2)$.
					Since  $\sigma(b_{i+1})=b_{i+2}$ by \Pref{EBG1} and $j=0$, $i\geq 2$ implies \begin{align*}
						\{s_{i,0},h_{i,0}\}\cup L_{i+2}^\P=\rhd\bigcup_{i'=2}^{i}\{W_{i'}^\P\colon\sigma(b_{i'})=g_{i'}\}\cup\bigcup_{i'\geq i+1}\{W_{i'}^\P\colon\sigma(b_{i'})=g_{i'}\}=L_2^\P
					\end{align*}
					and $i=1$ implies $\{s_{i,0},h_{i,0}\}\cup L_{i+2}^\P=\{s_{1,0},h_{1,0}\}\cup L_3^\P=\{s_{1,0},h_{1,0}\}\cup L_2^\P\rhd L_2^\P.$
					
					Thus let $\valustar_\sigma^\P(b_2)=R_2^\P$ and consider the case $\relbit{\sigma}\leq i+1$ first.
					Then, by \Cref{lemma: VV Lemma} and since $\sigma(b_{i+1})=b_{i+2}$, it holds that $L_{i+1}^\P=L_{i+2}^\P$ and thus $R_2^\P\lhd R_1^\P\lhd\{s_{i,0},h_{i,0}\}\cup L_{i+1}^\P.$ 
					Now assume $\relbit{\sigma}>i+1$, implying $\relbit{\sigma}\neq 1$ and $i+2\leq\relbit{\sigma}$.
					This implies $\relbit{\sigma}\neq\min\{i'\colon\sigma(b_{i'})=b_{i'+1}\}$ since \Pref{EBG1} implies $\sigma(b_{i+1})=b_{i+2}$.
					Thus $\sigma(b_{\relbit{\sigma}})=g_{\relbit{\sigma}}$ by \Cref{lemma: Traits of Relbit}.
					Then, the statement follows from $\{s_{i,0},h_{i,0}\}\cup L_{i+2}^{\P}\rhd L_{i+2}^{\P}$ and $L_{i+2}^{\P}\rhd R_2^{\P}$ which follows from \Cref{lemma: VV Lemma}~(2).

					Next, let $B_{i+2}^\P=R_{i+2}^\P$.
					We show that this results in a contradiction.
					First, $B_{i+2}^\P=R_{i+2}^\P$ implies $i+2<\relbit{\sigma}$ and $\sigma(b_{i+2})=g_{i+2}$.
					In particular we have $\relbit{\sigma}\geq 4$, implying $\relbit{\sigma}-1\geq 3$.
					But then \Pref{BR1} implies $\sigma(g_1)=F_{1,1}$ which implies $\sigma(b_2)=g_2$ by assumption.
					Now consider level $i$.
					Again, by \Pref{BR1}, $\sigma(g_i)=F_{i,1}$.
					Now, combining all of this and using \Pref{S2} yields $\sigma(s_{i,1})=h_{i,1}$.
					But then, since we have $\sigma(b_{i+1})=b_{i+2}$ by assumption, \Pref{B3} now implies $\sigmabar(g_{i+1})\neq\sigmabar(b_{i+2})$.
					Since $i+1<\relbit{\sigma}-1$, \Pref{BR1} now implies $\sigma(g_{i+1})=F_{i+1,1}$, i.e., we have $\sigmabar(g_{i+1})=1$.
					But this now implies $\sigmabar(b_{i+2})=0$, i.e., $\sigma(b_{i+2})=b_{i+3}$ which is a contradiction since $\sigma(b_{i+2})=g_{i+2}$ by $B_{i+2}^\P=R_{i+2}^\P$.
				\item Let $j=1$.
				Then, $\sigma(b_{i+1})=g_{i+1}$, implying $\valustar_\sigma^\P(s_{i,j})=\{s_{i,1},h_{i,1}\}\cup B_{i+1}^\P$.
				We now show $\{s_{i,1},h_{i,1}\}\cup B_{i+1}^\P\rhd B_2^\P$ for all possible ``choices'' of $B_{i+1}^\P$ and $B_2^\P$.
				Let $B_{i+1}^\P=L_{i+1}^\P$ and $B_2^\P=L_2^\P$. 
				Then $\{s_{i,1},h_{i,1}\}\cup L_{i+1}^\P\rhd L_2^\P,$ so the statement holds.
				Now consider the case $B_2^\P=R_2^\P$, implying that $2<\relbit{\sigma}$.
				First assume that $\relbit{\sigma}\leq i$.
				Then, $\{s_{i,1},h_{i,1}\}\cup L_{i+1}^{\P}\rhd R_2^{\P}$ follows from \Cref{lemma: VV Lemma}~(3).
				It cannot happen that$\relbit{\sigma}>i$, as this would yield $\relbit{\sigma}\geq i+1$.
				But this is a contradiction as this would imply $B_{i+1}^{\P}=R_{i+1}^{\P}$ as we currently assume $B_{i+2}^{\P}=L_{i+2}^{\P}$.
				
				Now consider the case $B_{i+1}^\P=R_{i+1}^\P$ and $B_2^\P=R_2^\P$.
				Then $i+1<\relbit{\sigma}$, hence the statement follows from $\{s_{i,1},h_{i,1}\}\rhd\bigcup_{i'<i}W_{i'}$.
				Finally assume $B_{i+1}^\P=R_{i+1}^\P$ and $B_2^\P=L_2^\P$.
				Since $\relbit{\sigma}>i+1\geq 2$ it holds that $\sigma(b_2)=b_3$, implying $B_2^\P=B_3^\P$.
				Applying \Pref{B1} repeatedly thus yields $B_2^\P=B_{k}^\P=R_{k}^\P$ where $k=\min\{i'\in\{2,\dots,i+1\}:\sigma(b_{i'})=g_{i'}\}\leq i+1$.
				Thus, the statement follows from $\{s_{i,1},h_{i,1}\}\rhd\bigcup_{i'<i} W_{i'}$ resp. $\{s_{i,j},h_{i,1}\}\rhd\emptyset$.  \qedhere
			\end{enumerate}
	\end{enumerate}
\end{proof}

\subsection{Omitted proofs of \Cref{section: Improving switches technical}}

Here, we provide the formal proofs of all statements of \Cref{section: Improving switches technical} that have not been proven there.

\BothCCOpenForMDP*

\begin{proof}
To simplify notation let $j\coloneqq\indbit_{i+1}$.
Since both cycle centers are in the same state, it suffices to prove $\valustar_{\sigma}^\M(s_{i,j})>\valustar_{\sigma}^\M(s_{i,1-j})$.
By \Pref{USV1}$_i$ and \Pref{EV1}$_{i+1}$, $\valustar_{\sigma}^\M(s_{i,1-j})=\rew{s_{i,1-j}}+\valustar_{\sigma}^\M(b_1)$ and $\valustar_{\sigma}^\M(s_{i,j})=\rew{s_{i,j},h_{i,j}}+\valustar_{\sigma}^\M(b_{i+1})$.
Let $\relbit{\sigma}=1$.
Then, $\sigma(b_1)=b_2$ and thus $\valustar_{\sigma}^\M(b_1)=\valustar_{\sigma}^\M(b_2)=L_2^\M$ and $\valustar_{\sigma}^\M(b_{i+1})=L_{i+1}^\M$.
The statement then follows from $\rew{s_{i,j},h_{i,j}}>\sum_{\ell\leq i}W_{\ell}^{\M}+\rew{s_{i,1-j}}$.

Hence let $\relbit{\sigma}>1$, implying $\sigma(b_1)=g_1$.
We distinguish the following cases.
\begin{enumerate}
	\item Let $\valustar_{\sigma}^\M(b_1)=R_1^\M$ and $\valustar_{\sigma}^\M(b_{i+1})=R_{i+1}^\M$.
		This implies $i+1<\relbit{\sigma}$ and the statement thus again follows from $\rew{s_{i,j},h_{i,j}}>\sum_{\ell\leq i}W_{\ell}^{\M}+\rew{s_{i,1-j}}$. 
	\item Let $\valustar_{\sigma}^\M(b_1)=R_1^\M$ and $\valustar_{\sigma}^\M(b_{i+1})=L_{i+1}^\M$.
		\Pref{EV1}$_{i+1}$ implies \[\valu_{\sigma}^\M(s_{i,j})=\rew{s_{i,j},h_{i,j}}+L_{i+1}^\M>\sum_{\ell\leq i}W_{\ell}^\M+\rew{s_{i,1-j}}+L_{i+1}^\M.\]
		If $\relbit{\sigma}\leq i$, then \[\sum_{\ell\leq i}W_{\ell}^\M=\sum_{\ell\leq\relbit{\sigma}}W_{\ell}^\M+\sum_{\ell=\relbit{\sigma}+1}^{i}W_{\ell}^\M>\sum_{\ell<\relbit{\sigma}}W_{\ell}^\M+\sum_{\ell=\relbit{\sigma}+1}^{i}\{W_{\ell}^\M\colon\sigma(b_i)=g_i\},\] implying $\valu_{\sigma}^\M(s_{i,j})>\valu_{\sigma}^\M(s_{i,1-j})$.
		If $\relbit{\sigma}=i+1$, then $\sum_{\ell\leq i}W_{\ell}^\M=\sum_{\ell<\relbit{\sigma}}W_{\ell}^\M$, again implying $\valu_{\sigma}^\M(s_{i,j})>\valu_{\sigma}^\M(s_{i,1-j})$.
		Let $\relbit{\sigma}>i+1$.
		Then, by assumption, it needs to hold that $\sigma(b_{i+1})=b_{i+2}$.
		Thus $\relbit{\sigma}\neq\min\{i'\colon\sigma(b_{i'})=b_{i'+1}\}$, implying $\sigma(b_{\relbit{\sigma}})=g_{\relbit{\sigma}}$ by \Cref{lemma: Traits of Relbit}.
		The statement the follows from \begin{align*}
			\sum_{\ell\leq i} W_{\ell}^\M+L_{i+1}^\M&=\sum_{\ell\leq i} W_{\ell}^\M+\sum_{\ell=i+1}^{\relbit{\sigma}-1}\{W_{\ell}^\M\colon\sigma(b_{\ell})=g_{\ell}\}+W_{\relbit{\sigma}}^\M+L_{\relbit{\sigma}+1}^\M\\
				&>\sum_{\ell<\relbit{\sigma}}W_{\ell}^\M+L_{\relbit{\sigma}+1}^\M=R_1^\M.
		\end{align*}
	\item Let $\valustar_{\sigma}^\M(b_1)=\rew{g_k}+\sum_{\ell<k}W_{\ell}^\M+\valustar_{\sigma}^\M(b_2), k=\min\{i'\colon\nsigmabar(d_{i'})\}<\relbit{\sigma},$ and $\valustar_{\sigma}^\M(b_{i+1})=R_{i+1}^\M$.
		We show that these assumptions yield a contradiction.
		The second equality implies that $i+1<\relbit{\sigma}$.
		We now prove that $\sigma$ has \Pref{REL1} in any case, so assume that $i\geq\nsb$.
		This implies $\nsb<i+1<\relbit{\sigma}$ and thus $\nsb\neq\relbit{\sigma}$.
		Consequently, $\sigma$ cannot be a phase-$2$-strategy or phase-$3$-strategy for $\bit$ as it then had \Pref{REL2}, implying $\relbit{\sigma}=\nsb$.
		Therefore, by the definition of the phases, $\sigma$ has \Pref{REL1} in any case, so $\relbit{\sigma}=\min\{i'\colon\sigma(b_{i'})=b_{i'+1}\}$.
		Consequently, it holds that $\incorrect{\sigma}=\emptyset$.
		But then $i'<\relbit{\sigma}$ implies $\sigmabar(d_{i'})$ by \Cref{corollary: Simplified MDP Valuation}.
		This contradicts the characterization of $\valustar_{\sigma}^\M(b_1)$.
	\item Let $\valustar_{\sigma}^\M(b_1)=\rew{g_k}+\sum_{\ell<k}W_{\ell}^\M+\valustar_{\sigma}^\M(b_2), k=\min\{i'\colon\nsigmabar(d_{i'})\}<\relbit{\sigma},$ and $\valustar_{\sigma}^\M(b_{i+1})=L_{i+1}^\M$.
		Then \[\valustar_{\sigma}^\M(s_{i,j})=\rew{s_{i,j},h_{i,j}}+L_{i+1}^\M>\sum_{\ell\leq i}W_{\ell}^\M+L_{i+1}^\M+\rew{s_{i,1-j}}\geq\rew{s_{i,1-j}}+L_2^\M.\]
		If $\valustar_{\sigma}^\M(b_2)=L_2^\M$, the statement thus follows since  $\valustar_{\sigma}^\M(b_1)<\valustar_{\sigma}^\M(b_2)$ in this case.
		Thus assume $\valustar_{\sigma}^\M(b_2)=R_2^\M$, implying $\sigma(b_2)=g_2$ and $\relbit{\sigma}>2$.
		If $\sigma(b_{\relbit{\sigma}})=g_{\relbit{\sigma}}$, then $\rew{s_{i,1-j}}+L_2^\M>R_2^\M$, implying the statement.
		Hence assume $\sigma(b_{\relbit{\sigma}})=b_{\relbit{\sigma}+1}$, implying $\relbit{\sigma}=\min\{i'\colon\sigma(b_{i'})=b_{i'+1}\}$.
		In particular, $\sigma(b_{i'})=g_{i'}$ for all $i'\in\{1,\dots,\relbit{\sigma}-1\}$.
		Thus, since $\valustar_{\sigma}^\M(b_{i+1})=L_{i+1}^\M$, we need to have $i+1\geq\relbit{\sigma}$ and in particular $i\geq\relbit{\sigma}-1$.
		Then, the statement follows from \begin{align*}
			\valustar_{\sigma}^\M(s_{i,j})&>\sum_{\ell\leq i}W_{\ell}^\M+L_{i+1}^\M+\rew{s_{i,1-j}}\\
				&\geq\sum_{\ell<\relbit{\sigma}}W_{\ell}^\M+L_{\relbit{\sigma}+1}^\M+\rew{s_{i,1-j}}=R_2^\M+\rew{s_{i,1-j}}. 
		\end{align*}\qedhere
\end{enumerate}
\end{proof}

\NumericsOfEll*

\begin{proof}
Let $m\coloneqq\id_{j=0}\lastflip{\bit}{i+1}{}+\id_{j=1}\lastunflip{\bit}{i+1}{}\neq 0$.
Then, $\lastflip{\bit}{i}{\{(i+1,j)\}}\neq 0$ and we distinguish three cases.
\begin{enumerate}[after=\smallskip]
	\item Let $\bit_i=1$ and $\bit_{i+1}=1-j$.
		We prove $\lastflip{\bit}{i}{\{(i+1,j)\}}=\bit-2^{i}-\sum(\bit,i)$ and $m=\bit-2^{i-1}-\sum(\bit,i)$.
		By definition, $\bit'\coloneqq\lastflip{\bit}{i}{\{(i+1,j)\}}$ is the largest number smaller than $\bit$ such that $\nsb(\bit')=i$ and $\bit'_{i+1}=j$.
		Since $\bit_{i+1}=1-j$, subtracting $2^{i}$ switches bit $i+1$ and only bit $i+1$.
		By subtracting $\sum(\bit,i)$, all bits below bit $i$ that are equal to $1$ are set to $0$.
		Therefore, $\bit'=\bit-2^{i}-\sum(\bit,i)$.
		Note that $\bit'>0$.
		
		Assume $j=1$, implying $m=\lastflip{\bit}{i+1}{}\neq 0$.
		Since $m$ is the largest number smaller than $\bit$ with least significant set bit equal to $1$ being bit $i+1$ and since $\bit_{i+1}=1$, we have $m=\bit-\sum(\bit,i+1)=\bit-2^{i-1}-\sum(\bit,i).$
		Consequently, \begin{align*}
			\ell^{\bit}(i,j,k)&=\ceil{\frac{\bit-2^i-\sum(\bit,i)+1-k}{2}}+\bit-\bit+2^{i-1}+\sum(\bit,i)\\
				&=\ceil{\frac{\bit+\sum(\bit,i)+1-k}{2}}.
			\end{align*}
	\item By similar arguments, it can be shown that $\lastflip{\bit}{i}{\{(i+1,j)\}}=\bit-2^{i}-2^{i-1}-\sum(\bit,i)$ and $m=\bit-2^{i}-\sum(\bit,i)$ in this case, implying the statement analogously.
	\item By similar arguments, it can be shown that $\lastflip{\bit}{i}{\{(i+1,j)\}}=\bit-2^{i-1}-\sum(\bit,i)$ and $m=\bit-\sum(\bit,i)$ in this case, implying the statement analogously.
\end{enumerate}
If $\id_{j=0}\lastflip{\bit}{i+1}{}+\id_{j=1}\lastunflip{\bit}{i+1}{}=0$, the statement follows immediately.
\end{proof}

\ProgressForUnimportantCC*

\begin{proof}
Consider the first statement.
Assume $\lastunflip{\bit}{i+1}{}\neq\lastunflip{\bit+1}{i+1}{}$.
This can only occur if $\bit+1=\lastunflip{\bit+1}{i+1}{}$, implying $(\bit+1)_{i+1}=\dots=(\bit+1)_{1}=0$.
But this implies $\nsb(\bit+1)\geq i+2$, hence $\bit_1=\dots=\bit_i=\bit_{i+1}=1$.
Since $\bit_i=0\vee\bit_{i+1}\neq j$ by assumption, it thus needs to hold that $j=0$.
This proves that $\lastunflip{\bit}{i+1}{}\neq\lastunflip{\bit+1}{i+1}{}$ implies $j=0$.
In a similar way it can be proven that $\lastflip{\bit}{i+1}{}\neq\lastflip{\bit+1}{i+1}{}$ implies $j=1$.
Consequently, it is impossible that both $\lastunflip{\bit}{i+1}{}\neq\lastunflip{\bit+1}{i+1}{}$ and $\lastflip{\bit}{i+1}{}\neq\lastunflip{\bit+1}{i+1}{}$ hold.
If $\lastunflip{\bit}{i+1}{}\neq\lastunflip{\bit+1}{i+1}{}$, then it holds that $j=0$ and $\lastflip{\bit}{i+1}{}=\lastflip{\bit+1}{i+1}{}$.
If $\lastflip{\bit}{i+1}{}\neq\lastflip{\bit+1}{i+1}{}$, we have $j=1$ and $\lastunflip{\bit}{i+1}{}=\lastunflip{\bit+1}{i+1}{}$.
But this implies \[\id_{j=0}\lastflip{\bit}{i+1}{}-\id_{j=1}\lastunflip{\bit}{i+1}{}=\id_{j=0}\lastflip{\bit+1}{i+1}{}-\id_{j=1}\lastunflip{\bit+1}{i+1}{}.\]
Now let also $i\neq\nsb$.
It suffices to prove $\lastflip{\bit}{i}{\{(i+1,j)\}}=\lastflip{\bit+1}{i}{\{(i+1,j)\}}$.
But this follows directly since the choice of $i$ implies $\lastflip{\bit+1}{i}{\{(i+1,j)\}}\neq\bit+1$.
\end{proof}

\OccurrenceRecordsCycleVertices*

\begin{proof}
The first statement follows immediately since $\bit_{i}=0\vee\bit_{i+1}\neq j$ imply \[\occrec^{\canstrat}(d_{i,j,k},F_{i,j})=\min\left(\floor{\frac{\bit+1-k}{2}},\ell^{\bit}(i,j,k)+t_{\bit}\right)\leq\floor{\frac{\bit+1-k}{2}}\leq\floor{\frac{\bit+1}{2}}.\]
Consider the second statement and observe that it suffices to prove \begin{equation} \label{equation: OR Cycle Edges}
\occrec^{\canstrat}(d_{\nsb,j,k},F_{\nsb,j})=\floor{\frac{\bit+1-k}{2}}
\end{equation} for $k\in\{0,1\}$.
Let $\id_{j=0}\lastflip{\bit}{\nsb+1}{}-\id_{j=1}\lastunflip{\bit}{\nsb+1}{}=0$.
Then, by \Cref{lemma: Numerics Of Ell}, $\ell^{\bit}(\nsb,j,k)\geq\bit$.
In order to show \Cref{equation: OR Cycle Edges}, it thus suffices to prove that either $\bit-1\geq\floor{(\bit+1-k)/2}$ or that the parameter $t_{\bit}=-1$ is not feasible.
Since it holds that $\bit-1\geq\floor{(\bit+1)/2}$ for $\bit\geq 2$, it suffices to show that $t_{\bit}=-1$ is not feasible for $\bit=0,1$.
By \Cref{table: Occurrence Records}, the parameter $-1$ can only be feasible if $\bit_1=1\wedge\nsb\neq 1$.
It is therefore not feasible for $\bit=0$.
Let $\bit=1$ and $\occrec^{\canstrat}(d_{\nsb,j,k},F_{\nsb,j})=\ell^{\bit}(\nsb,j,k)-1$.
Since $\bit+1=2$ is a power of two and since $\canstrat$ has Properties (\ref{property: OR1})$_{*,*,*}$ to (\ref{property: OR4})$_{*,*,*}$, \Pref{OR3}$_{\nsb,j,k}$ implies $\occrec^{\canstrat}(d_{\nsb,j,k},F_{\nsb,j})=\floor{(\bit+1-k)/2}$.
Consequently, $\occrec^{\canstrat}(d_{\nsb,j,k},F_{\nsb,j})=\floor{(\bit+1-k)/2}$.
Now let $\id_{j=0}\lastflip{\bit}{\nsb+1}{}-\id_{j=1}\lastunflip{\bit}{\nsb+1}{}\neq 0$.
Then, by the definition of $\nsb$ and $j$ and \Cref{lemma: Numerics Of Ell}, \begin{align*}
\ell^{\bit}(\nsb,j,k)&=\ceil{\frac{\bit+2^{\nsb-1}+\sum(\bit,\nsb)+1-k}{2}}\geq\ceil{\frac{\bit+2^{\nsb-1}+1-k}{2}}\geq\floor{\frac{\bit+1-k}{2}}+1.
\end{align*}
Since $-1,0$ and $1$ are the only feasible parameters, this implies \Cref{equation: OR Cycle Edges}.

Consider the third statement and let $i=1, j=1-\bit_{2}$.
Then, independent of whether $\bit_1=0$ or $\bit_1=1$, $\ell^{\bit}(i,j,k)\geq\ceil{(\bit-k)/2}=\floor{(\bit+1-k)/2}$ by \Cref{lemma: Numerics Of Ell}.
By the first statement and by \Pref{OR1}$_{i,j,k}$ and \Pref{OR2}$_{i,j,k}$, this implies $\sigma(d_{i,j,k})\neq F_{i,j}$ for both $k\in\{0,1\}$ as $\occrec^{\canstrat}(d_{i,j,k},F_{i,j})=\ell^{\bit}(i,j,k)+1$ otherwise.
Furthermore, this implies $\ell^{\bit}(i,j,0)=\floor{(\bit+1)/2}$.
Assume that $\occrec^{\sigma}(d_{i,j,0},F_{i,j})=\ell^{\bit}(i,j,0)-1$.
Then, $\occrec^{\canstrat}(d_{i,j,0},F_{i,j})\neq \floor{(\bit+1)/2}$.
Hence, by \Pref{OR4}$_{i,j,k}$, it holds that $\bit$ is odd and $i=\nsb(\bit+1)$.
But this contradicts $i=1$.

Consequently, $\occrec^{\canstrat}(d_{i,j,0},F_{i,j})=\ell^{\bit}(i,j,0)=\floor{(\bit+1)/2}$.
\end{proof}

\NumericsOfOR*

\begin{proof}
As a reminder, a binary number $\bit$ \emph{matches} the pair $(i,q)$ if $\bit_i=q$.
It matches a set $S$ if $\bit$ matches every $(i,q)\in S$.
Consider the first two statements.
By definition, it holds that  $\flips{\bit+1}{i}{}=\flips{\bit}{i}{}+\id_{i=\nsb}$.
Let $S_i\coloneqq \{(i,1),(i-1,0),\dots,(1,0)\}$. 
By definition, $\flips{\bit}{i}{}$ is the number of numbers smaller than or equal to $\bit$ matching $S_i$.
Since $2^{i-1}$ is the smallest number matching $S_i$, the statement follows if $\bit<2^{i-1}$.
Let $m_k$ denote the $k$-th number matching the scheme $S_i$.
Then $m_1=2^{i-1}$.
As only numbers ending on the subsequence $(1,0,\dots,0)$ of length $i$ match $S_i$, we have $m_k=(k-1)\cdot 2^{i}+2^{i-1}$.
Since $\flips{m_k}{i}{}=k$ by definition and \[\floor{\frac{m_k+2^{i-1}}{2^i}}=\floor{\frac{(k-1)\cdot 2^i+2^{i-1}+2^{i-1}}{2^i}}=\floor{\frac{k\cdot 2^k}{2^k}}=k,\] this implies $\flips{m_k}{i}=\floor{(m_k+2^{i-1})/2^i}$.
Let $\bit\in\bitset_n$ and $k\in\mathbb{N}$ such that $\bit\in[m_k,m_{k+1})$.
Then, by the definition of $\flips{\bit}{i}{}$, we have $\flips{\bit}{i}{}=k$.
In addition, \[\floor{\frac{\bit+2^{i-1}}{2^i}}\geq\floor{\frac{m_k+2^{i-1}}{2^i}}=\flips{m_k}{i}{}=k\] by the choice of $k$ and \[\floor{\frac{\bit+2^{i-1}}{2^i}}<\floor{\frac{m_{k+1}+2^{i-1}}{2^i}}=\flips{m_{k+1}}{i}{}=k+1.\]
Integrality thus implies $\floor{(\bit+2^{i-1})/2^i}r=k$, hence $\flips{\bit}{i}{}=k=\floor{(\bit+2^{i-1})/2^i}$.

Now let $i_1,i_2\in[n]$ with $i_1<i_2$ and $\bit\geq2^{i_1-1}$.
Then, $\flips{\bit}{i_1}{}=\floor{(\bit+2^{i_1-1})/2^{i_1}}$ and similarly $\flips{\bit}{i_2}{}=\floor{(\bit+2^{i_2-1})/2^{i_2}}$.
If $\bit=2^{i_1-1}$, then $\bit<2^{i_2-1}$, implying \[\flips{\bit}{i_1}{}=\floor{\frac{\bit}{2^{i_1}}+\frac{1}{2}}=\floor{\frac{2^{i_1-1}}{2^{i_1}}+\frac{1}{2}}=1>0=\floor{\frac{2^{i_1-1}}{2^{i_2}}+\frac{1}{2}}=\flips{\bit}{i_2}{}.\]
Thus let $\bit>2^{i_1-1}$.
Choose $k\in\mathbb{N}$ such that $k\cdot2^{i_1-1}<\bit\leq(k+1)2^{i_1-1}$.
Then \[\flips{\bit}{i_1}{}=\floor{\frac{\bit}{2^{i_1}}+\frac{1}{2}}>\floor{\frac{k\cdot2^{i_1-1}}{2^{i_1}}+\frac{1}{2}}=\floor{\frac{k+1}{2}}\geq\floor{\frac{k}{2}}\]and thus $\flips{\bit}{i_1}{}\geq\floor{k/2}+1$ by integrality.
In addition, $\bit\leq(k+1)2^{i_1-1}<(k+1)2^{i_2-1}$ implies \[\flips{\bit}{i_2}{}=\floor{\frac{\bit}{2^{i_2}}+\frac{1}{2}}<\floor{\frac{(k+1)2^{i_2-1}}{2^{i_2}}+\frac{1}{2}}=\floor{\frac{k+2}{2}}=\floor{\frac{k}{2}}+1,\]hence $\flips{\bit}{i_2}{}\leq\floor{k/2}$, implying the statement.

Now consider the third statement.
By definition, $\nsb$ is the least significant set bit of $\bit+1$.
Consequently, $\bit+1$ is dividable by $2^{\nsb-1}$, hence $k\in\mathbb{N}$ and in particular $\bit=k\cdot2^{\nsb-1}-1$.
Using \Cref{lemma: Numerics Of Ell} and $\frac{1}{2}-\frac{1}{2^{\nsb-x}}\in(0,\frac{1}{2})$, this implies \begin{align*}
	\flips{\bit}{\nsb-x}{}&=\floor{\frac{\bit+2^{\nsb-x-1}}{2^{\nsb-x}}}=\floor{\frac{k\cdot 2^{\nsb-1}-1+2^{\nsb-x-1}}{2^{\nsb-x}}}=\floor{\frac{k\cdot 2^{\nsb-1}}{2^{\nsb-x}}-\frac{1}{2^{\nsb-x}}+\frac{1}{2}}\\
		&=\floor{\frac{k\cdot 2^{\nsb-1}}{2^{\nsb-x}}}=\floor{k\cdot2^{x-1}}=k\cdot2^{x-1}. 
\end{align*}\qedhere
\end{proof}

\NorClosingCycleCenterInPhaseOne*

\begin{proof}
Since $F_{i,j}$ is open for $\sigma$, \Pref{ESC1} implies $\sigmabar(eb_{i,j})=\sigmaebar(eb_{i,j}), \sigmabar(eg_{i,j})=\sigmaebar(eg_{i,j})$ and $\sigmabar(d_{i,j})=\sigmaebar(d_{i,j})=0$.
Hence, $\sigma$ being well-behaved implies that $\sigmae$ is well-behaved.
By the same arguments, $\sigmae$ is a phase-1-strategy for $\bit$ and it suffices to prove $I_{\sigma[e]}=\mathfrak{D}^{\sigmae}$.

Consider the case $G_n=S_n$.
By \Pref{ESC1}, it holds that $\sigmabar(eb_{i,j})=\sigmaebar(eb_{i,j})$ and $\sigmabar(eg_{i,j})=\sigmaebar(eg_{i,j})$.
Since $\sigmae$ is a phase-1-strategy for $\bit$, also $\relbit{\sigma}=\relbit{\sigmae}$ by the choice of $e$.
Thus, $\valustar_{\sigma}^\P(F_{i,j})=\valustar_{\sigmae}^\P(F_{i,j})$ by \Cref{lemma: Exact Behavior Of Counterstrategy}.
In particular, the valuation of $F_{i,j}$ does not change.
Since $F_{i,j}$ is the only vertex that has an edge towards $d_{i,j,k}$, this implies that the valuation of no other vertex but $d_{i,j,k}$ changes, hence $I_{\sigmae}=I_{\sigma}\setminus\{e\}$ if $G_n=S_n$.

Consider the case $G_n=M_n$ and let $j=\sigmabar(g_i)$.
Then, $F_{i,\sigmabar(g_i)}$ is not closed with respect to either $\sigma$ or $\sigmae$.
Therefore, the valuations of $F_{i,j}$ and $g_i$ increase, but only by terms of size $o(1)$.
Now, \Pref{EV1}$_i$ and \Pref{EV2}$_i$ imply $\sigma(b_i)=b_{i+1}$,  hence $\sigma(s_{i-1,1})=b_1$ by \Pref{USV1}$_{i-1}$.
In particular, the valuation of no other vertex than $d_{i,j,k}, F_{i,j}, g_i$ and $h_{i-1,1}$ increases.
It is now easy to calculate that $(b_i,g_i), (s_{i-1,1},h_{i-1,1})\notin I_{\sigmae}$ as the change of the valuation of $g_i$ is only of size $o(1)$, implying the statement.


Let $j\neq\sigmabar(g_i)$ and let $t^{\rightarrow}\coloneqq g_1$ if $\bit$ is odd and $t^{\rightarrow}\coloneqq b_2$ if $\bit$ is even.
Then, $d_{i,j,k}$ and $F_{i,j}$ are the only vertices whose valuation increases by applying $e$.
Since $(g_i,F_{i,j})\notin I_{\sigma}$ by assumption, it thus suffices to prove $(g_i,F_{i,j})\notin I_{\sigmae}$.
By the choice of $e$, it holds that $\valu_{\sigmae}^\M(F_{i,j})=\frac{1-\e}{1+\e}\valu_{\sigmae}^\M(t^{\rightarrow})+\frac{2\e}{1+\e}\valu_{\sigma}^\M(s_{i,j})$.
First assume that $F_{i,1-j}$ is $t^{\rightarrow}$-open. 
Then, by \Cref{lemma: Both CC Open For MDP} and since $(g_i,F_{i,j})\notin I_{\sigma}$, we have $j=1-\bit_{i+1}$.
We prove that $\sigma(s_{i,j})=b_1$ (\Pref{USV1}$_i$), $\sigma(b_1)=t^{\rightarrow}$ (\Pref{EV1}$_1$ and \Pref{ESC1}) as well as $\sigma(s_{i,1-j})=h_{i,1-j}$ and $\sigmabar(b_{i+1})=1-j$ (\Pref{USV1}$_i$) imply the statement.
We have \begin{align*}
	\valu_{\sigmae}^\M(F_{i,1-j})-&\valu_{\sigmae}^\M(F_{i,j})\\
		&=\frac{\e(1-\e)}{1+\e}\valu_{\sigmae}^\M(t^{\rightarrow})+\e\valu_{\sigmae}^\M(s_{i,1-j})-\frac{2\e}{1+\e}\left(\rew{s_{i,j}}+\valu_{\sigmae}^\M(t^{\rightarrow})\right)\\
		&=\e\left(\valu_{\sigmae}^\M(s_{i,1-j})-\valu_{\sigmae}^\M(t^{\rightarrow})-\frac{2}{1+\e}\rew{s_{i,j}}\right)\\
		&>\e\left(\valu_{\sigmae}^\M(s_{i,1-j})-\valu_{\sigmae}^\M(t^{\rightarrow})-2N^{10}\right)\\
		&=\e\left(\rew{s_{i,1-j},h_{i,1-j}}+\valu_{\sigmae}^\M(b_{i+1})-\valu_{\sigmae}^\M(t^{\rightarrow})-2N^{10}\right).
\end{align*}
It thus suffices to prove $\rew{s_{i,1-j},h_{i,1-j}}+\valu_{\sigmae}^\M(b_{i+1})-\valu_{\sigmae}^\M(t^{\rightarrow})-2N^{10}\geq 0$.
We distinguish three cases.
\begin{enumerate}
	\item Let $t^{\rightarrow}=b_2$.
		Then, $\sigma(b_1)=b_2$ and $\relbit{\sigma}=1$.
		In particular, $\valu_{\sigmae}^\M(b_{i+1})=L_{i+1}^\M$ and $\valu_{\sigmae}^\M(t^{\rightarrow})=L_2^\M$.
		Consequently, \begin{align*}
			\langle s_{i,1-j},&h_{i,1-j}\rangle+\valu_{\sigmae}^\M(b_{i+1})-\valu_{\sigmae}^\M(t^{\rightarrow})-2N^{10}\\
				&=\rew{s_{i,1-j},h_{i,1-j}}+L_{i+1}^\M-L_2^\M-2N^{10}=\rew{s_{i,1-j},h_{i,1-j}}+L_{2,i}^\M-2N^{10}\\
				&>\rew{s_{i,1-j},h_{i,1-j}}-\sum_{\ell=2}^{i}W_{\ell}^\M-2N^{10}\\
				&\geq N^{2i+10}+N^8-\sum_{\ell=1}^{i}(N^{2\ell+10}-N^{2\ell+9}+N^{10})-2N^{10}\\
				&=N^{2i+10}+N^8-\frac{N^{2i+11}-N^{11}}{N+1}-(i+2)N^{10}\\				
				&>N^{2i+10}+N^8-\frac{N^{2i+11}-N^{11}}{N+1}-N^{11}.
		\end{align*}
		This term is larger than 0 if $(2N+1)N^{2i+2}+N+1>N^3(N+2)$ which holds since $i\geq 1$ and $N$ is sufficiently large.
	\item Let $t^{\rightarrow}=g_1$ and $\valu_{\sigmae}^\M(b_{i+1})=R_{i+1}^\M$.
		Then $\sigma(b_1)=g_1$ and $\valu_{\sigmae}^\M(t^{\rightarrow})=R_1^\M$.
		In particular, since $\sigma$ is a phase-1-strategy and $i+1<\relbit{\sigma}\wedge\sigma(b_{i+1})=g_{i+1}$ by assumption, it holds that $\bit_1=\dots=\bit_{i+1}=1$.
		This then implies \begin{align*}
		\langle s_{i,1-j},&h_{i,1-j}\rangle+\valu_{\sigmae}^\M(b_{i+1})-\valu_{\sigmae}^\M(t^{\rightarrow})-2N^{10}\\
			&=\rew{s_{i,1-j},h_{i,1-j}}+R_{i+1}^\M-R_1^\M-2N^{10}=\rew{s_{i,1},h_{i,1}}-\sum_{\ell=1}^{i}W_{\ell}^\M-2N^{10}\\
			&=\rew{s_{i,1},h_{i,1}}-\sum_{\ell=1}^{i}(N^{2\ell+10}-N^{2\ell+9}+N^8)-2N^{10}\\
			&>N^{2i+10}+N^8-\sum_{\ell=1}^{i}(N^{2\ell+10}-N^{2\ell+9})-N^{11}
		\end{align*}
		which is larger than 0 as shown above.
	\item Let $t^{\rightarrow}=g_1$ and $\valu_{\sigmae}^\M(b_{i+1})=L_{i+1}^\M$.
		It cannot hold that $i+1<\relbit{\sigmae}$ since this implies $\sigmae(b_{i+1})=g_{i+1}$ and thus $\valu_{\sigmae}^\M(b_{i+1})=R_{i+1}^\M$.
		Consequently, $i+1\geq\relbit{\sigmae}$.
		In addition, $\valu_{\sigmae}^\M(b_1)=R_1^\M$ as before.
		Consequently, \begin{align*}
			\rew{s_{i,1-j},h_{i,1-j}}+&\valu_{\sigmae}^\M(b_{i+1})-\valu_{\sigmae}^\M(t^{\rightarrow})-2N^{10}\\
				&=\rew{s_{i,1-j},h_{i,1-j}}+L_{i+1}^\M-R_1^\M-2N^{10}\\
				&=\rew{s_{i,1-j},h_{i,1-j}}-\hspace*{-2.75pt}\sum_{\ell=1}^{\relbit{\sigmae-1}}\hspace*{-2.75pt}W_{\ell}^\M-\hspace*{-5.5pt}\sum_{\ell=\relbit{\sigmae}+1}^{i}\hspace*{-5.5pt}\{W_{\ell}^\M\colon\sigmae(b_{\ell})=g_{\ell}\}-2N^{10}\\
				&>\rew{s_{i,1-j},h_{i,1-j}}-\sum_{\ell=1}^{i}W_{\ell}^\M-2N^{10}
		\end{align*}
		which is larger than 0 as proven before.
\end{enumerate}

This concludes the case that $F_{i,1-j}$ is $t^{\rightarrow}$-open.
If it is not $t^{\rightarrow}$-open, then it has to be closed or $t^{\rightarrow}$-halfopen by \Pref{ESC1}.
Assume that it is closed.
If $1-j=\bit_{i+1}$, then $\valustar_{\sigmae}^\M(F_{i,1-j})=\rew{s_{i,1-j},h_{i,1-j}}+\valustar_{\sigmae}^\M(b_{i+1})$ by Properties (\ref{property: USV1})$_i$ and (\ref{property: EV1})$_{i+1}$.
Since $\rew{s_{i,1-j},h_{i,1-j}}>\sum_{\ell\in[i]}W_{\ell}^\M$, this implies $\valustar_{\sigmae}^\M(F_{i,1-j})>\valustar_{\sigmae}(F_{i,j})=\valustar_{\sigmae}^\M(t^{\rightarrow}).$
If $j=\bit_{i+1}$, then the same properties imply $\valustar_{\sigmae}(F_{i,1-j})=\rew{s_{i,1-j}}+\valustar_{\sigmae}^\M(b_1)$.
Since $\valustar_{\sigmae}^\M(b_1)=\valustar_{\sigmae}^\M(t^{\rightarrow})=\valustar_{\sigmae}^\M(F_{i,j})$, this implies the statement.

Hence let $F_{i,1-j}$ be $t^{\rightarrow}$-halfopen.
Then $\valu_{\sigmae}^\M(F_{i,1-j})=\frac{1-\e}{1+\e}\valu_{\sigma}^\M(t^{\rightarrow})+\frac{2\e}{1+\e}\valu_{\sigmae}^\M(s_{i,1-j}).$
We prove that $\valu_{\sigmae}^\M(s_{i,1-j})>\valu_{\sigmae}^\M(t^{\rightarrow})$ in this case.
If $1-j\neq \bit_{i+1}$, then this follows from \Pref{USV1}$_i$ as $\valu_{\sigmae}^\M(s_{i,1-j})=\rew{s_{i,1-j}}+\valu_{\sigmae}^\M(t^{\rightarrow})$ in that case.
If $1-j=\bit_{i+1}$, then it holds that $\valu_{\sigmae}^\M(s_{i,1-j})=\rew{s_{i,1-h},h_{i,1-j}}+\valu_{\sigmae}^\M(b_{i+1})$ by  Properties (\ref{property: USV1})$_i$ and (\ref{property: EV1})$_{i+1}$.
The statement then follows since $\rew{h_{i,1-j}}>\sum_{\ell\in[i]}W_{\ell}^\M$.

Consequently, $\valu_{\sigmae}^{\M}(s_{i,1-j})>\valu_{\sigmae}^{\M}(t^{\rightarrow})$. 
This implies $\valu_{\sigmae}^\M(F_{i,1-j})>(1-\e)\valu_{\sigmae}^\M(t^{\rightarrow})+\e\valu_{\sigmae}^\M(s_{i,1-j})$ which yields $\valu_{\sigmae}^\M(F_{i,1-j})>\valu_{\sigmae}^\M(F_{i,j})$ as proven earlier.
\end{proof}

\ClosingWhenSelectorIsPointing*

\begin{proof}
Since $\sigma(g_i)=F_{i,j}$ and $j\neq\bit_{i+1}=\indbit_{i+1}^{\sigma}$, we have $\sigmae(b_i)=\sigma(b_i)=b_{i+1}$ by \Pref{EV2}$_{i}$.
This implies $i\geq\relbit{\sigma}$ by \Pref{REL1}. 
As $\relbit{\sigma}=\relbit{\sigmae}$, this implies that $\sigmae$ has \Pref{CC1}$_{i'}$ for all indices $i'$.
This further implies that $\sigmae$ has \Pref{ESC1},(\ref{property: EV1})$_{i'}$ and (\ref{property: USV1})$_{i'}$ for all $i'\in[n]$.
Furthermore, since $\sigma$ has all other properties defining a phase-$1$-strategy, $\sigmae$ has them as well.
As we do not perform changes within the cycle center $F_{i,\indbit_{i+1}^{\sigma}}$, also $\indbit^{\sigma}=\indbit^{\sigmae}\eqqcolon\indbit$.
Since $\sigma$ has \Pref{SVG}$_i$ and since the cycle center $F_{i,j}$ is not closed for $\sigma$ by the choice of $e$, we have $j=0$.
This implies that $\sigmae$ has \Pref{SVG}$_i$ as well.
Hence $\sigmae$ is a phase-$1$-strategy for $\bit$.

Proving that $\sigmae$ is well-behaved follows by re-evaluating Properties (\ref{property: D1}), (\ref{property: MNS4}), (\ref{property: MNS6}), (\ref{property: EG2}), (\ref{property: DN1}) and (\ref{property: DN2}).
This set of properties is sufficient as we do not need to verify properties where the conclusion might become true or the premise might become false since the implication is then already true.

\begin{enumerate}[align=right, leftmargin=1.75cm]
	\item[(\ref{property: D1})] By the premise of this property, $\sigmae(b_i)=g_i$, contradicting $\sigmae(b_i)=b_{i+1}$.
		\Pref{D2} holds by the same argument since $i<\relbit{\sigmae}$ implies $\sigmae(b_i)=g_i$.
	\item[(\ref{property: MNS4})] Since $\sigma$ is well-behaved, this only needs to be reevaluated if $i=\minnegsige{s}$.
		Since $i\neq 1$ cannot occur by assumption, let $i>1$.
		Then, $1<\minnegsige{s}\leq\minnegsige{g}<\minsige{b}$.
		Thus, in particular $\sigmae(b_1)=b_{2}, \sigmae(g_1)=F_{1,1}$ and $\sigmae(s_{1,1})=h_{1,1}$.
		By \Pref{USV1}$_1$ and \Pref{EV1}$_2$, $\sigmae(b_2)=g_2$, implying $\minsige{b}=2$.
		But this contradicts the premise since $1<\minnegsige{s}<\minsige{b}$ implies $\minsige{b}\geq 3$.
	\item[(\ref{property: MNS6})] Since $\sigma$ is well-behaved, this only needs to be reevaluated if $i=\minnegsige{s}$.
		Since  $i=1$ cannot occur by assumption, let $i>1$.
		Then, $1<\minnegsige{s}\leq\minnegsige{g}<\minsige{b}$, implying the same contradiction as in the last case.
	\item[(\ref{property: EG2})] The cycle center $F_{i,j}$ is closed with respect to $\sigmae$, so the premise is incorrect.
	\item[(DN*)] Since the only cycle center in level $n$ is  $F_{n,0}$ and since we always have $\bit_{n+1}=0$ by definition, the choice of $j$ implies that we cannot have $i=n$.
\end{enumerate}


We next prove $I_{\sigmae}=\{(d_{i,j,k},F_{i,j})\colon\sigmae(d_{i,j,k})\neq F_{i,j}\}$.
The only vertices  that have an edge towards $F_{i,j}$ are $d_{i,j,*}$ and $g_i$.
Since closing $F_{i,j}$ increases its valuation, the valuation of these vertices might increase as well.
Since no player 0 vertex has an edge to either $d_{i,j,0}$ or $d_{i,j,1}$, no new improving switch involving these vertices can emerge.
However, the valuation of $g_i$ might increase due to $\sigma(g_i)=\sigmae(g_i)=F_{i,j}$.
We now prove that this increase does not create new improving switches and that all switches but $e$ that are improving for $\sigma$ stay improving for $\sigmae$.

It suffices to prove that $\sigmae(b_{i})=b_{i+1}$ and $\valu_{\sigmae}^\P(g_i)\unlhd\valu_{\sigmae}^\P(b_{i+1})$ as well as$\sigmae(s_{i-1,1})=b_1$ and $\valu_{\sigmae}^\P(h_{i-1,1})\unlhd\valu_{\sigmae}^\P(b_1)$. 
Since $\sigma(b_i)=\sigmae(b_i)=b_{i+1}$ and $i\geq\relbit{\sigmae}$, we have $\valustar_{\sigmae}^\P(b_{i+1})=L_{i+1}^\P$.
Since $\sigmaebar(d_{i,j})$ implies $\valustar_{\sigmae}^\P(F_{i,j})=\valustar_{\sigmae}^\P(s_{i,j})$ by \Cref{lemma: Exact Behavior Of Counterstrategy}, we have $\valustar_{\sigmae}^\P(g_i)=\{g_i,s_{i,j}\}\cup\valustar_{\sigmae}^\P(b_1)=\{g_i,s_{i,j}\}\cup B_1^\P$ by the choice of~$j$ and \Pref{USV1}$_i$.
Thus, $\{g_i,s_{i,j}\}\lhd\bigcup_{\ell\in[i]}W_{\ell}^\P$ and  $\{g_i,s_{i,j}\}\lhd\bigcup_{\ell\in[i]}\{W_{\ell}^\P\colon\sigmae(b_{\ell})=g_{\ell}\}$ yield $\valustar_{\sigmae}^\P(g_i)\lhd\valustar_{\sigmae}^\P(b_{i+1})$.
For the second statement, we observe that $\bit_i=0$ implies $\sigmae(s_{i-1,1})=b_1$ by \Pref{USV1}$_{i-1}$.
The second part then follows using similar calculations as before since $\valustar_{\sigmae}^\P(h_{i-1,1})=\{h_{i-1,1},g_i,s_{i,j}\}\cup\valustar_{\sigmae}^\P(b_1)$.
\end{proof}

\ClosingWhenSelectorIsNotPointing*

\begin{proof}
By similar arguments used in the proof of \Cref{lemma: Closing when selector is pointing}, Properties (\ref{property: ESC1}), (\ref{property: REL1}) and (\ref{property: USV1})$_{i'}$, (\ref{property: CC1})$_{i'}$, (\ref{property: EV1})$_{i'}$, (\ref{property: EV2})$_{i'}$ and (\ref{property: EV3})$_{i'}$ for $i'\in[n]$ are valid for $\sigmae$.
Consider \Pref{SVM}$_i$ and let $G_n=M_n$.
Since $1-j=\indbit^{\sigma}_{i+1}$, the premise of this property is incorrect, hence $\sigmae$ has \Pref{SVM}$_i$.
Consider \Pref{SVG}$_i$ and let $G_n=S_n$.
If $1-j=\indbit^{\sigma}_{i+1}=0$, then $\sigmae$ has \Pref{SVG}$_i$ as well.
Hence assume $1-j=\indbit^{\sigma}_{i+1}=1$.
Then, since $\sigma$ has \Pref{SVG}$_i$, it follows that $\sigmabar(d_{i,1})$, implying $\sigmaebar(d_{i,1})$.
Thus, $\sigmae$ has \Pref{SVG}$_i$ resp. \Pref{SVM}$_i$, implying that $\sigmae$ is a phase-1-strategy for $\bit$.

Since $\sigma(g_i)=F_{i,1-j}$, applying $e$ does not close the chosen cycle center.
It is thus not necessary to reevaluate the assumptions of \Cref{table: Well behaved properties} and thus, $\sigma$ being well-behaved implies that $\sigmae$ is well-behaved.
It hence remains to prove $I_{\sigmae}=\mathfrak{D}^{\sigmae}\cup\{(g_i,F_{i,j})\}$.

By \Pref{EV1}, $\sigma(b_i)=\sigmae(b_i)=b_{i+1}$ implies that $F_{i,1-j}$ is not closed with respect to both $\sigma$ and $\sigmae$.
Hence, $\valustar_{\sigmae}^*(F_{i,1-j})=\valustar_{\sigmae}^*(b_1)$ by \Cref{lemma: Cycle centers in Phase One}.
Since \Pref{USV1}$_{i}$ and the choice of $j$ imply $\sigmae(s_{i,j})=\sigma(s_{i,j})=b_1$, it holds that \[\valustar_{\sigmae}^*(F_{i,1-j})=\valustar_{\sigmae}(b_1)^*\lhd\ubracket{s_{i,j}}\oplus\valustar_{\sigmae}^*(b_1)=\valustar_{\sigmae}^*(F_{i,j}),\] implying $(g_i,F_{i,j})\in I_{\sigmae}$.

Since $\sigma(g_i)=F_{i,1-j}$, the only vertices whose valuations change by applying $e$ are $d_{i,j,0},d_{i,j,1}$ and $F_{i,j}$.
This implies that no new switches besides the switch $(g_i,F_{i,j})$ are created and that all improving switches for $\sigma$ but $e$ stay improving for $\sigmae$.
\end{proof}

\SwitchingSelectorInPhaseOne*

\begin{proof}
	Let $j\coloneqq 1-\bit_{i+1}=1-\indbit^{\sigma}_{i+1}$.
	Since $\sigma$ is a phase-1-strategy for $\bit$, $\bit_i=0$ implies $\sigma(b_i)=b_{i+1}$.
	Since no cycle center is closed when applying $e$, $\sigmae$ has Properties (\ref{property: ESC1}), (\ref{property: EV1})$_{i'}$, (\ref{property: EV3})$_{i'}$, (\ref{property: CC1})$_{i'}$ and (\ref{property: USV1})$_{i'}$ for all $i'\in[n]$. 
	Also, since $\sigmae(b_i)=b_{i+1}$, the premise of \Pref{EV2}$_i$ is incorrect with respect to $\sigmae$, hence it has the property for all indices.
	Since $\sigma$ has \Pref{REL1} and $\sigmae(b_i)=b_{i+1}$, it also has \Pref{REL1}.
	This implies $\incorrect{\sigmae}=\emptyset$, and thus $i\geq\relbit{\sigmae}=\relbit{\sigma}$.	
	Also, since $\sigmabar(d_{i,j})$ by assumption, $\sigmae$ has \Pref{SVG}$_i$ resp. \Pref{SVM}$_i$.

	By the choice of $e$, by $i\geq\relbit{\sigmae}$ and since $\sigma$ is well-behaved, it suffices to investigate Properties (\ref{property: B3}), (\ref{property: MNS4}), (\ref{property: MNS6}) and (\ref{property: EG4}) in order to prove that $\sigmae$ is well-behaved.
	\begin{enumerate}[align=right,leftmargin=1.75cm]
		\item[(\ref{property: B3})] Since $\sigma$ has \Pref{USV1}$_i$, the premise of this property is incorrect.
		\item[(\ref{property: MNS4})] Let the premise be correct, i.e., let $\relbit{\sigmae}=1\wedge\minnegsige{s}\leq\minnegsige{g}<\minsige{b}$.
			Let $i'\coloneqq\minnegsige{s}$.
			For the sake of a contradiction, let $i'=1$.
			Then $\sigmae(b_1)=b_2$ since $\relbit{\sigmae}=1$.
			If $\sigmae(g_1)=F_{1,1}$, then $\sigmae(s_{1,1})=b_1$.
			Thus, $\indbit_{2}=\bit_2=0$ by \Pref{USV1}$_1$.
			If also $\sigma(g_1)=F_{1,1}$, \Pref{SVG}$_1$ resp. \Pref{SVM}$_1$ would imply $\sigmabar(d_{1,1})$, contradicting \Pref{MNS4} for $\sigma$.
			If $\sigma(g_1)=F_{1,0}$, then $\minnegsig{g}=1$.
			If also $\minnegsige{g}=1$, then the statement follows by applying \Pref{MNS4} to $\sigma$.
			Otherwise, we need to have $e=(g_1,F_{1,1})$, contradicting the assumption $i\neq 1$.
			Hence consider the case $i'=\minnegsige{s}>1$.
			Then $1<\minnegsige{s}\leq\minnegsige{g}<\minsige{b}$, implying $\sigmae(g_1)=F_{1,1}, \sigmae(s_{1,1})=h_{1,1}$ and $\minsige{b}\geq 3$.
			By \Pref{USV1}$_{1}$, this  implies $\indbit_2=1$, hence $\sigmae(b_2)=g_2$ by \Pref{EV1}$_2$.
			But then $\minsige{b}=2$ which is a contradiction.
			Therefore the premise cannot be correct, hence the implication is correct.
		\item[(\ref{property: MNS6})] Assume the premise is correct, i.e., assume $\relbit{\sigmae}=1\wedge\minnegsige{s}<\minsige{b}\leq\minnegsige{g}$.
			Let $i'\coloneqq\minnegsige{s}$ and assume $i'=1$.
			Then $\sigmae(g_{1})=F_{1,1}, \sigmae(s_{1,1})=b_1$ and $\sigmae(b_1)=b_2$.
			If also $\sigma(g_1)=F_{1,1}$, then \Pref{SVG}$_1$ resp. \Pref{SVM}$_1$ would imply $\sigmabar(d_{1,1})$ as in the last case, contradicting \Pref{MNS6} for $\sigma$.
			However, since $\sigma(g_1)=F_{1,0}$ implies $e=(g_1,F_{1,1})$, this again contradicts the assumption $i\neq 1$.
			Hence consider the case $i'=\minnegsige{s}>1$.
			Then $1<\minnegsige{s}<\minsige{b}\leq\minnegsige{g}$ which implies the same contradiction that occurred when discussing \Pref{MNS4}.
		\item[(\ref{property: EG4})] By \Pref{ESC1}, the premise of this property is always incorrect, hence the implication is correct.
	\end{enumerate}
	Note that the other Properties (MNS*) do not need to be considered since their conclusion is correct for level $i$ by assumption.
	In addition, none of the properties (EBG*) needs to be checked due to \Pref{ESC1}.

	It remains to prove $I_{\sigmae}=I_{\sigma}\setminus\{e\}$.
	This follows by proving $\sigmae(b_{i})=b_{i+1}$ and $\valu_{\sigmae}^*(g_i)\prec\valu_{\sigmae}^*(b_{i+1})$ as well as $\sigmae(s_{i-1,1})=b_1$ and $\valu_{\sigmae}^*(h_{i-1,1})\prec\valu_{\sigmae}^*(b_1)$.
	This can be proven in the same way as it was proven in the proof of \Cref{lemma: Closing when selector is pointing}.
\end{proof}

\PhaseOneLowOR*

\begin{proof}
By transitivity, $\sigma\in\reach{\sigma_0}$.
Since $e=(d_{i,j,k},F_{i,j})\in I_{\sigma}$, the cycle center $F_{i,j}$ cannot be closed.
Hence, since $\sigma$ is a phase-1-strategy for $\bit$, either $\bit_i=0$ or $\bit_{i+1}\neq j$.
Consequently, exactly one of the following cases is true:
\begin{enumerate}
	\item $\sigma(d_{i,j,1-k})\neq F_{i,j}$
	\item $\sigma(d_{i,j,1-k})=F_{i,j}\wedge j\neq\bit_{i+1}\wedge\sigma(g_i)=F_{i,j}$ 
	\item $\sigma(d_{i,j,1-k})=F_{i,j}\wedge j\neq\bit_{i+1}\wedge\sigma(g_i)=F_{i,1-j}\wedge\sigma(b_i)\neq g_i$
	\item $\sigma(d_{i,j,1-k})=F_{i,j}\wedge j\neq\bit_{i+1}\wedge\sigma(g_i)=F_{i,1-j}\wedge\sigma(b_i)= g_i$
	\item $\sigma(d_{i,j,1-k})=F_{i,j}\wedge j=\bit_{i+1}$
\end{enumerate}

We prove that case four and five cannot occur.
Assume that the conditions of the fourth case were true.
Then, since $\sigma$ is a phase-1-strategy for $\bit$, \Pref{EV1}$_i$ and \Pref{EV2}$_i$, imply $\bit_{i}=1$ and $\bit_{i+1}=1-j$.
This also implies that $t_{\bit}=0$ is the only feasible parameter for $(d_{i,j,k},F_{i,j})$ and $(d_{i,j,1-k},F_{i,j})$.
Now, by assumption,  $\sigma(d_{i,j,1-k})=F_{i,j}$.
If $\canstrat(d_{i,j,1-k})=F_{i,j}$,  then \Pref{OR1}$_{i,j,k}$ and \Pref{OR2}$_{i,j,k}$ imply $\occrec^{\canstrat}(d_{i,j,1-k},F_{i,j})=\ell^{\bit}(i,j,1-k)+1$, contradicting that $t_{\bit}=0$ is the only feasible parameter.
Thus, assume $\canstrat(d_{i,j,1-k})\neq F_{i,j}$.
Then $(d_{i,j,1-k},F_{i,j})$ was applied during $\canstrat\to\sigma$ and in particular before $(d_{i,j,k},F_{i,j})$.
Thus, $\canstrat(d_{i,j,1-k})\neq F_{i,j}$ and \Pref{OR4}$_{i,j,1-k}$ implies $\occrec^{\canstrat}(d_{i,j,1-k},F_{i,j})=\floor{\frac{\bit+1}{2}}-1$.
Since $\bit_i=1\wedge\bit_{i+1}=1-j$, \Cref{lemma: Numerics Of Ell} implies \[\ell^{\bit}(i,j,1-k)=\ceil{\frac{\bit+\sum(\bit,i)+1-k}{2}}\geq\floor{\frac{\bit+1}{2}}.\]
As $t_{\bit}=0$ is the only feasible parameter, it thus needs to hold that \[\occrec^{\canstrat}(d_{i,j,1-k},F_{i,j})=\floor{\frac{\bit+1}{2}}-1=\floor{\frac{\bit+1-(1-k)}{2}}.\]
This implies $k=0$ and that $\bit$ is odd.
But then $\occrec^{\sigma}(d_{i,j,k},F_{i,j})=\floor{(\bit+1)/2}$, contradicting the assumptions.

Consider the fifth case.
Then, $\bit_i=0$ as $j=\bit_{i+1}$.
If $\canstrat(d_{i,j,1-k})=\sigma(d_{i,j,1-k})=F_{i,j}$, then $i\neq\nsb$ by the definition of a canonical strategy.
If $\canstrat(d_{i,j,1-k})\neq F_{i,j}$, then the switch $(d_{i,j,1-k},F_{i,j})$ was applied during $\canstrat\to\sigma$.
This implies $\occrec^{\sigma}(d_{i,j,1-k},F_{i,j})<\floor{(\bit+1)/2}$.
By \Cref{lemma: Occurrence Records Cycle Vertices}, there can be at most one improving switch in level $\nsb$ with an occurrence record strictly smaller than $\floor{(\bit+1)/2}$.
This implies $i\neq\nsb$.
Since $i=1$ would imply $i=\nsb$ due to $\bit_i=0$, we thus have $i\geq 2$.

Assume $\canstrat(d_{i,j,1-k})=F_{i,j}$.
Then $\occrec^{\sigma}(d_{i,j,1-k},F_{i,j})=\ell^{\bit}(i,j,1-k)+1<\floor{(\bit+1)/2}$ by \Pref{OR2}$_{i,j,1-k}$ and \Pref{OR1}$_{i,j,-1-k}$, hence $\ell^{\bit}(i,j,1-k)<\floor{(\bit+1)/2}-1$.
However, since $\bit_i=0$ and $\bit_{i+1}=j$, \Cref{lemma: Numerics Of Ell} implies that either $\ell^{\bit}(i,j,1-k)\geq\bit$ or \[\ell^{\bit}(i,j,1-k)=\ceil{\frac{\bit+2^{i-1}+\sum(\bit,i)+1-(1-k)}{2}}\geq\floor{\frac{\bit+3}{2}}=\floor{\frac{\bit+1}{2}}+1\]which is a contradiction in either case.
The case $\canstrat(d_{i,j,1-k})\neq F_{i,j}$ yields the same contradiction developed for case four.

Thus one of the first three listed cases needs to be true.
In the first resp. third case, we can apply \Cref{lemma: Not closing cycle center in phase one} resp. \ref{lemma: Closing when selector is not pointing} to prove the part of the statement regarding the improving switches.
In order to apply \Cref{lemma: Closing when selector is pointing}, we need to prove that the conditions of the second case can only occur if $G_n=S_n$ and $i\neq1$.

Thus assume that the conditions of the second case were true.
Assume $i=1$.
Then, since $\sigma(g_1)=F_{1,1-\bit_2}$, wit holds that $\bit_1=0$ by \Pref{EV1}$_1$ and \Pref{EV2}$_1$.
By the choice of $j$ and \Cref{lemma: Numerics Of Ell}, this implies $\ell^{\bit}(i,j,k)=\floor{(\bit+1-k)/2}$.
By \Pref{OR3}$_{i,j,k}$ and \Pref{OR4}$_{i,j,k}$, it thus needs to hold that $\occrec^{\canstrat}(d_{i,j,k},F_{i,j})=\ell^{\bit}(i,j,k)=\floor{(\bit+1)/2}-1$.
This can only happen if $k=1$ and if $\bit$ is odd, contradicting $\bit_1=0$.
Consequently, $i\neq 1$.
Proving that the conditions can only occur if $G_n=S_n$ can be done by proving that we have $(g_i,F_{i,1-j})\in I_{\sigma}$ if $G_n=M_n$, contradicting $I_{\sigma}=\mathfrak{D}^{\sigma}$.
As proving this is rather tedious, we omit this part here.

It remains to show that there is a feasible parameter $t_{\bit+1}$ for $\bit+1$ such that \[\occrec^{\sigmae}(d_{i,j,k},F_{i,j})=\min\left (\left \lfloor\frac{(\bit+1)+1-k}{2}\right \rfloor,\ell^{\bit+1}(i,j,k)+t_{\bit+1}\right ).\]
Since $\occrec^{\sigma}(e)=\occrec^{\canstrat}(e)$, we have $\occrec^{\sigmae}(e)=\occrec^{\canstrat}(e)+1$.
Also, there is a parameter $t_{\bit}$ feasible for $\bit$ such that $\occrec^{\sigma}(e)=\min(\floor{(\bit+1-k)/2},\ell^{\bit}(i,j,k)+t_{\bit})=\floor{(\bit+1)/2}-1$  by the choice of $e$.
Consequently, $\occrec^{\sigmae}(e)=\floor{(\bit+1)/2}$.
We distinguish two cases.
\begin{enumerate}
	\item $i=\nsb\wedge j=\bit_{\nsb+1}$.
		Since we have one of the first three cases discussed earlier, this implies $\sigma(d_{i,j,1-k})\neq F_{i,j}$.
		Moreover, $\bit$ needs to be odd since both cycle edges of $F_{\nsb,\bit_{\nsb+1}}$ have an occurrence record of $\floor{(\bit+1)/2}$ if $\bit$ is even.
		Consequently, $\nsb>1$.
		Thus, by the choice of $e$ and \Cref{lemma: Occurrence Records Cycle Vertices}, it holds that $k=1$.
		It therefore suffices to show $\floor{(\bit+1)/2}=\ceil{\lastflip{\bit+1}{i}{\{(i+1,j)\}}/2}.$
		This however follows immediately from the choice of $i$ and $j$ and the fact that $\bit$ is odd.
	\item $i\neq\nsb\vee j\neq\bit_{i+1}$.
		This implies $\bit_i=0\vee j\neq\bit_{i+1}$, hence $(\bit+1)_{i}=0\vee(\bit+1)_{i+1}\neq j$.
		We thus need to show that there is a parameter $t_{\bit+1}$ feasible for $\bit+1$ such that \[\floor{\frac{\bit+1}{2}}=\min\left(\floor{\frac{(\bit+1)+1-k}{2}}, \ell^{\bit+1}(i,j,k)+t_{\bit+1} \right).\]
		By \Cref{lemma: Progress for unimportant CC}, $\ell^{\bit}(i,j,k)+1=\ell^{\bit+1}(i,j,k)$.
		We distinguish the following cases.
		\begin{enumerate}
			\item Let $\floor{(\bit+1)/2}-1=\floor{(\bit+1-k)/2}$.
				This implies $k=1$ and $\bit\bmod2=1$.
				Consequently, $\occrec^{\sigmae}(e)=\floor{(\bit+1)/2}=\floor{(\bit+1+1-k)/2}$.
				It remains to define a feasible parameter $t_{\bit+1}$.
				Since $\occrec^{\canstrat}(e)=\floor{(\bit+1-k)/2}$, there is a feasible  $t_\bit$ for $\bit$ such that $\floor{(\bit+1-k)/2}\leq\ell^{\bit}(i,j,k)+t_{\bit}$.
				Since $\canstrat(d_{i,j,k})\neq F_{i,j}$ due to $e\in I_{\canstrat}$, \Pref{OR2}$_{i,j,k}$ implies $t_{\bit}\neq 1$.
				Hence we can choose $t_{\bit+1}=0$ as $\occrec^{\sigmae}(e)=\floor{(\bit+1-k)/2}+1\leq\ell^{\bit}(i,j,k)+1=\ell^{\bit+1}(i,j,k).$
			\item Let $\floor{(\bit+1)/2}-1=\ell^{\bit}(i,j,k)+t_{\bit}$ for some parameter $t_{\bit}$ feasible for $\bit$ but $\floor{(\bit+1)/2}-1\neq\floor{(\bit+1-k)/2}$.
				Then, \Pref{OR2}$_{i,j,k}$ implies $t_{\bit}\neq 1$.
				Consider the case $t_{\bit}=0$ first.
				Then $\occrec^{\sigmae}(e)=\ell^{\bit}(i,j,k)+1=\ell^{\bit+1}(i,j,k)$ and $\occrec^{\sigma}(e)=\floor{(\bit+1)/2}\leq\floor{(\bit+1+1-k)/2}$.
				Thus, choosing $t_{\bit+1}=0$ is a feasible choice giving the correct characterization.
				Thus consider the case $t_{\bit}=-1$.
				Then, by \Pref{OR3}$_{i,j,k}$, $\bit$ is odd and $k=0$.
				This then implies that $\occrec^{\sigmae}(e)=\ell^{\bit}(i,j,k)=\floor{(\bit+1)/2}=\floor{(\bit+1+1-k)/2}$ as well as $\occrec^{\sigmae}(e)=\ell^{\bit+1}(i,j,k)-1$.
				We thus choose $t_{\bit+1}=0$ which is a feasible choice, does not contradict \Pref{OR3} for $\bit+1$ and yields the desired characterization. \qedhere
		\end{enumerate}
\end{enumerate}
\end{proof}

\ClosingActiveCC*

\begin{proof}
We have $\nsb=\relbit{\sigma}$ as $\sigma$ has \Pref{REL1} and \Pref{EV1}$_{i'}$ for all $i'<\nsb$.
Also, $\relbit{\sigma}=\relbit{\sigmae}$ by the choice of $e$.
Since we do not close any cycle centers in any level below $\relbit{\sigmae}$, $\sigmae$ has \Pref{CC1}$_{i'}$ for all $i'\in[n]$.
\begin{enumerate}
	\item Since the cycle centers of levels $i>\nsb$ are not changed, $\indbit_{i}^{\sigmae}=\indbit_{i}^{\sigma}=\bit_i=(\bit+1)_i$ for all $i>\nsb$.
		Moreover, $\indbit_{i}^{\sigmae}=\sigmaebar(d_{\nsb,j})=1=(\bit+1)_{\nsb}$ by the definition of $\nsb$  and the choice of $e$.
		It remains to show $\indbit_{i}^{\sigmae}=0$ for all $i<\nsb$.
		This is proven by backwards induction.
		Hence let $i=\nsb-1$ and consider $\indbit_{i}^{\sigmae}=\sigmaebar(d_{i,\indbit_{i+1}^{\sigmae}})$.
		
		Since $\indbit_{i+1}^{\sigmae}=1$, we prove $\sigmaebar(d_{i-1,1})=0$.
		We have $\indbit_{\nsb-1}^{\sigma}=1$ and $\indbit_{\nsb}^{\sigma}=0$.
		Thus $\sigma(b_{\nsb-1})=g_{\nsb-1}$ by \Pref{EV1}$_{\nsb-1}$, so $0=\sigmabar(d_{\nsb-1,1-\indbit_{\nsb}^{\sigma}})=\sigmabar(d_{\nsb-1,1})=\sigmaebar(d_{\nsb-1,1})$ by \Pref{EV3}$_{\nsb-1}$. 
		
		Now consider some $i<\nsb-1$.
		By the induction hypotheses, $\indbit_{i+1}^{\sigmae}=0$.
		We hence prove $\sigmaebar(d_{i,\indbit_{i+1}^{\sigmae}})=\sigmaebar(d_{i,0})=0$.
		By the definition of~$\nsb$, $\indbit_{i}^{\sigma}=\indbit_{i+1}^{\sigma}=1$. 
		Hence,  $\sigma(b_{i})=g_{i}$ by \Pref{EV1}$_{i}$, implying $\sigmabar(d_{i,0})=0$ by \Pref{EV3}$_{i}$.
		
	\item We prove that $\sigmae$ has the listed properties.
		Since $\indbit_{i}^{\sigmae}=\indbit_{i}^{\sigma}$ for all $i>\nsb$, $\sigmae$ has \Pref{EV1}$_{i}$ for all $i\geq\nsb$.
		This also implies that it has \Pref{EV2}$_{i}$ and (\ref{property: USV1})$_{i}$ for all $i\geq\nsb$.
		In addition, it has \Pref{EV3}$_{i}$ for all $i\geq\nsb$ and thus in particular for all $i>\nsb$.
		As \Pref{REL1} does not consider cycle centers, it remains valid for $\sigmae$.
		Since $\indbit_1^{\sigma}=1$ if and only if $\indbit_{1}^{\sigmae}=0$ and since $\sigma$ has \Pref{ESC1}, $\sigmae$ has \Pref{ESC2} if $\nsb=0$.
		Thus $\sigmae$ has all properties for the bound $\relbit{\sigmae}$ if $\nsb>1$ resp. for the bound $1$ if $\nsb=1$ as specified in \Cref{table: Definition of Phases}. 
		
	\item Since $\sigma$ is well-behaved, it suffices to reevaluate Properties (\ref{property: MNS4}), (\ref{property: MNS6}), (\ref{property: DN1}) and (\ref{property: DN2}).
		\begin{enumerate}[align=right, leftmargin=1.75cm]
			\item[(\ref{property: MNS4})] By the choice of $e$, the premise of this property is true for $\sigmae$ if and only if it is true for $\sigma$.
				In particular, $\minnegsige{s}=\minnegsig{s}, \minnegsige{g}=\minnegsig{g}$ and $\minsig{b}=\minsige{b}$.
				In addition, $\relbit{\sigma}=\relbit{\sigmae}=1$ implies that we close the cycle center $F_{1,\bit_2}$.
				If $\minnegsig{s}\neq 1$, then the conclusion is correct for $\sigmae$ if and only if it is correct for~$\sigmae$, hence $\sigmae$ has \Pref{MNS4}.
				It thus suffices to consider the case $\minnegsige{s}=1$.
				Assume the conditions of the premise were fulfilled and let $j'\coloneqq\sigmabar(g_1)$.
				Then, by assumption, $\sigmaebar(s_{\minnegsige{s}})=\sigmaebar(s_{1,j'})=b_1$.
				Thus, by the choice of $j$, it follows that we do not close the cycle center $F_{1,j'}$.
				Hence, since $\sigmabar(eb_{\minnegsig{s}})\wedge\nsigmabar(eg_{\minnegsig{s}})$ by \Pref{MNS4}, also $\sigmaebar(eb_{\minnegsige{g}})\wedge\nsigmabar(eg_{\minnegsige{s}})$.
			\item[(\ref{property: MNS6})] This follows by the same arguments used for \Pref{MNS4}.
			\item[(\ref{property: DN1})] Since $i=n$ in this case, $\sigmae(b_1)=\sigma(b_1)=g_1$ by the definition of $\nsb$.
			\item[(\ref{property: DN2})] This statement only needs to be considered if $\nsigmabar(d_n)\wedge\sigmaebar(d_n)$, hence, only if $\nsb=n$.
				Then, $\indbit_{1}^{\sigma}=\dots=\indbit_{n-1}^{\sigma}=1$.	
				But then \Pref{EV1}$_{i}$ implies $\sigmae(b_{i})=g_{i}$ for all $i\leq n-1$.
		\end{enumerate}
		Since $\sigma$ is well-behaved, $\sigmae$ is thus well-behaved.

		\item We prove that $\sigma(g_\nsb)\neq F_{\nsb,j}$ and $\nsb=1$ imply $I_{\sigmae}=\mathfrak{D}^{\sigmae}\cup\{(g_{\nsb},F_{\nsb,j})\}.$

			We first prove $(g_\nsb,F_{\nsb,j})\in I_{\sigmae}$.
			Since $\nsb=\relbit{\sigmae}=1$ and by \Pref{ESC1} and \Pref{USV1}$_i$, either $\valustar_{\sigmae}^*(F_{\nsb,1-j})=\valustar_{\sigmae}^*(s_{\nsb,1-j})=\ubracket{s_{\nsb,1-j}}\oplus\valustar_{\sigmae}^*(b_1)$ or $\valustar_{\sigmae}^*(F_{\nsb,1-j})=\valustar_{\sigmae}^*(b_2)$.
			By \Pref{USV1}$_i$ and \Pref{EV1}$_{\nsb+1}$, it also holds that $\valustar_{\sigmae}^*(F_{\nsb,j})=\ubracket{s_{\nsb,j},h_{\nsb,j}}\oplus\valustar_{\sigmae}^*(b_{\nsb+1}).$
			The statement thus follows in either case since $\ubracket{h_{\nsb,j}}\succ\ubracket{s_{\nsb,1-j}}\oplus L_{1,\nsb}^*\succ L_{1,\nsb}^*$ and $\valustar_{\sigma}^*(b_1)=\valustar_{\sigma}^*(b_2)=L_{1,\nsb}^*\oplus L_{\nsb+1}^*$ as well as $\valustar_{\sigma}^*(b_{\nsb+1})=L_{\nsb+1}^*$, implying that the valuation of $F_{\nsb,j}$ is higher than the valuation of $F_{\nsb,1-j}$.

			Since $\sigma(g_\nsb)=F_{\nsb,1-j}$, the valuation of $g_\nsb$ does not change.
			Hence, only the valuations of the cycle vertices $d_{\nsb,j,0},d_{\nsb,j,1}$ can change.
			Since $F_{\nsb,j}$ is the only vertex with an edge to these vertices, the valuations of all other vertices remain the same.
			Thus,  all switches improving with respect to $\sigma$ but $e$ stay improving with respect to $\sigmae$ and no further improving switches are created. 
						
			Next, let $\sigma(g_\nsb)=F_{\nsb,j}$ and $\nsb=1$.
			We prove $I_{\sigmae}=\mathfrak{D}^{\sigmae}\cup\{(b_1,g_1)\}\cup\{(e_{*,*,*},g_1)\}$.
			We first prove $(b_1,g_1)\in I_{\sigmae}$.
			
			By \Pref{EV1}$_1$, $\sigmae(b_1)=b_2$, and it suffices to show $\valustar_{\sigmae}^*(g_1)\succ\valustar_{\sigmae}^*(b_2)$.
			Since $\relbit{\sigmae}=1$, we have $\valustar_{\sigmae}^*(b_2)=L_2^*$.
			Let $G_n=S_n$.
			We use \Cref{corollary: Complete Valuation Of Selection Vertices PG} to determine the valuation of $g_1$.
			We hence need to analyze $\lambda_1^\P$.
			If $\sigmae(b_2)=g_2$, then $\lambda_1^\P=1$.
			If $\sigmae(b_2)=b_3$, then $j=\bit_2=0$ by \Pref{EV1}$_{2}$ and thus $\sigmae(g_{j})=F_{j,0}$ by assumption.
			Thus, $\lambda_1^\P=1$ in either case.
			Consider the different cases listed in \Cref{corollary: Complete Valuation Of Selection Vertices PG}.
			Since $\sigmae(b_1)=b_2$, the first case cannot occur.
			In addition, since $\sigmae(g_1)=F_{1,j}$ and the cycle center $F_{1,j}$ is closed, the cases 2 to 5 cannot occur.
			Hence consider the sixth case.
			As before, $\sigmaebar(g_1)=j=\bit_{2}$ by assumption, implying $\sigmae(s_{1,\sigmaebar(g_1)})=\sigmae(s_{1,j})=h_{1,j}$ by \Pref{USV1}$_1$.
			Thus, the sixth case cannot occur.
			As a consequence, by applying either the seventh or eighth case of \Cref{corollary: Complete Valuation Of Selection Vertices PG}, \Pref{USV1}$_1$ implies $\valustar_{\sigmae}^\P(g_1)=W_1^\P\cup\valustar_{\sigmae}^\P(b_2)\rhd\valustar_{\sigmae}^\P(b_2)$ since $j=\sigmaebar(b_2)$.		
			This also implies that any edge $(e_{*,*,*},g_1)$ is an improving switch as claimed.
			
			Now consider the case $G_n=M_n$.
			We use \Cref{corollary: Complete Valuation Of Selection Vertices MDP} to evaluate $\valustar_{\sigma}^\M(g_1)$ and thus determine $\lambda_1^\M$.
			If $\sigmae(b_2)=b_3$, then $\lambda_1^\M=1$ by the same arguments used when analyzing $\lambda_1^\P$.
			Since $\sigmaebar(d_1)\wedge\sigmaebar(s_1)$ in this case, the conditions of the last case of \Cref{corollary: Complete Valuation Of Selection Vertices MDP} are fulfilled.
			Consequently, $\sigmae(b_2)=b_3$ implies that $\valustar_{\sigma}^\M(g_1)=W_1+\valustar_{\sigma}^\M(b_2)>\valustar_{\sigma}^\M(b_2)$.
			If $\sigmae(b_2)=g_2$, we have $\lambda_1^\M=2$.
			However, by \Cref{corollary: Complete Valuation Of Selection Vertices MDP}, case~1, $\valustar_{\sigma}^\M(g_1)=W_1+\valustar_{\sigma}^\M(b_2)>\valustar_{\sigma}^\M(b_2)$ holds also in this case.
			This again implies that any edge $(e_{*,*,*},g_1)$ is improving for $\sigmae$.
			
			We now show that no further improving switches are created and that existing improving switches remain improving.
			The only vertices having edges towards $g_1$ are the vertices $b_1$ and $e_{*,*,*}$.
			It thus suffices to show that the valuations of these vertices does not change.
			This however follows from $\sigmae(b_1)=b_2$ and $\sigmae(e_{i,j,k})\neq g_1$.
			
			It remains to show that $\sigmae$ is a phase-3-strategy for $\bit$ in either case.
			By the first two statements, it suffices to show that $\sigmae$ has \Pref{USV2}$_{i,\bit_{i+1}}$ for all $i<\nsb$.
			But, since $\nsb=1$, there is no such $i$.
			Also, by the definition of a pseudo phase-$3$-strategy, it directly follows that $\sigmae$ is a such a strategy if $\sigma(g_{\nsb})\neq F_{\nsb,b_{\nsb+1}}$.
		
		\item Since $\sigma$ is a phase-1-strategy for $\bit$, it follows that$\sigma(s_{i,\indbit^{\sigma}_{i+1}})=h_{i,\indbit^{\sigma}_{i+1}}$ and $\sigma(s_{i,1-\indbit^{\sigma}_{i+1}})=b_1$ by \Pref{USV1}$_i$ for all $i<\nsb$.
			As $\bit_i=\indbit^{\sigma}_i=1-\indbit^{\sigmae}_{i}=1-(\bit+1)_{i+1}$ for all $i\leq\nsb$, this implies that $\sigmae$ has \Pref{USV3}$_{i}$ for all $i<\nsb$.
			
			We prove that $\sigma(g_{\nsb})\neq F_{\nsb,j}$ and $\nsb>1$ imply $I_{\sigmae}=\mathfrak{D}^{\sigmae}\cup\{(g_{\nsb},F_{\nsb,j})\}$.
			We observe that either $\valustar_{\sigmae}^*(F_{\nsb,1-j})=\valustar_{\sigmae}^*(s_{\nsb,1-j})$ or $\valustar_{\sigmae}^*(F_{\nsb,1-j})=\valustar_{\sigmae}^*(g_1)$.
			In addition, $\bigoplus_{\ell\in[\nsb-1]}W_{\ell^*}\prec\ubracket{s_{\nsb,1-j}}\oplus\bigoplus_{\ell\in[\nsb-1]}W_{\ell}^*\prec\ubracket{h_{\nsb,j}}$.			
			The statement can thus be shown by the same arguments used in the case $\nsb>1$.
			 
			Let $\sigma(g_\nsb)=F_{\nsb,j}$ and $\nsb>1$. 
			We prove $I_{\sigmae}=\mathfrak{D}^{\sigmae}\cup\{(b_\nsb,g_\nsb)\}\cup\{(s_{\nsb-1,1},h_{\nsb-1,1})\}.$
			We first show that $(s_{\nsb-1,1},h_{\nsb-1,1})$ is improving for $\sigmae$.
			Since $\sigmae(s_{\nsb-1,1})=b_1$ by \Pref{USV1}$_{\nsb-1}$, we prove $\valustar_{\sigmae}^*(h_{\nsb-1,1})\succ\valustar_{\sigmae}^*(b_1)$.

			It holds that $\valustar_{\sigmae}^*(h_{\nsb-1,1})=\ubracket{h_{\nsb-1,1}}\oplus\valustar_{\sigmae}^*(g_\nsb)$.
			Since $F_{\nsb,j}$ is closed for $\sigmae$, Properties (\ref{property: USV1})$_{\nsb}$ and (\ref{property: EV1})$_{\nsb+1}$ imply $\valustar_{\sigmae}^*(g_\nsb)=W_\nsb^*\oplus\valustar_{\sigmae}^*(b_{\nsb+1})$.
			As it also holds that $\valustar_{\sigma}^*(b_{\nsb+1})=L_{\nsb+1}^*$, it hence follows that \[\valustar_{\sigmae}^*(h_{\nsb-1,1})=\ubracket{h_{\nsb-1,1}}\oplus W_\nsb^*\cup L_{\nsb+1}^*\succ\bigoplus_{i=1}^{\nsb-1}W_{i}^*\oplus L_{\nsb+1}^*=R_1^*=\valustar_{\sigmae}^*(b_1).\]
			Thus $(s_{\nsb-1,1},h_{\nsb-1,1})\in I_{\sigmae}$.
			Also, \[\valustar_{\sigmae}^*(g_\nsb)=W_\nsb^*\oplus\valustar_{\sigmae}^*(b_{\nsb+1})\rhd\valustar_{\sigmae}^*(b_{\nsb+1})=\valustar_{\sigmae}^*(b_\nsb)\] since $\sigmae(b_\nsb)=b_{\nsb+1}$, implying $(b_\nsb,g_\nsb)\in I_{\sigmae}$.

			We argue why no further improving switches are created and that existing improving switches remain improving.
			The only vertices with edges to $g_\nsb$ are $s_{\nsb-1,1}$ and $b_\nsb$.
			It thus suffices to show that their valuations does not change.
			But this follows from $\sigmae(b_\nsb)=b_{\nsb+1}$ and $\sigmae(s_{\nsb-1,1})=b_1$.
			
			It remains to prove that $\sigmae$ is a phase-2-strategy.
			By the first two statements, it suffices to show that there is some $i<\nsb$ such that \Pref{USV3}$_{i}$ and the negations of both \Pref{EV2}$_{i}$ and \Pref{EV3}$_{i}$ are fulfilled as $\nsb=\relbit{\sigmae}$.
			Choose any $i<\nsb$. 
			Then, by our previous arguments, $\sigmae$ has \Pref{USV3}$_i$.			
			We next show that $\sigmae$ does not have \Pref{EV2}$_{i}$. 
			This follows from $\indbit_{i+1}^{\sigma}=1$, \Pref{EV1}$_{i}$, \Pref{EV2}$_{i}$ (both applied to $\sigma$) and $1-\indbit_{i+1}^{\sigmae}=\indbit^{\sigma}_{i+1}$.
			We finally show that $\sigmae$ does not have \Pref{EV3}$_{i}$.
			But this also immediately follows from $1-\indbit_{i+1}^{\sigmae}=\indbit_{i+1}^{\sigma}$ and by applying \Pref{EV1}$_{i}$ and \Pref{EV3}$_{i}$ to $\sigma$. 
			By definition, this also implies that $\sigmae$ is a pseudo phase-2-strategy if $\sigma(g_{\nsb})\neq F_{\nsb,j}$.\qedhere
\end{enumerate}
\end{proof}

\TransitionToPhaseTwo*

\begin{proof}
Let $j\coloneqq\bit_{\nsb+1}$.
We prove that $\sigmae$ is a phase-$2$-strategy for $\bit$.
By the choice of $e$, $\indbit^{\sigmae}=\indbit^{\sigma}=\bit+1\eqqcolon\indbit$.
As $\sigma$ has \Pref{REL2}, $\nsb=\relbit{\sigma}$.
Since $e\in I_{\sigma}$ implies $\sigma(b_\nsb)=b_{\nsb+1}$ by \Pref{EV2}$_{\nsb}$, we have $\incorrect{\sigma}=\emptyset$ as~$\sigma$ has \Pref{REL1}.
By the choice of $e$ and $\sigmae(b_\nsb)=\sigma(b_\nsb)=b_{\nsb+1}$, this implies $\incorrect{\sigmae}=\incorrect{\sigma}=\emptyset$ and $\relbit{\sigmae}=\relbit{\sigma}=\nsb$.
Hence $\sigmae$ has Properties (\ref{property: REL1}) and (\ref{property: REL2}).
By the choice of $e$, $\sigmae(g_\nsb)=F_{\nsb,j}$.
Hence \Pref{EV2}$_{i}$ remains valid for all $i\geq \nsb$.
It remains to show that there is an $i<\nsb$ such that $\sigmae$ has \Pref{USV3}$_{i}$ but not \Pref{EV2}$_{i}$ and \Pref{EV3}$_{i}$.
Since $\sigma$ is a pseudo phase-$2$-strategy for $\bit$, there is such an index fulfilling these conditions with respect to $\sigma$.
This index also fulfills these conditions with respect to $\sigmae$.
As~$\sigma$ being a pseudo phase-$2$-strategy implies that $\sigmae$ has the remaining properties, $\sigmae$ is a phase-$2$-strategy for~$\bit$.

Since $\sigma$ is well-behaved, $\relbit{\sigma}=\relbit{\sigmae}=\nsb\neq 1$ and a switch involving a selector vertex is applied we need to reevaluate the following properties.

\begin{itemize}[align=right, leftmargin=1.75cm]
	\item[(\ref{property: B3})] Assume that the premise was fulfilled by $\sigmae$.
		Then, by \Pref{USV1}$_\nsb$ and \Pref{EV1}$_{\nsb+1}$, $\sigmae(s_{\nsb,1})=h_{\nsb,1}$ implies $j=\indbit_{\nsb+1}=1$.
		Consequently, it holds that $\sigmae(b_{\nsb+1})=g_{\nsb+1}$, contradicting $\sigma(b_{\nsb+1})=b_{\nsb+2}$.
	\item[(\ref{property: EG4})] Since $\nsb>1$, the target of $g_1$ is not changed.
	\item[(EBG*)] Any premise requires a cycle center to escape towards both $g_1$ and $b_2$, contradicting \Pref{ESC2}.
	\item[(\ref{property: DN2})] Since $\sigmae$ is a pseudo phase-2$-$strategy for $\bit$ there is some $i$ such that \Pref{EV2}$_{i}$ is not fulfilled.
		This implies $\sigmae(b_{i})=g_{i}$.
\end{itemize}

We prove $I_{\sigmae}=\mathfrak{D}^{\sigmae}\cup\{(b_\nsb,g_\nsb),(s_{\nsb-1,1},h_{\nsb-1,1})\}$ and prove $(s_{\nsb-1,1},h_{\nsb-1,1})\in I_{\sigmae}$ first.
By \Pref{USV3}$_{\nsb-1}$, $\sigmae(s_{\nsb-1,1})= b_1$.
It thus suffices to prove $\valustar_{\sigmae}^*(h_{\nsb-1,1})\succ\valustar_{\sigmae}^*(b_1)$.
It holds that \[\valustar_{\sigmae}^*(h_{\nsb-1,1})=\ubracket{h_{\nsb-1,1}}\oplus\valustar_{\sigmae}^*(g_\nsb)=\ubracket{h_{\nsb-1,1}}\oplus W_\nsb^*\oplus\valustar_{\sigmae}^*(b_{\nsb+1})\] since $\sigmaebar(g_\nsb)=\indbit_{\nsb+1}$ and $\sigmae$ has \Pref{USV1}$_{\nsb}$.
Since $\relbit{\sigmae}=\nsb$, we also have that $\valustar_{\sigmae}^*(b_{\nsb+1})=L_{\nsb+1}^*$ and $\sigmaebar(b_{\relbit{\sigmae}})=b_{\relbit{\sigmae}+1}$.
The statement then follows since \Cref{corollary: Simplified MDP Valuation} implies $\valustar_{\sigma}^*(b_1)=R_1^*$.

We next show $(b_\nsb,g_\nsb)\in I_{\sigmae}$.
Since $\sigmae(b_\nsb)=b_{\nsb+1}$, we prove$\valustar_{\sigmae}^*(g_\nsb)\rhd\valustar_{\sigmae}^*(b_{\nsb+1})$.
This however follows since $\valustar_{\sigmae}^*(g_\nsb)=W_\nsb^*\oplus\valustar_{\sigmae}^*(b_{\nsb+1})$ as discussed previously.

It remains to show that improving switches remain improving and that no new improving switches are created.
By the choice of $e$, the valuation of $g_\nsb$ increases.
However, as discussed before, $\sigmae(b_\nsb)=b_{\nsb+1}$ and $\sigmae(s_{\nsb-1,1})=b_1$.
Since $b_\nsb$ and $s_{\nsb-1,1}$ are the only vertices that have an edge towards $g_\nsb$, the vertex $g_\nsb$ is the only vertex whose valuation changes when transitioning from $\sigma$ to $\sigmae$, implying the statement.
\end{proof}

\BeginningOfPhaseTwo*

\begin{proof}
	We first show that $\sigmae$ is a phase-$2$-strategy for $\bit$.
	Since the same set of cycle centers is closed for $\sigma$ and $\sigmae$, $\indbit^{\sigmae}=\indbit^{\sigma}=\bit+1\eqqcolon\indbit$.
	Thus \Pref{USV1}$_{i'}$ remains valid for all $i\geq\relbit{\sigma}$ and \Pref{CC1}$_{i}$ remains valid for all $i\in[n]$.
	We next show $\relbit{\sigma}=\relbit{\sigmae}$.
	By the choice of $e$,  $\sigmae(b_\nsb)=g_\nsb$.
	In addition, since $\sigma$ has \Pref{REL2}, $\relbit{\sigma}=\nsb$.
	Thus $\sigmae(b_{\nsb-1})=\sigma(b_{\nsb-1})=g_{\nsb-1}$ by \Cref{lemma: Traits of Relbit} as \Pref{REL1} applied to $\sigma$ implies $\incorrect{\sigma}=\emptyset$.
	Note that \Cref{lemma: Traits of Relbit} is applied to $\sigma$ which is well-behaved.
	Since $\sigma$ is well behaved and $\nsb-1<\relbit{\sigma}$ we have $\sigmae(g_{i-1})=\sigma(g_{i-1})=F_{i,0}$ by \Pref{BR1}.
	But then, since $\sigmaebar(b_{i})=\sigmabar(b_{i})$ for all $i\in[n], i\neq \nu$ and $\sigmaebar(g_{i})=\sigmabar(g_{i})$ for all $i\in[n]$, we have $\incorrect{\sigma}=\{\nu-1\}$.
	Since $\sigmaebar(g_{\nsb})=\sigmaebar(b_{\nsb+1})$ by \Pref{CC2}, it therefore follows that $\relbit{\sigmae}=\nsb$. 
	Thus, since $\sigma$ is phase-$2$-strategy, any statement regarding a level larger than $\nsb=\relbit{\sigma}=\relbit{\sigmae}$ remains valid.
	\Pref{EV1}$_{\nsb}$ and \Pref{EV2}$_{\nsb}$ follow directly from \Pref{CC2} and the choice of $e$.
	It remains to show that there is some $i<\relbit{\sigmae}$ such that \Pref{USV3}$_{i}$ as well as the negations of  both \Pref{EV2}$_{i}$ and \Pref{EV3}$_{i}$ hold.
	However, since $\sigma$ is a phase-$2$-strategy, there exists such an index for $\sigma$, so the same index can be used for $\sigmae$.

	Since we switched the target of $b_\nsb$ and $\nsb=\relbit{\sigmae}\neq 1$ we need to reevaluate the following assumptions to prove that $\sigmae$ is well-behaved.
	\begin{itemize}[align=right, leftmargin=1.75cm]
		\item[(\ref{property: S1})] Since $\sigma(g_\nsb)=F_{\nsb,\indbit_{\nsb+1}}$ by \Pref{EV2}$_{\nsb}$, the premise and the conclusion are true.
		\item[(\ref{property: S2})] This property only needs to be checked if $\relbit{\sigmae}=2$.
			Then, the only index for which the premise might become true is $i=1$. 
			But then, it cannot hold that $\sigmae(b_2)=g_2\wedge i>1$.
			Thus, the premise is either incorrect for $i=1$, implying that the implication is correct for $\sigmae$, or one of the other two conditions of the premise is true for $\sigmae$.
			But then, these conditions were also already true for $\sigma$, and hence $\sigmabar(s_i)=\sigmaebar(s_i)=1$ follows.
		\item[(\ref{property: B3})] As discussed earlier, $\sigmae(b_{\nsb-1})=g_{\nsb-1}$, hence the premise is incorrect.
		\item[(\ref{property: D1})] By \Pref{EV1}$_{\nsb}$ and \Pref{EV2}$_{\nsb}$, the conclusion is true, hence the implication.	
		\item[(\ref{property: D2})]  Again, this property only needs to be checked if $\relbit{\sigmae}=2$.
			But then, there is no $i\geq2$ with $i<\relbit{\sigmae}$, hence the premise is incorrect.
		\item[(\ref{property: EG5})] We only need to show that the premise is not true for $j=0$.
			It thus suffices to show that the cycle center $F_{\nsb-1,0}$ is closed.
			If $\relbit{\sigmae}>2$, then $\nsb-1>1$.
			By \Cref{lemma: Traits of Relbit}, it then holds that $\sigmae(b_{\nsb-1})=\sigma(b_{\nsb-1})=g_{\nsb-1}$.
			Hence, by Properties (\ref{property: D1}) and (\ref{property: BR1}), $\sigmaebar(d_{\nsb-1})=\sigmaebar(d_{\nsb-1,0})$.
			This in particular implies $\neg\sigmaebar(eb_{\nsb-1,0})$, so the premise is incorrect of $\relbit{\sigmae}>2$.
			Now consider the case $\relbit{\sigmae}=2$.
			Then, by the definition of a phase-$2$-strategy, the negation of \Pref{EV3}$_1$ holds.
			Thus, since $\indbit_{2}=1$ in this case, we have $\sigmaebar(d_{\nsb,1-\indbit_{2}})=\sigmaebar(d_{1,0})$.
	\end{itemize}
	
	We next prove that $I_{\sigmae}=\mathfrak{D}^{\sigmae}\cup\{(b_{\nsb-1},b_\nsb),(s_{\nsb-2,0},h_{\nsb-2,0}),(s_{\nsb-1,1},h_{\nsb-1,1})\}$ if $\nsb\neq 2$.
	We first show that $(b_{\nsb-1},b_\nsb)\in I_{\sigmae}$.
	Since $\sigmae(b_\nsb)=g_\nsb$ and $\sigmae(b_{\nsb-1})=g_{\nsb-1}$ it suffices to show $\valu_{\sigmae}^*(b_\nsb)\succ\valu_{\sigmae}^*(b_{\nsb-1})$. 
	This follows since $\valustar_{\sigmae}^*(b_\nsb)=L_{\nsb}^*\succ R_{\nsb-1}^*=\valustar_{\sigmae}^*(b_{\nsb-1})$ by  \Cref{lemma: VV Lemma}.
	
	We next show $(s_{\nsb-2,0},h_{\nsb-2,0})\in I_{\sigmae}$.
	By \Pref{USV3}$_{\nsb-2}$, $\sigmae(s_{\nsb-2,0})\neq h_{\nsb-2,0}$.
	Using $\sigmae(b_\nsb)=g_\nsb,\nsb=\relbit{\sigmae}, \indbit_{\nsb}=1$, \Pref{USV1}$_{\nsb}$, $(s_{\nsb-2,0},h_{\nsb-2,0})\in I_{\sigmae}$ follows from \begin{align*}
	\valustar_{\sigmae}^*(h_{\nsb-2,0})&=\ubracket{h_{\nsb-2,0}}\oplus W_{\nsb}^*\oplus L_{\nsb+1}^*\succ\bigoplus_{i<\nsb}W_{i}^*\cup L_{\nsb+1}^*=R_1^*=\valu_{\sigmae}^*(b_1).
	\end{align*}
	Using the same arguments yields $(s_{\nsb-1,1},h_{\nsb-1,1})\in I_{\sigmae}$.
	Since the valuation of all other vertices is unchanged, no other switch becomes improving and improving switches stay improving.
	
	We prove that $I_{\sigmae}=\mathfrak{D}^{\sigmae}\cup\{(b_1,b_2),(s_{1,1},h_{1,1})\}\cup\{(e_{*,*,*},g_1)\}$ if $\nsb=2$.
	All of the equations developed for the case $\nsb\neq 2$ are also valid for $\nsb=2$.
	In particular we have $\valustar_{\sigmae}^*(b_2)\succ\valustar_{\sigmae}^*(g_1)$ and $\sigmae(b_1)=g_1$, implying $(b_1,b_2)\in I_{\sigmae}$.
	In addition, we have $\sigmae(e_{i,j,k})=g_1$ for all $i\in[n]$ and $j,k\in\{0,1\}$, hence $(e_{i,j,k},b_2)\in I_{\sigmae}$ for these indices.
	By the usual arguments, no other new improving switches are created and improving switches stay improving (with the exception of $e$).
\end{proof}

\UpperSelectionVerticesInPhaseTwo*

\begin{proof}
	We first observe that $\sigma(s_{1,\indbit^{\sigma}_2})=b_1$ by \Pref{USV3}$_1$.	
	Since $\sigma$ has Properties (\ref{property: REL2}) and (\ref{property: EV1})$_{\relbit{\sigma}}$, it follows that $\sigma(b_{\relbit{\sigma}})=g_{\relbit{\sigma}}$.
	Thus, by \Cref{lemma: Traits of Relbit}, $\incorrect{\sigma}\neq\emptyset$.
	By the choice of $e$, $\indbit^{\sigma}=\indbit^{\sigmae}\eqqcolon\indbit, \relbit{\sigma}=\relbit{\sigmae}, \incorrect{\sigmae}=\incorrect{\sigma}\neq\emptyset$ and $\sigmae(b_{\relbit{\sigmae}})=g_{\relbit{\sigmae}}$.
	In particular, $\sigmae$ has Properties (\ref{property: EV1})$_{\relbit{\sigmae}}$ and (\ref{property: EV1})$_{i+1}$.
	Let $i\neq 1$.
	We prove that $\sigmae$ is a phase-2-strategy.	
	Since $i<\relbit{\sigmae}$, it suffices to check the special conditions of phase $2$ since all other properties of \Cref{table: Definition of Phases} remain valid for $\sigmae$.
	We show that the index $1$ fulfills these special conditions.
	Since $\relbit{\sigmae}\neq 1$, we have $\sigmae(b_1)=g_1$.
	As the choice of $i\neq 1$ implies $\relbit{\sigma}=\relbit{\sigmae}>2$, applying \Pref{BR1} to $\sigma$ yields $\sigmaebar(g_1)=\sigmabar(g_1)=1$.
	For the sake of a contradiction, assume that $\sigmae$ had \Pref{EV2}$_1$.
	Then, $1=\sigmabar(g_1)=\indbit_{2}$, implying $\nsb=\relbit{\sigma}=2$, contradicting the choice of $i$.
	Consequently, \Pref{EV2}$_1$ does not hold for $\sigmae$.
	Now, for the sake of contradiction, assume that $\sigmae$ had \Pref{EV3}$_1$.
	Then, since $\sigma(b_1)=g_1$ and $\nsb=\relbit{\sigma}>2$, the cycle center $F_{1,1-\indbit_2}=F_{1,1}$ is not closed.
	By \Pref{ESC2}, this implies $\sigmabar(eg_{1,1})\wedge\neg\sigmabar(eb_{1,1})$.
	Since $\sigmabar(g_1)=1$, \Pref{EG3} then implies $\sigmabar(s_1)=\sigmabar(s_{1,1})=1$.
	Consequently, by \Pref{EG5}, this implies $\sigmabar(b_2)=1$, so $\sigma(b_2)=g_2$.
	But then, $\sigma(b_2)=g_2\Leftrightarrow \relbit{\sigma}>2$ as both statements are true.
	Thus, since $\sigma(b_1)=g_1$, \Pref{D1} implies that $F_{1,\sigmabar(g_1)}=F_{1,1}$ is closed which is a contradiction.
	Hence $\sigmae$ does not have \Pref{EV3}$_1$.
	Finally, we have $\sigmae(s_{1,0})=\sigma(s_{1,0})=\sigma(s_{1,\indbit_{2}})=b_1$ by assumption and $\sigmaebar(s_{1,1})=\sigmabar(s_{1,1})=\sigmabar(s_{1,1-\indbit_2})=1$ by \Pref{S2}.
	Hence the index $1$ fulfills all of the special conditions of the definition of a phase-$2$-strategy, so $\sigmae$ is a phase-$2$-strategy for $\bit$.
	
	If $i=1$, then the assumptions imposed on $\sigma$ and the choice of $e$ directly imply that $\sigmae$ is a phase-$3$-strategy for $\bit$.

	We prove that $\sigmae$ is well-behaved.
	Note that $\sigmaebar(g_i)=\sigmabar(g_i)$ and thus, by \Pref{BR1}, $\sigmabar(g_i)=1$ if and only if $i\neq\relbit{\sigmae}-1$, implying $j=1-\sigmabar(g_i)$.
	By the usual arguments, it suffices to investigate the following properties.
	\begin{itemize}[align=right, leftmargin=1.75cm]
		\item[(\ref{property: B3})] We only need to consider this property if $j=1$, i.e., if $\indbit_{i+1}=1$.
			Since $i<\relbit{\sigmae}$ this implies that $i=\relbit{\sigmae}-1$.
			But then $\sigmae(b_{i+1})=\sigmae(b_{\relbit{\sigmae}})=g_{\relbit{\sigmae}}$, so the premise is incorrect.
		\item[(\ref{property: EG5})] Since $\sigmae$ fulfills \Pref{EV1}$_{i+1}$, $\sigmabar(b_{i+1})=\indbit_{i+1}=j.$
			Thus, the conclusion of \Pref{EG5} is correct, implying that the implication is correct.
\end{itemize}
	
	It remains to show $I_{\sigmae}=I_{\sigma}\setminus\{e\}$.
	The vertex $F_{i,j}$ is the only vertex that has an edge to $s_{i,j}$.
	Let $G_n=S_n$ first.
	Since $\valu_{\sigma}^\P(s_{i,j})\unlhd\valu_{\sigmae}^\P(s_{i,j})$, proving $\tau^{\sigma}(F_{i,j})\neq s_{i,j}$ implies $\tau^{\sigmae}(F_{i,j})\neq s_{i,j}$.		
	This then implies that the valuation of no other vertex than $s_{i,j}$ changes, implying $I_{\sigmae}=I_{\sigma}\setminus\{e\}$.
	
	For the sake of a contradiction, assume $\tau^{\sigma}(F_{i,j})=s_{i,j}$.
	Then, by \Cref{lemma: Exact Behavior Of Counterstrategy}, one of three cases holds.
	Since $\relbit{\sigmae}\neq 1$, it cannot hold that $[\sigmabar(eb_{i,j})\wedge\nsigmabar(eb_{i,j})\wedge\relbit{\sigma}=1]$.
	As $\sigmabar(d_{i,1-j})$ by \Pref{BR1} and assumption, \Pref{CC1}$_i$ implies $\nsigmabar(d_{i,j})$.
	Since $\sigmae$ has \Pref{ESC2}, $\sigmabar(eb_{i,j})\wedge\neg\sigmabar(eg_{i,j})\wedge[\relbit{\sigma}\neq 1\vee(\sigmabar(s_{i,j})\wedge\sigmabar(b_{i+1})\neq j)]$ also cannot hold.
	Consequently, by \Cref{lemma: Exact Behavior Of Counterstrategy}, $\tau^{\sigma}(F_{i,j})\neq s_{i,j}$.	
	
	Now let $G_n=M_n$.
	Again, as  $\sigmaebar(d_{i,1-j})$ by assumption, \Pref{CC2} implies that $F_{i,j}$ is not closed.
	By \Cref{lemma: Exact Behavior Of Random Vertex} and \Pref{ESC2}, this implies $\valustar_{\sigmae}^\M(F_{i,j})=\valustar_{\sigmae}^\M(g_1)$ and in particular $\valustar_{\sigmae}^\M(F_{i,j})\neq\valustar_{\sigmae}^\M(s_{i,j})$.
	The only vertices that have an edge to $F_{i,j}$ are $d_{i,j,0},d_{i,j,1}$ and $g_i$.
	We prove that $\sigma(d_{i,j,k})\neq F_{i,j}$ implies $\valu_{\sigma}^\M(F_{i,j})>\valu_{\sigma}^\M(e_{i,j,k})$, so $(d_{i,j,k},F_{i,j})\in I_{\sigma}$.	
	Since $\sigma(d_{i,j,k})=F_{i,j}$ implies $(d_{i,j,k},F_{i,j})\notin I_{\sigma}$, this then proves that $\sigma(d_{i,j,k})\neq F_{i,j}\Leftrightarrow(d_{i,j,k},F_{i,j})\in I_{\sigma}$.
	We then argue why the same arguments can be applied to $\sigmae$ which proves $(d_{i,j,k},F_{i,j})\in I_{\sigma}\Leftrightarrow(d_{i,j,k},F_{i,j})\in I_{\sigmae}$.
	
	Hence assume $\sigma(d_{i,j,k})\neq F_{i,j}$, implying $\sigma(d_{i,j,k})=e_{i,j,k}$.
	By \Pref{ESC2}, all escape vertices escape to $g_1$, hence $F_{i,j}$ is either $g_1$-open or $g_1$-halfopen.
	Also, $\sigma(s_{i,j})=b_1$ and $\sigma(b_1)=g_1$ imply $\valu_{\sigma}^\M(s_{i,j})=\rew{s_{i,j}}+\valu_{\sigma}^\M(g_1)$.
	Thus, \[\valu_{\sigma}^\M(F_{i,j})-\valu_{\sigma}^\M(e_{i,j,k})=q[\valu_{\sigma}^\M(s_{i,j})-\valu_{\sigma}^\M(g_1)],\] where the exact value of $q>0$ depends on whether $F_{i,j}$ is open or halfopen.
	But then, $\valu_{\sigma}^\M(s_{i,j})=\rew{s_{i,j}}+\valu_{\sigma}^\M(g_1)$ implies $\valu_{\sigma}^\M(F_{i,j})>\valu_{\sigma}^\M(e_{i,j,k})$.
	Since \[\valu_{\sigmae}^\M(s_{i,j})=\rew{s_{i,j}}+\valu_{\sigmae}^\M(h_{i,j})>\rew{s_{i,j}}+\valu_{\sigmae}^\M(g_1)\] as the edge $(s_{i,j},b_1)$ would otherwise be improving for $\sigmae$ which cannot happen, the same argument implies $\valu_{\sigmae}^\M(F_{i,j})>\valu_{\sigmae}^\M(e_{i,j,k})$.
		
	No vertex but $F_{i,j}$ has an edge towards $d_{i,j,k}$.
	Thus, although the valuation of $d_{i,j,k}$ increases due to the application of $e$, it is impossible to have an improving switch $(*,d_{i,j,k})$ for either $\sigma$ or $\sigmae$.
	Consequently, we do not need to consider this vertex when investigating whether new improving switches are created.

	It thus remains to prove $\sigma(g_i)=\sigmae(g_i)=F_{i,1-j}$ and $(g_{i},F_{i,j})\notin I_{\sigma},I_{\sigmae}$.
	Once this statement is proven, combining all of the previous statements yields $I_{\sigmae}=I_\sigma\setminus\{e\}$.
	Since \Pref{BR1}$_i$ and the choice of $j$ imply that $\sigmae(g_i)=F_{i,1-j}$, it suffices to prove $\valu_{\sigmae}^\M(F_{i,1-j})>\valu_{\sigmae}^\M(F_{i,j})$.
	Since $\sigmae(b_1)=g_1$, the assumption $\sigmaebar(d_{i'})$ for all $i'<\relbit{\sigmae}$ implies $\valustar_{\sigmae}^\M(F_{i,j})=\valustar_{\sigmae}^\M(g_1)=R_1^\M.$
	This furthermore yields $\valustar_{\sigmae}^\M(F_{i,1-j})=\valustar_{\sigmae}^\M(s_{i,1-j})$.
	\Pref{USV2}$_{i,1-j}$ implies that $\sigmae(s_{i,1-j})=h_{i,1-j}$.
	If $j=\indbit_{i+1}=1$, then \[\valustar_{\sigmae}^\M(F_{i,1-j})=\ubracket{s_{i,1-j},h_{i,1-j}}+\valustar_{\sigmae}^\M(b_{i+2}).\]
	But this implies $i=\nsb-1$, so \[\valustar_{\sigmae}^\M(F_{i,1-j})=\ubracket{s_{\nsb-1,1-j},h_{\nsb-1,1-j}}+\valustar_{\sigmae}^\M(b_{\nsb+1})>\sum_{\ell=1}^{\nsb-1}W_{\ell}^{\M}+L_{\nsb+1}^{\M}=R_1^\M.\]
	If $j=\indbit_{i+1}=0$, then $i<\nsb-1$ and $\valustar_{\sigmae}^\M(F_{i,1-j})=\ubracket{s_{i,1-j},h_{i,1-j}}+\valustar_{\sigmae}^\M(g_{i+1})$.
	Using \Cref{lemma: Valuation of g if level small}, this implies \begin{align*}
		\valustar_{\sigmae}^\M(F_{i,1-j})&=\ubracket{s_{i,1-j},h_{i,1-j}}+R_{i+1}^\M\\
			&=\ubracket{s_{i,1-j},h_{i,1-j}}+\sum_{\ell=i+1}^{\nsb-1}W_{\ell}^{\M}+L_{\nsb+1}^{\M}>\sum_{\ell=1}^{\nsb-1}W_{\ell}^{\M}+L_{\nsb+1}^{\M}=R_1^\M.
		\end{align*}\qedhere
\end{proof}

\ResettingEntryVertices*

\begin{proof}
	By the choice of $e$ and by assumption, $\sigmae$ has \Pref{USV3}$_{i'}$ for all $i'<i$.
	In particular, $\sigma(s_{i-1,0})=\sigma(s_{i-2,0})=b_1$ since $i<\relbit{\sigma}=\nsb$.  
	
	We prove that $\sigmae$ is a phase-$2$-strategy for $\bit$.
	Since $e=(b_i,b_{i+1})$ and $i<\relbit{\sigma}=\nsb$, it suffices to prove $\relbit{\sigma}=\relbit{\sigmae}$ and that there is an index $i'<i$ fulfilling the special conditions of a phase-$2$-strategy.
	Since $\relbit{\sigma}>i, \relbit{\sigma}=\nsb$ and since $\sigma$ has \Pref{EV1}$_{\relbit{\sigma}}$, $\sigma(b_{\relbit{\sigma}})=g_{\relbit{\sigma}}$.
	Hence, by \Cref{lemma: Traits of Relbit}, $\incorrect{\sigma}\neq\emptyset$.
	We show $i=\max\{i'\in\incorrect{\sigma}\}$.
	By the choice of $e$, $\sigmabar(b_i)=1$.
	If $i+1=\relbit{\sigma}$, then $\sigmabar(b_{i+1})=\sigmabar(b_{\relbit{\sigma}})=1$ and $\sigmabar(g_i)=0$ by \Pref{BR1}.
	If $i+1<\relbit{\sigma}$, then $i+1\leq \relbit{\sigma}-1$, so $\sigmabar(b_{i+1})=\indbit^{\sigma}_{i+1}=(\bit+1)_{i+1}=0$ by \Pref{EV1}$_{i+1}$ and $\sigmabar(g_i)=1$ by \Pref{BR1}.
	In either case $\sigmabar(g_i)\neq\sigmabar(b_{i+1})$, hence $i\in\incorrect{\sigma}$.
	For any $i'\in\{i+1,\dots\nsb-1\}$, \Pref{EV1}$_{i'}$ and $\relbit{\sigma}=\nsb$ imply $\sigmabar(b_{i'})=0$.
	Thus $i=\max\{i'\in\incorrect{\sigma}\}$.
	We now prove $i-1\in\incorrect{\sigmae}$ since this suffices to prove $\relbit{\sigmae}=\relbit{\sigma}$ as $\sigmae(b_{i'})=b_{i'+1}$ for all $i'\in\{i,\dots,\relbit{\sigma}-1\}$.
	
	By \Pref{EV1}$_{i-1}$, it holds that $\sigmaebar(b_{i-1})=\sigmabar(b_{i-1})=1$ since \Pref{B2} would imply $\sigma(b_i)=b_{i+1}$ otherwise.
	By \Pref{BR1} and $i-1<\relbit{\sigma}-1$, it follows that $\sigmaebar(g_{i-1})=\sigmabar(g_{i-1})=1$.
	Also, $\sigmaebar(b_i)=0$ by the choice of $e$.
	Hence $i-1\in\incorrect{\sigmae}$, implying $i-1=\max\{i'\in\incorrect{\sigmae}\}$.
	Consequently, $\relbit{\sigmae}=\relbit{\sigma}=\nsb$, so $\sigmae$ has \Pref{REL2}.
	
	We show that $i-1$ fulfills the special conditions of \Cref{table: Definition of Phases} for phase-$2$-strategies.
	As shown previously, $\sigmaebar(b_{i-1})=1$ and $\sigmaebar(g_{i-1})=1=1-\indbit_{i}$.
	Thus \Pref{EV2}$_{i-1}$ does not hold for $\sigmae$.
	If $i>2$, \Cref{lemma: Traits of Relbit} implies $\sigma(b_2)=g_2$ as $i=\max\{i'\in\incorrect{\sigma}\}$.
	If $i=2$ then $\sigma(b_2)=g_2$ since $(b_2,b_3)\in I_{\sigma}$.
	Thus, by applying \Pref{D2} to $\sigma$, it follows that $\sigmabar(d_{i-1})=\sigmaebar(d_{i-1})=\sigmaebar(d_{i-1,1-\indbit_{i}})=1$.
	Thus $\sigmae$ fulfills the negation of \Pref{EV3}$_{i-1}$.
	Finally, $\sigmae$ also has \Pref{USV3}$_{i-1}$ by assumption.
	Thus the index $i-1$ fulfills the special conditions of \Cref{table: Definition of Phases}, so $\sigmae$ is phase-$2$-strategy for $\bit$.

	Since $\sigma$ is a phase-2-strategy, it suffices to check the following properties:
	\begin{enumerate}[align=right, leftmargin=1.75cm]
		\item[(\ref{property: B1})] If $i<\relbit{\sigmae}-1$, then $i+1<\relbit{\sigmae}$.
			Since $\sigma$ has \Pref{EV1}$_{i+1}$ by assumption, \Pref{REL2} implies $\sigmae(b_{i+1})=\sigma(b_{i+1})=b_{i+2}$.
		\item[(\ref{property: B3})] For this property, it might happen that either the premise becomes true with respect to $\sigmae$ or that it is true while the conclusion becomes false for $\sigmae$.
			Consider the first case first.
			Then, $\sigmae(s_{i-1,1})=h_{i-1,1}$ and $\sigmae(b_i)=b_{i+1}$.
			However, since $i=\max\{i'\in\incorrect{\sigmae}\}$, we have $\sigmaebar(g_i)=\sigmabar(g_i)\neq\sigmabar(b_{i+1})=\sigmaebar(b_{i+1})$, hence the conclusion is true as well.
			Now assume that the premise is correct for $\sigmae$ while the conclusion became false by the application of $e$.
			Then $\sigmae(s_{i-2,1})=h_{i-2,1}$ and $\sigmae(b_{i-1})=b_{i}$.
			But this cannot happen since $\sigmae(b_{i-1})=g_{i-1}$ as proven earlier.
		\item[(\ref{property: EG5})] Since $\sigmaebar(b_i)=0$ we need to show $\sigmaebar(d_{i-1,1})=1$ which was already shown earlier.
	\end{enumerate}

	We now show that $i\neq 2$ implies $I_{\sigmae}=(I_{\sigma}\setminus\{e\})\cup\{(b_{i-1},b_i),(s_{i-2,0},h_{i-2,0})\}$ and that $i=2$ implies $I_{\sigmae}=(I_{\sigma}\setminus\{e\})\cup\{(b_1,b_2)\}\cup\{(e_{*,*,*},g_1)\}$.
	By \Cref{lemma: Traits of Relbit}, $\sigmae(b_{i-1})=g_{i-1}$ and also, by assumption, $\sigma(s_{i-2,0})=\sigmae(s_{i-2,0})=b_1$ if $i\neq 2$.
	We hence need to show $\valustar_{\sigmae}^*(b_i)\succ\valustar_{\sigmae}^*(g_{i-1})$ and $\ubracket{h_{i-2,0}}\oplus\valustar_{\sigmae}^*(b_i)\succ\valustar_{\sigmae}^*(b_1)$.
	This in particular implies that any edge $(e_{*,*,*},b_2)$ is improving for $\sigmae$ if $i=2$ since $\sigma(e_{*,*,*})=g_1$ by \Pref{ESC2}.
	Since either $i+1<\relbit{\sigmae}$ and $\sigmae(b_{i+1})=b_{i+2}$ by \Pref{B1} or $i+1=\relbit{\sigmae}$, we have $\valustar_{\sigmae}^*(b_{i})=L_{i+1}^*$ since $\sigmae(b_i)=b_{i+1}$.
	By assumption, $\sigmabar(d_{i'})=\sigmaebar(d_{i'})=1$ for all $i'<\relbit{\sigmae}$.
	Thus, $\valustar_{\sigmae}^*(g_{i-1})=R_{i-1}^*$ by \Cref{lemma: Valuation of g if level small}.		
	Therefore, $\sigmae(b_{\relbit{\sigmae}})=g_{\relbit{\sigmae}}$ implies $\valustar_{\sigmae}^*(b_i)\succ\valustar_{\sigmae}^*(g_{i-1})$ by \Cref{lemma: VV Lemma}.
	Since $\valustar_{\sigmae}^*(b_1)=R_1^*$, \Cref{lemma: VV Lemma} further implies $\ubracket{h_{i-2,0}}\oplus\valustar_{\sigmae}^*(b_i)\succ\valustar_{\sigmae}^*(b_1)$ if $i\neq 2$.
	Thus $(b_{i-1},b_i),(s_{i-2,0},h_{i-2,0})\in I_{\sigmae}$ if $i\neq 2$ and $(b_1,b_2),(e_{*,*,*},b_2)\in I_{\sigma}$ if $i=2$.
\end{proof}

\PossibleBeginningOfPhaseThree*

\begin{proof}
Let $j\coloneqq\bit_{\nsb+1}$.
We first show $\relbit{\sigmae}=\relbit{\sigma}$.
Since $\sigma$ is a pseudo phase-3-strategy for $\bit$, $\relbit{\sigma}=\nsb=1$.
Hence, $\sigma(b_1)=b_2$, implying $\sigmae(b_1)=b_2$, so $\relbit{\sigmae}=1$.
Note that this implies that $\sigmae$ has \Pref{CC2} as the cycle center $F_{\nsb,j}$ is closed due to $\indbit^{\sigma}=\indbit^{\sigmae}=\bit+1\eqqcolon\indbit$.

We next show that $\sigmae$ is a phase-3-strategy for $\bit$.
The only properties other than \Pref{CC2} involving $e=(g_\nsb,F_{\nsb,j})$ are Properties (\ref{property: REL1}) and (\ref{property: EV2})$_1$.
These do not need to be fulfilled for a phase-3-strategy so, $\sigma$ being a pseudo phase-$3$-strategy implies that $\sigmae$ is a phase-$3$-strategy for $\bit$.

We now show that $\sigmae$ is well-behaved.
Since $\sigma$ is a well-behaved pseudo phase-$3$-strategy for $\bit$, $\relbit{\sigmae}=1$ and by the choice of $e$, it suffices to investigate the following properties:
\begin{itemize}[align=right, leftmargin=1.75cm]
	\item[(\ref{property: MNS1})] Since the cycle center $F_{1,j}$ is closed with respect to $\sigmae$ due to $\indbit^{\sigmae}=\bit+1$ and since $j=\sigmaebar(g_i)$, the conclusion of this property is true for $\sigmae$.
	\item[(\ref{property: MNS2})] If $\sigmae(b_2)=g_2$, then $\minsige{b}=2$ and there cannot be an index fulfilling the conditions of the premise.
		If $\sigmae(b_2)=b_3$, then $\indbit_2=j=0$ by \Pref{EV1}$_2$, implying $\sigmae(g_1)=F_{1,0}$.
		But then $\minnegsige{g}=1$, hence there cannot be an index such that the conditions of the premise are fulfilled.
	\item[(\ref{property: MNS3})] This follows by the same arguments used for \Pref{MNS2}.
	\item[(\ref{property: MNS4})] By the choice of $e$, the definition of $j$ and \Pref{EV1}$_2$, we either have $\minnegsige{g}=1\wedge\minnegsige{b}>2$ or $\minnegsige{g}>1\wedge\minsige{b}=2$.
		Since the second case contradicts the conditions of the premise, assume $\minnegsige{g}=1\wedge\minsige{b}>2$.
		Then, $\sigmae(b_2)=b_3$, hence $j=\indbit^{\sigmae}_2=0$, implying $\sigmae(g_1)=F_{1,0}$.
		In addition, $\minnegsige{s}=1$ as $\minnegsige{s}\leq\minnegsige{g}$.
		Thus, by the definition of $\minnegsige{s}$, we have $\sigmae(s_{1,0})=b_1$.
		But this contradicts \Pref{USV1}$_1$ as this implies $\sigmae(s_{1,\indbit_2})=\sigmae(s_{1,0})=h_{1,0}$.
	\item[(\ref{property: MNS5})] This follows by the same arguments used for \Pref{MNS2}.
	\item[(\ref{property: MNS6})] By the choice of $e$, the definition of $j$ and \Pref{EV1}$_2$, we either have $\minnegsige{g}=1\wedge\minnegsige{b}>2$ or $\minnegsige{g}>1\wedge\minsige{b}=2$.
		Since the first case contradicts the conditions of the premise, assume $\minnegsige{g}>1\wedge\minsige{b}=2$.
		Then, $\sigmae(b_2)=g_2$, hence $j=\indbit_2=1$, implying $\sigmae(g_1)=F_{1,1}$.
		In addition, $\minnegsige{s}=1$ as $\minnegsige{s}<\minsige{b}=2$.
		Thus, by the definition of $\minnegsige{s}$, we have $\sigmae(s_{1,1})=b_1$.
		But this contradicts \Pref{USV1}$_1$ as this implies $\sigmae(s_{1,\indbit^{\sigmae}_2})=\sigmae(s_{1,1})=h_{1,1}$.
	\item[(\ref{property: EG4})] Since the conclusion is true for $\sigmae$ by the choice of $e$, the implication is true.
	\item[(\ref{property: EBG2})] Since $\sigmaebar(g_1)=\sigmaebar(b_2)=\indbit_2=j$, $\relbit{\sigmae}=1$ and \Pref{USV1}$_1$ together imply $\sigmaebar(s_1)=1$.
	\item[(\ref{property: EBG3})] Since $\sigmaebar(g_1)=\sigmaebar(b_2)=\indbit_2=j$, $\nsb=1$ and $\indbit=\bit+1$ imply $\sigmaebar(d_1)=1$.
\end{itemize}

It thus remains to show that $I_{\sigmae}=(I_{\sigma}\setminus\{e\})\cup\{(b_1,g_1)\}\cup\{(e_{*,*,*},g_1)\}$.
Since $\sigmae(b_1)=b_2$ and $\relbit{\sigmae}=1$, this can be shown by using the same arguments used in the proof of \Cref{lemma: Closing active CC}~(4).
\end{proof}

\EscapeVerticesPhaseThree*

\begin{proof}
	Since a phase-3-strategy does not need to fulfill \Pref{ESC1} or (\ref{property: ESC2}), $\sigma$ being a phase-3-strategy for $\bit$ implies that $\sigmae$ is a phase-3-strategy for $\bit$.	
	We prove that $\sigmae$ is well-behaved.
	By assumption, $\sigmae(d_{i,j,k})=\sigma(d_{i,j,k})=F_{i,j}$.
	Hence $F_{i,j}$ escapes towards $g_1$ resp. $b_2$ with respect to $\sigma$ if and only if it escapes towards the same vertex with respect to $\sigmae$.
	Since there are no other conditions on escape vertices except the escape of cycle centers in \Cref{table: Well behaved properties}, $\sigmae$ is well-behaved since $\sigma$ is well-behaved.

	It remains to show the statements related to the improving switches.
	We first prove that $\sigma(d_{i,j,1-k})=e_{i,j,1-k}$ implies $I_{\sigmae}=(I_{\sigma}\setminus\{e\})\cup\{(d_{i,j,k},e_{i,j,k})\}$.
	As $d_{i,j,k}$ is the only vertex having an edge to $e_{i,j,k}$, it suffices to prove $(d_{i,j,k},e_{i,j,k})\in I_{\sigmae}$.
	We distinguish two cases.
	\begin{enumerate}
		\item Let $\relbit{\sigmae}=1$. 
			Then $\valustar_{\sigmae}^*(d_{i,j,k})=\valustar_{\sigmae}^*(F_{i,j})$ and $\valustar_{\sigmae}^*(e_{i,j,k})=\valustar_{\sigmae}^*(g_1)$.
			Moreover, $\valustar_{\sigmae}^*(g_1)=W_1^*\oplus\valustar_{\sigmae}^*(b_2)$ and $\valustar_{\sigmae}^*(b_2)=L_2^*$  by \Cref{lemma: Valuations in Phase Three}.
			We thus prove $\valustar_{\sigmae}^*(g_1)\succ\valustar_{\sigmae}^*(F_{i,j})$ and distinguish two further cases.
			\begin{enumerate}
				\item Let $\sigmae(e_{i,j,1-k})=g_1$.
					Since $\sigmae(d_{i,j,1-k})=e_{i,j,1-k}$ and $\sigmae(d_{i,j,k})=F_{i,j}$ we then have $\sigmaebar(eg_{i,j})\wedge\nsigmaebar(eb_{i,j})$ and, by assumption, $\relbit{\sigmae}=1$.
					Let $G_n=S_n$.
					Then, \Cref{lemma: Exact Behavior Of Counterstrategy} implies $\valustar_{\sigmae}^\P(F_{i,j})=\valustar_{\sigmae}^\P(s_{i,j})$.
					Since \Pref{EG1} implies $\sigmae(s_{i,j})=b_1$ and $\sigmae(b_1)=b_2$ follows from $\relbit{\sigmae}=1$, it holds that $\valustar_{\sigmae}^\P(F_{i,j})=\{s_{i,j}\}\cup\valustar_{\sigmae}^\P(b_2)$.
					Hence, by \Cref{lemma: Valuations in Phase Three}, $\valustar_{\sigmae}^\P(g_1)\rhd\valustar_{\sigmae}^\P(F_{i,j}),$ so $(d_{i,j,k},e_{i,j,k})\in I_{\sigmae}$.
					Now let $G_n=M_n$.
					Then, since $F_{i,j}$ is $g_1$-halfopen, \begin{align*}
						\valu_{\sigmae}^\M(g_1)-\valu_{\sigmae}^\M(F_{i,j})=\frac{2\e}{1+\e}[\valu_{\sigma}^\M(g_1)-\valu_{\sigma}^\M(s_{i,j})],
					\end{align*}
					so it suffices to prove $\valu_{\sigmae}^\M(s_{i,j})<\valu_{\sigmae}^\M(g_1)$.
					This follows by the previous arguments since $\valu_{\sigmae}^\M(s_{i,j})=\rew{s_{i,j}}+\valu_{\sigmae}^\M(b_2)<W_1^\M+\valu_{\sigmae}^\M(b_2)=\valu_{\sigmae}^\M(g_1).$
			
				\item Now assume $\sigmae(e_{i,j,1-k})=b_2$.
					By the same arguments used in case 1(a) this implies $\sigmaebar(eb_{i,j})\wedge\nsigmaebar(eg_{i,j})$.
					Let $G_n=S_n$ and consider \Cref{lemma: Exact Behavior Of Counterstrategy}.
					Either the conditions of case four or of case five are then fulfilled.
					If the conditions of case four are true, then $\valustar_{\sigmae}^\P(F_{i,j})=\valustar_{\sigmae}^\P(b_2)$.
					Since $\valustar_{\sigmae}^\P(g_1)=W_1^\P\cup\valustar_{\sigmae}^\P(b_2)$, this implies $\valustar_{\sigmae}^\P(g_1)\rhd\valustar_{\sigmae}^\P(F_{i,j}).$
					For the sake of a contradiction, assume that the conditions of case five were true.
					Then $\sigmaebar(s_{i,j})$ and $\sigmaebar(b_{i+1})\neq j$.
					But then \Pref{USV1}$_i$ implies $j=\indbit^{\sigmae}_{i+1}$ and \Pref{EV1}$_{i+1}$ implies $j\neq\indbit_{i+1}$ which is a contradiction.
					If $G_n=M_n$, then the statement follows directly since $\valustar_{\sigmae}^\M(F_{i,j})=\valustar_{\sigmae}^\M(b_2)$.
			\end{enumerate}
		\item Let $\relbit{\sigmae}\neq 1$.
			Then, $\valustar_{\sigmae}^*(e_{i,j,k})=\valustar_{\sigmae}^*(b_2)$ and we prove $\valustar_{\sigmae}^*(F_{i,j})\prec\valustar_{\sigmae}^*(b_2)$.
			\begin{enumerate}
				\item Assume $\sigmae(e_{i,j,1-k})=b_2$.
					Then $\sigmaebar(eb_{i,j})\wedge\nsigmaebar(eg_{i,j})$ and $\relbit{\sigmae}\neq 1$ by assumption.
					From \Pref{EB1} and \Pref{EV1}$_{i+1}$, it follows that $j\neq\indbit_{i+1}$.
					Let $G_n=S_n$.
					By \Cref{lemma: Exact Behavior Of Counterstrategy}, $\valustar_{\sigmae}^\P(F_{i,j})=\valustar_{\sigmae}^\P(s_{i,j})$.
					Consider the case $\sigmae(s_{i,j})=b_1$ first.
					Then $\valustar_{\sigmae}^\P(F_{i,j})=\{s_{i,j}\}\cup\valustar_{\sigmae}^\P(b_1)=\{s_{i,j}\}\cup R_1^\P.$
					Since $\sigmae(b_{\relbit{\sigmae}})=g_{\relbit{\sigmae}}$ and $W_{\relbit{\sigmae}}^\P\rhd\bigcup_{\ell<\relbit{\sigmae}}W_{\ell}^\P\cup\{s_{i,j}\}$, this implies  \begin{align*}
						\valustar_{\sigmae}^\P(F_{i,j})&=\{s_{i,j}\}\cup R_1^\P\lhd W_{\relbit{\sigmae}}^\P\cup\bigcup_{\ell\geq\relbit{\sigmae}+1}\{W_{\ell}^\P\colon\sigmae(b_{\ell})=g_{\ell}\}\\
							&=L_{\relbit{\sigmae}}^{\P}=L_2^{\P}=\valustar_{\sigmae}^\P(b_2).
					\end{align*}
					Next let $\sigmae(s_{i,j})=h_{i,j}$.
					Then, since $j=1-\indbit_{i+1}$ and since $\sigmae$ has \Pref{USV1}$_{\ell}$ for all $\ell\geq\relbit{\sigmae}$, we must have $i<\relbit{\sigmae}$.
					Consequently, \[\valustar_{\sigmae}^\P(F_{i,j})=\{s_{i,j},h_{i,j}\}\cup\begin{cases}\valustar_{\sigmae}^\P(g_{i+1}),&j=1\\\valustar_{\sigmae}^\P(b_{i+2}), &j=0\end{cases}.\]
					We now focus on the case $j=1$ and continue considering the case $G_n=S_n$.
					Then,  $\valustar_{\sigmae}^\P(F_{i,j})=\{s_{i,j},h_{i,j}\}\cup\valustar_{\sigmae}^\P(g_{i+1})$ as \Pref{BR1} implies $i<\relbit{\sigmae}-1$. 
					We now determine $\valustar_{\sigmae}^\P(g_{i+1})$ using \Cref{corollary: Complete Valuation Of Selection Vertices PG}.
					By \Pref{BR1}, it holds that $\sigma(g_{i'})=F_{i',1}$ for all $i'<\relbit{\sigmae}-1$ as well as $\sigmae(g_{\relbit{\sigmae}-1})=F_{\relbit{\sigmae}-1,0}$.
					By assumption, we also have $\sigmae(s_{i',j'})=h_{i',j'}$ for all $i'<\relbit{\sigmae}$ and $j'\in\{0,1\}$.
					Since $\sigma(b_{i'})=b_{i'+1}$ for all $i'<\relbit{\sigmae}$, we obtain $\lambda_{i+1}^\P=\relbit{\sigmae}-1$.
					By \Pref{BR2}, $\nsigmaebar(eg_{\lambda_{i+1}^\P})$.
					But this implies that the conditions of the sixth case of \Cref{corollary: Complete Valuation Of Selection Vertices PG} are fulfilled, hence \begin{align*}
						\valustar_{\sigmae}^\P(F_{i,j})&=\{s_{i,j},h_{i,j}\}\cup\bigcup_{\ell=i+1}^{\relbit{\sigmae}-1}W_{\ell}^\P\cup\valustar_{\sigmae}^\P(b_{\relbit{\sigmae}+1})\\
							&\lhd W_{\relbit{\sigmae}}^\P\cup\valustar_{\sigmae}^\P(b_{\relbit{\sigmae}+1})=\valustar_{\sigmae}^\P(b_2).
						\end{align*}
					
					Let $G_n=M_n$ and $j=1$.
					It then suffices to prove $\valu_{\sigmae}^\M(s_{i,j})<\valu_{\sigmae}^\M(b_2)$.
					Consider the case $\sigmae(s_{i,j})=b_1$, implying $\valustar_{\sigmae}^\M(s_{i,j})=\rew{s_{i,j}}+\valustar_{\sigmae}^\M(b_1)$.
					If $\valustar_{\sigmae}^\M(b_1)=B_1^\M,$ then $\valustar_{\sigmae}^\M(b_1)=R_1^\M$ and the arguments follows by the same arguments used for the case $G_n=S_n$.
					If $\valustar_{\sigmae}^\M(b_1)=B_2^\M+\sum_{j'<k}W_{j'}^\M+\rew{g_k}$ where $k=\min\{i'\geq1\colon\nsigmaebar(d_{i'})\}<\relbit{\sigma}$, then \[\valustar_{\sigmae}^\M(s_{i,j})=\rew{s_{i,j},g_k}+\sum_{j'<k}W_{j'}^k+B_2^\M<B_2^\M=\valustar_{\sigmae}^\M(b_2).\]
					Let $\sigmae(s_{i,j})=h_{i,j}$.
					As $j=1$, it holds that $\valustar_{\sigmae}^\M(F_{i,j})=\rew{s_{i,j},h_{i,j}}+\valustar_{\sigmae}^\M(g_{i+1})$.
					We use \Cref{corollary: Complete Valuation Of Selection Vertices MDP} to evaluate $\valustar_{\sigmae}^\M(g_{i+1})$ and thus determine $\lambda_{i+1}^\M$.
					Assume there was some index $i'\in\{i+1,\dots,\relbit{\sigmae}-1\}$ such that $\sigmae(s_{i',\sigmabar(g_{i'})})=b_1$.
					Then, by assumption, $\sigmaebar(eb_{i'})\wedge\nsigmaebar(eg_{i'})$.
					As $\sigmae(b_{i'})=b_{i'+1}$ by \Pref{EV1}$_{i'}$, \Cref{corollary: Complete Valuation Of Selection Vertices MDP} then implies the statement.

					Hence assume there was no such index.
					By \Pref{BR1}, $\lambda_{i+1}^\M\leq\relbit{\sigmae}-1$.
					The case $\sigmaebar(d_{\lambda_{i+1}^\M})\wedge\nsigmaebar(s_{\lambda_{i+1}^\M})$ cannot happen by assumption.
					Also, by \Pref{BR2}, $\nsigmaebar(eg_{\lambda_{i+1}^\M})$.
					Consequently, either $\sigmaebar(eb_{\lambda_{i+1}^\M})\wedge\nsigmaebar(eg_{\lambda_{i+1}^\M})$ and the statement follows by the same arguments used before or $\sigmaebar(d_{\lambda_{i+1}^\M})$.
					This however implies $\lambda_{i+1}^\M=\relbit{\sigmae}-1$ and thus \begin{align*}
						\valustar_{\sigmae}^\M(s_{i,j})&=\rew{s_{i,j,h_{i,j}}}+\sum_{\ell=i+1}^{\relbit{\sigmae}-1}W_{\ell}^\M+\valustar_{\sigmae}^\M(b_{\relbit{\sigmae}+1})\\
							&<W_{\relbit{\sigmae}}^\M+\valustar_{\sigmae}^\M(b_{\relbit{\sigmae}+1})=\valustar_{\sigmae}^\M(b_2).
						\end{align*}

					This concludes the case $j=1$, hence let $j=0$.
					Then, $\indbit_{i+1}=1$, so $i=\relbit{\sigmae}-1$ and the statement follows from $\valustar_{\sigmae}^*(F_{i,j})=\ubracket{s_{\relbit{\sigmae}-1,0},h_{\relbit{\sigmae}-1,0}}\oplus L_{\relbit{\sigmae}+1}^*$ and $\ubracket{s_{\relbit{\sigmae}-1,0},h_{\relbit{\sigmae}-1,0}}\prec W_{\relbit{\sigmae}}^{*}$.
				\item Now assume $\sigmae(e_{i,j,1-k})=g_1$.
					Then $\sigmaebar(eg_{i,j})\wedge\nsigmaebar(eb_{i,j})\wedge\relbit{\sigmae}\neq 1$.
					Thus, by \Cref{lemma: Exact Behavior Of Random Vertex,lemma: Exact Behavior Of Counterstrategy}, $\valustar_{\sigmae}^*(F_{i,j})=\valustar_{\sigmae}^*(g_1).$
					But then \Cref{lemma: Valuations in Phase Three} implies $\valustar_{\sigma}^*(F_{i,j})\prec\valustar_{\sigma}^*(b_2)$.
			\end{enumerate}
	\end{enumerate}
	We now show that $\sigmae(d_{i,j,1-k})=F_{i,j}$ and $j\neq\indbit_{i+1}$ also imply $(d_{i,j,k},e_{i,j,k})\in I_{\sigmae}$ by proving \begin{equation} \label{equation: Valuation in Escape Vertices Phase Three} 
		\valustar_{\sigmae}^*(F_{i,j})\prec\valustar_{\sigmae}^*(e_{i,j,k}).
	\end{equation}
	In this case, $F_{i,j}$ is closed.
	Thus, by \Cref{lemma: Exact Behavior Of Random Vertex,lemma: Exact Behavior Of Counterstrategy}, $\valustar_{\sigmae}^*(F_{i,j})=\valustar_{\sigmae}^*(s_{i,j})$.
	If $\sigmae(s_{i,j})=b_1$, then this implies $\valustar_{\sigmae}^*(F_{i,j})=\ubracket{s_{i,j}}\oplus\valustar_{\sigmae}^*(b_1)$.
	If $\relbit{\sigmae}=1$, then $\valustar_{\sigmae}^*(e_{i,j,k})=\valustar_{\sigmae}^*(g_1)$.
	We then obtain (\ref{equation: Valuation in Escape Vertices Phase Three}) as in case 1(a).
	
	If $\relbit{\sigmae}\neq 1$, then $\valustar_{\sigmae}^*(e_{i,j,k})=\valustar_{\sigmae}^*(b_2)$ and $\valustar_{\sigmae}^*(F_{i,j})=\ubracket{s_{i,j}}\oplus\valustar_{\sigmae}^*(g_1)$.
	By the same arguments used in case~2~(b), this implies (\ref{equation: Valuation in Escape Vertices Phase Three}).
	Hence let $\sigmae(s_{i,j})=h_{i,j}$.
	By the choice of $j$, this implies that we need to have $\relbit{\sigmae}\neq 1$ due to \Pref{USV1}$_i$.
	But then, (\ref{equation: Valuation in Escape Vertices Phase Three}) can be shown by the same arguments used for case~2~(a).
	
	Finally, we show that $\sigmae(d_{i,j,1-k})=F_{i,j}$ and $j=\indbit_{i+1}$ imply $I_{\sigmae}=I_{\sigma}\setminus\{e\}$.
	As only the valuation of $e_{i,j,k}$ can increase, it suffices to show that $(d_{i,j,k},e_{i,j,k})\notin I_{\sigma},I_{\sigmae}$.
	As $F_{i,j}=F_{i,\indbit_{i+1}}$ is closed for $\sigma$, this implies $\indbit_i=1$ by the definition of the induced bit state.
	Thus, by \Pref{REL2}, it follows that $i\geq\nsb=\relbit{\sigma}=\relbit{\sigmae}$.
	Consequently, by \Pref{USV1}$_i$, \Pref{EV1}$_{i+1}$ and $\sigmae(d_{i,j,1-k})=F_{i,j}$, it follows that \[\valustar_{\sigmae}^*(F_{i,j})=\valustar_{\sigmae}^*(s_{i,j})=\ubracket{s_{i,j}}\oplus\valustar_{\sigmae}^*(h_{i,j})=\ubracket{s_{i,j},h_{i,j}}\oplus\valustar_{\sigmae}^*(b_{i+1}).\]
	We distinguish the following cases.
	\begin{enumerate}
		\item Let $\relbit{\sigmae}=1$.
			Then $\sigmae(e_{i,j,k})=g_1$, hence $\valustar_{\sigmae}^*(e_{i,j,k})=W_1^*\oplus\valustar_{\sigmae}^*(b_2)=W_1^*\oplus L_2^*$ by \Cref{lemma: Valuations in Phase Three}.
			In this case, $\relbit{\sigmae}=1$ implies \[\valustar_{\sigmae}^*(F_{i,j})=\ubracket{s_{i,j},h_{i,j}}\oplus L_{i+1}^*\succ\bigoplus_{i'=1}^{i}W_{i'}^*\oplus L_{i+1}^*\succeq W_1^*\oplus L_2^* =\valustar_{\sigmae}^*(e_{i,j,k})\] if $i\neq 1$ whereas $i=1$ implies \[\valustar_{\sigmae}^*(F_{i,j})=\ubracket{s_{1,j},h_{1,j}}\oplus\valustar_{\sigmae}^*(b_2)\succ W_1^*\oplus\valustar_{\sigmae}^*(b_2)=\valustar_{\sigmae}^*(e_{i,j,k}).\]
			Hence $\valustar_{\sigmae}^*(F_{i,j})\succ\valustar_{\sigmae^*}(e_{i,j,k})$ in either case, so $(d_{i,j,k},e_{i,j,k})\notin I_{\sigmae}$.
		\item Let $\relbit{\sigmae}\neq 1$.
			Then, $\sigmae(e_{i,j,k})=b_2$, hence $\valustar_{\sigmae}^*(e_{i,j,k})=\valustar_{\sigmae}^*(b_2)=L_2^*$ by \Cref{lemma: Valuations in Phase Three}.
			If $\valustar_{\sigmae}^*(b_{i+1})=L_{i+1}^*$, then the statement follows from \[\valustar_{\sigmae}^*(F_{i,j})=\ubracket{s_{i,j},h_{i,j}}\oplus L_{i+1}^*\succ\bigoplus_{\ell=2}^{i}W_{\ell}\oplus L_{i+1}^*\succ L_2^*.\]
			If $\valustar_{\sigmae}^*(b_{i+1})=R_{i+1}^*$, then $\sigmae(b_{i+1})=g_{i+1}$ and $i+1<\relbit{\sigmae}$.
			But then, by \Pref{EV1}$_{i+1}$, $\indbit_{i+1}=1$ contradicting $\relbit{\sigmae}=\nsb$. 
	\end{enumerate}
	
	Since these arguments can be applied to $\sigma$ analogously, $(d_{i,j,k},e_{i,j,k})\notin I_{\sigma},I_{\sigmae}$ and the statement follows.
\end{proof}

\PhaseThreeWellBehaved*

\begin{proof}
We first show that $\sigmae$ is a phase-3-strategy for $\bit$.
If $F_{i,j}$ is halfopen, then the application of $e$ can only influence Properties (\ref{property: ESC1}) and (\ref{property: ESC2}).
In that case, there is nothing to show as $\sigmae$ does not need to fulfill these properties.
If $F_{i,j}$ is closed for $\sigma$, then $j\neq\indbit^{\sigma}_{i+1}$ by assumption.
As $j\neq\indbit^{\sigmae}_{i+1}$, \Pref{EV1}$_i$ remains valid.
Consequently, $\sigmae$ is a phase-3-strategy for $\bit$ and in particular $\relbit{\sigmae}=\nsb=\relbit{\sigma}$.
We thus skip the upper index $\sigma$ resp. $\sigmae$ when referring to the induced bit status as $\bit+1=\indbit^{\sigma}=\indbit^{\sigmae}$.

We prove that $\sigmae$ is well-behaved.
We need to consider all properties related to escape vertices where the premise might become true or where the conclusion might become false.
\begin{itemize}[align=right, leftmargin=1.75cm]
	\item[(\ref{property: BR2})] This property only needs to be checked if the conclusion becomes false.
		Thus, assume $i<\relbit{\sigmae}$, implying $\relbit{\sigmae}>1$.
		But then, $e=(e_{i,j,k},b_2)$, hence it cannot happen that the conclusion becomes false by applying the switch $e$.
	\item[(\ref{property: D2})] By \Pref{REL2},$\sigmae(b_2)=g_2$ implies $\relbit{\sigmae}\leq 2$, hence the premise cannot be fulfilled.
	\item[(\ref{property: MNS1})] Assume $\relbit{\sigmae}=1\wedge\minsige{b}\leq\minnegsige{s},\minnegsige{g}$ and $G_n=S_n$.
		By \Pref{REL2}, it holds that $\relbit{\sigmae}=\nsb=1$.
		Since $\sigmae$ has \Pref{CC2}, this implies that the cycle center $d_{1,\sigmaebar(g_1)}$ is closed and thus $\nsigmaebar(eb_1)$.
	\item[(\ref{property: MNS2})] Assume $\relbit{\sigmae}=1$ and let $i<\minnegsige{g}<\minnegsige{s},\minsige{b}$.
		If $\sigmae(b_2)=g_2$, then $\minsige{b}=2$ by definition and there cannot be an index fulfilling the conditions of the premise.
		If $\sigmae(b_2)=b_3$, then $\indbit_2=0$ by \Pref{EV1}$_2$, implying $\sigmae(g_1)=F_{1,0}$ by \Pref{CC2}.
		But then $\minnegsige{g}=1$, hence there cannot be an index such that the conditions of the premise are fulfilled.
	\item[(\ref{property: MNS3})] This follows by the same arguments used for \Pref{MNS2}.
	\item[(\ref{property: MNS4})] Assume $\relbit{\sigmae}=1\wedge\minnegsige{s}\leq\minnegsige{g}<\minsige{b}$.
		If $\minsige{b}=2$, then $\minnegsige{s}=\minnegsige{g}=1$, implying $\sigmae(g_1)=F_{1,0}$ and $\sigmae(s_{1,0})=b_1$.
		But this contradicts \Pref{CC2} since \Pref{EV1}$_2$ implies $\sigmaebar(g_1)=\indbit_2=\sigmaebar(b_2)=1$.
		Thus assume $\minsige{b}>3$, implying $\sigmae(b_2)=b_3$.
		Then, by \Pref{CC2}, $\sigmae(g_1)=F_{1,0}$ and $\sigmae(s_{1,0})=h_{1,0}$ by \Pref{USV1}$_1$.
		But this implies $1=\minnegsige{g}<\minnegsige{s}$, contradicting the premise.
	\item[(\ref{property: MNS5})] Assume $\relbit{\sigmae}=1\wedge i<\minnegsige{s}<\minsige{b}\leq\minnegsige{g}$.
		Since no such $i$ can exist if $\minsige{b}=2$, assume $\minsige{b}>2$.
		This in particular implies $\sigmae(b_2)=b_3$, hence $\sigmae(g_1)=F_{1,0}$ by \Pref{CC2}.
		This implies $\minnegsige{g}=1$, contradicting the assumption.
	\item[(\ref{property: MNS6})] If $\minsige{b}>2$, then the same arguments used for \Pref{MNS5} can be applied again.
		Hence assume $\minsige{b}=2$, implying $\sigmae(b_2)=g_2$.
		Then, by \Pref{CC2}, $\sigmae(g_1)=F_{1,1}$ and $\sigmae(s_{1,1})=h_{1,1}$.
		But this implies $\minnegsige{s}>2$, contradicting $\minnegsige{s}<\minsige{b}$.
	\item[(\ref{property: EG1})] This property only needs to be considered if its premise is incorrect for $\sigma$ but correct for $\sigmae$.
		Therefore, since $\relbit{\sigma}=\relbit{\sigmae}=1$ implies that the switch $(e_{i,j,k},g_1)$ was applied, we need to have $\nsigmabar(eg_{i,j})\wedge\nsigmabar(eb_{i,j})$.
		This implies $\sigmabar(d_{i,j})$, hence $j\neq\indbit_{i+1}$ by assumption, implying $\sigmaebar(s_{i,j})=b_1$ by \Pref{USV1}$_i$.
	\item[(\ref{property: EG2})] If $\relbit{\sigmae}=1$, then $\sigmaebar(d_1)=\sigmaebar(d_{\relbit{\sigmae}})$ by \Pref{CC2} and the choice of $e$, so the implication is correct.
	\item[(\ref{property: EG3})] The premise of this property can only become true if $\relbit{\sigmae}=1$ as we need to apply the switch $(e_{i,j,k},g_1)$.
		Thus, \Pref{CC2} implies that $\sigmaebar(d_1), \sigmaebar(g_1)=\indbit_{i+1}$ and hence, by \Pref{USV1}$_{1}$, $\sigmaebar(s_1)$.
	\item[(\ref{property: EG4})] By \Pref{EV1}$_2$, it holds that $\sigmaebar(b_2)=\indbit_{2}$ and by \Pref{CC2} and the premise we have $\sigmaebar(g_{\relbit{\sigmae}})=\sigmaebar(g_1)=\indbit_{2}$.
	\item[(\ref{property: EG5})] The premise of this implication cannot become correct since $\sigmaebar(eg_{i,j})\wedge\neg\sigmaebar(eb_{i,j})$ imply $\relbit{\sigmae}=1$.
	\item[(\ref{property: EB1})] Assume that the conditions of the premise were fulfilled.
		Then $\sigmaebar(eb_{i,j})\wedge\neg\sigmaebar(eg_{i,j})$.
		If also $\sigmabar(eb_{i,j})\wedge\neg\sigmabar(eg_{i,j})$, the statement follows since $\sigma$ is well-behaved.
		Hence suppose that this is not the case.
		Then, $\sigmabar(d_{i,j})$, implying that $j\neq\indbit_{i+1}=\sigmaebar(b_{i+1})$ by assumption and \Pref{EV1}$_{i+1}$.
	\item[(\ref{property: EB2})] If the premise is true for $\sigma$, then there is nothing to prove.
		Hence assume that it is incorrect for $\sigma$.	
		Then, either $\sigmabar(eb_{i,0})\wedge\sigmabar(eg_{i,0})$ or $\sigmabar(d_{i,0})$.
		In the first case, \Pref{EBG1} (applied to $\sigma$) yields $\sigmabar(b_{i+1})=0$.
		This is a contradiction to \Pref{EB1} (applied to $\sigmae$) since it implies $\sigmabar(b_{i+1})=\sigmaebar(b_{i+1})\neq 0$.
		Consequently, $\sigmabar(d_{i,0})$, implying $0\neq\indbit_{i+1}$ and thus $1=\indbit_{i+1}$.		
		We now show that $i<\relbit{\sigmae}$, implying $\relbit{\sigmae}=i+1$ since $\relbit{\sigmae}=\nsb$ by \Pref{REL2}.
		For the sake of a contradiction assume $i\geq\relbit{\sigmae}$.
		Then, by \Pref{USV1}$_{i}$ and \Pref{EV1}$_{i+1}$, $\sigmaebar(b_{i+1})=\indbit_{i+1}=0$.
		However, by \Pref{EB1}, also $\sigmaebar(b_{i+1})\neq 0$ which is a contradiction.	
	\item[(\ref{property: EB3})] As before, there is nothing to prove if the premise is also correct for $\sigma$.
		By the same arguments used for \Pref{EB2}, we can deduce that assuming $\sigmabar(eb_{i,j})\wedge\sigmabar(eg_{i,j})$ yields a contradiction.
		Hence, $\sigmabar(d_{i,j})$, implying $j\neq\indbit_{i+1}^{\sigmae}$ by assumption.
		If $\relbit{\sigmae}>2$, then $\sigmae(b_2)=b_3$ follows by \Pref{EV1}$_2$.
		Hence assume $\relbit{\sigmae}=2$, implying $i\geq\relbit{\sigmae}$.
		Consequently, $\sigmae$ has \Pref{USV1}$_i$, implying $j=\indbit_{i+1}$ since $\sigmae(s_{i,j})=h_{i,j}$ by assumption.
		This is however a contradiction to the choice of $j$.
	\item[(\ref{property: EB4})] If the premise is true for $\sigma$, then there is nothing to prove.
		Hence assume that it is incorrect for $\sigma$.
		By the same arguments used earlier, we deduce $\sigmabar(d_{i,1})$, implying $1\neq\indbit_{i+1}$ and thus $0=\indbit_{i+1}$.		
		We now show that $i<\relbit{\sigmae}$, implying $\relbit{\sigmae}>i+1$ since $\relbit{\sigmae}=\nsb$ by \Pref{REL2}.
		Towards a contradiction assume $i\geq\relbit{\sigmae}$.
		Then, by \Pref{USV1}$_{i}$ and \Pref{EV1}$_{i+1}$, $\sigmaebar(b_{i+1})=\indbit_{i+1}=0$.
		However, by \Pref{EB1}, also $\sigmaebar(b_{i+1})\neq 0$ which is a contradiction.	
	\item[(\ref{property: EB6})] By \Pref{REL2}, $\relbit{\sigmae}>2$ implies $\nsb>2$.
		Hence $\indbit_{2}=0$, implying $\sigmae(b_2)=b_3$ by \Pref{EV1}$_2$.
	\item[(\ref{property: EB5})] By \Cref{lemma: Traits of Relbit} it suffices to show $\incorrect{\sigmae}\neq\emptyset$.
		Consider the case $\relbit{\sigmae}>2$ first.
		Then $\sigmae(b_2)=b_3$ by \Pref{EB6}, implying $\relbit{\sigmae}\neq \min\{i'\colon\sigmae(b_{i'})=b_{i'+1}\}$.
		Hence $\incorrect{\sigmae}\neq\emptyset$ in this case.
		Now consider the case $\relbit{\sigmae}=2$.
		Then, by assumption, $\sigmae(b_1)=g_1$ and $\sigmae(b_2)=g_2$ by \Pref{REL2} and \Pref{EV1}$_2$.
		Furthermore, by \Pref{BR1}, $\sigmaebar(g_1)=F_{1,0}$, hence $\incorrect{\sigmae}\neq\emptyset$.
	\item[(\ref{property: EBG1})] If $i\geq\relbit{\sigmae}$, then \Pref{USV1}$_i$ and $\sigmae(s_{i,j})=h_{i,j}$ imply that $j=\indbit_{i+1}$.
		Hence $j=\sigmaebar(b_{i+1})$ by \Pref{EV1}$_{i+1}$.
		Thus assume $i<\relbit{\sigmae}$, implying $\relbit{\sigmae}>1$.
		Therefore, the switch $(e_{i,j,k},b_2)$ was applied, implying $\sigmabar(eg_{i,j})\wedge\nsigmabar(eb_{i,j})$.
		But then, $\sigmaebar(b_{i+1})=\sigmabar(b_{i+1})=j$ by \Pref{EG5}. 
	\item[(\ref{property: EBG2})] By assumption and \Pref{EV2}$_2$, it follows that $\sigmaebar(g_1)=\sigmaebar(b_2)=\indbit_{2}$.
		Thus, $\sigmaebar(s_1)=\sigmaebar(s_{1,\sigmaebar(g_1)})=\sigmaebar(s_{1,\indbit_{2}})$ and hence, by either \Pref{USV2}$_{1,\indbit_2}$ or \Pref{USV1}, $\sigmaebar(s_1)=1$.
	\item[(\ref{property: EBG3})] If $\relbit{\sigmae}=1$, then $\sigmaebar(d_1)$ follows from \Pref{CC2}.
		Thus assume $\relbit{\sigmae}>1$, implying $\sigmae(b_1)=g_1$.
		Towards a contradiction, assume that the cycle center $F_{1,\indbit_{2}}=F_{1,\sigmaebar(g_1)}$ was not closed.
		Since the game is a sink game resp. weakly unichain, the cycle center cannot escape towards $g_1$ since player 1 could then create a cycle in $S_n$ resp. since $M_n$ would not have the weak unichain condition.
		Thus, by assumption, $\sigmaebar(eb_{1,\sigmaebar(g_1)})\wedge\nsigmaebar(eg_{1,\sigmaebar(g_1)})\wedge\sigmae(b_1)=g_1$.
		But then, \Pref{EB1} implies $\sigmaebar(g_1)\neq\sigmaebar(b_2)$, contradicting the assumption.
	\item[(\ref{property: EBG4})] The assumptions $\sigmaebar(b_2)=g_2$ and $\sigmaebar(g_1)=F_{1,0}$ imply $\relbit{\sigmae}=2$ if $\sigmae(b_1)=g_1$.
		If this is not the case, we have $\sigmaebar(b_1)=b_2$, implying $\relbit{\sigmae}=1$.
	\item[(\ref{property: EBG5})] By assumption $\sigmaebar(g_1)=F_{1,1}$ and $\sigmaebar(b_2)=b_3$.
		Towards a contradiction assume $\relbit{\sigmae}=2$.
		Then we need to have $\sigmae(b_1)=g_1$ and $\relbit{\sigmae}=\min\{i'\colon\sigmae(b_{i'})=b_{i'+1}\}$.
		But then $\sigmaebar(b_1)\wedge\sigmaebar(g_1)\neq\sigmaebar(b_2)$, implying $\incorrect{\sigmae}\neq\emptyset$, contradicting \Cref{lemma: Traits of Relbit}.
	\item[(\ref{property: D1})] Assume $i\neq 1$.
		Then, by \Pref{EV1}$_i$, $\sigmae(b_i)=g_i$ implies $\indbit_i=1$ and $\sigmaebar(g_i)=\indbit_{i+1}$ by \Pref{EV2}$_i$.
		Since we only open inactive cycle centers, there is nothing to show in this case.
		Hence assume $i=1$.
		If $\relbit{\sigmae}=1$, then $\sigmae(b_1)=b_2$, hence the premise is incorrect.
		Thus assume $\relbit{\sigmae}>1$.
		This implies $t^\to=b_2$.
		In particular, $\sigmaebar(eb_{i,j})\wedge\neg\sigmaebar(eg_{i,j})$, hence $\relbit{\sigmae}>2$ implies $\sigmae(b_2)=b_3$ by \Pref{EB6}.
		Thus, $\relbit{\sigmae}>2\not\Leftrightarrow\sigmae(b_2)=g_2$, so the premise is incorrect.
\end{itemize}
Hence $\sigmae$ is well-behaved. 
\end{proof}

\PGCycleVerticesPhaseThree*

\begin{proof}
By \Cref{lemma: Phase 3 Well-Behaved}, it suffices to prove $I_{\sigmae}=I_{\sigma}\setminus\{e\}$.
By construction, $d_{i,j,k}$ is the only vertex having an edge to $e_{i,j,k}$ and $F_{i,j}$ is the only vertex having an edge towards $d_{i,j,k}$.
It thus suffices to show that the valuation of $F_{i,j}$ does not change when applying $e$, so we prove $\valu_{\sigma}^{\P}(F_{i,j})=\valu_{\sigmae}^{\P}(F_{i,j}).$
Since $\sigma(d_{i,j,k})\neq e_{i,j,k}$ implies that we cannot have $\sigmabar(eg_{i,j})\wedge\sigmabar(eb_{i,j})$, it suffices to distinguishing the following cases.
\begin{itemize}
	\item Let $\sigmabar(d_{i,j})$, implying $\valustar_{\sigmae}^\P(F_{i,j})=\valustar_{\sigmae}^\P(s_{i,j})$ by \Cref{lemma: Exact Behavior Of Counterstrategy}.
		If $t^\to=b_2$, then $\sigmaebar(eb_{i,j})\wedge\neg\sigmaebar(eg_{i,j})$ and $\relbit{\sigmae}\neq 1$.
		If $t^\to=g_1$, then $\sigmaebar(eg_{i,j})\wedge\neg\sigmaebar(eb_{i,j})$ and $\relbit{\sigmae}=1$.
		This however yields $\tau^{\sigmae}(F_{i,j})=s_{i,j}$ in either case by \Cref{lemma: Exact Behavior Of Counterstrategy}, implying the statement.
	\item Let $\sigmabar(eg_{i,j})\wedge\neg\sigmabar(eb_{i,j})$.
		Since the target of $F_{i,j}$ does not change when applying~$e$ if $\sigmaebar(eg_{i,j})\wedge\neg\sigmaebar(eb_{i,j})$ we can assume $\sigmaebar(eg_{i,j})\wedge\sigmaebar(eb_{i,j})$.
		This implies $t^{\rightarrow}=b_2$, so $\relbit{\sigmae}\neq 1$, implying $\tau^{\sigmae}(F_{i,j})=g_1$ by \Cref{lemma: Exact Behavior Of Counterstrategy}.
		It now suffices to show $\sigmaebar(g_1)\neq\sigmaebar(b_2)$.
		For the sake of contradiction, assume $\sigmaebar(g_1)=\sigmaebar(b_2)$.
		Then, by \Pref{EBG3}, $F_{1,\sigmaebar(g_1)}$ is closed.
		Thus, by the definition of $\indbit$, it holds that$1=\sigmaebar(d_{1,\indbit_{2}})=\indbit_{1}$, implying $\nsb=1$.
		But this contradicts $\relbit{\sigmae}\neq 1$ by \Pref{REL2}.
	\item Let $\sigmabar(eb_{i,j})\wedge\neg\sigmaebar(eg_{i,j})$.
		By the same arguments used in the last case we can assume $\sigmaebar(eb_{i,j})\wedge\sigmaebar(eg_{i,j})$.
		Hence, $t^{\rightarrow}=g_1$ and $\relbit{\sigma}=\relbit{\sigmae}=1$.
		By \Pref{USV1}$_i$ and \Pref{EV1}$_{i+1}$, we have $\sigmabar(s_{i,j})=b_1$ if $j\neq\indbit_{i+1}$ or $\sigmabar(b_{i+1})=j$ if $j=\indbit_{i+1}$.
		In either case, $\tau^{\sigma}(F_{i,j})=b_2$ by \Cref{lemma: Exact Behavior Of Counterstrategy}.
		It hence suffices to show $\sigmaebar(g_1)=\sigmaebar(b_2)$.
		This however follows from \Pref{CC2} since $\relbit{\sigmae}=1$.\qedhere
\end{itemize}
\end{proof}

\MDPPhaseThreeLevelIsSet*

\begin{proof}
By \Cref{lemma: Phase 3 Well-Behaved}, it suffices to prove $I_{\sigmae}=I_{\sigma}\setminus\{e\}$.

By \Pref{REL2}, $\nsb=\relbit{\sigma}=\relbit{\sigmae}$, so in particular $i\geq\relbit{\sigmae}$.
By \Pref{EV2}$_i$, we have $\sigmabar(g_i)=1-j=\indbit_{i+1}$.
We thus begin by showing $(g_i,F_{i,j})\notin I_{\sigma},I_{\sigmae}$.
This is done by showing \begin{equation} \label{equation: MDP Phase 3 Level is set}
\begin{split}
\valu_{\sigma}(F_{i,1-j})&=\valu_{\sigmae}(F_{i,1-j}),\\
\valu_{\sigmae}(F_{i,1-j})&>\valu_{\sigmae}(F_{i,j})
\end{split}
\end{equation} which suffices as $\valu_{\sigmae}(F_{i,j})\geq\valu_{\sigma}(F_{i,j})$.
Since $\indbit_i=1$ implies $i\geq\nsb$, we have $\sigmabar(d_{i,1-j})$ by \Pref{EV1}$_i$ and $\sigma(s_{i,1-j})=h_{i,1-j}$ by \Pref{USV1}$_i$.
By \Pref{EV1}$_{i+1}$, this yields \[\valu_{\sigma}(F_{i,1-j})=\rew{s_{i,1-j},h_{i,1-j}}+\valu_{\sigma}(b_{i+1})=\rew{s_{i,1-j},h_{i,1-j}}+L_{i+1}.\]
This implies $\valu_{\sigma}(F_{i,1-j})=\valu_{\sigmae}(F_{i,1-j})$.
Since $F_{i,j}$ is not closed for $\sigmae$, it either escapes to $t^{\rightarrow}$ or is mixed.
Consequently, it either holds that $\valustar_{\sigmae}(F_{i,j})=\valustar_{\sigmae}(t^{\rightarrow})$ or $\valustar_{\sigmae}(F_{i,j})=\frac{1}{2}\valustar_{\sigmae}(t^{\rightarrow})+\frac{1}{2}\valustar_{\sigmae}(t^{\leftarrow})$.
As $\valustar_{\sigmae}(t^{\rightarrow})>\valustar_{\sigmae}(t^{\leftarrow})$ by \Cref{lemma: Valuations in Phase Three}, we have $\valustar_{\sigmae}(F_{i,j})\leq\valustar_{\sigmae}(t^{\rightarrow})$.
Let $t^{\rightarrow}=b_2$.
Since $\sigma(b_2)=g_2$  implies $\nsb=2$ by \Pref{EV1}$_2$, we have $\valustar_{\sigmae}(b_2)=L_2$ in any case.
Consequently, \[\valustar_{\sigmae}(F_{i,j})\leq\valustar_{\sigmae}(b_2)=L_2=L_{2,i}+L_{i+1}<\rew{s_{i,1-j},h_{i,1-j}}+L_{i+1}=\valustar_{\sigmae}(F_{i,1-j}).\]
Let $t^{\rightarrow}=g_1$.
Then, $\valustar_{\sigmae}(g_1)=W_1+\valustar_{\sigmae}(b_2)$ by \Cref{lemma: Valuations in Phase Three}.
Consequently, \begin{align*}
\valustar_{\sigmae}(F_{i,j})&\leq W_1+\valustar_{\sigmae}(b_2)=W_1+L_2=W_1+L_{2,i}+L_{i+1}\\
	&<\rew{s_{i,1-j},h_{i,1-j}}+L_{i+1}=\valustar_{\sigmae}(F_{i,1-j}).
\end{align*}
Thus, $(g_{i},F_{i,j})\notin I_{\sigma},I_{\sigmae}$ as claimed.

We now consider the cycle edges of $F_{i,j}$.
First, $(d_{i,j,k},F_{i,j})\notin I_{\sigmae}$ as $(d_{i,j,k},e_{i,j,k})$ was just applied.
If $\sigma(d_{i,j,1-k})=F_{i,j}$, then also $(d_{i,j,1-k},F_{i,j})\notin I_{\sigma},I_{\sigmae}$, implying $I_{\sigmae}=I_{\sigma}\setminus\{e\}$ since the valuation of no further vertex changes due to $\sigma(g_i)=F_{i,1-j}$.
Hence assume $\sigma(d_{i,j,1-k})=e_{i,j,1-k}$.
We prove $(d_{i,j,1-k},F_{i,j})\in I_{\sigma}\Leftrightarrow(d_{i,j,1-k},F_{i,j})\in I_{\sigmae}$ which suffices to prove the statement since no other vertex but $F_{i,j}$ has an edge to $d_{i,j,1-k}$.

Let $v=\sigma(e_{i,j,1-k})=\sigmae(e_{i,j,1-k})$.
We prove that $i\geq\relbit{\sigmae}$ and $\sigma(g_i)=\sigmae(g_i)=F_{i,1-j}$ imply $\valu_{\sigma}(v)=\valu_{\sigmae}(v)$.
If $v=b_2$, then $\valu_{\sigma}(v)=\valu_{\sigmae}(v)=L_2$ by \Cref{lemma: Valuations in Phase Three}.
If $v=g_1$, then either \[\valustar_{\sigma}(v)=\valustar_{\sigmae}(v)=\rew{g_{\ell}}+\sum_{i'\in[\ell-1]}W_{i'}+\valustar_{\sigma}(b_2)\] where $\ell=\min\{i'\geq 1\colon\nsigmabar(d_{i'})\}<\relbit{\sigma}$ or $\valustar_{\sigma}(v)=\valustar_{\sigmae}(v)=R_1$.
We furthermore observe that $\sigmabar(g_i)=\sigmaebar(g_i)\neq j$ implies $\min\{i'\geq 1\colon\nsigmabar(d_{i'})\}=\min\{i'\geq1\colon\nsigmaebar(d_{i'})\}$.
Thus, if we have $\valu_{\sigma}(F_{i,j})>\valu_{\sigma}(v)$, then also $\valu_{\sigmae}(F_{i,j})>\valu_{\sigmae}(v)$, so \[(d_{i,j,1-k},F_{i,j})\in I_{\sigma}\Rightarrow(d_{i,j,1-k},F_{i,j})\in I_{\sigmae}.\]
Hence consider the case $(d_{i,j,1-k},F_{i,j})\notin I_{\sigma}$, implying $\valu_{\sigma}(F_{i,j})\leq\valu_{\sigma}(v)$.
Since $\sigma(d_{i,j,k})=F_{i,j}$ and $\sigma(d_{i,j,1-k})=e_{i,j,1-k}$, we have \begin{align*}
\valu_{\sigma}(F_{i,j})-\valu_{\sigma}(v)=\frac{2\e}{1+\e}[\rew{s_{i,j}}+\valu_{\sigma}(b_1)-\valu_{\sigma}(v)]\leq 0.
\end{align*}
Thus, $\valu_{\sigma}(b_1)+\rew{s_{i,j}}\leq \valu_{\sigma}(v)$, hence $\valu_{\sigma}(b_1)<\valu_{\sigma}(v)$.
Since $\valu_{\sigma}(t^{\leftarrow})<\valu_{\sigma}(t^{\rightarrow})$ by \Cref{lemma: Valuations in Phase Three}, this implies that $v=t^{\rightarrow}$ and $\sigma(b_1)=t^{\leftarrow}$ have to hold.
As it holds that $\valu_{\sigma}(t^{\leftarrow})=\valu_{\sigmae}(t^{\leftarrow})$, this then implies \begin{align*}
\valu_{\sigmae}(F_{i,j})-\valu_{\sigmae}(v)&=(1-\e)\valu_{\sigmae}(t^{\rightarrow})+\e\valu_{\sigmae}(s_{i,j})-\valu_{\sigmae}(t^{\rightarrow})\\
	&=\e[\rew{s_{i,j}}+\valu_{\sigmae}(b_1)-\valu_{\sigmae}(t^{\rightarrow})]\leq 0,
\end{align*}
hence $(d_{i,j,1-k},F_{i,j})\notin I_{\sigmae}$.
\end{proof}

\ImprovingSwitchinOtherCCPhaseThree*

\begin{proof}
Since $F_{i,j}$ is $t^{\leftarrow}$-halfopen, the choice of $e$ implies $\sigma(d_{i,j,1-k})=e_{i,j,1-k}$.
Consequently, by \Cref{lemma: Phase 3 Well-Behaved}, it suffices to prove $I_{\sigmae}=I_{\sigma}\setminus\{e\}$.

Since $\sigma(g_i)=F_{i,1-j}$, the application of $e$ can only increase the valuation of $F_{i,j}, d_{i,j,0}$ and $d_{i,j,1}$.
In addition, and since there are no player 0 vertices $v$ with $(v,d_{i,j,*})\in E_0$.
It thus suffices to prove $(g_i,F_{i,j})\notin I_{\sigma}, I_{\sigmae}$.
This however follows from \Cref{lemma: Valuations in Phase Three} since \[\valustar_{\sigma}(F_{i,1-j})=\valustar_{\sigma}(t^{\rightarrow})>\valustar_{\sigma}(t^{\leftarrow})=\valustar_{\sigma}(F_{i,j})\]and \[\valustar_{\sigmae}(F_{i,1-j})=\valustar_{\sigmae}(t^{\rightarrow})>\frac{1}{2}\valustar_{\sigmae}(t^{\rightarrow})+\frac{1}{2}\valustar_{\sigmae}(t^{\leftarrow})=\valustar_{\sigmae}(F_{i,j}).\]
\end{proof}

\NoCCClosedPhaseThree*

\begin{proof}
In both cases, $e=(d_{i,j,k},e_{i,j,k})$ is applied within a $t^{\leftarrow}$-halfopen cycle center.
This implies $\sigma(d_{i,j,1-k})=e_{i,j,1-k}$, hence, by \Cref{lemma: Phase 3 Well-Behaved}, it suffices to prove $I_{\sigmae}=I_{\sigma}\setminus\{e\}$.

Let both cycle centers be $t^{\leftarrow}$-halfopen for $\sigma$ and let $j\coloneqq\sigmabar(g_i)=\indbit_{i+1}$.
By \Cref{lemma: Both CC Open For MDP}, $\valu_{\sigma}(F_{i,j})>\valu_{\sigma}(F_{i,1-j})$ and by \Cref{lemma: Valuations in Phase Three}, also $\valu_{\sigmae}(F_{i,j})>\valustar_{\sigmae}(F_{i,1-j})$.
Thus, $(g_i,F_{i,1-j})\notin I_{\sigma},I_{\sigmae}$.
Note that \Cref{lemma: Both CC Open For MDP} can be applied since $i\geq\relbit{\sigma}+1=\nsb+1$ by \Pref{REL2} and since it has \Pref{USV1}$_i$ and \Pref{EV1}$_{i+1}$.
Due to the application of the switch $e$, the valuation of $g_i$ increases. 
We prove that this does not create new improving switches.
We thus first prove $\sigma(b_i)\neq g_i$ and $(b_i,g_i)\notin I_{\sigma}, I_{\sigmae}$.
Since $i\geq\relbit{\sigma}+1$, no cycle centers being closed implies $\indbit_i=0$, hence $\sigma(b_i)=\sigmae(b_i)=b_{i+1}$ by \Pref{EV1}$_i$.
Furthermore, if $\relbit{\sigma}>1$ and thus $t^{\rightarrow}=b_2$, then \begin{align*}
	\valustar_{\sigmae}(b_{i+1})&=L_{i+1}>\rew{g_i}+\sum_{\ell=1}^{i-1}W_{\ell}+L_{i+1}\geq\rew{g_i}+L_2=\rew{g_i}+\valustar_{\sigmae}(b_2)=\valustar_{\sigmae}(g_i)
\end{align*}
and, by \Cref{lemma: Valuations in Phase Three}, $\valustar_{\sigma}(b_{i+1})>\rew{g_i}+\valustar_{\sigma}(g_1)=\valustar_{\sigma}(g_i)$ follows by the same estimation.
Consequently, $(b_i,g_i)\notin I_{\sigma}, I_{\sigmae}$ if $\relbit{\sigma}>1$.
The statement follows by a similar argument if $\relbit{\sigma}=1$ by $\valustar_{\sigma}(g_1)=W_1+\valustar_{\sigma}(b_2)$.

We now prove that $\sigma(s_{i-1,1})=b_1$ and $(s_{i-1,1},h_{i-1,1})\notin I_{\sigma}, I_{\sigmae}$.
It cannot happen that $i=1$ due to $i\geq\relbit{\sigma}+1$, hence we do not need to consider a possible increase of the valuation of $g_1$.
Since $i\geq\relbit{\sigma}+1$ implies $i-1\geq\relbit{\sigma}$ and since $\indbit_{i+1}=0$, \Pref{USV1}$_{i-1}$ implies $\sigma(s_{i-1,1})=\sigmae(s_{i-1,1})=b_1$.
It remains to prove $\valustar_{\sigma}(b_1)>\valustar_{\sigma}(h_{i-1,1})$ and $\valustar_{\sigmae}(b_1)>\valustar_{\sigmae}(h_{i-1,1})$.
We only prove the second statement since it implies the first statement due to $\valustar_{\sigmae}(b_1)=\valustar_{\sigma}(b_1)$ and $\valustar_{\sigmae}(h_{i-1,1})>\valustar_{\sigma}(h_{i-1,1})$.
If $\relbit{\sigmae}=1$, then $\sigmae(b_1)=b_2, t^{\rightarrow}=g_1$ and $\valustar_{\sigmae}(g_1)=W_1+\valustar_{\sigmae}(b_2)$.
Consequently, \begin{align*}
	\valustar_{\sigmae}(b_1)&=\valustar_{\sigmae}(b_2)>\rew{h_{i-1,1},g_i}+W_1+\valustar_{\sigmae}(b_2)\\
		&=\rew{h_{i-1,1},g_i}+\valustar_{\sigmae}(g_1)=\valustar_{\sigmae}(h_{i-1,1}).
\end{align*}
If $\relbit{\sigmae}>1$, then $\sigmae(b_1)=g_1$ and $t^{\rightarrow}=b_2$.
The statement can then be shown by similar arguments as $i>\relbit{\sigma}$ implies $\rew{g_i}<\sum_{\ell\in[i-1]}W_{\ell}$.

This concludes the case that both cycle centers of level $i$ are $t^{\leftarrow}$-open.
Consider the case that $F_{i,\indbit_{i+1}}$ is mixed and that $F_{i,1-\indbit_{i+1}}$ is $t^{\leftarrow}$-halfopen for $\sigma.$
Then,  $j=\indbit_{i+1}=1-j$ implies that no other edge but $(g_i,F_{i,j})$ can become improving.
However, after the application of $e$, both cycle centers are mixed.
Hence, $\valu_{\sigmae}(F_{i,\indbit_{i+1}})>\valu_{\sigmae}(F_{i,1-\indbit_{i+1}})$ by \Cref{lemma: Both CC Open For MDP}.
\end{proof}

\MDPPhaseThreeOpenClosedCycleCenter*

\begin{proof}
By \Cref{lemma: Phase 3 Well-Behaved}, it suffices to prove the statements related to the set of improving switches.
We first prove $(g_i,F_{i,*})\notin I_{\sigma},I_{\sigmae}$.
Let $\indbit_{i}=0$ first.
Then, $\sigma(g_i)=F_{i,j}$ by assumption~1., implying $(g_i,F_{i,j})\notin I_{\sigma},I_{\sigmae}$.
We thus prove $(g_i,F_{i,1-j})\notin I_{\sigma},I_{\sigmae}$.
By assumption~1., $\sigmabar(d_{i,j})$, hence $\valustar_{\sigma}(F_{i,j})=\valustar_{\sigma}(s_{i,j})$ and $\valustar_{\sigmae}(F_{i,j})=\valustar_{\sigmae}(t^{\rightarrow})$ by \Cref{lemma: Exact Behavior Of Random Vertex}.
Since $F_{i,1-j}$ is $t^{\leftarrow}$-halfopen by assumption,  $\valustar_{\sigmae}(F_{i,1-j})<\valustar_{\sigmae}(F_{i,j})$ since $\valu_{\sigmae}(t^{\leftarrow})<\valu_{\sigmae}(t^{\rightarrow})$ by \Cref{lemma: Valuations in Phase Three}.
Consequently, $(g_i,F_{i,1-j})\notin I_{\sigmae}$ and it remains to prove $\valustar_{\sigma}(F_{i,1-j})<\valustar_{\sigma}(F_{i,j})$.

If $\sigma(s_{i,j})=b_1$, then the equivalences $\sigma(b_1)=b_2\Leftrightarrow\relbit{\sigma}=1\Leftrightarrow t^{\leftarrow}=b_2$ implies  \[\valustar_{\sigma}(F_{i,j})=\rew{s_{i,j}}+\valustar_{\sigma}(b_1)=\rew{s_{i,j}}+\valustar_{\sigma}(t^{\leftarrow})>\valu_{\sigma}(t^{\leftarrow})=\valu_{\sigma}(F_{i,1-j}).\]
Hence assume $\sigma(s_{i,j})=h_{i,j}$, implying $i<\relbit{\sigma}$ by \Pref{USV1}$_i$ as $j=1-\indbit_{i+1}$. 
This implies $\relbit{\sigma}>1$, so $t^{\leftarrow}=g_1$ and $t^{\rightarrow}=b_2$.
Let $i=\relbit{\sigma}-1$.
Then $j=1-\indbit_{\nsb}=0$, hence $\valustar_{\sigma}(F_{i,j})=\rew{s_{i,j},h_{i,j}}+\valustar_{\sigma}(b_{\nsb+1})$.
By assumption~2, $\sigmabar(d_{i'})$ for all $i'<i$.
Consequently, \begin{align*}
\valustar_{\sigma}(F_{i,1-j})&=\valustar_{\sigma}(g_1)=R_1<\rew{s_{\nsb-1,j},h_{\nsb-1,j}}+ L_{\nsb+1}=\valustar_{\sigma}(F_{i,j}).
\end{align*}

Let $i<\relbit{\sigma}-1$, implying $j=1$ and thus $\valustar_{\sigma}(F_{i,j})=\rew{s_{i,j},h_{i,j}}+\valustar_{\sigma}(g_{i+1})$.
Since $\sigmabar(eb_{i+1})$ by assumption~2. and $\nsigmabar(eg_{i+1})$ by \Pref{BR2} and thus in particular $\nsigmabar(d_{i+1})$, we have $\valustar_{\sigma}(g_{i+1})=\rew{g_{i+1}}+\valustar_{\sigma}(b_2)$ by \Cref{corollary: Complete Valuation Of Selection Vertices MDP}.
Furthermore, as $\sigmabar(d_{i'})$ for all $i'<i$, it holds that $\valustar_{\sigma}(g_1)=\sum_{\ell<i}W_{\ell}+\rew{g_{i+1}}+\valustar_{\sigma}(b_2)$.
Consequently, as $\valustar_{\sigma}(g_1)=\valustar_{\sigma}(F_{i,1-j}),$ it follows that\begin{align*}
	\valustar_{\sigma}(F_{i,j})\hspace*{-1pt}=\hspace*{-1pt} \rew{s_{i,j},h_{i,j},g_{i+1}}+\valustar_{\sigma}(b_2)\hspace*{-1pt}>\hspace*{-1pt}\sum_{\ell\in[i]}W_{\ell}^\M+\rew{g_{i+1}}+\valustar_{\sigma}(b_2)\hspace*{-1pt}=\hspace*{-1pt}\valustar_{\sigma}(F_{i,1-j}).
\end{align*}
Hence, $\valustar_{\sigma}(F_{i,1-j})<\valustar_{\sigma}(F_{i,j})$ holds in any case, so $(g_i,F_{i,1-j})\notin I_{\sigma},I_{\sigmae}$. 
As also $(g_i,F_{i,j})\notin I_{\sigma},I_{\sigmae}$, this proves that $\indbit_{i}=0$ implies $(g_i,F_{i,*})\notin I_{\sigma},I_{\sigmae}$.

Now let $\indbit_i=1$.
Then, by \Pref{REL2}, $i\geq\relbit{\sigma}=\relbit{\sigmae}=\nsb$.
By \Pref{CC2}, $\sigma(g_i)=F_{i,\indbit_{i+1}}=F_{i,1-j}$.
We hence prove $(g_i,F_{i,j})\notin I_{\sigma},I_{\sigmae}$.
Since $\indbit_i=1$, the cycle center $F_{i,1-j}$ is closed.
By \Pref{USV1}$_i$, \Pref{EV1}$_{i+1}$ and since $i\geq\relbit{\sigma},\relbit{\sigmae}$ we thus obtain $\valustar_{\sigma}(F_{i,1-j})=\rew{s_{i,1-j},h_{i,1-j}}+\valustar_{\sigma}(b_{i+1})=\rew{s_{i,1-j},h_{i,1-j}}+L_{i+1}$ and the corresponding equality for $\valustar_{\sigmae}(F_{i,1-j})$.

As $F_{i,j}$ is closed with respect to $\sigma$ and escapes towards $t^{\rightarrow}$ with respect to $\sigmae$, \Pref{USV1}$_i$ yields $\valustar_{\sigma}(F_{i,j})=\rew{s_{i,j}}+\valustar_{\sigma}(t^{\leftarrow})$ and $\valustar_{\sigmae}(F_{i,j})=\valustar_{\sigmae}(t^{\rightarrow})$.
It is now easy to see that $\rew{s_{i,1-j},h_{i,1-j}}>\sum_{\ell=\in[i]}W_{\ell}$ implies $\valustar_{\sigma}(F_{i,1-j})>\valustar_{\sigma}(F_{i,j})$ and $\valustar_{\sigmae}(F_{i,1-j})>\valustar_{\sigmae}(F_{i,j})$.
Therefore, $(g_i,F_{i,j})\notin I_{\sigma},I_{\sigmae}$ if $\indbit_i=1$.
As also $(g_i,F_{i,1-j})\notin I_{\sigma},I_{\sigmae}$ in this case due to $\sigma(g_i)=F_{i,1-j}$, this proves $(g_i,*)\notin I_{\sigma},I_{\sigmae}$ in any case.
By the choice of $e$, we also have $(d_{i,j,k},F_{i,j})\notin I_{\sigma},I_{\sigmae}$ and $\sigma(d_{i,j,1-k})=F_{i,j}$ as we assume $F_{i,j}$ to be closed with respect to $\sigma$.
Thus also $(d_{i,j,1-k},F_{i,j})\notin I_{\sigma},I_{\sigmae}$.

If $\indbit_i=1$, the increase of the valuation of $F_{i,j}$ can only have an immediate effect on the vertices $g_i, d_{i,j,0}$ and $d_{i,j,1}$.
However, as $\indbit_i=1$ implies $\sigma(g_i)=F_{i,1-j}$ and since there are no player 0 vertices edges towards $d_{i,j,*}$, we immediately obtain $I_{\sigmae}=I_{\sigma}\setminus\{e\}$.
We thus only consider the case $\indbit_i=0$ for the remainder of this proof, implying $\sigma(g_i)=F_{i,j}$.

Since $\sigma(g_i)=F_{i,j}$, the valuation of $g_i$ increases due to the increase of the valuation of $F_{i,j}$.
We investigate how this increase influences the set of improving switches.
We first prove that $i\neq 1$ implies $\sigma(b_i)=b_{i+1}$ and $(b_i,g_i)\notin I_{\sigma},I_{\sigmae}$.
If $i=1$, then $\relbit{\sigma}>1$ as $\indbit_i=0$ by assumption, implying $\sigma(b_1)=g_1$ and thus $(b_i,g_i)\notin I_{\sigma},I_{\sigmae}$.
Hence let $i\neq 1$.
Then, $\sigma(b_i)=b_{i+1}$ by \Pref{EV1}$_i$.
If $\sigma(b_{i+1})=g_{i+1}$, then $\indbit_{i+1}=1$ and thus $i+1\geq\relbit{\sigma}=\nsb$, implying $\valustar_{\sigma}(b_{i+1})=L_{i+1}$ in any case.
The same holds for $\sigmae$, so in particular $\valustar_{\sigma}(b_{i+1})=\valustar_{\sigmae}(b_{i+1})$.
It hence suffices to prove $\valustar_{\sigmae}(b_{i+1})>\valustar_{\sigmae}(g_i)$ as $\valustar_{\sigmae}(g_i)\geq\valustar_{\sigma}(g_i)$.

By the choice of $e$ and our assumptions, $\valustar_{\sigmae}(g_i)=\rew{g_i}+\valustar_{\sigmae}(t^{\rightarrow})$.
If $t^{\rightarrow}=b_2$, then $\valustar_{\sigmae}(t^{\rightarrow})=\valustar_{\sigmae}(b_2)$.
If $t^{\rightarrow}=g_1$, then $\valustar_{\sigmae}(t^{\rightarrow})=\valustar_{\sigmae}(g_1)=W_1+\valustar_{\sigmae}(b_2)$ by \Cref{lemma: Valuations in Phase Three} as $\relbit{\sigmae}=1$.
This in particular yields \begin{align*}
	\valustar_{\sigmae}(g_i)&=\rew{g_i}+\valustar_{\sigmae}(t^{\rightarrow})\leq\rew{g_i}+W_1+\valustar_{\sigmae}(b_2)=\rew{g_i}+W_1+L_2\\
		&=\rew{g_i}+W_1+L_{2,i-1}+L_{i+1}<L_{i+1}=\valustar_{\sigmae}(b_{i+1})
\end{align*}
as $\sigma(b_i)=b_{i+1}$.
Thus $(b_i,g_i)\notin I_{\sigma},I_{\sigmae}$ for all $i$.

Now let $i>\relbit{\sigmae}$.
We prove that this implies $\sigma(s_{i-1,1})=b_1, \valu_{\sigma}(b_1)>\valu_{\sigma}(h_{i-1,1})$ and $\valu_{\sigmae}(b_1)>\valu_{\sigmae}(h_{i-1,1})$.
When proving these statements, we will also prove that $\valu_{\sigma}(b_1)=\valu_{\sigmae}(b_1), \valu_{\sigma}(g_1)=\valu_{\sigmae}(g_1)$ and $\valu_{\sigma}(b_2)=\valu_{\sigmae}(b_2)$.
We then argue why this suffices to prove the statement for $\relbit{\sigmae}=1$ and then consider the case $\relbit{\sigmae}>1$ and $i<\relbit{\sigmae}$.
It is not necessary to consider the case $i=\relbit{\sigmae}$ as $\indbit_i=0$.

Since $\indbit_{i}=0$ and $i>\relbit{\sigmae}$ implies $i-1\geq\relbit{\sigmae}$, \Pref{USV1}$_{i-1}$ implies $\sigma(s_{i-1,1})=\sigmae(s_{i-1,1})=b_1$.
This implies that the valuation of no further vertex than $h_{i-1,1}$ and the vertices discussed previously can change when transitioning from $\sigma$ to $\sigmae$.
None of these vertices are part of the valuation of $b_1,g_1$ and $b_2$ since $i>\relbit{\sigmae}, \sigma(b_i)=\sigmae(b_i)=b_{i+1}$ and $\sigma(s_{i-1,1})=\sigma(s_{i-1,1})=b_1$, implying that their valuations do not change.
In particular, $\valu_{\sigma}(b_1)=\valu_{\sigmae}(b_1)$.
If we can show $(s_{i-1,1},b_1)\notin I_{\sigma},I_{\sigmae}$, this thus proves $I_{\sigmae}=I_{\sigma}\setminus\{e\}$ for the case $\relbit{\sigmae}=1$.
Since $\valu_{\sigma}(b_1)=\valu_{\sigmae}(b_1)$, it suffices to prove $\valu_{\sigmae}(b_1)>\valu_{\sigmae}(h_{i-1,1})$ as $\valu_{\sigmae}(h_{i-1,1})\geq\valu_{\sigma}(h_{i-1,1})$.
Consider the case $\relbit{\sigmae}=1$ first, implying $t^{\rightarrow}=g_1$ and $\valustar_{\sigmae}(g_1)=W_1+\valustar_{\sigmae}(b_2)$ by \Cref{lemma: Valuations in Phase Three}.
Thus, since $\valustar_{\sigmae}(g_i)=\rew{g_i}+\valustar_{\sigmae}(t^{\rightarrow})$ and $\valustar_{\sigmae}(b_2)=\valu_{\sigmae}(b_1)$ since $\sigmae(b_1)=b_2$ due to $\relbit{\sigmae}=1$, it follows that
\begin{align*}
\valustar_{\sigmae}(h_{i-1,1})&=\rew{h_{i-1,1},g_i}+\valustar_{\sigmae}(g_i)=\rew{h_{i-1,1},g_i}+W_1+\valustar_{\sigmae}(b_2)<\valustar_{\sigmae}(b_1).
\end{align*}
Consider the case $\relbit{\sigmae}>1$ next.
Then $t^{\rightarrow}=b_2$ and $\sigmae(b_1)=g_1$.
Now, either $\valustar_{\sigmae}(b_1)=R_1$ or $\valustar_{\sigmae}(b_1)=g_{i'}+\sum_{\ell<i'}W_{\ell}+\valustar_{\sigmae}(b_2)$ where $i'=\min\{\ell\geq 1\colon\nsigmaebar(d_{i'})\}<\relbit{\sigmae}$.
In the first case, $i>\relbit{\sigmae}$ implies \begin{align*}
	\valustar_{\sigmae}(h_{i-1,1})&=\rew{h_{i-1,1},g_i}+\valustar_{\sigmae}(b_2)<L_{\relbit{\sigmae}+1}<R_1=\valustar_{\sigmae}(b_1).
\end{align*}
In the second case, $i>\relbit{\sigmae}>i'$ implies \[\valustar_{\sigmae}(h_{i-1,1})=\rew{h_{i-1,1},g_i}+\valustar_{\sigmae}(b_2)<\rew{g_{i'}}+\sum_{\ell<i'}W_{\ell}+\valustar_{\sigmae}(b_2)=\valustar_{\sigmae}(b_1).\]
Hence $i>\relbit{\sigmae}$ implies $(s_{i-1,1},h_{i-1,1})\notin I_{\sigma},I_{\sigmae},$ proving the statement for $\relbit{\sigmae}=1$.

It remains to investigate the case $i<\relbit{\sigmae}$, implying $\relbit{\sigmae}>1$ and $t^{\rightarrow}=b_2, t^{\leftarrow}=g_1$.
In this case, opening the cycle center $F_{i,j}$ changes the valuation of $g_1$.
Since $\relbit{\sigmae}>1$ implies $\sigma(b_1)=g_1$, this also changes the valuation of $b_1$ and of possibly every vertex that has an edge to either one of these vertices.
These are in particular upper selection vertices, escape vertices and cycle centers.
We begin by observing that \begin{align*}
	\valu_{\sigma}(F_{i,j})&=\valu_{\sigma}(s_{i,j})=\rew{s_{i,j},h_{i,j}}+\begin{cases}\valu_{\sigma}(b_{i+2}), &j=0\\\valu_{\sigma}(g_{i+1}), &j=1\end{cases}\\
		&=\rew{s_{i,j},h_{i,j}}+\begin{cases}\valu_{\sigma}(b_{i+2}), &i=\relbit{\sigmae}-1\\\valu_{\sigma}(g_{i+1}), &i<\relbit{\sigmae}-1\end{cases}\\
	\valu_{\sigmae}(F_{i,j})&=\frac{1-\e}{1+\e}\valu_{\sigmae}(b_2)+\frac{2\e}{1+\e}\valu_{\sigmae}(s_{i,j})\\
	\valustar_{\sigma}(g_1)&=\begin{cases}R_1, &i=\relbit{\sigmae}-1\\\rew{g_{i+1}}+\sum_{\ell<i+1}W_{\ell}+\valu_{\sigma}(b_2), &i<\relbit{\sigmae}-1\end{cases},\\
	\valustar_{\sigmae}(g_1)&=\rew{g_i}+\sum_{\ell<i}W_{\ell}+\valu_{\sigmae}(b_2).
\end{align*}
Furthermore, $\valu_{\sigma}(b_2)=\valu_{\sigmae}(b_2)=L_2$ and $i\neq 1$ implies $\sigma(b_i)=\sigmae(b_i)=b_{i+1}$.
Note that we have $\valustar_{\sigma}(g_1)<\valustar_{\sigma}(b_2)$ and $\valustar_{\sigmae}(g_1)<\valustar_{\sigmae}(b_2)$ by \Cref{lemma: Valuations in Phase Three}.
We begin by investigating upper selection vertices and prove that $(s_{i,j},b_1)$ is improving for $\sigmae$.
Since $\sigma(s_{i,j})=\sigmae(s_{i,j})=h_{i,j}$ by assumption, it suffices to prove $\valu_{\sigmae}(h_{i,j})<\valu_{\sigmae}(b_1)$.
Consider the case $i=\relbit{\sigmae}-1$ first, implying $j=0$.
Then, since \Pref{EV1}$_{i'}$ implies $\sigma(b_{i'})=b_{i'+1}$ for all $i'\in\{2,\dots,\relbit{\sigmae}-1\}$, \begin{align*}
\valustar_{\sigmae}(h_{i,j})&=\rew{h_{i,j}}+\valu_{\sigmae}(b_{i+2})<\rew{g_i}+\sum_{\ell<i}W_{\ell}+W_{i+1}+L_{i+2}\\
	&=\rew{g_i}+\sum_{\ell<i}W_{\ell}+L_{i+1}=\rew{g_i}+\sum_{\ell<i}W_{\ell}+L_2=\valustar_{\sigmae}(g_1)=\valustar_{\sigmae}(b_1).
\end{align*}
Therefore, $(s_{i,j},b_1)\in I_{\sigmae}$ if $i=\relbit{\sigmae}-1$.
Consider the case $i<\relbit{\sigmae}-1$, implying $j=1$.
Then, since $\sigmabar(eb_{i+1})\wedge\nsigmabar(eg_{i+1})$ by assumption~2 and \Pref{BR2}, \begin{align*}
	\valustar_{\sigmae}(h_{i,j})&=\rew{h_{i,j}}+\valustar_{\sigmae}(g_{i+1})=\rew{h_{i,j},g_{i+1}}+\valustar_{\sigmae}(b_2)\\
		&<\rew{g_i}+\sum_{\ell<i}W_{\ell}+\valustar_{\sigmae}(b_2)=\valustar_{\sigmae}(g_1),
\end{align*}
hence $(s_{i,j},b_1)\in I_{\sigmae}$ if $i<\relbit{\sigmae}$.
It remains to prove that no further improving switch is created.

First, we prove that for all $i'\in[n]$ and $j'\in\{0,1\}$ with $(i',j')\neq(i,j)$, $\sigma(s_{i',j'})=h_{i',j'}$ implies $(s_{i',j'},b_1)\notin I_{\sigma},I_{\sigmae}$.
Hence let $i',j'$ be such a pair of indices and $i'\geq\relbit{\sigmae}$ first.
Then, $\valustar_{\sigma}(h_{i',j'})=\rew{h_{i',j'}}+\valustar_{\sigma}(b_{i'+1})=\rew{h_{i',j'}}+L_{i'+1}$ follows from \Pref{USV1}$_{i'}$ and \Pref{EV1}$_{i'+1}$ and $\valustar_{\sigmae}(h_{i',j'})=\rew{h_{i',j'}}+L_{i'+1}$ follows analogously.
This implies $\valustar_{\sigma}(h_{i',j'}),\valustar_{\sigmae}(h_{i',j'})>\valustar_{\sigma}(b_2)=\valustar_{\sigmae}(b_2)$.
Since $\valustar_{\sigma}(g_1)<\valustar_{\sigma}(b_2)$ as well as $\valustar_{\sigmae}(g_1)<\valustar_{\sigmae}(b_2)$ by \Cref{lemma: Valuations in Phase Three}, $\sigma(s_{i',j'})=h_{i',j'}$ thus implies $(s_{i',j'},b_1)\notin I_{\sigma},I_{\sigmae}$ for $i'\geq\relbit{\sigmae}$.
Next let $i'<\relbit{\sigmae}$ and $i<i'<\relbit{\sigmae}$.
Then, by assumption~3, $j'=\indbit_{i'+1}$, so $\valu_{\sigmae}(h_{i',j'})=\rew{h_{i',j'}}+\valu_{\sigmae}(b_{i'+1})=\rew{h_{i',j'}}+L_{i'+1}$ by \Pref{EV1}$_{i'+1}$.
Furthermore, \[\valustar_{\sigmae}(b_1)=\valustar_{\sigmae}(g_1)=\rew{g_i}+\sum_{\ell<i}W_{\ell}+\valustar_{\sigmae}(b_2)=\rew{g_i}+\sum_{\ell<i}W_{\ell}+L_2.\]
Since $\indbit_{2}=\dots=\indbit_{i'}=0$ due to $i'<\relbit{\sigmae}=\nsb$ and $\rew{g_i}+\sum_{\ell<i}W_{\ell}<0$, this implies \[\valustar_{\sigmae}(b_1)<L_2=L_{i'+1}<\rew{h_{i',j'}}+L_{i'+1}=\valustar_{\sigmae}(h_{i',j'})\] and thus $(s_{i',j'},b_1)\notin I_{\sigmae}$.
The same calculations also yield $(s_{i',j'},b_1)\notin I_{\sigma}$.

Now let $i'<i<\relbit{\sigmae}$.
If $j'=\indbit_{i'+1}$, then the same arguments used previously can be applied again.
Hence let $j'=1-\indbit_{i'+1}$.
Since $i'<i<\relbit{\sigmae}$ implies $i'<\relbit{\sigmae}-1$, we have $\indbit_{i'+1}=0$ and thus $j'=1$.
By assumption~2, the cycle center $F_{\ell,\sigmaebar(g_{\ell})}$ is closed and $\sigmaebar(s_{\ell})$ for all $\ell<i$.
Since $\ell<i<\relbit{\sigmae}$ furthermore implies $\sigmae(g_{\ell})=F_{\ell,1}$ by \Pref{BR1}, we have $\lambda_{i'+1}=i$.
Consequently, as the last case of \Cref{corollary: Complete Valuation Of Selection Vertices MDP} is fulfilled, \begin{align*}
	\valustar_{\sigmae}(h_{i',j'})&=\rew{h_{i',j'}}+\valustar_{\sigmae}(g_{i'+1})=\rew{h_{i',j'}}+\rew{g_i}+\sum_{\ell=i'+1}^{i-1}W_{\ell}+\valustar_{\sigmae}(b_2)\\
		&>\sum_{\ell=1}^{i'}W_{\ell}+\rew{g_i}+\sum_{\ell=i'+1}^{i-1}W_{\ell}+\valustar_{\sigmae}(b_2)\\
		&=\rew{g_i}+\sum_{\ell<i}W_{\ell}+\valustar_{\sigmae}(b_2)=\valustar_{\sigmae}(b_1),
\end{align*}
implying $(s_{i',j'},b_1)\notin I_{\sigmae}$.
The same arguments can also be used to show $(s_{i',j'},b_1)\notin I_{\sigma}$.
If $i'=i<\relbit{\sigmae}$, then $j'=\indbit_{i'+1}$ as we consider indices $(i',j')\neq(i,j)$, implying the statement by the same arguments as before.

We next investigate escape vertices $e_{i',j',k'}$.
If $\sigma(e_{i',j',k'})=\sigmae(e_{i',j',k'})=b_2$, then \Cref{lemma: Valuations in Phase Three} implies $\valustar_{\sigma}(b_2)>\valustar_{\sigma}(g_1)$, hence $(e_{i',j',k'},g_1)\notin I_{\sigma}$.
Analogously, $(e_{i',j',k'},g_1) \notin ,I_{\sigmae}$ and, by definition, $(e_{i',j',k'},b_2)\notin I_{\sigma},I_{\sigmae}$.
Using the same arguments, it follows that $\sigma(e_{i',j',k'})=g_1$ implies  $(e_{i',j',k'},b_2)\in I_{\sigma},I_{\sigmae}$ as well as $(e_{i',j',k'},g_1)\notin I_{\sigma}, I_{\sigmae}$.
In particular, we have $(e_{i',j',k'},*)\in I_{\sigma}\Leftrightarrow (e_{i',j',k'},*)\in I_{\sigmae}$.

We next investigate the selector vertices $g_{i'}$.
We do not need to consider the case $i'=i$ as we already proved $(g_{i}, *)\notin I_{\sigma}, I_{\sigmae}$.
Consider the case $\indbit^{\sigma}_{i'}=1$ first, implying $i'\geq\relbit{\sigma}>i$ by \Pref{REL2}.
Since $i'\geq\relbit{\sigma}>1$, we have $\sigmabar(g_{i'})=\indbit_{i'+1}, \sigmabar(d_{i'})$ and $\sigmabar(s_{i'})$.
We thus prove $\valu_{\sigma}(F_{i',\indbit_{i'+1}})\geq\valu_{\sigma}(F_{i',1-\indbit_{i'+1}})$.
By the previously shown properties, \[\valu_{\sigma}(F_{i',\indbit_{i'+1}})=\rew{s_{i',\indbit_{i'+1}},h_{i',\indbit_{i'+1}}}+\valu_{\sigma}(b_{i'+1}).\]
Now, either $\valustar_{\sigma}(F_{i',1-\indbit_{i'+1}})=\frac{1}{2}\valustar_{\sigma}(g_1)+\frac{1}{2}\valustar_{\sigma}(b_2)$ or $\valustar_{\sigma}(F_{i',1-\indbit^{\sigma}_{i'+1}})=\valustar_{\sigma}(b_2)$ by assumption~4.
As $\valu_{\sigma}(b_2)>\valu_{\sigma}(g_1)$ by \Cref{lemma: Valuations in Phase Three}, it suffices to consider the second case.
The statement then follows since \begin{align*}
	\valustar_{\sigma}(F_{i',\indbit_{i'+1}})&=\rew{s_{i',\indbit_{i'+1}},h_{i',\indbit_{i'+1}}}+\valustar_{\sigma}(b_{i'+1})=\rew{s_{i',\indbit_{i'+1}}, h_{i',\indbit_{i'+1}}}+L_{i'+1}\\
		&\geq L_{2,i'}+L_{i'+1}=L_2=\valustar_{\sigmae}(b_2)\geq\valustar_{\sigma}(F_{i',1-\indbit_{i'+1}}).
\end{align*}
Consequently, $(g_{i'},*)\notin I_{\sigma}$.
As $i'\neq i$, these estimations can also be applied to $\sigmae$, implying $(g_{i'},*)\notin I_{\sigmae}$.
Hence $\indbit^{\sigma}_{i'}=1$ implies $(g_{i'},*)\notin I_{\sigma},I_{\sigmae}$.

Next assume $\indbit_{i'}=0$ and $i'>i$.
Then, by assumption~4, either [$\sigmabar(g_{i'})=\indbit_{i'+1}$ and $F_{i',0},F_{i',1}$ are mixed] or [$\sigmabar(g_{i'})=1-\indbit_{i'+1}$, $F_{i',1-\indbit_{i'+1}}$ is $b_2$-open and $F_{i',\indbit^{\sigma}_{i'+1}}$ is mixed].
In the first case, both cycle centers are in the same state with respect to both $\sigma$ and $\sigmae$.
Consequently, it suffices to prove $\valu_{\sigma}(s_{i',\indbit_{i'+1}})>\valu_{\sigma}(s_{i',1-\indbit_{i'+1}})$.
But this follows as $\sigma(s_{i',1-\indbit_{i'+1}})=b_1, \sigma(s_{i',\indbit_{i'+1}})=h_{i',\indbit_{i'+1}}$ and $\sigma(b_{\relbit{\sigma}})=g_{\relbit{\sigma}}$.
As these arguments also apply to $\sigmae$, it follows that $(g_{i'},*)\notin I_{\sigma},I_{\sigmae}$.
In the second case, the argument follows since $\valustar_{\sigma}(b_2)>\valustar_{\sigma}(g_1)$ by \Cref{lemma: Valuations in Phase Three}.
By the same argument, $(g_{i'},*)\notin I_{\sigmae}$.
This concludes the case $\indbit_{i'}=0$ and $i'>i$.

Consider the case $\indbit_{i'}=0\wedge i'<i$ next.
Then, since $i'<i\leq\relbit{\sigma}-1$, \Pref{BR1} implies $\sigmabar(g_{i'})=1$.
We thus prove $\valustar_{\sigma}(F_{i',1})>\valustar_{\sigma}(F_{i',0})$.
By assumption~2, it holds that $\sigma(s_{i',0})=h_{i',0}$ and $\sigma(s_{i',1})=h_{i',1}$.
Since $F_{i'\sigmabar(g_{i'})}$ is closed by assumption~2, this implies $\valustar_{\sigma}(F_{i',1})=\rew{s_{i',1},h_{i',1}}+\valustar_{\sigma}(g_{i'+1}).$
Using \Cref{corollary: Complete Valuation Of Selection Vertices MDP}, it follows that $\valustar_{\sigma}(g_{i'+1})=\sum_{\ell=i'+1}^{\relbit{\sigma}-1}W_{\ell}+\valustar_{\sigma}(b_{\relbit{\sigma}+1})$.
By assumption~2., the other cycle center $F_{i',0}$ of level $i'$ is $g_1$-halfopen.
Let $i=\relbit{\sigma}-1$, implying $\valustar_{\sigma}(g_1)=R_1$.
Then \begin{align*}
	\valustar_{\sigma}(F_{i',1})&=\rew{s_{i',1},h_{i',1}}+\valustar_{\sigma}(g_{i'+1})=\rew{s_{i',1},h_{i',1}}+\sum_{\ell=i'+1}^{\relbit{\sigma}-1}W_{\ell}+\valustar_{\sigma}(b_{\relbit{\sigma}+1})\\
		&>\sum_{\ell=1}^{i'}W_{\ell}+\sum_{\ell=i'+1}^{\relbit{\sigma}-1}W_{\ell}+\valustar_{\sigma}(b_{\relbit{\sigma}+1})=R_1=\valustar_{\sigma}(g_1)=\valustar_{\sigma}(F_{i',0}).
\end{align*}

Let $i<\relbit{\sigma}-1$, implying $\valustar_{\sigma}(g_1)=\rew{g_{i+1}}+\sum_{\ell<i+1}W_{\ell}+\valustar_{\sigma}(b_2)$.
Then, due to $\sigmabar(eb_{i+1})\wedge\nsigmabar(eg_{i+1})$, it follows that $\valustar_{\sigma}(g_{i'+1})=\rew{g_{i+1}}+\sum_{\ell=i'+1}^{i}W_{\ell}+\valustar_{\sigma}(b_2)$.
Consequently, \begin{align*}
	\valustar_{\sigma}(F_{i',1})&=\rew{s_{i',1}, h_{i',1}}+\valustar_{\sigma}(g_{i'+1})=\rew{s_{i',1},h_{i',1},g_{i+1}}+\sum_{\ell=i'+1}^{i}W_{\ell}+\valustar_{\sigma}(b_2)\\
		&>\rew{g_{i+1}}+\sum_{\ell=1}^{i'}W_{\ell}+\sum_{\ell=i'+1}^{i}W_{\ell}+\valustar_{\sigma}(b_2)=\valustar_{\sigma}(g_1)=\valustar_{\sigma}(F_{i',0}).
\end{align*}
Hence, $(g_{i'},*)\notin I_{\sigma}$.
Since it holds that $\valustar_{\sigmae}(g_{i'+1})=\rew{g_i}+\sum_{\ell=i'+1}^{i-1}W_{\ell}+\valustar_{\sigma}(b_2)$ and $\valustar_{\sigmae}(g_1)=\rew{g_i}+\sum_{\ell<i}W_{\ell}+\valustar_{\sigmae}(b_2)$, the same calculation can be used to obtain $\valustar_{\sigmae}(F_{i',1})>\valustar_{\sigmae}(F_{i',0})$.
Hence, also $(g_{i'},*)\notin I_{\sigmae}$.

This covers all cases, hence $(g_{i'},F_{i',*})\notin I_{\sigma},I_{\sigmae}$ for any index $i'\in[n]$.

We next consider entry vertices $b_{i'}$ for $i'\in[n]$.
First of all, since $\sigma(b_1)=\sigmae(b_1)=g_1$ and $\valu_{\sigma}(b_2)>\valu_{\sigma}(g_1)$ as well as $\valu_{\sigmae}(b_2)>\valu_{\sigmae}(g_1)$ by \Cref{lemma: Valuations in Phase Three}, we have $(b_1,b_2)\in I_{\sigma},I_{\sigmae}$.
Thus consider some edge $(b_{i'},b_{i'+1})$ for $i'\neq 1$.
Since $(b_{i'},b_{i'+1})\notin I_{\sigma},I_{\sigmae}$ if $\sigma(b_{i'})=b_{i'+1}$, assume $\sigma(b_{i'})=g_{i'}$.
Then $\indbit_{i'}=1$, implying $i'\geq\relbit{\sigma}$.
This directly implies $\valu_{\sigma}(b_{i'})=L_{i'}=W_{i'}+L_{i'+1}>L_{i'+1}=\valu_{\sigma}(b_{i'+1}),$ hence $(b_{i'},b_{i'+1})\notin I_{\sigma}$.
The same argument can be used to prove $(b_{i'},b_{i'+1})\notin I_{\sigmae}$.

Next, consider some edge $(b_{i'},g_{i'})$.
If $\sigma(b_{i'})=g_{i'}$, then $(b_{i'},g_{i'})\notin I_{\sigma},I_{\sigmae}$, hence assume $\sigma(b_{i'})=b_{i'+1}$, implying $i'>1$.
By \Pref{EV1}$_{i'}$, it then holds that $\indbit_{i'}=0$.
Let $i'>i$.
By assumption~4, the cycle center $F_{i',\sigmabar(g_{i'})}$ is then either mixed or $b_2$-open.
Since $\valu_{\sigma}(b_2)>\valu_{\sigma}(g_1)$, it suffices to consider the case that it is $b_2$-open.
Thus,\[\valustar_{\sigma}(g_{i'})\leq\rew{g_{i'}}+\valu_{\sigma}(b_2)=\rew{g_{i'}}+L_2=\rew{g_{i'}}+L_{2,i'-1}+L_{i'+1}<L_{i'+1}=\valustar_{\sigma}(b_{i'+1}),\]hence $(b_{i'},g_{i'})\notin I_{\sigma}$ and, by the same arguments, also $(b_{i'},g_{i'})\notin I_{\sigmae}$.
Now consider the case $i'<i$.
Then, by assumption~2, the cycle center $F_{i',\sigmabar(g_{i'})}$ is closed.
Depending on whether $i=\relbit{\sigma}-1$ or $i<\relbit{\sigma}-1$, we then have \[\valustar_{\sigma}(g_{i'})=\sum_{\ell=i'}^{\relbit{\sigma}-1}W_{\ell}+\valustar_{\sigma}(b_{\relbit{\sigma}+1})\quad \text{ or } \quad\valustar_{\sigma}(g_{i'})=\rew{g_{i+1}}+\sum_{\ell=i'}^{i}W_{\ell}+\valustar_{\sigma}(b_2).\]
The statement follows in either case since $i'<i<\relbit{\sigma}$ implies $\valustar_{\sigma}(b_{i'+1})=\valustar_{\sigma}(b_2)$ and since $\sigma(b_{\relbit{\sigma}})=g_{\relbit{\sigma}}$.
As the same arguments can be applied to $\sigmae$, this implies $(b_{i'},g_{i'})\notin I_{\sigma}, I_{\sigmae}$.
The case $i'=i$ can be shown by similar arguments since $\valustar_{\sigma}(b_{i+1})=L_2$ and in particular $W_{\relbit{\sigma}}\subset L_2$ due to $\sigma(b_{\relbit{\sigma}})=g_{\relbit{\sigma}}$.
Hence, $(b_{i'},*)\in I_{\sigma}\Leftrightarrow(b_{i'},*)\in I_{\sigmae}$ for all $i'\in[n]$.

We next consider upper selection vertices $s_{i',j'}$ for arbitrary $i',j'$.
We already proved that $(i',j')\neq(i,j)$ implies $(s_{i',j'},b_1)\notin I_{\sigma},I_{\sigmae}$.
We thus only prove $(s_{i',j},h_{i',j'})\notin I_{\sigma},I_{\sigmae}$ for arbitrary $i',j'$.
This is immediate if $\sigma(s_{i',j'})=h_{i',j'}$, so let $\sigma(s_{i',j'})=b_1$.
This implies $i'>i$ due to assumption~2 as well as $\valu_{\sigma}(s_{i',j'})=\rew{s_{i',j'}}+\valu_{\sigma}(g_1)$ and $j'=1-\indbit_{i'+1}$.
We now distinguish several cases.
First assume $i=\relbit{\sigma}-1$, implying $\valustar_{\sigma}(g_1)=R_1$.
Also, since $i'>i=\relbit{\sigma}-1$, we have $i'\geq\relbit{\sigma}$.
Assume $\indbit_{i'+1}=0$, implying $j'=1-\indbit_{i'+1}=1$.
Thus, $\valustar_{\sigma}(h_{i',j'})=\rew{h_{i',j'}}+\valustar_{\sigma}(g_{i'+1})$.
By assumption~2, $F_{i'+1,\sigmabar(g_{i'+1})}$ is either mixed or $b_2$-open.
Since the valuation of the cycle center is larger if it is $b_2$-open, it suffices to consider this case.
Using $\indbit_{i'+1}=0$, we thus have \begin{align*}
	\valustar_{\sigma}(h_{i',j'})&\leq\rew{h_{i',j'},g_{i'+1}}+\valustar_{\sigma}(b_2)=\rew{h_{i',j'},g_{i'+1}}+L_2\\
		&<L_{i'+2}=L_{i'+1}\leq L_{\relbit{\sigma}+1}<\sum_{\ell<\relbit{\sigma}}W_{\ell}+L_{\relbit{\sigma}+1}=\valustar_{\sigma}(g_1),
\end{align*} implying $(s_{i',j'},h_{i',j'})\notin I_{\sigma}$ if $\indbit_{i'+1}=0$.
Consider the case $\indbit_{i'+1}=1$, implying $j'=0$.
Then \begin{align*}
	\valustar_{\sigma}(h_{i',j'})&<L_{i'+1}\leq L_{\relbit{\sigma}+1}<\sum_{\ell<\relbit{\sigma}}W_{\ell}+L_{\relbit{\sigma}+1}=\valustar_{\sigma}(g_1)
\end{align*}
implying $(s_{i',j'},h_{i',j'})\notin I_{\sigma}$ if $\indbit_{i'+1}=1$.
This concludes the case $i=\relbit{\sigma}-1$.
Hence assume $i<\relbit{\sigma}-1$, implying $\valustar_{\sigma}(g_1)=\rew{g_{i+1}}+\sum_{\ell<i+1}W_{\ell}+\valustar_{\sigma}(b_2)$.
Note that it thus might happen that $i'\leq\relbit{\sigma}$.
Consider the case $\indbit_{i'+1}=0$.
By the same arguments used for the case $i=\relbit{\sigma}-1$, we then have $\valustar_{\sigma}(h_{i',j'})<L_{i'+1}.$
If $i'>\relbit{\sigma}$, then $i+1<\relbit{\sigma}$ thus implies \begin{align*}
	\valustar_{\sigma}(h_{i',j'})&<L_{i'+1}\leq L_{\relbit{\sigma}+1}<\rew{g_{i+1}}+\sum_{\ell<i+1}W_{\ell}+W_{\relbit{\sigma}}+L_{\relbit{\sigma}+1}\\
		&=\rew{g_{i+1}}+\sum_{\ell<i+1}W_{\ell}+\valustar_{\sigma}(b_2)=\valustar_{\sigma}(g_1),
\end{align*}
hence $(s_{i',j'},h_{i',j'})\notin I_{\sigma}$.
If $i'<\relbit{\sigma}$, then $i<i'$ implies \begin{align*}
	\valustar_{\sigma}(h_{i',j'})&\leq\rew{h_{i',j'},g_{i'+1}}+L_{2,i'}+L_{i'+1}=\rew{h_{i',j'},g_{i'+1}}+L_{\relbit{\sigma}}\\
		&<\rew{g_{i+1}}+\sum_{\ell<i+1}W_{\ell}+\valustar_{\sigma}(b_2)=\valustar_{\sigma}(g_1),
\end{align*}
hence $(s_{i',j'},h_{i',j'})\notin I_{\sigma}$.
Thus consider the case $\indbit_{i'+1}=1$ , implying $i'\geq\relbit{\sigma}-1$.
Since we consider the case $i<\relbit{\sigma}-1$, it holds that $j=1-\indbit_{i'+1}=0$, implying \begin{align*}
\valustar_{\sigma}(h_{i',j'})&=\rew{h_{i',j'}}+\valustar_{\sigma}(b_{i'+2})=\rew{h_{i',j'}}+L_{i'+2}\\
	&<\rew{g_{i+1}}+\sum_{\ell<i+1}W_{\ell}+W_{i'+1}+L_{i'+2}\\
	&=\rew{g_{i+1}}+\sum_{\ell<i+1}W_{\ell}+L_{i'+1}\leq\rew{g_{i+1}}+\sum_{\ell<i+1}W_{\ell}+L_{\relbit{\sigma}}=\valustar_{\sigma}(g_1),
\end{align*}
hence $(s_{i',j'},h_{i',j'})\notin I_{\sigma}$.
Thus, under all circumstances, $(s_{i',j'},h_{i',j'})\notin I_{\sigma}$.

We now prove $(s_{i',j'},h_{i',j'})\notin I_{\sigmae}$. 
We have $\valustar_{\sigmae}(g_1)=\rew{g_i}+\sum_{\ell<i}W_{\ell}+\valustar_{\sigmae}(b_2)$.
Let $\indbit_{i'+1}=0$ first.
Then, using the same arguments used for $\sigma$ as well as $i'>i$, we obtain \begin{align*}
	\valustar_{\sigmae}(h_{i',j'})&=\rew{h_{i',j'}}+\valustar_{\sigmae}(g_{i'+1})<\rew{g_i}+\sum_{\ell<i}W_{\ell}+\valustar_{\sigmae}(b_2)=\valustar_{\sigmae}(g_1),
\end{align*}
so $(s_{i',j'},h_{i',j'})\notin I_{\sigmae}$.
Thus let $\indbit_{i'+1}=1$.
Then, by the same arguments used before and $i'>i$, it follows that\begin{align*}
	\valustar_{\sigmae}(h_{i',j'})&=\rew{h_{i',j'}}+\valustar_{\sigmae}(b_{i'+2})=\rew{h_{i',j'}}+L_{i'+2}<\rew{g_i}+\sum_{\ell<i}W_{\ell}+W_{i'+1}+L_{i'+2}\\
		&=\rew{g_{i}}+\sum_{\ell<i}W_{\ell}+L_{i'+1}\leq\rew{g_{i}}+\sum_{\ell<i}W_{\ell}+L_2=\valustar_{\sigmae}(b_2),
\end{align*}
implying $(s_{i',j'},h_{i',j'})\notin I_{\sigmae}$.
We thus have $(s_{i',j'},*)\notin I_{\sigma},I_{\sigmae}$ for all  indices $i',j'$ with the exception of the edge $(s_{i,j},b_1)$.

Since there are no indices $i',j',k'$ besides $i,j,k$ such that $(d_{i',j',k'},e_{i',j',k'})\in I_{\sigma}$ by assumption, it suffices to prove $(d_{i',j',k'},e_{i',j',k'})\notin I_{\sigmae}$ for all such indices.
The statement follows if $\sigma(d_{i',j',k'})=e_{i',j',k'}$, hence assume $\sigma(d_{i',j',k'})=F_{i',j'}$.
Let $\sigma(e_{i',j',k'})=b_2$.
Then, since $\valu_{\sigma}(b_2)=\valu_{\sigmae}(b_2)=L_2$, the valuation of $e_{i',j',k'}$ does not increase by the application of $e$.
As $(d_{i',j',k'},e_{i',j',k'})\notin I_{\sigma}$ implies $\valu_{\sigma}(F_{i',j'})\geq\valu_{\sigma}(e_{i',j',k'})$ and the valuation of $F_{i',j'}$ can only increase, this implies $(d_{i',j',k'},e_{i',j',k'})\notin I_{\sigmae}$.
Thus let $\sigma(e_{i',j',k'})=g_1$, implying $\valustar_{\sigmae}(e_{i',j',k'})=\rew{g_i}+\sum_{\ell<i}W_{\ell}+\valustar_{\sigmae}(b_2)$.
Assume that $F_{i',j'}$ is not closed with respect to $\sigma$.
Then the assumption $\sigmae(d_{i',j',k'})=F_{i',j'}$ implies that the cycle center is halfopen with respect to $\sigma$.
Due to the assumptions of this lemma, it is easy to see that this implies that $F_{i',j'}$ is $g_1$-halfopen with respect to both $\sigma$ and $\sigmae$ and that we either have $i'=i$ and $j'=1-j$ or $i'<i<\relbit{\sigmae}$ and $j'=1-\sigmaebar(g_{i'})$.
In both cases, we have $j'=1-\sigmaebar(g_{i'})=\indbit_{i'+1}$ by \Pref{BR1} and $\sigmae(s_{i',j'})=h_{i',j'}$ by \Pref{USV2}$_{i',\indbit_{i'+1}}$.
Consequently, by \Pref{EV1}$_{i'+1}$ and $i'+1\leq i+1\leq \relbit{\sigmae}$, we obtain\begin{align*}
	\valustar_{\sigmae}(s_{i',j'})&=\rew{s_{i',j'},h_{i',j'}}+\valustar_{\sigmae}(b_{i'+1})=\rew{s_{i',j'},h_{i',j}}+\valustar_{\sigmae}(b_{\relbit{\sigma}})>\valustar_{\sigmae}(b_2),
\end{align*}
implying $\valu_{\sigmae}(F_{i',j'})>\valu_{\sigmae}(g_1)$ as $\valustar_{\sigmae}(b_2)>\valustar_{\sigmae}(g_1)$.
Consequently, it holds that $(d_{i',j',k'},e_{i',j',k'})\notin I_{\sigmae}$.
Now let $F_{i',j'}$ be closed with respect to $\sigma$.
Consider the case $\indbit_{i'}=1\wedge\indbit_{i'+1}=j'$ first, implying $i'>i$ as $i<\relbit{\sigmae}=\nsb$ by assumption.
Then, \Pref{USV1}$_{i'}$ implies \begin{align*}
	\valustar_{\sigmae}(F_{i',j'})&=\rew{s_{i',j'}, h_{i',j'}}+\valu_{\sigmae}(b_{i'+1})=\rew{s_{i',j'},h_{i',j'}}+L_{i'+1}>L_{2,i'}+L_{i'+1}\\
		&>\rew{g_i}+\sum_{\ell<i}W_{\ell}+L_{2,i'}+L_{i'+1}=\rew{g_i}+\sum_{\ell<i}W_{\ell}+L_2=\valustar_{\sigmae}(g_1),
\end{align*}
hence $(d_{i',j',k'},e_{i',j',k'})\notin I_{\sigmae}$.
Next assume $\indbit_{i'}=1\wedge\indbit_{i'+1}\neq j'$.
Then, \Pref{USV1}$_{i'}$ implies $\sigmae(s_{i',j'})=b_1$ as $i'\geq\relbit{\sigmae}=\nsb$.
Since $F_{i',j'}$ is closed, we then have \[\valustar_{\sigmae}(F_{i',j'})=\valustar_{\sigmae}(s_{i',j'})=\rew{s_{i',j'}}+\valustar_{\sigmae}(g_1)>\valustar_{\sigmae}(g_1),\]hence $(d_{i',j',k'},e_{i',j',k'})\notin I_{\sigmae}$.
Since $\indbit_{i'}=0\wedge\indbit_{i'+1}=j'$ is impossible if $F_{i',j'}$ is closed, it remains to consider the case $\indbit_{i'}=0\wedge\indbit_{i'+1}\neq j'$.
By assumption~4, we then need to have $i'\leq i<\relbit{\sigmae}$ as well as $\sigmae(s_{i',j'})=h_{i',j'}$.
Consider the case $i'=i$ first, implying $j'=j$.
As we just applied the switch $(d_{i,j,k},e_{i,j,k})$, it is clear that this switch is not improving for $\sigmae$.
Hence consider $(d_{i,j,1-k},e_{i,j,1-k})$.
We have \[\valu_{\sigmae}(F_{i,j})=\frac{1-\e}{1+\e}\valu_{\sigmae}(b_2)+\frac{2\e}{1+\e}\valu_{\sigmae}(s_{i,j}).\]
If $\sigma(e_{i,j,1-k})=g_1$, then $(d_{i,j,1-k},e_{i,j,1-k})\notin I_{\sigmae}$ due to $\valustar_{\sigmae}(b_2)>\valustar_{\sigmae}(g_1)$.
It cannot happen that $\sigma(e_{i,j,1-k})=b_2$ since this would imply $(d_{i,j,1-k},e_{i,j,1-k})\in I_{\sigma}$, contradicting the assumption.
The reason is that $W_{\relbit{\sigma}}$ is not part of the valuation of~$F_{i,j}$ which results in $\valustar_{\sigma}(F_{i,j})=\valustar_{\sigma}(s_{i,j})<\valustar_{\sigma}(b_2)=\valustar_{\sigma}(e_{i,j,1-k})$.
Hence $(d_{i,j,1-k},e_{i,j,1-k})\notin I_{\sigmae}$ and we consider the case $i'<i$ next.
Then, since $\sigmaebar(d_{i'})$ by assumption~2 and since $i'<i<\relbit{\sigmae}$, we need to have $j'=\sigmaebar(g_{i'})=1-\indbit_{i'+1}=1$.
Consequently, \begin{align*}
	\valustar_{\sigmae}(F_{i',j'})&=\valustar_{\sigmae}(s_{i',1})=\rew{s_{i',1},h_{i',1}}+\valustar_{\sigmae}(g_{i'+1})\\
		&=\rew{s_{i',1},h_{i',1}}+\rew{g_i}+\sum_{\ell=i'+1}^{i-1}W_{\ell}+\valustar_{\sigmae}(b_2)\\
		&>\sum_{\ell=1}^{i'}W_{\ell}+\sum_{\ell=i'+1}^{i-1}W_{\ell}+\rew{g_i}+\valustar_{\sigmae}(b_2)=\valustar_{\sigmae}(g_1).
\end{align*}
Thus, if $\sigma(e_{i',j',k'})=g_1$, then $(d_{i',j',k'},e_{i',j',k'})\notin I_{\sigmae}$.
Since \[\valustar_{\sigma}(F_{i',j'})\leq\valustar_{\sigmae}(F_{i',j'})<\valustar_{\sigmae}(b_2)=\valustar_{\sigma}(b_2),\] also $\sigma(e_{i',j',k'})=g_1$ has to hold since $\sigma(e_{i',j',k'})=b_2$ implies $(d_{i',j',k'},e_{i',j',k'})\in I_{\sigma}$, contradicting our assumption.
Consequently, $(d_{i',j',k'},e_{i',j',k'})\notin I_{\sigmae}$ for all indices $i',j',k'$.

It remains to consider edges $(d_{*,*,*},F_{*,*})$.
We prove that $(d_{i',j',k'},F_{i',j'})\in I_{\sigma}\Leftrightarrow (d_{i',j',k'},F_{i',j'})\in I_{\sigmae}$.
If $\sigma(d_{i',j',k'})=F_{i',j'}$, then $(d_{i',j',k'},F_{i',j'})\notin I_{\sigma},I_{\sigmae}$ and the statement follows.
Note that also $(d_{i,j,k},F_{i,j})\notin I_{\sigma},I_{\sigmae}$.
Hence fix some indices $i',j',k'$ with $\sigma(d_{i',j',k'})=e_{i',j',k'}$.
Then, the cycle center $F_{i',j'}$ is not closed with respect to $\sigma$ or $\sigmae$, implying $\indbit_{i'}=0\vee \indbit _{i'+1}\neq j'$.
Consider the case $\indbit_{i'}=0$ first and assume $i'>i$.
By assumption~4., $F_{i',j'}$ is ether mixed or $b_2$-open.
Consider the case that $F_{i',j'}$ is mixed.
Then \begin{align*}
	\valustar_{\sigma}(F_{i',j'})&=\frac{1}{2}\valustar_{\sigma}(b_2)+\frac{1}{2}\valustar_{\sigma}(g_1), \\
	\valustar_{\sigmae}(F_{i',j'})&=\frac{1}{2}\valustar_{\sigmae}(b_2)+\frac{1}{2}\valustar_{\sigmae}(g_1).
\end{align*}
If $\sigma(e_{i',j',k'})=g_1$, then $\valustar_{\sigma}(e_{i',j',k'})<\valustar_{\sigma}(F_{i',j'})$ and $\valustar_{\sigmae}(e_{i',j',k'})<\valustar_{\sigmae}(F_{i',j'})$ by \Cref{lemma: Valuations in Phase Three}, implying $(d_{i',j',k'},F_{i',j'})\in I_{\sigma},I_{\sigmae}$.
If $\sigma(e_{i',j',k'})=b_2$, then \Cref{lemma: Valuations in Phase Three} implies $\valustar_{\sigma}(e_{i',j',k'})>\valustar_{\sigma}(F_{i',j'})$ and $\valustar_{\sigmae}(e_{i',j',k'})>\valustar_{\sigmae}(F_{i',j'})$ and it follows that $(d_{i',j',k'},F_{i',j'})\notin I_{\sigma},I_{\sigmae}$.
Next assume that $F_{i',j'}$ is $b_2$-open.
If $\indbit_{i'+1}\neq j'$, then assumption~3 implies $\sigma(s_{i',j'})=b_1$ and thus $\valu_{\sigma}(s_{i',j'})=\rew{s_{i',j'}}+\valu_{\sigma}(g_1)$.
By \Cref{lemma: Valuations in Phase Three}, it holds that $\rew{s_{i',j'}}+\valu_{\sigma}(g_1)<\valu_{\sigma}(b_2)$, implying $\valu_{\sigma}(F_{i',j'})<\valu_{\sigma}(e_{i',j',k'})$.
Since the same holds for $\sigmae$, this implies $(d_{i',j',k'},F_{i',j'})\notin I_{\sigma},I_{\sigmae}$.
If $\indbit_{i'+1}=j'$, then \Pref{EV1}$_{i'+1}$ implies that $\valu_{\sigma}(s_{i',j'})=\rew{s_{i',j'},h_{i',j'}}+\valu_{\sigma}(b_{i'+1})$ and thus $\valu_{\sigma}(s_{i',j'})>\valu_{\sigma}(b_2)$ as $\rew{s_{i',j'},h_{i',j'}}>L_{2,i'}$.
Since the same holds for $\sigmae$, we thus have $(d_{i',j',k'},F_{i',j'})\in I_{\sigma},I_{\sigmae}$.

Next, assume $\indbit_{i'+1}=0$ and $i'\leq i<\relbit{\sigmae}$.
Then, $F_{i',j'}$ is closed for $\sigma$ if $j'=\sigmabar(g_{i'})$ and $g_1$-halfopen if $j'=1-\sigmabar(g_{i'})$ by either assumption~1 or assumption~2.
Since we assume $\sigma(d_{i',j',j'})=e_{i',j',k'}$, it suffices to consider the second case.
Then $\sigma(e_{i',j',k'})=g_1$ and \[\valu_{\sigma}(F_{i',j'})=\frac{1-\e}{1+\e}\valu_{\sigma}(g_1)+\frac{2\e}{1+\e}\valu_{\sigma}(s_{i',j'}).\]
Now, by \Pref{BR1}, $\sigmabar(g_{i'})=1-\indbit_{i'+1}$, implying $1-\sigmabar(g_{i'})=j'=\indbit_{i'+1}$.
Consequently, by \Pref{USV2}$_{i'}$ and \Pref{EV1}$_{i'+1}$, we have \[\valu_{\sigma}(s_{i',j'})=\rew{s_{i',j'},h_{i',j'}}+\valu_{\sigma}(b_{i'+1})>\valu_{\sigma}(b_2)>\valu_{\sigma}(g_1).\]
This implies $\valu_{\sigma}(F_{i',j'})>\valu_{\sigma}(g_1)=\valu_{\sigma}(e_{i',j',k'})$, hence $(d_{i',j',k'},F_{i',j'})\in I_{\sigma}$.
Since the same arguments can be used for $\sigmae$, we also obtain $(d_{i',j',k'},F_{i',j'})\in I_{\sigmae}$.

Finally, consider the case $\indbit_{i'}=1\wedge\indbit_{i'+1}\neq j'$.
Then $i'\geq\relbit{\sigma}>i$, hence $F_{i',j'}$ is either mixed or $b_2$-open by assumption~5.
We however already showed that this implies $(d_{i',j',k'},F_{i',j'})\in I_{\sigma}\Leftrightarrow (d_{i',j',k'},F_{i',j'})\in I_{\sigmae}$.
\end{proof}

\USVInPhaseThree*

\begin{proof}
By the choice of $e$, $\sigmae$ is a phase-3-strategy for $\bit$ with $\sigmae\in\reach{\sigma_0}$.
Consequently, by the choice of $i$ and \Pref{REL2}, we have $\nsb=\relbit{\sigma}=\relbit{\sigmae}>1$.
We prove that $\sigmae$ is well-behaved.
Since $\relbit{\sigmae}>1$ and since $\sigmaebar(eb_{i,j})\wedge\nsigmaebar(eg_{i,j})$, we only need to reevaluate \Pref{S2}.
We show that the premise of this property cannot be fulfilled.
By \Pref{BR1}, we have $\sigmabar(g_i)=j$, hence $\nsigmaebar(d_i)$.
As $\relbit{\sigmae}>1$ implies $\sigmae(b_1)=g_1$, assume $\sigmae(b_2)=g_2$.
Then, by \Pref{BR1}, $\sigmaebar(g_1)=0\neq 1=\sigmaebar(b_1)$ and $\sigmaebar(b_2)$.
Consequently, $1\in\incorrect{\sigmae}$.
However, since $\sigmaebar(b_{i'})$ implies $\sigmaebar(g_{i'})=\sigmaebar(b_{i'+1})$ by \Pref{EV2}$_{i'}$ for $i'>1$, this implies $\incorrect{\sigmae}=\{1\}$ and thus $\relbit{\sigmae}=2$.
But then, $i=1$, hence the premise of the property cannot be fulfilled.

It remains to show that $I_{\sigmae}=I_{\sigma}\setminus\{e\}$.
Applying $e=(s_{i,j},b_1)$ increases the valuation of~$F_{i,j}$.
However, since $\sigmaebar(eb_{i,j})\wedge\nsigmaebar(eg_{i,j})$, the valuation is only changed by terms of order $o(1)$.
It is now easy but tedious to prove that the increase of the valuation of $F_{i,j}$ by terms of order $o(1)$ neither creates further improving switches nor makes improving switches unimproving.
This implies the statement.
\end{proof}

\OpenCompletelyInPhaseThree*

\begin{proof}
Since $F_{i,j}$ is $t^{\rightarrow}$-halfopen and by the choice of $e$, it holds that $\sigma(d_{i,j,1-k})=e_{i,j,1-k}$.
We can thus apply \Cref{lemma: Phase 3 Well-Behaved} and only need to prove $I_{\sigmae}=I_{\sigma}\setminus\{e\}$.

Since $F_{i,j}$ is $t^{\rightarrow}$-halfopen with respect to $\sigma$ and $t^{\rightarrow}$-open with respect to $\sigmae$, the valuation of $F_{i,j}$ only changes by terms of order $o(1)$ when applying the switch $e$.
It is easy but tedious to verify that this implies $I_{\sigmae}=I_{\sigma}\setminus\{e\}$.
\end{proof}

\StartOfPhaseFourValu*

\begin{proof}
	We first show $\relbit{\sigmae}=1$.
	Since $\sigma$ is a phase-$3$-strategy for $\bit$, it has \Pref{EV1}$_i$ and \Pref{EV2}$_i$ for all $i>1$.
	This implies $i\notin\incorrect{\sigma}$ for all $i>1$ and thus, by the choice of $e$, $i\notin\incorrect{\sigmae}$ for all $i>1$.
	Since $\sigmae(b_1)=b_2$, also $1\notin\incorrect{\sigmae}$, hence $\incorrect{\sigmae}=\emptyset$, implying $\relbit{\sigmae}=\min\{i\colon\sigmae(b_i)=b_{i+1}\}=1$.
	
	We now show that $\sigmae$ is a phase-$4$-strategy for $\bit$.
	By the choice of $e$ and the induced bit state, $\indbit^{\sigma}=\indbit^{\sigmae}=\bit+1\eqqcolon\indbit$.
	Since $\sigma$ is a phase-$3$-strategy and by the choice of $e$ it suffices to show that $\sigma$ has Properties (\ref{property: EV1})$_1$, (\ref{property: EV2})$_2$, (\ref{property: EV3})$_1$, (\ref{property: EV3})$_{\nsb}$, (\ref{property: CC2}) and (\ref{property: REL1}).
	Furthermore, we need to show that there is an index $i<\nsb$ with $\sigma(s_{i,1-\indbit_{i+1}})=h_{i,1-\indbit_{i+1}}$.
	
	First, $\sigmae$ has \Pref{EV3}$_{\nsb}$ since $\sigmabar(d_{i,j})\Leftrightarrow\indbit_{i}=1\wedge\indbit_{i+1}=j$.
	Second, the special condition as well as \Pref{ESC4}$_{i,j}$ and \Pref{ESC5}$_{i,j}$ are fulfilled for the relevant indices by assumption as well.
	In addition, $0=\indbit_1=\sigmaebar(d_{1,\indbit_2})$ by definition.
	Thus, since $\sigmae(b_1)=b_2$, $\sigmae$ has \Pref{EV1}$_1$ and consequently also \Pref{EV2}$_1$ and \Pref{EV3}$_1$.
	In addition, $\sigmae$ has \Pref{CC2} since $\sigma$ is a phase-$3$-strategy.
	Since $\incorrect{\sigmae}=\emptyset$, it also has \Pref{REL1}.
	Hence, $\sigmae$ is a phase-$4$-strategy.
	
	We show that $\sigmae$ is well-behaved.
	Since $\relbit{\sigmae}=1$ and since the target of the vertex $b_1$ changed when transitioning to $\sigmae$, the following assumptions need to be reevaluated.
	\begin{itemize}[align=right, leftmargin=1.75cm]
		\item[(\ref{property: S1})] Let $i\geq\relbit{\sigmae}=1$ and $\sigmae(b_i)=g_i$.
			By \Pref{EV1}$_i$ this implies $i\geq\nsb$ and thus, by \Pref{USV1}$_i$, $\sigmaebar(s_i)$.
		\item[(\ref{property: D1})] $\sigmae(b_i)=g_i$ implies $\sigmaebar(d_i)$ by \Pref{EV1}$_i$ and \Pref{EV2}$_i$.
		\item[(\ref{property: MNS1})] Assume that the premise was correct.
			Then, by \Cref{lemma: Config implied by Aeb}, $\minsige{b}=2$.
			This implies $\sigmae(b_2)=\sigma(b_2)=g_2$ and thus in particular $\relbit{\sigma}=\nsb=2$.
			But then, by \Pref{BR1} (applied to $\sigma)$, this implies $\sigma(g_1)=\sigmae(g_1)=F_{1,0}$ and thus $\minnegsige{g}=1$.
			This however contradicts the premise.
		\item[(\ref{property: MNS2})] Assume there was some $i<\minnegsige{g}<\minnegsige{s},\minsige{b}$ and that $\nsigmaebar(b_{\minnegsige{g}+1})$.
			Then $\minsige{b}=\nsb=\relbit{\sigma}$, implying $\minnegsige{g}=\minsige{b}-1=\relbit{\sigma}-1$ by \Pref{BR1}.
			But then $\sigmae(b_{\minnegsige{g}+1})=\sigmae(b_{\minsige{b}})=\sigmae(b_{\nsb})=g_{\nsb}$ by \Pref{EV1}$_{\nsb}$, contradicting the assumption.
		\item[(\ref{property: MNS3})] Assume there was some index $i<\minnegsige{s}\leq\minnegsige{g}<\minsige{b}$ and let $\ell\coloneqq\minnegsige{s}$.
			Then, $\sigmae(s_{\ell,\sigmaebar(g_{\ell})})=b_1$ and $\ell<\minsige{b}=\nsb$.
			But this contradicts the assumption that $\sigmae(s_{i',j'})=h_{i',j'}$ for all $i'<\nsb,j'\in\{0,1\}$.
			This argument also applies to Properties (\ref{property: MNS4}), (\ref{property: MNS5}) and (\ref{property: MNS6}), hence $\sigmae$ has all of these properties.
		\item[(EG*)] It can easily be checked that for all indices $i\in[n], j\in\{0,1\}$ not listed in either of the sets $S_1$ or $S_2$, $\sigmabar(d_{i,j})$ and thus $\sigmaebar(d_{i,j})$ holds.
			Hence $\sigmaebar(eg_{i,j})\implies\sigmaebar(eb_{i,j})$, so the premise of any of any of the assumptions (EG*) is incorrect.
		\item[(\ref{property: DN1})] By \Pref{EV1}$_n$, $\sigmaebar(d_n)$ implies $\sigmaebar(b_n)$.
		\item[(\ref{property: DN2})] We only need to consider this assumption if $\minnegsige{g}=n$.
			Since $\relbit{\sigma}\neq 1$, this implies $\sigmaebar(g_i)=\sigmabar(g_i)=1$ for all $i\in[n-1]$ by \Pref{BR1} (applied to $\sigma$).
			Thus, by assumption, $i\neq\relbit{\sigma}-1$ for all of those $i$, hence $n=\relbit{\sigma}-1$.
			But this implies $\relbit{\sigma}=n+1$, contradicting the definition of $\relbit{\sigma}$.
	\end{itemize}
	
	We now show the statements regarding the improving switches.
	First, $(s_{\nsb-1,0},b_1)\in I_{\sigmae}$ follows since $\sigma(s_{\nsb-1,0})=h_{\nsb-1,0}$ by assumption and since $\relbit{\sigmae}=1$ and $\sigmae(b_1)=b_2$ imply \begin{align*}
		\valustar_{\sigmae}^\P(b_1)=L_1^\P=L_{\nsb}^\P=W_{\nsb}^\P\cup L_{\nsb+1}^\P\rhd\{h_{\nsb-1,0}\}\cup L_{\nsb+1}^\P=\valustar_{\sigmae}^\P(h_{\nsb-1,0}).
	\end{align*}
	
	We now show $(s_{i,1},b_1)\in I_{\sigmae}$ for all $i\leq\nsb-2$.
	Fix some $i\leq\nsb-2$.
	Then, since $i+1<\nsb=\relbit{\sigma}$, \begin{align*}
		\valustar_{\sigma}^\P(h_{i,1})&=\{h_{i,1}\}\cup\valustar_{\sigma}^\P(g_{i+1})=\{h_{i,1}\}\cup R_{i+1}^\P\\
			&\rhd \bigcup_{i'=1}^{i}W_{i'}^\P\cup\bigcup_{i'=i+1}^{\relbit{\sigma}-1}W_{i'}^\P\cup L_{\relbit{\sigma}+1}^\P=R_1^\P=\valustar_{\sigma}^\P(b_1).
	\end{align*}
	Since $\sigma(s_{i,1})=\sigmae(s_{i,1})=h_{i,1}$, the statement then follows from \[\valustar_{\sigmae}^\P(b_1)=L_1^\P=L_{\relbit{\sigma}}^\P=W_{\relbit{\sigma}}^\P\cup L_{\relbit{\sigma}+1}^\P\rhd\bigcup_{i'=i+1}^{\relbit{\sigma}-1}W_{i'}^\P\cup\{h_{i,1}\}\cup L_{\relbit{\sigma}+1}^\P=\valustar_{\sigmae}^\P(h_{i,1}).\]
	
	We now show that the edges contained in the sets $X_0$ and $X_1$ are improving switches if $\indbit^{\sigma}$ is not a power of 2 and that no other edge is an improving switch otherwise.
	We distinguish the following cases.
	\begin{enumerate}
		\item Let $\indbit=2^k$ for some $k\in\mathbb{N}$, implying $\valustar_{\sigmae}^{\P}(b_1)=L_1^{\P}=W_{\relbit{\sigma}}^{\P}$.
			By applying the improving switch $e=(b_1,b_2)$ the valuation of $b_1$ increased.
			The only vertices with edges towards $b_1$ are upper selection vertices.
			We hence show that for any vertex~$s_{i,j}$, one of the following statements is true:
			\begin{enumerate}
				\item $\sigmae(s_{i,j})=h_{i,j}$ and $\valu_{\sigmae}^\P(h_{i,j})\unrhd\valu_{\sigmae}^\P(b_1)$.
				\item $\sigmae(s_{i,j})=h_{i,j}$ and $(s_{i,j},b_1)\in I_{\sigmae}$.
				\item $\sigmae(s_{i,j})=b_1$ and $\tau^{\sigma}(F_{i,j}),\tau^{\sigmae}(F_{i,j})\neq s_{i,j}$.
			\end{enumerate}
			We distinguish the following cases:
			\begin{itemize}
				\item \boldall{$i\leq\nsb-2$ and $j=0$:} Then, $\sigmae(s_{i,0})=h_{i,0}$.
					Also, $\sigmae(h_{i,0})=b_{i+2}$, so $i+2\leq \relbit{\sigma}$ implies $\valustar_{\sigmae}^\P(h_{i,0})=\{h_{i,0}\}\cup\valustar_{\sigmae}^\P(b_{i+2})=\{h_{i,0}\}\cup W_{\relbit{\sigma}}^\P\rhd \valustar_{\sigmae}^\P(b_1).$
				\item \boldall{$i\leq\nsb-2$ and $j=1$:} As proven before, all of these edges are improving.
				\item \boldall{$i=\nsb-1$ and $j=0$:} As proven before, $(s_{\relbit{\sigma}-1,0},b_1)$ is improving for $\sigmae$.
				\item \boldall{$i=\nsb-1$ and $j=1$:} By assumption, it holds that $\sigmae(s_{i,1})=h_{i,1}$.
					Thus, by the choice of~$i$, $\valustar_{\sigmae}^\P(h_{i,1})=\{h_{i,1}\}\cup\valustar_{\sigmae}^\P(g_{\relbit{\sigmae}})=\{h_{i,1}\}\cup W_{\relbit{\sigmae}}^\P\rhd\valustar_{\sigmae}^\P(b_1).$
				\item \boldall{$i=\nsb$ and $j=0$:} Since $\indbit^{\sigma}=2^k$, this then implies $\sigmae(s_{\nsb,0})=h_{\nsb,0}$ by \Pref{USV1}$_{\nsb}$.
					But then, $\valustar_{\sigmae}^\P(h_{\nsb,0})=\{h_{\nsb,0}\}\rhd W_{\nsb}^\P=\valustar_{\sigmae}^\P(b_1)$.
				\item \boldall{$i=\nsb$ and $j=1$:} Then, by \Pref{USV1}$_{\nsb}$, $\sigmae(s_{\nsb,1})=b_1$.
					We need to show $\tau^{\sigma}(F_{i,j})\neq s_{i,j}$ and $\tau^{\sigmae}(F_{i,j})\neq s_{i,j}$.
					This is done by showing that the first, second and fifth case of \Cref{lemma: Exact Behavior Of Counterstrategy} cannot occur.
					The first case cannot occur since $j=1=1-\indbit_{i+1}=1-\indbit_{i+1}$ and both $\sigma$ and $\sigmae$ have \Pref{EV3}$_{\nsb}$.
					The second case cannot occur with respect to both $\sigma$ and $\sigmae$ since there is no cycle center $F_{i,j}$ with $\sigmabar(eg_{i,j})\wedge\nsigmabar(eb_{i,j})$ by \Pref{ESC4}$_{i,j}$ and \Pref{ESC5}$_{i,j}$.
					The fifth case cannot occur for $\sigmae$ since $\relbit{\sigmae}=1$ and $\sigmae(s_{i,j})=b_1$.
					It can also not occur for $\sigma$ since $\indbit=2^k$ implies $\sigmabar(eb_{\relbit{\sigma},1})\wedge\sigmabar(eg_{\relbit{\sigma},1})$ by \Pref{ESC5}$_{\relbit{\sigma},1}$.
				\item \boldall{$i>\nsb$ and $j=0$:} Since $\indbit_{i'}=0$ for all $i'\neq\nsb$, $i>\nsb$ implies $\indbit_i=\indbit_{i+1}=0$.
					Hence, by \Pref{USV1}$_{i}$, $\sigmae(s_{i,0})=h_{i,0}$ and consequently $\valustar_{\sigmae}^\P(h_{i,0})=\{h_{i,0}\}\rhd W_{\relbit{\sigma}}^\P=\valustar_{\sigmae}^\P(b_1).$
				\item \boldall{$i>\nsb$ and $j=1$:} Then $\sigmae(s_{i,j})=\sigmae(s_{i,1})=b_1$ by \Pref{USV1}$_{i}$, hence it suffices to show $\tau^{\sigma}(F_{i,1}),\tau^{\sigmae}(F_{i,1})\neq s_{i,1}$.
					This is again proven by showing that the first, second and fifth case of \Cref{lemma: Exact Behavior Of Counterstrategy} cannot be fulfilled.
					Since $\indbit=2^k$ for some $k\in\mathbb{N}$ by assumption, $m=\max\{i\colon\sigmae(b_i)=g_i\}=\nsb$.
					Hence, by \Pref{ESC5}$_{i,1}$ (resp. by assumption), we have $\sigmaebar(eb_{i,1})\wedge\sigmaebar(eg_{i,j})$ and $\sigmabar(eb_{i,j})\wedge\sigmabar(eg_{i,j})$.
					Consequently, either the sixth or the seventh case of \Cref{lemma: Exact Behavior Of Counterstrategy} is true, both implying $\tau^{\sigmae}(F_{i,1}),\tau^{\sigma}(F_{i,1})\neq s_{i,1}$.
			\end{itemize}
		\item Now assume that there is no $k\in\mathbb{N}$ such that $\indbit=2^k$.
			We prove $X_0,X_1\subseteq I_{\sigmae}$ and that the edges contained in $I_{\sigmae}$ according to the lemma are indeed improving.
			Fix some $k\in\{0,1\}$.
			We prove $X_k\subseteq I_{\sigmae}$.
			\begin{itemize}
				\item We first show $(d_{i,j,k},F_{i,j})\in I_{\sigmae}$ where $i=\relbit{\sigma}=\nsb$ and $j=1-\indbit_{i+1}$.
					By assumption, $\sigma(d_{i,j,k})=\sigmae(d_{i,j,k})\neq F_{i,j}$.
					Hence $\sigmae(d_{i,j,k})=\sigma(d_{i,j,k})=e_{i,j,k}$ and it suffices to show $\valu_{\sigmae}^\P(F_{i,j})\rhd\valu_{\sigmae}^\P(e_{i,j,k})$.
					Since $\sigmae(d_{i,j,k})=e_{i,j,k}$, \Pref{ESC4}$_{i,j}$ implies $\sigmae(e_{i,j,k})=b_2$, so $\valu_{\sigmae}^\P(e_{i,j,k})=\{e_{i,j,k}\}\cup\valu_{\sigmae}^\P(b_2)$.
					Since $\neg\sigmaebar(s_{i,j})$ by the choice of $j$ and \Pref{USV1}$_i$, $\relbit{\sigmae}=1$, $\sigmaebar(eb_{i,j})$ and $\nsigmaebar(eg_{i,j})$, \Cref{lemma: Exact Behavior Of Counterstrategy} thus implies \[\valu_{\sigmae}^\P(F_{i,j})=\{F_{i,j},d_{i,j,k'},e_{i,j,k'}\}\cup\valu_{\sigmae}^\P(b_2)\rhd\{e_{i,j,k}\}\cup\valu_{\sigmae}^\P(b_2)=\valu_{\sigmae}^\P(e_{i,j,k})\]for some $k'\in\{0,1\}$.
					Hence $(d_{i,j,k},F_{i,j})\in I_{\sigmae}$.
					
				\item Let $i\in\{\nsb+1,\dots,m-1\}$ with $\indbit_i=0$ and $j=1-\indbit_{i+1}$.
					We prove $\sigmae(d_{i,j,k})\neq F_{i,j}$ and $\valu_{\sigmae}^{\P}(F_{i,j})\rhd\valu_{\sigmae}^{\P}(e_{i,j,k})$.
					However, since $\sigmaebar(eb_{i,j})\wedge\neg\sigmaebar(eg_{i,j})$ this can be shown by the same arguments used before.
					
			\end{itemize}				
					We prove that no other edge becomes an improving switch.
					Let $(i,j)$ be a pair of indices for which the edge $(s_{i,j},b_1)$ does not become improving for $\sigmae$.
					By our assumptions on $I_{\sigma}$, it then suffices to prove that one of the following three cases is true.
					\begin{enumerate}
						\item $\sigmae(s_{i,j})=h_{i,j}$ and $\valu_{\sigmae}^{\P}(h_{i,j})\unrhd\valu_{\sigmae}^{\P}(b_1)$ or
						\item $\sigmae(s_{i,j})=b_1$ and $\tau^{\sigma}(F_{i,j}),\tau^{\sigmae}(F_{i,j})\neq s_{i,j}$ or
						\item $\sigmae(s_{i,j})=b_1, j=1-\sigmabar(g_i)$ and $\valu_{\sigmae}^{\P}(F_{i,1-j})>\valu_{\sigmae}^{\P}(F_{i,j})$.
					\end{enumerate}
					We distinguish the following cases:
					\begin{itemize}
						\item \boldall{$i\leq\nsb-1$ and $j\in\{0,1\}$:} Then, the statement follows by the same arguments used for the corresponding cases for $\indbit=2^k,k\in\mathbb{N}$.
						\item \boldall{$i=\nsb$ and $j=\indbit_{\nsb+1}$:} Then, $\sigmae(s_{i,j})=h_{i,j}$ by \Pref{USV1}$_{\nsb}$.
							Hence, by \Pref{EV1}$_{i+1}$ and since $\{h_{\nsb,j}\}\rhd L_{1,\nsb}$, \begin{align*}
								\valustar_{\sigmae}^\P(h_{i,j})&=\{h_{\nsb,j}\}\cup\valustar_{\sigmae}^\P(b_{\nsb+1})=\{h_{\nsb,j}\}\cup L_{\nsb+1}\rhd L_1^\P=\valustar_{\sigmae}^\P(b_1).
							\end{align*}
						\item \boldall{$i=\nsb$ and $j=1-\indbit_{\nsb+1}$:} Then, $\sigmae(s_{i,j})=b_1$ by \Pref{USV1}$_{\nsb}$ and $\sigmae(g_i)=F_{i,1-j}$ by \Pref{EV2}$_{\nsb}$.
							We prove $\valustar_{\sigmae}^{\P}(F_{i,1-j})\rhd\valustar_{\sigmae}^{\P}(F_{i,j})$.
							Note that we do not need to consider the cycle vertices here as we proved that the corresponding edges become improving for $\sigmae$.
							Since $(i,j)\in S_1$,  $\sigmae$ has \Pref{ESC4}$_{i,j}$.
							Thus, $\sigmaebar(eb_{i,j})\wedge\nsigmaebar(eg_{i,j})\wedge\relbit{\sigmae}=1$, implying $\valustar_{\sigmae}(F_{i,j})=\valustar_{\sigmae}(b_2)$.
							Since $F_{i,1-j}$ is closed by \Pref{EV1}$_{\nsb}$, $\sigmae(s_{i,1-j})=h_{i,1-j}$ by \Pref{USV1}$_i$, \Pref{EV1}$_{i+1}$ and the choice of $i$ imply \begin{align*}
								\valustar_{\sigmae}^\P(F_{i,1-j})&=\{s_{i,1-j},h_{i,1-j}\}\cup\valustar_{\sigmae}^\P(b_{i+1})\rhd W_{i}^\P\cup\valustar_{\sigmae}^\P(b_{i+1})\\
									&=\valustar_{\sigmae}^\P(b_i)=\valustar_{\sigmae}^\P(b_{\nsb})=\valustar_{\sigmae}^\P(b_2)=\valustar_{\sigmae}^{\P}(F_{i,j}).
							\end{align*}
						\item \boldall{$i\in\{\nsb+1,\dots,m-1\}, \indbit_i=0$ and $j=\indbit_{i+1}$:} Again, $\sigmae(s_{i,j})=h_{i,j}$ by \Pref{USV1}$_{i}$  in this case.
							By \Pref{EV1}$_{i+1}$ we then obtain \begin{align*}
								\valustar_{\sigmae}^\P(h_{i,j})&=\{h_{i,j}\}\cup\valustar_{\sigmae}^\P(b_{i+1})=\{h_{i,j}\}\cup L_{i+1}^\P\rhd L_1^\P=\valustar_{\sigmae}^\P(b_1).
							\end{align*}
						\item \boldall{$i\in\{\nsb+1,\dots,m-1\}, \indbit_i=0$ and $j=1-\indbit_{i+1}$:} Then, by \Pref{USV1}$_i$, we have $\sigmae(s_{i,j})=b_1$.
							In addition, $(i,j)\in S_1$ and $(i,1-j)\in S_2$, implying $\sigmabar(eb_{i,j})\wedge\nsigmabar(eg_{i,j})$ as well as $\sigmabar(eb_{i,1-j})\wedge\nsigmabar(eg_{i,1-j})$.
							By \Pref{EBG3} and since only cycle centers $F_{i',\indbit_{i'+1}}$ are closed by assumption, $\nsb>1$ implies $\sigmabar(g_1)\neq\sigmabar(b_2)$.
							Consequently, by \Cref{lemma: Exact Behavior Of Counterstrategy} and since player~$1$ always chooses the vertex minimizing the valuation, \[\valustar_{\sigmae}^\P(F_{i,j})=\valustar_{\sigmae}^\P(b_2)>\valustar_{\sigmae}^\P(g_1)=\valustar_{\sigmae}^{\P}(F_{i,1-j}).\]
							By our assumptions on $I_{\sigma}$, this implies that it holds that $\sigma(g_i)=\sigmae(g_i)=F_{i,j}$.
							We thus prove $\sigmae(b_i)=b_{i+1}$ and $\valustar_{\sigmae}^\P(b_{i+1})>\valustar_{\sigmae}^\P(g_i)$ to prove $(b_i,g_{i})\notin I_{\sigma}$ and $\sigmae(s_{i-1,1})=b_1$ and $\valustar_{\sigmae}^{\P}(b_1)>\valustar_{\sigmae}^{\P}(h_{i-1,1})$ to prove $(s_{i-1,1},h_{i-1,1})\notin I_{\sigmae}$.
							
							First, $\sigmae(b_i)=b_{i+1}$ follows from \Pref{EV1}$_i$ whereas $\sigmae(s_{i-1,1})=b_1$ follows from $i-1\geq\nsb$ and \Pref{USV1}$_{i-1}$.
							Since we need to analyze $\valustar_{\sigmae}^\P(g_i)$ using \Cref{corollary: Complete Valuation Of Selection Vertices PG}, we determine $\lambda_i^\P$.
							However, since $\sigmae(s_{i,j})=b_1$ and $\sigmaebar(g_i)=j$, this lemma implies $\valustar_{\sigmae}^{\P}(g_i)=\{g_i\}\cup\valustar_{\sigmae}^{\P}(g_1)$.
							Since the conditions of the third case of \Cref{lemma: Valuation of g1} are fulfilled (by \Pref{BR1} applied to $\sigma$, \Pref{EV1}$_{i'}$ for $i'\leq\nsb$ and our assumption), \begin{align*}
								\valustar_{\sigmae}^\P(g_i)&=\{g_i\}\cup\valustar_{\sigmae}^\P(g_1)=\{g_i\}\cup\bigcup_{i'=1}^{\nsb-1}W_{i'}^\P\cup\valustar_{\sigma}^\P(b_{\nsb+1})\\
									&=\{g_i\}\cup\bigcup_{i'=1}^{\nsb-1}W_{i'}^\P\cup L_{\nsb+1,i-1}^\P+L_{i+1}^\P\lhd L_{i+1}^\P=\valustar_{\sigmae}^\P(b_{i+1}), \\
								\valustar_{\sigmae}^\P(h_{i-1,1})&=\{h_{i-1,1}, g_i\}\hspace*{-1pt}\cup\hspace*{-1pt}\valustar_{\sigmae}^\P(g_1)\hspace*{-1pt}\lhd\hspace*{-1pt}\valustar_{\sigmae}^\P(g_1)\hspace*{-1pt}\lhd\hspace*{-1pt}\valustar_{\sigmae}^\P(b_2)=\valustar_{\sigmae}^\P(b_1).
							\end{align*}
						\item \boldall{$i\in\{\nsb+1,\dots,m-1\}, \indbit_i=1$ and $j=\indbit_{i+1}$:} As before, $\sigmae(s_{i,j})=h_{i,j}$ by \Pref{USV1}$_{i}$.
							The statement thus follows by the same arguments used before.
						\item \boldall{$i\in\{\nsb+1,\dots,m-1\}, \indbit_i=1$ and $j=1-\indbit_{i+1}$:} Then, $\sigmae(s_{i,j})=b_1$ by \Pref{USV1}$_{i}$.
							We prove $\tau^{\sigma}(F_{i,j}),\tau^{\sigmae}(F_{i,j})\neq s_{i,j}$.
							By \Pref{ESC5}$_{i,j}$, both $\sigmaebar(eg_{i,j})\wedge\sigmaebar(eb_{i,j})$ and $\sigmabar(eg_{i,j})\wedge\sigmabar(eb_{i,j})$ hold.
							Hence, by \Cref{lemma: Exact Behavior Of Counterstrategy}, $\tau^{\sigma}(F_{i,j}), \tau^{\sigmae}(F_{i,j})\neq s_{i,j}$.
						\item \boldall{$i\geq m$ and $j=\indbit_{i+1}$:} 
							By the choice of $i$, we then have $\indbit_i=0$.
							For $i\neq n$, the statements follows similar to the last cases.
							For $i=n$, we have \[\valustar_{\sigmae}^{\P}(h_{i,0})=\{h_{n,0}\}\rhd \bigcup_{i'\geq 1}\{W_{i'}^{\P}\colon\sigmae(b_{i'})=g_{i'}\}=L_1^{\P}=\valustar_{\sigmae}^{\P}(b_1).\]
						\item \boldall{$i\geq m$ and $j=1-\indbit_{i+1}$:} Then, it holds that $\sigmae(s_{i,j})=b_1$ by \Pref{USV1}$_i$.
							Hence we need to show $\tau^{\sigma}(F_{i,j}),\tau^{\sigmae}(F_{i,j})\neq s_{i,j}$.
							However, this follows immediately from \Cref{lemma: Exact Behavior Of Counterstrategy} since \Pref{ESC5} implies $\sigmabar(eb_{i,j})\wedge\sigmabar(eg_{i,j})$ as well as $\sigmaebar(eb_{i,j})\wedge\sigmaebar(eg_{i,j})$. \qedhere
					\end{itemize}
	\end{enumerate}	
\end{proof}

\PhaseThreeToFiveMDP*

\begin{proof}
We begin by proving that $\sigmae$ is a phase-$5$-strategy.
Since $\indbit^{\sigmae}=\indbit^{\sigma}=\bit+1\eqqcolon\indbit$ and $\nsb>1$, $\sigmae$ has Properties (\ref{property: EV1})$_i$, (\ref{property: EV2})$_i$ and (\ref{property: EV3})$_i$ for all $i\in[n]$.
Also, $\sigmae$ does not have \Pref{ESC1} as it has \Pref{ESC5}$_{i,j}$ for all $(i,j)\in S_2$ and $S_2\neq\emptyset$.
Therefore, as $\sigmae$ has \Pref{USV1}$_i$ for all $i\in[n]$ by assumption, it is a phase-$5$-strategy for $\bit$.
We next prove that $\sigmae$ is well-behaved.
Since $\relbit{\sigma}\neq 1$ but $\relbit{\sigmae}=1$ as $\incorrect{\sigmae}=\emptyset$ due to the choice of $e$, we need to reevaluate the following properties.
\begin{enumerate}[align=right, leftmargin=1.75cm]
	\item[(\ref{property: S1})] By Properties (\ref{property: USV1})$_i$ and (\ref{property: EV2})$_i$, $\sigmae(b_i)=g_i$ implies $\sigmaebar(s_i)$ for all $i\geq 1$.
	\item[(\ref{property: D1})] By Properties (\ref{property: EV1})$_i$ and (\ref{property: EV2})$_i$, $\sigmae(b_i)=g_i$ implies $\sigmaebar(d_i)$.
	\item[(\ref{property: MNS2})] Assume there was some $i<\minnegsige{g}<\minnegsige{s},\minsige{b}$.
		Then $1<\minnegsige{g}$, implying $\sigmae(g_1)=F_{1,1}$.
		By the choice of $i$, it holds that $\minsige{b}\geq 3$, hence $\sigmae(b_2)=b_3$.
		But then, \Pref{USV1}$_1$ implies $\sigmae(s_{1,\sigmaebar(g_1)})=\sigmae(s_{1,1})=b_1$, contradicting $\minnegsige{g}<\minnegsige{s}$.
	\item[(\ref{property: MNS3})] Assume there was some $i<\minnegsige{s}\leq\minnegsige{g}<\minsige{b}$.
		Then $1<\minnegsige{s}\leq\minnegsige{g}$, implying $\sigmae(g_1)=F_{1,1}$ and $\sigmae(s_{1,1})=h_{1,1}$.
		Hence, by \Pref{USV1}$_1$, it holds that $\sigmae(b_2)=g_2$ and thus $\minsige{b}=2$.
		But this is a contradiction as the premise implies $\minsige{b}>3$.
	\item[(\ref{property: MNS4})] If $\minnegsige{s}>1$, then the same arguments used for \Pref{MNS3} can be used again.
		Hence let $\minnegsige{s}=1$.
		Then, $\sigmae(s_{1,\sigmabar(g_1)})=b_1$.
		In particular, \Pref{USV1}$_1$ implies $\sigmabar(g_1)\neq\indbit_{2}$, hence $\sigmaebar(g_1)=1-\indbit^{\sigmae}_2$.
		But then, by \Pref{ESC4}$_{1,1-\indbit_2}$, we have $\sigmaebar(eb_{\minnegsige{s}})\wedge\nsigmaebar(eg_{\minnegsige{s}})$.
	\item[(\ref{property: MNS5})] Assuming that there was some $i<\minnegsige{s}<\minsige{b}\leq\minnegsige{b}$ yields the same contradiction devised for \Pref{MNS3}.
	\item[(\ref{property: MNS6})] If $\minnegsige{s}>1$, then the same arguments used for \Cref{property: MNS3} can be used to show that the premise cannot hold.
		Hence assume $\minnegsige{s}=1$.
		But then, the same arguments used for \Pref{MNS4} can be used to prove the statement.
	\item[(EG*)] It is easy to verify that each cycle center is either closed, escapes only to $b_2$ or to both $b_2$ and $g_1$.
		In particular, there is no cycle center $F_{i,j}$ with $\sigmabar(eg_{i,j})\wedge\nsigmabar(eb_{i,j})$.
	\item[(\ref{property: DN1})] By \Pref{EV1}$_n$, $\sigmaebar(d_n)$ implies $\sigmae(b_n)=g_n$.
	\item[(\ref{property: DN2})] We only need to consider this assumption if $\minnegsige{g}=n$.
		Since $\relbit{\sigma}\neq 1$, this implies $\sigmaebar(g_i)=\sigmabar(g_i)=1$ for all $i\in\{1,\dots,n-1\}$ by \Pref{BR1} (applied to $\sigma$).
		Thus, by assumption, $i\neq\relbit{\sigma}-1$ for all of those $i$, hence $n=\relbit{\sigma}-1$.
		But this implies $\relbit{\sigma}=n+1$, contradicting \Pref{REL2} for $\sigma$.
\end{enumerate}

It remains to prove the statement regarding the improving switches.
We observe that $\valu_{\sigmae}^\M(g_1)<\valustar_{\sigmae}^\M(b_2)$ since $\sigma(e_{i,j,k})=g_1$ implies $(e_{i,j,k},b_2)\in I_{\sigma}$.

Let $i<\nsb, j\coloneqq 1-\indbit_{i+1}$ and $k\in\{0,1\}$.
We prove  $(d_{i,j,k},F_{i,j})\in I_{\sigmae}$.
By assumption, the cycle center $F_{i,j}$ is open, so in particular $\sigmae(d_{i,j,k})\neq F_{i,j}$.
By the choice of $i$ and~$j$, \Pref{ESC4}$_{i,j}$ implies $\sigmaebar(eb_{i,j})\wedge\nsigmaebar(eg_{i,j})$.
Consequently, by \Cref{lemma: Exact Behavior Of Random Vertex}, $\valu_{\sigmae}^\M(F_{i,j})=(1-\e)\valu_{\sigmae}^\M(b_2)+\e\cdot\valu_{\sigmae}^\M(s_{i,j})$.
Since $\sigmae(b_1)=b_2$, the choice of $j$ and \Pref{USV1}$_i$ imply $\valu_{\sigmae}^\M(s_{i,j})=\rew{s_{i,j}}+\valu_{\sigmae}^{\M}(b_2)$.
Thus, $(d_{i,j,k},F_{i,j})\in I_{\sigmae}$ follows from  $\valu_{\sigmae}^\M(F_{i,j})=(1-\e)\valu_{\sigmae}^\M(b_2)+\e\valu_{\sigmae}^\M(s_{i,j})>\valu_{\sigmae}^{\M}(b_2)=\valu_{\sigmae}^\M(e_{i,j,k}).$

We  prove that $X_0,X_1$ are improving for $\sigmae$ if $\indbit$ is not a power of two.
Fix $k\in\{0,1\}$ and let $i\coloneqq\nsb, j\coloneqq 1-\indbit_{\nsb+1}$.
We begin by proving $(d_{i,j,k},F_{i,j})\in I_{\sigmae}$.
By the choice of $j$ and our assumptions, $\sigmae(d_{i,j,k})\neq F_{i,j}$.
In addition \Pref{ESC4}$_{i,j}$ implies that $\sigmaebar(eb_{i,j})\wedge\nsigmaebar(eg_{i,j})$.
Since this implies $\valu_{\sigmae}(F_{i,j})=(1-\e)\valu_{\sigmae}^\M(b_2)+\e\valu_{\sigmae}^\M(s_{i,j})$ as well as $\valu_{\sigmae}^\M(e_{i,j,k})=\valu_{\sigmae}^\M(b_2)$, it suffices to prove $\valu_{\sigmae}^\M(s_{i,j})>\valu_{\sigmae}^\M(b_2)$.
This however follows directly since \Pref{USV1}$_i$ and the choice of $j$ imply $\valu_{\sigmae}^\M(s_{i,j})=\rew{s_{i,j}}\cup\valu_{\sigmae}^\M(b_2)$.
By applying the same arguments, we also obtain $(d_{i,j,k},F_{i,j})\in I_{\sigmae}$ for $i\in\{\nsb+1,\dots,m-1\}$ with $\indbit_i=0$ and $j=1-\indbit_{i+1}$ as $(i,j)\in S_1$ for these indices.

We now prove that no further improving switch is created.
Note that no additional improving switches $(d_{i,j,k},F_{i,j})$ but the ones discussed earlier are created in any case.
The reason is that the only indices $(i,j)$ with $(i,j)\in S_1$ are $i<\nsb$ and $j=1-\indbit_{i+1}$ if $\indbit$ is a power of 2.
All other indices $(i,j)$ are contained in $S_2$ since $\nsb=m$.
Consequently, $\valustar_{\sigmae}^{\M}(F_{i,j})=\frac{1}{2}\valustar_{\sigmae}^\M(b_2)+\valustar_{\sigmae}^\M(g_1)<\valustar_{\sigmae}^\M(b_2)$.
By the same argument, no further improving switch $(d_{i,j,k},F_{i,j})$ besides the ones discussed earlier is created for the case that $\indbit$ is not a power of 2.

The application of $e$ increases the valuation of the vertex $b_1$.
The only vertices that have an edge towards $b_1$ are upper selection vertices $s_{i,j}$.
As we fully covered the cycle vertices, it now suffices to prove that the following statements hold:
\begin{enumerate}
	\item If $\sigmae(s_{i,j})=h_{i,j}$, then $(s_{i,j},b_1)\notin I_{\sigma}, I_{\sigmae}$.
	\item If $\sigmae(s_{i,j})=b_1$ and $\sigmaebar(g_i)\neq j$, then $(g_i,F_{i,j})\notin I_{\sigma}, I_{\sigmae}$.
	\item If $\sigmae(s_{i,j})=b_1$ and $\sigmaebar(g_i)=j$, then $\valu_{\sigmae}^\M(g_i)-\valu_{\sigma}^\M(g_i)\in o(1)$.
\end{enumerate}

Consider the case $\sigmae(s_{i,j})=h_{i,j}$ first.
Then, by \Pref{USV1}$_i$, $j=\indbit_{i+1}$.
Consequently, since $\sigmae(b_1)=b_2$ and since $\rew{h_{i,j}}>\sum_{\ell\in[i]}W_{\ell}^\M$, \Pref{EV1}$_{i+1}$ yields \begin{align*}
\valustar_{\sigmae}^\M(h_{i,j})&=\rew{h_{i,j}}+\valustar_{\sigmae}^\M(b_{i+1})=\rew{h_{i,j}}+L_{i+1}^\M> L_{1,i}^\M+L_{i+1}^\M=L_1^\M=\valustar_{\sigmae}^\M(b_1).
\end{align*}
Since $\valustar_{\sigma}(h_{i,j})=\valustar_{\sigmae}(h_{i,j})$ and $\valustar_{\sigmae}^\M(b_1)>\valustar_{\sigma}^\M(b_1)$, this implies $(s_{i,j},b_1)\notin I_{\sigma},I_{\sigmae}$.
Now consider the case $\sigmae(s_{i,j})=b_1$ and $\sigmaebar(g_i)\neq j$.
Then, by \Pref{USV1}$_{i}$, it holds that $j=1-\indbit_{i+1}$ and thus, by \Pref{EV1}$_{i+1}$ \begin{align*}
	\valustar_{\sigmae}^\M(s_{i,j})&=\rew{s_{i,j}}+\valustar_{\sigmae}^\M(b_1)=\rew{s_{i,j}}+L_{\nsb}^{\M},\\
	\valustar_{\sigmae}^\M(h_{i,1-j})&=\rew{h_{i,1-j}}+\valustar_{\sigmae}^\M(b_{i+1})=\rew{h_{i,1-j}}+L_{i+1}^\M>\rew{s_{i,j}}+L_{\nsb}^{\M}.	
\end{align*}
We prove that this implies $\valu_{\sigmae}^\M(F_{i,1-j})>\valu_{\sigmae}^\M(F_{i,j})$ in any case.
Let $F_{i,1-j}$ be closed.
Then, $\indbit_i=1\wedge\indbit_{i+1}=1-j$.
Consequently, \[\valustar_{\sigmae}^\M(F_{i,1-j})=\rew{s_{i,1-j},h_{i,1-j}}+\valustar_{\sigmae}^\M(b_{i+1})=\rew{s_{i,1-j},h_{i,1-j}}+L_{i+1}^\M>L_2^\M=\valustar_{\sigmae}^\M(b_2).\]
Since either $\valustar_{\sigmae}^\M(F_{i,j})=\valustar_{\sigmae}^\M(b_2)$ or $\valustar_{\sigmae}^\M(F_{i,j})=\frac{1}{2}\valustar_{\sigmae}(b_2)+\frac{1}{2}\valustar_{\sigmae}^\M(g_1)$ and $\valustar_{\sigmae}^\M(g_1)<\valustar_{\sigmae}^\M(b_2)$, this implies the statement.
Thus assume that $F_{i,1-j}$ is not closed, implying $\indbit_{i}=0$.
For the sake of a contradiction, assume $i<\nsb$.
Then, $\sigmaebar(g_i)=1-j=\indbit_{i+1}$.
However, since $\nsb=\relbit{\sigma}$, applying \Pref{BR1} to $\sigma$ implies $\sigmabar(g_i)=\sigmaebar(g_i)=1-\indbit_{i+1}$ which is a contradiction.
Since $F_{i,1-j}$ is not closed, it suffices to consider the case $i>\nsb$.
If $i<m$, then $(i,1-j)\in S_1$ and $(i,j)\in S_2$.
Then, by Properties (\ref{property: ESC4})$_{i,1-j}$, (\ref{property: ESC5})$_{i,j}$ and (\ref{property: EV1})$_{i+1}$, we have \[\valustar_{\sigmae}^\M(F_{i,1-j})=\valustar_{\sigmae}^\M(b_2)\quad\text{and}\quad\valustar_{\sigmae}^\M(F_{i,j})=\frac{1}{2}\valustar_{\sigmae}^\M(b_2)+\frac{1}{2}\valustar_{\sigmae}^\M(g_1),\] implying the statement.
If $i>m$, then $(i,j),(i,1-j)\in S_2$ and the statement follows from $\valustar_{\sigmae}^\M(h_{i,1-j})>\valustar_{\sigmae}^\M(s_{i,j})$.
As it holds that $\valustar_{\sigma}^\M(s_{i,j})<\valustar_{\sigmae}^\M(s_{i,j})$ and $\valustar_{\sigma}^\M(h_{i,1-j})=\valustar_{\sigmae}^\M(h_{i,1-j})$ we thus have $(g_i,F_{i,j})\notin I_{\sigma},I_{\sigmae}$.

Finally, assume $\sigmae(s_{i,j})=b_1$ and $\sigmaebar(g_i)=j$.
Since $\valu_{\sigmae}^\M(g_i)-\valu_{\sigma}^\M(g_i)\geq 0$, we prove that this difference is smaller than 1.
By \Pref{USV1}$_i$, $j=1-\indbit_{i+1}$.
By the assumptions of the lemma, this implies that $F_{i,j}$ is neither closed with respect to $\sigma$ nor to $\sigmae$.
Consequently, either $\valu_{\sigmae}^\M(g_i)-\valu_{\sigma}^\M(g_i)=\e[\valu_{\sigmae}^\M(s_{i,j})-\valu_{\sigma}^\M(s_{i,j})]$ or $\valu_{\sigmae}^\M(g_i)-\valu_{\sigma}^\M(g_i)=\frac{2\e}{1+\e}[\valu_{\sigmae}^\M(s_{i,j})-\valu_{\sigma}^\M(s_{i,j})]$.
In either case, the difference is smaller than 1 by the choice of $\e$.
\end{proof}

\PhaseThreeToFive*

\begin{proof}
	As usual, the choice of $e$ implies $\indbit^{\sigma}=\indbit^{\sigmae}=\bit+1\eqqcolon\indbit$.
	
	We begin by proving $\relbit{\sigmae}=u$.
	Since $\sigma$ is a phase-$3$-strategy for $\bit$, it has \Pref{EV1}$_i$ and \Pref{EV2}$_i$ for all $i>1$.
	This implies $i\notin\incorrect{\sigma}$ and thus $i\notin\incorrect{\sigmae}$ for all $i>1$.
	Since $\sigmae(b_1)=g_1$, it suffices to show $\sigmaebar(g_i)=\sigmaebar(b_{i+1})$, implying $\incorrect{\sigmae}=\emptyset$ and $\relbit{\sigmae}=u$.
	This however follows directly as~$\sigma$ has \Pref{CC2}.
	
	We now show that $\sigmae$ is a phase-$5$-strategy for $\bit$.
	Since $\sigma$ is a phase-$3$-strategy for $\bit$ and $e=(b_1,b_2)$, it suffices to show that the $\sigma$ has Properties (\ref{property: EV1})$_1$, (\ref{property: EV2})$_1$, (\ref{property: EV3})$_1$ and (\ref{property: CC2}).
	Note that $\nsb=1$ implies $S_3\neq\emptyset$, implying that $\sigmae$ does not have \Pref{ESC1}.
	By definition, $1=\indbit_1=\sigmaebar(d_{1,\indbit_2})$.
	Thus, since $\sigmae(b_1)=g_1$, $\sigmae$ has \Pref{EV1}$_1$.
	It also has \Pref{EV2}$_1$ and \Pref{CC2} since $\sigma$ has \Pref{CC2}.
	Since $\nsigmabar(d_{1,1-\indbit_2})$ by assumption, also $\nsigmaebar(d_{1,1-\indbit_2})$, hence $\sigmae$ has \Pref{EV3}$_1$.
	Thus, $\sigmae$ is a phase-$5$-strategy for $\bit$.

	We prove that $\sigmae$ is well-behaved.
	Since $e=(b_1,g_1), \relbit{\sigma}=1$ and $\relbit{\sigmae}=u>1$, it suffices to investigate the following properties.
	\begin{itemize}[align=right, leftmargin=1.75cm]
		\item[(\ref{property: S2})] Let $i<\relbit{\sigmae}$.
			Then, since $\relbit{\sigmae}=u$, we have $\sigmae(b_i)=g_i$.
			Consequently, by \Pref{EV1}$_i$, \Pref{EV2}$_i$ and \Pref{USV1}$_i$, $\sigmaebar(d_i)$ and $\sigmaebar(s_i)$.
		\item[(\ref{property: B1})] As $\relbit{\sigmae}=u$, the premise can never be correct as $i<\relbit{\sigmae}-1$ implies $\sigmae(b_i)=g_i$.
		\item[(\ref{property: B2})] This again holds since $\relbit{\sigmae}=u$.
		\item[(\ref{property: BR1})] Let $i\coloneqq\relbit{\sigmae}-1$.
			Then $\sigmae(b_i)=g_i$ and $\sigmae(b_{i+1})=b_{i+2}$.
			Thus, by Properties (\ref{property: EV1})$_i$ and (\ref{property: EV1})$_{i+1}$ as well as \Pref{EV2}$_i$, we have $\sigmae(g_i)=F_{i,0}$.
			For $i<\relbit{\sigmae}-1$, we have $\sigmae(g_i)=F_{i,1}$ as we then have $\sigmae(b_{i+1})=g_{i+1}$.
		\item[(\ref{property: BR2})] Since $i<\relbit{\sigmae}$ implies $\sigmae(b_i)=g_i$, \Pref{EV1}$_i$ implies $\sigmaebar(d_i)$ and thus $\nsigmaebar(eg_i)$.
		\item[(\ref{property: D2})] This follows by the same argument used for \Pref{BR2}.
		\item[(\ref{property: EG5})] By \Pref{USV1}$_i$, $\sigmaebar(s_{i,j})$ implies $\sigmaebar(b_{i+1})=j$.
		\item[(EB*)] Any pair of indices $i\in[n], j\in\{0,1\}$ either fulfills \Pref{ESC5}$_{i,j}$, \Pref{ESC3}$_{i,j}$ or $\sigmaebar(d_{i,j})$.
			Hence, there are no indices such that $\sigmaebar(eg_{i,j})\wedge\neg\sigmaebar(eb_{i,j})$, so the premise of any of the Properties (EB*) is always incorrect.
		\item[(\ref{property: EBG4})] Since $\relbit{\sigmae}>1$ implies $\sigmae(b_1)=g_1$, it is impossible that both $\sigmae(g_1)=F_{1,0}$ and $\sigmae(b_2)=g_2$.
		\item[(\ref{property: EBG5})] Since $\relbit{\sigmae}>1$ implies $\sigmae(b_1)=g_1$, it is impossible that both $\sigmae(g_1)=F_{1,1}$ and $\sigmae(b_2)=b_3$.
	\end{itemize}
	It remains to show that \[I_{\sigmae}=(I_{\sigma}\setminus\{e\})\cup\bigcup_{\substack{i=u+1\\\indbit_i=0}}^{m-1}\{(d_{i,1-\indbit_{i+1},0},F_{i,1-\indbit_{i+1}}), (d_{i,1-\indbit_{i+1},1},F_{i,1-\indbit_{i+1}})\}.\]
	Let $i\in\{u+1,\dots,m-1\}, \indbit_i=0, j=1-\indbit_{i+1}$ and $k\in\{0,1\}$.
	We prove $(d_{i,j,k},F_{i,j})\in I_{\sigmae}$.
	By \Pref{ESC3}$_{i,j}$, it holds that $\sigmaebar(eg_{i,j})\wedge\neg\sigmaebar(eb_{i,j})$.
	In addition, $\sigma(d_{i,j,k})\neq F_{i,j}$.
	It thus suffices to show $\valu_{\sigmae}^*(e_{i,j,k})\prec\valu_{\sigmae}^*(F_{i,j})$.

	Consider the case $G_n=S_n$ first.
	Since $\sigmae(d_{i,j,k})=e_{i,j,k}$, \Pref{ESC3}$_{i,j}$ implies $\valu_{\sigmae}^\P(e_{i,j,k})=\{e_{i,j,k}\}\cup\valu_{\sigmae}^\P(g_1)$.
	Now, by \Cref{lemma: Exact Behavior Of Counterstrategy}, we obtain \[\valu_{\sigmae}^\P(F_{i,j})=\{F_{i,j},d_{i,j,k'},e_{i,j,k'}\}\cup\valu_{\sigmae}^\P(g_1)\] for some $k'\in\{0,1\}$ as $\relbit{\sigmae}\neq 1$.
	This however implies the statement for $G_n=S_n$ as the priority of $F_{i,j}$ is even and larger than the priorities of both $d_{i,j,k'},e_{i,j,k'}$.

	Consider the case $G_n=M_n$.
	Then $\valu_{\sigmae}^\M(F_{i,j})=(1-\e)\valu_{\sigmae}^\M(g_1)+\e\valu_{\sigmae}^\M(s_{i,j})$, it therefore suffices to prove $\valu_{\sigmae}^\M(s_{i,j})>\valu_{\sigmae}^\M(g_1)$.
	This however follows directly as \Pref{USV1}$_i$, the choice of $j$ and $\sigmae(b_1)=g_1$ imply $\valu_{\sigmae}^\M(s_{i,j})=\rew{s_{i,j}}+\valu_{\sigmae}^\M(g_1)$.
	
	It remains to show that no other improving switch is created and that switches that are improving for $\sigma$ are improving for $\sigmae$ (with the exception of $e$).
	By applying $e=(b_1,g_1)$, the valuation of $b_1$ increases.
	The only vertices that have an edge towards $b_1$ are the vertices $s_{i,j}, i\in[n], j\in\{0,1\}$.
	To show that no other improving switch is created and that switches that are improving for $\sigma$ are also improving for $\sigmae$, it suffices to show that one of the following holds for all $i\in[n], j\in\{0,1\}$ not considered earlier:
	\begin{enumerate}
		\item $\sigmae(s_{i,j})=h_{i,j}$ implies $\valu_{\sigmae}^*(h_{i,j})\succ\valu_{\sigmae}^*(b_1)$.
		\item $\sigmae(s_{i,j})=b_1$ and $(i,j)\notin S_4$ implies $(d_{i,j,k},F_{i,j})\in I_{\sigma}\Leftrightarrow(d_{i,j,k},F_{i,j})\in I_{\sigmae}$.
		\item If $G_n=S_n$, then $\sigmae(s_{i,j})=b_1$ implies either $\tau^{\sigma}(F_{i,j})=\tau^{\sigmae}(F_{i,j})\neq s_{i,j}$ or $\sigmae(g_i)=F_{i,j}$ and $\sigmae(b_i)=b_{i+1}\wedge(b_i,g_i)\notin I_{\sigmae}$ as well as $\sigmae(s_{i-1,1})=b_1\wedge(s_{i-1,1},h_{i-1,1})\notin I_{\sigmae}$.
		\item If $G_n=M_n$, then $\sigmae(s_{i,j})=b_1$ and $\sigmaebar(g_i)=1-j$ implies $\valu_{\sigmae}^\M(F_{i,1-j})>\valu_{\sigmae}^\M(F_{i,j})$.
		\item If $G_n=M_n$, then $\sigmae(s_{i,j})=b_1$ and $\sigmaebar(g_i)=j$ implies $\valu_{\sigmae}^\M(g_i)-\valu_{\sigma}^\M(g_i)\in(0,1)$.
	\end{enumerate}
	We now prove these statements one after another.
	\begin{enumerate}
		\item Fix indices $i,j$ with $\sigmae(s_{i,j})=h_{i,j}$.
			Then, by \Pref{USV1}$_i$, $j=\indbit_{i+1}$.
			Consequently, by \Pref{EV1}$_{i+1}$, $\valustar_{\sigmae}^*(h_{i,j})=\ubracket{h_{i,j}}\oplus\valustar_{\sigmae}^*(b_{i+1})$.
			Since $\sigmae$ has \Pref{EV1}$_{i'}$ and \Pref{EV2}$_{i'}$ for all $i'<\relbit{\sigmae}$, there is no $i'<\relbit{\sigmae}$ with $\nsigmaebar(d_{i'})$.
			Consequently, by \Cref{lemma: Valuation of g if level small}, $\valustar_{\sigmae}^*(g_1)=R_1^*$ in any case.
			Now, since $\ubracket{h_{i,j}}\succ\bigoplus_{i'\leq i}W_{i'}^*$, this implies $\valustar_{\sigmae}^*(h_{i,j})\succ\valustar_{\sigmae}^*(g_1)$ for both possible cases $\valustar_{\sigmae}^*(b_{i+1})=L_{i+1}^*$ and $\valustar_{\sigmae}^*(b_{i+1})=R_{i+1}^*$.
		\item Consider some edge $(d_{i,j,k},F_{i,j})$ for which we did not prove that $(d_{i,j,k},F_{i,j})\in I_{\sigmae}$, i.e., assume $(i,j)\notin S_4$.
			We show that $(d_{i,j,k},F_{i,j})\in I_{\sigma}\Leftrightarrow(d_{i,j,k},F_{i,j})\in I_{\sigmae}$.
			Let $(d_{i,j,k},F_{i,j})\in I_{\sigma}$.
			Then, by our assumptions on $I_{\sigma}$, $\sigma(e_{i,j,k})=b_2$.
			This implies that $\sigmabar(eg_{i,j})\wedge\sigmabar(eb_{i,j})$ has to hold and, due to $(d_{i,j,k},F_{i,j})\in I_{\sigma}$ and $\sigmabar(g_1)=\sigmabar(b_2)$, \begin{align*}
				\valustar_{\sigma}^\M(F_{i,j})&=\frac{1}{2}\valustar_{\sigma}^\M(g_1)+\frac{1}{2}\valustar_{\sigma}^\M(b_2)>\valustar_{\sigma}^\M(b_2)=\valustar_{\sigma}^\M(e_{i,j,k}), \\
				\valu_{\sigma}^\P(F_{i,j})&=\{F_{i,j},d_{i,j,k'},e_{i,j,k'}\}\cup\valu_{\sigma}^\P(b_2)\rhd\{e_{i,j,k'}\}\cup\valu_{\sigma}^\P(b_2)=\valu_{\sigma}^\P(e_{i,j,k}).
			\end{align*}
			But then, $\valustar_{\sigma}^*(e_{i,j,k})=\valustar_{\sigmae}^*(e_{i,j,k})$ and $\valustar_{\sigmae}^*(g_1)\succeq\valustar_{\sigma}^*(g_1)$, the same estimation holds for $\sigmae$, implying $(d_{i,j,k},F_{i,j})\in I_{\sigmae}$.
			Now let $(d_{i,j,k},F_{i,j})\notin I_{\sigma}$, implying $\sigma(e_{i,j,k})=\sigmae(e_{i,j,k})=g_1$.
			If $\sigma(d_{i,j,k})=F_{i,j}$, then there is nothing to show hence assume $\sigma(d_{i,j,k})=e_{i,j,k}$.
			This implies $\sigmabar(eg_{i,j})$ and $\sigmaebar(eg_{i,j})$.
			Assume $\sigmabar(eb_{i,j})$.
			Then, using the same estimations used for the case $(d_{i,j,k},F_{i,j})\in I_{\sigma}$, we can show $(d_{i,j,k},F_{i,j})\notin I_{\sigmae}$.
			Thus assume $\nsigmabar(eb_{i,j})$.
			Then, $\sigmabar(eg_{i,j})\wedge\nsigmabar(eb_{i,j})$, implying $(i,j)\in S_4$.
			This however contradicts our choice of $i$ and $j$, proving the statement.		
		\item Let $G_n=S_n$ and $\sigma(s_{i,j})=\sigmae(s_{i,j})=b_1$.
			Then, by \Pref{USV1}$_i$, $j=1-\indbit_{i+1}$.
			By our assumptions, $F_{i,j}$ is thus either mixed or $g_1$-open.
			Consider the case that it is mixed first.
			Then, by \Pref{EV2}$_1$, $\sigmaebar(g_1)=\sigmaebar(b_2)$, implying $\sigmabar(g_1)=\sigmabar(b_2)$ by the choice of $e$.
			Consequently, by \Cref{lemma: Exact Behavior Of Counterstrategy}, $\valustar_{\sigma}^\P(F_{i,j})=\valustar_{\sigma}^\P(b_2)$ and $\valustar_{\sigmae}^\P(F_{i,j})=\valustar_{\sigmae}^\P(b_2)$.
			Thus, $\tau^{\sigma}(F_{i,j}),\tau^{\sigmae}(F_{i,j})\neq s_{i,j}$.

			Now assume that $F_{i,j}$ is $g_1$-open.
			Then, by our assumptions on the cycle centers, $(i,j)\in S_4$, implying $i\in\{u+1,\dots,m-1\}$ with $\indbit_i=0$ and $j=1-\indbit_{i+1}$.
			We prove $\valu_{\sigma}^\P(F_{i,j})\rhd\valu_{\sigma}^\P(F_{i,1-j})$, implying $\sigma(g_i)=\sigmae(g_i)=F_{i,j}$ by our assumptions on $I_{\sigma}$.
			We then prove that $\sigmae(b_i)=b_{i+1}$ and $(b_i,g_i)\notin I_{\sigmae}$ as well as $\sigmae(s_{i-1,1})=b_1$ and $(s_{i-1,1},h_{i-1,1})\notin I_{\sigmae}$ (if $i>1$), implying that the valuation of no further vertex changes, proving the statement.	
			
			Since $1-j=\indbit_{i+1}$, we have $(i,1-j)\in S_3$ by assumption.
			Consequently, by \Pref{ESC5}$_{i,j}$, it holds that $\sigmabar(eg_{i,1-j})\wedge\sigmabar(eb_{i,1-j})$.
			As pointed out earlier, this implies $\valustar_{\sigma}^\P(F_{i,1-j})=\valustar_{\sigma}^\P(b_2)$.
			Since $\relbit{\sigma}=1$, \Cref{lemma: Exact Behavior Of Counterstrategy} yields $\valustar_{\sigma}^\P(F_{i,j})=\{s_{i,j}\}\cup\valustar_{\sigmae}^\P(b_2)$.
			This implies $\valustar_{\sigma}^\P(F_{i,j})\rhd\valustar_{\sigmae}^\P(F_{i,1-j})$.
			We hence need to have $\sigma(g_i)=\sigmae(g_i)=F_{i,j}$ by our assumptions on $I_{\sigma}$.
			
			By $\indbit_{i}=0$, \Pref{EV1}$_{i}$ implies $\sigmae(b_i)=b_{i+1}$.
			Since $\nsb=1$ implies $\indbit_1=1$, we have $i>1$.
			Thus, \Pref{USV1}$_{i-1}$ implies $\sigmae(s_{i-1,1})=b_1$.
			Consequently, $\valustar_{\sigmae}^\P(b_{i+1})=L_{i+1}^\P$.
			Since $\sigmaebar(eg_{i,j})\wedge\nsigmaebar(eb_{i,j})\wedge\relbit{\sigmae}>1$, \Cref{lemma: Exact Behavior Of Counterstrategy} implies $\valustar_{\sigmae}^\P(F_{i,j})=\valustar_{\sigmae}^\P(g_1)=R_1^\P$.
			Consequently, since $i\geq\relbit{\sigmae}$ and $\sigmae(b_i)=b_{i+1}$ imply $R_{i}^\P=L_{i}^\P=L_{i+1}^\P$, \begin{align*}
				\valustar_{\sigmae}^\P(g_i)&=\rew{g_i}+R_1^\P=\rew{g_i}+R_{1,i-1}^\P+R_i^\P\\
					&=\rew{g_i}+R_{1,i-1}^\P+L_{i+1}^\P\lhd L_{i+1}^\P=\valustar_{\sigmae}^\P(b_{i+1}).
			\end{align*}
			As it holds that $\valustar_{\sigmae}^\P(b_1)\rhd\valustar_{\sigmae}^\P(b_{i+1})$ and $\valustar_{\sigmae}^\P(h_{i-1,1})=\rew{h_{i-1,1}}\cup\valustar_{\sigmae}^\P(g_i)$, a similar estimation yields $\valustar_{\sigmae}^\P(h_{i-1,1})\lhd\valustar_{\sigmae}^\P(b_1)$.
			Consequently, $(b_i,g_i)\notin I_{\sigmae}$ as well as $(s_{i-1,1},h_{i-1,1})\notin I_{\sigmae}$.
		\item Let $G_n=M_n$ and $\sigmae(s_{i,j})=b_1$ and $\sigmaebar(g_{i})=1-j$.
			Then, by \Pref{USV1}$_i$, $j=1-\indbit_{i+1}$ and $\sigmae(s_{i,1-j})=h_{i,1-j}$.
			Assume that $F_{i,1-j}$ is closed.
			Then, by \Pref{EV1}$_{i+1}$, we have $\valustar_{\sigmae}^\M(F_{i,1-j})=\rew{s_{i,1-j},h_{i,1-j}}+\valustar_{\sigmae}^\M(b_{i+1})$.
			Since it is not possible that $F_{i,j}$ is closed, either \[\valustar_{\sigmae}^\M(F_{i,j})=\valustar_{\sigmae}^\M(g_1)\quad\text{ or }\quad\valustar_{\sigmae}^\M(F_{i,j})=\frac{1}{2}\valustar_{\sigmae}^\M(g_1)+\frac{1}{2}\valustar_{\sigmae}^\M(b_2).\]
			As $\valustar_{\sigmae}^\M(b_2)<\valustar_{\sigmae}^\M(g_1)$ and $\sum_{\ell\in[i]}W_{\ell}^\M<\rew{s_{i,1-j},h_{i,1-j}}$, this implies that we have $\valustar_{\sigmae}^\M(F_{i,1-j})>\valustar_{\sigmae}^\M(F_{i,j})$ in any case.

			Thus assume that $F_{i,1-j}$ is not closed.
			Since $1-j=\indbit_{i+1}$, it then follows that $(i,1-j)\in S_3$.
			If $(i,j)\in S_3$, then both cycle centers are in the same state.
			Since $\sigmae$ has \Pref{USV1}$_i$, \Pref{EV1}$_{i+1}$ and since $i\geq\nsb=1$, the statement thus follows from \Cref{lemma: Both CC Open For MDP}.

			Thus assume $(i,j)\in S_4$.
			Then, by \Pref{ESC5}$_{i,1-j}$ and \Pref{ESC4}$_{i,j}$, we have $\valustar_{\sigma}^\M(F_{i,1-j})=\frac{1}{2}\valustar_{\sigma}(g_1)+\frac{1}{2}\valustar_{\sigma}(b_2)$ and $\valustar_{\sigma}^\M(F_{i,j})=\valustar_{\sigma}(g_1)$.
			But then, $\valustar_{\sigma}^\M(F_{i,j})>\valustar_{\sigma}^\M(F_{i,1-j})$, implying $(g_i,F_{i,1-j})\in I_{\sigma}$, contradicting our assumption on $I_{\sigma}$.
		\item Let $G_n=M_n$ and $\sigmae(s_{i,j})=b_1$ and $\sigmaebar(g_i)=j$.
			Then, \Pref{USV1}$_i$ implies $j=1-\indbit_{i+1}$.
			In particular, $F_{i,1-j}$ is not closed.
			Thus, either $\sigmabar(eg_{i,j})\wedge\nsigmabar(eb_{i,j})$ or $\sigmabar(eg_{i,j})\wedge\sigmabar(eb_{i,j})$.
			In the first case, $\valu_{\sigma}^\M(F_{i,j})=(1-\e)\valu_{\sigma}^\M(g_1)+\e\valu_{\sigma}^\M(s_{i,j})$.
			But then, since $\valu_{\sigma}^\M(g_1)=\valu_{\sigmae}^\M(b_1)=R_1^\M$, this implies $\valu_{\sigmae}^\M(F_{i,j})-\valu_{\sigma}^\M(F_{i,j})\in(0,1)$.
			If $\sigmabar(eg_{i,j})\wedge\sigmabar(eb_{i,j})$, the statement follows analogously since $\valu_{\sigma}^\M(b_2)=\valu_{\sigmae}^\M(b_2)$.
	\end{enumerate}
\end{proof}

\PhaseFourComplete*

\begin{proof}
We first note that $\relbit{\sigma}=\relbit{\sigmae}=1$ since $\sigma$ is a phase-$4$-strategy for $\bit$ and by the choice of $e$.
Furthermore, by the choice of $e$, $\indbit^{\sigma}=\indbit^{\sigmae}=\bit+1\eqqcolon\indbit$.
We first prove that $\sigmae$ is well-behaved.
By the choice of $e$, we need to reevaluate the following properties:
\begin{itemize}[align=right, leftmargin=1.75cm]
	\item[(\ref{property: S1})] $\sigmae(b_i)=g_i$ implies $\indbit_i=1$ by \Pref{EV1}$_i$, hence $i\geq\nsb$.
	\item[(\ref{property: MNS1})] By \Cref{lemma: Config implied by Aeb}, $\minsige{b}\leq\minnegsige{s},\minnegsige{g}$ implies $\minsige{b}=\nsb=2$.
		But then, by assumption~3., $\sigmaebar(g_1)=F_{1,0}$, implying $\minnegsige{g}=1$ and thus contradicting the premise.
	\item[(\ref{property: MNS2})] By assumption~3, $\minnegsige{g}=\nsb-1$.
		By the choice of $i$ and $j$, we have $\minnegsige{s}\leq\nsb-1$.
		Thus $\minnegsige{s}\leq\minnegsige{g}$ and the premise of this property cannot be fulfilled.
	\item[(\ref{property: MNS4})] The conclusion is always true since $i'<\nsb$ implies $(i',1-\indbit_{i'+1})\in S_1$, implying that $\nsigmaebar(eb_{i'})\wedge\nsigmaebar(eg_{i'})$.
	\item[(\ref{property: MNS6})] See \Pref{MNS4}.
	\item[(\ref{property: EG3})] For every pair of indices $i,j$, either $\sigmabar(d_{i,j})$ or $\sigmabar(eb_{i,j})\wedge\nsigmabar(eg_{i,j})$ or $\sigmabar(eb_{i,j})\wedge\sigmabar(eg_{i,j})$ by either $\indbit=\bit+1$ or \Pref{ESC4}$_{i,j}$ resp. (\ref{property: ESC5})$_{i,j}$.
		Consequently, the premise is incorrect for $\sigmae$,
	\item[(\ref{property: EBG2})] By assumption~3 and \Pref{EV1}$_2$, $\sigmaebar(g_1)=\sigmabar(g_1)\neq\sigmabar(b_2)=\sigmaebar(b_2)$, hence the premise is incorrect.
\end{itemize}

We next prove that $\sigma$ is a phase-$4$-strategy if there is an index $i'<i$ such that $(s_{i',1-\indbit^{\sigma}_{i'+1}},b_1)\in I_{\sigma}$ and a phase-$5$-strategy otherwise.
By the definition of the phases, it suffices to prove that $\sigmae$ has \Pref{USV1}$_{\ell}$ for all $\ell\in[n]$ if there is no such index.
Hence assume that no such index exists and let $i'<i$ as there is nothing to prove if $i'\geq i$ and let $j'\coloneqq 1-\indbit_{i'+1}$.
Then, since $(s_{i',j'},b_1)\notin I_{\sigma}$, assumption~4 implies $\sigma(s_{i',j'})=h_{i',j'}$ and $\valu_{\sigma}^\P(h_{i',j'})\rhd\valu_{\sigma}^\P(b_1)$.
It now suffices to prove that this cannot happen, hence we prove that $\sigma(s_{i',j'})=h_{i',j'}$ implies $\valu_{\sigma}^\P(h_{i',j'})\lhd\valu_{\sigma}^\P(b_1)$.

Since $i'<i<\nsb$, we have $i'\leq\nsb+2$, implying $j'=1-\indbit_{i'+1}=1$.
Consequently, $\valustar_{\sigma}^\P(h_{i',j'})=\{h_{i',j'}\}\cup\valustar_{\sigma}^\P(g_{i'+1})$.
Since $\sigma$ has \Pref{USV1}$_{i'+1}$ by assumption~1 and $\sigmabar(g_{i'+1})=1-\indbit_{i'+1}$ by assumption~3, we have $\sigma(s_{i'+1,\sigmabar(g_{i'+1})})=b_1$.
This implies $\lambda_{i'+1}^\P=i'+1$.
Since $i'+1<\nsb$, the cycle center $F_{i'+1,\sigmabar(g_{i'+1})}$ cannot be closed by assumption~2.
As also $\nsigmabar(s_{i'+1})$ and $\nsigmabar(b_{i'+1})$ by \Pref{EV1}$_{i'+1}$, \Cref{corollary: Complete Valuation Of Selection Vertices PG}, implies that either \begin{align*}
\valustar_{\sigma}^\P(g_{i'+1})&=\{g_{i'+1}\}\cup\valustar_{\sigma}^\P(g_1)\quad\text{or}\\
\valustar_{\sigma}^\P(g_{i'+1})&=\{g_{i'+1}\}\cup\valustar_{\sigma}^\P(b_2)\quad\text{or}\\
\valustar_{\sigma}^\P(g_{i'+1})&=\{g_{i'+1},s_{i'+1,\sigmabar(g_{i'+1})}\}\cup\valustar_{\sigma}^\P(b_1).
\end{align*}
As $\relbit{\sigmae}=1$ implies $\sigmae(b_1)=b_2$, the statement directly follows in the last two cases.
In the first case, it follows from $\valustar_{\sigmae}^\P(g_1)\lhd\valustar_{\sigmae}^\P(b_2)$ which can be shown by using \Cref{lemma: Valuations in Phase Three,lemma: Config implied by Aeb} and assumption~3. 

We now show that $(d_{i,j,0},F_{i,j}), (d_{i,j,1},F_{i,j})$ are improving for $\sigmae$.
Let $k\in\{0,1\}$.
It suffices to show $\valu_{\sigmae}^\P(F_{i,j})\rhd\valu_{\sigmae}^\P(e_{i,j,k})$ since $\sigmae(d_{i,j,k})=e_{i,j,k}$ by assumption.
By \Pref{ESC4}$_{i,j}$, we have $\sigmaebar(eb_{i,j})\wedge\neg\sigmaebar(eg_{i,j})$.
Since $\sigmae(d_{i,j,k})=e_{i,j,k}$, this implies $\sigmae(e_{i,j,k})=b_2$.
Hence, by \Cref{lemma: Exact Behavior Of Counterstrategy}, \[\valu_{\sigmae}^\P(F_{i,j})=\{F_{i,j},d_{i,j,k},e_{i,j,k}\}\cup\valu_{\sigmae}^\P(b_2)\rhd\{e_{i,j,k}\}\cup\valu_{\sigmae}^\P(b_2)=\valu_{\sigmae}^\P(e_{i,j,k}).\]
We now explain how we prove $I_{\sigmae}=(I_{\sigma}\setminus\{e\})\cup\{(d_{i,j,0},F_{i,j}), (d_{i,j,1},F_{i,j})\}$.
Applying the switch $e$ increases the valuation of $F_{i,j}$.
By the choice of $j$ and assumption~3, the valuation of $g_i$ increases as well.
We thus begin by showing $\sigmae(b_i)=b_{i+1}$ and $(b_i,g_i)\notin I_{\sigma},I_{\sigmae}$.
However, applying the switch $e$ also increases the valuation of several vertices contained in levels below level $i$.
To be precise, since $\sigma(g_{\ell})=F_{\ell,1}, \tau^{\sigma}(F_{\ell,1})=s_{\ell,1}$ and $\sigma(s_{\ell,1})=h_{\ell,1}$ for all $\ell<i$, the valuation of all of these vertices $g_{\ell}$ and $F_{\ell,1}$ increases.
We thus show that the following statements hold:
\begin{enumerate}
	\item $\sigmae(b_{\ell})=b_{\ell+1}$ and $(b_{\ell},g_{\ell})\notin I_{\sigma},I_{\sigmae}.$
	\item The edges $(d_{\ell,1,0}, F_{\ell,1})$ and $(d_{\ell,1,1},F_{\ell,1})$ are not improving for $\sigmae$.
\end{enumerate}
Since also the valuation of $g_1$ increases, we also prove that $\sigma(e_{i',j',k'})=b_2$ implies $(e_{i',j',k'},g_1)\notin I_{\sigma},I_{\sigmae}$ for any indices $i',j',k'$.
This then proves the statement as $\sigmae(b_1)=b_2$ due to $\relbit{\sigmae}=1$, implying that the valuation of no further vertex can change.

First, since $i<\nsb$, \Pref{EV1}$_i$ implies $\sigmae(b_i)=b_{i+1}$.
By the choice of $i$ and $j$, \Pref{ESC4}$_{i,j}$ implies $\sigmabar(eb_{i,j})\wedge\nsigmabar(eg_{i,j})$.
As $\relbit{\sigmae}=1$ and $\sigmae(s_{i,j})=b_1$ by the choice of $e$, \Cref{corollary: Complete Valuation Of Selection Vertices PG} implies \[\valustar_{\sigmae}^\P(g_i)=\{g_i,s_{i,\sigmabar(g_i)}\}\cup\valustar_{\sigmae}^\P(b_2)<\valustar_{\sigmae}^\P(b_2)=\valustar_{\sigmae}^\P(b_{i+1})\]as $i+1\leq\nsb$.
Since $\valustar_{\sigma}^\P(b_{i+1})=\valustar_{\sigmae}^\P(b_{i+1})=L_{i+1}^\P$ and $\valustar_{\sigma}^\P(g_i)\leq\valustar_{\sigmae}^\P(g_i)$, this implies $(b_i,g_i)\notin I_{\sigma},I_{\sigmae}$.
Now, for any $\ell<i$, $\sigmae(b_{\ell})=b_{\ell+1}$ follows also by \Pref{EV1}$_{\ell}$ and $(b_{\ell},g_{\ell})\notin I_{\sigma},I_{\sigmae}$ follows as \[\valustar_{\sigmae}^\P(g_{\ell})=\bigcup_{i'=\ell}^{i-1}W_{\ell}^\P\cup\{g_{i},s_{i,\sigmabar(g_i)}\}\cup\valustar_{\sigmae}^\P(b_2)\lhd\valustar_{\sigmae}^\P(b_2).\]
Similarly, as $\valustar_{\sigmae}^\P(F_{\ell,1})=\valustar_{\sigmae}(g_{\ell})\setminus\{g_{\ell}\}$ and $\valustar_{\sigmae}^\P(e_{\ell,1,k})=\valustar_{\sigmae}^\P(b_2)$ by \Pref{ESC4}$_{\ell,1}$, the same estimation yields $(d_{\ell,1,0},F_{\ell,1}), (d_{\ell,1,1}, F_{\ell,1})\notin I_{\sigma},I_{\sigmae}$.
Finally, if $\sigmae(e_{i',j',k'})=b_2$, then the same estimation implies $(e_{i',j',k'},g_1)\notin I_{\sigma},I_{\sigmae}$.
Consequently, $I_{\sigmae}=(I_{\sigma}\setminus\{e\}\cup\{(d_{i,j,0},F_{i,j}), (d_{i,j,1}, F_{i,j})\}$.
\end{proof}

\PhaseFiveEscapeEasy*

\begin{proof}
As usual, $\indbit^{\sigma}=\indbit^{\sigmae}\eqqcolon\indbit$ by the choice of $e$.
Since $\sigma$ is a phase-$5$-strategy, it has \Pref{REL1}.
Thus, $\relbit{\sigma}=\min\{i'\colon\sigma(b_{i'})=b_{i'+1}\}$.
Also, by the choice of $e$, $\relbit{\sigmae}=\relbit{\sigma}$.
We first discuss some statements that will be used several times during this proof.

Since $\sigma$ is well-behaved, $\nsb>1$ implies that there is no cycle center $F_{i',j'}$ with $\sigmabar(eg_{i',j'})\wedge\nsigmabar(eb_{i',j'})$.
More precisely, for the sake of a contradiction, assume there was such a cycle center.
By Properties (\ref{property: EG2}) and (\ref{property: EG3}), this implies $\sigmabar(d_1)$ and $\sigmabar(s_1)$.
Thus, by \Pref{USV1}$_1$, the cycle center $F_{1,\indbit^{\sigma}_2}$ is closed.
This however contradicts $\nsb>1$.
Similarly, it is easy to verify that $\nsb=1$ implies that there is no cycle center $F_{i',j'}$ with $\sigmabar(eb_{i',j'})\wedge\nsigmabar(eg_{i',j'})$.
As we assume $\sigmabar(eb_{i,j})\wedge\sigmabar(eg_{i,j})$, the choice of $e$ implies \[\valu_{\sigmae}^\P(F_{i,j})=\{F_{i,j},d_{i,j,*},e_{i,j,*}\}\cup\valu_{\sigmae}^\P(t^{\rightarrow})\] and \[\valu_{\sigmae}^\M(F_{i,j})=(1-\e)\valu_{\sigmae}^\M(t^{\rightarrow})+\e\cdot\valu_{\sigmae}^\M(s_{i,j}).\]

Using these observations, we now prove the statements of the lemma.
\begin{enumerate}
	\item If there are indices $i',j',k', (i,j,k)\neq (i',j',k')$ such that $(e_{i',j',k'},t^{\rightarrow})\in I_{\sigma}$, then $\sigmae$ cannot have \Pref{ESC1}.
		This implies that $\sigmae$ is a phase-$5$-strategy for $\bit$.
		If there is an index $i'$ such that $\sigma$ does not have \Pref{SVG}$_{i'}$ resp. (\ref{property: SVM})$_{i'}$, then $\sigmae$ can also not have the corresponding property.
		Consequently, due to the special condition of phase 5, $\sigmae$ is a phase-$5$-strategy for $\bit$.
	\item We next show that the strategy $\sigmae$ is well-behaved.
		Depending on $\nsb$, we thus need to investigate the following properties:
		\begin{itemize}[align=right, leftmargin=1.75cm]
			\item[(\ref{property: MNS1})] Assume that the premise of this property is correct.
				Then, by \Cref{lemma: Config implied by Aeb}, $\minsige{b}=2$ and consequently $\sigmae(g_1)=F_{1,1}\wedge\sigmae(s_{1,1})=h_{1,1}$.
				But this contradicts that $\nsb=2$ implies $\sigmae(g_1)=\sigma(g_1)=F_{1,0}$ if $G_n=S_n$.
			\item[(\ref{property: MNS2})] Assume that the premise of this property is correct.
				Then $\sigmae(g_1)=F_{1,1}$ and $\sigmae(s_{1,1})=h_{1,1}$.
				But then, by \Pref{USV1}$_1$, $\sigmae(b_2)=g_2$, implying $\minsige{b}=2$, contradicting the assumption.
			\item[(\ref{property: MNS3})] This follows by the same arguments used for \Pref{MNS2}.
			\item[(\ref{property: MNS4})] We only need to investigate this property if $i=\minnegsige{s}=\minnegsig{s}$ and if the premise is true for $\sigmae$.
				But then, the premise was already true for $\sigma$, implying $\sigmabar(eb_{\minnegsig{s}})\wedge\nsigmabar(eg_{\minnegsig{s}})$.
				Since $\relbit{\sigma}=1$ implies $\nsb>1$, we apply an improving switch $(e_{i,j,k},b_2)$.
				But then, also $\sigmaebar(eb_{\minnegsige{s}})\wedge\nsigmaebar(eg_{\minnegsige{s}})$.
			\item[(\ref{property: MNS5})] This follows by the same arguments used for \Pref{MNS2}.
			\item[(\ref{property: MNS6})] This follows by the same arguments used for \Pref{MNS4}.
			\item[(\ref{property: EG1})] By assumption, $\sigmabar(eb_{i,j})\wedge\sigmabar(eg_{i,j})$.
				In order to have $\sigmaebar(eg_{i,j})\wedge\neg\sigmaebar(eb_{i,j})\wedge\relbit{\sigmae}=1$ we thus need to have applied a switch $(e_{i,j,*},g_1)$.
				This however implies $\relbit{\sigmae}\neq 1$, contradicting the premise.
			\item[(\ref{property: EG2})] Follows by the same arguments.
			\item[(\ref{property: EG4})] Follows by the same arguments.
			\item[(\ref{property: EG3})] Assume the premise was true.
				Since $\sigmabar(eg_{i,j})\wedge\sigmabar(eb_{i,j})$, we need to have $\relbit{\sigma}\neq 1$.
				But then, $\sigmabar(b_1)=1$ by \Pref{EV1}$_1$, implying $\sigmabar(s_1)$ by \Pref{EV2}$_1$ and \Pref{USV1}$_1$.
			\item[(\ref{property: EG5})] If the premise is true, then $\sigmabar(s_{i,j})$.
				Hence, $j=\indbit_{i+1}=\sigmabar(b_{i+1})$ by \Pref{USV1}$_{i}$ and \Pref{EV1}$_{i+1}$.
				Since $\sigmabar(s_{i,j})=\sigmaebar(s_{i,j})$ and $\sigmabar(b_{i+1})=\sigmaebar(b_{i+1})$ by the choice of $e$, the statement follows.
			\item[(EB*)] Assume $\sigmaebar(eb_{i,j})\wedge\neg\sigmaebar(eg_{i,j})$.
				Since $\sigmabar(eb_{i,j})\wedge\sigmabar(eg_{i,j})$, we need to have applied a switch $(e_{i,j,*},b_2)$.
				But this implies $\relbit{\sigma}=1$ and thus $\sigma(b_1)=b_2$, contradicting all of the premises.
			\item[(EBG*)] Since $\sigmabar(eb_{i,j})\wedge\sigmabar(eg_{i,j})$ by assumption, it is impossible to have $\sigmaebar(eb_{i,j})\wedge\sigmaebar(eg_{i,j})$ after applying $e$.
				Hence the premise of any of these properties is incorrect.
		\end{itemize}
	
	\item Assume that there are no indices $i',j',k'$ such that $(e_{i',j',k'},t^{\rightarrow})\in I_{\sigma}$ and that $\sigma$ has \Pref{SVG}$_{i'}$/(\ref{property: SVM})$_{i'}$ for all $i'\in[n]$
		We prove that $\sigmae$ is a phase-$1$-strategy for $\bit+1$.
		Since $\sigma$ is a phase-$5$-strategy for $\bit$ and by the choice of $e$, we have $\indbit^{\sigma}=\bit+1=\indbit^{\sigmae}$.
		Also, by our assumptions and the definition of phase-$1$-strategies and phase-$5$-strategies, it suffices to show that $\sigmae$ has \Pref{ESC1}.
		
		Consider the case $\nsb>1$, implying $\relbit{\sigmae}=1$.
		Then, by assumption, there are no further indices $i',j',k'$ such that $(e_{i',j',k'},b_2)\in I_{\sigma}$.
		Thus, for all these indices, it either holds that $\sigma(e_{i',j',k'})=b_2$ or $\sigma(e_{i',j',k'})=g_1\wedge\valu_{\sigma}^*(g_1)\succeq\valu_{\sigma}^*(b_2)$.
		It suffices to prove that the second case cannot occur.
		We do so by proving
		\begin{equation} \label{equation: Valuation Phase Five If NSB>1}
			\nsb>1\implies\valu_{\sigma}^*(g_1)\prec\valu_{\sigma}^*(b_2)
		\end{equation}
		and showing that the arguments also apply to $\sigmae$.
		Consider the different cases listed in \Cref{lemma: Valuation of g1}.		
		If the conditions of either the first or the ``otherwise'' case are fulfilled, then the statement follows.
		It thus suffices to prove that the conditions of the second and third case cannot be fulfilled.
		
		Assume that the conditions of the second case were fulfilled.		
		Since $\relbit{\sigma}=1$ implies $\sigma(b_1)=b_2$ and thus $\minsig{b}\geq 2$, we then have $\sigma(g_1)=F_{1,1}$ and $\sigma(s_{1,1})=h_{1,1}$.
		But then \Pref{USV1}$_1$ yields $\nsb=2$ which is a contradiction if $G_n=S_n$.
		We thus need to have $G_n=M_n$ and $\sigmabar(d_1)$.
		But then, the cycle center $F_{1,\indbit_2}$ is closed, implying $\indbit_1=1$ by definition. 
		This however contradicts $\nsb>1$.
		Thus the conditions of the second case cannot be fulfilled.
		
		Assume that the conditions of the third case were fulfilled.
		Let $\minnegsig{g}>1$, implying $\minsig{b}>2$.
		Then $\sigma(g_1)=F_{1,1}$ and $\sigma(s_{1,1})=h_{1,1}$.
		Then \Pref{USV1}$_1$ implies $\minsig{b}=2$ which is a contradiction.
		Thus let $\minnegsig{g}=1$.
		Then $\sigma(g_1)=F_{1,0}$ and $\sigma(s_{1,0})=h_{1,0}$.
		Consequently, by \Pref{USV1}$_1$, $\indbit_{2}=0$, so in particular $\nsigmabar(b_{\minnegsig{g}+1})$.
		For the sake of a contradiction, assume $\nsigmabar(eb_1)$.
		Since $\nsb>1$ and $\indbit_2=0$ imply that $F_{1,0}$ cannot be closed, it thus needs to hold that $\sigmabar(eg_{1,0})\wedge\nsigmabar(eb_{1,0})$.
		But then, \Pref{EG1} implies $\sigma(s_{1,0})=b_1$, contradicting $\sigma(s_{1,0})=h_{1,0}$.
		Thus the conditions of the third case of \Cref{lemma: Valuation of g1} cannot be fulfilled, implying \Cref{equation: Valuation Phase Five If NSB>1}.
		Note that these arguments can also be applied to $\sigmae$, hence the same statement holds for $\sigmae$.
		
		Now consider the case $\nsb=1$, implying $\relbit{\sigma}\neq 1$ and $\sigma(b_1)=g_1$.
		Then, by assumption, there are no further indices $i',j',k'$ such that $(e_{i',j',k'},g_1)\in I_{\sigma}$.
		Thus, either $\sigma(e_{i',j',k'})=g_1$ or $\sigma(e_{i',j',k'})=b_2$ and $\valu_{\sigma}^*(b_2)\succeq\valu_{\sigma}^*(g_1)$ by the choice of $e$.
		We now show that the second case cannot occur by proving
		\begin{equation} \label{equation: Valuation Phase Five If NSB=1}			
			\nsb=1\implies\valu_{\sigma}^*(b_2)\prec\valu_{\sigma}^*(g_1).
		\end{equation}
		It holds that $\valustar_{\sigma}^*(g_1)=R_1^*$ as $i'<\relbit{\sigma}$ implies $\sigmabar(d_{i'})$ by \Cref{corollary: Simplified MDP Valuation}.
		Consider the case  $\sigma(b_2)=g_2$ first.
		Then, $\relbit{\sigma}>2$, implying $\valustar_{\sigma}^*(b_2)=R_2^*$.		
		Hence $\valustar_{\sigma}^*(b_2)=R_2^*\prec R_1^*=\valustar_{\sigma}^*(g_1)$.
		
		Now let $\sigma(b_2)=b_3$.
		Then, $\relbit{\sigma}=2$ and $\valustar_{\sigma}^*(b_2)=L_2^*$, so \Cref{lemma: VV Lemma} implies $\valustar_{\sigmae}^*(b_2)\prec R_1^*=\valustar_{\sigmae}^*(g_1)$.
		Therefore, $\valustar_{\sigma}^*(b_2)\prec\valustar_{\sigma}^*(g_1)$ in any case, contradicting $\valu_{\sigma}^*(b_2)\succeq\valu_{\sigma}^*(g_1)$.
		Note again that the same arguments apply to $\sigmae$.
	\item[$*$] Before we prove the remaining aspects, we prove that $(d_{i,j,1-k},F_{i,j})\in I_{\sigmae}$ in any case.
		As we assume $\sigmabar(eb_{i,j})\wedge\sigmabar(eg_{i,j})$, it holds that $\sigmae(d_{i,j,1-k})=e_{i,j,1-k}$.
		It thus suffices to show $\valu_{\sigmae}^*(F_{i,j})\succ\valu_{\sigmae}^*(e_{i,j,1-k})$.
		By assumption $\sigma(e_{i,j,1-k})=\sigmae(e_{i,j,1-k})=t^{\rightarrow}$.
		If $G_n=S_n$, the statement thus follows from \[\valu_{\sigmae}^\P(F_{i,j})=\{F_{i,j},d_{i,j,*},e_{i,j,*}\}\cup\valu_{\sigmae}^\P(t^{\rightarrow})\rhd\{e_{i,j,1-k}\}\cup\valu_{\sigmae}^\P(t^{\rightarrow})=\valu_{\sigmae}^\P(e_{i,j,1-k}).\]
		If $G_n=M_n$, we then have $\valu_{\sigmae}^\M(F_{i,j})=(1-\e)\valu_{\sigmae}^\M(t^{\rightarrow})+\e\valu_{\sigmae}^\M(s_{i,j})$.
		To prove the statement, it thus suffices to prove $\valu_{\sigmae}^\M(s_{i,j})>\valu_{\sigmae}^\M(t^{\rightarrow})$.

		We only consider the case $\nsb>1$, the case $\nsb=1$ follows analogously.
		In this case, $t^{\rightarrow}=b_2$.
		If $j=\indbit_{i+1}$, then $\valustar_{\sigmae}^\M(s_{i,j})=\rew{s_{i,j},h_{i,j}}+\valustar_{\sigmae}^\M(b_{i+1})$ by \Pref{EV1}$_i$.
		The statement thus follows since $\rew{h_{i,j}}>\sum_{\ell\in[i]}W_{\ell}^\M$.
		If $j\neq\indbit_{i+1}$, then $\relbit{\sigmae}=1$ implies $\sigmae(b_1)=b_2$  and thus $\valustar_{\sigmae}^\M(s_{i,j})=\rew{s_{i,j}}+\valustar_{\sigmae}^\M(b_1)=\rew{s_{i,j}}+\valustar_{\sigmae}^\M(b_2)$.
		
	\item We prove that $(g_i,F_{i,j})$ is improving for $\sigmae$ if and only if the corresponding conditions are fulfilled.
		Assume that the corresponding conditions are fulfilled.
		We distinguish the following cases.
		\begin{enumerate}
			\item The cycle center $F_{i,1-j}$ cannot be closed as either $\sigmabar(eb_{i,1-j})$ or $\sigmabar(eg_{i,1-j})$.
			\item Let $F_{i,1-j}$ be $t^{\rightarrow}$-open.
				Then, if $G_n=S_n$, the statement follows since $j=0$ and \begin{align*}
				\valu_{\sigmae}^\P(F_{i,0})&=\{F_{i,0},d_{i,j,*},e_{i,j,*}\}\cup\valu_{\sigmae}^\P(t^{\rightarrow})\\
					&\rhd\{F_{i,1},d_{i,1,*},e_{i,1,*}\}\cup\valu_{\sigmae}^\P(t^{\rightarrow})=\valu_{\sigmae}^\P(F_{i,1}).
				\end{align*}
				If $G_n=M_n$, then the statement follows by \Cref{lemma: Both CC Open For MDP} since $j=\indbit_{i+1}$ and $F_{i,j}$ is also $t^{\rightarrow}$-open.
			\item Let $F_{i,1-j}$ be $t^{\rightarrow}$-halfopen.
				If $G_n=S_n$, then the statement follows analogously to the last case.
				If $G_n=M_n$, then the statement follows by an easy but tedious calculation.
			\item Let $F_{i,1-j}$ be mixed.
				Then, by either Equation~(\ref{equation: Valuation Phase Five If NSB>1}) or Equation~(\ref{equation: Valuation Phase Five If NSB=1}), it holds that $\valustar_{\sigmae}^*(F_{i,1-j})\preceq\valustar_{\sigmae}^*(t^{\rightarrow})$.
				We can thus use the same arguments used in one of the last two cases to prove the statement.
		\end{enumerate}
		As we proved that it is not possible that any cycle center escapes only to $t^{\leftarrow}$ at the beginning of this proof, these are all cases that need to be covered.
		Hence, if all the stated conditions are fulfilled, then the edge $(g_i,F_{i,j})$ is an improving switch for $\sigmae$.
		We now prove that $(g_i,F_{i,j})$ is not improving for $\sigmae$ if any of these conditions is not fulfilled, proving the claimed equivalence.
		We consider the different conditions one after another.
		\begin{enumerate}
			\item Let $\indbit_i=1$.
				Then, since $\sigmae$ is a phase-$5$-strategy for $\bit$, it holds that $\sigmae(b_i)=g_i$ and $\sigmae(g_i)=F_{i,\indbit_{i+1}}$.
				Furthermore, this cycle center is then closed.
				Since $\sigmabar(eb_{i,j})\wedge\sigmabar(eg_{i,j})$ by assumption, we then need to have $j=1-\indbit_{i+1}$ and consequently $\sigmaebar(g_i)=1-j$.
				By Properties (\ref{property: USV1})$_i$ and (\ref{property: EV1})$_{i+1}$, this thus yields $\valustar_{\sigmae}^*(F_{i,1-j})=\ubracket{s_{i,1-j},h_{i,1-j}}\oplus\valustar_{\sigmae}^*(b_{i+1}).$
				Since $\valustar_{\sigmae}^*(F_{i,j})=\valustar_{\sigmae}^*(t^{\rightarrow})$, the statement thus follows from $\ubracket{s_{i,j},h_{i,j}}\succ\bigoplus_{\ell\leq i}W_{\ell}$.
			\item Let $\sigmabar(g_i)=0$ resp. $\sigmabar(g_i)=\indbit_{i+1}$ depending whether we consider $G_n=S_n$ or $G_n=M_n$.
				Due to the first case, we may assume $\indbit_i=0$.
				We furthermore may assume $j=1$ resp. $j=1-\indbit_{i+1}$ if $G_n=S_n$ resp. $G_n=M_n$ as we otherwise already have $(g_i,F_{i,j})\notin I_{\sigmae}$ by definition.
				We prove $\valu_{\sigmae}^*(F_{i,1-j})\preceq\valu_{\sigmae}^*(F_{i,j})$ by distinguishing the different possible states of $F_{i,1-j}$.
				\begin{enumerate}
					\item Let $F_{i,1-j}$ be closed.
						Since $\indbit_i=0$, this implies $1-j=1-\indbit_{i+1}$.
						Thus, by \Pref{USV1}$_i$, \begin{align*}
							\valustar_{\sigmae}^*(F_{i,1-j})&=\ubracket{s_{i,1-j}}\oplus\valustar_{\sigmae}^*(b_1)\succ\valustar_{\sigmae}^*(b_1)=\valustar_{\sigmae}^*(F_{i,j}).
						\end{align*}
					\item Let $F_{i,1-j}$ be $t^{\rightarrow}$-open or $t^{\rightarrow}$-halfopen.
						Since $1-j=0$ resp. $1-j=\indbit_{i+1}$, the same arguments used when proving that $(g_i,F_{i,j})\in I_{\sigmae}$ can be applied if the corresponding conditions are fulfilled to obtain $\valu_{\sigmae}^*(F_{i,1-j})\succeq\valu_{\sigmae}^*(F_{i,j})$ in either case.
					\item Let $F_{i,1-j}$ be mixed.
						Then, $\sigmaebar(eb_{i,1-j})\wedge\sigmaebar(eg_{i,1-j})$ by the choice of $e$.
						However, in any context and for any $\nsb$, this contradicts the assumptions of the lemma.
				\end{enumerate}
				By the observations made at the beginning of this proof, these are all cases that can occur.
			\item Let $j=1$ resp. $j=1-\indbit_{i+1}$.
				Due to the first two cases, we may assume $\indbit_{i+1}=0$ and $\sigmaebar(g_i)=1$ resp. $\sigmaebar(g_i)=1-\indbit_{i+1}$.
				But this implies $(g_i,F_{i,j})\notin I_{\sigmae}$ by the definition of an improving switch.
			\item We only discuss the last condition for $\nsb>1$ as the statement follows for $\nsb=1$ analogously.
				Hence let $\nsigmabar(eb_{i,1})$ resp. $\nsigmaebar(eb_{i,1-\indbit_{i+1}})$.
				Due to the last cases, we may assume $\indbit_i=0$, $\sigmabar(g_i)=1$ resp. $\sigmabar(g_i)=1-\indbit_{i+1}$ and $j=0$ resp. $j=\indbit_{i+1}$.
				By the observations made at the beginning of the proof, we cannot have $\sigmaebar(eg_{i,1})$ resp. $\sigmaebar(eg_{i,1-\indbit_{i+1}})$.
				This implies that we need to have $\sigmaebar(d_{i,1})$ resp. $\sigmaebar(d_{i,1-\indbit_{i+1}})$, implying the statement in either case.
				More precisely, we then either have \[\valustar_{\sigmae}^*(F_{i,1-j})=\ubracket{s_{i,1-j}}\oplus\valustar_{\sigmae}^*(b_1)=\ubracket{s_{i,1-j}}\oplus\valustar_{\sigmae}^*(b_2)\succ\valustar_{\sigmae}^*(F_{i,j})\] or \[\valustar_{\sigmae}^*(F_{i,1-j})=\ubracket{s_{i,1-j},h_{i,1-j}}\oplus\valustar_{\sigmae}^*(b_{i+1})\succ\valustar_{\sigmae}^*(b_2)=\valustar_{\sigmae}^*(F_{i,j}).\]
		\end{enumerate}
		Thus, if any of the given conditions is not fulfilled, then $(g_i,F_{i,j})\notin I_{\sigmae}$.
		Consequently, $(g_i,F_{i,j})\in I_{\sigmae}$ if and only if the stated conditions are fulfilled.

		It remains to show that no other improving switches are created in any case.
		By assumption, $\sigmabar(eb_{i,j})\wedge\sigmabar(eg_{i,j})$, implying $\sigmae(d_{i,j,*})=e_{i,j,*}$.
		Since the application of~$e$ increases the valuation of $F_{i,j}$, we begin by proving $(d_{i,j,k},F_{i,j})\in I_{\sigma},I_{\sigmae}$.
		This however follows easily by Equations~(\ref{equation: Valuation Phase Five If NSB>1}) and (\ref{equation: Valuation Phase Five If NSB=1}) since \begin{align*}
			\valustar_{\sigma}^*(d_{i,j,k})=\valustar_{\sigma}^*(t^{\leftarrow})&\prec\valustar_{\sigma}^*(t^{\rightarrow})=\valustar_{\sigma}^*(F_{i,j}),\\
			\valu_{\sigmae}^\P(d_{i,j,k})=\{e_{i,j,k}\}\cup\valu_{\sigmae}^\P(t^{\rightarrow})&\lhd\{F_{i,j}, e_{i,j,k^*}, d_{i,j,k^*}\}\cup\valu_{\sigmae}^\P(t^{\rightarrow})=\valu_{\sigmae}^\P(F_{i,j})\\,
			\valu_{\sigmae}^\M(d_{i,j,k})=\valu_{\sigmae}^\M(t^{\rightarrow})&<(1-\e)\valu_{\sigmae}^\M(t^{\rightarrow})+\e\valu_{\sigmae}^\M(s_{i,j})=\valu_{\sigmae}^\M(F_{i,j})
		\end{align*}
		for some $k^*\in\{0,1\}$ since $\valu_{\sigmae}^\M(s_{i,j})>\valu_{\sigmae}^\M(b_1)$.
		Moire precisely, if $j=1-\indbit_{i+1}$, then this follows directly since we have $\valu_{\sigmae}(s_{i,j})=\rew{s_{i,j}}+\valu_{\sigmae}^\M(b_1)$ in that case.
		If $j=\indbit_{i+1}$ then this follows as $\valu_{\sigmae}^\M(s_{i,j})=\rew{s_{i,j},h_{i,j}}+\valu_{\sigmae}^\M(b_{i+1})$ in that case.
		
		We now consider the possible change of the valuation of $g_i$.
		If $\sigmae(g_i)\neq j$, then the edge $(g_i,F_{i,j})$ is the only edge (besides the edges $(d_{i,j,*},F_{i,j})$ that we already considered) that might become improving for $\sigmae$.
		We however already completely described the conditions under which this edge becomes improving.
		Hence consider the case $\sigmae(g_i)=j$.
		We first observe that we cannot have $i=1$ since $\sigmabar(eg_{i,j})$ would then contradict the fact that $G_n$ is a sink game resp. weakly unichain.
		We prove that we then have $\sigmae(b_i)=b_{i+1}$ and $\valu_{\sigmae}^*(b_{i+1})\succeq\valustar_{\sigmae}^*(g_i)$ as well as $\sigmae(s_{i-1,1})=b_1$ and $\valu_{\sigmae}^*(b_1)\succeq\valu_{\sigmae}^*(h_{i-1,1})$.
		By assumption and the choice of $e$, the cycle center $F_{i,j}$ is not closed for $\sigmae$.
		If $j=\indbit_{i+1}$, this implies $\sigmae(b_i)=b_{i+1}$ and $\sigmae(s_{i-1,1})=b_1$ by \Pref{EV1}$_i$ resp. \Pref{USV1}$_{i-1}$.
		If $j=1-\indbit_{i+1}$, then we need to have $\indbit_i=0$ since \Pref{EV1}$_i$ and \Pref{EV2}$_i$ would imply $\sigmae(g_i)=\indbit_{i+1}=1-j$ otherwise.
		Thus $\sigmae(b_i)=b_{i+1}$ and $\sigmae(s_{i-1,1})=b_1$ in any case.
		We now prove $\valu_{\sigmae}^*(b_{i+1})\succeq\valu_{\sigmae}^*(g_i)$.
		Since $\sigmae(b_{i})=b_{i+1}$ implies $i\geq\relbit{\sigmae}$ by \Pref{REL1}, we have $\valustar_{\sigmae}^*(b_{i+1})=L_{i+1}^*$.
		If $\nsb=1$, then $\sigmaebar(eg_{i,j})\wedge\nsigmaebar(eb_{i,j})\wedge\relbit{\sigmae}\neq 1$.
		Then, by \Cref{lemma: Exact Behavior Of Random Vertex} resp. \Cref{lemma: Exact Behavior Of Counterstrategy}, $\valustar_{\sigmae}^*(F_{i,j})=\valustar_{\sigmae}^*(g_1)=R_1^*.$
		This in particular implies $\valustar_{\sigmae}^*(g_i)=\ubracket{g_i}\oplus R_1^*$.
		The statement thus follows directly since $\ubracket{g_i}\prec\bigoplus_{\ell=1}^{i-1}W_{\ell}^*$.
		If $\nsb>1$, then $\sigmaebar(eb_{i,j})\wedge\nsigmaebar(eg_{i,j})\wedge\relbit{\sigmae}=1$.
		Since Properties (\ref{property: EV1})$_{i+1}$ and (\ref{property: USV1})$_i$ imply that either $\sigmae(s_{i,j})=b_1$ or $\sigmae(b_{i+1})=j$, \Cref{lemma: Exact Behavior Of Random Vertex} resp. \Cref{lemma: Exact Behavior Of Counterstrategy} thus imply $\valustar_{\sigmae}^*(g_i)=\ubracket{g_i}\oplus\valustar_{\sigmae}^*(b_2).$
		Then, the statement again follows since $\ubracket{g_i}\prec\bigoplus_{\ell=1}^{i-1}W_{\ell}^*$.
		Therefore $\valustar_{\sigmae}^*(b_{i+1})\prec\valustar_{\sigmae}^*(g_i)$ in any case.
		Since $\relbit{\sigmae}=1\Leftrightarrow\sigmae(b_1)=b_2\Leftrightarrow\nsb>1$, the same arguments imply \[\valustar_{\sigmae}^*(h_{i-1,1})=\ubracket{h_{i-1,1},g_i}\oplus\valustar_{\sigmae}^*(b_1)\prec\valustar_{\sigmae}^*(b_1).\]\qedhere
	\end{enumerate}
\end{proof}

\PhaseFiveCycleVertices*

\begin{proof}
	As in the last proofs, we have $\indbit^{\sigma}=\indbit^{\sigmae}\eqqcolon\indbit$ by the choice of $e$.
	Let $j\coloneqq 1-\indbit_{i+1}$.
	We begin by showing that $\sigmae$ is a phase-$5$-strategy for $\bit$.
	If $\sigma(d_{i,j,1-k})\neq F_{i,j}$, then the same cycle centers are closed with respect to $\sigma$ and $\sigmae$.
	In this case, $\sigma$ being a phase-$5$-strategy immediately implies that $\sigmae$ is a phase-$5$-strategy.
	Thus assume $\sigma(d_{i,j,1-k})=F_{i,j}$.
	Then, the cycle center $F_{i,j}=F_{i,1-\indbit_{i+1}}$ is closed with respect to $\sigma$ but not with respect to $\sigmae$.
	It thus sufficient to investigate \Pref{EV3}$_i$ and \Pref{CC1}$_i$.
	Since, by assumption, $\sigmae(b_i)=\sigma(b_i)=b_{i+1}$, \Pref{EV3}$_i$ remains valid.
	As $\sigma$ has \Pref{REL1}, $\sigmae$ also has \Pref{REL1}.
	This implies $\relbit{\sigmae}=\min\{i'\colon\sigmae(b_{i'})=b_{i'+1}\}$.
	Thus, $i\geq\relbit{\sigmae}$, so $\sigmae$ has \Pref{CC1}$_i$.
	Therefore, $\sigmae$ is a phase-$5$-strategy for $\bit$.
	
	We now show that $\sigmae$ is well-behaved.
	If $\sigma(d_{i,j,1-k})\neq F_{i,j}$, this follows immediately since $\sigma$ is well-behaved.
	Hence assume $\sigma(d_{i,j,1-k})=F_{i,j}$ and note that we have $i\geq\relbit{\sigmae}$ as argued earlier.
	Since we close a cycle center $F_{i,\sigmabar(g_i)}$ with $i\geq\relbit{\sigmae}$, we investigate the following properties.
	\begin{itemize}[leftmargin=1.75cm]
		\item[(\ref{property: MNS4})] Since $\sigma$ is well-behaved, this property only needs to be reevaluated if $i=\minnegsige{s}$.
			Since the case $i=1$ cannot occur by assumption, assume $i>1$.
			This implies $1<\minnegsige{s}\leq\minnegsige{g}<\minsige{b}$.
			Thus, in particular $\sigmae(b_1)=b_{2}, \sigmae(g_1)=F_{1,1}$ and $\sigmae(s_{1,1})=h_{1,1}$.
			By \Pref{USV1}$_1$ and \Pref{EV1}$_2$, $\sigmae(b_2)=g_2$, implying $\minsige{b}=2$.
			But this contradicts the premise since $1<\minnegsige{s}<\minsige{b}$ implies $\minsige{b}\geq 3$.
		\item[(\ref{property: MNS6})] Since $\sigma$ is well-behaved, this only needs to be reevaluated if $i=\minnegsige{s}$.
		Since $i\neq1$ by assumption, assume $i>1$.
		Then, $1<\minnegsige{s}\leq\minnegsige{g}<\minsige{b}$, implying the same contradiction as in the last case.
		\item[(DN*)] Since there is no cycle center $F_{n,1}$ by construction, we cannot have $i=n$.
	\end{itemize}
	Consequently, $\sigmae$ is well-behaved.
	
	It remains to show that $I_{\sigmae}=I_{\sigma}\setminus\{e\}$.
	We distinguish three different cases.
	\begin{itemize}
		\item The cycle center $F_{i,j}$ is closed with respect to $\sigmae$.
			Then, since $j=1-\indbit_{i+1}$, we have $\sigmae(s_{i,j})=b_1$ by \Pref{USV1}$_i$, implying $\valustar_{\sigmae}^*(F_{i,j})=\ubracket{s_{i,j}}\oplus\valustar_{\sigmae}^*(b_1)$.
			The only vertices that have an edge towards $F_{i,j}$ are $d_{i,j,0},d_{i,j,1}$ and $g_i$. 
			Since $F_{i,j}$ is closed for $\sigmae$ and $\sigmabar(g_i)=j$ by assumption, the valuation of these vertices might change when applying $e$.
			However, $\sigmae(d_{i,j,k})=F_{i,j}$ implies that $(d_{i,j,0},F_{i,j}),(d_{i,j,1},F_{i,j})\notin I_{\sigmae}$.
			Since no player 0 vertex has an edge to $d_{i,j,*}$, consider the vertex $g_i$.		
			The only vertices having an edge towards $g_i$ are $b_i$ and $h_{i-1,1}$.
			Since $\sigmae(b_i)=b_{i+1}$ by assumption and $\sigmae(s_{i-1,1})=b_1$ by \Pref{USV1}$_i$, it suffices to show $(b_i,g_i),(s_{i-1,1},h_{i-1,1})\notin I_{\sigmae}$.
			We begin by showing $(b_i,g_i)\notin I_{\sigmae}$.
			It suffices to show $\valustar_{\sigmae}^*(b_{i+1})\rhd\valustar_{\sigmae}^*(g_i)$.
			As mentioned before, we have $\valustar_{\sigmae}^*(g_i)=\ubracket{g_i,s_{i,j}}\oplus\valustar_{\sigmae}^*(b_1)$.
			If $\relbit{\sigmae}=1$, the statement follows since \begin{align*}
				\valustar_{\sigmae}^*(g_i)&=\ubracket{g_i,s_{i,j}}\oplus\valustar_{\sigmae}^*(b_1)=\ubracket{g_i,s_{i,j}}\oplus L_1^*\prec L_{i+1}^*=\valustar_{\sigmae}^*(b_{i+1}).
			\end{align*}
			Now consider the case $\relbit{\sigmae}\neq 1$, implying $\valustar_{\sigmae}^*(b_1)=R_1^*$ since $i'<\relbit{\sigmae}$ implies $\sigmabar(d_{i'})$ by \Cref{corollary: Simplified MDP Valuation} and \Pref{REL1}.
			Since $\sigmae(b_i)=b_{i+1}$ implies $i+1>\relbit{\sigmae}$, we have $\valustar_{\sigmae}^*(b_{i+1})=L_{i+1}^*$.
			Consequently, \begin{align*}
				\valustar_{\sigmae}^*(g_i)&=\ubracket{g_{i},s_{i,j}}\oplus R_1^*\prec\bigoplus_{i'\geq i+1}\{W_{i'}^*\colon\sigmae(b_{i'})=g_{i'}\}=L_{i+1}^*=\valustar_{\sigmae}^*(b_{i+1}).
			\end{align*}		 
			It remains to prove $(s_{i-1,1},h_{i-1,1})\notin I_{\sigmae}$ by showing $\valu_{\sigma}^*(b_1)\succ\valu_{\sigmae}^*(h_{i-1,1})$.
			This however follows by $\valustar_{\sigmae}^*(h_{i-1,1})=\ubracket{h_{i-1,1},g_i,s_{i,j}}\oplus\valustar_{\sigmae}^*(b_1)\prec\valustar_{\sigmae}^*(b_1)$.
		\item The cycle center $F_{i,j}$ is not closed and $\sigmaebar(eb_{i,j})$.
			Then, since $\sigmae(d_{i,j,k})=F_{i,j}$, we have $\sigmaebar(eb_{i,j})\wedge\nsigmaebar(eg_{i,j})$.
			Consider the case that $G_n=S_n$.
			Then, by \Cref{lemma: Exact Behavior Of Counterstrategy}, either $\valustar_{\sigmae}^\P(F_{i,j})=\valustar_{\sigmae}^\P(s_{i,j})$ or $\valustar_{\sigmae}^\P(F_{i,j})=\valustar_{\sigmae}^\P(b_2)$.
			In the first case we can use the same arguments as before to prove $(b_i,g_i),(s_{i-1,1},h_{i-1,1})\notin I_{\sigmae}$.
			Hence consider the second case.
			Then, by \Cref{lemma: Exact Behavior Of Counterstrategy}, we need to have $\relbit{\sigmae}=1$.
			But then, it follows that $\valustar_{\sigmae}^\P(g_i)=\{g_i\}\cup\valustar_{\sigmae}^\P(b_1)=\{g_i\}\cup L_1^\P\lhd L_{i+1}^\P=\valustar_{\sigmae}^\P(b_{i+1})$ and $\valustar_{\sigmae}^\P(h_{i-1,1})=\{h_{i-1,1},g_1\}\cup\valustar_{\sigmae}^\P(b_1)\lhd\valustar_{\sigmae}^\P(b_1),$ so $(b_i,g_i),(s_{i-1,1},h_{i-1,1})\notin I_{\sigmae}$.
			Now consider the case $G_n=M_n$.
			Since $F_{i,j}$ is $b_2$-halfopen, $\valustar_{\sigmae}^\M(F_{i,j})=\valustar_{\sigmae}^\M(b_2)$.
			It suffices to prove $\relbit{\sigmae}=1$ as we can then apply the same arguments used in the case $G_n=S_n$.
			But this follows from \Pref{EB5} as $\sigmaebar(eb_{i,j})\wedge\nsigmaebar(eg_{i,j})\wedge\relbit{\sigmae}\neq 1$ implies $\sigmae(b_{\relbit{\sigmae}})=g_{\relbit{\sigmae}}$, contradicting \Pref{REL1}.
		\item The cycle center $F_{i,j}$ is not closed and $\sigmaebar(eg_{i,j})$.
			Similar to the last case, we then have $\sigmaebar(eg_{i,j})\wedge\nsigmaebar(eg_{i,j})$.
			Consider the case $G_n=S_n$.
			Then,  by \Cref{lemma: Exact Behavior Of Counterstrategy}, either $\valustar_{\sigmae}^\P(F_{i,j})=\valustar_{\sigmae}^\P(s_{i,j})$ or $\valustar_{\sigmae}^\P(F_{i,j})=\valustar_{\sigmae}^\P(g_1)$.
			Since the second case implies $\relbit{\sigmae}\neq 1$, similar arguments as the ones used previously can be used to show $(b_i,g_i),(s_{i-1,1},h_{i-1,1})\notin I_{\sigmae}$ in both cases.
			If $G_n=M_n$, then $\valustar_{\sigmae}^\M(F_{i,j})=\valustar_{\sigmae}^\M(g_1)$ and it again suffices to prove $\relbit{\sigmae}\neq 1$.
			This follows from \Pref{EG1} and \Pref{EG2} since these properties would imply that the cycle center $F_{1,\indbit_2}$ was closed.
			Then, \Pref{EV1}$_1$ would imply $\sigmae(b_1)=g_1$, contradicting $\relbit{\sigmae}=1$.\qedhere
	\end{itemize}
\end{proof}

\PhaseFiveSelectorVertices*

\begin{proof}
	We first show that $\sigmae$ is a phase-$1$-strategy for $\bit+1$ resp. a phase-$5$-strategy for $\bit$.
	We observe that we have $\indbit^{\sigma}=\indbit^{\sigmae}\eqqcolon\indbit$ as the status of no cycle center or entry vertex is changed.
	Since we change the target of a selector vertex with $\indbit_i=0$, it suffices to check Properties (\ref{property: REL1}), (\ref{property: CC2}), (\ref{property: EV2})$_i$ and (\ref{property: SVG})$_i$/(\ref{property: SVM})$_i$.
	
	It is immediate that $\sigmae$ has \Pref{REL1} as $\indbit_i=0$.
	To prove that it has \Pref{CC2} assume $i=\nsb$.
	But this  implies $\indbit=1$, contradicting again the assumption.
	By definition, $\sigmae$ has \Pref{ESC1} if and only if there are no indies $i',j',k'$ with $\sigma(e_{i',j',k'})\neq t^{\rightarrow}$.
	Thus, if there are no such indices and if $\sigmae$ has Property (\ref{property: SVG})$_{i'}$/(\ref{property: SVM})$_{i'}$ for all $i'\in[n]$, then $\sigmae$ is a phase-$1$-strategy for $\bit+1$.
	Otherwise, it is  a phase-$5$-strategy for $\bit$.
	In particular, $\relbit{\sigma}=\relbit{\sigmae}=\min\{i'\colon\sigmae(b_{i'})=b_{i'+1}\}$.
	
	We next show that $\sigmae$ is well-behaved.
	Since we change the target of a selector vertex and $i\neq n$, we need to investigate the following assumptions:
	\begin{itemize}[leftmargin=1.75cm]
		\item[(\ref{property: S2})] Since $\indbit_i=0$ implies $i\geq\relbit{\sigmae}$, it cannot hold that $i<\relbit{\sigmae}$.
		\item[(\ref{property: D2})] This follows by the same argument.
		\item[(\ref{property: B3})] Since $\sigmae$ has \Pref{USV1}$_i$ and \Pref{EV1}$_{i+1}$, the premise of this assumption is always incorrect.
		\item[(\ref{property: BR1})] Since $\indbit_i=0$ implies $i\geq\relbit{\sigmae}$, it cannot hold that $i<\relbit{\sigmae}$.
		\item[(\ref{property: MNS1})] If the premise is correct for both $\sigma$ and $\sigmae$, then $\sigmae$ has this property as $\sigma$ has it.
			The implication is also fulfilled if the premise is incorrect for $\sigmae$.
			Hence assume that the premise is correct for $\sigmae$ but incorrect for $\sigma$.
			Since $\sigmae$ is well-behaved, \Cref{lemma: Config implied by Aeb} implies $\minsige{b}=2$.
			Thus, $\sigmae(b_2)=g_2$, hence $\sigmae(s_{1,1})=h_{1,1}$ and $\sigmae(s_{1,0})=b_1$.
			As we assume that the premise is incorrect for $\sigma$, the choice of $e$ implies $\minnegsig{g}=1$ and thus $e=(g_1,F_{1,1})$.
			We thus need to have $\sigma(g_1)=F_{1,0}$ and $\valu_{\sigma}^\P(F_{1,1})\rhd\valu_{\sigma}^\P(F_{1,0})$.
			We show that this cannot be true.

			Since we have $\relbit{\sigmae}=\relbit{\sigma}=1$ and $\sigmae(b_2)=g_2$, the cycle center $F_{1,1}$ cannot be closed.
			Consequently, by assumption, $\sigmabar(eb_{1,1})\wedge\nsigmabar(eg_{1,1})$.
			Thus, \[\valu_{\sigma}^\P(F_{1,1})=\{F_{1,1},d_{1,1,k^*},e_{1,1,k^*}\}\cup\valu_{\sigma}^\P(b_2)\] for some $k^*\in\{0,1\}$ by \Cref{lemma: Exact Behavior Of Counterstrategy}.
			If also $\sigmabar(eb_{1,0})\wedge\nsigmaebar(eg_{1,0})$, then this also yields $\valu_{\sigma}^\P(F_{1,0})=\{F_{1,0},d_{1,0,k^*},e_{1,0,k^*}\}$ for some $k^*\in\{0,1\}$.
			The statement then follows since $\Omega(F_{1,0})>\Omega(F_{1,1})$ and since the priority of $F_{1,0}$ is even.
			If this is not the case, then $F_{1,0}$ is closed by assumption.
			But this implies \[\valustar_{\sigma}^\P(F_{1,0})=\{s_{1,0}\}\cup\valustar_{\sigma}^\P(b_2)>\valustar_{\sigma}^\P(b_2)=\valustar_{\sigma}^\P(F_{1,1}).\]

		\item[(\ref{property: MNS2})] Assume $\relbit{\sigmae}=1$, let $i'<\minnegsige{g}<\minnegsige{s},\minsige{b}$ and let $G_n=S_n$ imply $\nsigmaebar(b_{\minnegsige{g}+1})$.
			Then $\sigmae(b_2)=b_3$ since $\sigmae(b_2)=g_2$ implies $\minsige{b}=2$, contradicting the premise.
			Consequently, $\indbit_{2}=0$.
			However, since $1<\minnegsige{g}<\minnegsige{s}$, it holds that $\sigmae(g_1)=F_{1,1}$ and $\sigmae(s_{1,1})=h_{1,1}$.
			But this implies $\indbit_2=1$ by \Pref{USV1}$_1$ and \Pref{EV1}$_2$ which is a contradiction.
		\item[(\ref{property: MNS3})] If the premise is true, then $\indbit_2=0$ since we need to have $\sigmae(b_2)=b_3$.
			But, since $1<\minnegsige{s}\leq\minnegsige{g}$ implies $\sigmae(g_1)=F_{1,1}$ and $\sigmae(s_{1,1})=h_{1,1}$, we also have $\indbit_2=1$ which is a contradiction.
		\item[(\ref{property: MNS4})] Let $\relbit{\sigmae}=1$ and $\minnegsige{s}\leq\minnegsige{g}<\minsige{b}$.
			If $\minnegsige{s}>1$, then the same arguments used for proving that $\sigmae$ has \Pref{MNS3} can be used to prove that $\sigmae$ has \Pref{MNS4} as follows.
			Thus assume $1=\minnegsige{s}$, implying that we have $\sigmae(s_{1,\sigmaebar(g_1)})=b_1$.
			In particular, by \Pref{USV1}$_1$ and \Pref{EV2}$_2$, it holds that $\sigmaebar(g_1)\neq\sigmaebar(b_2)=\indbit_2$.
			If $\minnegsig{s}=\minnegsige{s}$ and $\minnegsig{g}=\minnegsige{g}$, then the statement follows by applying \Pref{MNS4} to~$\sigma$.
			Thus assume $\minnegsige{s}\neq\minnegsig{s}$.
			Then $\sigmae(s_{1,\sigmabar(g_1)})=h_{1,\sigmabar(g_1)}$.
			But this implies $e=(g_i,F_{i,j})=(g_1,F_{1,1-\indbit^{\sigmae}_{i+1}})$ and thus, by assumption, $\sigmaebar(eb_{1,j})\wedge\nsigmaebar(eg_{i,j})$.
			Hence assume $\minnegsige{s}=\minnegsig{s}$ and $\minnegsige{g}\neq\minnegsig{g}$.
			If $\minnegsig{g}<\minsig{b}$, then the statement follows since we can again apply \Pref{MNS4} to $\sigma$.
			Thus assume $\minnegsig{g}\geq\minsig{b}$.
			But then, $\minnegsig{s}<\minsig{b}\leq\minnegsig{g}$, hence $\minnegsige{s}=\minnegsig{s}$ and applying \Pref{MNS6} to $\sigma$ imply $\sigmaebar(eb_{\minnegsige{s}})\wedge\sigmaebar(eg_{\minnegsige{g}})$.
		\item[(\ref{property: MNS5})] If the premise is true, then $1<\minnegsige{s}<\minsige{b}\leq\minnegsige{g}$.
			In particular, $\sigmae(g_1)=F_{1,1}$ and $\sigmae(s_{1,1})=h_{1,1}$.
			By \Pref{EV1}$_2$, this implies $\sigmae(b_2)=g_2$ and thus $\minsige{b}=2$.
			This however is a contradiction since the premise implies $\minsige{b}\geq 3$.
		\item[(\ref{property: MNS6})] If $\minnegsige{s}>1$, then the same arguments used for \Pref{MNS5} can be used to prove that the premise cannot be correct.
			Hence assume $\minnegsige{s}=1$, implying $\sigmae(g_1)=F_{1,1}$ and $\sigmae(s_{1,1})=b_1$.
			This in particular implies $1=\sigmaebar(g_1)\neq\indbit^{\sigmae}_2=0$ and thus $\sigma(b_2)=\sigmae(b_2)=b_3$ by \Pref{USV1}$_1$ and \Pref{EV1}$_2$.
			If $\minnegsig{s}=\minnegsige{s}$ and $\minnegsig{g}=\minnegsige{g}$, then the statement follows by applying \Pref{MNS6} to $\sigma$.
			Assume $\minnegsig{s}\neq\minnegsige{s}$.
			Then $\sigma(s_{1,\sigmabar(g_1)})=h_{1,\sigmabar(g_1)}$ and thus $\sigmabar(g_1)=\indbit^{\sigmae}_{2}$.
			But then, $e=(g_i,F_{i,j})=(g_1,F_{1,1-\indbit^{\sigmae}_2})$, so $\sigmaebar(eb_{i,j})\wedge\nsigmaebar(eg_{i,j})$ by assumption.
			Thus assume $\minnegsig{s}=\minnegsige{s}$ and $\minnegsig{g}\neq\minnegsige{g}$.
			Since the statement follows by applying \Pref{MNS6} to $\sigma$ if $\minnegsig{s}<\minsig{b}\leq\minnegsig{g}$, assume $\minnegsig{g}<\minsig{b}$.
			But then, $1=\minnegsig{s}\leq\minnegsig{g}<\minsig{b}$.
			Since applying an improving switch in level 1 implies $\minnegsig{s}\neq\minnegsige{s}$, we have $\sigmabar(g_1)=\sigmaebar(g_1)$.
			But then the statement follows by applying \Pref{MNS4} to $\sigma$.
		\item[(EG*)] Since $\relbit{\sigmae}=1$ implies $\nsb>1$, we have $\sigmaebar(eb_{i,j})\wedge\nsigmaebar(eg_{i,j})$.
			Hence the premise of any of the Properties (\ref{property: EG1}) to (\ref{property: EG4}) is incorrect.
			Note that we do not need to validate \Pref{EG5}. 
		\item[(EBG*)] By assumption, we cannot have $\sigmaebar(eb_{i,j})\wedge\sigmaebar(eg_{i,j})$, hence the premise of any of these assumptions is incorrect.
	\end{itemize}

	Hence $\sigmae$ is a well-behaved strategy.
	
	It remains to show that $I_{\sigmae}=I_{\sigma}\setminus\{e\}$.
	Since we apply the improving switch $e=(g_{i},F_{i,j})$, the valuation of $g_i$ increases.
	If $i\neq 1$, then there are only two vertices that have an edge to $g_i$, namely $b_i$ and $h_{i-1,1}$.
	However, if $i=1$, then also the valuation of escape vertices and hence cycle centers might be influenced.
	We prove that $\sigmae(b_i)=b_{i+1}\wedge (b_i,g_i)\notin I_{\sigmae}$ for all $i\in[n]$ and $\sigmae(s_{i-1,1})=b_1\wedge(s_{i-1,1},h_{i-1,1})\notin I_{\sigmae}$ if $i>1$.
	We then discuss the case $i=1$ at the end of this proof.

	Thus let $i\in[n]$.	
	Since $\indbit_i=0$ and by \Pref{EV1}$_i$, it holds that $\sigmae(b_i)=b_{i+1}$.
	It thus suffices to prove $\valu_{\sigmae}^*(b_{i+1})\succ\valu_{\sigmae}^*(g_i)$.
	We distinguish the following cases.
	\begin{enumerate}
		\item Let $\relbit{\sigmae}=1$.
			Then $\valustar_{\sigmae}^*(b_{i+1})=L_{i+1}^*$.
			By assumption, $\sigmaebar(eb_{i,j})\wedge\neg\sigmaebar(eg_{i,j})$.
			Thus, depending on whether $G_n=S_n$ or $G_n=M_n$, \Cref{lemma: Exact Behavior Of Counterstrategy} and \Pref{USV1}$_{i}$ respectively \Cref{lemma: Exact Behavior Of Random Vertex} imply $\valustar_{\sigmae}^*(F_{i,j})=\valustar_{\sigmae}^*(b_2)$.
			Consequently, \begin{align*} 
				\valustar_{\sigmae}^*(g_i)&=\ubracket{g_i}\oplus\valustar_{\sigmae}^*(b_2)=\ubracket{g_i}\oplus L_2^*=\ubracket{g_i}\oplus L_{2,i-1}^*\oplus L_{i+1}^*\prec L_{i+1}^*.
			\end{align*}
		\item Let $\relbit{\sigmae}\neq 1$. 
			Since $\indbit_i=0$, it cannot hold that $\valustar_{\sigmae}^*(b_{i+1})=R_{i+1}^*$.
			This implies that $\valustar_{\sigmae}^*(b_{i+1})=L_{i+1}^*$ and $i\geq\relbit{\sigmae}$.
			By assumption, $\sigmaebar(eg_{i,j})\wedge\nsigmaebar(eb_{i,j})$.
			Thus, by \Cref{lemma: Exact Behavior Of Random Vertex} resp. \ref{lemma: Exact Behavior Of Counterstrategy}, $\valustar_{\sigmae}^*(g_i)=\ubracket{g_i}\oplus\valustar_{\sigmae}^*(g_1)$.
			Note that $\valustar_{\sigmae}^*(b_1)=R_1^*$ in any case by \Cref{corollary: Simplified MDP Valuation}.
			Thus, by \Pref{USV1}$_i$ and since $i\geq\relbit{\sigmae}$, \begin{align*}
				\valustar_{\sigmae}^*(g_i)&=\ubracket{g_i,s_{i,j}}\oplus\valustar_{\sigmae}^*(b_1)=\ubracket{g_i,s_{i,j}}\oplus R_1^*\\
					&=\ubracket{g_i,s_{i,j}}\oplus\bigoplus_{i'=1}^{\relbit{\sigmae}-1}W_{i'}^*\oplus L_{\relbit{\sigmae}+1,i-1}^*\oplus L_{i+1}^*\prec L_{i+1}^*=\valustar_{\sigmae}^*(b_{i+1}).
			\end{align*}
	\end{enumerate}
	Thus $\valustar_{\sigmae}^*(g_i)\prec\valustar_{\sigmae}^*(b_{i+1})$ in any case, implying $(b_i,g_i)\notin I_{\sigmae}$.
	
	We prove that $i\neq 1$ implies  $\sigmae(s_{i-1,1})\neq h_{i-1,1}$ and $(s_{i-1,1},h_{i-1,1})\notin I_{\sigmae}$.
	The first statement follows since $\indbit_{i}=0$ and \Pref{USV1}$_{i-1}$ imply $\sigmae(s_{i-1,1})=b_1$.
	It thus remains to prove $\valu_{\sigmae}^*(b_1)\succ\valu_{\sigmae}^*(h_{i-1,1})$.
	We again distinguish the following cases.
	\begin{enumerate}
		\item Assume $\relbit{\sigmae}=1$ first.
			Then $\valustar_{\sigmae}^*(b_1)=L_1^*$.
			By assumption, $\sigmae(eb_{i,j})\wedge\neg\sigmae(eg_{i,j})$.
			By \Pref{USV1}$_i$ and since $\relbit{\sigmae}=1$ implies $\sigmae(b_1)=b_2$, we then have \[\valustar_{\sigmae}^*(h_{i-1,1})=\{h_{i-1,1},g_i\}\oplus\valustar_{\sigmae}^*(b_2)\prec\valustar_{\sigmae}^*(b_2)=\valustar_{\sigmae}^*(b_1).\]
		\item Let $\relbit{\sigmae}\neq 1$.
			Then $\sigmae(b_1)=g_1$, implying $\valustar_{\sigmae}^*(b_1)=\valustar_{\sigmae}^*(g_1)=R_1^*$ by \Cref{corollary: Simplified MDP Valuation}.
			By assumption, $\sigmaebar(eg_{i,j})\wedge\nsigmaebar(eb_{i,j})$.
			Thus,  by \Cref{lemma: Exact Behavior Of Random Vertex} resp. \Cref{lemma: Exact Behavior Of Counterstrategy}, $\valustar_{\sigmae}^*(h_{i-1,1})=\ubracket{h_{i-1,1},g_i}\oplus\valustar_{\sigmae}^*(g_1)\prec\valustar_{\sigmae}^*(b_1).$
	\end{enumerate}
	
	It remains to discuss the case $i=1$.
	Since $\indbit_i=0$ by assumption, we then have $\relbit{\sigmae}=1$.
	In particular, for any cycle center $F_{i',j'}$, either $\sigmaebar(d_{i',j'})$ or $\sigmaebar(eb_{i',j'})\wedge\nsigmaebar(eg_{i',j'})$ by assumption.
	Thus, the valuation of no cycle center is increased as this could only happen if $\sigmaebar(eg_{i',j'})$.
	Moreover, there is $d_{i',j',k'}$ with $\sigmae(d_{i',j',k'})=e_{i',j',k'}$ and $\sigmae(e_{i',j',k'})=g_1$.
	It thus suffices to prove that $\sigmaebar(d_{i',j'})\wedge\sigmae(e_{i',j',k'})=g_1$ implies $(d_{i',j',k'},e_{i',j',k'})\notin I_{\sigma},I_{\sigmae}$ and that $\sigmae(e_{i',j',k'})=b_2$ implies $(e_{i',j',k'},g_1)\notin I_{\sigma},I_{\sigmae}$.
	
	Consider the second statement first.
	It suffices to prove $\valustar_{\sigmae}^*(b_2)\succ\valustar_{\sigmae}^*(g_1)$.
	If $\sigmaebar(eb_1)\wedge\nsigmaebar(eg_1)$, then this follows since $\valustar_{\sigmae}^*(g_1)=\ubracket{g_1}\oplus\valustar_{\sigmae}^*(b_2)$ in that case.
	If $\sigmaebar(d_1)$, then we need to have $\sigmaebar(g_i)=1-\indbit_2$ due to $\relbit{\sigmae}=1$ and $\nsb>1$.
	But then, the statement follows since $\valustar_{\sigmae}^*(g_1)=\ubracket{g_1,s_{1,\indbit_{2}}}\oplus\valustar_{\sigmae}^*(b_2)$.
	Since the same arguments hold for $\sigma$, the statement follows.
	Thus consider some cycle center $F_{i',j'}$ closed with respect to $\sigmae$.
	It suffices to show $\valustar_{\sigmae}(F_{i',j'})\succ\valustar_{\sigmae}^*(g_1)$.
	If $j'=1-\indbit_{i'+1}$, then the statement follows since $\valustar_{\sigmae}^*(F_{i',j'})=\ubracket{s_{i',j'}}\oplus\valustar_{\sigmae}^*(b_2)$ in this case and since $\valustar_{\sigmae}^*(b_2)\succ\valustar_{\sigmae}^*(g_1)$ as proven before.
	If $j'=\indbit_{i'+1}$, then $\valustar_{\sigmae}^*(F_{i',j'})=\ubracket{s_{i',j'},h_{i',j'}}\oplus\valustar_{\sigmae}^*(b_{i'+1})$ and the statement follows since $\valustar_{\sigmae}^*(F_{i',j'})\succ\valustar_{\sigmae}^*(b_2)$.
\end{proof}

\subsection{Omitted proofs of \Cref{section: Final Formal Proofs}}

The following statements are claims that are used within proofs of the statements in \Cref{section: Final Formal Proofs}.
Each claim thus refers to the notation used in the corresponding proof, and this notation is not restated here.

\SelectorVerticesInPhaseOne*

\begin{proof}
Consider the first phase-$1$-strategy $\sigma$ such that after applying an improving switch $e=(d_{i,j',k},F_{i,j'})$ to $\sigma$, the edge $(g_i,F_{i,j'})$ becomes improving for $\sigmae$.
Then, $\applied{\canstrat}{\sigmae}\subseteq\mathbb{D}^1$.
Furthermore, $\sigmae$ is a phase-$1$-strategy for~$\bit$ by \Cref{lemma: Phase 1 Low OR} and $I_{\sigmae}=\mathfrak{D}^{\sigmae}\cup\{(g_{i},F_{i,j'})\}$.
Moreover, $F_{i,j'}$ is closed for $\sigmae$ and $\occrec^{\canstrat}(g_i,F_{i,j'})\leq\occrec^{\canstrat}(d_{i,j',k},F_{i,j'})$ by \Cref{table: Occurrence Records}.
Since $(d_{i,j',k},F_{i,j'})$ minimized the occurrence record for $\sigma$, the switch $(g_i,F_{i,j'})$ minimizes the occurrence record for $\sigmae$.
By the tie-breaking rule, this switch is thus applied next.
Since $e\in I_{\canstrat}^{<\maxocc}$, \Cref{lemma: Occurrence Records Cycle Vertices} implies that it cannot happen the cycle center $F_{1,1-\bit_2}$ was closed by applying~$e$, so $i\neq 1$.
It is easy to verify that the other conditions of row 4 of \Cref{table: Phase 1 Switches} hold as well, since $(g_i,F_{i,j'})$ would not have become an improving switch otherwise.
Thus, by row 4 of \Cref{table: Phase 1 Switches}, the strategy $\sigma'$ obtained by applying $(g_i,F_{i,j'})$ to $\sigmae$ is a well-behaved phase-$1$-strategy for $\bit$ with $\sigma'\in\reach{\sigma_0}$ and $I_{\sigma'}=\mathfrak{D}^{\sigma'}$.
This proves that the first switch of the type $(g_*,F_{*,*})$ is applied immediately when it becomes an improving switch.
The same arguments can however also be applied for any edge $(g_*,F_{*,*})$ that becomes improving.
\end{proof}

\NoCycleCenterOpenPhaseOne*

\begin{proof}
Assume there were indices $i\in[n], j'\in\{0,1\}$ with $\sigma(d_{i,j',1})\neq F_{i,j'}$ and let $e\coloneqq(d_{i,j',1},F_{i,j'})$.
Then $e\in I_{\sigma}$, so $\occrec^{\sigma}(e)=\maxocc$.
Since $\sigma(d_{i,j',1})\neq F_{i,j'}$ implies that  $F_{i,j'}$ is not closed, $\bit_i=0\vee\bit_{i+1}\neq j'$ as $\sigma$ is a phase-$1$-strategy for $\bit$.
As $e$ was not applied during $\canstrat\to\sigma$, this yields \[\occrec^{\canstrat}(e')=\occrec^{\sigma}(e')=\min\left (\floor{\frac{\bit+1-k}{2}}, \ell^{\bit}(i,j',k)+t_{\bit}\right )\] for a feasible $t_{\bit}$ for $\bit$.
In particular, $\occrec^{\sigma}(e)\leq\floor{(\bit+1-1)/2}=\floor{\bit/2}$.
But this is a contradiction, since $\occrec^{\sigma}(e)=\maxocc>\floor{\bit/2}$ since $\bit$ is odd.
\end{proof}

\OcccurrenceRecordOfCycleVerticesPhaseOne*

\begin{proof}
Consider some fixed indices $i\in[n],j,k\in\{0,1\}$ such that $(d_{i,j,k},F_{i,j})\in\applied{\canstrat}{\sigma^{(3)}}$.
As argued previously, this edge is contained in a cycle center $F_{i,j}$ which is open for~$\canstrat$.
If the occurrence record of one of the cycle edges of $F_{i,j}$ is $\maxocc-1$, then the application of $(d_{i,j,k},F_{i,j})$ is described by \Cref{lemma: Phase 1 Low OR} and we do not need to consider it here.
Also, due to the tie-breaking rule, we do not apply an improving switch contained in halfopen cycle centers (with the exception of $F_{\nsb,\bit_{\nsb+1}}$) as we only consider switches contained in~$I_{\canstrat}^{\maxocc}$.
We may thus assume that $F_{i,j}$ is open with respect to $\canstrat$ and that both cycle edges have an occurrence record of $\maxocc$.

We now distinguish between several possible indices.
Consider the case $i\neq\nsb$ or $i=\nsb\wedge j\neq\bit_{i+1}$ first.
By the tie-breaking rule, we then need to have $k=0$ as the edge $e\coloneqq(d_{i,j,0},F_{i,j})$ is then applied as improving switch.
Let $\sigma$ denote the strategy in which $e$ is applied.
Since $\bit$ is even, $\occrec^{\sigmae}(e)=\maxocc+1=\floor{(\bit+2)/2}$.
It thus suffices to show that there is a parameter~$t_{\bit+1}$ feasible for $\bit+1$ such that \begin{equation} \tag{$\star$}
\floor{\frac{\bit}{2}}+1\leq\ell^{\bit+1}(i,j,k)+t_{\bit+1}.
\end{equation}
By the choice of $i$ and $j$, \Cref{lemma: Progress for unimportant CC} implies $\ell^{\bit}(i,j,k)+1=\ell^{\bit+1}(i,j,k)$.
Therefore, $\occrec^{\sigma}(e)+1\leq\ell^{\bit}(i,j,k)+t_{\bit}+1\leq\ell^{\bit+1}(i,j,k)+t_{\bit}$ for some $t_{\bit}$ feasible for~$\bit$.
Since $\bit$ is even, \Pref{OR4}$_{i,j,0}$ implies $\occrec^{\canstrat}(e)\neq \ell^{\bit}(i,j,k)-1$.
In addition, by \Pref{OR2}$_{i,j,0}$, it holds that $t_{\bit}\neq 1$ as this would imply $\canstrat(d_{i,j,0})=F_{i,j}$, contradicting our assumption.
Consequently, $t_{\bit}=0$, implying that  $t_{\bit+1}=0$ is a feasible parameter that yields $\star$.

Consider the case $i=\nsb$ and $j=\bit_{\nsb+1}$ next.
Then, both switches $(d_{i,j,*},F_{i,j})$ are applied.
Using \Cref{lemma: Numerics Of Ell}, it is easy to verify that $\occrec^{\canstrat}(d_{i,j,k},F_{i,j})=\floor{(\bit+1-k)/2}$ for both $k\in\{0,1\}$.
Also, by the tie-breaking rule, $F_{i,j}$ is closed once there are no more open cycle centers.
In particular, both cycle edges of $F_{i,j}$ are then applied and their application is described by row 1 resp.~5 of \Cref{table: Phase 1 Switches}.
Let $\sigma$ denote the strategy obtained after closing~$F_{i,j}$.
Then, by definition and the choice of $i$ and $j$, it holds that  $\bit+1=\lastflip{\bit+1}{\nsb}{\{(\nsb+1,j)\}}$.
Since \[\ceil{\frac{(\bit+1)+1-k}{2}}=\floor{\frac{\bit+1+2-k}{2}}=\floor{\frac{\bit+1-k}{2}}+1,\] this then implies \[\occrec^{\sigma}(d_{i,j,k},F_{i,j})=\ceil{\frac{\lastflip{\bit+1}{\nsb}{\{(\nsb+1,j)\}}+1-k}{2}}\] as required.
\end{proof}

\USVPhaseTwoIfNSBTwo*

\begin{proof}
As a reminder, the current strategy is denoted by $\sigmae$ and the set of improving switches for $\sigmae$ is given by $I_{\sigmae}=\mathfrak{D}^{\sigmae}\cup\{(b_1,b_2),(s_{1,1},h_{1,1})\}\cup\{(e_{*,*,*},b_2)\}$.
By \Cref{table: Occurrence Records}, $\occrec^{\sigmae}(e_{*,*,*},b_2)=\floor{\bit/2},\occrec^{\sigmae}(b_1,b_2)=\flips{\bit}{1}{}-1$ and $\occrec^{\sigmae}(s_{1,1},h_{1,1})=\flips{\bit}{2}{}$.
Since~$\bit$ is odd, \Cref{lemma: Numerics Of OR} implies $\flips{\bit}{1}{}-1=\floor{\frac{\bit}{2}}$.
Consequently, $\occrec^{\sigmae}(b_1,b_2)=\occrec^{\sigmae}(e_{*,*,*},b_2).$
If $\bit=1$, then $\flips{\bit}{2}{}=\floor{\frac{\bit+2}{4}}=0=\floor{\frac{\bit}{2}}$ and $\occrec^{\sigmae}(s_{1,1},h_{1,1})=\occrec^{\sigmae}(b_1,b_2)$.
In this situation $(s_{1,1},h_{1,1})$ is applied next as this is exactly the exception described in which the tie-breaking rule behaves differently, see \Cref{definition: Tie-Breaking exponential}.
If $\bit>1$, then $\nsb=2$ implies $\bit\geq 5$.
But this implies $\occrec^{\sigmae}(s_{1,1},h_{1,1})<\occrec^{\sigmae}(b_1,b_2)$, so $(s_{1,1},h_{1,1})$ is applied next. 
Consequently, $e'\coloneqq(s_{1,1},h_{1,1})$ is applied next in any case.

We  now prove that the requirements of row 2 \Cref{table: Phase 2 Switches} are fulfilled. 
Since $\relbit{\sigmae}=\nsb=2$, we show the following statements:
\begin{enumerate}
	\item \boldall{$\sigmaebar(d_1):$} No switch of the type $(d_{*,*,*},e_{*,*,*})$ was applied during $\canstrat\to\sigma^{(2)}$ by \Cref{lemma: Extended Reaching phase 2}.
		Also, no such switch or switch of the type $(g_{*},F_{*,*,})$ was applied during $\sigma^{(2)}\to\sigmae$.
		Thus, by \Cref{lemma: Extended Reaching phase 2}, $\sigmaebar(d_{1})$.
	\item \boldall{$\sigmae$ has \Pref{USV3}$_{1}$:} Since $\nsb=2$, we have $\indbit^{\sigmae}_2=1$.
		Since we did not apply any improving switch of the type $(s_{*,*},*)$ during $\canstrat\to\sigmae$, the statement then follows by applying \Pref{USV1}$_1$ to~$\canstrat$.
	\item \boldall{$\sigmae$ has \Pref{EV2}$_2$ and \Pref{CC2}:} 
		We already argued that $\sigmae$ has these properties when applying the statement described by row 1 of \Cref{table: Phase 2 Switches}.
\end{enumerate}		
Thus, all requirements of row 2 of \Cref{table: Phase 2 Switches} are met.

To simplify notation ,we denote the strategy obtained by applying $e'=(s_{1,1},h_{1,1})$ to $\sigmae$ again by $\sigma$.
Then, $\sigma$ is a phase-3-strategy for $\bit$ with $\sigma\in\reach{\sigma_0}$ and \[I_{\sigma}=I_{\sigmae}\setminus\{e'\}=\mathfrak{D}^{\sigma}\cup\{(b_1,b_2)\}\cup\{(e_{*,*,*},b_2)\}.\]
Since $\bit$ is odd, \[\occrec^{\sigma}(e')=\flips{\bit}{2}{}+1=\floor{\frac{\bit+2}{4}}+1=\floor{\frac{(\bit+1)+2}{4}}=\flips{\bit+1}{2}{}\] and \Cref{table: Occurrence Records} describes the occurrence record of $(s_{1,1},h_{1,1})$ with respect to $\bit+1$.
Since we did not apply any improving switch $(g_*,F_{*,*})$ or $(d_{*,*,*},e_{*,*,*})$, the conditions on cycle centers in levels below $\nsb$ hold for $\sigma^{(3)}$ as they held for $\sigma^{(2)}$.
Therefore,~$\sigma$ is a strategy as described by the respective rows of \Cref{table: Properties at start of phase,table: Switches at start of phase}.
\end{proof}

\SecondSelectionSwitchInPhaseTwo*

\begin{proof}
By the definition of $\nsb$,  there is a number $k\in\mathbb{N}$ such that $\bit=k\cdot2^{\nsb-1}-1$.
By \Cref{table: Occurrence Records}, \Cref{lemma: Numerics Of OR} and using $\nsb>2$, we obtain the following:
\allowdisplaybreaks
\begin{align*}
	\occrec^{\sigma}(s_{\nsb-1,1},h_{\nsb-1,1})&=\flips{\bit}{\nsb}{}=\floor{\frac{k\cdot2^{\nsb-1}-1+2^{\nsb-1}}{2^{\nsb}}}=\floor{\frac{k}{2}}	\\[1em]
	\occrec^{\sigma}(s_{\nsb-2,0},h_{\nsb-2,0})&=\flips{\bit}{\nsb-1}{}-1=k\cdot 2^{0}-1=k-1	\\[1em]
	\occrec^{\sigma}(s_{\nsb-3,0},h_{\nsb-3,0})&=\occrec^{\sigma}(b_{\nsb-2},b_{\nsb-1})=\flips{\bit}{\nsb-2}{}-1=k\cdot 2^{2-1}-1=2k-1\\[1em]
	\occrec^{\sigma}(e_{*,*,*},b_2)&=\floor{\frac{\bit}{2}}=\floor{\frac{k\cdot2^{\nsb-1}-1}{2}}=\floor{k\cdot 2^{\nsb-2}-\frac{1}{2}}=2^{\nsb-2}k-1.
\end{align*}						
If $k>2$, then $(s_{\nsb-1,1},h_{\nsb-1,1})$ is the unique improving switch minimizing the occurrence records.
If $k\leq 2$, then the occurrence records of $(s_{\nsb-1,1},h_{\nsb-,1})$ and $(s_{\nsb-2,0},h_{\nsb-2,0})$ are identical and lower than the occurrence record of every other improving switch.
Since the tie-breaking rule applies improving switches at selection vertices contained in higher levels first, $(s_{\nsb-1,1},h_{\nsb-1,1})$ is also applied first then.
Consequently, $e\coloneqq (s_{\nsb-1,1},h_{\nsb-1,1})$ is applied next in any case.

We prove that $\sigma$ fulfills the conditions of row 2 of \Cref{table: Phase 2 Switches}.
By our previous arguments, it suffices to prove that $\sigma$ has \Pref{USV3}$_{\nsb-1}$.
As $\indbit_{\nsb}=1$, this  follows since $(s_{\nsb-1,1},h_{\nsb-1,1})\in I_{\sigma}$ and since $\sigma$ has \Pref{USV2}$_{\nsb-1,0}$ by the definition of a phase-$2$-strategy.
By our previous arguments and row 2 of \Cref{table: Phase 2 Switches}, $\sigmae$ then has Properties~(\ref{property: USV2})$_{\nsb-1,1}$, (\ref{property: CC2}), (\ref{property: EV1})$_\nsb$ and (\ref{property: USV3})$_{i,1-\indbit_{i+1}}$ for all $i<\nsb-1$.
Furthermore, \[I_{\sigmae}=\mathfrak{D}^{\sigmae}\cup\{(s_{\nsb-2,0},h_{\nsb-2,0}),(b_{\nsb-2},b_{\nsb-1}),(s_{\nsb-3,0},h_{\nsb-3,0})\}\] if $\nsb-1>2$ and $\nsb>2$ implies \[I_{\sigmae}=\mathfrak{D}^{\sigmae}\cup\{(e_{*,*,*},b_2)\}\cup\{(b_1,b_2),(s_{1,0},h_{1,0})\}.\]
Also note that $\nsb>2$ implies $\occrec^{\sigmae}(s_{\nsb-1,1},h_{\nsb-1,1})=\flips{\bit}{\nsb}{}+1=\flips{\bit+1}{\nsb}{}$ by \Cref{lemma: Numerics Of OR}, so \Cref{table: Occurrence Records} specifies its occurrence record with respect to $\bit+1$.
\end{proof}

\EasyCVSwitchesPhaseThree*

\begin{proof}
As usual, we define $t^{\leftarrow}\coloneqq\{g_1,b_2\}\setminus\{t^{\rightarrow}\}$.
We consider the case $G_n=S_n$ and $j=\indbit_{i+1}$ first.
Then, since $\sigma$ is a phase-3-strategy for $\bit$, it cannot happen that $\sigmabar(d_{i,j})$ as this would imply $\indbit_i=1$.
In particular, $\sigmabar(d_{i,j})$ thus implies $\indbit_{i+1}\neq j$.
But then, by \Cref{lemma: PG Cycle Vertices Phase Three}, $\sigmae$ is a well-behaved phase-3-strategy for $\bit$ with $I_{\sigmae}=I_{\sigma}\setminus\{e\}$.

Hence consider the case $G_n=M_n$ and $\indbit_i=1$.
Then, by assumption, $\indbit_{i+1}=1-j$, so the statement follows by \Cref{lemma: MDP Phase 3 Level Is Set}.
It thus suffices to consider the case $\indbit^{\sigma}_i=0$, implying $i>\nsb$ by assumption.
We remind here that $\relbit{\sigma}=\nsb$ by \Pref{REL2} and distinguish two cases.

\boldall{$1:$ Let $\indbit_{i+1}=j$.}
We prove that the application of $e$ is then described by \Cref{lemma: Improving Switch in Other CC Phase Three} or \Cref{lemma: No CC Closed Phase Three}.
We begin by proving that $F_{i,j}$ is $t^{\leftarrow}$-halfopen and then discuss the possible states of $F_{i,1-j}$.

Since $\indbit_i=0\wedge\indbit_{i+1}=j$, the cycle center $F_{i,j}$ was not closed for $\sigma^{(3)}$.
In particular, as the choice of $e$ implies $\sigma(d_{i,j,k})=\sigma^{(3)}(d_{i,j,k})=F_{i,j}$,  \Cref{corollary: No Open CC In Phase Three} and $\applied{\sigma^{(3)}}{\sigma}\cap\mathbb{D}^1=\emptyset$ imply $\sigma(d_{i,j,1-k})=e_{i,j,1-k}$.
As improving switches were applied according to Zadeh's pivot rule and our tie-breaking rule, this implies $(e_{i,j,1-k},t^{\rightarrow})\notin\applied{\sigma^{(3)}}{\sigma}$.
Consequently, $\sigma(d_{i,j,1-k})=e_{i,j,1-k}\wedge\sigma(e_{i,j,1-k})=t^{\leftarrow}$, so $F_{i,j}$ is $t^{\leftarrow}$-halfopen with respect to $\sigma$.

We now discuss the possible states of $F_{i,1-j}$.
First, $F_{i,1-j}$ cannot be $t^{\leftarrow}$-open for $\sigma$ as this would imply that it is also $t^{\leftarrow}$-open for $\sigma^{(3)}$, contradicting \Cref{corollary: No Open CC In Phase Three}.

Also, $F_{i,1-j}$ cannot be closed as it would then be the unique closed cycle center in level $i$.
Then, the tie-breaking rule would have applied some switch $(e_{i,1-j,*},t^{\rightarrow})$.
But this would have made the corresponding edge $(d_{i,1-j,*},e_{i,1-j,*})$ improving by \Cref{lemma: Escape Vertices Phase Three}.
Furthermore, this switch would then already have been applied, contradicting the assumption that $F_{i,1-j}$ was closed.

Let, for the sake of contradiction, $F_{i,1-j}$ be mixed.
Then, $\sigma(d_{i,1-j,*})=e_{i,1-j,*}$ as well as $\sigma(e_{i,1-j,k'})=t^{\rightarrow}$ and $\sigma(e_{i,1-j,1-k'})=t^{\leftarrow}$ for some $k'\in\{0,1\}$.
This implies that $F_{i,1-j}$ was $t^{\leftarrow}$-halfopen with respect to $\sigma^{(3)}$ and that $(e_{i,1-j,k'},t^{\rightarrow})\in\applied{\sigma^{(3)}}{\sigma}$.
Hence, this switch was ranked higher by the tie-breaking rule.
But this is a contradiction as the tie-breaking rule ranks switches contained in $F_{i,\indbit^{\sigma}_{i+1}}=F_{i,j}$ higher if both cycle centers are $t^{\leftarrow}$-halfopen.

It is also immediate that $F_{i,1-j}$ cannot be $t^{\rightarrow}$-halfopen as the tie-breaking rule would then choose the edge $(e_{i,1-j,k'},t^{\rightarrow})$ with $\sigma(e_{i,1-j,k'})=t^{\leftarrow}$ as next improving switch.

Now assume that $F_{i,1-j}$ is $t^{\rightarrow}$-open.
We show that this implies that $F_{i,1-j}$ was closed at the end of phase 1. 
As the cycle center $F_{i,1-j}$ is $t^{\rightarrow}$-open and $\sigma^{(3)}(e_{i,1-j,*})=t^{\leftarrow}$, this implies $(e_{i,1-j,0},t^{\rightarrow}), (e_{i,1-j,1},t^{\rightarrow})\in\applied{\sigma^{(3)}}{\sigma}$.
As all improving switches $(e_{*,*,*},t^{\rightarrow})$ have the same occurrence records, this implies that the tie-breaking rule ranked $(e_{i,1-j,0},t^{\rightarrow})$ and $(e_{i,1-j,1},t^{\rightarrow})$ higher than $(e_{i,j,k},t^{\rightarrow})$.
However, since $j=\indbit_{i+1}$, this can only happen if $F_{i,1-j}$ was closed with respect to $\sigma^{(3)}$. 
If $F_{i,1-j}$ was not closed for $\canstrat$, then \Cref{corollary: Selection Vertices In Phase One} and $\applied{\sigma ^3}{\sigma}\subseteq\mathbb{D}^0\cup\mathbb{E}^1$ imply $\sigma^{(3)}(g_i)=\sigma(g_i)=F_{i,1-j}$.
If it was closed for $\canstrat$, then $\canstrat(g_i)=F_{i,1-j}$ by \Cref{definition: Canonical Strategy MDP}.
Moreover, by the choice of $j$ and $i$ and \Cref{corollary: Selection Vertices In Phase One}, it is not possible that the cycle center $F_{i,j}$ was closed during phase 1.
Consequently, also $\sigma(g_i)=F_{i,1-j}$.
Thus, the statement follows by \Cref{lemma: Improving Switch in Other CC Phase Three}.

Finally, assume that $F_{i,1-j}$ is $t^{\leftarrow}$-halfopen.
Then, since $\applied{\sigma^{(3)}}{\sigma}\subseteq\mathbb{E}^1\cup\mathbb{D}^0$, this implies that $F_{i,1-j}$ was $t^{\leftarrow}$-halfopen for $\sigma^{(3)}$.
In particular, this implies that no cycle center of level~$i$ was closed during phase~1.
But this implies $\sigma(g_i)=\sigma^{(3)}(g_i)=\canstrat(g_i)=F_{i,\indbit_{i+1}}=F_{i,j}$ by \Cref{corollary: Selection Vertices In Phase One}.
Since $i\geq\nsb+1=\relbit{\sigma}+1$ by assumption, \Cref{lemma: No CC Closed Phase Three} implies the statement.
This concludes the case $j=\indbit_{i+1}$.

\boldall{$2:$ Let $1-\indbit_{i+1}=j$.}
We investigate $F_{i,j}$ first.
As $j=1-\indbit_{i+1}$, it is possible that $F_{i,j}$ was closed with respect to $\sigma^{(3)}$.
Depending on whether or not improving switches corresponding to $F_{i,j}$ were applied during $\sigma^{(3)}\to\sigma$, the cycle center is either (a) closed, (b) $t^{\rightarrow}$-halfopen or (c) $t^{\leftarrow}$-halfopen for $\sigma$.
Consider the cycle center $F_{i,1-j}$.
It cannot be closed as $1-j=\indbit_{i+1}$ and $\indbit_i=0$.
If $F_{i,1-j}$ was $t^{\leftarrow}$-open with respect to $\sigma$, then the assumption $\applied{\sigma^{(3)}}{\sigma}\subseteq\mathbb{E}^1\cup\mathbb{D}^0$ implies that it was $t^{\leftarrow}$-open with respect to $\sigma^{(3)}$, contradicting \Cref{corollary: No Open CC In Phase Three}.
If $F_{i,1-j}$ was $t^{\rightarrow}$-open, then $(e_{i,1-j,k'},t^{\rightarrow})\in\applied{\sigma^{(3)}}{\sigma}$ for both $k'\in\{0,1\}$.
This implies that $\sigma^{(3)}(d_{i,1-j,k'})=F_{i,1-j}$, hence $F_{i,1-j}$ was closed with respect to $\sigma^{(3)}$.
But this is not possible as $1-j=\indbit_{i+1}$ and $i>\nsb$ then imply $\indbit_i=1$, contradicting that $F_{i,1-j}$ is $t^{\rightarrow}$-open.
By the same argument, $F_{i,1-j}$ cannot be $t^{\rightarrow}$-halfopen for $\sigma$.

Now assume that $F_{i,1-j}$ is mixed.
Then, $(e_{i,1-j,k'},b_2), (d_{i,1-j,k'},e_{i,1-j,k'})\in\applied{\sigma_3}{\sigma}$ for some $k'\in\{0,1\}$.
This implies that $(e_{i,1-j,k'},t^{\rightarrow})$ precedes $(e_{i,j,k},t^{\rightarrow})$ within the tie-breaking rule.
Consequently, $F_{i,j}$ cannot be closed or $t^{\rightarrow}$-halfopen and is hence $t^{\leftarrow}$-halfopen.
Furthermore, this implies that $F_{i,1-j}=F_{i,\indbit_{i+1}}$ was also $t^{\leftarrow}$-halfopen for $\sigma^{(3)}$.
Therefore, no cycle center of level $i$ was closed at the end of phase 1.
Thus, by \Cref{corollary: Selection Vertices In Phase One}, $\sigma(g_i)=\sigma^{(3)}(g_i)=\canstrat(g_i)=F_{i,1-j}$.
The statement thus follows by \Cref{lemma: No CC Closed Phase Three}.

Next, assume that $F_{i,1-j}$ is $g_1$-halfopen.
Then $F_{i,j}$ cannot be $g_1$-halfopen since the tie-breaking rule would then choose to apply an improving switch involving $F_{i,1-j}$ as $1-j=\indbit_{i+1}$.
Thus consider the case that $F_{i,j}$ is closed.
We show that we can apply \Cref{lemma: MDP Phase 3 Open Closed Cycle Center}.
By assumption, $I_{\sigma}\cap\mathbb{D}^0=\{e\},$ hence there is no other improving switch $(d_{*,*,*}, e_{*,*,*})$.
As $\indbit_i=0$ and since $F_{i,1-j}$ is $t^{\leftarrow}$-halfopen, we also need to prove $\sigma(g_i)=F_{i,j}$.
This however follows by \Cref{corollary: Selection Vertices In Phase One} if $F_{i,j}$ is closed during phase~1 resp. \Cref{definition: Canonical Strategy MDP} if it was already closed with respect to $\canstrat$.
Since $\sigma$ is a phase-3-strategy, $i'>i>\nsb$ implies $\sigma(s_{i,1-\indbit_{i'+1}})=b_1$ by \Pref{USV1}$_i$.
Now, let $i'>i$ and $\indbit_{i'}=0$.
Then, due to the tie-breaking rule, all improving switches $(e_{i',*,*},b_2)\in \mathbb{E}^1$ have already been applied.
Since $\indbit_{i'}=0$, the cycle center $F_{i',\indbit_{i'+1}}$ cannot have been closed with respect to $\sigma^{(3)}$.
If both cycle centers of level $i'$ were $t^{\leftarrow}$-halfopen for $\sigma^{(3)}$, then they are mixed for $\sigma$, and, in addition, $\sigma(g_{i'})=\sigma^{(3)}(g_{i'})=\canstrat(g_{i'})=F_{i',\indbit_{i'+1}}$.
If the cycle center $F_{i',1-\indbit_{i'+1}}$ is closed for $\sigma^{(3)}$, then $F_{i',\indbit_{i'+1}}$ can only be $t^{\leftarrow}$-halfopen for $\sigma^{(3)}$.
Consequently, by \Cref{corollary: Selection Vertices In Phase One} resp. \Cref{definition: Canonical Strategy MDP} and our previous arguments, this implies $\sigmabar(g_{i'})=1-\indbit_{i'+1}$.
Furthermore, $F_{i',1-\indbit_{i'+1}}$ is then $t^{\rightarrow}$-open and $F_{i',\indbit_{i'+1}}$ is $t^{\rightarrow}$-halfopen (for $\sigma$).
Similarly, if $i'>i$ and $\indbit_{i'+1}=1$, then $F_{i',1-\indbit_{i'+1}}$ is $t^{\rightarrow}$-open if it was closed for $\sigma^{(3)}$ and mixed if it was $t^{\leftarrow}$-halfopen.
Hence, all requirements of \Cref{lemma: MDP Phase 3 Open Closed Cycle Center} are met and the statement follows since $i>\nsb$.

Finally, assume that $F_{i,j}$ is $t^{\rightarrow}$-halfopen.
If we can prove that $\sigma(g_i)=F_{i,j}$, then the statement follows by \Cref{lemma: Open Completely In Phase Three}.
This however follows immediately since $F_{i,j}$ can only be $t^{\rightarrow}$-halfopen if it was closed with respect to $\sigma^{(3)}$, implying $\sigma(g_i)=F_{i,j}$ by the same statements used several times before.
\end{proof}

\SigmabarInPhaseThree*

\begin{proof}
The definition of the sets $S_1$ to $S_4$ implies that $\indbit_i=0\vee\indbit_{i+1}\neq j$ for all of the relevant indices.
We thus begin by considering some fixed but arbitrary indices $i\in[n],j,k\in\{0,1\}$ with $\indbit_{i}=0\vee\indbit_{i+1}\neq j$.
Then, due to the previous application of the improving switches during phase $3$, it holds that $(e_{i,j,k},t^{\rightarrow})\in\applied{\sigma^{(3)}}{\sigma}$ if and only if $\sigma^{(3)}(d_{i,j,k})=F_{i,j}$.
Thus, $(e_{i,j,k},t^{\rightarrow})\notin\applied{\sigma^{(3)}}{\sigma}$ if and only if $\sigma^{(3)}(d_{i,j,k})=e_{i,j,k}$.
Since~$e_{i,j,k}$ has an outdegree of 2 by construction, this implies that $\sigma(e_{i,j,k})=t^{\leftarrow}$ if and only if $\sigma^{(3)}(d_{i,j,k})=e_{i,j,k}$.
In particular, due to $\indbit_i=0\vee\indbit_{i+1}\neq j$, the switch $(d_{i,j,k}, e_{i,j,k})$ was then also applied.
Hence, if there is a $k'\in\{0,1\}$ with $\sigma^{(3)}(d_{i,j,k'})=e_{i,j,k'}$, then $\sigmabar(eg_{i,j})$ if $\nsb>1$ resp. $\sigmabar(eb_{i,j})$ if $\nsb=1$.

Now, consider some fixed indices $i\in[n],j\in\{0,1\}$ and the corresponding cycle center $F_{i,j}$.
Since every cycle center is closed or escapes to $t^{\rightarrow}$ with respect to $\sigma$, either $\sigmabar(eb_{i,j})\wedge\sigmabar(eg_{i,j})$ or $\sigmabar(eb_{i,j})\wedge\nsigmabar(eg_{i,j})$ or $\sigmabar(d_{i,j})$ if $\nsb>1$.
Similarly, if $\nsb=1$, either $\sigmabar(eg_{i,j})\wedge\sigmabar(eb_{i,j})$ or $\sigmabar(eg_{i,j})\wedge\nsigmabar(eb_{i,j})$ or $\sigmabar(d_{i,j})$.
Consequently, $\sigmabar(eb_{i,j})\wedge\sigmabar(eg_{i,j})$ holds if and only if there is a$k\in\{0,1\}$ such that that $\canstrat(d_{i,j,k})\neq F_{i,j}$ and $(d_{i,j,k},F_{i,j})$ was \emph{not} applied during phase~$1$.
By \Cref{lemma: Extended Reaching phase 3}, all improving switches of the type $(d_{*,*,*},F_{*,*})$ not applied in phase~1 had $\occrec^{\canstrat}(d_{*,*,*},F_{*,*})=\maxocc$.
By \Cref{corollary: No Open CC In Phase Three}, it thus suffices to prove that there is a $k\in\{0,1\}$ with $\occrec^{\canstrat}(d_{i,j,k},F_{i,j})=\maxocc$ to prove $\sigmabar(eb_{i,j})\wedge\sigmabar(eg_{i,j})$.
Analogously, to prove $\sigmabar(eb_{i,j})\wedge\nsigmabar(eg_{i,j})$ resp. $\sigmabar(eg_{i,j})\wedge\nsigmabar(eb_{i,j})$, it suffices to show that $F_{i,j}$ was closed at the end of phase~1.

Let $\nsb>1$.
Let $m=\max\{i\colon\sigma(b_i)=g_i\}$ and $u=\min\{i\colon\sigma(b_i)=b_{i+1}\}$.
\begin{enumerate}
	\item We prove that $\occrec^{\canstrat}(d_{i,j,0},F_{i,j})=\maxocc$ for all $(i,j)\in S_2$.
		\begin{itemize}
			\item Let $i\leq\nsb-1, j=\indbit_{i+1}$ and $k\in\{0,1\}$.
				Then, $\bit_{i+1}\neq(\bit+1)_{i+1}=\indbit_{i+1}$ by the choice of $i$.
				In particular, $j\neq\bit_{i+1}$.
				Thus, there is a feasible $t_{\bit}$ for $\bit$ with \[\occrec^{\canstrat}(d_{i,j,k},F_{i,j})=\min\left(\floor{\frac{\bit+1-k}{2}},\ell^{\bit}(i,j,k)+t_{\bit}\right ).\]
				However, the choice of $i$ implies $\bit_i=1$ and thus $t_{\bit}=0$ is the only feasible parameter.
				It thus suffices to show $\ell^{\bit}(i,j,0)\geq\maxocc$.
				Since $\bit_i=1$ and $j\neq\bit_{i+1}$, this follows from \Cref{lemma: Numerics Of Ell}.
			\item Let $i\in\{\nsb+1,\dots,m\},\indbit_i=1$ and $j=1-\indbit_{i+1}$.
				Since $i>\nsb$ implies $\indbit_{i}=\bit_i$ and $\indbit_{i+1}=\bit_{i+1}$, we can deduce $\ell^{\bit}(i,j,0)\geq\maxocc$ as in the previous case.
			\item Let $i\in\{\nsb,\dots,m-1\}\wedge\indbit_i=0$ and $j=\indbit_{i+1}$.
				Since $i+1>\nsb$ implies that we have $\indbit_{i+1}=\bit_{i+1}, \bit_{\nsb-1}=1$ and $\nsb\geq 2$, we obtain $\ell^{\bit}(i,j,0)>\maxocc+1$ as \begin{align*}
					\ell^{\bit}(i,j,0)&=\ceil{\frac{\bit+2^{i-1}+\sum(\bit,i)+1}{2}}\geq\ceil{\frac{\bit+2^{i-1}+2^{\nsb-2}+1}{2}}\\
						&\geq\ceil{\frac{\bit+2^{\nsb-1}+2^{\nsb-2}+1}{2}}\geq\ceil{\frac{\bit+4}{2}}=\floor{\frac{\bit+5}{2}}>\floor{\frac{\bit+1}{2}}+1.
				\end{align*}

				Thus, $\ell^{\bit}(i,j,0)+t_{\bit}>\maxocc$ for every $t_{\bit}$ feasible for $\bit$, implying the statement.
			\item Let $i>m$ and $j\in\{0,1\}$.
				Then, $\lastflip{\bit}{i+1}{}=\lastunflip{\bit}{i+1}{}=0$ since $\bit'_i=0$ for all $\bit'\leq\bit$.
				Hence, by \Cref{lemma: Numerics Of Ell}, $\ell^{\bit}(i,j,k)\geq \bit$.
				Consequently, $\occrec^{\canstrat}(d_{i,j,0},F_{i,j})=\maxocc$
			\item Let $\bit+1=2^{l}$ for some $l\in\mathbb{N}$.
				Then $\nsb=l+1$ and $\bit_{\nsb}=0$.
				This implies $\lastflip{\bit}{\nsb}{\{(\nsb,1)\}}=\lastflip{\bit}{\nsb+1}{}=\lastunflip{\bit}{\nsb+1}{}=0$ and consequently $\occrec^{\canstrat}(d_{i,j,0},F_{i,j})=\maxocc$.
		\end{itemize}
	\item We prove that either $\canstratbar(d_{i,j})$ or $\occrec^{\canstrat}(d_{i,j,k},F_{i,j})<\maxocc$ for both $k\in\{0,1\}$ holds for all $(i,j)\in S_1$.
		\begin{itemize}
			\item Let $i\leq\nsb-1$ and $j=1-\indbit_{i+1}$.
				Then $\bit_i=1$ and $j=1-\indbit_{i+1}=\bit_{i+1}$.
				Hence $F_{i,j}$ was closed with respect to $\canstrat$.
			\item Let $i\in\{\nsb,\dots,m-1\}, \indbit_{i}=0, j=1-\indbit_{i+1}$ and $k\in\{0,1\}$.
			 	Then $\bit_i=\indbit_{i}=0, \indbit_{i+1}=\bit_{i+1}$ and $\indbit_i=0$ implies $i\neq\nsb$.
			 	In particular, $\nsb\leq i-1$ and $\bit_\nsb=0$.
			 	Using \Cref{lemma: Numerics Of Ell}, this implies $\ell^{\bit}(i,j,k)\leq\floor{(\bit+1-k)/2}-1$.
			 	Rearranging this yields $\occrec^{\canstrat}(d_{i,j,k},F_{i,j})\leq\ell^{\bit}(i,j,1)+1$.
			 	If this inequality is strict, the statement follows.
			 	If the inequality is tight, then $\canstrat(d_{i,j,k})=F_{i,j}$ by \Pref{OR2}$_{i,j,k}$ and thus $\occrec^{\canstrat}(d_{i,j,k},F_{i,j})<\maxocc$ by \Pref{OR1}$_{i,j,k}$.
			\item Assume that there is no $l\in\mathbb{N}$ with $\bit+1=2^l$ and let $i=\nsb$ and $j=1-\bit_{\nsb+1}$.
				Since $\bit$ is odd, \Pref{OR3}$_{i,j,0}$ implies $\occrec^{\canstrat}(d_{i,j,0},F_{i,j})<\maxocc.$
				For $k=1$, $\bit$ being odd implies $\occrec^{\canstrat}(d_{i,j,1},F_{i,j})\leq\floor{(\bit+1-1)/2}<\maxocc$.
		\end{itemize}
\end{enumerate}

We now consider the case $\nsb=1$, implying $\bit_i=(\bit+1)_i$ for all $i>1$.
\begin{enumerate}
	\item We prove that $\occrec^{\canstrat}(d_{i,j,0},F_{i,j})=\maxocc$ for all $(i,j)\in S_3$.
		\begin{itemize}
			\item Let $i\in[u]$ and $j=1-\indbit_{i+1}$.
				By the definition of $u$, it holds that $\indbit_i=\bit_i=1$ if $i<u\wedge i\neq 1$ and $\bit_i=0$ if $i=u\vee i=1$.
				In either case, $j=1-\bit_{i+1}$.
				Hence, in the first case, \Cref{lemma: Numerics Of Ell} implies \[\ell^{\bit}(i,j,0)=\ceil{\frac{\bit+\sum(\bit,i)+1}{2}}\geq\ceil{\frac{\bit+1}{2}}\geq\floor{\frac{\bit+1}{2}}.\]
				This implies $\occrec^{\canstrat}(d_{i,j,0},F_{i,j})=\maxocc$ since~$-1$ is not a feasible parameter as~$\bit$ is even.
				Consider the second case, implying \[\ell^{\bit}(i,j,0)=\ceil{\frac{\bit-2^{i-1}+\sum(\bit,i)+1}{2}}.\]
				If $i=1$, then $\ell^{\bit}(i,j,0)=\ceil{\bit/2}=\maxocc$.
				If $i=u$, then $\bit_{i'}=1$ for all indices $i'\in\{2,\dots,u-1\}$ and $\bit_1=0$.
				This implies $\ell^{\bit}(i,j,0)=\maxocc$, and hence the statement since $\bit$ is even.
			\item  Let $i\in\{u+1,\dots,m\}, \indbit_i=1$ and $j=1-\indbit_{i+1}$.
				Then $i\geq 2, \bit_i=1$ and $j=1-\bit_{i+1}$.
				Thus, $\occrec^{\canstrat}(d_{i,j,0},F_{i,j})=\maxocc$ follows by the same arguments used in the last case.
			\item Let $i\in\{u+1,\dots,m-1\}, \indbit_i=0$ and $j=\indbit_{i+1}$.
				Then $i\geq2$ as well as $\bit_i=0$ and $j=\bit_{i+1}$ and $\occrec^{\canstrat}(d_{i,j,0},F_{i,j})=\maxocc$ follows from \Cref{lemma: Numerics Of Ell} and  \[\ell^{\bit}(i,j,0)=\ceil{\frac{\bit+2^{i-1}+\sum(\bit,i)+1}{2}}\geq\ceil{\frac{\bit+3}{2}}\geq\floor{\frac{\bit+1}{2}}+1.\]
			\item Let $i>m$ and $j\in\{0,1\}$.
				Then, $\id_{j=0}\lastflip{\bit}{i+1}{}+\id_{j=1}\lastunflip{\bit}{i+1}{}=0$ by the definition of $m$.
				Hence, by \Cref{lemma: Numerics Of Ell}, $\ell^{\bit}(i,j,k)\geq\bit$.
				This implies $\occrec^{\canstrat}(d_{i,j,0},F_{i,j})=\maxocc$.
			\item Finally consider the pair $(u,\indbit^{\sigma}_{u+1})$. 
				Then, by definition, $\indbit_u=0$ and $\indbit_{u+1}=\bit_{u+1}$.
				If $u>1$, the statement follows as in the third case.
				The case $u=1$ is not possible since $\nsb=1$.
			\end{itemize}
	\item We prove that $\occrec^{\canstrat}(d_{i,j,k},F_{i,j})<\maxocc$ for both $k\in\{0,1\}$ for all $(i,j)\in S_4$.
		First, $(i,j)\in S_4$ implies $i\in\{u+1,\dots,m-1\}, \indbit_i=0$ and $j=1-\indbit_{i+1}$.
		Since $i>u$ implies $i>1$, we have $\bit_i=0$ and $j=1-\bit_{i+1}$.
		Consequently, by \Cref{lemma: Numerics Of Ell},\begin{align*}
			\ell^{\bit}(i,j,k)&=\ceil{\frac{\bit-2^{i-1}+\sum(\bit,i)+1-k}{2}}\\
				&=\ceil{\frac{\bit-2^{i-1}+\sum_{l=2}^{u-1}2^{l-1}+\sum_{l=u+1}^{i-1}\bit_l2^{l-1}+1-k}{2}}.
		\end{align*}
		We prove that this implies $\ell^{\bit}(i,j,k)<\maxocc$, implying the statement as we then either have  $\occrec^{\canstrat}(d_{i,j,k},F_{i,j})<\maxocc$ or $\occrec^{\canstrat}(d_{i,j,k},F_{i,j})=\ell^{\bit}(i,j,k)+1$.
		If $u=2$, then \begin{align*}
			\ell^{\bit}(i,j,k)&=\ceil{\frac{\bit-2^{i-1}+\sum_{l=3}^{i-1}\bit_l2^{l-1}+1-k}{2}}\leq\ceil{\frac{\bit-2^{i-1}+2^{i-1}-4+1-k}{2}}\\
				&=\ceil{\frac{\bit-3-k}{2}}=\ceil{\frac{\bit-1-k}{2}}-1\leq\floor{\frac{\bit+1-k}{2}}-1\\
				&\leq\floor{\frac{\bit+1}{2}}-1<\floor{\frac{\bit+1}{2}}.
		\end{align*}
		If $u>2$, then \begin{align*}
			\ell^{\bit}(i,j,k)&\leq\ceil{\frac{\bit-2^{i-1}+\sum_{l=2}^{u-1}2^{l-1}+\sum_{l=u+1}^{i-1}2^{l-1}+1-k}{2}}\\
				&=\ceil{\frac{\bit-2^{i-1}+2^{u-1}-2+2^{i-1}-2^u+1-k}{2}}=\ceil{\frac{\bit-2^{u-1}+1-k}{2}}\\
				&\leq\ceil{\frac{\bit-4+1-k}{2}}<\floor{\frac{\bit+1}{2}}.
		\end{align*}\qedhere
\end{enumerate}
\end{proof}

\ApplicationCycleEdgesInPhaseThreeMDP*

\begin{proof}
As a reminder, we have $i=\nsb-1, j=1-\bit_{i+1}$ and $k\in\{0,1\}$.
There are no other indices $i',j',k'$ with $(d_{i',j',k'},e_{i',j',k'})\in I_{\sigmae}$.
Also, since no such switch was applied previously in any level below level $i$, the cycle center $F_{i,j}$ is closed for $\sigmae$ as it was closed for $\sigma^{(3)}$ by \Cref{lemma: Extended Reaching phase 3}.
As $i<\nsb$ and $\bit_i=1\wedge\bit_{i+1}\neq\indbit^{\sigmae}_{i+1}$, \Cref{definition: Canonical Strategy MDP} implies that $\canstrat(g_i)=F_{i,j}$.
By the same arguments used when discussing the case $G_n=S_n$ resp. \Cref{claim: Sigmabar in Phase Three}, it can be proven that $F_{i,1-j}$ was not closed during phase~1 as $(i,1-j)\in S_2$.
Consequently, $\sigmae(g_i)=F_{i,j}$ follows from \Cref{corollary: Selection Vertices In Phase One}.
By the tie-breaking rule, no improving switch involving $F_{i,1-j}$ was applied yet.
Therefore, $\smash{\sigmae(e_{i,1-j,*})=\sigma^{(3)}(e_{i,1-j,*})=g_1}$ as well as $\smash{\sigmae(d_{i,1-j,*})=\sigma^{(3)}(d_{i,1-j,*})}$.
By \Cref{corollary: No Open CC In Phase Three}, $F_{i,1-j}$ cannot be open for $\sigma^{(3)}$, so it is not open for $\sigmae$.
Therefore, as $\indbit_i=0$ and $1-j=\indbit_{i+1}$, it is $g_1$-halfopen.
Thus, the first requirement of \Cref{lemma: MDP Phase 3 Open Closed Cycle Center} is met.

By \Cref{lemma: Extended Reaching phase 3} and since $(s_{i',*},h_{i',*})\notin\applied{\sigma^{(3)}}{\sigmae}$ for all $i'<\nsb$, it follows that we have $\sigmae(s_{i',*})=\sigma^{(3)}(s_{i',*})=h_{i',*}$ for all $i'<\nsb$.
Furthermore, $i'<\nsb$ implies $\bit_{i'}=1$ and no improving switch $(d_{*,*,*},e_{*,*,*})$ below level $\nsb$ was applied yet.
Consequently, $\sigmaebar(d_{i'})$ for all $i'<\nsb$.
Now consider some cycle center $F_{i',j'}$ where $i'<i$ and $j'=1-\sigmaebar(g_{i'})$.
We prove that $F_{i',j'}$ is $g_1$-halfopen.
The cycle center $\smash{F_{i',\indbit_{i'+1}}}$ is not closed while $\smash{F_{i',1-\indbit_{i'+1}}}$ is closed due to $1-\indbit_{i'}=\bit_{i'}$.
Thus, by \Cref{corollary: Selection Vertices In Phase One} and the same arguments used before, $\sigmaebar(g_{i'})=\canstratbar(g_{i'})=1-\indbit_{i'+1}$ and, in particular, $j'=\indbit_{i'+1}$.
However, by \Cref{corollary: No Open CC In Phase Three} and the tie-breaking rule, this implies that $F_{i',j'}$ is $g_1$-halfopen as before.
Thus, the second requirement is met.

The third requirement is met as $i'>i=\nsb-1$ and since $\sigmae$ has \Pref{USV1}$_{i'}$.

Consider the fourth requirement.
Let $i'>i$ and $\indbit_{i'}=0$.
Then, due to the tie-breaking rule, all improving switches $(e_{i',j',k'},b_2)$ with $\sigma^{(3)}(d_{i',j',k'})=F_{i',j'}$ have already been applied.
Since $\indbit_{i'}=0$, $F_{i',\indbit_{i'+1}}$ cannot have been closed for $\sigma^{(3)}$.
If both cycle centers of level $i'$ were $g_1$-halfopen for $\sigma^{(3)}$, then they are mixed for $\sigma$, and $\sigma(g_{i'})=\sigma^{(3)}(g_{i'})=\canstrat(g_{i'})=F_{i',\indbit_{i'+1}}$.
If $F_{i',1-\indbit_{i'+1}}$ is closed for $\sigma^{(3)}$, then $F_{i',\indbit_{i'+1}}$ can only be $g_1$-halfopen for $\sigma^{(3)}$.
Consequently, by \Cref{corollary: Selection Vertices In Phase One} resp. \Cref{definition: Canonical Strategy MDP}, $\sigmabar(g_{i'})=1-\indbit_{i'+1}$.
Furthermore, $F_{i',1-\indbit_{i'+1}}$ is then $b_2$-open and $F_{i',\indbit_{i'+1}}$ is $b_2$-halfopen (for $\sigma$).
Thus, the fourth requirement is met.

By the same argument, if $i'>i$ and $\indbit_{i'+1}=1$, then $\smash{F_{i',1-\indbit_{i'+1}}}$ is $b_2$-open if it was closed for~$\sigma^{(3)}$ and mixed if it was $g_1$-halfopen.
Thus, the fifth and final requirement is met.
\end{proof}

\FirstSwitchPhaseFour*

\begin{proof}
We first consider the case that $\bit+1$ is a power of two, implying $\bit=2^{\nsb-1}-1$.
We distinguish four kinds of improving switches.
\begin{enumerate}
	\item Let $e=(s_{\nsb-1,0},b_1)$.
		Then, $\occrec^{\sigma}(e)=0$ follows from \[\occrec^{\sigma}(e)=\flips{\bit}{\nsb}{}=\floor{\frac{\bit+2^{\nsb-1}}{2^{\nsb}}}<\floor{\frac{2^{\nsb-1}+2^{\nsb-1}}{2^{\nsb}}}.\]
	\item Let $e=(s_{i,1},b_1)$ for $i\leq\nsb-2$.
		Then, $\occrec^{\sigma}(e)=\flips{\bit}{i+1}{}-j\cdot\bit_{i+1}=\floor{(\bit+2^{i})/2^{i+1}}-1.$
		If $i=\nsb-2$, then $\nsb\geq 3$ and \[\occrec^{\sigma}(e)=\floor{\frac{2^{\nsb-1}-1+2^{\nsb-2}}{2^{\nsb-1}}}-1=\floor{1+\frac{2^{\nsb-2}-1}{2^{\nsb-1}}}-1=0.\]
		If $i\leq\nsb-3$, then $\nsb\geq 4$ and \[\occrec^{\sigma}(e)\geq\floor{\frac{2^{\nsb-1}-1+2^{\nsb-3}}{2^{\nsb-2}}}-1=\floor{2+\frac{2^{\nsb-3}-1}{2^{\nsb-2}}}-1=1.\]
	\item Let $e=(d_{i,j,k},F_{i,j})$ for some indices $i\in[n], j,k\in\{0,1\}$ with $\sigma(e_{i,j,k})=g_1$.
		Then $\occrec^{\sigma}(e)=\maxocc\geq 1$.
	\item Let $e=(e_{i,j,k},b_2)$ for some indices $i\in[n], j,k\in\{0,1\}$ with $\sigma(e_{i,j,k})=g_1$.
		Then, $\occrec^{\sigma}(e)=\floor{\bit/2}\geq 1$ if $\bit>1$ and $\occrec^{\sigma}(e)=0$ if $\bit=1$.
\end{enumerate}
Thus, $(s_{\nsb-1,0},b_1)$ and $(s_{\nsb-2,1},b_1)$ both minimize the occurrence record if $\bit>1$.
If $\bit=1$, then all switches $(e_{i,j,k},b_2)$ with $i\in[n], j,k\in\{0,1\}$ and $\sigma(e_{i,j,k})=g_1$ also minimize the occurrence record.
Due to the tie-breaking rule, $(s_{\nsb-1,0},b_1)$ is thus applied next in either case.

Now consider the case that $\bit+1$ is not a power of two.
Then $\bit\geq2^{\nsb}+2^{\nsb-1}-1$ and $\bit\geq 6$, implying $\floor{(\bit+2)/4}<\maxocc$ and $\floor{(\bit+2)/4}<\floor{\bit/2}$.
We prove that $(s_{\nsb-1,0},b_1)$ minimizes the occurrence record.
\begin{enumerate}
	\item Let $e=(s_{\nsb-1,0},b_1)$.
		Then, $\occrec^{\sigma}(e)=\flips{\bit}{\nsb}{}=\floor{(\bit+2^{\nsb-1})/2^{\nsb}}\leq\floor{(\bit+2)/4}$ as $\bit\geq6$ implies $\nsb\geq 2$.
	\item Let $e=(d_{i,j,k},F_{i,j})$ with $i\in[n],j,k\in\{0,1\}$ and $\sigma(e_{i,j,k})=g_1$.
		Then $\occrec^{\sigma}(e)=\maxocc$, implying that $\occrec^{\sigma}(e)>\occrec^{\sigma}(s_{\nsb-1,0},b_1)$.
	\item Let $e=(e_{i,j,k},b_2)$ with $i\in[n],j,k\in\{0,1\}$ and $\sigma(e_{i,j,k})=g_1$.
		Then $\occrec^{\sigma}(e)=\floor{\bit/2}$ by \Cref{table: Occurrence Records}, implying that $\occrec^{\sigma}(e)>\occrec^{\sigma}(s_{\nsb-1,0},b_1)$.
	\item Let $e=(s_{i,1},b_1)$ with $i\leq\nsb-2$.
		Then, $\occrec^{\sigma}(e)=\flips{\bit}{i+1}{}-\bit_{i+1}=\flips{\bit}{i+1}{}-1$ by \Cref{table: Occurrence Records}.
		Hence, $\occrec^{\sigma}(e)=\flips{\bit}{i+1}{}-1>\flips{\bit}{\nsb}{}-1>\occrec^{\sigma}(s_{\nsb-1,0},b_1)-1$ by \Cref{lemma: Numerics Of OR}.
		Thus, by integrality, $\occrec^{\sigma}(e)\geq\occrec^{\sigma}(s_{\nsb-1,0},b_1)$.
	\item Let $e=(d_{i,j,k},F_{i,j})$, with $i=\nsb,j=1-\indbit_{i+1}$ and $k\in\{0,1\}$.
		By the definition of a canonical strategy, $\canstrat(d_{i,j,k})\neq F_{i,j}$.
		Hence $\occrec^{\canstrat}(e)=\min(\floor{\frac{\bit+1-k}{2}},\ell^{\bit}(i,j,k)+t_\bit)$, where $t_{\bit}$ is feasible for $\bit$.
		Since $\bit_i=\bit_{\nsb}=0$ and $\indbit_{i+1}=\bit_{i+1}$, \Cref{lemma: Numerics Of Ell} then implies \begin{align*}
			\ell^{\bit}(i,j,k)&=\ceil{\frac{\bit-2^{i-1}+\sum(\bit,i)+1-k}{2}}\\
				&=\ceil{\frac{\bit-2^{i-1}+2^{i-1}-1+1-k}{2}}=\floor{\frac{\bit+1-k}{2}}.
		\end{align*}
		Hence, by \Pref{OR3}$_{i,j,k}$, $\occrec^{\canstrat}(e)=\ell^{\bit}(i,j,k)-1=\maxocc-1$ for $k=0$ and $\occrec^{\canstrat}(e)=\ell^{\bit}(i,j,k)=\maxocc-1$ for $k=1$.
		If $k=1$, \Cref{corollary: Switches With Low OR In Phase One} implies that the edge~$e$ was applied during phase~1.
		Consequently, $\occrec^{\sigma}(e)=\maxocc>\occrec^{\sigma}(s_{\nsb-1,0},b_1)$.
	\item Let $e=(d_{i,j,k},F_{i,j})$ with $i\in\{\nsb+1,\dots,m-1\}, \indbit_i=0, j=1-\indbit_{i+1}$ and $k\in\{0,1\}$.
		By the choice of $i$, it then follows that $\bit_i=0$ and $j=1-\bit_{i+1}$.
		If we have $\occrec^{\canstrat}(e)=\floor{(\bit+1-k)/2}$, then it either holds that $\occrec^{\canstrat}(e)=\maxocc-1$ or $\occrec^{\canstrat}(e)=\maxocc.$
		In both cases, $\bit\geq 6$ implies that $\occrec^{\sigma}(s_{\nsb-1,0},b_1)\leq\occrec^{\canstrat}(e)\leq\occrec^{\sigma}(e)$.		
		Thus assume $\occrec^{\canstrat}(e)=\ell^{\bit}(i,j,k)+t_{\bit}$ for some $t_{\bit}$ feasible for $\bit$ and $\occrec^{\canstrat}(e)\neq\floor{(\bit+1-k)/2}$.
		By \Cref{lemma: Numerics Of Ell}, \begin{align*}
			\ell^{\bit}(i,j,k)&=\ceil{\frac{\bit-2^{i-1}+\sum(\bit,i)+1-k}{2}}\geq\ceil{\frac{\bit-2^{i-1}+2^{\nsb-1}-k}{2}}\\
				&=\floor{\frac{\bit-2^{i-1}+2^{\nsb-1}+1-k}{2}}.
		\end{align*}
		We prove that $\ell^{\bit}(i,j,k)>\floor{(\bit+2)/4}$, implying $\occrec^{\sigma}(e)\geq\occrec^{\sigma}(s_{\nsb-1,0},b_1)$.
		We begin by observing \begin{align*}
			\ell^{\bit}(i,j,k)&=\floor{\frac{\bit-2^{i-1}+2^{\nsb-1}+1-k}{2}}\\
				&=\floor{\frac{2\bit-2^i+2^{\nsb}+2-2k}{4}}\geq\floor{\frac{2\bit-2^i+2^{\nsb}}{4}}.
		\end{align*}
		By the choice of $i$ and $j$, the cycle center $F_{i,j}$ was closed at least once during some previous transition.
		But, since $\bit_i=0$, the cycle center was also opened again later.
		This implies $\bit\geq2^{i-1}+2^{i-1}+2^{\nsb-1}-1=2^{i}+2^{\nsb-1}-1$.		
		Thus, \begin{align*}
			2\bit-2^i+2^\nsb-[\bit+2]&=\bit-2^i+2^{\nsb}-2\\		
				&\geq2^{i}+2^{\nsb-1}-1-2^{i}+2^{\nsb}-2\\
				&=2^{\nsb}+2^{\nsb-1}-3\geq 4+1-3=2.
		\end{align*}
		Since $2\bit-2^{i}+2^{\nsb}$ is even and larger than 0 and since $\bit$ being odd implies $\bit+2$ being odd, this difference is at least 3.
		It is easy to show that, in general, $x$ being even and larger than 0, $y$ being odd and $x-y\geq 3$ implies $\floor{x/4}>\floor{y/4}$.
		This yields $\ell^{\bit}(i,j,k)>\floor{(\bit+2)/4}$.
\end{enumerate}

It remains to prove that \Cref{lemma: Phase Four Complete} describes the application of $e$.
Since $\sigma$ is a phase-4-strategy and since $i'>i=\nsb-1$ implies $i'\geq\nsb$, $\sigma$ has \Pref{USV1}$_{i'}$ for all $i'>i$.
By \Cref{lemma: Extended Reaching phase 4 or 5}, it follows that $\sigma$ also meets the other requirements of \Cref{lemma: Phase Four Complete}.
\end{proof}

\BoundedORBeginningPhaseFive*

\begin{proof}
By \Cref{lemma: Extended Reaching phase 5}, the set of improving switches can be partitioned as follows:
\begin{enumerate}
	\item Let $e=(d_{i,j,k},F_{i,j})$ resp. $e=(e_{i,j,k},g_1)$ with $i\in[n],j,k\in\{0,1\}$ and $\sigma(e_{i,j,k})=b_2$.
		Then $\occrec^{\sigma}(e)=\maxocc$ resp.  $\occrec^{\sigma}(e)=\occrec^{\canstrat}(e)=\ceil{\bit/2}=\maxocc$ by \Cref{lemma: Extended Reaching phase 5}.
	\item Let $e=(d_{i,j,k},F_{i,j})$ with $\indbit_i=0, i\in\{u+1,\dots,m-1\}, j=1-\indbit_{i+1}$ and $k\in\{0,1\}$.
		Then, $\bit_i=0$ and $j=1-\bit_{i+1}$ since $i\geq u+1>1$ and $\nsb=1$.
		In addition, $\bit_1=0$ and, due to $i>u$, there is at least one $l\in\{2,\dots,i-1\}$ with $(\bit+1)_l=\bit_l=0$.
		Consequently, \Cref{lemma: Numerics Of Ell} yields \[\ell^{\bit}(i,j,k)=\ceil{\frac{\bit-2^{i-1}+\sum(\bit,i)+1-k}{2}}\leq\ceil{\frac{\bit-3-k}{2}}=\floor{\frac{\bit-k}{2}}-1.\]
		Since there is a $t_{\bit}$ feasible for $\bit$,  $\occrec^{\canstrat}(e)=\min(\floor{(\bit+1-k)/2},\ell^{\bit}(i,j,k)+t_{\bit})$.
		We thus distinguish the following cases.
		\begin{enumerate}
			\item Let $\occrec^{\canstrat}(e)=\ell^{\bit}(i,j,k)+1$.
				Then, by \Pref{OR2}$_{i,j,k}$, $\canstrat(d_{i,j,k})=F_{i,j}$ and $e$ was not applied during $\canstrat\to\sigma$ as switches of this type were only applied during phase $1$ so far.
				Consequently, $\occrec^{\sigma}(e)=\occrec^{\canstrat}(e)<\maxocc$ by  \Pref{OR1}$_{i,j,k}$.
			\item Let $\occrec^{\canstrat}(e)=\ell^{\bit}(i,j,k)$.
				Then, $\occrec^{\canstrat}(e)\leq\occrec^{\sigma}(e)\leq\floor{(\bit-k)/2}-1<\maxocc$ as well as $\canstrat(d_{i,j,k})\neq F_{i,j}$ by \Pref{OR2}$_{i,j,k}$.
				Using \Pref{OR4}$_{i,j,k}$, this implies $\occrec^{\canstrat}(e)=\maxocc-1$.
				Hence, by \Cref{corollary: Switches With Low OR In Phase One}, $e$ was applied during phase 1.
				Consequently, $\occrec^{\sigma}(e)=\occrec^{\canstrat}(e)+1=\maxocc$.
			\item The case $\occrec^{\canstrat}(e)=\ell^{\bit}(i,j,k)-1$ cannot occur since the parameter $t_{\bit}=-1$ is not feasible as $\bit$ is even.
			\item Let $\occrec^{\canstrat}(e)=\floor{(\bit+1-k)/2}$ but $\floor{(\bit+1-k)/2}\neq\ell^{\bit}(i,j,k),\ell^{\bit}(i,j,k)+1$.
					This implies that we need to have $\floor{(\bit+1-k)/2}<\ell^{\bit}(i,j,k)$ since it holds that  \[\occrec^{\canstrat}(e)=\min(\floor{(\bit+1-k)/2},\ell^{\bit}(i,j,k)+t_{\bit}).\]
					But this is a contradiction since $\ell^{\bit}(i,j,k)\leq\floor{(\bit-k)/2}-1$.
		\end{enumerate}
\end{enumerate}
Since the only improving switches with an occurrence record lower than $\maxocc$ are the switches described in case 2.a), the second part of the statement follows.
\end{proof}

\RowOneInPhaseFive*

\begin{proof}
We currently consider the first phase-$5$-strategy $\sigma$ as described by \Cref{lemma: Extended Reaching phase 5} and an improving switch $e=(d_{i,j,k},F_{i,j})$ with $i\in\{u+1,\dots,m-1\}, \indbit_i=0, j=1-\indbit_{i+1}$ and  $k\in\{0,1\}$ as well as $\canstrat(d_{i,j,k})=F_{i,j}$.
We have to prove $\sigma(b_i)=b_{i+1}, j=1-\indbit_{i+1}, \sigmabar(g_i)=1-\indbit_{i+1}$ and $i\neq 1$.
The first two statements follow directly since $\sigma$ is a phase-$5$-strategy and $\indbit_i=0$ as well as by the choice of $j$.
Also, $i\neq 1$ follows from $i\geq u+1>1$.
It thus suffices to show $\sigmabar(g_i)=1-\indbit_{i+1}$.

For the sake of a contradiction, let $\sigma(g_i)=F_{i,\indbit_{i+1}}$.
Since $\indbit_i=0$ and $\nsb=1$, it holds that $\bit_i=(\bit+1)_i=0$, implying $i\neq\nsb$.
By \Cref{corollary: Selection Vertices In Phase One}, the only improving switch from a selector vertex towards the active cycle center of a level that can be performed during phase~1 is $(g_{\nsb},F_{\nsb,\bit_{\nsb+1}})$.
This implies $(g_i,F_{i,\indbit_{i+1}})\notin\applied{\canstrat}{\sigma}$, hence $\canstrat(g_i)=F_{i,\indbit_{i+1}}$.
If $\canstrat(d_{i,j,1-k})=F_{i,j}$, then $F_{i,j}=F_{i,1-\indbit_{i+1}}$ was closed at the beginning of phase~1 as $\canstrat(d_{i,j,k})=F_{i,j}$.
But this implies $\canstrat(g_i)=F_{i,1-\indbit_{i+1}}$ by the definition of a canonical strategy which is a contradiction.
Thus let $\canstrat(d_{i,j,1-k})\neq F_{i,j}$, implying that we have $\occrec^{\canstrat}(d_{i,j,1-k},F_{i,j})\neq \ell^{\bit}(i,j,1-k)+1$.
Then, by the same arguments used when proving \Cref{claim: Bounded OR Beginning Phase 5}, it follows that $\ell^{\bit}(i,j,1-k)\leq\floor{(\bit-(1-k))/2}-1$.
Also, by these arguments, it cannot happen that $\occrec^{\canstrat}(d_{i,j,1-k})=\floor{(\bit+1-(1-k))/2}\neq\ell^{\bit}(i,j,1-k)$.
Since the parameter $t_{\bit}=-1$ is not feasible, we thus have \[\occrec^{\canstrat}(d_{i,j,1-k},F_{i,j})=\ell^{\bit}(i,j,1-k)\leq\floor{\frac{\bit-(1-k)}{2}}-1<\floor{\frac{\bit+1}{2}}.\]
But this implies that $(d_{i,j,1-k},F_{i,j})$ was applied in phase~1 by \Cref{corollary: Switches With Low OR In Phase One}.
Hence, $F_{i,j}$ was closed in phase 1.
But then, by \Cref{corollary: Selection Vertices In Phase One}, $(g_i,F_{i,j})$ became improving during phase 1 and was thus applied.
This implies $\sigma(g_i)=F_{i,j}=F_{i,1-\indbit_{i+1}}$, contradicting the assumption.
Consequently, $\sigma(g_i)=F_{i,1-\indbit_{i+1}}$.
\end{proof}

\EscapeVerticesPhaseFiveNSBOne*

\begin{proof}
First, we show that $F_{i,j}$ is mixed.
Since $e=(e_{i,j,k},g_1)\in I_{\sigma}$ implies $(d_{i,j,k},F_{i,j})\in I_{\sigma}$ by \Cref{equation: Characterization NSB==1}, we have $\sigmabar(eb_{i,j})$.
In particular, $F_{i,j}$ is not closed, so $\indbit_{i}=0\vee\indbit_{i+1}\neq j$.
Consequently, $(i,j)\in S_3$ or $(i,j)\in S_4$. 
By \Cref{lemma: Extended Reaching phase 5}, $\sigmabar(eb_{i,j})$, and as no improving switch $(e_{*,*,*},b_2)$ was applied during $\sigma^{(5)}\to\sigma$, we need to have $(i,j)\in S_3$, implying the statement.
We now prove that $j=1$ implies $\nsigmabar(eb_{i,1-j})$ if $G_n=S_n$.
Since $j=1$, we need to prove $\nsigmabar(eb_{i,0})$.
If $F_{i,0}$ is closed, then the statement follows.
If $F_{i,0}$ is not closed, then $\indbit_i=0\vee\indbit_{i+1}\neq j$ as $F_{i,1}$ cannot be closed by the choice of $e$.
Consequently, $(i,0)\in S_3$ or $(i,0)\in S_4$.
In the second case, $\nsigmabar^{(5)}(eb_{i,0})$ by \Cref{lemma: Extended Reaching phase 5} and the statement follows as no improving switch $(e_{*,*,*},b_2)$ was applied during $\sigma^{(5)}\to\sigma$.
Consider the case $(i,0)\in S_3$.
Then, by \Cref{lemma: Extended Reaching phase 5}, $F_{i,0}$ and $F_{i,1}$ are mixed with respect to $\sigma^{(5)}$.
Thus, as we consider the case $G_n=S_n$, the tie-breaking rule must have applied the improving switches $(e_{i,0,*},g_1)$ prior to $(e_{i,1,k},g_1)$, implying the statement.
Note that the statement \enquote{$j=1-\indbit_{i+1}\implies\nsigmabar(eb_{i,1-j})$ if $G_n=M_n$} follows by the same arguments and since the tie-breaking rule applies improving switches $(e_{i,\indbit_{i+1},*},g_1)$ first.
\end{proof}

\FulfillingORFiveIfNSBIsOne*

\begin{proof}
Let $\sigma$ denote the strategy before the application of the switch $(e_{i,j,k},g_1)$ and let $\sigmae$ denote the strategy obtained after the application of this switch.
By the same arguments used in the proof of \Cref{claim: Escape Vertices Phase Five NSB One}, it follows that $F_{i,j}$ is mixed with respect to $\sigma$.
By the characterization of $I_{\sigma}$ given in \Cref{equation: Characterization NSB==1}, it holds that $(d_{i,j,k},F_{i,j})\in I_{\sigma}$, implying $\sigma(d_{i,j,k})\neq F_{i,j}$ and $\sigmae(d_{i,j,k})\neq F_{i,j}$.
Since~$\sigma$ is a phase-5-strategy for $\bit$, this furthermore implies $\indbit_i=0\vee\indbit_{i+1}\neq j$.

Assume that $\tilde{e}$ was applied previously in this transition.
It is not possible that $\tilde{e}$ was applied during phase 5 as this would imply $\sigma(d_{i,j,1-k})=F_{i,j}$, contradicting that $F_{i,j}$ is mixed with respect to $\sigma$.
Consequently, $\tilde{e}$ was applied during phase~1.
We thus need to have $\canstrat(d_{i,j,1-k})\neq F_{i,j}$ and $\occrec^{\canstrat}(\tilde{e})\in\{\maxocc-1,\maxocc\}$ by \Pref{OR4}$_{i,j,1-k}$.
This in particular implies $\occrec^{\sigmae}(\tilde{e})\in\{\maxocc, \maxocc+1\}$ and hence the statement.

Now assume that $\tilde{e}$ was not applied previously in this transition, implying $\occrec^{\sigmae}(\tilde{e})=\occrec^{\canstrat}(\tilde{e})$.
Let $\canstrat(d_{i,j,1-k})\neq F_{i,j}$.
Then, by \Pref{OR4}$_{i,j,1-k}$ and \Cref{corollary: Switches With Low OR In Phase One}, it follows that $\occrec^{\canstrat}(\tilde{e})=\maxocc$.
Thus let $\canstrat(d_{i,j,1-k})=F_{i,j}$.
Then, by Properties (\ref{property: OR1})$_{i,j,1-k}$ and (\ref{property: OR2})$_{i,j,1-k}$, it holds that $\occrec^{\canstrat}(\tilde{e})=\ell^{\bit}(i,j,1-k)+1<\maxocc.$
We now prove that this yields a contradiction.
\begin{enumerate}
	\item Let, for the sake of a contradiction, $\indbit_i=1\wedge\indbit_{i+1}\neq j$.
		Assume $\bit_i=0\wedge\bit_{i+1}\neq j$.
		Since $\nsb=1$, this implies $i=\nsb=1$.
		Thus, as $\bit$ is even, \[\ell^{\bit}(i,j,1-k)=\ceil{\frac{\bit-2^0+\sum(\bit,i)+1-(1-k)}{2}}=\floor{\frac{\bit+1}{2}}=\maxocc\] by \Cref{lemma: Numerics Of Ell}.
		But then, $\occrec^{\canstrat}(\tilde{e})=\ell^{\bit}(i,j,1-k)+1=\maxocc+1$ which is a contradiction.
		Assuming $\bit_i=1\wedge\bit_{i+1}\neq j$ also results in a contradiction since \[\ell^{\bit}(i,j,1-k)=\ceil{\frac{\bit+\sum(\bit,i)+1-(1-k)}{2}}\geq\ceil{\frac{\bit+k}{2}}\geq\floor{\frac{\bit+1}{2}}.\]
	\item Let, for the sake of a contradiction, $\indbit_i=0\wedge\indbit_{i+1}=j$.
		Then also $\bit_i=0\wedge\bit_{i+1}=j$ and $i\geq 2$ since $\nsb=1$.
		Then, \Cref{lemma: Numerics Of Ell} implies \[\ell^{\bit}(i,j,1-k)\geq\ceil{\frac{\bit+3-(1-k)}{2}}\geq\floor{\frac{\bit+1}{2}}+1.\]
		This yields a contradiction as before.
	\item Let, for the sake of a contradiction, $\indbit_i=0\wedge\indbit_{i+1}\neq j$.
		This implies that $i\geq u$.
		If $i>m$, then \Cref{lemma: Numerics Of Ell} implies $\ell^{\bit}(i,j,1-k)\geq \bit$, contradicting that we have $\occrec^{\canstrat}(\tilde{e})=\ell^{\bit}(i,j,1-k)+1<\maxocc$.
		We hence may assume $i\in\{u,\dots,m-1\}$.
		If $i\neq u$, then \Cref{lemma: Extended Reaching phase 5} implies $(d_{i,j,1-k},F_{i,j})\in I_{\sigma^{(5)}}$.
		But then, the switch was applied  during $\canstrat\to\sigmae$, contradicting the assumption.
		Hence let $i=u$.
		Then, $\indbit_{i'}=1$ for all $i'<u=i$.
		Consequently, \begin{align*}
			\ell^{\bit}(i,j,1-k)&=\ceil{\frac{\bit-2^{i-1}+\sum(\bit,i)+1-(1-k)}{2}}\\
				&=\ceil{\frac{\bit-2^{i-1}+2^{i-1}-1-1+1-1+k}{2}}=\ceil{\frac{\bit-2+k}{2}}\\
				&=\floor{\frac{\bit-1+k}{2}}.
		\end{align*}
		But this implies $\occrec^{\canstrat}(d_{i,j,1-k})=\ell^{\bit}(i,j,1-k)+1=\maxocc$ and contradicts \Pref{OR1}$_{i,j,1-k}$.\qedhere
\end{enumerate}
\end{proof}

\SVGSVMInPhaseFive*

\begin{proof}
Consider some arbitrary but fixed index $i\in[n]$.
If $\indbit_i=1$, then the statements follow from the definition of a phase-5-strategy.
If $\indbit_i=0$ and $(g_i,F_{i,j})\in\applied{\sigma^{(5)}}{\sigma}$, then this follows from \Cref{corollary: Selection Vertices Phase 5}. 
Thus, let $\indbit_i=0$ and $(g_i,F_{i,j})\notin\applied{\sigma^{(5)}}{\sigma}$, implying $i\neq 1$ since $\nsb=1$.
We now prove the following statement.
If $\sigmabar(g_i)=1$ resp. $\sigmabar(g_i)=1-\indbit_{i+1}$ and $\nsigmabar(d_{i,1})$  resp. $\nsigmabar(d_{i,1-\indbit_{i+1}})$  then $(g_i,F_{i,0})\in I_{\sigma}$ resp. $(g_i,F_{i,\indbit_{i+1}})\in I_{\sigma}$.
This is sufficient to prove the statement as $I_{\sigma}\cap\mathbb{G}=\emptyset$.

Thus, let $j\coloneqq 0$ (if $G_n=S_n$) resp. $j\coloneqq\indbit_{i+1}$ (if $G_n=M_n$) and assume $\nsigmabar(d_{i,1-j})$.
It suffices to prove $\valu_{\sigma}^*(F_{i,j})\succ\valu_{\sigma}^*(F_{i,1-j})$.
Since  $\sigma(e_{i',j',k'})=g_1$ for all $(i',j',k')\neq(1,\indbit_2,k), i\neq 1$ and $\relbit{\sigma}=u\neq 1$, the two cycle centers $F_{i,*}$ are either closed or escape only to $g_1$.

Consider the case $G_n=S_n$.
If both cycle centers escape towards $g_1$, then the statement follows from \begin{align*}
\valu_{\sigma}^\P(F_{i,0})&=\{F_{i,0},d_{i,0,*},e_{i,0,*},b_1\}\cup\valu_{\sigma}^\P(g_1)\\
	&\rhd\{F_{i,1},d_{i,1,*},e_{i,1,*},b_1\}\cup\valu_{\sigma}^\P(g_1)=\valu_{\sigma}^\P(F_{i,1}).
\end{align*}
Since we currently consider the case $\indbit_i=0$, only $F_{i,1-\indbit_{i+1}}$ can be closed, so assume this is the case.
Assume that $0=1-\indbit_{i+1}$, so $j=1-\indbit_{i+1}$.
Then, by \Pref{USV1}$_i$ and $\sigma(b_1)=g_1$, the statement follows since $\valu_{\sigma}^\P(F_{i,0})=\{s_{i,0},b_1\}\cup\valu_{\sigma}^\P(g_1)$ and $\valu_{\sigma}^\P(F_{i,1})=\{F_{i,1},d_{i,1,k},e_{i,1,k},b_1\}\cup\valustar_{\sigma}^\P(g_1)$ for some $k\in\{0,1\}$.
Thus assume that $0=\indbit_{i+1}$, so $j=\indbit_{i+1}$.
Then, the cycle center $F_{i,1-j}=F_{i,1}$ is closed, contradicting the assumption developed at the beginning of the proof.

Consider the case $G_n=M_n$ next and note that we thus have $j=\indbit_{i+1}$ from now on.
If both cycle centers of level $i$ are $g_1$-open or $g_1$-halfopen, then the statement follows by \Cref{lemma: Both CC Open For MDP} since $i>\nsb=1$. 
Thus consider the case that $F_{i,j}$ is $g_1$-open and that $F_{i,1-j}$ is $g_1$-halfopen.
By assumption, $j=\indbit_{i+1}$, implying $\sigma(s_{i,j})=h_{i,j}$ and $\sigma(s_{i,1-j})=b_1$ by \Pref{USV1}$_i$.
Thus, by \Pref{EV1}$_{i+1}$, \[\valu_{\sigma}^\M(F_{i,j})=(1-\e)\valu_{\sigma}^\M(g_1)+\e\left[\rew{s_{i,j},h_{i,j}}+\valu_{\sigma}^\M(b_{i+1})\right]\] and \[\valu_{\sigma}^\M(F_{i,1-j})=\valu_{\sigma}^{M}(g_1)+\frac{2\e}{1+\e}\rew{s_{i,1-j}}.\]
To prove $\valu_{\sigma}^\M(F_{i,j})>\valu_{\sigma}^\M(F_{i,1-j})$, it thus suffices to prove \[\rew{s_{i,j},h_{i,j}}+\valu_{\sigma}^\M(b_{i+1})-\valu_{\sigma}^\M(g_1)-\frac{2}{1+\e}\rew{s_{i,1-j}}>0.\]
This can be shown by an easy but tedious calculation using $\valu_{\sigma}^\M(g_1)=R_1^\M$, $\indbit_{i}=0$, $i+1>\relbit{\sigma}$, and  $\valu_{\sigma}^\M(b_{i+1})=L_{i+1}^\M$.
Now let $F_{i,j}$ be $g_1$-halfopen and $F_{i,1-j}$ be $g_1$-open.
Then, by the same arguments used before, \[\valu_{\sigma}^{\M}(F_{i,j})=\frac{1-\e}{1+\e}\valu_{\sigma}^{\M}(g_1)+\frac{2\e}{1+\e}[\rew{s_{i,j},h_{i,j}}+\valu_{\sigma}^{\M}(b_{i+1})]\] and \[\valu_{\sigma}^{\M}(F_{i,1-j})=\valu_{\sigma}^{\M}(g_1)+\e\rew{s_{i,1-j}}.\]
It thus suffices to prove \[\rew{s_{i,j},h_{i,j}}+\valu_{\sigma}^{\M}(b_{i+1})-\valu_{\sigma}^{\M}(g_1)-\frac{1+\e}{2}\rew{s_{i,1-j}}>0\]which follows analogously.
Since only $F_{i,1-\indbit_{i+1}}=F_{i,1-j}$ can be closed in level $i$, the statement then follows by the same argument used for the case $G_n=S_n$.
\end{proof}

\ImprovingSwitchesForBPlusOneForNSBOne*

\begin{proof}
To simplify notation, let $\sigma\coloneqq\sigma_{\bit+1}$.
Consider the strategy $\sigma^{(5)}$.
Using the characterization of the strategy that was obtained after having applied all switches $(d_{i,j,k},F_{i,j})$ with an occurrence record smaller than $\maxocc$ (see Equation~(\ref{equation: Characterization NSB==1})), we obtain \begin{align*}
I_{\sigma}&=\{(d_{i,j,*},F_{i,j})\colon\sigma^{(5)}(e_{i,j,*})=b_2\}\\
	&\quad\cup\bigcup_{\substack{i=u+1\\\indbit^{\sigma}_i=0}}^{m-1}\left\{e=(d_{i,1-\indbit^{\sigma}_{i+1},*},F_{i,1-\indbit^{\sigma}_{i+1}})\colon\occrec^{\sigma^{(5)}}(e)=\floor{\frac{\bit+1}{2}}\right\}.
\end{align*}
In particular, $I_{\sigma}\subseteq\{(d_{i,j,k},F_{i,j})\colon\sigma(d_{i,j,k})\neq F_{i,j}\}$ and every improving switch has an occurrence record of at least $\maxocc$.
To prove $\{(d_{i,j,k},F_{i,j})\colon\sigma(d_{i,j,k})\neq F_{i,j}\}\subseteq I_{\sigma}$, let $e\coloneqq(d_{i,j,k},F_{i,j})$ with $\sigma(d_{i,j,k})\neq F_{i,j}$.
It suffices to show $\valu_{\sigma}^*(F_{i,j})\succ\valu_{\sigma}^*(e_{i,j,k})$.
\Pref{ESC1} and $\nsb=1$ imply $\sigmabar(eg_{i,j})\wedge\nsigmabar(eb_{i,j})$.
Furthermore, \Pref{REL1} yields $\relbit{\sigma}=\min\{i'\colon\sigma(b_{i'})=b_{i'+1}\}\neq 1$.
This implies $\valu_{\sigma}^\P(F_{i,j})=\{F_{i,j}\}\cup\valu_{\sigma}^\P(e_{i,j,k})$, implying the statement if $G_n=S_n$.
If $G_n=M_n$, it suffices to prove $\valu_{\sigma}^{\M}(s_{i,j})>\valu_{\sigma}^{\M}(g_1)$ as this implies $\valu_{\sigma}^{\M}(F_{i,j})>\valu_{\sigma}^{\M}(g_1)$.
Since $\sigma(d_{i,j,k})\neq F_{i,j}$, either $\indbit_i=0$ or $\indbit_{i+1}\neq j$.
In the second case, \Pref{USV1}$_i$ implies $\sigma(s_{i,j})=b_1$ and the statement follows since $\valu_{\sigma}^{\M}(s_{i,j})=\rew{s_{i,j}}+\valu_{\sigma}^{\M}(g_1)$ due to $\sigma(b_1)=g_1$.
Thus let $\indbit_i=0\wedge\indbit_{i+1}=j$.
Then, the statement follows since $\valu_{\sigma}^{\M}(s_{i,j})=\rew{s_{i,j},h_{i,j}}+\valu_{\sigma}^{\M}(b_{i+1})$ by \Pref{EV1}$_{i+1}$ and $\rew{s_{i,j},h_{i,j}}>\sum_{\ell<i}\rew{g_\ell,s_{\ell,\sigmabar(g_{\ell})},h_{\ell,\sigmabar(g_{\ell})}}$.
\end{proof}

\ORBeginningPhaseFiveForNSBLargerOne*

\begin{proof}
We consider each cell of the table individually.
We also observe that  $\sigma(e_{i,j,k})=g_1$ implies $(e_{i,j,k},b_2)\in I_{\sigma}$, it holds that $\valu_{\sigma}^*(g_1)\prec\valu_{\sigma}^*(b_2)$ .
\begin{enumerate}
	\item Let $e=(d_{i,j,k},F_{i,j})$ with $\sigma(e_{i,j,k})=g_1$.
		Then, $\occrec^{\sigma}(e)=\occrec^{\canstrat}(e)=\maxocc$ by \Cref{lemma: Extended Reaching phase 5}.
	\item Let $e=(e_{i,j,k},b_2)$ with $\sigma(e_{i,j,k})=g_1$.
		Then, $\occrec^{\sigma}(e)=\occrec^{\canstrat}(e)=\maxocc-1$ by \Cref{table: Occurrence Records}.
	\item Let $e=(d_{\nsb,j,k},F_{\nsb,j})$ with $j\coloneqq 1-\indbit_{\nsb+1}$ for some $k\in\{0,1\}$.
		This edge is only an improving switch if $\bit+1$ is not a power of~two.
		Note that this in particular implies $\id_{j=0}\lastflip{\bit}{\nsb+1}{}+\id_{j=1}\lastunflip{\bit}{\nsb+1}{}\neq 0$.
		Since $\bit_{\nsb}=0\wedge\bit_{\nsb+1}\neq j$, \Cref{lemma: Numerics Of Ell} thus implies  \begin{align*}
			\ell^{\bit}(\nsb,j,k)&=\ceil{\frac{\bit-2^{\nsb-1}+\sum(\bit,\nsb)+1-k}{2}}=\ceil{\frac{\bit-2^{\nsb-1}+2^{\nsb-1}-1+1-k}{2}}\\
				&=\ceil{\frac{\bit-k}{2}}=\floor{\frac{\bit+1-k}{2}}=\floor{\frac{\bit+1}{2}}-k.
		\end{align*}
		Since $\bit+1$ is not a power of two, the parameter $t_{\bit}=-1$ is not feasible by \Pref{OR3}$_{i,j,k}$.
		Hence $\occrec^{\canstrat}(d_{\nsb,j,k},F_{\nsb,j})=\maxocc-k$.
		Then, \Cref{corollary: Switches With Low OR In Phase One} implies that $(d_{\nsb,j,1},F_{\nsb,j})$ was applied during $\canstrat\to\sigma$.
		Consequently, for both $k\in\{0,1\}$, it holds that $\occrec^{\sigma}(d_{\nsb,j,k},F_{\nsb,j})=\maxocc$.
	\item Let $e=(d_{i,j,k},F_{i,j})$ with $i\in\{\nsb+1,\dots,m-1\}, \indbit_i=0, j\coloneqq 1-\indbit_{i+1}$ and $k\in\{0,1\}$.
		This edge is only an improving switch if $\bit+1$ is not a power of two.
		Since $i>\nsb$, $\indbit_i=0$ implies $\bit_i=0\wedge\bit_{i+1}\neq j$.
		Also, $i<m$ implies $\id_{j=0}\lastflip{\bit}{i+1}{}+\id_{j=1}\lastunflip{\bit}{i+1}{}\neq 0$ since $j=1-\indbit_{i+1}$ and $\bit\geq 1$ by the choice of $i$.
		Since $\bit_{\nsb}=0$, this yields \begin{align*}
			\ell^{\bit}(i,j,k)&=\ceil{\frac{\bit-2^{i-1}+\sum(\bit,i)+1-k}{2}}\leq\ceil{\frac{\bit-2^{i-1}+2^{i-1}-1-2^{\nsb-1}+1-k}{2}}\\
				&=\ceil{\frac{\bit-2^{\nsb-1}-k}{2}}\leq\ceil{\frac{\bit-2-k}{2}}=\floor{\frac{\bit-1-k}{2}}\leq\floor{\frac{\bit-1}{2}}\leq\maxocc-1.
		\end{align*}
		There are two cases.
		If $\canstrat(d_{i,j,k})=F_{i,j}$, then $\occrec^{\canstrat}(d_{i,j,k},F_{i,j})=\ell^{\bit}(i,j,k)+1\leq\maxocc-1$ by \Pref{OR1}$_{i,j,k}$.
		If $\canstrat(d_{i,j,k})\neq F_{i,j}$, then $\occrec^{\canstrat}(d_{i,j,k},F_{i,j})=\ell^{\bit}(i,j,k)\leq\maxocc-1$.
		In the first case, $e$ was not applied during phase 1 and $\occrec^{\canstrat}(e)=\occrec^{\sigma}(e)\leq\maxocc-1$.
		In the second case, $\occrec^{\canstrat}(e)=\maxocc-1$ by \Pref{OR4}$_{i,j,k}$.
		Then, $e$ was applied during phase 1, implying $\occrec^{\sigma}(d_{i,j,k},F_{i,j})=\maxocc$.		
	\item Let $e=(d_{i,j,k},F_{i,j})$ with $i\leq\nsb-1$ and $j\coloneqq 1-\indbit_{i+1}$.
		Then, bit $i$ and bit $i+1$ switched during $\canstrat\to\sigma^{(5)}$.
		In particular, $F_{i,j}$ was closed with respect to $\canstrat$ and consequently $(d_{i,j,k},F_{i,j})\notin\applied{\canstrat}{\sigma}$.
		Hence, by \Cref{table: Occurrence Records},  \begin{align*}
			\occrec^{\sigma}(e)=\occrec^{\canstrat}(e)&=\ceil{\frac{\lastflip{\bit}{i}{\{(i+1,j)\}}+1-k}{2}}=\ceil{\frac{\bit-\sum(\bit,i)+1-k}{2}}\\
				&=\ceil{\frac{\bit-2^{i-1}+1+1-k}{2}}=\floor{\frac{\bit-2^{i-1}+3-k}{2}}.
		\end{align*}
		We now distinguish several cases.
		\begin{itemize}
			\item For $i=1$, $\occrec^{\sigma}(e)=\floor{(\bit+2-k)/2}=\maxocc$ independent of $k$.
			\item For $i=2$, $\occrec^{\sigma}(e)=\floor{(\bit+1-k)/2}$, so $\occrec^{\sigma}(e)=\maxocc-k$.
			\item For $i=3$, $\occrec^{\sigma}(e)=\floor{(\bit-1-k)/2}$, so $\occrec^{\sigma}(e)=\maxocc-1$ if $k=0$ and $\occrec^{\sigma}(e)=\maxocc-2$ if $k=1$.
			\item For $i>3$, it is easy to see that the occurrence record is always strictly smaller than $\maxocc-1$.
		\end{itemize}
\end{enumerate}
\end{proof}

\ApplicationOfPhaseThreeSwitches*

\begin{proof}
As a reminder, $e=(d_{i,j,k},F_{i,j})$ being a type $3$ switch implies that we either have $i<\nsb-1, j=1-\indbit_{i+1}$ or $i\in\{\nsb+1,\dots,m-1\}, \indbit_i=0, j\coloneqq 1-\indbit_{i+1}$.
In the second case, $\canstrat(d_{i,j,k})=F_{i,j}$ holds as well.
Since it is easy to verify that $i\neq 1$ and $\sigma(b_i)=b_{i+1}$ (for example by the arguments used in the proof of \Cref{claim: OR Beginning Phase 5 For NSB>1}), we only show $\sigmabar(g_i)=1-\indbit_{i+1}$.
By \Cref{lemma: Extended Reaching phase 5}, this holds for all $i\leq\nsb-1$.
It thus suffices to prove this for $i\in\{\nsb+1,\dots,m-1\}\wedge\indbit_{i}=0$.
We show the statement by proving that $\sigmabar(g_i)=\indbit_{i+1}$ implies $(g_i,F_{i,1-\indbit_{i+1}})\in I_{\sigma}$, contradicting the characterization of $I_{\sigma}$ given in Equation~(\ref{equation: IS Beginning Phase 5 For NSB>1}).

Since $j=1-\indbit_{i+1}$, it suffices to prove $\valu_{\sigma}^{*}(F_{i,j})\succ\valu_{\sigma}^*(F_{i,1-j})$. 
We have $(i,j)\in S_1$ and $(i,1-j)\in S_2$.
Thus, by \Cref{lemma: Extended Reaching phase 5}, $\sigmabar(eb_{i,j})\wedge\nsigmabar(eg_{i,j})$ as well as $\sigmabar(eb_{i,1-j})\wedge\sigmabar(eg_{i,1-j})$.
Also, by the choice of $j$ and \Pref{USV1}$_i$, $\sigma(s_{i,j})=b_1$.
Thus, by \Cref{lemma: Exact Behavior Of Counterstrategy,lemma: Exact Behavior Of Random Vertex}, $\valustar^*_{\sigma}(F_{i,j})=\valustar^*_{\sigma}(b_2)$ regardless of whether $G_n=S_n$ or $G_n=M_n$.
Also, since $\nsb\geq2$, $\sigmabar(g_1)=1-\indbit_2\neq\sigmabar(b_2)$ by \Cref{lemma: Extended Reaching phase 5}.
Thus, if $G_n=S_n$, then \Cref{lemma: Exact Behavior Of Counterstrategy} implies $\valustar^{\P}_{\sigma}(F_{i,1-j})=\valustar^{\P}_{\sigma}(g_1)\lhd\valustar_{\sigma}^{\P}(b_2)=\valustar_{\sigma}^{\P}(F_{i,j})$ as player~1 minimizes the valuation.
If $G_n=M_n$, then $\valustar_{\sigma}^\M(F_{i,1-j})=\frac{1}{2}\valustar_{\sigma}^\M(g_1)+\frac{1}{2}\valustar_{\sigma}^{\M}(b_2)$, hence the statement follows since $\valustar_{\sigma}^\M(g_1)<\valustar_{\sigma}^\M(b_2)$.
\end{proof}

\ApplicationOfTypeTwoSwitchesEscape*

\begin{proof}
Let $i\in[n], j,k\in\{0,1\}$ and let $e\coloneqq (e_{i,j,k},b_2)$ be an improving switch.
We begin by proving that the cycle center $F_{i,j}$ is mixed.
Since only improving switches of type~$3$ were applied so far during phase $5$, $\sigma(e_{i,j,k})=g_1$ implies $\sigma(d_{i,j,k})=e_{i,j,k}$.
Consequently, we have $\sigmabar(eg_{i,j})$.
In particular, $F_{i,j}$ is not closed, so $\indbit_i=0\vee\indbit_{i+1}=j$.
Thus, either $(i,j)\in S_1$ or $(i,j)\in S_2$.
By \Cref{lemma: Extended Reaching phase 5}, $\sigmabar(eg_{i,j})$ and as no switch $(e_{*,*,*}, g_1)$ was applied during $\sigma^{(5)}\to\sigma$, we need to have $(i,j)\in S_2$, implying that $F_{i,j}$ is mixed.

We next prove that $j=1$ resp. $j=1-\indbit_{i+1}$ (depending on whether $G_n=S_n$ or $G_n=M_n$) implies $\nsigmabar(eg_{i,1-j})$.
Consider the case $G_n=S_n$ and thus $j=1$ first.
We prove $\nsigmabar(eg_{i,0})$.
If $F_{i,0}$ is closed, then the statement follows.
If it is not closed, then $\indbit_i=0\vee\indbit_{i+1}\neq 0$. 
Consequently, either $(i,0)\in S_1$ or $(i,0)\in S_2$.
In the first case, $\nsigmabar(eg_{i,0})$ follows from \Cref{lemma: Extended Reaching phase 5} as no improving switch $(e_{*,*,*}, b_2)$ was applied during $\sigma^{(5)}\to\sigma$, so assume $(i,0)\in S_2$.
Then, by the same lemma, both cycle centers $F_{i,0}, F_{i,1}$ were mixed for $\sigma^{(5)}$.
Thus, as we consider the case $G_n=S_n$, the tie-breaking rule must have applied the improving switches $(e_{i,0,*},b_2)$ prior to $(e_{i,j,k},b_2)$, implying $\nsigmabar(eg_{i,0})$.		
If $G_n=M_n$, then $\nsigmabar(eg_{i,1-j})$ follows by the same arguments as the tie-breaking rule applied the improving switches $(e_{i,\indbit_{i+1},*},b_2)$ first.
Finally, as no improving switch $(g_*,F_{*,*})$ was applied during $\sigma^{(5)}\to\sigma$, $\nsb=2$ implies $\sigma(g_1)=F_{1,0}$ if $G_n=S_n$ by \Cref{lemma: Extended Reaching phase 5}.
Thus, all requirements of \Cref{lemma: Phase Five Escape Easy} are met.
\end{proof}

\FulfillingORFiveIfNSBIsNotOne*

\begin{proof}
To simplify the notation, let $\sigma$ denote the strategy obtained after the application of $(e_{i,j,k},b_2)$ and let $e\coloneqq(d_{i,j,1-k},F_{i,j})\in I_{\sigma}$.
We prove that $\sigma$ has \Pref{OR4}$_{i,j,1-k}$ by proving\begin{equation} \label{equation: Fulfilling OR Five}
\occrec^{\sigma}(d_{i,j,1-k},F_{i,j})\in\left \{\floor{\frac{\bit+1+1}{2}}-1,\floor{\frac{\bit+1+1}{2}}\right \}.
\end{equation}
The second statement is shown along the way.

Consider the case $e\in\applied{\canstrat}{\sigma}$.
Since $e\in\applied{\sigma^{(5)}}{\sigma}$ would imply $\sigma(d_{i,j,1-k})=F_{i,j}$, we need to have $e\in\applied{\canstrat}{\sigma^{(5)}}$.
This implies that the switch was applied during phase 1 as well as $\canstrat(d_{i,j,1-k})\neq F_{i,j}$ and $\occrec^{\canstrat}(e)\in\{\maxocc-1,\maxocc\}$.
The only improving switches of type $(d_{*,*,*},F_{*,*})$ with an occurrence record of $\maxocc$ applied in phase 1 are the cycle edges of $F_{\nsb,\indbit_{\nsb+1}}$.
Consequently, $\occrec^{\canstrat}(e)=\maxocc-1$ as $F_{\nsb,\indbit_{\nsb+1}}$ is closed and its cycle edges cannot become improving switches.
Hence \[\occrec^{\sigma}(e)=\occrec^{\canstrat}(e)+1\in\left \{\floor{\frac{\bit+1+1}{2}}-1, \floor{\frac{\bit+1+1}{2}}\right \},\] proving both parts of the statement.

Consider the case  $e\notin\applied{\canstrat}{\sigma}$ next.
Since the switch was not applied, this then  implies $\occrec^{\canstrat}(e)=\occrec^{\sigma}(e)$.
We distinguish two cases.
\begin{enumerate}
	\item Consider the case $\canstrat(d_{i,j,1-k})=F_{i,j}$ first.
		Then, \Pref{OR1}$_{i,j,1-k}$ implies $\occrec^{\canstrat}(e)=\ell^{\bit}(i,j,1-k)+1\leq\maxocc-1$.
		Assume $\occrec^{\canstrat}(d_{i,j,1-k},F_{i,j})<\maxocc-1$.
		This implies $\ell^{\bit}(i,j,1-k)\leq\maxocc-3$ by integrality.
		Since $\ell^{\bit}(i,j,0)$ and $\ell^{\bit}(i,j,1)$ differ by at most $1$, it follows that $\ell^{\bit}(i,j,k)\leq\maxocc-2$.
		This implies that $\occrec^{\canstrat}(d_{i,j,k},F_{i,j})\leq\maxocc-1$.
		But this is a contradiction since $(e_{i,j,k},b_2)\in I_{\sigma^{(5)}}$ implies $\occrec^{\canstrat}(d_{i,j,k},F_{i,j})=\maxocc$ by \Cref{lemma: Extended Reaching phase 5}.
		Hence the statement follows from $\occrec^{\canstrat}(e)=\occrec^{\sigma}(e)=\maxocc-1$.
	\item Let $\canstrat(d_{i,j,1-k})\neq F_{i,j}$.
		Then, $\occrec^{\sigma}(e)=\occrec^{\canstrat}(e)=\maxocc$ by \Pref{OR4}$_{i,j,1-k}$ and \Cref{corollary: Switches With Low OR In Phase One}.
		This already implies the first of the two statements. 
		In addition, either $\canstrat(d_{i,j,k})=F_{i,j}$ and $\occrec^{\canstrat}(d_{i,j,k},F_{i,j})\leq\maxocc-1$ or $\canstrat(d_{i,j,k})\neq F_{i,j}$ and $\occrec^{\canstrat}(d_{i,j,k},F_{i,j})=\maxocc-1$.
		If none of these were true, then $F_{i,j}$ would be open at the end of phase~1, contradicting \Cref{corollary: No Open CC For Phase Two}.
		This proves the second part of the statement. \qedhere
\end{enumerate}
\end{proof}

\SelectionVerticesPhaseFiveNSBLargerOne*

\begin{proof}
Let $e=(g_i,F_{i,j})$.
By \Cref{lemma: Phase Five Escape Easy}, $e\in I_{\sigma}$ if and only if $\indbit_i=0, \sigmabar(eb_{i,1-j})$ and $[j=0\wedge\sigmabar(g_i)=1]$ if $G_n=S_n$ resp. $[j=\indbit_{i+1}\wedge\sigmabar(g_i)=1-\indbit_{i+1}]$ if $G_n=M_n$.
Let, for the sake of contradiction, $(g_i,F_{i,j})\in\applied{\canstrat}{\sigma}$.
The conditions on $j$ and $\sigmabar(g_i)$ imply $(g_i,F_{i,j})\notin\applied{\sigma^{(5)}}{\sigma}$.
Since $\indbit_i=0$ implies $i\neq\nsb$, also $(g_i,F_{i,j})\neq(g_\nsb,F_{\nsb,*})$.
Thus, by \Cref{lemma: Extended Reaching phase 5}, $\bit_i=0\wedge\bit_{i+1}\neq j$.
Consequently, $0=\bit_i=\indbit_{i+1}=(\bit+1)_{i+1}$ and $j=1-\bit_{i+1}$.
Since all bits below level $\nsb$ have $\bit_i=1\wedge(\bit+1)_i=0$, this implies $i>\nsb$.
Therefore, $\bit_{i+1}=(\bit+1)_{i+1}=1-j$ and in particular $j=1-\indbit_{i+1}$
This is a contradiction if $G_n=M_n$ as $j=\indbit_{i+1}$.
Hence consider the case $G_n=S_n$.
Then, $j=1-\indbit_{i+1}=0$, implying $\indbit_{i+1}=1$.
Thus, $i\in\{\nsb+1,\dots,m-1\}, \indbit_{i}=0$ and $j=1-\indbit_{i+1}$, implying $(i,j)\in S_1$.
Therefore, $\sigmabar^5(eb_{i,j})\wedge\nsigmabar^5(eg_{i,j})$, contradicting $(e_{i,j,k},b_2), (d_{i,j,k},e_{i,j,k})\in I_{\sigma^{(5)}}$.
Thus, $(g_i,F_{i,j})\notin\applied{\canstrat}{\sigma^{(5)}}$, implying the statement.
\end{proof}

\PropertySVAtEndOfPhaseFive*

\begin{proof}
As a reminder, $\sigma$ is the strategy obtained after the application of the switch $(e_{1,\indbit_2,k},b_2)$ resp. after the application of the switch $(g_1,F_{1,\indbit_2})$ if it becomes improving.
First consider some $i\geq 2$.
If $\indbit_i=1$, then $\sigma$ has Property (SV*)$_i$ as it is a phase-5-strategy.
If $\indbit_i=0$ and $(g_i,F_{i,j})\in\applied{\sigma^{(5)}}{\sigma}$, then this follows from \Cref{corollary: Selection Vertices Phase 5}. 
Thus, let $\indbit_i=0$ and $(g_i,F_{i,j})\notin\applied{\sigma^{(5)}}{\sigma}$.
For the sake of a contradiction, assume that $\sigma$ does not have Property (SV*)$_1$.
Then, $\sigmabar(g_i)=1$ resp. $\sigmabar(g_i)=1-\indbit$ (depending on whether $G_n=S_n$ or $G_n=M_n$) and $\nsigmabar(d_{i,1})$ resp. $\nsigmabar(d_{i,1-\indbit_{i+1}})$.
To simplify notation, let $j\coloneqq 0$ resp. $j\coloneqq\indbit_{i+1}$.
We show that we then have $\valu_{\sigma}^*(F_{i,j})\succ\valu_{\sigma}^*(F_{i,1-j})$, implying $(g_i,F_{i,j})\in I_{\sigma}$.
But this is a contradiction as any improving switch of this kind is applied immediately and $i\geq 2$ implies that the application of $(e_{1,\indbit^{\sigma}_2,k},b_2)$ cannot have unlocked this switch.

As the last improving switch of the type $(e_{*,*,*},b_2)$ was just applied, any cycle center is either closed or escapes to $b_2$.
We first consider the case $G_n=S_n$.
Since $\sigma$ is a phase-5-strategy for $\bit$, it has \Pref{USV1}$_i$ and \Pref{EV1}$_{i+1}$.
Consequently, either $\sigmabar(b_{i+1})=j$ or $\nsigmabar(s_{i,j})$
If both cycle centers of level $i$ escape towards $b_2$, then the statement follows since \[\valu_{\sigma}^P(F_{i,j})=\{F_{i,0},e_{i,0,*},d_{i,0,*}\}\cup\valu_{\sigma}^P(b_2)\rhd\{F_{i,1},e_{i,1,*},d_{i,1,*}\}\cup\valu_{\sigma}^P(b_2)=\valu_{\sigma}^P(F_{i,1-j})\] by \Cref{lemma: Exact Behavior Of Counterstrategy}.
Since $\indbit_i=0$, only $F_{i,1-\indbit_{i+1}}$ can be closed in level~$i$.
Let this cycle center be closed.
If $j=1-\indbit_{i+1}=0$, then \Pref{USV1}$_i$ and $\sigma(b_1)=b_2$ implies $\valustar_{\sigma}^P(F_{i,0})=\{s_{i,0}\}\cup\valustar_{\sigma}^P(b_2)$ and the statement follows from $\valustar_{\sigma}^{\P}(F_{i,1})=\valustar_{\sigma}^{\P}(b_2)$.
If $j=\indbit_{i+1}=0$, then  $F_{i,1-j}=F_{i,1}$ is closed,  contradicting the assumption $\nsigmabar(d_{i,1})$.

Consider the case $G_n=M_n$.
If both cycle centers are $b_2$-open or $b_2$-halfopen, then the statement follows by \Cref{lemma: Both CC Open For MDP} since $\sigma$ has \Pref{REL1}.
If $F_{i,j}$ is $b_2$-open and $F_{i,1-j}$ is $b_2$-halfopen, then the statement follows by an easy but tedious calculation.
Thus consider the case that $F_{i,j}$ is $b_2$-halfopen and that $F_{i,1-j}$ is $b_2$-open.
Then, by the choice of $j$ and \Pref{USV1}$_i$, \begin{align*}
	\valu_{\sigma}^\M(F_{i,j})&=\frac{1-\e}{1+\e}\valu_{\sigma}^\M(b_2)+\frac{2\e}{1+\e}\valu_{\sigma}^\M(s_{i,j})\\
		&=\frac{1-\e}{1+\e}\valu_{\sigma}^\M(b_2)+\frac{2\e}{1+\e}[\rew{s_{i,j},h_{i,j}}+\valu_{\sigma}^\M(b_{i+1})]\\
	\valu_{\sigma}^\M(F_{i,1-j})&=(1-\e)\valu_{\sigma}^\M(b_2)+\e\valu_{\sigma}^\M(s_{i,1-j})=\valu_{\sigma}^\M(b_2)+\e\rew{s_{i,1-j}}, \\
	\valu_{\sigma}^\M(F_{i,j})-\valu_{\sigma}^\M(F_{i,1-j})&=\frac{2\e}{1-\e}(\rew{s_{i,j},h_{i,j}}+\valu_{\sigma}^\M(b_{i+1}))-\frac{2\e}{1+\e}\valu_{\sigma}^\M(b_2)-\e\rew{s_{i,j}}\\
		&=\e\left[\frac{2}{1+\e}(\rew{s_{i,j},h_{i,j}}+\valu_{\sigma}^\M(b_{i+1})-\valu_{\sigma}^\M(b_2))-\rew{s_{i,1-j}}\right]
\end{align*}
It thus suffices to show that the last term is larger than zero which follows easily from $\indbit_i=0$.

In level $i$, only $F_{i,1-\indbit_{i+1}}=F_{i,1-j}$ can be closed.
Then, the statement follows by the same argument used for the case $G_n=S_n$.

We now consider Property (SV*)$_1$.
Assume that $(g_1,F_{1,\indbit_2})$ does not become improving when applying $(e_{1,\indbit_2,k},b_2)$.
Then, by \Cref{lemma: Phase Five Escape Easy}, we need to have $\sigmabar(g_i)=\indbit_{i+1}$ if $G_n=M_n$.
Consider the case $G_n=S_n$.
If $\indbit_{2}=0$, then \Cref{lemma: Phase Five Escape Easy} implies that we need to have $\sigmabar(g_1)=0$.
If $\indbit_2=1$, then $\nsb=2$.
But this implies $\canstrat(g_1)=F_{1,0}$ since the cycle center $F_{1,0}$ was then closed with respect to $\canstrat$.
For this reason, the switch $(g_1,F_{1,1})$ was not applied during phase 1.
Since a switch involving a selection vertex $gi$ can only be applied during phase 5 if $\sigmabar(g_i)=1$ by \Cref{lemma: Phase Five Escape Easy}, the switch cannot have been applied during phase 5.
Consequently, $\sigma(g_1)=\canstrat(g_1)=F_{1,0}$
Thus, $\sigma$ has Property (SV*)$_1$.
If the edge $(g_1,F_{1,\indbit_2})$ becomes an improving switch, then the strategy obtained after applying it has Property (SV*)$_i$ by \Cref{corollary: Selection Vertices Phase 5}.
\end{proof}

\ImprovingSwitchesOfNewCanstrat*

\begin{proof}
Let $\sigma^{(5)}$ denote the phase-5-strategy of \Cref{lemma: Extended Reaching phase 5} with $\sigma\in\reach{\sigma^{(5)}}$.
We first observe that $\indbit^{\sigma}=\indbit^{\sigma^{(5)}}$, so the upper index can be omitted.
It is easy to verify that $I_{\sigma}$ can be partitioned as \begin{align*}
I_{\sigma}&=\{(d_{i,j,*},F_{i,j})\colon\sigma^{(5)}(e_{i,j,*})=g_1\}\cup\{(d_{\nsb,1-\indbit_{\nsb+1},*},F_{\nsb,1-\indbit_{\nsb+1}})\}\\
	&\cup\left \{e=(d_{i,1-\indbit_{i+1},*},F_{i,1-\indbit_{i+1}})\colon i\in\{\nsb+1,\dots,m-1\}, \indbit_i=0, \occrec^{\sigma_5}(e)=\maxocc-1\right \}\\
	&\cup\left \{e=(d_{i,1-\indbit_{i+1},*},F_{i,1-\indbit_{i+1}})\colon i<\nsb, \occrec^{\sigma^{(5)}}(e)\geq\maxocc-1\right \},
\end{align*} if $\bit+1$ is not a power of two.
A similar partition can be derived if $\bit+1$ is a power of two.
In particular, $I_{\sigma}\subseteq\{(d_{i,j,k},F_{i,j})\colon\sigma(d_{i,j,k})\neq F_{i,j}\}$.
We prove that $e=(d_{i,j,k},F_{i,j})$ implies $e\in I_{\sigma}$ if $\sigma(d_{i,j,k})\neq F_{i,j}$.
		
If $\sigma^{(5)}(e_{i,j,k'})=g_1$ for some $k'\in\{0,1\}$, then $e\in I_{\sigma}$ as one of the cycle edges of $F_{i,j}$ is improving for $\sigma^{(5)}$ while the other becomes improving after applying $(e_{i,j,k'},b_2)$.
Thus let $\sigma^{(5)}(e_{i,j,*})=b_2$, implying $\nsigmabar^{(5)}(eg_{i,j})$.
Then, by \Cref{lemma: Extended Reaching phase 5}, either $\sigmabar^{(5)}(d_{i,j})$ or $\sigmabar^{(5)}(eb_{i,j})\wedge\nsigmabar^{(5)}(eg_{i,j})$.
In the first case, $\indbit_i=1\wedge\indbit_{i+1}=j$ by \Cref{lemma: Extended Reaching phase 5}.
But this implies $\sigmabar(d_{i,j})$ since $\sigma$ is a phase-5-strategy for $\bit$ and thus has \Pref{EV1}$_i$.
This however contradicts $\sigma(d_{i,j,k})\neq F_{i,j}$.
Hence, assume that $\sigmabar^{(5)}(eb_{i,j})\wedge\nsigmabar^{(5)}(eg_{i,j})$.
Then, by \Cref{lemma: Extended Reaching phase 5}, $(i,j)\in S_1$.
We distinguish three cases.
\begin{enumerate}
	\item Let $(i,j)\in\{(i,1-\indbit_{i+1})\colon i\leq\nsb-1\}$.
		If $\occrec^{\sigma^{(5)}}(e)<\maxocc-1$, then $e$ was an improving switch of type 3 for $\sigma^{(5)}$ and thus applied during phase~5.
		But this contradicts $\sigma(d_{i,j,k})\neq F_{i,j}$ since no switch $(d_{*,*,*},e_{*,*,*})$ is applied during phase~5.
		This implies $(i,j)\in\{(i,1-\indbit_{i+1})\colon i\leq\nsb-1, \occrec^{\sigma^{(5)}}(e)\geq\maxocc-1\}$, hence $e\in I_{\sigma}$.
	\item Let $(i,j)\in \{(i,1-\indbit_{i+1})\colon i\in\{\nsb+1,\dots,m-1\}, \indbit_{i}=0\}$ which can only occur if $\bit+1$ is not a power of 2.
		As proved when discussing $I_{\sigma^{(5)}}$, we then either have $\canstrat(d_{i,j,k})=F_{i,j}$, implying $\occrec^{\canstrat}(d_{i,j,k},F_{i,j})\leq\maxocc-1$ or $\canstrat(d_{i,j,k})\neq F_{i,j}$ and $\occrec^{\canstrat}(d_{i,j,k},F_{i,j})=\maxocc-1$.
		Consider the first case.
		If the inequality is strict, the switch was applied previously during phase 5, yielding a contradiction.
		Otherwise, $(d_{i,j,k},F_{i,j})\in I_{\sigma}$.
		In the second case, the switch was applied during phase~1, hence it was a switch of type 1 during phase~5, also implying $(d_{i,j,k},F_{i,j})\in I_{\sigma}$.
	\item Finally, let $i=\nsb\wedge j=1-\indbit_{\nsb+1}$ which only needs to be considered if $\bit+1$ is not a power of 2.
		In this case we however have $e\in I_{\sigma^{(5)}}$, implying $e\in I_{\sigma}$. 
\end{enumerate}
Thus, $e\in I_{\sigma}$ in all case, proving the statement.
\end{proof}

\NewCanstratIsCanstrat*

\begin{proof}
Consider \Cref{definition: Canonical Strategy,definition: Canonical Strategy MDP}. 
As $\sigma$ is a phase-$5$-strategy for $\bit$, it holds that $\indbit=\bit+1$.
Thus, condition~1 follows since $\sigma(e_{*,*,*})=b_2$ and $\nsb>1$.
This also implies that conditions 2(a), 2(c), 3(a) and 3(b) are fulfilled as $\sigma$ has \Pref{EV1}$_*$ and \Pref{EV2}$_*$.

Consider condition 2(b) and let $i\in[n]$.
Since $(\bit+1)_i=1$ implies that $F_{i,(\bit+1)_{i+1}}$ is closed, we prove that $F_{i,1-(\bit+1)_{i+1}}$ is not closed.
Let $j\coloneqq 1-(\bit+1)_{i+1}$.
Then, by \Cref{lemma: Extended Reaching phase 5}, $\sigma^{(5)}(d_{i,j,*})=e_{i,j,*}$ and it suffices to prove $(d_{i,j,0},F_{i,j})\notin\applied{\sigma^{(5)}}{\sigma}$.
As such a switch is applied during $\sigma^{(5)}\to\sigma$ if and only if it is of type 3 by \Cref{corollary: Improving Switches Of Type 3}, we prove 
\begin{equation} \label{equation: High OR at end of phase 5 nsb>1}
\occrec^{\sigma^{(5)}}(d_{i,j,0},F_{i,j})\geq\floor{\frac{\bit+1}{2}}-1.
\end{equation}
This follows if $\id_{j=0}\lastflip{\bit}{i+1}{}+\id_{j=1}\lastunflip{\bit}{i+1}{}=0$ since this implies $\ell^{\bit}(i,j,k)\geq \bit$.
Thus suppose that this term is not 0.
Then, since $\bit_1=1$ and by the choice of $i$ and $j$, \[\ell^{\bit}(i,j,k)=\ceil{\frac{\bit+\sum(\bit,i)+1-k}{2}}\geq\ceil{\frac{\bit+2-k}{2}}=\floor{\frac{\bit+1-k}{2}}+1.\] 
But this implies Inequality (\ref{equation: High OR at end of phase 5 nsb>1}) since $\ell^{\bit}(i,j,0)\geq\floor{(\bit+1)/2}+1$.

Consider condition 3(c) and let $i\in[n]$ and  $j\coloneqq 1-(\bit+1)_{i+1}$.
It is easy to prove that $\sigma$ has condition~3(c) since $(\bit+1)_i=0$ and $F_{i,j}$ being closed imply $\valu_{\sigma}^*(F_{i,j})\succ\valu_{\sigma}^*(F_{i,1-j})$.
Since $F_{i,j}$ is closed,  $\valustar_{\sigma}^*(F_{i,j})=\ubracket{s_{i,j}}\oplus\valustar_{\sigma}^*(b_2)$ by \Pref{USV1}$_i$, \Cref{lemma: Exact Behavior Of Counterstrategy} and $\sigma(b_1)=b_2$.
Since $F_{i,1-j}$ cannot be closed due to the choice of $j$ and $(\bit+1)_i=0$, we have $\sigmabar(eb_{i,1-j})\wedge\nsigmabar(eg_{i,1-j})$.
Consequently, $\valustar_{\sigma}^*(F_{i,1-j})=\valustar_{\sigma}^*(b_2)$ since $\sigma(b_1)=b_2$.
But this implies $\valustar_{\sigma}^*(F_{i,1-j})\prec\valustar_{\sigma}^*(F_{i,j})$.
Next, consider condition 3(d) and consider a level $i$ with $(\bit+1)_i=0$.
Let $j\coloneqq 0$ resp. $j\coloneqq\indbit_{i+1}$ depending on whether $G_n=S_n$ or $G_n=M_n$.
We prove that $\valu_{\sigma}^*(F_{i,j})\succ\valu_{\sigma}^*(F_{i,1-j})$ if none of the two cycle centers is closed.
If $G_n=M_n$, this either follows from \Cref{lemma: Both CC Open For MDP} since $\sigma$ has \Pref{REL1} or by an easy but tedious calculation.
If $G_n=S_n$, this follows since $\Omega(F_{i,0})>\Omega(F_{i,1})$ and as these priorities are even.

\Pref{USV1} implies that $\sigma$ fulfills conditions 4 and 5 for all indices.
Finally, consider condition~6 and let $i=\ell(\bit+2), j=\indbit_{\ell(\bit+2)+1}$.
By the same argument used for condition~3(c), it suffices to prove $\occrec^{\sigma}(d_{i,j,k},F_{i,j})\geq\floor{(\bit+1)/2}-1$ for both $k\in\{0,1\}$.
This however follows from $\ell(\bit+1)=1, \indbit_2=1-\bit_2$ and $\bit_1=1$ by \[\ell^{\bit}(1,1-\bit_2,k)=\ceil{\frac{\bit+\sum(\bit,1)+1-k}{2}}=\ceil{\frac{\bit+1-k}{2}}=\floor{\frac{\bit+2-k}{2}}=\floor{\frac{\bit+1}{2}},\]implying the statement.
Hence, $\sigma$ is a canonical strategy for $\bit+1$.
\end{proof}

\OccRecForBPlusTwo*

\begin{proof}
Assume $\bit+2=2^l$ for some $l\in\mathbb{N}$.
Then, the choice of $i$ and $j$ implies $\bit+2=2^{i-1}$ and $j=1$.
In particular $\lastunflip{\bit}{i+1}{}=0$, implying  $\ell^{\bit}(i,j,k)\geq\bit$ by \Cref{lemma: Numerics Of Ell}.
Consequently, $\occrec^{\canstrat}(d_{i,j,k},F_{i,j})=\floor{(\bit+1-k)/2}=\maxocc$ since $\bit+1$ is odd.
In addition, $\occrec^{\canstrat}(d_{i,j,k},F_{i,j})\neq\ell^{\bit}(i,j,k)+1$, hence $\canstrat(d_{i,j,k})\neq F_{i,j}$ as $\canstrat$ has the canonical properties.

Thus assume that $\bit+2$ is not a power of 2.
Since $\bit$ is even and by the choice of $i$, it holds that $\bit_1=0$ and $\bit_2=\dots=\bit_{i-1}=1$.
In particular $\id_{j=0}\lastflip{\bit}{i+1}{}+\id_{j=1}\lastunflip{\bit}{i+1}{}\neq 0$.
Hence, by \Cref{lemma: Numerics Of Ell}, \begin{align*}
\ell^{\bit}(i,j,k)&=\ceil{\frac{\bit-2^{i-1}+\sum(\bit,i)+1-k}{2}}=\ceil{\frac{\bit-2^{i-1}+2^{i-1}-2+1-k}{2}}\\
	&=\ceil{\frac{\bit-1-k}{2}}=\floor{\frac{\bit-k}{2}}.
\end{align*}
Furthermore, $\bit$ being even and \Pref{OR4}$_{i,j,0}$ implies $\occrec^{\canstrat}(d_{i,j,k},F_{i,j})\neq\ell^{\bit}(i,j,k)-1$.
Hence $\occrec^{\canstrat}(d_{i,j,0},F_{i,j})=\floor{(\bit+1)/2}$ and $\ell^{\bit}(i,j,1)=\floor{(\bit-1)/2}=\maxocc-1$.
Also, $\canstrat(d_{i,j,k})\neq F_{i,j}$ for both $k\in\{0,1\}$ since $\occrec^{\canstrat}(d_{i,j,k},F_{i,j})=\ell^{\bit}(i,j,k)+1<\maxocc$ otherwise, contradicting the previous arguments.
\end{proof}

\ContrapositionORThree*

\begin{proof}
As a remainder, we currently consider a canonical strategy $\sigma$ for $\bit+1$ with $I_{\sigma}=\mathfrak{D}^{\sigma}$.
We prove the statements one after another.
\begin{enumerate}
	\item We distinguish several cases.
		\begin{enumerate}
			\item Let $\bit_i=1\wedge\bit_{i+1}=j$.
				This implies $i\neq 1$ since $i=1$ contradicts the choice of $j$.
				Also $\bit\geq 4$ for the same reason.
				Let \[\id_{j=0}\lastflip{\bit+1}{i+1}{}+\id_{j=1}\lastunflip{\bit+1}{i+1}{}=0.\]
				Then $\ell^{\bit+1}(i,j,k)-1\geq\bit$ by \Cref{lemma: Numerics Of Ell}.
				Since $\occrec^{\sigma}(e)\leq\occrec^{\canstrat}(e)+1$, we then have $\occrec^{\sigma}(e)\leq\maxocc+1$.
				Since $\bit\geq 4$, this implies $\occrec^{\sigma}(e)<\ell^{\bit+1}(i,j,k)-1$, implying the statement.
				Let $\id_{j=0}\lastflip{\bit+1}{i+1}{}+\id_{j=1}\lastunflip{\bit+1}{i+1}{}\neq 0$ and observe \[\occrec^{\canstrat}(e)=\ceil{\frac{\lastflip{\bit}{i}{\{(i+1,j)\}}+1-k}{2}}=\ceil{\frac{\bit-\sum(\bit,i)+1-k}{2}}.\]
				We distinguish two more cases.
				\begin{enumerate}
					\item Let $(\bit+1)_i=1$, implying $(\bit+1)_{i+1}=j$.
						Then $e\notin\applied{\canstrat}{\sigma}$, and consequently $\occrec^{\canstrat}(e)=\occrec^{\sigma}(e)$.
						It is easy to verify that \[\id_{j=0}\lastflip{\bit+1}{i+1}{}+\id_{j=1}\lastunflip{\bit+1}{i+1}{}=\bit+1-\sum(\bit+1,i)-2^{i-1}-2^{i-1}\] in this case.
						Hence, by the definition of $\ell^{\bit+1}(i,j,k)$, we have \begin{align*}
							\ell^{\bit+1}(i,j,k)&=\ceil{\frac{\lastflip{\bit+1}{i}{\{(i+1,j)\}}+1-k}{2}}-\sum(\bit+1,i)+2^i\\
								&=\occrec^{\sigma_{\bit+1}}(e)+\sum(\bit+1,i)+2^i\geq\occrec^{\sigma_{\bit+1}}(e)+4,
						\end{align*} so $\occrec^{\sigma_{\bit+1}}(e)\leq\ell^{\bit+1}(i,j,k)-4$, implying \Cref{equation: Occrec not ell-1}.
					\item Assume $(\bit+1)_i=0$, implying $i<\nsb$ and thus $(\bit+1)_{i+1}\neq\bit_i=j$.
						Then $\bit_1=\dots=\bit_i=1$ and $(\bit+1)_1=\dots=(\bit+1)_{i}=0$.
						Hence, by \Cref{lemma: Numerics Of OR}, \begin{align*}
							\ell^{\bit+1}(i,j,k)&=\ceil{\frac{\bit+1-2^{i-1}+\sum(\bit+1,i)+1-k}{2}}\\
								&=\ceil{\frac{\bit-(2^{i-1}-1)+1-k}{2}}\\
								&=\ceil{\frac{\bit-\sum(\bit,i)+1-k}{2}}=\occrec^{\canstrat}(e).
						\end{align*}
						This implies $\occrec^{\sigma}(e)\geq\occrec^{\canstrat}(e)=\ell^{\bit+1}(i,j,k)$, and thus \Cref{equation: Occrec not ell-1}.
				\end{enumerate}
			\item Let $\bit_i=1\wedge\bit_{i+1}\neq j$.
				Since $j=(\bit+2)_{i+1}$, this implies $(\bit+2)_{i+1}\neq\bit_{i+1}$.
				Hence, bit $i+1$ was switched when transitioning from $\canstrat$ to $\sigma_{\bit+2}$.
				In one of the two transitions, the first bit switched from 0 to 1 and this bit was the only bit that was switched in this transition.
				Thus, either $[i<\ell(\bit+1)$ and $\ell(\bit+2)=1]$ or $[\ell(\bit+1)=1$ and $i<\ell(\bit+2)]$.
				Consider $[i<\ell(\bit+1)$ and $\ell(\bit+2)=1]$ first.
				Since $\bit_i=1$ and $\bit_{i+1}\neq j$, \Cref{lemma: Numerics Of Ell} implies that it holds that $\ell^{\bit}(i,j,k)=\ceil{(\bit+\sum(\bit,i)+1-k)/2}.$
				Now, since $i<\ell(\bit+1)$, we have $\bit_l=1$ for all $l<i$, implying $\ell^{\bit}(i,j,k)=\ceil{(\bit+2^{i-1}-k)/2}.$
				If $i=1$, then \[\ell^{\bit}(i,j,k)=\ceil{\frac{\bit+1-k}{2}}\geq\floor{\frac{\bit+1-k}{2}}.\]
				This implies \Cref{equation: Occrec is high} since the only feasible tolerance for $i=1$ is 0.
				If $i>1$, then \Cref{equation: Occrec is high} follows from\[\ell^{\bit}(i,j,k)\geq\ceil{\frac{\bit+2-k}{2}}=\floor{\frac{\bit+3}{2}}=\floor{\frac{\bit+1-k}{2}}+1.\]
				
				Thus, $\occrec^{\canstrat}(e)=\floor{(\bit+1-k)/2}\neq\ell^{\bit}(i,j,k)+1$.
				Consequently, by \Pref{OR1}$_{i,j,k}$, $\canstrat(d_{i,j,k})\neq F_{i,j}$ for both $k\in\{0,1\}$.
				Since $\ell(\bit+1)>i\geq 1$ implies that $\bit$ is odd, this yields $\occrec^{\canstrat}(d_{i,j,1},F_{i,j})<\occrec^{\canstrat}(d_{i,j,0},F_{i,j})$.
				Combining these implies that $(d_{i,j,1},F_{i,j})$ is applied during phase 1 of $\canstrat\to\sigma_{\bit+1}$, so \begin{align*}
				\occrec^{\sigma}(d_{i,j,1},F_{i,j})&=\floor{\frac{\bit}{2}}+1=\floor{\frac{\bit+1}{2}}=\floor{\frac{\bit+1+1-1}{2}}\quad\text{ and}\\
				\occrec^{\sigma}(d_{i,j,0},F_{i,j})&=\floor{\frac{\bit+1}{2}}=\floor{\frac{\bit+1+1-0}{2}}.
				\end{align*}				
				
				Next let $\ell(\bit+1)=1$ and $i<\ell(\bit+2).$
				Since $\ell(\bit+1)=1$ implies $\bit_1=0$, $\bit_i=1$ implies $i>1$.
				In addition, $i<\ell(\bit+2)$ implies $\bit_{i'}=1$ for all $i'\in\{2,\dots,i-1\}$.
				Consequently, as in the last case, \begin{align*}
					\ell^{\bit}(i,j,k)&=\ceil{\frac{\bit+\sum(\bit,i)+1-k}{2}}=\ceil{\frac{\bit+2^{i-1}-2+1+k}{2}}\\
						&\geq\ceil{\frac{\bit+1-k}{2}}\geq\floor{\frac{\bit+1-k}{2}}.
				\end{align*}
				Since $\bit$ is even, $t_\bit=-1$ is not a feasible parameter for $\bit$.
				This implies that $\occrec^{\canstrat}(e)=\floor{(\bit+1-k)/2}$ and in particular $\occrec^{\canstrat}(e)\neq\ell^{\bit}(i,j,k)+1$.
				Thus, \Pref{OR1}$_{i,j,k}$ implies $\canstrat(d_{i,j,k})\neq F_{i,j}$ for both $k\in\{0,1\}$.
				Since $\bit$ is even, $\occrec^{\canstrat}(d_{i,j,0},F_{i,j})=\occrec^{\canstrat}(d_{i,j,1},F_{i,j})$.
				Hence, as discussed previously, the switch $(d_{i,j,0},F_{i,j})$ is applied during $\canstrat\to\sigma$.
				Thus, \Cref{equation: Occrec is high} follows from \begin{align*}
					\occrec^{\sigma_{\bit+1}}(d_{i,j,0},F_{i,j})&=\floor{\frac{\bit+1}{2}}+1=\floor{\frac{\bit+1+1-0}{2}}\quad\text{ and}\\
					\occrec^{\sigma_{\bit+1}}(d_{i,j,1},F_{i,j})&=\floor{\frac{\bit}{2}}=\floor{\frac{\bit+1}{2}}=\floor{\frac{\bit+1+1-1}{2}}.
				\end{align*}
			\item Let $\bit_i=0$ and $(\bit+1)_i=1\wedge j=(\bit+1)_{i+1}$, implying $i=\ell(\bit+1)$. 
				Hence, as the occurrence record of $e$ with respect to $\sigma$ is described by \Cref{table: Occurrence Records}, \begin{align*}
				\occrec^{\sigma}(e)&=\ceil{\frac{\lastflip{\bit+1}{i}{\{(i+1,j)\}}+1-k}{2}}\\
					&=\ceil{\frac{\bit+1+1-k}{2}}=\floor{\frac{\bit+1-k}{2}}+1.
				\end{align*}
				Since $\id_{j=0}\lastflip{\bit+1}{i+1}{}+\id_{j=1}\lastunflip{\bit+1}{i+1}{}<\bit+1-2^i$, this implies $\ell^{\bit+1}(i,j,k)>\occrec^{\sigma}(e)+2$ and consequently \Cref{equation: Occrec not ell-1}.
			\item Let $\bit_i=0$ and $(\bit+1)_{i}=1\wedge j\neq(\bit+1)_{i+1}$.
				Then, $i=\ell(\bit+1)$.
				Since $j=(\bit+2)_{i+1}$ by assumption, the bit with index $i+1$ has switched when transitioning from $\bit+1$ to $\bit+2$.
				This is however only possible if $i=\ell(\bit+1)=1$.
				As this also implies $\bit_{i+1}=(\bit+1)_{i+1}\neq j$, this implies \[
					\ell^{\bit}(i,j,k)=\ceil{\frac{\bit-2^{i-1}+\sum(\bit,i)+1-k}{2}}=\ceil{\frac{\bit-1+1-k}{2}}=\floor{\frac{\bit+1-k}{2}}\]
				\Pref{OR3}$_{i,j,k}$ applied to $\canstrat$ thus implies $\occrec^{\canstrat}(e)=\floor{(\bit+1-k)/2}$.
				By the same arguments used in the earlier cases, we devise  $(d_{i,j,0},F_{i,j})\in\applied{\canstrat}{\sigma}$.
				But, similar to the previous cases, this implies \Cref{equation: Occrec is high}.
			\item Let $\bit_i=0$ and $(\bit+1)_i=0$.
				Then, $i>1$ and $\bit_{i+1}=(\bit+1)_{i+1}=(\bit+2)_{i+1}$, hence $j=\bit_{i+1}$.
				Thus, by \Cref{lemma: Numerics Of Ell}, \begin{align*}
					\ell^{\bit}(i,j,k)&=\ceil{\frac{\bit+2^{i-1}+\sum(\bit,i)+1-k}{2}}\geq\ceil{\frac{\bit+2^{i-1}+1-k}{2}}\\
						&\geq\ceil{\frac{\bit+2+1-k}{2}}=\ceil{\frac{\bit+1-k}{2}}+1\geq\floor{\frac{\bit+1-k}{2}}+1.
				\end{align*}
				Therefore $\floor{(\bit+1-k)/2}\leq\ell^{\bit}(i,j,k)-1$, implying $\occrec^{\canstrat}(e)=\floor{(\bit+1-k)/2}$.
				This implies \Cref{equation: Occrec is high} by using the same arguments as in the last cases.
		\end{enumerate}

	\item Since $i\neq\ell(\bit+2)$, it is not possible that $(\bit+1)_{i+1}=0\wedge(\bit+2)_{i+1}=1$.
		It thus suffices to investigate the following cases.
		\begin{enumerate}
			\item Let $\bit_i=0\wedge(\bit+1)_{i}=1$, i.e., $i=\ell(\bit+1)=\nsb$.
				Then $\bit_{i+1}=(\bit+1)_{i+1}$ and $(\bit+1)_{i+1}=(\bit+2)_{i+1}$ if and only if $i\neq 1$.
				Consider the case $i\neq 1$ first.
				Then, $j\neq\bit_{i+1}$, hence $j=1-\bit_{i+1}$.
				Since $i=\ell(\bit+1)$, \Cref{lemma: Numerics Of Ell} then implies \begin{align*}
					\ell^{\bit}(i,j,k)&=\ceil{\frac{\bit-2^{i-1}+\sum(\bit,i)+1-k}{2}}=\ceil{\frac{\bit-2^{\nsb-1}+2^{\nsb-1}-1+1-k}{2}}\\
						&=\floor{\frac{\bit+1-k}{2}}.
				\end{align*}
				Now, $\occrec^{\canstrat}(e)=\floor{(\bit+1-k)/2}$ or $\occrec^{\canstrat}(e)=\ell^{\bit}(i,j,k)-1\neq\floor{(\bit+1-k)/2}$.
				Consider the first case.
				Then $\canstrat(d_{i,j,k})\neq F_{i,j}$ by \Pref{OR2}$_{i,j,k}$.
				Using the same arguments used when proving the first statement, this implies $\occrec^{\sigma}(d_{i,j,k},F_{i,j})=\floor{(\bit+1+1-k)/2}$.
				Consider the second case.
				By our previous calculation and by \Pref{OR3}$_{i,j,k}$, $\occrec^{\canstrat}(d_{i,j,k},F_{i,j})=\maxocc-1$ and $k=0$.
				Also, by \Pref{OR2}$_{i,j,k}$, $\canstrat(d_{i,j,0})\neq F_{i,j}$.
				Thus, by \Cref{corollary: Switches With Low OR In Phase One}, $e\in\applied{\canstrat}{\sigma}$.
				Hence $\occrec^{\sigma}(e)=\floor{(\bit+1)/2}=\floor{(\bit+2)/2}$ since $\bit$ is odd.
				Consequently, $\occrec^{\sigma}(e)=\floor{(\bit+1+1-0)/2}$, so \Cref{equation: Occrec is high} holds.
				
				Now consider the case $i=1$, i.e., $(\bit+1)_{i+1}\neq(\bit+2)_{i+1}$.
				Then $j=\bit_{i+1}$, hence \[\ell^{\bit}(i,j,k)=\ceil{\frac{\bit+2^{i-1}+\sum(\bit,i)+1-k}{2}}=\ceil{\frac{\bit+2-k}{2}}=\floor{\frac{\bit+1-k}{2}}+1.\]
				Thus, $\ell^{\bit}(i,j,k)-1=\floor{(\bit+1-k)/2}$, implying $\occrec^{\canstrat}(e)=\floor{(\bit+1-k)/2}$.
				Since $F_{1,j}$ is the cycle center that is closed during the transition from $\canstrat$ to $\sigma_{\bit+1}$, the switch $(d_{1,j,k},F_{1,j})$ is applied for both $k$.
				Since $\bit$ is even, \Cref{equation: Occrec is high} follows from $\occrec^{\sigma}(e)=\floor{(\bit+1-k)/2}+1=\floor{(\bit+1+1-k)/2}.$
			\item Let $\bit_i=0\wedge(\bit+1)_{i}=0$.
				Since $i\neq\ell(\bit+2)$, we have $(\bit+2)_{i}=0$.
				This implies $\bit_{i+1}=(\bit+1)_{i+1}=(\bit+2)_{i+1}$, so $j=1-\bit_{i+1}$ and $i>2$.
				Thus, by \Cref{lemma: Numerics Of Ell}, \[\ell^{\bit}(i,j,k)=\ceil{\frac{\bit-2^{i-1}+\sum(\bit,i)+1-k}{2}}.\]
				Since $(\bit+1)_i=0$ implies $i\neq\ell(\bit+1)$, we either have $\occrec^{\canstrat}(d_{i,j,k},F_{i,j})\neq\ell^{\bit}(i,j,k)-1$ or $\occrec^{\canstrat}(d_{i,j,k},F_{i,j})=\floor{(\bit+1-k)/2}$ by \Pref{OR3}$_{i,j,k}$.
				Assume $\occrec^{\canstrat}(d_{i,j,k},F_{i,j})=\floor{(\bit+1-k)/2}$ and let $\ell^{\bit}(i,j,k)+1=\floor{(\bit+1-k)/2}$ first.
				Then, since $\ell^{\bit}(i,j,k)+1=\ell^{\bit+1}(i,j,k)$ by \Cref{lemma: Numerics Of OR}, we obtain \[\occrec^{\sigma}(e)\geq\occrec^{\canstrat}(d_{i,j,k},F_{i,j})=\ell^{\bit}(i,j,k)+1=\ell^{\bit+1}(i,j,k).\]
				Hence $\occrec^{\sigma}(d_{i,j,k},F_{i,j})\neq\ell^{\bit+1}(i,j,k)-1$.
				Now let $\occrec^{\canstrat}(d_{i,j,k},F_{i,j})\neq\ell^{\bit}(i,j,k)+1$.
				Then, by \Pref{OR2}$_{i,j,k}$, $\canstrat(d_{i,j,k})\neq F_{i,j}$.
				Since $\occrec^{\canstrat}(e)=\floor{(\bit+1-k)/2}$, we can apply the same arguments used when discussing previous cases to obtain \Cref{equation: Occrec is high}.
				
				Let $\occrec^{\canstrat}(d_{i,j,k},F_{i,j})\neq\ell^{\bit}(i,j,k)-1\neq\floor{(\bit+1-k)/2}$ as we can apply the same arguments used before otherwise.
				Thus, either $\occrec^{\canstrat}(d_{i,j,k},F_{i,j})=\ell^{\bit}(i,j,k)$ or $\occrec^{\canstrat}(d_{i,j,k},F_{i,j})=\ell^{\bit}(i,j,k)+1$.
				Since $\ell^{\bit+1}(i,j,k)=\ell^{\bit}(i,j,k)+1$, the statement follows directly if we consider the case $\occrec^{\canstrat}(d_{i,j,k},F_{i,j})=\ell^{\bit}(i,j,k)+1$.
				Hence assume $\occrec^{\canstrat}(d_{i,j,k},F_{i,j})=\ell^{\bit}(i,j,k)$.
				Then, \Pref{OR2}$_{i,j,k}$ implies $\canstrat(d_{i,j,k})\neq F_{i,j}$.
				As $\occrec^{\canstrat}(d_{i,j,k},F_{i,j})<\floor{(\bit+1-k)/2}$, $e$ is applied when transitioning from $\canstrat$ to $\sigma$.
				Thus $\occrec^{\sigma}(e)=\occrec^{\canstrat}(e)+1=\ell^{\bit+1}(i,j,k)$, implying \Cref{equation: Occrec not ell-1}.
			\item Let $\bit_i=1\wedge j=1-\bit_{i+1}$.
				Then, \Cref{lemma: Numerics Of Ell} implies\[\ell^{\bit}(i,j,k)=\ceil{\frac{\bit+\sum(\bit,i)+1-k}{2}}\geq\floor{\frac{\bit+1-k}{2}}.\]
				Since  $i\neq\ell(\bit+1)$ as $\bit_i=1$, this yields $\occrec^{\canstrat}(e)=\floor{(\bit+1-k)/2}$ by \Pref{OR3}$_{i,j,k}$.
				In particular, $\canstrat(d_{i,j,k})\neq F_{i,j}$, so the same arguments used previously yield \Cref{equation: Occrec is high}.
			\item Let $\bit_{i}=1\wedge j=\bit_{i+1}$, implying $\occrec^{\canstrat}(e)=\ceil{(\lastflip{\bit}{i}{\{(i+1,j)\}}+1-k)/2}$.
				Since $j\neq(\bit+2)_{i+1}$, bit $i+1$ switched.
				As $i\neq\ell(\bit+2)$ and $\bit_i=1$ yields $i\neq 1$, this implies $i\leq\nsb-1$.
				In particular, $\occrec^{\canstrat}(e)=\ceil{(\bit-2^{i-1}+2-k)/2}$.
				Furthermore, we then have $(\bit+1)_{i+1}=0\wedge(\bit+1)_{i+1}\neq j$.
				Hence, by \Cref{lemma: Numerics Of Ell}, \begin{align*}
					\ell^{\bit+1}(i,j,k)&=\ceil{\frac{\bit+1-2^{i-1}+\sum(\bit+1)+1-k}{2}}\\
						&=\ceil{\frac{\bit-2^{i-1}+2-k}{2}}=\occrec^{\canstrat}(e).
				\end{align*}
				Thus, $\occrec^{\sigma}(e)\geq\occrec^{\canstrat}(e)=\ell^{\bit+1}(i,j,k)$, so $\occrec^{\sigma}(e)\neq\ell^{\bit+1}(i,j,k)-1$.
		\end{enumerate}

	\item Since $i=\ell(\bit+2)$, we have $(\bit+1)_{i}=0$.
		This further implies $(\bit+1)_{i+1}=(\bit+2)_{i+1}$, so $j\neq (\bit+1)_{i+1}$.
		Thus, by \Cref{lemma: Numerics Of Ell}, \begin{align*}
			\ell^{\bit+1}(i,j,k)&=\ceil{\frac{\bit+1-2^{i-1}+\sum(\bit+1,i)+1-k}{2}}\\
				&=\ceil{\frac{\bit+1-2^{\nsb-1}+2^{\nsb-1}-1+1-k}{2}}\\
					&=\ceil{\frac{\bit+1-k}{2}}=\floor{\frac{\bit+1+1-k}{2}}.
		\end{align*}
		Since $\bit+1$ is even, the parameter $t_{\bit+1}=-1$ is not feasible.
		This implies \Cref{equation: Occrec is high}.
		
	\item By the choice of $i$, there is no number $\bit'\leq\bit+1$ with $i=\ell(\bit')$.
		Consequently, it holds that $\lastflip{\bit+1}{i+1}{}=\lastunflip{\bit+1}{i+1}{}=0$.
		Thus, by \Cref{lemma: Numerics Of Ell}, $\ell^{\bit+1}(i,j,k)\geq\bit+1>\floor{(\bit+1+1-k)/2}$, implying \Cref{equation: Occrec is high}.
		
	\item Since $k=1$, it suffices to show $\occrec^{\sigma}(d_{i,j,1},F_{i,j})=\maxocc$.
		Since $\bit$ is even, we have $\ell(\bit+1)=1$, hence $\bit_i=0$.
		As shown in the proof of \Cref{claim: OccRec for B Plus Two}, this implies $\ell^{\bit}(i,j,k)=\floor{(\bit-1)/2}=\maxocc-1$.
		Consequently, by \Cref{claim: OccRec for B Plus Two}, it holds that $\occrec^{\canstrat}(d_{i,j,1},F_{i,j})=\ell^{\bit}(i,j,k)<\occrec^{\canstrat}(d_{i,j,0},F_{i,j})$.
		Since \Pref{OR2}$_{i,j,1}$ now implies $\canstrat(d_{i,j,1})\neq F_{i,j}$, this implies $(d_{i,j,1},F_{i,j})\in\applied{\canstrat}{\sigma}$.
		But then, the statement follows since we then have $\occrec^{\sigma}(d_{i,j,1},F_{i,j})=\maxocc$.
\end{enumerate}
\end{proof}

\OROfNotAppliedEdges*

\begin{proof}
Let $i\in[n], j,k\in\{0,1\}$.
We distinguish four cases.
\begin{enumerate}
	\item Let $\bit_i=1\wedge\bit_{i+1}=j$.
		Since $\nsb=1$, this implies $(\bit+1)_{i}=1\wedge(\bit+1)_{i+1}=j.$
		Hence, $F_{i,j}$ is closed for both $\canstrat$ and $\sigma$ and the switch was not applied during $\canstrat\to\sigma$, implying $i\neq 1$.
		Consequently, $\lastflip{\bit}{i}{\{(i+1,j)\}}=\lastflip{\bit+1}{i}{\{(i+1,j)\}}$, implying the statement. 
	\item Let $\bit_i=0\wedge\bit_{i+1}\neq j$.					
		Consider the case $(\bit+1)_i=0\wedge(\bit+1)_{i+1}\neq j$, implying $i\neq 1$.
		Let $\id_{j=0}\lastflip{\bit}{i+1}{}+\id_{j=1}\lastunflip{\bit}{i+1}{}=0$, implying $\ell^{\bit}(i,j,k)\geq\bit$ by \Cref{lemma: Numerics Of Ell}.
		Since $\occrec^{\canstrat}(d_{i,j,k},F_{i,j})\leq\maxocc$, this implies $\occrec^{\sigma}(e)=\occrec^{\canstrat}(e)=\maxocc$ independent of $k$ since $\bit+1$ is odd.
		Note that this implies $\canstrat(d_{i,j,*})\neq F_{i,j}$ by \Pref{OR1}$_{i,j,*}$.
		Consequently, both $(d_{i,j,0},F_{i,j})$ and $(d_{i,j,1},F_{i,j})$ could have been applied during phase~1.
		However, due to the tie-breaking rule, only $(d_{i,j,0},F_{i,j})$ was applied during phase 1.
		It thus suffices to investigate $e\coloneqq(d_{i,j,1},F_{i,j})$.
		Since $\ell^{\bit+1}(i,j,1)\geq\bit$ by \Cref{lemma: Numerics Of Ell} and $\occrec^{\sigma}(e)=\occrec^{\canstrat}(e)$, we thus obtain \[
			\occrec^{\sigma}(e)=\floor{\frac{\bit+1}{2}}=\floor{\frac{(\bit+1)+1-k}{2}}=\min\left (\floor{\frac{(\bit+1)+1-k}{2}},\ell^{\bit+1}(i,j,k)\right ),
		\]
		hence choosing $t_{\bit+1}=0$ yields the correct description of the occurrence record.
		
		Let $\id_{j=0}\lastflip{\bit}{i+1}{}+\id_{j=1}\lastunflip{\bit}{i+1}{}\neq 0$.
		Then, by \Cref{lemma: Numerics Of Ell} and since $\bit_1=0$, \begin{align}\label{equation: Edges Not Switched NSBOne}
			\ell^{\bit}(i,j,k)&\leq\ceil{\frac{\bit-2^{i-1}+2^{i-1}-1-1+1-k}{2}}=\maxocc-k.
		\end{align}
		If $\canstrat(d_{i,j,k})\neq F_{i,j}$ and $\occrec^{\canstrat}(e)<\maxocc$, then $(d_{i,j,k},F_{i,j})$ was applied in phase 1 by \Cref{corollary: Switches With Low OR In Phase One}. 
		By \Pref{OR1}$_{i,j,k}$, $\canstrat(d_{i,j,k})=F_{i,j}$ and $\occrec^{\canstrat}(e)=\maxocc$ is not possible.
		Consider the case $\canstrat(d_{i,j,k})=F_{i,j}$ and $\occrec^{\canstrat}(d_{i,j,k},F_{i,j})<\maxocc$.
		We show that~$e$ was then applied during phase~5.
		By \Cref{corollary: Cycle Edges Phase Five If NSB One}, we need to show $i>u$ and $i<m$.
		If $i\neq u$, the first statement follows as $\bit_i=0$ and $(\bit+1)_i=0$.
		If $i=u$, Inequality~(\ref{equation: Edges Not Switched NSBOne}) is tight, contradicting $\canstrat(d_{i,j,k})=F_{i,j}$ since we then had $\occrec^{\canstrat}(d_{i,j,k},F_{i,j})=\ell^{\bit}(i,j,k)+1\geq\maxocc$.
		Assume, for the sake of contradiction, $i>m$.
		Then $\bit<2^{i-1}$, hence $\lastflip{\bit}{i+1}{}=\lastunflip{\bit}{i+1}{}=0$.
		Consequently, by \Cref{lemma: Numerics Of Ell}, $\ell^{\bit}(i,j,k)\geq\bit$, contradicting Properties (\ref{property: OR1})$_{i,j,k}$ and (\ref{property: OR2})$_{i,j,k}$. 
		Therefore, $e$ was applied during phase 5 and we do not consider it here.
		It thus suffices to consider the case $\canstrat(d_{i,j,k})\neq F_{i,j}$ and $\occrec^{\canstrat}(e)=\maxocc$.
		It then holds that $\occrec^{\canstrat}(e)=\ell^{\bit}(i,j,k)+t_{\bit}$ for some feasible $t_{\bit}$ due to Inequality~(\ref{equation: Edges Not Switched NSBOne}).
		This implies $\occrec^{\sigma}(e)=\ell^{\bit}(i,j,1)+1$ if $k=1$, contradicting $\canstrat(d_{i,j,1})\neq F_{i,j}$.
		Thus $k=0$ and $\occrec^{\canstrat}(e)=\ell^{\bit}(i,j,0)$.
		But this implies that Inequality~(\ref{equation: Edges Not Switched NSBOne}) is an equality.
		Consequently, $i=\ell(\bit+2)$.
		Also, $\id_{j=0}\lastflip{\bit}{i+1}{}+\id_{j=1}\lastunflip{\bit}{i+1}{}\neq 0$ by assumption, hence $\bit+2$ is not a power of two.
		Hence, by \Pref{OR3}$_{i,j,0}$, the parameter $t_{\bit+1}=-1$ is feasible.
		Since $\ell^{\bit+1}(i,j,0)=\ell^{\bit}(i,j,0)+1$ by \Cref{lemma: Progress for unimportant CC} and $\maxocc=\floor{(\bit+1+1-0)/2}-1$, this parameter describes the occurrence record with respect to $\sigma$.
		Hence, \[\occrec^{\sigma}(e)=\min\left (\floor{\frac{\bit+1+1-0}{2}}, \ell^{\bit+1}(i,j,0)+t_{\bit+1}\right)<\floor{\frac{\bit+1+1-0}{2}}\] for $t_{\bit+1}=-1$, so the occurrence record is correctly described by \Cref{table: Occurrence Records}.
		This concludes the case $(\bit+1)_i=0\wedge(\bit+1)_{i+1}\neq j$.
		
		Consider the case $(\bit+1)_{i}=1$ and $(\bit+1)_{i+1}\neq j$ next, implying $i=\nsb=1$.
		If $\id_{j=0}\lastflip{\bit}{2}{}+\id_{j=1}\lastunflip{\bit}{2}{}=0$, we can use the same arguments used for the case $(\bit+1)_{i}=0\wedge(\bit+1)_{i+1}\neq j$.
		Hence let $\id_{j=0}\lastflip{\bit}{2}{}+\id_{j=1}\lastunflip{\bit}{2}{}\neq 0$.
		Then by \Cref{lemma: Numerics Of Ell}, $\ell^{\bit}(1,j,k)=\maxocc$ for both choices of $k\in\{0,1\}$.
		Since the parameter $t_{\bit}=-1$ is not feasible as $\bit$ is even and choosing $t_{\bit}=1$ violates \Cref{lemma: Occurrence Records Cycle Vertices}, it thus holds that $\occrec^{\sigma}(e)=\ell^{\bit}(i,j,k)=\floor{(\bit+1)/2}$ for both $k\in\{0,1\}$.
		In particular, $\canstrat(d_{i,j,k})\neq F_{i,j}$ for both $k\in\{0,1\}$.
		Hence, by the tie-breaking rule, $(d_{i,j,0},F_{i,j})$ is applied during phase 1.
		Consequently, $(d_{i,j,1},F_{i,j})$ is not applied during phase 1 and the same arguments used previously can be used to show that choosing $t_{\bit+1}=0$ is feasible and implies the desired characterization.
		This concludes the case $(\bit+1)_i=1\wedge(\bit+1)_{i+1}\neq j$.
		Since only the first bit switches during $\canstrat\to\sigma$, this also concludes the case $\bit_i=0\wedge\bit_{i+1}\neq j$.
		
	\item Let $\bit_i=0\wedge\bit_{i+1}=j$.
		Since only the first bit switches, it suffices to consider $i\neq 1$ and $(\bit+1)_i=0\wedge(\bit+1)_{i+1}=j$.
		As before, if $\id_{j=0}\lastflip{\bit}{i+1}{}+\id_{j=1}\lastunflip{\bit}{i+1}{}=0$, then the statement follows directly.
		Hence assume $\id_{j=0}\lastflip{\bit}{i+1}{}+\id_{j=1}\lastunflip{\bit}{i+1}{}\neq 0.$
		Then, $\ell^{\bit}(i,j,k)\geq\ceil{(\bit+2+1-k)/2}\geq\maxocc+1$ by \Cref{lemma: Numerics Of Ell}.
		Since the parameter $t_{\bit}=-1$ is not feasible, this implies $\occrec^{\canstrat}(d_{i,j,*},F_{i,j})=\maxocc$ and $\canstrat(d_{i,j,*})\neq F_{i,j}$.
		By the same arguments used before,  $(d_{i,j,1},F_{i,j})$ is not applied and its occurrence record with respect to $\sigma$ is described by \Cref{table: Occurrence Records} when interpreted for $\bit+1$.
		
	\item Finally, consider the case $\bit_i=1\wedge\bit_{i+1}\neq j$.
		Since only the first bit switches, this implies $i\neq 1$ and $(\bit+1)_i=1$ and $(\bit+1)_{i+1}\neq j$.
		It is easy to see that this enables us to use the same arguments used previously. \qedhere
\end{enumerate}
\end{proof}

\InequalityisMetOrOneEdgeIsSwitched*

\begin{proof}
By Equation~(\ref{equation: Minima are the same}), at most one of the two edges of the cycle center $F_{i,j}$ is switched.
We distinguish the following cases.
\begin{enumerate}
	\item Let $F_{i,j}$ be open for $\sigma_{\bit-1}$.
		Then, one of the two cycle edges is applied during phase 1 of $\sigma_{\bit-1}\to\canstrat$ since no cycle center is open at the end of phase 1 by \Cref{corollary: No Open CC For Phase Two} resp. \ref{corollary: No Open CC In Phase Three}.
	\item Let $F_{i,j}$ be closed for $\sigma_{\bit-1}$.
		Then,  since $\bit_i=0$, either $(\bit-1)_i=1\wedge(\bit-1)_{i+1}=j$ or $(\bit-1)_i=0\wedge(\bit-1)_{i+1}\neq j$.
		Consider the first case.
		This case can only happen if $i<\ell(\bit)$, additionally implying $j\neq\bit_{i+1}$.
		In addition, we then have \begin{align*}
			\occrec^{\sigma_{\bit-1}}(d_{i,j,k},F_{i,j})&=\floor{\frac{\lastflip{\bit-1}{i}{\{(i+1,j)\}}-k}{2}}+1\\
				&=\floor{\frac{\bit-1-2^{i-1}+1-k}{2}}+1\\
				&=\floor{\frac{\bit-2^{i-1}-k}{2}}+1=\floor{\frac{\bit-2^{i-1}+2-k}{2}}.
		\end{align*} for $k\in\{0,1\}$.
		For $i=2$, this implies \begin{align*}
			\occrec^{\sigma_{\bit-1}}(d_{i,j,0},F_{i,j})&=\floor{\frac{\bit-2+2-0}{2}}=\floor{\frac{\bit}{2}}\\
			\occrec^{\sigma_{\bit-1}}(d_{i,j,1},F_{i,j})&=\floor{\frac{\bit-2+2-1}{2}}=\floor{\frac{\bit-1}{2}}\geq\floor{\frac{\bit}{2}}-1
		\end{align*} and thus the statement.
		For $i\geq 3$ it is easy to verify that this implies that the occurrence record of at least one of the cycle edges is so low that the corresponding edge is applied as improving switch during phase 5 of $\sigma_{\bit-1}\to\canstrat$.
		This concludes the first case, hence assume $(\bit-1)_i=0\wedge(\bit-1)_{i+1}\neq j.$
		Then, $F_{i,j}$ being closed implies $\occrec^{\sigma_{\bit-1}}(d_{i,j,k},F_{i,j})=\ell^{\bit-1}(i,j,k)+1\leq\floor{\bit/2}-1$ and $\sigma_{\bit-1}(d_{i,j,k})=F_{i,j}$ for both $k\in\{0,1\}$.
		If this inequality is met with equality for both $k$, then the statement follows as the occurrence record of the edges is sufficiently high.
		If the inequality is strict for at least one $k$, then the corresponding switch is applied during phase 5 of the transition $\sigma_{\bit-1}\to\canstrat$.
	\item Let $F_{i,j}$ be halfopen for $\sigma_{\bit-1}$.
		Then, for some $k\in\{0,1\}$, $\sigma_{\bit-1}(d_{i,j,k})=F_{i,j}$ as well as $\occrec^{\sigma_{\bit-1}}(d_{i,j,k},F_{i,j})=\ell^{\bit-1}(i,j,k)+1\leq\floor{\bit/2}-1$.
		Furthermore, $\sigma_{\bit-1}(d_{i,j,1-k})=F_{i,j}$ and $\occrec^{\sigma_{\bit-1}}(d_{i,j,k},F_{i,j})\in\{\floor{\bit/2}-1, \floor{\bit/2}\}$.
		If $\occrec^{\sigma_{\bit-1}}(d_{i,j,k},F_{i,j})=\floor{\bit/2}-1$, then the edge is applied as an improving switch and the statement follows.
		Hence assume $\occrec^{\sigma_{\bit-1}}(d_{i,j,k},F_{i,j})=\floor{\bit/2}$.
		Then, due to $\sigma_{\bit-1}(d_{i,j,k})\neq F_{i,j}$, this implies $\ell^{\bit-1}(i,j,1-k)+1\neq\floor{\bit/2}$.
		However, $\ell^{\bit-1}(i,j,k)\leq\floor{\bit/2}-2$ since it then holds that $\occrec^{\sigma_{\bit-1}}(d_{i,j,k},F_{i,j})=\ell^{\bit-1}(i,j,k)+1\leq\floor{\bit/2}-1$.
		But, since $\ell^{\bit-1}(i,j,k)$ and $\ell^{\bit-1}(i,j,1-k)$ differ by at most one, this implies $\ell^{\bit-1}(i,j,1-k)\leq\floor{\bit/2}-1$.
		But this is a contradiction to $\occrec^{\sigma_{\bit-1}}(d_{i,j,1-k},F_{i,j})=\floor{\bit/2}$. \qedhere
\end{enumerate}
\end{proof}

\NeedToHaveZero*

\begin{proof}
Let, for the sake of contradiction, $(\bit-1)_i=1$.
Then, as $\bit_i=0$, we have $i<\ell(\bit)$ and consequently $(\bit-1)_{i+1}\neq\bit_{i+1}$.
It further implies that $\bit$ is even.
Then, since $e\in\applied{\sigma_{\bit-1}}{\canstrat}$ by assumption, the switch $e=(g_i,F_{i,j})$ was applied during phase 5 of $\sigma_{\bit-1}\to\canstrat$.
This implies $j=0$ if $G_n=S_n$ resp. $j=\bit_{i+1}=1-(\bit-1)_{i+1}$ if $G_n=M_n$.
Consider the case $G_n=M_n$ first.
Then, since $(\bit-1)_i=1\wedge j=1-(\bit-1)_{i+1}$ imply $\ell^{\bit-1}(i,j,k)\geq \floor{(\bit-k)/2}+1$ for $k\in\{0,1\}$ by \Cref{lemma: Numerics Of Ell}, we obtain  $\occrec^{\sigma_{\bit-1}}(d_{i,j,k},F_{i,j})=\floor{(\bit-k)/2}$. 

In addition, since $\occrec^{\sigma_{\bit-1}}(d_{i,j,k},F_{i,j})\neq \ell^{\bit-1}(i,j,k)+1$ for both $k\in\{0,1\}$, $F_{i,j}$ is then open with respect to $\sigma_{\bit-1}$ by \Pref{OR2}$_{i,j,*}$.
This implies that $(d_{i,j,1},F_{i,1})$ is applied during phase~1 of $\sigma_{\bit-1}\to\canstrat$.
But then $\min_{k\in\{0,1\}}\occrec^{\sigma_{\bit-1}}(d_{i,j,k},F_{i,j})<\min_{k\in\{0,1\}}\occrec^{\canstrat}(d_{i,j,k},F_{i,j}),$ contradicting Equation~(\ref{equation: Minima are the same}).
Now consider the case $G_n=S_n$.
If $j=0=1-(\bit-1)_{i+1}$, then the statement follows by the same arguments.
This is the case if and only if $i<\ell(\bit)-1$, so let $i=\ell(\bit)-1$.
By \Cref{definition: Canonical Strategy}, this implies $\sigma_{\bit-1}(g_i)=F_{i,0}$.
Since $F_{i,0}$ is then closed during phase 1 of the transition $\sigma_{\bit-1}\to\canstrat$ and since $(g_i,F_{i,1})$ cannot be applied during phase 5 in $S_n$, this is a contradiction.
\end{proof}

\SelectorSwitchAppliedInPhaseOne*

\begin{proof}
By \Cref{corollary: Selection Vertices In Phase One}, there is some index $k\in\{0,1\}$ such that $\sigma_{\bit-1}(d_{i,j,k})=F_{i,j}$ as well as $\occrec^{\sigma_{\bit-1}}(d_{i,j,k},F_{i,j})=\ell^{\bit-1}(i,j,k)+1\leq\floor{\bit/2}-1$.
Moreover, it holds that $\occrec^{\sigma_{\bit-1}}(d_{i,j,1-k},F_{i,j})=\floor{\bit/2}-1$ and $(d_{i,j,1-k},F_{i,j})$ is applied during phase 1 of $\sigma_{\bit-1}\to\canstrat$.
By \Cref{equation: Minima are the same}, $(d_{i,j,k},F_{i,j})$ is not applied as improving switch during $\sigma_{\bit-1}\to\canstrat$.
It is easy to verify that this implies that we need to have $\occrec^{\sigma_{\bit-1}}(d_{i,j,k},F_{i,j})=\floor{\bit/2}-1$ as well as $\ell(\bit)>1$.
This implies that $\bit$ is even and, since $(\bit-1)_i=0$ and that we cannot have $i=2$.
In particular, the first statement holds and we have $\ell^{\bit-1}(i,j,k)=\floor{\bit/2}-2$.
Since $\ell^{\bit-1}(i,j,1)\leq\ell^{\bit-1}(i,j,0)$, this implies  $k=1$ as $\occrec^{\sigma_{\bit-1}}(d_{i,j,1-k},F_{i,j})\leq\floor{\bit/2}-2$ held otherwise. 
This implies $\floor{(\bit-2)/2}=\floor{(\bit-2^{i-1}+\sum(\bit-1,i)+1)/2}$ as \begin{align*}
\floor{\frac{\bit-2}{2}}&=\floor{\frac{\bit}{2}}-1=\ell^{\bit-1}(i,j,1)+1=\ceil{\frac{\bit-1-2^{i-1}+\sum(\bit-1,i)}{2}}+1\\
	&=\floor{\frac{\bit-2^{i-1}+\sum(\bit-1,i)+2}{2}}=\floor{\frac{\bit-2^{i-1}+\sum(\bit-1,i)+1}{2}}.
\end{align*}
where the last equality follows since $i\geq 3$ and since $\sum(\bit-1,i)$ is odd as $\bit$ is even.
Since $\bit$ is even and $i\neq 1$, the nominators on both sides are then even.
But this implies that the nominators have to be equal.
Since $\sum(\bit,i)=\sum(\bit-1,i)+1$, this implies $\sum(\bit,i)=2^{i-1}-2$.
More precisely,
\begin{align*}
	\sum(\bit,i)&=\sum(\bit-1,i)+1=\bit-2-\bit+2^{i-1}-1+1=2^{i-1}-2,
\end{align*} implying the second statement.

We now prove that it is not possible that $(g_i,F_{i,j})$ was applied during both $\sigma_{\bit-2}\to\sigma_{\bit-1}$ and $\sigma_{\bit-3}\to\sigma_{\bit-2}$, implying (c).
First note that $i\geq 3$ implies $\bit-3\geq\bar{\bit}$, i.e., $\bit-3$ is indeed  \enquote{contributing} to $\mathfrak{N}(\bar{\bit},\bit-1)$.
Since the statement follows if $(g_i,F_{i,j})\notin\applied{\sigma_{\bit-2}}{\sigma_{\bit-1}}$, assume that this was the case.
Since we assume that $(g_i,F_{i,j})$ is applied during phase 1 of $\sigma_{\bit-1}\to\canstrat$ and since $i\neq\ell(\bit)$ and $i\neq\ell(\bit-1)$, it is not possible that $(g_i,F_{i,j})$ was applied during phase 5 of $\sigma_{\bit-2}\to\sigma_{\bit-1}$.
It was thus applied during phase 1 of $\sigma_{\bit-2}\to\sigma_{\bit-1}$.
Since~$\bit$ is even, this implies $\occrec^{\sigma_{\bit-2}}(d_{i,j,k},F_{i,j})\leq\floor{(\bit-2+1)/2}-1=\floor{\bit/2}-2$ for both $k\in\{0,1\}$.
Now, for the sake of contradiction, assume $(g_i,F_{i,j})\in\applied{\sigma_{\bit-3}}{\sigma_{\bit-2}}$.
By the same argument used previously, it is not possible that this switch was applied during phase 5 of that transition.
It thus needs to be applied during phase 1.
This is an immediate contradiction if $(\bit-3)_i=1$.
Thus assume $(\bit-3)_i=0$.
Then, since $(g_i,F_{i,j})$ was applied during phase 1 of $\sigma_{\bit-3}\to\sigma_{\bit-2}$, there is some $k\in\{0,1\}$ such that $\sigma_{\bit-3}(d_{i,j,1-k})\neq F_{i,j}$ as well as  $\occrec^{\sigma_{\bit-3}}(d_{i,j,1-k},F_{i,j})=\floor{(\bit-3+1)/2}-1=\floor{\bit/2}-2.$
Since this switch is then applied during phase 1, this implies $\occrec^{\sigma_{\bit-2}}(d_{i,j,1-k},F_{i,j})=\floor{\bit/2}-1$ which is a contradiction.
\end{proof}

\SelectorSwitchAppliedInPhaseFive*

\begin{proof}
We remind here that we have $\bit_i=0$ and $(\bit-1)_i=0$, implying $\bit_{i+1}=(\bit-1)_{i+1}$.
By the conditions describing under which circumstances an improving switch $(g_i,F_{i,j})$ becomes improving in phase 5 (see \Cref{lemma: Phase Five Escape Easy}) and the assumption, we have $j=0$ if $G_n=S_n$ resp. $j=\indbit^{\sigma}_{i+1}=\bit_{i+1}=(\bit-1)_{i+1}$ if $G_n=M_n$.
By \Cref{definition: Canonical Strategy} resp. \ref{definition: Canonical Strategy MDP}, $\sigma_{\bit-1}(g_i)=F_{i,1-j}$ then implies that $F_{i,1-j}$ has to be closed with respect to $\sigma_{\bit-1}$.
As $(\bit-1)_i=0$, it also implies that $1-j=1-(\bit-1)_{i+1}$ needs to hold for both $M_n$ and $S_n$
But this implies that \[\occrec^{\sigma_{\bit-1}}(d_{i,1-j,k},F_{i,1-j})=\ell^{\bit-1}(i,1-j,k)+1\leq\floor{\frac{\bit}{2}}-1\] for both $k\in\{0,1\}$.
We show this implies $i\neq 2$ as assuming $\occrec^{\sigma_{\bit-1}}(d_{2,1-j,0},F_{2,1-j})=\ell^{\bit-1}(2,1-j,0)+1$ contradicts $\occrec^{\sigma_{\bit-1}}(d_{2,1-j,0},F_{2,1-j})\leq\floor{\bit/2}-1$.

Assume $i=2$.
Then, since $(\bit-1)_2=\bit_2=0$, we need to have $\bit_1=1$, so $\bit$ is odd.
Hence, $\ell^{\bit-1}(2,1-j,0)+1=\floor{(\bit-2+1)/2}=\floor{\bit/2}.$
Also, $\floor{(\bit-2)/2}=\floor{\bit/2}-1$ due to the parity of $\bit$ and thus $\floor{(\bit-2)/2}\geq \floor{(\bit-2^{i-1}+\sum(\bit-1,i)+1)/2}$, which contradicts the previously given inequality. 
Using the same arguments used when proving \Cref{claim: Selector Switch Applied in Phase One}, this implies $\sum(\bit,i)\leq 2^{i-1}-2$.

It remains to prove that $(g_i,F_{i,j})\in\applied{\sigma_{\bit-2}}{\sigma_{\bit-1}}$ implies $(g_i,F_{i,j})\notin\applied{\sigma_{\bit-3}}{\sigma_{\bit-2}}$.
As we have $\sigmabar_{\bit-1}(g_i)=1-j$, the switch $(g_i,F_{i,j})$ cannot have been applied during phase 5 of $\sigma_{\bit-2}\to\sigma_{\bit-1}$.
It was thus applied during phase 1 of $\sigma_{\bit-2}\to\sigma_{\bit-1}$.
Assume $\bit-2=\bar{\bit}$ which can happen if $\bit$ is odd.
But then, $(g_i,F_{i,j})$ was not applied during phase 1 of $\sigma_{\bit-2}\to\sigma_{\bit-1}$ as we then have $(\bit-2)_i=1$.
Thus assume $\bar{\bit}\leq\bit-3$ and that $(g_i,F_{i,j})$ was applied during phase 1 of $\sigma_{\bit-2}\to\sigma_{\bit-1}$.
Towards a contradiction, assume that $(g_i,F_{i,j})$ is applied during $\sigma_{\bit-3}\to\sigma_{\bit-2}$.
Since we apply the same switch during phase~1 of $\sigma_{\bit-2}\to\sigma_{\bit-1}$ and $i\neq\ell(\bit-1)$ due to $(\bit-1)_i=0$, the switch must have been applied during phase~1 of $\sigma_{\bit-3}\to\sigma_{\bit-2}$.
If this was not the case, i.e., if it was applied during phase 5, we had $\sigma_{\bit-2}(g_i)=F_{i,j}$ and it would not be possible to apply $(g_i,F_{i,j})$.
Note that this implies $(\bit-3)_i=0$ and in particular $(\bit-3)_{i+1}=j$.		
Then, since $(g_i,F_{i,j})$ is applied during both $\sigma_{\bit-3}\to\sigma_{\bit-2}$ and $\sigma_{\bit-2}\to\sigma_{\bit-1}$, the improving switch $(g_i,F_{i,1-j})$ has to be applied in between.
This switch can only be applied during phase 1 of $\sigma_{\bit-2}\to\sigma_{\bit-1}$ since $j=(\bit-3)_{i+1}=\bit_{i+1}$ and since $(g_i,F_{i,j})$ was applied during phase 5 of $\sigma_{\bit-1}\to\canstrat$.
But this is a contradiction as we apply $(g_i,F_{i,j})$ during phase 1 of that transition and $i\neq\ell(\bit-2)$.
\end{proof}

%
%
%
%
%
%
%
%
%

\end{document}